# SPECIAL SET LINEAR ALGEBRA AND SPECIAL SET FUZZY LINEAR ALGEBRA

**W. B. Vasantha Kandasamy**
e-mail: **vasanthakandasamy@gmail.com**
web: **http://mat.iitm.ac.in/home/wbv/public_html/**
**www.vasantha.in**

**Florentin Smarandache**
e-mail: **smarand@unm.edu**

**K Ilanthenral**
e-mail: **ilanthenral@gmail.com**

**2009**

# SPECIAL SET LINEAR ALGEBRA AND SPECIAL SET FUZZY LINEAR ALGEBRA

W. B. Vasantha Kandasamy
Florentin Smarandache
K Ilanthenral

**2009**

# CONTENTS

# DEDICATION

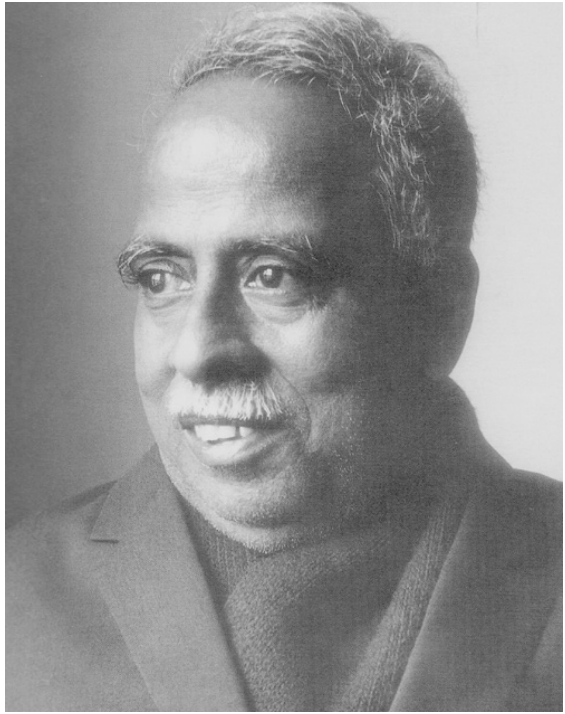

(15-09-1909 to 03-02-1969)

*We dedicate this book to late Thiru C.N.Annadurai (Former Chief Minister of Tamil Nadu) for his Centenary Celebrations. He is fondly remembered for legalizing self respect marriage, enforcing two language policy and renaming Madras State as TamilNadu. Above all he is known for having dedicated his rule to Thanthai Periyar.*

# PREFACE

This book for the first time introduces the notion of special set linear algebra and special set fuzzy linear algebra. This is an extension of the book set linear algebra and set fuzzy linear algebra. These algebraic structures basically exploit only the set theoretic property, hence in applications one can include a finite number of elements without affecting the systems property. These new structures are not only the most generalized structures but they can perform multi task simultaneously; hence they would be of immense use to computer scientists.

    This book has five chapters. In chapter one the basic concepts about set linear algebra is given in order to make this book a self contained one. The notion of special set linear algebra and their fuzzy analogue is introduced in chapter two. In chapter three the notion of special set semigroup linear algebra is introduced. The concept of special set n- vector spaces, n greater than or equal to three is defined and their fuzzy analogue is given in chapter four. The probable applications are also mentioned.

The final chapter suggests 66 problems. Our thanks are due to Dr. K. Kandasamy for proof-reading this book. We also acknowledge our gratitude to Kama and Meena for their help with corrections and layout.

W.B.VASANTHA KANDASAMY
FLORENTIN SMARANDACHE
K.ILANTHENRAL

**Chapter One**

# BASIC CONCEPTS

In this chapter we just introduce the notion of set linear algebra and fuzzy set linear algebra. This is mainly introduced to make this book a self contained one. For more refer [60].

**DEFINITION 1.1:** *Let S be a set. V another set. We say V is a set vector space over the set S if for all $v \in V$ and for all $s \in S$; vs and $sv \in V$.*

*Example 1.1:* Let $V = \{1, 2, \ldots, \infty\}$ be the set of positive integers. $S = \{2, 4, 6, \ldots, \infty\}$ the set of positive even integers. V is a set vector space over S. This is clear for $sv = vs \in V$ for all $s \in S$ and $v \in V$.

It is interesting to note that any two sets in general may not be a set vector space over the other. Further even if V is a set vector space over S then S in general need not be a set vector space over V.

For from the above example 1.1 we see V is a set vector space over S but S is also a set vector space over V for we see for every $s \in S$ and $v \in V$, $vs = sv \in S$. Hence the above example is both important and interesting as one set V is a set vector space another set S and vice versa also hold good inspite of the fact $S \neq V$.

Now we illustrate the situation when the set V is a set vector space over the set S. We see V is a set vector space over the set S and S is not a set vector space over V.

***Example 1.2:*** Let V = {$Q^+$ the set of all positive rationals} and S = {2, 4, 6, 8, … , $\infty$}, the set of all even integers. It is easily verified that V is a set vector space over S but S is not a set vector space over V, for $\frac{7}{3} \in V$ and $2 \in S$ but $\frac{7}{3}.2 \notin S$. Hence the claim.

Now we give some more examples so that the reader becomes familiar with these concepts.

***Example 1.3:*** Let $M_{2\times 2} = \left\{ \begin{pmatrix} a & b \\ c & d \end{pmatrix} \;\middle|\; a, b, c, d \in \right.$ (set of all positive integers together with zero)} be the set of 2 × 2 matrices with entries from N. Take S = {0, 2, 4, 6, 8, … , $\infty$}. Clearly $M_{2\times 2}$ is a set vector space over the set S.

***Example 1.4:*** Let V = {$Z^+ \times Z^+ \times Z^+$} such that $Z^+$ = {set of positive integers}. S = {2, 4, 6, …, $\infty$}. Clearly V is a set vector space over S.

***Example 1.5:*** Let V = {$Z^+ \times Z^+ \times Z^+$} such that $Z^+$ is the set of positive integers. S = {3, 6, 9, 12, … , $\infty$}. V is a set vector space over S.

***Example 1.6:*** Let $Z^+$ be the set of positive integers. $pZ^+ = S$, p any prime. $Z^+$ is a set vector space over S.

**Note:** Even if p is replaced by any positive integer n still $Z^+$ is a set vector space over $nZ^+$. Further $nZ^+$ is also a set vector space

over $Z^+$. This is a collection of special set vector spaces over the set.

***Example 1.7:*** Let Q[x] be the set of all polynomials with coefficients from Q, the field of rationals. Let S = {0, 2, 4, …, ∞}. Q[x] is a set vector space over S. Further S is not a set vector space over Q[x].

Thus we see all set vector spaces V over the set S need be such that S is a set vector space over V.

***Example 1.8:*** Let R be the set of reals. R is a set vector space over the set S where S = {0, 1, 2, … , ∞}. Clearly S is not a set vector space over R.

***Example 1.9:*** $\mathbb{C}$ be the collection of all complex numbers. Let $Z^+ \cup \{0\}$ = {0, 1, 2, … , ∞} = S. $\mathbb{C}$ is a set vector space over S but S is not a set vector space over $\mathbb{C}$.

At this point we propose a problem.
 Characterize all set vector spaces V over the set S which are such that S is a set vector space over the set V.

Clearly $Z^+ \cup \{0\}$ = V is a set vector space over S = p $Z^+ \cup \{0\}$, p any positive integer; need not necessarily be prime. Further S is also a set vector space over V.

***Example 1.10:*** Let V = $Z_{12}$ = {0, 1, 2, … , 11} be the set of integers modulo 12. Take S = {0, 2, 4, 6, 8, 10}. V is a set vector space over S. For s ∈ S and v ∈ V, sv ≡ vs (mod 12).

***Example 1.11:*** Let V = $Z_p$ = {0, 1, 2, … , p – 1} be the set of integers modulo p (p a prime). S = {1, 2, … , p – 1} be the set. V is a set vector space over S.

In fact V in example 1.11 is a set vector space over any proper subset of S also.

This is not always true, yet if V is a set vector space over S. V need not in general be a set vector space over every proper subset of S. Infact in case of $V = Z_p$ every proper subset $S_i$ of S = {1, 2, …, p – 1} is such that V is a set vector space over $S_i$.

**Note:** It is important to note that we do not have the notion of linear combination of set vectors but only the concept of linear magnification or linear shrinking or linear annulling. i.e., if $v \in V$ and s is an element of the set S where V is the set vector space defined over it and if sv makes v larger we say it is linear magnification.

***Example 1.12:*** If $V = Z^+ \cup \{0\}$ = {0, 1, 2, …, ∞} and S = {0, 2, 4, …, ∞} we say for s = 10 and v = 21, sv = 10.21 is a linear magnification.

Now for v = 5 and s = 0, sv = 0.5 = 0 is a linear annulling. Here we do not have the notion of linear shrinking.

So in case of any modeling where the researcher needs only linear magnification or annulling of the data he can use the set vector space V over the set S.

*Notation:* We call the elements of the set vector space V over the set S as set vectors and the scalars in S as set scalars.

***Example 1.13:*** Suppose $Q^+ \cup \{0\}$ = {all positive rationals} = $Q^* = V$ and
$$S = \left\{0, 1, \frac{1}{2}, \frac{1}{2^2}, \frac{1}{2^3}, \ldots\right\}.$$

Then V is a set vector space over S. We see for $6 \in Q^* = V$ and $\frac{1}{2} \in S$ $\frac{1}{2}.6 = 3$. This is an instance of linear shrinking. Suppose s = 0 and $v = \frac{21}{3}$ then sv = 0. $\frac{21}{3} = 0$ is an instance of linear annulling.

Now we see this vector space is such that there is no method to get a linearly magnified element or the possibility of linear magnification.

Now we have got a peculiar instance for $1 \in S$ and $1.v = v$ for every $v \in V$ linear neutrality or linearly neutral element. The set may or may not contain the linearly neutral element. We illustrate yet by an example in which the set vector space which has linear magnifying, linear shrinking, linear annulling and linearly neutral elements.

***Example 1.14:*** Let $Q^+ \cup \{0\} = V = \{$set of all positive rationals with zero$\}$.

$$S = \left\{0, 1, \frac{1}{3}, \frac{1}{3^2}, ..., 2, 4, 6, ..., \infty\right\}$$

be the set. V is a set vector space over the set S. We see V is such that there are elements in S which linearly magnify some elements in V; for instance if $v = 32$ and $s = 10$ then $sv = 10.32 = 320$ is an instance of linear magnification. Consider $v = 30$ and $s = \frac{1}{3}$ then $sv = \frac{1}{3} \times 30 = 10 \in V$. This is an instance of linear shrinking. Thus we have certain elements in V which are linearly shrunk by elements of the set S.

Now we have $0 \in S$ such that $0\,v = 0$ for all $v \in V$ which is an instance of linear annulling. Finally we see $1 \in S$ is such that $1.v = v$ for all $v \in V$; which is an instance of linearly neutral. Thus we see in this set vector space V over the set S given in example 1.14, all the four properties holds good. Thus if a researcher needs all the properties to hold in the model he can take them without any hesitation. Now we define yet another notion called linearly normalizing element of the set vector space V. Suppose $v \in V$ is a set vector and $s \in S$ is a set scalar and $1 \in V$ which is such that $1.v = v.1 = v$ and $s.1 = 1.s = s$ for all $v \in V$ and for all $s \in S$; we call the element scalar s in S to be a linearly normalizing element of S if we have v ∈ V such that sv = vs = 1.

For example in the example 1.14, we have for 3 ∈ V and $\frac{1}{3}$ ∈ S; $\frac{1}{3}$ 3 = 1 ∈ S. Thus $\frac{1}{3}$ is a linearly normalizing element of S.

It is important to note that as in case of linearly annulling or linearly neutral element the scalar need not linearly normalize every element of the set vector space V. In most cases an element can linearly normalize only one element.

Having seen all these notions we now proceed on to define the new notion of set vector subspace of a set vector space V.

**DEFINITION 1.2:** *Let V be a set vector space over the set S. Let W ⊂ V be a proper subset of V. If W is also a set vector space over S then we call W a set vector subspace of V over the set S.*

We illustrate this by a few examples so that the reader becomes familiar with this concept.

*Example 1.15:* Let V = {0, 1, 2, …, ∞} be the set of positive integers together with zero. Let S = {0, 2, 4, 6, …, ∞}, the set of positive even integers with zero. V is a set vector space over S. Take W = {0, 3, 6, 9, 12, …, ∞} set of all multiples of 3 with zero. W ⊆ V; W is also a set vector space over S. Thus W is a set vector subspace of V over S.

*Example 1.16:* Let Q[x] be the set of all polynomials with coefficients from Q; the set of rationals. Q[x] is a set vector space over the set S = {0, 2, 4, … , ∞}. Take W = {0, 1, …, ∞} the set of positive integers with zero. W is a set vector space over the set S. Now W ⊆ Q ⊆ Q[x]; so W is a set vector subspace of V over the set S.

***Example 1.17:*** Let $V = Z^+ \times Z^+ \times Z^+ \times Z^+$ a set of vector space over the set $Z^+ = \{0, 1, 2, \ldots, \infty\}$. Let $W = Z^+ \times Z^+ \times \{0\} \times \{0\}$, a proper subset of V. W is also a set vector space over $Z^+$, i.e., W is a set vector subspace of V.

***Example 1.18:*** Let $V = 2Z^+ \times 3Z^+ \times Z^+$ be a set, V is a set vector space over the set $S = \{0, 2, 4, \ldots, \infty\}$. Now take $W = 2Z^+ \times \{0\} \times 2Z^+ \subseteq V$; W is a set vector subspace of V over the set S.

***Example 1.19:*** Let

$$V = M_{3\times 2} = \left\{ \begin{pmatrix} a & d \\ b & e \\ c & f \end{pmatrix} \;\middle|\; a, b, c, d, e, f \in Z^+ \cup \{0\} \right\}$$

be the set of all $3 \times 2$ matrices with entries from the set of positive integers together with zero. V is a set vector space over the set $S = Z^+ \cup \{0\}$. Take

$$W = \left\{ \begin{pmatrix} a & 0 \\ b & 0 \\ c & 0 \end{pmatrix} \;\middle|\; a, b, c \in Z^+ \cup \{0\} \right\} \subseteq V;$$

W is a set vector subspace of V over the set S.

Now having defined set vector subspaces, we proceed on to define the notion of zero space of a set vector space V over the set S. We as in the case of usual vector spaces cannot define set zero vector space at all times. The set zero vector space of a set vector space V exists if and only if the set vector space V over S has $\{0\}$ in V i.e., $\{0\}$ is the linearly annulling element of S or $0 \in V$ and $0 \notin S$ in either of the two cases we have the set zero subspace of the set vector space V.

***Example 1.20:*** Let $Z^+ = V = \{1, 2, \ldots, \infty\}$ be a set vector space over the set $S = \{2, 4, 6, \ldots, \infty\}$; V is a set vector space but $0 \notin V$, so V does not have a set vector zero subspace.

It is interesting to mention here that we can always adjoin the zero element to the set vector space V over the set S and this does not destroy the existing structure. Thus the element $\{0\}$ can always be added to make the set vector space V to contain a set zero subspace of V.

We leave it for the reader to prove the following theorem.

**THEOREM 1.1:** *Let V be a set vector space over the set S. Let $W_1, \ldots, W_n$ be n proper set vector subspaces of V over S. Then $\bigcap_{i=1}^{n} W_i$ is a set vector subspace of V over S. Further $\bigcap_{i=1}^{n} W_i = \phi$ can also occur, if even for a pair of set vector subspaces $W_i$ and $W_j$ of V we have $W_i \cap W_j = \phi\, (i \neq j)$.*

*Clearly even if $0 \in V$ then also we cannot say $\bigcap_{i=1}^{n} W_i = \{0\}$ as 0 need not be present in every set vector subspace of V.*

We illustrate the situation by the following example.

***Example 1.21:*** Let $Z^+ = \{1, 2, \ldots, \infty\}$ be a set vector space over the set $S = \{2, 4, 6, \ldots, \infty\}$. Take

$$W_1 = \{2, \ldots, \infty\},$$
$$W_2 = \{3, 6, \ldots, \infty\}, \ldots,$$

and

$$W_p = \{p, 2p, 3p, \ldots, \infty\}.$$

We see $W_i \cap W_j \neq \phi$ for every $i, j\, (i \neq j)$.

Will $\bigcap W_i = \phi$ if $i = 1, 2, \ldots, \infty$ ?

Will $\bigcap_{i=1}^{n} W_i \neq \phi$ if $i = 1, 2, \ldots, n$; $n < \infty$?

**Example 1.22:** Let $V = \{1, 2, \ldots, \infty\}$ be a set vector space over $S = \{2, , \ldots, \infty\}$. $W = \{2, 2^2, \ldots, 2^n \ldots\}$ is set vector subspace of V over the set S.

$W_1 = \{3, 3^2, 3^3, \ldots, \infty\}$ is a proper subset of V but $W_1$ is not a set vector subspace of V as $W_1$ is not a set vector space over S as $2.3 = 6 \notin W_1$ for $2 \in S$ and $3 \in W_1$. Thus we by this example show that in general every proper subset of the set V need not be a set vector space over the set S.

Now having seen the set vector subspaces of a set vector space we now proceed to define yet another new notion about set vector spaces.

**DEFINITION 1.3:** *Let V be a set vector space over the set S. Let T be a proper subset of S and W a proper subset of V. If W is a set vector space over T then we call W to be a subset vector subspace of V over T.*

We first illustrate this situation by examples before we proceed on to give more properties about them.

**Example 1.23:** Let $V = Z^+ \cup \{0\} = \{0, 1, 2, \ldots, \infty\}$ be a set vector space over the set $S = \{2, 4, 6, \ldots, \infty\}$. Consider $W = \{3, 6, 9, \ldots, \} \in V$.
Take $T = \{6, 12, 18, 24, \ldots\} \subseteq S$. Clearly W is a set vector space over T; so W is a subset vector subspace over T.

**Example 1.24:** Let $V = Z^+ \times Z^+ \times Z^+$ be a set vector space over $Z^+ = S$. Take $W = Z^+ \times \{0\} \times \{0\}$ a proper subset of V and $T = \{2, , 6, \ldots, \infty\} \subseteq Z^+ = S$. Clearly W is a set vector space over T i.e., W is a subset vector subspace over T.

Now we show all proper subsets of a set vector space need not be a subset vector subspace over every proper subset of S. We illustrate this situation by the following examples.

*Example 1.25:* Let $V = Z^+ \times Z^+ \times Z^+$ be a set vector space over $Z^+ = S = \{0, 1, 2, ..., \infty\}$. Let $W = \{3, 3^2, ..., \infty\} \times \{5, 5^2, ...\} \times \{0\} \subseteq V$. Take $T = \{2, 4, 6, ...\} \subseteq S$.

Clearly W is not a set vector space over T. That is W is not a subset vector subspace of V over T. Thus we see every subset of a set vector space need not in general be a subset vector subspace of the set vector space V over any proper subset T of the set S.

We illustrate this concept with some more examples.

*Example 1.26:* Let $V = 2Z^+ \times 3Z^+ \times 5Z^+ = \{(2n, 3m, 5t) \mid n, m, t \in Z^+\}$; V is a set vector space over the set $S = Z^+ = \{1, 2, ..., \infty\}$. Take $W = \{2, 2^2, ..., \infty\} \times \{3, 3^2, ...\} \times \{5, 5^2, ...\} \subseteq V$. W is not a subset vector subspace over the subset $T = \{2, 4, ..., \infty\} \subseteq Z^+$. Further W is not even a set vector subspace over $Z^+ = S$. Take $W = \{2, 2^2, ...\} \times \{0\} \times \{0\}$ a proper subset of V. Choose $T = \{2, 2^3, 2^5, 2^7\}$. W is a subset vector subspace of V over the subset $T \subseteq S$.

*Example 1.27:* Let

$$V = \left\{ \begin{pmatrix} a & b \\ c & d \end{pmatrix} \middle| a, b, c, d \in Z^+ \cup \{0\} = \{0, 1, 2, ...\} \right\}$$

be the set of all $2 \times 2$ matrices with entries from $Z^+ \cup \{0\}$. V is a set vector space over the set $S = Z^+ = \{1, 2, ..., \infty\}$. Take

$$W = \left\{ \begin{pmatrix} x & y \\ z & w \end{pmatrix} \middle| x, y, z, w \in 2Z^+ = \{2, 4, 6, ...\} \right\} \subseteq V.$$

W is a subset vector subspace over the subset $T = \{2, 4, ..., \infty\}$.

*Example 1.28:* Let

$$V = \left\{ \begin{pmatrix} a & b & c & d \\ e & f & g & h \end{pmatrix} \middle| \ a, b, c, d, e, f, g \text{ and } h \in Z \right\}$$

be the set of all $2 \times 4$ matrices with entries from the set of integers. V is a set vector space over the set $S = Z^+ = \{1, 2, \ldots, \infty\}$. Take

$$W = \left\{ \begin{pmatrix} a & b & c & d \\ e & f & g & h \end{pmatrix} \middle| \ a, b, c, d, e, f, g \text{ and } h \in Z^+ \right\} \subseteq V.$$

Clearly W is a subset vector subspace over $T = \{2, 4, \ldots, \infty\} \subseteq S = Z^+$.

Now as in case of vector spaces we can have the following theorem.

**THEOREM 1.2:** *Let V be a set vector space over the set S. Suppose $W_1, \ldots, W_n$ be a set of n subset vector subspaces of V over the same subset T of S. Then either $\bigcap_{i=1}^{n} W_i = \phi$ or $\bigcap_{i=1}^{n} W_i$ is a subset vector subspace over T or if each $W_i$ contains 0 then $\bigcap_{i=1}^{n} W_i = \{0\}$, the subset vector zero subspace of V over T.*

The proof is left as an exercise for the reader as the proof involves only simple set theoretic techniques.

*Note:* We cannot say anything when the subset vector subspaces of V are defined over different subsets of S. We illustrate this situation by some more examples.

*Example 1.29:* Let $V = Z^+ \times Z^+ \times Z^+$ be a set vector space over a set $S = Z^+ = \{1, 2, \ldots, \infty\}$. $W = \{2, 4, 6, 8, \ldots\} \times \{\phi\} \times \{\phi\}$ is a subset vector subspace of V over the subset $T = \{2, 4, 8, 16, \ldots\}$. $W_1 = \{3, 6, 9, \ldots\} \times \phi \times \phi \subset V$, is a subset subvector space over the subset $T_1 = \{3, 3^2, \ldots\}$. Clearly we cannot define $W \cap$ $W_1$ for we do not have even a common subset over which it can be defined as $T \cap T_1 = \phi$.

From this example a very natural question is if $T \cap T_1 \neq \phi$ and if $W_1 \cap W$ is not empty can we define some new structure. For this we make the following definition.

**DEFINITION 1.4:** *Let V be a set vector space over the set S. Suppose $W_i$ is a subset vector subspace defined over the subset $T_i$ of S for i = 1, 2, ..., n; n < $\infty$ and if $W = \cap W_i \neq \phi$ and $T = \cap T_i \neq \phi$; then we call W to be a sectional subset vector sectional subspace of V over T.*

*Note:* We call it a sectional for every subset vector subspace contributes to it.

We give illustration of the same.

***Example 1.30:*** Let $V = Z^+ \times Z^+ \times Z^+ \times Z^+$ be a set vector space over the set $Z^+ = \{1, 2, ..., \infty\}$. Let $W_1 = \{2, 4, 6, ...\} \times \{2, 4, 6, ...\} \times Z^+ \times Z^+$ be a subset vector subspace over $T_1 = \{2, 4, ..., \infty\}$. $W_2 = \{Z^+\} \times \{2, 4, 6, ...\} \times Z^+ \times Z^+$ be a subset vector subspace over $T_2 = \{2, 2^2, ..., \infty\}$. Let $W_3 = \{2, 2^2, ..., \infty\} \times \{Z^+\} \times \{2, 4, 6, ...\} \times \{2, 2^2, ..., \infty\}$ be a subset vector subspace over $T_3 = \{2, 2^3, 2^5, ..., \infty\}$. Consider $W_1 \cap W_2 \cap W_3 = \{2, 2^2, ..., \infty\} \times \{2, 4, 6, ...\} \times \{2, 4, 6, ..., \infty\} \times \{2, 2^2, ..., \infty\} = W$. Now $T = T_1 \cap T_2 \cap T_3 = \{2, 2^3, 2^5, ...\}$; W is a sectional subset vector sectional subspace of V over T.

We see a sectional subset vector sectional subspace is a subset vector subspace but a sectional subset vector sectional subspace in general is not a subset vector subspace.

We prove the following interesting theorem.

**THEOREM 1.3:** *Every sectional subset vector sectional subspace W of set vector space V over the set S is a subset vector subspace of a subset of S but not conversely.*

*Proof:* Let V be the set vector space over the set S. $W = \bigcap_{i=1}^{n} W_i$ be a sectional subset vector sectional subspace over $T = \bigcap_{i=1}^{n} T_i$ where each $W_i$ is a subset vector subspace over $T_i$ for $i = 1, 2, \ldots, n$ with $W_i \neq W_j$ and $T_i \neq T_j$; for $i \neq j$, $1 \leq i, j \leq n$. We see W is a sectional subset vector subspace over $T_i$ ; $i = 1, 2, \ldots, n$.

We illustrate the converse by an example.

*Example 1.31:* Let $V = \{2, 2^2, \ldots\}$ be a set vector space over the set $S = \{2, 2^2\}$. $W_1 = \{2^2, 2^4, \ldots\}$ is a subset vector subspace of V over the set $T_1 = \{2\}$ and $W_2 = \{2^3, 2^6, \ldots, \infty\}$ is a subset vector subspace of over the set $T_2 = \{2^2\}$. Now $W_1 \cap W_2 \neq \phi$ but $T_1 \cap T_2 = \phi$. We do not and cannot make $W_1$ or $W_2$ as sectional subset vector sectional subspace.

That is why is general for every set vector space V over the set S we cannot say for every subset vector subspace $W(\subset V)$ over the subset $T(\subset S)$ we can find atleast a subset vector subspace $W_1(\subset V)$ over the subset $T_1(\subset S)$ such that $W_1 \cap W_2 \neq \phi$ and $T \cap T_1 \neq \phi$.

Now before we define the notion of basis, dimension of a set S we mention certain important facts about the set vector spaces.

1. All vector spaces are set vector spaces but not conversely.

The converse is proved by giving counter examples.
Take $V = Z^+ = \{0, 1, 2, \ldots, \infty\}$, V is a set vector space over the set $S = \{2, 4, 6, \ldots, \infty\}$. We see V is not an abelian group under addition and S is not a field so V can never be a vector space over S. But if we have V to be a vector space over the field F we see V is a set vector space over the set F as for every $c \in F$ and $v \in V$ we have $cv \in V$. Thus every vector space is a set vector space and not conversely.

2. All semivector spaces over the semifield F is a set vector space but not conversely.

We see if V is a semivector space over the semifield F then V is a set and F is a set and for every $v \in V$ and $a \in F$; $av \in V$ hence V is trivially a set vector space over the set F.

However take $\{-1, 0, 1, 2, \ldots, \infty\} = V$ and $S = \{0, 1\}$ V is a set vector space over S but V is not even closed with respect to '+' so V is not a group. Further the set $S = \{0, 1\}$ is not a semifield so V is not a semivector space over S. Hence the claim.

Thus we see the class of all set vector spaces contains both the collection of all vector spaces and the collection of all semivector spaces. Thus set vector spaces happen to be the most generalized concept.

Now we proceed on to define the notion of generating set of a set vector space V over the set S.

**DEFINITION 1.5:** *Let V be a set vector space over the set S. Let $B \subseteq V$ be a proper subset of V, we say B generates V if every element v of V can be got as sb for some $s \in S$ and $b \in B$. B is called the generating set of V over S.*

***Example 1.32:*** Let $Z^+ \cup \{0\} = V$ be a set vector space over the set $S = \{1, 2, 4, \ldots, \infty\}$. The generating set of V is $B = \{0, 1, 3, 5, 7, 9, 11, 13, \ldots, 2n + 1, \ldots\}$ B is unique. Clearly the cardinality of B is infinite.

***Examples 1.33:*** Let $V = Z^+ \cup \{0\}$ be a set vector space over the set $S = \{1, 3, 3^2, \ldots, \infty\}$. $B = \{0, 1, 2, 4, 5, 6, 7, 8, 10, 11, 12, 13, 14, 15, 16, 17, 18, 19, 20, 21, 22, 23, 24, 25, 26, 28, \ldots, 80, 82, \ldots, 242, 244, \ldots, 3^6 - 1, \ldots, 3^6 + 1, \ldots, \infty\} \subseteq V$ is the generating set of V.

***Example 1.34:*** Let $V = \{0, 3, 3^2, \ldots, 3^n, \ldots\}$ be a set vector space over $S = \{0, 1, 3\}$. $B = \{0, 3\}$ is a generating set for $3.0 = 0$, $3.1 = 3$, $3.3 = 3^2$, $3.3^2 = 3^3$ so on.

*Example 1.35:* Let

$$V = \left\{ \begin{pmatrix} a & b \\ c & d \end{pmatrix} \middle| \, a, b, c, d \in Z^+ = \{1, 2, \ldots, \infty\} \right\}$$

be the set vector space over the set $S = Z^+$. The generating set of V is infinite.

Thus we see unlike in vector space or semivector spaces finding the generating set of a set vector space V is very difficult.
When the generating set is finite for a set vector space V we say the set vector space is finite set hence finite cardinality or finite dimension, otherwise infinite or infinite dimension.

*Example 1.36:* Let $Z^+ = V = \{1, 2, \ldots, \infty\}$ be a set vector space over the set $S = \{1\}$. The dimension of V is infinite.

*Example 1.37:* Let $V = \{1, 2, \ldots, \infty\}$ be a set vector space over the set $S = \{1, 2, \ldots, \infty\} = V$. Then dimension of V is 1 and B is uniquely generated by $\{1\}$. No other element generates V.

Thus we see in case of set vector space we may have only one generating set. It is still an open problem to study does every set vector space have one and only one generating subset?

*Example 1.38:* Let $V = \{2, 4, 6, \ldots\}$ be a set vector space over the set $S = \{1, 2, 3, \ldots\}$. $B = \{2\}$ is the generating set of V and dimension of V is one.

Thus we have an important property enjoyed by set vector spaces. We have a vector space V over a field F of dimension one only if V = F, but we see in case of set vector space V, dimension V is one even if V ≠ S. The example 1.38 gives a set vector space of dimension one where V ≠ S.

*Example 1.39:* Let $V = \{2, 4, 6, \ldots, \infty\}$ be a set vector space over the set $S = \{2, 2^2, \ldots, \infty\}$. $B = \{2, 6, 10, 12, 14, 18, 20, 22, 24, \ldots, 30, 34, \ldots\}$ is a generating set of V.

Now in case of examples 1.38 and 1.39 we see V = {2, 4, …, ∞} but only the set over which they are defined are different so as in case of vector spaces whose dimension is dependent on the field over which it is defined are different so as in case of vector spaces whose dimension is dependent on the field over which it is defined so also the cardinality of the generating set of a vector space V depends on the set over which V is defined. This is clear from examples 1.38 and 1.39.

Now having defined cardinality of set vector spaces we define the notion of linearly dependent and linearly independent set of a set vector space V over the set S.

**DEFINITION 1.6:** *Let V be a set vector space over the set S. B a proper subset of V is said to be a linearly independent set if x, y ∈ B then x ≠ sy or y ≠ s'x for any s and s' in S. If the set B is not linearly independent then we say B is a linearly dependent set.*

We now illustrate the situation by the following examples.

*Example 1.40:* Let $Z^+$ = V = {1, 2, ..., ∞} be a set vector space over the set S = {2, 4, 6, …, ∞}. Take B = {2, 6, 12} ⊂ V; B is a linearly dependent subset for 12 = 6.2, for 6 ∈ S. B = {1, 3} is a linearly independent subset of V.

*Example 1.41:* Let V = {0, 1, 2, …, ∞} be a set vector space over the set S = {3, $3^2$, …}. B = {1, 2, 4, 8, 16, …} is a linearly independent subset of V. As in case of vector spaces we can in case of set vector spaces also say a set B which is the largest linearly dependent subset of V? A linearly independent subset B of V which can generate V, then we say B is a set basis of V or the generating subset of V and cardinality of B gives the dimension of V.

**DEFINITION 1.7:** *Let V and W be two set vector spaces defined over the set S. A map T from V to W is said to be a set linear transformation if*
$$T(v) = w$$

*and*
$$T(sv) = sw = sT(v)$$

*for all $v \in V$, $s \in S$ and $w \in W$.*

***Example 1.42:*** *Let $V = Z^+ = \{0, 1, 2, \ldots, \infty\}$ and $W = \{0, 2, 4, \ldots, \infty\}$ be set vector spaces over the set $S = \{0, 2, 2^2, 2^3, \ldots\}$.*

$$T : V \to W$$
$$T(0) = 0,$$
$$T(1) = 2,$$
$$T(2) = 4,$$
$$T(3) = 6$$

and so on. $T(2^n p) = 2^n (2p)$ for all $p \in V$. Thus T is a set linear transformation from V to W. As in case of vector spaces we will not be always in a position to define the notion of null space of T in case of set vector spaces. For only if $0 \in V$ as well as W; we will be in a position to define null set of a set linear transformation T.

Now we proceed on to define the notion of set linear operator of a set vector space V over the set S.

**DEFINITION 1.8:** *Let V be a set vector space over the set S. A set linear transformation T from V to V is called the set linear operator of V.*

***Example 1.43:*** *Let $V = Z^+ = \{1, 2, \ldots, \infty\}$ be a set vector space over the set $S = \{1, 3, 3^2, \ldots, \infty\}$. Define T from V to V by $T(x) = 2x$ for every $x \in V$. T is a set linear operator on V.*

Now it is easy to state that if V is a set vector space over the set S and if $O_S(V)$ denotes the set of all set linear operators on V then $O_S(V)$ is also a set vector space over the same set S. Similarly if V and W are set vector spaces over the set S and $T_S(V, W)$ denotes the set of all linear transformation from V to W then $T_S(V,W)$ is also a set vector space over the set S. For if we want to prove some set V is a set vector space over a set S it is enough if we show for every $s \in S$ and $v \in V$; $sv \in V$. Now in case of $O_S(V)$ = {set of all set linear operators from V to V}. $O_S(V)$ is a set and clearly for every $s \in S$ and for every $T \in O_S(V)$, $sT \in O_s(V)$ i.e., sT is again a set linear operator of V. Hence $O_S(V)$ is a set vector space over the set S.

Likewise if we consider the set $T_S(V, W)$ = {set of all set linear transformations from V to W}, then the set $T_S(V, W)$ is again a set vector space over S. For if $s \in S$ and $T \in T_S(V, W)$ we see sT is also a set linear operator from V to W. Hence the claim. Now we can talk about the notion of invertible set linear transformation of the set vector spaces V and W defined over the set S.

Let T be a set linear transformation from V into W. We say T is set invertible if there exist a set linear transformation U from W into V such that UT and TU are set identity set maps on V and W respectively. If T is set invertible, the map U is called the set inverse of T and is unique and is denoted by $T^{-1}$.

Further more T is set invertible if and only if (1) T is a one to one set map that is $T\alpha = T\beta$ implies $\alpha = \beta$. (2) T is onto that is range of T is all of W.

We have the following interesting theorem.

**THEOREM 1.4:** *Let V and W be two set vector spaces over the set S and T be a set linear transformation from V into W. If T is invertible the inverse map $T^{-1}$ is a set linear transformation from W onto V.*

*Proof:* Given V and W are set vector spaces over the set S. T is a set linear transformation from V into W. When T is a one to one onto map, there is a uniquely determined set inverse map $T^{-1}$ which set maps W onto V such that $T^{-1}T$ is the identity map on V and $TT^{-1}$ is the identity function on W.
    Now what we want to prove is that if a set linear transformation T is set invertible then $T^{-1}$ is also set linear.

Let x be a set vector in W and c a set scalar from S. To show
$$T^{-1}(cx) = cT^{-1}(x).$$

Let $y = T^{-1}(x)$ that is y be the unique set vector in V such that Ty = x. Since T is set linear T(cy) = cTy. Then cy is the unique set vector in V which is set by T into cx and so

$$T^{-1}(cx) = cy = c\, T^{-1}(x)$$

and $T^{-1}$ is set linear !

Suppose that we have a set invertible set linear transformation T from V onto W and a set invertible set linear transformation U from W onto Z. Then UT is set invertible and $(UT)^{-1} = T^{-1}U^{-1}$.

To obtain this it is enough if we verify $T^{-1}U^{-1}$ is both left and right set inverse of UT.

Thus we can say in case of set linear transformation; T is one to one if and only if $T\alpha = T\beta$ if and only if $\alpha = \beta$.

Since in case of set linear transformation we will not always be in a position to have zero to be an element of V we cannot define nullity T or rank T. We can only say if $0 \in V$, V a set vector space over S and $0 \in W$, W also a set vector space over S then we can define the notion of rank T and nullity T, where T is a set linear transformation from V into W.

We may or may not be in a position to have results of linear transformation from vector spaces.

Also the method of representing every vector space V over a field F of dimension n as $V \cong \underbrace{F \times ... \times F}_{n-\text{times}}$ may not be feasible in case of set vector spaces. Further the concept of representation of set linear transformation as a matrix is also not possible for all set vector spaces. So at this state we make a note of the inability of this structure to be always represented in this nice form.

Next we proceed on to define the new notion of set linear functionals of a set vector space V over a set S.

**DEFINITION 1.9:** *Let V be a set vector space over the set S. A set linear transformation from V onto the set S is called a set linear functional on V, if*

$$f : V \to S$$
$$f(c\alpha) = cf(\alpha)$$

*for all $c \in S$ and $\alpha \in V$.*

***Example 1.44:*** Let $V = \{0, 1, 2, ..., \infty\}$ and $S = \{0, 2, 4, ..., \infty\}$. $f : V \to S$ defined by $f(0) = 0$, $f(1) = 2$, $f(2) = 4$, ..., is a set linear functional on V.

***Example 1.45:*** Let $V = Z^+ \times Z^+ \times Z^+$ be a set vector space over $Z^+$. A set linear functional $f : V \to Z^+$ defined by $f(x, y, z) = x + y + z$.

Now we proceed onto define the new notion of set dual space of a set vector space V.

**DEFINITION 1.10:** *Let V be a set vector space over the set S. Let L(V,S) denote the set of all linear functionals from V into S, then we call L(V,S) the set dual space of V. Infact L(V, S) is also a set vector space over S.*

The study of relation between set dimension of V and that of L(V, S) is an interesting problem.

***Example 1.46:*** Let $V = \{0, 1, 2, ..., \infty\}$ be a set vector space over the set $S = \{0, 1, 2, ..., \infty\}$. Clearly the set dimension of V over S is 1 and this has a unique generating set $B = \{1\}$. No other proper subset of B can ever generate V.

What is the dimension of L(V, S)?

We cannot define the notion of set hyperspace in case of set vector space using set linear functionals. We can define the concept of set annihilator if and only if the set vector space V and the set over which it is defined contains the zero element, otherwise we will not be in a position to define the set annihilator.

**DEFINITION 1.11:** *Let V be a set vector space over the set S, both V and S has zero in them. Let A be a proper subset of the set V, the set annihilator of A is the set $A^o$ of all set linear functionals f on V such that $f(\alpha) = 0$ for every $\alpha$ in A.*

**DEFINITION 1.12:** *A fuzzy vector space (V, $\eta$) or $\eta V$ is an ordinary vector space V with a map $\eta : V \to [0, 1]$ satisfying the following conditions;*

1. $\eta(a + b) \geq \min \{\eta(a), \eta(b)\}$
2. $\eta(-a) = \eta(a)$
3. $\eta(0) = 1$
4. $\eta(ra) \geq \eta(a)$

*for all a, b, $\in V$ and $r \in F$ where F is the field.*

We now define the notion of set fuzzy vector space or $V_\eta$ or $V\eta$ or $\eta V$.

**DEFINITION 1.13:** *Let V be a set vector space over the set S. We say V with the map $\eta$ is a fuzzy set vector space or set fuzzy vector space if $\eta: V \to [0, 1]$ and $\eta(ra) \geq \eta(a)$ for all $a \in V$ and $r \in S$. We call $V_\eta$ or $V\eta$ or $\eta V$ to be the fuzzy set vector space over the set S.*

We now illustrate this situation by the following example.

*Example 1.47:* Let V = {(1 3 5), (1 1 1), (5 5 5), (7 7 7), (3 3 3), (5 15 25), (1 2 3)} be set which is a set vector space over the set S = {0, 1}.
  Define a map $\eta: V \to [0, 1]$ by

$$\eta(x, y, z) = \left(\frac{x + y + z}{50}\right) \in [0, 1]$$

for $(x, y, z) \in V$. $V_\eta$ is a fuzzy set vector space.

*Example 1.48:* Let $V = Z^+$ the set of integers. $S = 2Z^+$ be the set. V is a set vector space over S. Define $\eta: V \to [0, 1]$ by, for every $v \in V$; $\eta(v) = \dfrac{1}{v}$. $\eta V$ is a set fuzzy vector space or fuzzy set vector space.

*Example 1.49:* Let $V = \{(a_{ij}) \mid a_{ij} \in Z^+; 1 \leq i, j \leq n\}$ be the set of all $n \times n$ matrices with entries from $Z^+$.

Take $S = 3Z^+$ to be the set. V is a fuzzy set vector space where $\eta: V \to [0, 1]$ is defined by

$$\eta(A = (a_{ij})) = \begin{cases} \dfrac{1}{5|A|} & \text{if } |A| \neq 0 \\ 1 & \text{if } |A| = 0. \end{cases}$$

$V\eta$ is the fuzzy set vector space.

The main advantage of defining set vector spaces and fuzzy set vector spaces is that we can include elements x in the set vector spaces V provided for all $s \in S$, $sx \in V$ this cannot be easily done in usual vector spaces. Thus we can work with the minimum number of elements as per our need and work with them by saving both time and money.

We give yet some more examples.

*Example 1.50:* Let $V = 2Z^+ \times 5Z^+ \times 7Z^+$ be a set vector space over the set $Z^+$; with $\eta: V \to [0, 1]$ defined by

$$\eta((x, y, z)) = \dfrac{1}{x + y + z}$$

makes, $\eta V$ a fuzzy set vector space.

Now we define the notion of set fuzzy linear algebra.

**DEFINITION 1.14:** *A set fuzzy linear algebra (or fuzzy set linear algebra) (V, η) or ηV is an ordinary set linear algebra V with a map such η: V → [0, 1] such that η(a + b) ≥ min (η(a), η(b)) for a, b ∈ V.*

*Since we know in the set vector space V we merely take V to be a set but in case of the set linear algebra V we assume V is closed with respect to some operation usually denoted as '+' so the additional condition η(a + b) ≥ min (η(a), η(b)) is essential for every a, b ∈ V.*

We illustrate this situation by the following examples.

*Example 1.51:* Let $V = Z^+[x]$ be a set linear algebra over the set $S = Z^+$; η: V → [0, 1].

$$\eta(p(x)) = \begin{cases} \dfrac{1}{\deg(p(x))} \\ 1 \quad \text{if } p(x) \text{ is a constant.} \end{cases}$$

Clearly Vη is a set fuzzy linear algebra.

*Example 1.52:* Let

$$V = \left\{ \begin{pmatrix} a & b \\ c & d \end{pmatrix} \middle| a, b, d, c \in Z^+ \right\}$$

be set linear algebra over $2Z^+ = S$. Define

$$\eta\left(\begin{pmatrix} a & b \\ c & d \end{pmatrix}\right) = \begin{cases} \dfrac{1}{|ad - bc|} & \text{if } ad \neq bc \\ 0 & \text{if } ad = bc \end{cases}$$

for every a, b, c, d ∈ $Z^+$. Clearly Vη is a fuzzy set linear algebra.

***Example 1.53:*** Let $V = Z^+$ be a set linear algebra over $Z^+$. Define $\eta : V \to [0, 1]$ as $\eta(a) = \dfrac{1}{a}$. $V\eta$ is a fuzzy set linear algebra.

Now we proceed onto define the notion of fuzzy set vector subspace and fuzzy set linear subalgebra.

**DEFINITION 1.15:** *Let V be a set vector space over the set S. Let $W \subset V$ be the set vector subspace of V defined over S. If $\eta: W \to [0, 1]$ then $W_\eta$ is called the fuzzy set vector subspace of V.*

We illustrate this by the following example.

***Example 1.54:*** Let $V = \{(1\ 1\ 1), (1\ 0\ 1), (0\ 1\ 1), (0\ 0\ 0), (1\ 0\ 0)\}$ be a set vector space defined over the set $\{0, 1\}$. Define $\eta : V \to [0, 1]$ by

$$\eta(x\ y\ z) = \dfrac{(x + y + z)}{9} \pmod{2}.$$

So that
$$\eta\,(0\ 0\ 0) = 0$$
$$\eta\,(1\ 1\ 1) = \dfrac{1}{9}$$
$$\eta\,(1\ 0\ 1) = 0$$
$$\eta\,(1\ 0\ 0) = \dfrac{1}{9}$$
$$\eta\,(0\ 1\ 1) = 0.$$

$V\eta$ is a set fuzzy vector space. Take $W = \{(1\ 1\ 1), (0\ 0\ 0), (0\ 1\ 1)\} \subset V$. W is a set vector subspace of V. $\eta: W \to [0, 1]$.

$$\eta\,(0\ 0\ 0) = 0$$
$$\eta\,(111) = \dfrac{1}{9}$$
$$\eta\,(011) = 0.$$

$W_\eta$ is the fuzzy set vector subspace of V.

***Example 1.55:*** Let V = {(111), (1011), (11110), (101), (000), (0000), (0000000), (00000), (1111111), (11101), (01010), (1101101)} be a set vector space over the set S = {0,1}.

Let W = {(1111111), (0000000), (000), (00000), (11101), (01010) (101)} $\subset$ V. Define $\eta: W \to [0, 1]$ by

$$\eta(x_1, x_2, \ldots, x_r) = \frac{1}{12}.$$

$\eta W$ is a fuzzy set vector subspace.

We now proceed on to define the notion of fuzzy set linear subalgebra.

**DEFINITION 1.16:** *Let V be a set linear algebra over the set S. Suppose W is a set linear subalgebra of V over S. Let $\eta : W \to [0, 1]$, $\eta W$ is called the fuzzy set linear subalgebra if $\eta (a + b) \geq \min \{\eta (a), \eta (b)\}$ for $a, b, \in W$.*

We give some examples of this new concept.

***Example 1.56:*** Let $V = Z^+ \times Z^+ \times Z^+$ be a set linear algebra over the set $S = 2Z^+$. $W = Z^+ \times 2Z^+ \times 4Z^+$ is a set linear subalgebra over the set $S = 2Z^+$. Define $\eta: W \to [0, 1]$

$$\eta (a\ b\ c) = 1 - \frac{1}{a+b+c}.$$

Clearly $\eta (x, y) \geq \min \{\eta (x), \eta (y)\}$ where $x = (x_1, x_2, x_3)$ and $y = (y_1, y_2, y_3)$; $x, y \in W$. $W\eta$ is a fuzzy set linear subalgebra.

***Example 1.57:*** Let

$$V = \left\{ \begin{pmatrix} a & b \\ c & d \end{pmatrix} \middle| a,b,c,d, \in Z^+ \right\}.$$

V is a set linear algebra over the set $S = \{1, 3, 5, 7\} \subseteq Z^+$. Let

$$W = \left\{ \begin{pmatrix} a & a \\ a & a \end{pmatrix} \middle| a \in Z^+ \right\}$$

be the set linear subalgebra of V. Define $\eta: W \to [0, 1]$ by

$$\eta \begin{pmatrix} a & a \\ a & a \end{pmatrix} = 1 - \frac{1}{a}.$$

$W_\eta$ or $W\eta$ is a set fuzzy linear subalgebra.

For more about set linear algebra and set fuzzy linear algebra please refer [60].

**Chapter Two**

# SPECIAL SET VECTOR SPACES AND FUZZY SPECIAL SET VECTOR SPACES AND THEIR PROPERTIES

In this chapter we for the first time introduce the notion of special set vector spaces and fuzzy special set vector spaces. This chapter has four sections. Section one introduces the notion of special set vector spaces and describes some of their properties. In section two special set vector bispaces are introduced and their properties are studied. Section three generalizes the notion of special set vector bispaces to special set vector n-spaces. The final section introduces the notion of fuzzy special set vector spaces and fuzzy special set n-vector spaces.

## 2.1 Special Set Vector Spaces and their Properties

In this section for the first time we define the notion of special set vector space and special set linear algebra and describe some of their properties.

**DEFINITION 2.1.1**: *Let $V = \{S_1, S_2, ..., S_n\}$ be a set of collection of sets, such that each $S_i$ is distinct i.e., $S_i \not\subseteq S_j$ or $S_j \not\subseteq S_i$ if $i \neq j$; $1 \leq i, j \leq n$. Suppose P is any set which is nonempty such that each $S_i$ is a set vector space over the set P; for each $i = 1, 2, ..., n$ then we call $V = \{S_1, S_2, ..., S_n\}$ to be a special set vector space over the set P.*

We illustrate this situation by some simple examples.

**Examples 2.1.1:** Let $V = (V_1, V_2, V_3, V_4)$ where

$$V_1 = \{(1\ 1\ 1\ 1), (0\ 0\ 0\ 0), (1\ 1\ 0\ 0), (0\ 0\ 1\ 1),$$
$$(1\ 1), (0\ 1), (1\ 0), (0\ 0)\},$$

$$V_2 = \left\{ \begin{pmatrix} a & a & a \\ a & a & a \end{pmatrix}, \begin{pmatrix} a & a \\ a & a \\ a & a \end{pmatrix} \middle| a \in Z_2 = \{0, 1\} \right\},$$

$$V_3 = \{Z_2 \times Z_2 \times Z_2, Z_2 \times Z_2 \times Z_2 \times Z_2 \times Z_5\}$$

and

$$V_4 = \left\{ \begin{pmatrix} 1 & 0 & 0 \\ 0 & 1 & 0 \\ 0 & 0 & 0 \end{pmatrix}, \begin{pmatrix} 1 & 0 & 0 \\ 0 & 1 & 0 \\ 0 & 0 & 0 \end{pmatrix}, \begin{pmatrix} 1 & 1 & 1 \\ 0 & 0 & 0 \\ 1 & 1 & 1 \end{pmatrix}, \begin{pmatrix} 1 & 0 & 0 & 0 \\ 1 & 1 & 0 & 0 \\ 0 & 0 & 1 & 1 \\ 0 & 0 & 0 & 1 \end{pmatrix}, \begin{pmatrix} 1 & 1 & 1 & 1 \\ 0 & 1 & 1 & 1 \\ 0 & 0 & 1 & 1 \\ 0 & 0 & 0 & 1 \end{pmatrix}, \begin{pmatrix} 0 & 0 & 0 & 0 \\ 0 & 0 & 0 & 0 \\ 0 & 0 & 0 & 0 \\ 0 & 0 & 0 & 0 \end{pmatrix} \right\}.$$

Clearly V is a special set vector space over the set $P = \{0, 1\}$.

**Example 2.1.2:** Let $V = \{V_1, V_2, V_3, V_4, V_5\}$ where

$V_1 = \{$all polynomial of degree 3 and all polynomials of even degree with coefficients from $Z^+ \cup \{0\}\}$,

$$V_2 = \{3Z^+ \cup \{0\} \times 5Z^+ \cup \{0\} \times 7Z^+ \cup \{0\}, 8Z^+ \cup \{0\} \times 9Z\},$$

$$V_3 = \left\{ \begin{pmatrix} a & a \\ a & a \end{pmatrix}, \begin{pmatrix} a & a & a \\ a & a & a \\ a & a & a \end{pmatrix} \middle| a \in Z^0 = Z^+ \cup \{0\} \right\}$$

$$V_4 = \left\{ \begin{pmatrix} a & b & c & d \\ 0 & e & f & g \\ 0 & 0 & h & i \\ 0 & 0 & 0 & j \end{pmatrix}, \begin{pmatrix} a & 0 & 0 & 0 \\ b & e & 0 & 0 \\ c & f & g & 0 \\ d & p & h & c \end{pmatrix} \middle| \right.$$

$$a, b, c, d, e, f, g, h, i, p, j \in 2Z^+ \cup \{0\}\}$$

and

$$V_5 = \left\{ \begin{bmatrix} a_1 & a_2 \\ a_3 & a_4 \\ a_5 & a_6 \\ a_7 & a_8 \end{bmatrix}, \begin{bmatrix} a_1 \\ a_2 \\ a_3 \\ a_4 \\ a_5 \\ a_6 \end{bmatrix} \middle| a_i \in 3Z^+ \cup \{0\}, 1 \leq i \leq 8 \right\}.$$

V is a special set vector space over the set $S = Z^0 = Z^+ \cup \{0\}$.

*Example 2.1.3:* $V = \{Z_{10}, Z_{15}, Z_{11}, Z_{19}, Z_{22}, Z_4\} = \{V_1, V_2, V_3, V_4, V_5, V_6\}$ is a special set vector space over the set $S = \{0, 1\}$.

*Example 2.1.4:* Let $V = \{V_1, V_2, V_3, V_4\}$ where $V_1 = \{(1\ 1\ 1\ 1\ 1), (0\ 0\ 0\ 0\ 0), (1\ 1\ 1\ 0\ 0), (1\ 0\ 1\ 0\ 1), (0\ 1\ 1\ 0\ 1), (1\ 1\ 1), (0\ 0\ 0), (1\ 0\ 0), (0\ 0\ 1)\}$, $V_2 = \{(1\ 1), (0\ 0), (1\ 1\ 1\ 1\ 1\ 1), (0\ 0\ 0\ 0\ 0\ 0), (0\ 1\ 0\ 1\ 0\ 1), (1\ 0\ 1\ 0\ 1\ 0), (1\ 1\ 1), (0\ 0\ 0)\}$, $V_3 = \{(1\ 1\ 1\ 1), (0\ 0\ 0\ 0), (1\ 1\ 0\ 1), (0\ 1\ 1\ 1), (1\ 1\ 1\ 1\ 1\ 1), (0\ 0\ 0\ 0\ 0\ 0\ 0), (1\ 1\ 0\ 0\ 0\ 1\ 0)\}$ and $V_4 = \{(1\ 1\ 1\ 1\ 1\ 1\ 1), (0\ 0\ 0\ 0\ 0\ 0), (1\ 1\ 1), (0\ 0\ 0), (1\ 1\ 1\ 1), (0\ 0\ 0\ 0), (1\ 1), (0\ 0)\}$ be a special set vector space over the set $S = Z_2 = \{0, 1\}$. Note we see $S_i \not\subseteq S_j$ and $S_j \not\subseteq S_i$ if $i \neq j$ but $S_i \cap S_j \neq \phi$, that is in general is nonempty if $i \neq j$, $1 \leq i, j \leq 4$.

Now having seen special set vector spaces now we proceed onto define the notion of special set vector subspace.

**DEFINITION 2.1.2:** *Let $V = (V_1, V_2, ..., V_n)$ be a special set vector space over the set S. Take $W = (W_1, W_2, ..., W_n) \subseteq (V_1, V_2, ..., V_n)$ such that $W_i \subset V_i$, $1 \leq i \leq n$ and $W_i \neq W_j$ and $W_j \not\subseteq W_i$, if $i \neq j$, $1 \leq i, j \leq n$. Further if for every $s \in S$ and $w_i \in W_i$, $sw_i \in W$ for $i = 1, 2, ..., n$ then we define $W = (W_1, ..., W_n)$ to be the special set vector subspace of V over the set S. Clearly every subset of V need not in general be a special set vector subspace of V.*

We illustrate this by some simple examples.

***Example 2.1.5:*** Let $V = (V_1, V_2, V_3, V_4, V_5)$ where

$$V_1 = \left\{ \begin{pmatrix} a & a & a & a & a \\ a & a & a & a & a \end{pmatrix} \middle| a \in Z_{10} \right\},$$

$$V_2 = \left\{ \begin{pmatrix} a & a & a \\ a & a & a \\ a & a & a \end{pmatrix}, \begin{pmatrix} a & a \\ a & d \end{pmatrix} \middle| a, d \in Z_{10} \right\},$$

$V_3 = \{Z_{10} \times Z_{10} \times Z_{10}, (1\ 1\ 1\ 1\ 1), (0\ 0\ 0\ 0\ 0), (2\ 2\ 2\ 2\ 2), (6\ 6\ 6\ 6\ 6\ 6), (3\ 3\ 3\ 3\ 3), (4\ 4\ 4\ 4\ 4), (5\ 5\ 5\ 5\ 5), (7\ 7\ 7\ 7\ 7), (8\ 8\ 8\ 8\ 8), (9\ 9\ 9\ 9\ 9)\}$,

$$V_4 = \left\{ \begin{bmatrix} a & a \\ a & a \\ a & a \\ a & a \\ a & a \end{bmatrix}, \begin{bmatrix} a_1 \\ a_2 \\ a_3 \\ a_4 \\ a_5 \end{bmatrix} \middle| a, b, a_i \in Z_{10}, 1 \leq i \leq 5 \right\}$$

and

$V_5 = \{$all polynomials of even degree and all polynomials of degree less than or equal to three with coefficients from $Z_{10}\}$ be a special set vector space over the set $Z_{10}$. Take $W = (W_1, W_2, \ldots, W_5)$ with $W_i \subseteq V_i$, $1 \leq i \leq 5$; where

$$W_1 = \left\{ \begin{pmatrix} a & a & a & a & a \\ a & a & a & a & a \end{pmatrix} \middle| a \in \{0, 2, 4, 6, 8\} \right\} \subseteq V_1,$$

$$W_2 = \left\{ \begin{pmatrix} a & a \\ a & d \end{pmatrix} \middle| a, d \in Z_{10} \right\} \subseteq V_2,$$

$$W_3 = \{Z_{10} \times Z_{10} \times Z_{10}\} \subseteq V_3,$$

$$W_4 = \left\{ \begin{bmatrix} a_1 \\ a_2 \\ a_3 \\ a_4 \\ a_5 \end{bmatrix} \middle| a_i \in Z_{10},\ 1 \leq i \leq 5 \right\} \subseteq V_4$$

and $W_5 = \{$all polynomials of even degree with coefficients from $Z_{10}\} \subseteq V_5$. $W = (W_1, W_2, \ldots, W_5) \subseteq (V_1, V_2, \ldots, V_5) = V$ is clearly a special set vector subspace of $V$ over the set $Z_{10}$.

We give yet another example.

*Example 2.1.6:* Let $V = (V_1, V_2, V_3, V_4)$ where
$$V_1 = \{Z^+[x]\},$$

$$V_2 = \{Z^+ \times 2Z^+ \times 3Z^+ \times 5Z^+, Z^+ \times Z^+\},$$

$$V_3 = \left\{ \begin{pmatrix} a & b \\ c & d \end{pmatrix}, \begin{pmatrix} a_1 & a_2 & a_3 \\ a_4 & a_5 & a_6 \\ a_7 & a_8 & a_9 \end{pmatrix} \middle| a_i \in Z^+, 1 \leq i \leq 9 \right\}$$

and

$$V_4 = \left\{ \begin{bmatrix} a & a & a & a & a \\ a & a & a & a & a \end{bmatrix}, \begin{bmatrix} a & a & a \\ a & a & a \\ a & a & a \\ a & a & a \\ a & a & a \\ a & a & a \\ a & a & a \end{bmatrix} \,\middle|\, a \in Z^+ \right\}$$

be a special set vector space over the set $S = 2Z^+$. Take $W = (W_1, W_2, W_3, W_4) \subseteq (V_1, V_2, V_3, V_4)$ ($W_i \subseteq V_i$; $i = 1, 2, 3, 4$) where $W_1 = \{$all polynomials of degree less than or equal to 6 with coefficients from $Z^+\} \subseteq V_1$, $W_2 = \{Z^+ \times Z^+\} \subseteq V_2$,

$$W_3 = \left\{ \begin{pmatrix} a_1 & a_2 & a_3 \\ a_4 & a_5 & a_6 \\ a_7 & a_8 & a_9 \end{pmatrix} \,\middle|\, a_i \in Z^+\ 1 \leq i \leq 9 \right\} \subseteq V_3$$

and

$$W_4 = \left\{ \begin{bmatrix} a & a & a \\ a & a & a \\ a & a & a \\ a & a & a \\ a & a & a \\ a & a & a \\ a & a & a \end{bmatrix} \,\middle|\, a \in Z^+ \right\} \subseteq V_4$$

to be the subset of V. It is easily verified that $W = (W_1, W_2, W_3, W_4) \subseteq (V_1, V_2, V_3, V_4)$ is a special set vector subspace of V.

Now we give yet another example.

*Example 2.1.7:* Let $V = (V_1, V_2, V_3, V_4, V_5)$ where $V_1 = \{Z_9\}$, $V_2 = \{Z_{14}\}$, $V_3 = \{Z_{16}\}$ $V_4 = \{Z_{18}\}$ and $V_5 = \{Z_8\}$ be a special set vector space over the set $S = \{0, 1\}$. Take $W = (W_1, W_2, W_3, W_4, W_5)$ where $W_i \subseteq V_i$, $i = 1, 2, 3, 4, 5$ and $W_1 = \{0, 3, 6\} \subseteq V_1$, $W_2 = \{0, 2, 4, 6, 8, 10, 12\} \subseteq Z_{14} = V_2$, $W_3 = \{0, 4, 8, 12\} \subseteq V_3$, $W_4 = \{0, 9\} \subseteq V_4$ and $W_5 = \{0, 2, 4, 6\} \subseteq V_5$. Clearly $W = (W_1, W_2, W_3, W_4, W_5) \subseteq V$ is a special set vector subspace of V.

One of the natural question which is pertinent is that, will every proper subset of a special set vector space be a special set vector subspace of V? The answer is no and we prove this by a simple example.

*Example 2.1.8:* Let $V = (V_1, V_2, V_3, V_4)$ where

$$V_1 = \{(1\ 1\ 1\ 1\ 1\ 1), (0\ 0\ 0\ 0\ 0\ 0), (1\ 1\ 1\ 1\ 1), (0\ 0\ 0\ 0\ 0),\\ (0\ 0\ 0\ 0), (1\ 1\ 1\ 1)\},$$

$$V_2 = \left\{ \begin{pmatrix} a & b & c \\ d & e & f \\ g & h & i \end{pmatrix} \middle| a, b, c, d, e, f, g, h, i \in Z_5 \right\},$$

$$V_3 = \left\{ \begin{bmatrix} a & a & a & a & a \\ b & b & b & b & b \end{bmatrix} \middle| a, b \in Z_{10} \right\}$$

and

$$V_4 = \left\{ \begin{bmatrix} a & a & a \\ a & a & a \\ a & a & a \\ a & a & a \\ a & a & a \\ a & a & a \\ a & a & a \end{bmatrix} \middle| a \in Z_7 \right\}$$

be a special set vector space over the set $S = \{0, 1\}$. Take $W = (W_1, W_2, W_3, W_4)$ where $W_1 = \{(1\ 1\ 1\ 1\ 1\ 1), (1\ 1\ 1\ 1), (0\ 0\ 0\ 0\ 0)\} \subseteq V_1$,

$$W_2 = \left\{ \begin{pmatrix} a & a & a \\ a & a & a \\ a & a & a \end{pmatrix} \middle| a \in Z_5 \setminus \{0\} \right\} \subseteq V_2,$$

$$W_3 = \left\{ \begin{bmatrix} 1 & 1 & 1 & 1 & 1 & 1 \\ 2 & 2 & 2 & 2 & 2 & 2 \end{bmatrix}, \begin{bmatrix} 3 & 3 & 3 & 3 & 3 & 3 \\ 4 & 4 & 4 & 4 & 4 & 4 \end{bmatrix}, \right.$$
$$\begin{bmatrix} 5 & 5 & 5 & 5 & 5 & 5 \\ 2 & 2 & 2 & 2 & 2 & 2 \end{bmatrix}, \begin{bmatrix} 7 & 7 & 7 & 7 & 7 & 7 \\ 1 & 1 & 1 & 1 & 1 & 1 \end{bmatrix},$$
$$\left. \begin{bmatrix} 2 & 2 & 2 & 2 & 2 & 2 \\ 9 & 9 & 9 & 9 & 9 & 9 \end{bmatrix} \right\} \subseteq V_3$$

and

$$W_4 = \left\{ \begin{bmatrix} a & a \\ a & a \\ a & a \\ a & a \\ a & a \end{bmatrix} \middle| a \in \{1, 2, 5\} \subseteq Z_7 \right\} \subseteq V_4.$$

Clearly $W = (W_1, W_2, W_3, W_4)$ is a proper subset of V but W is not a special set vector subspace of V. For take (1 1 1 1), (1 1 1 1 1 1) ∈ $W_1$; 0(1 1 1 1) = (0 0 0 0) ∉ $W_1$ and 0(1 1 1 1 1 1) = (0 0 0 0 0 0) ∉ $W_1$ hence $W_1$ is not a set vector subspace of $V_1$ over S = {0, 1}. Also take

$$\begin{pmatrix} a & a & a \\ a & a & a \\ a & a & a \end{pmatrix} \in W_2, \; 0\begin{pmatrix} a & a & a \\ a & a & a \\ a & a & a \end{pmatrix} = \begin{pmatrix} 0 & 0 & 0 \\ 0 & 0 & 0 \\ 0 & 0 & 0 \end{pmatrix} \notin W_2.$$

So $W_2$ is also not a set vector subspace over the set S = {0, 1}.
Likewise in $W_3$ we see

$$\begin{bmatrix} 1 & 1 & 1 & 1 & 1 & 1 \\ 2 & 2 & 2 & 2 & 2 & 2 \end{bmatrix} \in W_3$$

but

$$0 \cdot \begin{bmatrix} 1 & 1 & 1 & 1 & 1 & 1 \\ 2 & 2 & 2 & 2 & 2 & 2 \end{bmatrix} = \begin{bmatrix} 0 & 0 & 0 & 0 & 0 & 0 \\ 0 & 0 & 0 & 0 & 0 & 0 \end{bmatrix} \notin W_3.$$

So $W_3$ too is not a set vector subspace of $V_3$. Finally

$$\begin{bmatrix} a & a \\ a & a \\ a & a \\ a & a \\ a & a \end{bmatrix} \in W_4$$

but

$$0 \begin{bmatrix} a & a \\ a & a \\ a & a \\ a & a \\ a & a \end{bmatrix} = \begin{bmatrix} 0 & 0 \\ 0 & 0 \\ 0 & 0 \\ 0 & 0 \\ 0 & 0 \end{bmatrix}$$

is not in $W_4$. Thus $W_4$ too is not a set vector subspace of $V_4$. So $W = (W_1, W_2, W_3, W_4) \subseteq (V_1, V_2, V_3, V_4)$ is a proper subset of $V$ but is not a special set vector subspace of $V$.

Now we want to make a mention that in the subset $W = (W_1, W_2, W_3, \ldots, W_n)$ even if one of the $W_i \subseteq V_i$ is not a set vector subspace of $V_i$, $1 \le i \le n$, then also $W = (W_1, W_2, \ldots, W_n) \subseteq (V_1, V_2, \ldots V_n)$ will not be a special set vector subspace of $V$.

Now we proceed onto define the notion of generating special subset of a special set vector space $V$ over the set $S$ over which $V$ is defined.

**DEFINITION 2.1.3:** *Let $V = (V_1, V_2, \ldots, V_n)$ be a special set vector space over the set $S$. Let $(P_1, P_2, \ldots, P_n \subseteq V = (V_1, V_2, \ldots, V_n)$ where each $P_i \subseteq V_i$ is such that $P_i$ generates $V_i$ over the set $S$; for $1 \le i \le n$. Then we call $P$ to be the special generating subset of $V$. However it may so happen that at times for some $P_i \subseteq V_i$; $P_i = V_i, i \in \{1, 2, \ldots, n\}$.*

We now give an example of a special generating subset of V.

***Example 2.1.9:*** Let $V = (V_1, V_2, V_3, V_4, V_5)$ where

$$V_1 = \{Z^+ \cup \{0\}\},$$

$$V_2 = \left\{ \begin{pmatrix} a & a & a \\ a & a & a \\ a & a & a \end{pmatrix} \middle| a \in Z^+ \cup \{0\} \right\},$$

$$V_3 = \left\{ \begin{bmatrix} a & a \\ a & a \\ a & a \\ a & a \\ a & a \end{bmatrix} \middle| a \in Z^+ \cup \{0\} \right\},$$

$$V_4 = \{(a\ a\ a\ a\ a), (a\ a\ a\ a\ a\ a\ a) \mid a \in Z^+ \cup \{0\}\}$$

and

$$V_5 = \{\text{all polynomials of the form}$$
$$n(1 + x + x^2 + x^3 + x^4 + x^5) \mid n \in Z^+ \cup \{0\}\}.$$

Clearly $V = (V_1, V_2, V_3, V_4, V_5)$ is a special set vector space over the set $S = Z^+ \cup \{0\}$. Take $P = (P_1, P_2, P_3, P_4, P_5) \subseteq V = (V_1, V_2, V_3, V_4, V_5)$; where $P_1 = \{1\}$,

$$P_2 = \begin{pmatrix} 1 & 1 & 1 \\ 1 & 1 & 1 \\ 1 & 1 & 1 \end{pmatrix}, P_3 = \begin{bmatrix} 1 & 1 \\ 1 & 1 \\ 1 & 1 \\ 1 & 1 \\ 1 & 1 \end{bmatrix},$$

$P_4 = \{(1\ 1\ 1\ 1\ 1), (1\ 1\ 1\ 1\ 1\ 1\ 1)\}$ and $P_5 = \{(1 + x + x^2 + x^3 + x^4 + x^5)\}$. It is easily verified that P is the special generating subset of V over the set $S = Z^+ \cup \{0\}$. Further we say in this case the special set vector space V is finitely generated over S. The n-cardinality of V is denoted by $(|P_1|, |P_2|, ..., |P_n|)$ is called the n-dimension of V.

Hence $|P| = (1, 1, 1, 1, 1)$ is the 5-dimension. However it is pertinent to mention as in case of vector spaces the special set vector spaces are also depend on the set over which it is defined. For if we take in the above example instead of $S = Z^+ \cup \{0\}$, $T = \{0, 1\}$, still $V = (V_1, V_2, V_3, V_4, V_5)$ is a special set vector space over $T = \{0, 1\}$ however now the special generating subset of V over T is infinite. Thus the n-dimension becomes infinite even if one of the generating subset of V say $P_i$ of $(P_1, P_2, ..., P_n) = P$ is infinite; we say the n-dimension of V is infinite.

Now we give yet another example.

*Example 2.1.10:* Let $V = (V_1, V_2, V_3, V_4, V_5)$ where $V_1 = (Z_3 \times Z_3)$, $V_2 = \{Z_2 \times Z_2 \times Z_2\}$, $V_3 = Z_{11}$, $V_4 = \{S \times S \times S \mid S = \{0, 2, 4\} \subseteq Z_6\}$ and $V_5 = Z_7$, be a special set vector space over the set $S = \{0, 1\}$.

Take $P = (P_1, P_2, P_3, P_4, P_5)$ where

$P_1$ = $\{(1\ 1), (1\ 2), (2\ 1), (2\ 2), (0\ 1), (0\ 2), (2\ 0), (1\ 0)\} \subseteq V_1$,
$P_2$ = $\{(1\ 1\ 1), (1\ 1\ 0), (0\ 1\ 1), (1\ 0\ 1), (1\ 0\ 0), (0\ 1\ 0), (0\ 0\ 1)\} \subseteq V_2$,
$P_3$ = $\{1, 2, 3, 4, 5, 6, 7, 8, 9, 10\} \subseteq V_3$,
$P_4$ = $\{(0\ 2\ 4), (0\ 4\ 2), (2\ 2\ 2), (4\ 4\ 4), (0\ 2\ 2), (0\ 4\ 4), (4\ 2\ 0), (2\ 4\ 0), (2\ 0\ 4), (4\ 0\ 2), (2\ 0\ 2), (2\ 2\ 0), (4\ 0\ 4), (4\ 4\ 0), (2\ 2\ 4), (2\ 4\ 2), (4\ 2\ 2), (2\ 4\ 4), (4\ 2\ 4), (4\ 4\ 2), (0\ 0\ 2), (0\ 0\ 4), (0\ 2\ 0), (0\ 4\ 0), (4\ 0\ 0), (2\ 0\ 0)\} \subseteq V_4$ and
$P_5$ = $\{1, 2, 3, 4, 5, 6\} \subseteq Z_7$.

Clearly $P = (P_1, P_2, P_3, P_4, P_5)$ is a special generating subset of V over $S = \{0, 1\}$. Further V is finite 5-dimensional. The 5-dimension of V is $|P| = (|P_1|, |P_2|, |P_3|, |P_4|, |P_5|) = (8, 7, 10, 26, 6)$.

We now give an example of an infinite dimensional special set vector space.

**Example 2.1.11:** Let $V = (V_1, V_2, V_3, V_4)$ where

$V_1 = \{Z[x];$ all polynomials in the variable x with coefficients from $Z\}$,

$$V_2 = (Z^+ \cup \{0\}) \times (Z^+ \cup \{0\}) \times (Z^+ \cup \{0\}),$$

$$V_3 = \left\{ \begin{pmatrix} a_1 & a_2 & a_3 \\ a_4 & a_5 & a_6 \\ a_7 & a_8 & a_9 \end{pmatrix} \middle| a_i \in Z^+ \cup \{0\}, 1 \leq i \leq 9 \right\}$$

and

$$V_4 = \left\{ \begin{bmatrix} a & a & a & a & a \\ a & a & a & a & a \end{bmatrix} \middle| a \in Z^+ \cup \{0\} \right\}$$

be a special set vector space over the set $S = \{0, 1\}$. Now consider the special generating subset $P = (P_1, P_2, P_3, P_4)$ of $V = (V_1, V_2, V_3, V_4)$ where $P_1 = \{$every non zero polynomial $p(x)$ in $Z[x]\} \subseteq Z[x]$ and $P_1 = [Z \setminus \{0\}][x]$, $P_2 = \{\langle x\ y\ z \rangle\} \subseteq V_2$ and $x, y, z \in Z^+$, that is $P_2 = V_2 \setminus \{(0\ 0\ 0)\}$,

$$P_3 = \left\{ \begin{pmatrix} a & b & c \\ d & e & f \\ g & h & i \end{pmatrix} \middle| a, b, c, d, e, f, g, h, i \in Z^+ \cup \{0\} \right\} \setminus$$

$$\left\{ \begin{pmatrix} 0 & 0 & 0 \\ 0 & 0 & 0 \\ 0 & 0 & 0 \end{pmatrix} \right\} \setminus \left\{ \begin{pmatrix} 0 & 0 & 0 \\ 0 & 0 & 0 \\ 0 & 0 & 0 \end{pmatrix} \right\} \subseteq V_3$$

and

$$P_4 = \left\{ \begin{bmatrix} a & a & a & a & a \\ a & a & a & a & a \end{bmatrix} \middle| a \in Z^+ \right\} \subseteq V_4.$$

We see every subset $P_i$ of $V_i$ is infinite. Thus the special set vector space is only infinitely generated and n-dimension of V is $(\infty, \infty, \infty, \infty)$. If we replace the set $S = \{0, 1\}$ by the set $S = Z^+$

∪{0} still the special set vector space is of infinite dimension but not the special generating subset of V given by $T = (T_1, T_2, T_3, T_4)$ is such that $T_1$ is infinite set of polynomials, $T_2$ too is an infinite subset of $V_2$, $T_3$ too is an infinite subset of $V_3$ but different from $P_3$ and

$$T_4 = \left\{ \begin{pmatrix} 1 & 1 & 1 & 1 & 1 \\ 1 & 1 & 1 & 1 & 1 \end{pmatrix} \right\}.$$

Thus $|T| = (|T_1|, |T_2|, |T_3|, |T_4|) = (\infty, \infty, \infty, 1)$ but it is important and interesting to observe that the infinite of $T_i$ is different from the infinite of $P_i$, $1 \le i \le 4$.

Let us now define the notion of special set linear algebra.

**DEFINITION 2.1.4:** *Let $V = (V_1, V_2, ..., V_n)$ be a special set vector space over the set S. If at least one of the $V_i$ of V is a set linear algebra over S then we call V to be a special set linear algebra; $i \in \{(1, 2, ..., n)\}$.*

*We see every special set linear algebra is a special set vector space and not conversely; i.e., a special set vector space in general is not a special set linear algebra.*

We now illustrate this definition by some examples.

*Example 2.1.12:* Let $V = (V_1, V_2, V_3, V_4)$ where

$V_1 = Z_5 \times Z_5 \times Z_5$, $V_2 = \{(a\ a\ a\ a), (a\ a\ a) \mid a \in Z^+ \cup \{0\}\}$,

$$V_3 = \left\{ \begin{pmatrix} a & b \\ c & d \end{pmatrix} \middle| a,b,c,d \in \{0,1\} \right\}$$

and

$$V_4 = \left\{ \begin{pmatrix} a & a \\ a & a \\ a & a \\ a & a \\ a & a \end{pmatrix}, \begin{bmatrix} a & a & a & a & a \\ a & a & a & a & a \end{bmatrix} \middle| a,b,c \in \{0,1\} \right\}$$

be a special set linear algebra over the set S = {0 1}. We see $V_1$ and $V_3$ are set linear algebras over the set S = {0, 1}. However on $V_4$ and $V_2$ we cannot define any compatible operation. So $V_2$ and $V_4$ are not set linear algebras over S = {0, 1} they are only just set vector spaces over S.

We give yet another example.

***Example 2.1.13:*** Let V = ($V_1$, $V_2$, $V_3$, $V_4$, $V_5$) where

$$V_1 = \left\{ \begin{pmatrix} a & b \\ c & d \end{pmatrix} \middle| a,b,c,d \in \{2, 3, 5, 7, 9, 14, 18, 22, 15\} \right\},$$

$V_2$ = {(1 1 1 1 1), (0 0 0 0 0), (1 1 1), (0 0 0), (1 0 1), (1 1 1 1 1), (0 0 0 0), (1 1 0 0), (1 0 0 1)},

$$V_3 = \left\{ \begin{pmatrix} a & a & a \\ a & a & a \end{pmatrix} \middle| a \in Z^+ \cup \{0\} \right\}$$

and

$$V_4 = \left\{ \begin{bmatrix} a \\ a \\ a \\ a \\ a \end{bmatrix} \middle| a \in \{0,1\} \right\}$$

be a special set linear algebra over the set S = {0, 1}. Clearly $V_1$ and $V_2$ are not set linear algebras over the set S = {0, 1}, but $V_3$ and $V_4$ are set linear algebras over the set S = {0, 1}.

Now we proceed onto define some substructures of these special set linear algebras.

**DEFINITION 2.1.5:** *Let V = ($V_1$, $V_2$, ..., $V_n$) be a special set linear algebra over the set S. If W = ($W_1$, $W_2$, ..., $W_n$) ⊆ ($V_1$, $V_2$, ..., $V_n$) with $W_i$ ⊆ $V_i$ where at least one of the $W_i$ ⊆ $V_i$ is a set* *linear subalgebra of a set linear algebra over S and if W = ($W_1$, $W_2$, ..., $W_n$) ⊆ ($V_1$, $V_2$, ..., $V_n$) is also a special set vector subspace of V over S then we call W to be a special set linear subalgebra of V over S.*

We now illustrate this by a simple example.

***Example 2.1.14:*** Let V = ($V_1$, $V_2$, $V_3$, $V_4$, $V_5$) where $V_1$ = {$Z_3$[x] | all polynomials of degree less than or equal to 4 with coefficients from $Z_3$}, $V_2$ = {(a a a a), (a a) | a ∈ $Z_5$},

$$V_3 = \left\{ \begin{pmatrix} a_1 & a_2 \\ a_3 & a_4 \\ a_5 & a_6 \\ a_7 & a_8 \\ a_9 & a_{10} \end{pmatrix}, \begin{bmatrix} a & a & a & a \\ a & a & a & a \end{bmatrix} \middle| a, a_i \in Z_7; 1 \leq i \leq 10 \right\},$$

$$V_4 = \left\{ \begin{pmatrix} a & b & c \\ d & e & f \end{pmatrix}, \begin{pmatrix} a & a \\ a & a \\ a & a \end{pmatrix} \middle| a,b,c,d,e,f \in Z_{11} \right\}$$

and

$$V_5 = \left\{ \begin{pmatrix} a & b \\ c & d \end{pmatrix} \middle| a,b,c,d \in Z_6 \right\}$$

be a special set linear algebra over the set S = {0, 1}. Take W = ($W_1$, $W_2$, ..., $W_5$) with $W_i$ ⊆ $V_i$, 1 ≤ i ≤ 5, where $W_1$ = {0, (1 + x + $x^2$ + $x^3$ + $x^4$), 2(1 + x + $x^2$ + $x^3$ + $x^4$)} ⊆ $V_1$, $W_2$ = {(a a a a) | a ∈ $Z_5$}, ⊆ $V_2$,

$$W_3 = \left\{ \begin{bmatrix} a & a & a & a \\ a & a & a & a \end{bmatrix} \middle| a \in Z_7 \right\} \subseteq V_3,$$

$$W_4 = \left\{ \begin{pmatrix} a & b & c \\ d & e & f \end{pmatrix} \middle| a,b,c,d,e,f \in Z_{11} \right\} \subseteq V_4$$

and

$$W_5 = \left\{ \begin{pmatrix} a & a \\ a & a \end{pmatrix} \middle| a \in Z_6 \right\} \in V_5.$$

W is clearly a special set vector subspace of V as well as W is a special set linear subalgebra of V over S = {0 1}.

We give yet another example.

***Example 2.1.15:*** Let $V = (V_1, V_2, V_3, V_4)$ where

$$V_1 = \{Z_7 \times Z_7 \times Z_7\},$$

$$V_2 = \left\{ \begin{pmatrix} a & b \\ c & d \end{pmatrix} \middle| a,b,c,d \in Z_7 \right\},$$

$$V_3 = \{(a\ a\ a\ a\ a\ a), (a\ a\ a\ a\ a) \mid a \in Z_7\}$$

and

$$V_4 = \left\{ \begin{bmatrix} a & a \\ a & a \\ a & a \\ a & a \end{bmatrix}, \begin{bmatrix} a & a & a \\ a & a & a \end{bmatrix} \middle| a \in Z_7 \right\}$$

be a special set linear algebra over the set $S = Z_7$. Take $W = (W_1, W_2, W_3, W_4)$ where $W_i \subseteq V_i$, $1 \leq i \leq 4$ and

$$W_1 = \{(a\ a\ a) \mid a \in Z_7\} \subseteq V_1,$$

$$W_2 = \left\{ \begin{pmatrix} a & a \\ a & a \end{pmatrix} \middle| a \in Z_7 \right\} \subseteq V_2,$$

$$W_3 = \{(a\ a\ a\ a\ a\ a) \mid a \in Z_7) \subseteq V_3$$

and

$$W_4 = \left\{ \begin{bmatrix} a & a & a \\ a & a & a \end{bmatrix} \bigg| a \in Z_7 \right\}$$

is a special set linear subalgebra of V over $Z_7$.

Now we define yet another new algebraic substructure of V.

**DEFINITION 2.1.6:** *Let $V = (V_1, V_2, ..., V_n)$ be a special set linear algebra over the set S. If $W = (W_1, W_2, ..., W_n) \subseteq (V_1, V_2, ..., V_n)$ with each $W_i \subseteq V_i$ is only a set vector subspace of $V_i$ and never a set linear subalgebra even if one $V_i$ is a set linear algebra $1 \le i \le n$, then we call W to be a pseudo special set vector subspace of the special set linear algebra V over S.*

We illustrate this by the following example.

*Example 2.1.16:* Let $V = (V_1, V_2, ..., V_n)$, where

$$V_1 = \{Z^o[x] \text{ where } Z^o = Z^+ \cup \{0\}\},$$

$$V_2 = \{(a\ a\ a\ a\ a), (a\ a\ a\ a\ a\ a) \mid a \in Z^+ \cup \{0\}\},$$

$$V_3 = \left\{ \begin{pmatrix} a & b \\ c & d \end{pmatrix}, \begin{pmatrix} a & a & a & a \\ b & b & b & b \end{pmatrix} \bigg| a,b,c,d \in Z^+ \cup \{0\} \right\}$$

and

$$V_4 = \left\{ \begin{pmatrix} a & b & c & d \\ 0 & e & f & g \\ 0 & 0 & h & i \\ 0 & 0 & 0 & j \end{pmatrix} \bigg| a,b,c,d,e,f,g,h,i,j \in Z^+ \cup \{0\} \right\}$$

be a special set linear algebra over the set $S = Z^o = Z^+ \cup \{0\}$. Clearly $V_1, V_4$ are set linear algebras over S. Take $W = (W_1, W_2, W_3, W_4) \subseteq V$; where $W_1 = \{n(5x^3 + 3x + 1), n(7x^7 + 8x^6 + 5x + 7), s(3x + 7), t(9x^3 + 49x + 1)\} \subseteq V_1$; s, t, n $\in Z^+ \cup \{0\}\}$.

$W_2 = \{(2a, 2a, 2a, 2a, 2a), (5a, 5a, 5a, 5a, 5a) \mid a \in Z^+ \cup \{0\}\} \subseteq V_2$,

$$W_3 = \left\{ \begin{pmatrix} a & a & a & a \\ a & a & a & a \end{pmatrix} \middle| a, b \in 2Z^+ \cup \{0\} \right\} \subseteq V_3$$

and

$$W_4 = \left\{ n \begin{pmatrix} 1 & 1 & 1 & 1 \\ 0 & 1 & 1 & 2 \\ 0 & 0 & 1 & 2 \\ 0 & 0 & 0 & 2 \end{pmatrix} \middle| n \in Z^+ \cup \{0\} \right\} \subseteq V_4$$

is a special set vector subspace of V which is not a special set linear subalgebra of V; so W is a pseudo special set vector subspace of V.

Now having given example of pseudo special set vector subspace we make a mention of the following result. A pseudo special set vector subspace of V is never a special set linear subalgebra of V.

Now we proceed onto define the notion of special set linear transformation of V onto W, V and W are special set vector spaces over the same set S.

**DEFINITION 2.1.7:** *Let $V = (V_1, V_2, ..., V_n)$ be a special set vector space over the set S. $W = (W_1, W_2, ..., W_n)$ be another special set vector space over the same set S. A special set linear map $T = (T_1, ..., T_n)$ where $T_i: V_i \to W_i$; $1 \leq i \leq n$ such that $T_i(\alpha v) = \alpha T_i(v)$ for all $\alpha \in S$ and $v \in V_i$, $1 \leq i \leq n$ is called the special set linear transformation of V into W. S $Hom_S(V, W) = \{Hom_S(V_1, W_1), ..., Hom_S(V_n, W_n)\}$ where $Hom_S(V_i, W_i)$ denotes the set of all set linear transformations of the vector space $V_i$ into the vector space $W_i$, $1 \leq i \leq n$.*

It is easy to verify $SHom_S(V, W)$ is again a special set vector space over the set S.

We shall illustrate this by some examples.

***Example 2.1.17:*** Let $V = (V_1, V_2, V_3, V_4)$ and $W = (W_1, W_2, W_3, W_4)$ be special set vector spaces over the set $S = \{0, 1\}$, here

$$V_1 = \{Z_2 \times Z_2 \times Z_2\},$$

$$V_2 = \left\{ \begin{pmatrix} a & b \\ c & d \end{pmatrix} \middle| a, b, c, d \in \{0, 1\} \right\},$$

$V_3 = \{Z_2[x]\}$ all polynomials of degree less than or equal to 5, with coefficients from $Z_2 = \{0, 1\}\}$ and

$$V_4 = \left\{ \begin{bmatrix} a_1 & a_2 & a_3 & a_7 \\ a_4 & a_5 & a_6 & a_8 \end{bmatrix} \middle| a_1, \ldots, a_8 \in Z_2 \right\}$$

are set vector spaces of the special set vector space $V$ over the set $S = \{0, 1\}$. Now

$$W_1 = \left\{ \begin{pmatrix} a & b \\ 0 & c \end{pmatrix} \middle| a, b, c \in Z_2 \right\},$$

$$W_2 = Z_2 \times Z_2 \times Z_2 \times Z_2,$$

$$W_3 = \left\{ \begin{pmatrix} a_1 & a_2 & a_3 \\ a_4 & a_5 & a_6 \end{pmatrix} \middle| a_i \in Z_2, 1 \leq i \leq 6 \right\}$$

and

$$W_4 = \left\{ \begin{bmatrix} a_1 & a_5 \\ a_2 & a_6 \\ a_3 & a_7 \\ a_4 & a_8 \end{bmatrix} \middle| a_i \in Z_2, 1 \leq i \leq 8 \right\}$$

are set vector space of the special set vector space $W$ over the set $S = \{0, 1\}$. Define a special set linear map $T = (T_1, T_2, T_3, T_4)$ from $V = (V_1, V_2, V_3, V_4)$ into $W = (W_1, W_2, W_3, W_4)$ such that $T_i : V_i \to W_i$; $1 \leq i \leq 4$ where

$T_1: V_1 \to W_1$ is defined as $T_1(a\ b\ c) = \begin{pmatrix} a & b \\ 0 & c \end{pmatrix}$,

$T_2: V_2 \to W_2;\ T_2 \begin{pmatrix} a & b \\ 0 & c \end{pmatrix} = (a, b, c, d)$,

$T_3: V_3 \to W_3;$

$T_3(a_0 + a_1x + a_2x^2 + a_3x^3 + a_4x^4 + a_5x^5) = \begin{pmatrix} a_0 & a_1 & a_2 \\ a_5 & a_4 & a_5 \end{pmatrix}$

and

$T_4: V_4 \to W_4$ as $T_4\left(\begin{bmatrix} a_1 & a_2 & a_3 & a_4 \\ a_5 & a_6 & a_7 & a_8 \end{bmatrix}\right) = \begin{bmatrix} a_1 & a_5 \\ a_2 & a_6 \\ a_3 & a_7 \\ a_4 & a_8 \end{bmatrix}$.

It is easily verified that $T = (T_1, T_2, T_3, T_4)$ is a special set linear transformation of V into W.

We give yet another example.

***Example 2.1.18:*** Let $V = (V_1, V_2, V_3, V_4, V_5)$ be a special set vector space over the set $S = \{0, 1\}$ where $V_1 = \{Z_7 \times Z_7\}$, $V_2 = Z_{14} \times Z_{14} \times Z_{14}$; $V_3 = Z_{19}$, $V_4 = \{(a\ a\ a\ a\ a) \mid a \in Z_{19}\}$ and $V_5 = \{Z_7 \times Z_{12} \times Z_{14}\}$ is a special set vector space over the set $S = \{0\ 1\}$. Let $W = (W_1, W_2, W_3, W_4, W_5)$ where

$W_1 = \{(a\ a) \mid a \in Z_7\}$,

$W_2 = \left\{ \begin{pmatrix} a & a \\ b & c \end{pmatrix} \middle| a, b, c \in Z_{14} \right\}$,

$W_3 = Z_{19} \times Z_{19}$,

$$W_4 = \left\{ \begin{pmatrix} a & b & c \\ 0 & d & e \end{pmatrix} \middle| \ a, b, c, d, e \in Z_{19} \right\}$$

and

$$W_5 \{(a\ a\ a) \mid a \in Z_7\}$$

be the special set vector space over the set $S = \{0, 1\}$. Define $T = (T_1, T_2, T_3, T_4, T_5)$ from $V = (V_1, V_2, V_3, V_4, V_5)$ into $W = (W_1, W_2, W_3, W_4, W_5)$ by

$$T_1: V_1 \to W_1;\ \text{by } T_1(a, b) = (a\ a),$$

$$T_2: V_2 \to W_2;\ T_2(a\ b\ c) = \begin{pmatrix} a & a \\ b & c \end{pmatrix};$$

$$T_3: V_3 \to W_3\ ;\ \text{by } T_3(a) = (a, a),$$

$$T_4: V_4 \to W_4;\ T_4(a\ a\ a\ a\ a) = \begin{pmatrix} a & a & a \\ 0 & a & a \end{pmatrix}$$

and

$$T_5: V_5 \to W_5\ \text{as}\ T_5(a\ b\ c) = (a\ a\ a).$$

It is clear that $T = (T_1, T_2, T_3, T_4, T_5)$ is a special set linear transformation of V into W.

Now we proceed onto define the notion of special set linear operator on V.

**DEFINITION 2.1.8:** *Let $V = (V_1\ V_2, \ldots, V_n)$ be a special vector space over the set S. A special set linear map $T = (T_1, T_2, \ldots, T_n)$ from $V = (V_1, V_2, V_3, \ldots, V_n)$ into $V = (V_1, V_2, V_3, \ldots, V_n)$ is defined by $T_i: V_i \to V_i$, $i = 1, 2, \ldots, n$ such that $T_i(\alpha V_i) = \alpha T(V_i)$ for every $\alpha \in S$ and $V_i \in V_i$. We call $T = (T_1, \ldots, T_n)$ to be a special set linear operator on V. If we denote by $SHom_s (V, V) = \{Hom_s (V_1, V_1), Hom_s (V_2, V_2), \ldots, Hom_s (V_n, V_n)\}$ then it can be verified $SHom_s(V, V)$ is a special set linear algebra over S under the composition of maps.*

We illustrate this by some simple examples.

***Example 2.1.19:*** Let $V = \{V_1, V_2, V_3, V_4\}$ where $V_1 = Z_3 \times Z_3 \times Z_3$,

$$V_2 = \left\{ \begin{pmatrix} a & b & c & d \\ e & f & g & h \end{pmatrix} \middle| a, b, c, d, e, f, g, h \in Z_5 \right\},$$

$$V_3 = \left\{ \begin{pmatrix} a_1 & b_2 \\ a_3 & a_4 \end{pmatrix} \middle| a_i \in Z_7, 1 \le i \le Z_4 \right\}$$

and

$$V_4 = \left\{ \begin{pmatrix} a & b & c & d \\ 0 & e & f & g \\ 0 & 0 & h & i \\ 0 & 0 & 0 & e \end{pmatrix} \middle| a, b, c, d, e, f, g, h, i, 1 \in Z_9 \right\}$$

be a special set vector space over the set $S = \{0\ 1\}$. Define a special set linear operator $T = (T_1, T_2, T_3, T_4) : V = (V_1, V_2, V_3, V_4) \to (V_1, V_2, V_3, V_4) = V$ where

$T_i: V_i \to V_i$ such that $T_1(a\ b\ c) = (c\ b\ a)$,

$T_2: V_2 \to V_2$ as $T_2 \begin{pmatrix} a & b & c & d \\ e & f & g & g \end{pmatrix} = \begin{pmatrix} e & f & g & h \\ a & b & c & d \end{pmatrix}$,

$T_3: V_3 \to V_3$ defined by $T_3 \begin{pmatrix} a_1 & a_2 \\ a_3 & a_4 \end{pmatrix} = \begin{pmatrix} a_1 & a_4 \\ a_2 & a_3 \end{pmatrix}$

and

$T_4: V_4 \to V_4$ is such that $T_4 \begin{pmatrix} a & b & c & d \\ 0 & e & f & g \\ 0 & 0 & h & i \\ 0 & 0 & 0 & 1 \end{pmatrix} = \begin{pmatrix} a & b & 0 & 0 \\ 0 & e & f & 0 \\ 0 & 0 & 0 & h \\ 0 & 0 & 0 & 1 \end{pmatrix}.$

It is easily verified that $T = (T_1, T_2, T_3, T_4)$ is a special set linear operator on V. Define another special linear operator $P = (P_1, P_2, P_3, P_4)$ from V to V as $P_i: V_i \to V_i$; $i = 1, 2, 3, 4$ given by

$P_1: V_1 \to V_1$ $P_1$ (a b c) = (a a a),

$P_2: V_2 \to V_2$ given by $P_2 \begin{pmatrix} a & b & c & d \\ e & f & g & h \end{pmatrix} = \begin{pmatrix} a & b & c & d \\ 0 & 0 & 0 & 0 \end{pmatrix}$,

$P_3: V_3 \to V_3$ defined by $P_3 \begin{pmatrix} a_1 & a_2 \\ a_3 & a_4 \end{pmatrix} = \begin{pmatrix} a_2 & a_2 \\ a_2 & a_2 \end{pmatrix}$

and

$P_4: V_4 \to V_4$ by $P_4 \begin{pmatrix} a & b & c & d \\ 0 & e & f & g \\ 0 & 0 & h & i \\ 0 & 0 & 0 & 1 \end{pmatrix} = \begin{pmatrix} d & g & i & 1 \\ 0 & c & f & h \\ 0 & 0 & b & e \\ 0 & 0 & 0 & a \end{pmatrix}$.

$P = (P_1, P_2, P_3, P_4)$ is a special set linear operator on V. Define

$$P \circ T = (P_1, P_2, P_3, P_4) \circ (T_1, T_2, T_3, T_4)$$
$$= (P_1 \circ T_1, P_2 \circ T_2, P_3 \circ T_3, P_4 \circ P_4)$$

as
$P_1 \circ T_1: V_1 \to V_1$ as
$$P_1 \circ T_1(a\ b\ c) = P_1(c\ b\ a) = (c\ c\ c).$$

$P_2 \circ T_2: V_2 \to V_2$ given by
$$P_2 \circ T_2 \begin{pmatrix} a & b & c & d \\ e & f & g & h \end{pmatrix} = P_2 \begin{pmatrix} e & f & g & h \\ a & b & c & d \end{pmatrix} = \begin{pmatrix} e & f & g & h \\ 0 & 0 & 0 & 0 \end{pmatrix}.$$

$P_3 \circ T_3: V_3 \to V_3$ defined by
$$P_3 \circ T_3 \begin{pmatrix} a_1 & a_2 \\ a_3 & a_4 \end{pmatrix} = P_3 \begin{pmatrix} a_1 & a_4 \\ a_2 & a_3 \end{pmatrix} = \begin{pmatrix} a_4 & a_4 \\ a_4 & a_4 \end{pmatrix}$$

and
$P_4 \circ T_4: V_4 \to V_4$ gives

$$P_4 T_4 \begin{pmatrix} a & b & c & d \\ 0 & e & f & g \\ 0 & 0 & h & i \\ 0 & 0 & 0 & 1 \end{pmatrix} = P_4 \begin{pmatrix} a & b & 0 & 0 \\ 0 & e & f & 0 \\ 0 & 0 & 0 & h \\ 0 & 0 & 0 & 1 \end{pmatrix} = \begin{pmatrix} 0 & 0 & h & 1 \\ 0 & 0 & f & 0 \\ 0 & 0 & b & e \\ 0 & 0 & 0 & a \end{pmatrix}.$$

It is easily verified P o T = ($P_1$ o $T_1$, $P_2$ o $T_2$, $P_3$ o $T_3$, $P_4$ o $T_4$) is a special set linear operator on V. Now find

$$T \text{ o } P = (T_1 \text{ o } P_1, T_2 \text{ o } P_2, T_3 \text{ o } P_3, T_4 \text{ o } P_4): V \to V.$$

$$T_1 \text{ o } P_1(a\ b\ c) = T_1\ (a\ a\ a) = (a\ a\ a) \neq P_1 \text{ o } T_1.$$

$T_2$ o $P_2$: $V_2 \to V_2$ is such that

$$T_2 \text{ o } P_2 \begin{pmatrix} a & b & c & d \\ e & f & g & h \end{pmatrix} = T_2 \begin{pmatrix} a & b & c & d \\ 0 & 0 & 0 & 0 \end{pmatrix} = \begin{pmatrix} 0 & 0 & 0 & 0 \\ a & b & c & d \end{pmatrix}$$

$$\neq P_2 \text{ o } T_2.$$

$T_3$ o $P_3$: $V_3 \to V_3$ is defined by

$$T_3 \text{ o } P_3 \begin{pmatrix} a_1 & a_2 \\ a_3 & a_4 \end{pmatrix} = T_3 \begin{pmatrix} a_2 & a_2 \\ a_2 & a_2 \end{pmatrix} = \begin{pmatrix} a_2 & a_2 \\ a_2 & a_2 \end{pmatrix}.$$

Thus $T_3$ o $P_3 \neq P_3$ o $T_3$ and now

$T_4$ o $P_4$: $V_4 \to V_4$ given by

$$T_4 \text{ o } P_4 \begin{pmatrix} a & b & c & d \\ 0 & e & f & g \\ 0 & 0 & h & i \\ 0 & 0 & 0 & 1 \end{pmatrix} T_4 \begin{pmatrix} d & g & i & 1 \\ 0 & c & f & g \\ 0 & 0 & b & e \\ 0 & 0 & 0 & a \end{pmatrix} = \begin{pmatrix} d & g & 0 & 0 \\ 0 & c & f & 0 \\ 0 & 0 & 0 & e \\ 0 & 0 & 0 & a \end{pmatrix}.$$

Here also $T_4$ o $P_4 \neq P_4$ o $T_4$. Thus we see T o P ≠ P o T but both T o P and P o T are special set linear operator on V.

Thus we can prove $SHom_s(V, V) = (Hom_s(V_1, V_1), …, Hom_s(V_n, V_n))$ is a special set linear algebra over the set S.

Now we give yet another example of a special set linear operator on V.

*Example 2.1.20:* Let $V = \{V_1, V_2, V_3, V_4\}$ where

$$V_1 = Z_9 \times Z_9 \times Z_9 \times Z_9,$$

$$V_2 = \left\{ \begin{pmatrix} a & b \\ a & a \end{pmatrix} \middle| a, b, c, d \in Z_8 \right\},$$

$$V_3 = \left\{ \begin{bmatrix} a & b & c & d & e \\ f & g & h & i & j \end{bmatrix} \middle| a, b, c, d, e, f, g, h \in Z_6 \right\}$$

and $V_4 = Z_{12} \times Z_{12} \times Z_{12}$ be a special set linear algebra over the set $S = \{0, 1\}$. Define $T = (T_1, T_2, T_3, T_4): V \to V = (V_1, V_2, V_3, V_4)$ where $T_i: V_i \to V_i$ where $T_1(a\ b\ c\ d) = (b\ c\ d\ a)$,

$$T_2 \begin{pmatrix} a & b \\ a & a \end{pmatrix} = \begin{pmatrix} d & b \\ b & d \end{pmatrix}$$

where $T_2: V_2 \to V_2$, $T_3: V_3 \to V_3$ is such that

$$T_3 \left( \begin{bmatrix} a & b & c & d & e \\ f & g & h & i & j \end{bmatrix} \right) = \begin{bmatrix} a & 0 & c & 0 & e \\ f & 0 & h & 0 & j \end{bmatrix}$$

and $T_4: V_4 \to V_4$ is defined as $T_4\ (a\ b\ c) = (a\ c\ a)$.

Clearly $T = (T_1, T_2, T_3, T_4)$ is a special set linear operator on V. Now we find out $T \circ T = (T_1 \circ T_1, T_2 \circ T_2, T_3 \circ T_3, T_4 \circ T_4)$ from $V \to V = (V_1, V_2, V_3, V_4)$ where $T_i \circ T_i; V_i \to V_i; 1 \le i \le 4$

$T_1 \circ T_1 \circ T_1 \circ T_1 (a\ b\ c\ d)$
$= T_1 \circ T_1 \circ T_1 (b\ c\ d\ a)$
$= T_1 \circ T_1 (c\ d\ a\ b)$
$= T_1 (d\ a\ b\ c)$

$\quad = \quad$ (a d c d).

Thus $T_1 \circ T_1 \circ T_1 \circ T_1 = I_1: V_1 \to V_1$. $I_1$ the identity special linear operator on V.

$$T_2 \circ T_2 \begin{pmatrix} a & b \\ c & d \end{pmatrix} = T_2 \begin{pmatrix} d & b \\ b & d \end{pmatrix} = \begin{pmatrix} d & b \\ b & d \end{pmatrix}.$$

Thus $T_2 \circ T_2 = T_2$ on $V_2$.

$$T_3 \circ T_3 \begin{bmatrix} a & b & c & d & e \\ f & g & h & i & j \end{bmatrix}$$
$$= T_3 \begin{bmatrix} a & 0 & c & 0 & e \\ f & 0 & h & 0 & j \end{bmatrix} = \begin{bmatrix} a & 0 & c & 0 & e \\ f & 0 & h & 0 & j \end{bmatrix}.$$

Thus $T_3 \circ T_3 = T_3$ on $V_3$. Finally $T_4 \circ T_4 \circ T_4$ (a b c) = $T_4 \circ T_4$ (a c a) = $T_4$ (a a a) = (a a a). Thus $T_4 \circ T_4 \circ T_4 = T_4 \circ V_4$. Further we see $T \circ T = (T_1 \circ T_1, T_2 \circ T_2, T_3 \circ T_3, T_4 \circ T_4)$ is yet another special set linear operator on V.

Now we proceed onto define the notion of special set idempotent operator on $V = (V_1, V_2, \ldots, V_n)$.

**DEFINITION 2.1.9:** *Let $V = (V_1, V_2, \ldots, V_n)$ be a special set vector space over the set S. Suppose $T = (T_1, T_2, \ldots, T_n): V \to V = (V_1, V_2, \ldots, V_n)$ where $T_i: V_i \to V_i$, $i \leq i \leq n$ be a special set linear operator on V and if $T \circ T = (T_1 \circ T_1, T_2 \circ T_2, \ldots, T_n \circ T_n) = (T_1, T_2, \ldots, T_n)$ then we call T to be a special set linear idempotent operator on V.*

We illustrate this by a simple example.

***Example 2.1.21:*** Let $V = (V_1, V_2, V_3, V_4, V_5)$ be a special set vector space over the set $S = \{0, 1\}$ where $V_1 = \{Z_2 \times Z_2 \times Z_2 \times Z_2\}$,

$$V_2 = \left\{ \begin{bmatrix} a & b \\ c & d \end{bmatrix} \middle| a, b, c, d \in Z_2 \right\},$$

$$V_3 = \left\{ \begin{pmatrix} a_1 & a_2 & a_3 & a_4 & a_5 \\ a_6 & a_7 & a_8 & a_9 & a_{10} \end{pmatrix} \middle| a_i \in Z_2; 1 \leq i \leq 10 \right\},$$

$$V_4 = \left\{ \begin{bmatrix} a_1 & a_6 \\ a_2 & a_7 \\ a_3 & a_8 \\ a_4 & a_9 \\ a_5 & a_{10} \end{bmatrix} \middle| a_i \in Z_2; 1 \leq i \leq 10 \right\}$$

and
$V_5 = \{(1\ 1\ 1\ 1\ 1\ 1\ 1), (0\ 0\ 0\ 0\ 0\ 0\ 0), (1\ 1\ 1\ 1), (0\ 0\ 0\ 0),$
$(0\ 0\ 0), (1\ 1\ 1), (1\ 1\ 1\ 1\ 1), (0\ 0\ 0\ 0\ 0), (1\ 1), (0\ 0)\}.$

Define $T = (T_1, T_2, T_3, T_4, T_5) : V = (V_1, V_2, V_3, V_4, V_5) \to V = (V_1, V_2, \ldots, V_5)$ with $T_i : V_i \to V_i$, $1 \leq i \leq 5$ such that each $T_i$ is a set linear operator of the set vector space $V_i$, $i = 1, 2, \ldots, 5$. $T_1 : V_1 \to V_1$ is such that here $T_1(a\ b\ c\ d) = (a\ a\ a\ a)$, $T_1 \circ T_1 (a\ b\ c\ d) = T_1 (a\ a\ a\ a) = (a\ a\ a\ a)$. Thus $T_1 \circ T_1 = T_1$ on $V_1$ i.e., $T_1$ is a set idempotent linear operator on $V_1$.

$T_2 : V_2 \to V_2$, with

$$T_2 \begin{pmatrix} a & b \\ c & d \end{pmatrix} = \begin{pmatrix} d & d \\ d & d \end{pmatrix}$$

so that

$$T_2 \circ T_2 \begin{pmatrix} a & b \\ c & d \end{pmatrix} = T_2 \begin{pmatrix} d & d \\ d & d \end{pmatrix}.$$

Hence $T_2 \circ T_2 = T_2$ on $V_2$ hence $T_2$ is a set idempotent linear operator on $V_2$. $T_3 : V_3 \to V_3$ defined by

$$T_3 = \begin{pmatrix} a_1 & a_2 & a_3 & a_4 & a_5 \\ a_6 & a_7 & a_8 & a_9 & a_{10} \end{pmatrix} = T_3 \begin{pmatrix} a_2 & a_2 & a_3 & a_3 & a_3 \\ a_7 & a_7 & a_8 & a_8 & a_8 \end{pmatrix}.$$

Thus $T_3 \circ T_3 = T_3$ on $V_3$ i.e. $T_3$ is also a set idempotent linear operator on $V_3$. Finally $T_4: V_4 \to V_4$ defined by $T_4(a_1, \ldots, a_i) = (a_1, \ldots, a_i)$ for every $(a_1, \ldots, a_i) \in V_4$. We see $T_4 \circ T_4 = T_4$, i.e., $T_4$ is also a set idempotent operator on $V_4$. Thus $T \circ T = (T_1 \circ T_1, T_2 \circ T_2, T_3 \circ T_3, T_4 \circ T_4) = (T_1, T_2, T_3, T_4) = T$. Hence $T: V \to V$ is a special set idempotent linear operator on V.

Now we proceed onto define the notion of pseudo special set linear operator on a special set vector space.

**DEFINITION 2.1.10:** *Let $V = (V_1, V_2, \ldots, V_n)$ be a special set vector space over the set S. Let $T = (T_1, T_2, \ldots, T_n): V = (V_1, V_2, \ldots, V_n) \to V = (V_1, V_2, \ldots, V_n)$ be a special set linear transformation from V to V where $T_i = V_i \to V_j$ where $i \neq j$ for at least one $T_i$; $1 \le i, j \le n$. Then we call $T = (T_1, \ldots, T_n)$ to be a pseudo special set linear operator on V. Clearly $SHom_S(V, V) = \{Hom_s(V_{i_1}, V_{i_2}) \ldots Hom_s(V_{i_n}, V_{i_j})\}$ is only a special set vector space over S and not a special set linear algebra in general.*

We illustrate this by an example.

*Example 2.1.22:* Let $V = (V_1, V_2, V_3, V_4)$ be a special set vector space over the set $S = Z^+ \cup \{0\}$; where $V_1 = \{S \times S \times S \times S\}$,

$$V_2 = \left\{ \begin{pmatrix} a & b & c & d \\ e & f & g & h \end{pmatrix} \middle| a,b,c,d,e,f,g,h \in Z^+ \cup \{0\} \right\};$$

$V_3 = \{[Z^+ \cup \{0\}] [x]$ set of all polynomials of degree less than or equal to three$\}$

and

$$V_4 = \left\{ \begin{bmatrix} a & e \\ b & f \\ c & g \\ d & n \end{bmatrix} \middle| a,b,c,d,e,f,g,h \in Z^+ \cup \{0\} \right\}.$$

Define a special map $T = (T_1, T_2, T_3, T_4): V = (V_1, V_2, V_3, V_4) \to (V_1, V_2, V_3, V_4)$ by $T_i: V_i \to V_j$ for at least one $i \neq j$, $1 \leq i, j \leq 4$. Let

$$T_1: V_1 \to V_3$$
$$T_2: V_2 \to V_4$$
$$T_3: V_3 \to V_2$$
$$T_4: V_4 \to V_1$$

be defined on $V_i$'s as follows:

$$T_1(a\ b\ c\ d) = a + bx + cx^2 + dx^3.$$

$$T_2 \begin{pmatrix} a & b & c & d \\ e & f & g & h \end{pmatrix} = \begin{bmatrix} a & e \\ b & f \\ c & g \\ d & h \end{bmatrix};$$

$$T_3((a_0 + a_1x + a_2x^2 + a_3x^3)) = \begin{pmatrix} a_0 & a_1 & a_2 & a_3 \\ 0 & 0 & 0 & 0 \end{pmatrix}$$

and

$$T_4 \left( \begin{bmatrix} a & e \\ b & f \\ c & g \\ d & h \end{bmatrix} \right) = (a\ b\ c\ d).$$

$T = (T_1, T_2, T_3, T_4)$ is a pseudo special set linear operator on V.

Now we proceed onto define special direct sum in special set vector spaces.

**DEFINITION 2.1.11:** *Let $V = (V_1, ..., V_n)$ be a special set vector space over the set S. If each $V_i$ of V can be written as $V_i = W_1^i \oplus ... \oplus W_{n_i}^i$ where $W_t^i \cap W_k^i = \phi$ or zero element if $t \neq k$ where $i_t, i_k \in \{i_1, i_2, ..., i_n\}$ and each $W_k^i$ is a set vector subspace of $V_i$. If this is true for each $i = 1, 2, ..., n$ then*

$$V = (W_1^1 \oplus \ldots \oplus W_{n_1}^1, W_1^2 \oplus \ldots \oplus W_{n_2}^2, \ldots, W_1^n \oplus \ldots \oplus W_{n_n}^n)$$

*is called the special set direct sum of the special set vector space $V = (V_1, V_2, \ldots, V_n)$.*

We first illustrate this situation by an example.

***Example 2.1.23:*** Let $V = (V_1, V_2, V_3, V_4)$ where

$$V_1 = \{(1\ 1\ 1\ 1), (0\ 1), (1\ 0), (0\ 0), (0\ 0\ 0\ 0)\},$$

$$V_2 = \left\{ \begin{pmatrix} a & a \\ a & a \end{pmatrix}, \begin{bmatrix} a & a \\ a & a \\ a & a \\ a & a \end{bmatrix} \middle| a \in \{0, 1\} \right\},$$

$$V_3 = \{Z_2 \times Z_2 \times Z_2\}$$

and

$$V_4 = \{1 + x^2, 0, 1 + x + x^2, 1 + x, x^3 + 1, x^2 + x, x^3 + x, x^3 + x^2, 1 + x^3 + x^2 + x\}$$

be a special set vector space over the set $S = \{0, 1\}$. Write

$$V_1 = \{(1\ 1\ 1\ 1), (0\ 0\ 0\ 0)\} \oplus \{(0\ 1), (1\ 0), (0\ 0)\} = W_1^1 \oplus W_2^1.$$

$$V_2 = \left\{ \begin{pmatrix} a & a \\ a & a \end{pmatrix} \right\} \oplus \left\{ \begin{bmatrix} a & a \\ a & a \\ a & a \\ a & a \end{bmatrix} \right\} = W_1^2 \oplus W_2^2,$$

$$V_3 = \{Z_2 \times \{0\} \times Z_2 \times \{0\}\} \oplus \{\{0\} \times \{0\} \times \{0\} \times Z_2\} \oplus \{\{0\} \times Z_2 \times \{0\} \times \{0\}\} = W_1^3 \oplus W_2^3 \oplus W_3^3$$

and

$$V_4 = \{1 + x^3 + x^2 + x, 1 + x^3\ x^2 + x, 0) + \{x^2 + x + 1, 1 + x^2, 1 + x, 0\} + \{0, x^3 + x, x^3 + x^2\} = W_1^4 \oplus W_2^4 \oplus W_3^4$$

with $W_j^4 \cap W_K^4 = (0)$, if $i \neq j$.

Now
$$V = (V_1, V_2, V_3, V_4)$$
$$= \{W_1^1 \oplus W_2^1, W_1^2 \oplus W_2^2, \ W_1^3 \oplus W_2^3 \oplus W_3^3, W_1^4 \oplus W_2^4 \oplus W_3^4\}$$

is the special direct sum of special set vector subspace of V. However it is interesting to note that the way of writing V as a special direct sum is not unique.

Now we define special set projection operator of V.
Further if
$$V = (W_1^1 \oplus W_2^1, W_1^2 \oplus W_2^2, W_1^3 \oplus W_2^3 \oplus W_3^3, \ W_1^4 \oplus W_2^4 \oplus W_3^4)$$
then W any special set vector subspace of V is
$$(W_{i_1}^1, W_{i_2}^2, \ldots, W_{i_n}^n) \subseteq (V_1, V_2, \ldots, V_n),$$
i.e.; $W_{i_t}^t \subseteq V_t$; $1 \leq t \leq n$. So any way of combinations of $W_{i_t}^t$ will give special set vector subspaces of V.

**DEFINITION 2.1.12:** *Let $V = (V_1, V_2, \ldots, V_n)$ be a special set vector space over the set S. Let $W = (W_1, W_2, \ldots, W_n) \subseteq (V_1, V_2, \ldots, V_n) = V$ where $W_i \subseteq V_i$; $1 \leq i \leq n$ be a special set vector subspace of V. Define $T = (T_1, T_2, \ldots, T_n)$ a special linear operator on V such that $W_i$ is invariant under $T_i$ for $i = 1, 2, \ldots, n$; i.e., $T_i : V_i \to V_i$ with $T_i(W_i) \subseteq W_i$ then T is called the special set projection of V.*

We illustrate this by one example.

*Example 2.1.24:* Let $V = (V_1, V_2, V_3, V_4)$ where

$$V_1 = \{Z_5 \times Z_5 \times Z_5 \times Z_5\},$$

$$V_2 = \left\{ \begin{pmatrix} a & b \\ c & d \end{pmatrix} \middle| a,b,c,d \in Z_7 \right\},$$

$$V_3 = \left\{ \begin{pmatrix} a & a & a & a \\ a & a & a & a \end{pmatrix} \middle| a \in Z_{12} \right\}$$

and

$$V_4 = \left\{ \begin{bmatrix} a_1 & a_2 & a_3 \\ a_4 & a_5 & a_6 \\ a_7 & a_8 & a_9 \\ a_{10} & a_{11} & a_{12} \end{bmatrix} \middle| a_i \in Z_{10}, 1 \leq i \leq 12 \right\}$$

be a special set vector space over the set $S = \{0, 1\}$. Take $W = (W_1, W_2, W_3, W_4) \subseteq V$ where

$$W_1 = \{Z_5 \times Z_5 \times \{0\} \times \{0\}\} \subseteq V_1,$$

$$W_2 = \left\{ \begin{pmatrix} a & a \\ a & a \end{pmatrix} \middle| a \in Z_7 \right\} \subseteq V_2,$$

$$W_3 = \left\{ \begin{pmatrix} a & a & a & a \\ a & a & a & a \end{pmatrix} \middle| a \in \{0, 2, 4, 6, 8, 10\} \right\} \subseteq V_3$$

and

$$W_4 \left\{ \begin{bmatrix} a & a & a & a \\ a & a & a & a \\ a & a & a & a \\ a & a & a & a \end{bmatrix} \middle| a \in Z_{10} \right\} \subseteq V_4$$

be a special set vector subspace of V over S. Define a special set linear operator T on V by $T = (T_1, T_2, T_3, T_4): V = (V_1, V_2, V_3, V_4) \to V = (V_1, V_2, V_3, V_4)$ such that $T_i: V_i \to V_i$, $i = 1, 2, 3, 4$. $T_1: V_1 \to V_1$ is such that $T_1(a\ b\ c\ d) = \{(a\ b\ 0\ 0)\}$, $T_2: V_2 \to V_2$ is such that

$$T_2 \begin{pmatrix} a & b \\ c & d \end{pmatrix} = \begin{pmatrix} a & a \\ a & a \end{pmatrix},$$

$T_3: V_3 \to V_3$ is defined by

$$T_3 \begin{pmatrix} a & a & a & a \\ a & a & a & a \end{pmatrix} = \left\{ \begin{pmatrix} b & b & b & b \\ b & b & b & b \end{pmatrix} \middle| b \in \{0, 2, 4, 6, 8, 10\} \right\}$$

and

$T_4: V_4 \to V_4$ is such that

$$T_4\left(\begin{bmatrix} a_1 & a_2 & a_3 \\ a_4 & a_5 & a_6 \\ a_7 & a_8 & a_9 \\ a_{10} & a_{11} & a_{12} \end{bmatrix}\right) = \begin{bmatrix} a_1 & a_1 & a_1 \\ a_1 & a_1 & a_1 \\ a_1 & a_1 & a_1 \\ a_1 & a_1 & a_1 \end{bmatrix}$$

such that $a_1 \in Z_{10}$. $T = (T_1, T_2, T_3, T_4)$ is easily seen to be a special set linear operator on V. Further it can be verified that T is a special set projection operator on V. For T o T = $(T_1$ o $T_1$, $T_2$ o $T_2$, $T_3$ o $T_3$, $T_4$ o $T_4)$ where $T_i$ o $T_i: V_i \to V_I$ given by $T_1$ o $T_1$ (a b c d) = $T_1$ (a b 0 0) = (a b 0 0), i.e.; $T_1$ o $T_1 = T_1$ on $V_1$, i.e., $T_1$ is a set projection operator on $V_1$. We see

$T_2$ o $T_2 : V_2 \to V_2$ gives

$$T_2 \text{ o } T_2 \left(\begin{pmatrix} a & b \\ c & d \end{pmatrix}\right) = T_2\left(\begin{pmatrix} a & a \\ a & a \end{pmatrix}\right) = \begin{pmatrix} a & a \\ a & a \end{pmatrix},$$

i.e., $T_2$ o $T_2 = T_2$ is again a set projection operator on $V_2$. Consider $T_3$ o $T_3 : V_3 \to V_3$ given by

$$T_3 \text{ o } T_3 \left(\begin{pmatrix} a & a & a & a \\ a & a & a & a \end{pmatrix}\right) = T_3\left(\begin{pmatrix} b & b & b & b \\ b & b & b & b \end{pmatrix}\right)$$

such that $b \in \{0, 2, 4, 6, 8, 10\}$ and

$$T_3\left(\begin{pmatrix} b & b & b & b \\ b & b & b & b \end{pmatrix}\right) = \begin{pmatrix} b & b & b & b \\ b & b & b & b \end{pmatrix}$$

then by making $T_3$ is set projector operator on $V_3$; i.e., $T_3$ o $T_3 = T_3$. Finally $T_4$ o $T_4: V_4 \to V_4$ gives

$$T_4 \circ T_4 \left( \begin{bmatrix} a_1 & a_2 & a_3 \\ a_4 & a_5 & a_6 \\ a_7 & a_8 & a_9 \\ a_{10} & a_{11} & a_{12} \end{bmatrix} \right) = T_4 \begin{bmatrix} a_1 & a_1 & a_1 \\ a_1 & a_1 & a_1 \\ a_1 & a_1 & a_1 \\ a_1 & a_1 & a_1 \end{bmatrix} = \begin{bmatrix} a_1 & a_1 & a_1 \\ a_1 & a_1 & a_1 \\ a_1 & a_1 & a_1 \\ a_1 & a_1 & a_1 \end{bmatrix}$$

there by making $T_4 \circ T_4 = T_4$; thus $T_4$ is a set projection operator on $V_4$.

Thus $T = (T_1, T_2, T_3, T_4)$ is such that $T \circ T = (T_1 \circ T_1, T_2 \circ T_2, T_3 \circ T_3, T_4 \circ T_4)$ leaving $W = (W_1, W_2, W_3, W_4)$ to be invariant under T. $T(W) = (T_1(W_1), T_2(W_2), T_3(W_3), T_4(W_4)) \subseteq (W_1, W_2, W_3, W_4) = W \subseteq (V_1, V_2, V_3, V_4) = V$. We see $T \circ T = T (T_1, T_2, T_3, T_4)$.

Interesting results in this direction can be obtained by interested readers.

## 2.2 Special Set Vector Bispaces and their Properties

In this section we introduce the notion of special set vector bispaces and study some of their properties. Now we proceed onto define the notion of special set vector bispaces.

**DEFINITION 2.2.1:** *Let $V = (V_1 \cup V_2)$ where $V_1 = (V_1^1, V_2^1, \ldots, V_n^1)$ and $V_2 = (V_1^2, V_2^2, \ldots, V_{n_2}^2)$ are distinct special set vector spaces over the same set S i.e., $V_1 \not\subseteq V_2$ or $V_2 \not\subseteq V_1$. Then we call $V = V_1 \cup V_2$ to be a special set bivector space or special set vector bispace over the set S.*

We give some examples of special set vector bispaces.

*Example 2.2.1:* Let

$$V = V_1 \cup V_2 = \{ V_1^1, V_2^1, V_3^1, V_4^1 \} \cup \{ V_1^2, V_2^2, V_3^2, V_4^2, V_5^2 \}$$

where
$$V_1^1 = \{Z_2 \times Z_2\},$$

$$V_2^1 = \left\{ \begin{pmatrix} a & b \\ c & d \end{pmatrix} \middle| a,b,c,d \in Z_5 \right\},$$

$$V_3^1 = \left\{ \begin{bmatrix} a & a & a & a & a \\ a & a & a & a & a \end{bmatrix} \middle| a \in Z_7 \right\}$$

and

$$V_4^1 = \left\{ \begin{bmatrix} a & a & a & a & a & a \\ a & a & a & a & a & a \end{bmatrix} \middle| a \in Z_{10} \right\}$$

so that $V_1 = (V_1^1, V_2^1, V_3^1, V_4^1)$ is a special set vector space over the set $S = \{0, 1\}$. Here in $V_2 = (V_1^2, V_2^2, V_3^2, V_4^2, V_5^2)$ we have $V_1^2 = Z_7 \times Z_7 \times Z_7 \times Z_7$, $V_2^2 = \{Z_5[x]$ is set all polynomials in the variable x of degree less than or equal to 4$\}$,

$$V_3^2 = \left\{ \begin{pmatrix} a_1 & a_2 & a_3 \\ a_4 & a_5 & a_6 \\ a_7 & a_8 & a_9 \end{pmatrix} \middle| a_i \in Z_9; 1 \le i \le 9 \right\},$$

$$V_4^2 = \{(a\ a\ a), (a\ a), (a\ a\ a\ a\ a) \mid a \in Z_{11}\}$$

and

$$V_5^2 = \left\{ \begin{bmatrix} a_1 & a_2 \\ a_3 & a_4 \\ a_5 & a_6 \\ a_7 & a_8 \end{bmatrix}, \begin{pmatrix} a & a \\ a & a \end{pmatrix} \middle| a, a_i \in Z_8, 1 \le i \le 8 \right\}$$

is a special set vector space over the set $S = \{0, 1\}$. Thus $V = V_1 \cup V_2 = (V_1^1, V_2^1, V_3^1, V_4^1) \cup (V_1^2, V_2^2, V_3^2, V_4^2, V_5^2)$ is a special set vector bispace over the set $S = \{0, 1\}$.

*Example 2.2.2:* Let $V = V_1 \cup V_2 = (V_1^1, V_2^1, V_3^1) \cup (V_1^2, V_2^2, V_3^2)$ where

$$V_1^1 = \left\{ \begin{pmatrix} a & b \\ c & d \end{pmatrix} \middle| a,b,c,d \in S^o = Z^+ \cup \{0\} \right\},$$

$$V_2^1 = \{S^o \times S^o\},$$

$$V_3^1 = \left\{ \begin{bmatrix} a \\ b \\ c \\ d \end{bmatrix} \middle| a,b,c,d \in S^o = Z^+ \cup \{0\} \right\},$$

$$V_1^2 = \{(a\ a\ a\ a\ a), (a\ a\ a\ a) \mid a \in Z^+ \cup \{0\} = S^o\},$$

$$V_2^2 = \{S^o \times 2S^o \times 3S^o \times 4S^o\}$$

and

$$V_3^2 = \left\{ \begin{bmatrix} a \\ a \\ a \\ a \end{bmatrix}, \begin{bmatrix} a \\ a \\ a \\ a \\ a \\ a \end{bmatrix} \middle| a \in S^o \right\}.$$

Clearly $V_1$ and $V_2$ are special vector spaces over the same set $S^o$. Thus $V = V_1 \cup V_2$ is a special set vector bispace over the set $S^o$.

Now we proceed onto define the notion of special set linear bialgebra over the set S.

**DEFINITION 2.2.2:** *Let $V = V_1 \cup V_2$ be a special set vector bispace over the set S. If both $V_1$ and $V_2$ are special set linear algebras on the set S, then we call $V = V_1 \cup V_2$ to be a special set linear bialgebra.*

It is important to mention here certainly all special set linear bialgebras are special ser vector bispaces however all special set vector bispaces in general are not special set linear bialgebras.

We first illustrate this by a simple example.

***Example 2.2.3:*** Let
$$V = V_1 \cup V_2 = (V_1^1, V_2^1, V_3^1) \cup (V_1^2, V_2^2, V_3^2, V_4^2)$$
where
$$V_1^1 = \{(a\ a\ a\ a), (a\ a) | a \in Z_2 = \{0\ 1\}\},$$

$$V_2^1 = \left\{ \begin{pmatrix} a & a \\ a & a \\ a & a \end{pmatrix}, \begin{pmatrix} a & a & a \\ a & a & a \end{pmatrix} \middle| a \in Z_2\{0,1\} \right\}$$

and
$$V_3^1 = \{0, x^2 + 1, x^3 + x + 1, x + 1, x^4 + x^3 + x + 1,$$
$$x^2 + x^5 + x^4 + 1, x^7 + 1\}.$$

Clearly $V_1 = (V_1^1, V_2^1, V_3^1)$ is not a special set linear algebra for none of the $V_i^1$'s are set linear algebras; i =1, 2, 3,. Now even if $V_2 = (V_1^2, V_2^2, V_3^2, V_4^2)$ is a special set linear algebra over the set $\{0, 1\}$, still $V = V_1 \cup V_2$ is not a special set linear bialgebra over the set $\{0, 1\}$.

In view of this example we have the following new notion.

**DEFINITION 2.2.3:** *Let $V = V_1 \cup V_2$ where $V_1 = (V_1^1, V_2^1, \ldots, V_{n_1}^1)$ is a special set linear algebra over the set S and $V_2 = (V_1^2, V_2^2, \ldots, V_{n_2}^1)$ is a special set vector space over the set S which is not a special set linear algebra. We call $V = V_1 \cup V_2$ to be a quasi special set linear bialgebra.*

Now we proceed onto describe the substructures.

**DEFINITION 2.2.4:** *Let*
$$V = (V_1 \cup V_2) = (V_1^1, V_2^1, \ldots, V_{n_1}^1) \cup (V_1^2, V_2^2, \ldots, V_{n_2}^2)$$
*where both $V_1$ and $V_2$ are special set vector spaces over the set S. Thus $V = V_1 \cup V_2$ is a special set bivector space over the set S. Now consider a proper biset* $W = W_1 \cup W_2 = (W_1^1, W_2^1, \ldots, W_{n_1}^1) \cup (W_1^2, W_2^2, \ldots, W_{n_2}^2) \subseteq V_1 \cup V_2$
such that $W_t^i \subseteq V_t^i$; $i = 1, 2$ and $1 \leq t \leq n_1$ or $n_2$ and each $W_t^i$ is a set vector subspace of $V_t^i$, $i = 1, 2$ and $t = 1, 2, \ldots n_1$ and $t = 1, 2, \ldots, n_2$, we call $W_1 \cup W_2 = (W_1^1, \ldots, W_{n_1}^1) \cup (W_1^2, \ldots, W_{n_2}^2)$ as a special set vector bisubspace of $V = V_1 \cup V_2$.

We illustrate this situation by some simple examples.

*Example 2.2.4:* Let
$$V = V_1 \cup V_2 = \{V_1^1, V_2^1, V_3^1\} \cup \{V_1^2, V_2^2, V_3^2, V_4^2\}$$
where
$$V_1^1 = \{Z_2 \times Z_2 \times Z_2\},$$
$$V_2^1 = \left\{ \begin{pmatrix} a & b \\ c & d \end{pmatrix} \middle| a, b, c, d \in Z_2 \right\}$$
and
$$V_3^1 = \left\{ \begin{pmatrix} a & a & a & a \\ a & a & a & a \end{pmatrix}, \begin{pmatrix} a & a & a \\ a & a & a \end{pmatrix} \middle| a \in Z_2 \right\}$$

is a special set vector space over the set $S = \{0, 1\}$.

$$V_2^1 = \{Z_9 \times Z_9 \times Z_9\},$$

$$V_2^2 = \left\{ \begin{pmatrix} a & b & c \\ c & d & f \end{pmatrix} \middle| a, b, e, c, d, f \in Z_6 \right\},$$

$$V_3^2 = \left\{ \begin{bmatrix} a \\ a \\ a \\ a \\ a \end{bmatrix}, \begin{bmatrix} a \\ a \\ a \\ a \\ a \\ a \end{bmatrix} \middle| a \in Z_{12} \right\}$$

and

$$V_4^2 = \left\{ \begin{pmatrix} a & a & a & a & a \\ a & a & a & a & a \end{pmatrix}, \begin{pmatrix} a & a & a & a & a \\ a & a & a & a & a \end{pmatrix} \middle| a \in Z_{10} \right\}.$$

Clearly $V_2 = \{V_1^2, V_2^2, V_3^2, V_4^2\}$ is a special set vector space over the same set $S = \{0, 1\}$. We see $V = V_1 \cup V_2$ is a special set vector bispace over the set $S = \{0, 1\}$. Take $W = W_1 \cup W_2 = (W_1^1, W_2^1, W_3^1) \cup (W_1^2, W_2^2, W_3^2, W_4^2) \subseteq V_1 \cup V_2$ where $W_1^1 \subseteq V_1^1$ and

$$W_1^1 = \{Z_2 \times Z_2 \times \{0\}\} \subseteq V_1^1,$$

$$W_2^1 = \left\{ \begin{pmatrix} a & a \\ a & a \end{pmatrix} \middle| a \in Z_2 \right\} \subseteq V_2^1, \text{ and}$$

$$W_3^1 = \left\{ \begin{pmatrix} a & a & a \\ a & a & a \end{pmatrix} \middle| a \in Z_2 \right\} \subseteq V_3^1.$$

Thus $W_1 = (W_1^1, W_2^1, W_3^1) \subseteq (V_1^1, V_2^1, V_3^1) = V_1$ is a special set vector subspace of $V_1$. Consider $W_2 = (W_1^2, W_2^2, W_3^2, W_3^2)$ where

$$W_1^2 = \{Z_9 \times Z_9 \times \{0\}\} \subseteq V_1^2,$$

$$W_2^2 = \left\{ \begin{pmatrix} a & a & a \\ a & a & a \end{pmatrix} \middle| a \in Z_6 \right\} \subseteq V_2^2,$$

$$W_3^2 = \left\{ \begin{bmatrix} a \\ a \\ a \\ a \\ a \end{bmatrix} \middle| a \in Z_{12} \right\} \subseteq V_3^2$$

and

$$W_4^2 = \left\{ \begin{pmatrix} a & a & a & a & a & a \\ a & a & a & a & a & a \end{pmatrix} \middle| a \in Z_{10} \right\} \subseteq V_4^2.$$

Thus $W_2 = (W_1^2, W_2^2, W_3^2, W_4^2) \subseteq (V_1^2, V_2^2, V_3^2, V_4^2)$ is a special set vector subspace of $V_2$. Thus $W_1 \cup W_2 = (W_1^1, W_2^1, W_3^1)$ $\cup (W_1^2, W_2^2, W_3^2, W_4^2) \subseteq V_1 \cup V_2 = V$ is a special set vector bisubspace of $V$.

We give yet another example.

**Example 2.2.5:** Let $V = V_1 \cup V_2 = (V_1^1, V_2^1, V_3^1, V_4^1) \cup (V_1^2, V_2^2, V_3^2, V_4^2)$ be a special set vector bispace of $V$ over the set $S = Z^+ \cup \{0\} = Z^o$. Here $V_1 = (V_1^1, V_2^1, V_3^1, V_4^1)$ where

$$V_1^1 = \left\{ \begin{pmatrix} a & b \\ c & d \end{pmatrix} \middle| a, b, c, d \in Z^o \right\},$$

$$V_2^1 = \{Z^o \times Z^o \times Z^o\},$$

$$V_3^1 = \{(a\ a\ a\ a),\ (a\ a),\ (a\ a\ a\ a\ a\ a) \mid a \in Z^o\}$$

and

$V_4^1 = \{Z^o[x]$ set of all polynomials of degree less than or equal to four with coefficients from $Z^o\}$.

$V_1 = (V_1^1, V_2^1, V_3^1, V_4^1)$ is a special set vector space over the set $S = Z^o$. Now

$$V_2 = (V_1^2, V_2^2, V_3^2, V_4^2)$$

where

$$V_1^2 = \left\{ \begin{pmatrix} a_1 & a_2 & a_3 \\ a_4 & a_5 & a_6 \\ a_7 & a_8 & a_9 \end{pmatrix} \middle| a_i \in Z^o, 1 \le i \le 9 \right\},$$

$V_2^2 = \{n(x^2 + 1), n(x + 3x^2 + 9x^7 + 1), n(x^9 + 3x^6 + 6x^3 + 9) \mid n$ a positive integer $n \in Z^o\}$,

$$V_3^2 = \left\{ \begin{pmatrix} a \\ a \\ a \end{pmatrix}, \begin{bmatrix} a \\ a \\ a \\ a \\ a \end{bmatrix} \middle| a \in Z^o \right\}$$

and
$$V_4^2 = \{Z^o \times Z^o \times Z^o \times Z^o \times Z^o\},$$

$V_2 = (V_1^2, V_2^2, V_3^2, V_4^2)$ is a special set vector space over the set $Z^o$. Clearly $V = (V_1^2, V_2^2, V_3^2, V_4^2)$ is a special set vector bispace over the set $Z^o$. Take
$$W = W_1 \cup W_2$$
$$= \{W_1^1, W_2^1, W_3^1, W_4^1\} \cup \{W_1^2, W_2^2, W_3^2, W_4^2\} \subseteq V_1 \cup V_2$$
with
$$W_1^1 = \left\{ \begin{pmatrix} a & a \\ a & a \end{pmatrix} \middle| a \in Z^o \right\} \subseteq V_1^1,$$
$$W_2^1 = Z^o \times Z^o \times \{0\} \subseteq V_2^1,$$
$$W_3^1 = \{(a\ a\ a\ a), (a\ a) \mid a \in Z^o\} \subseteq V_3^1$$
and
$$W_4^1 = \{n(1+x+x^2+x^3+x^4) \mid n \in Z^o\} \subseteq V_4^1;$$
$W_1 = (W_1^1, W_2^1, W_3^1, W_4^1) \subseteq V_1$ is a special set vector subspace of $V_1$ over the set $Z^o$.

$$W_2 = (W_1^2, W_2^2, W_3^2, W_4^2) \subseteq V_2$$
where
$$W_1^2 = \left\{ \begin{pmatrix} a & a & a \\ a & a & a \\ a & a & a \end{pmatrix} \middle| a \in Z^o \right\} \subseteq V_1^2,$$
$$W_2^2 = \{n(x^2 + 1) \mid n \in Z^o) \subseteq V_2^2,$$

$$W_3^2 = \left\{ \begin{bmatrix} a \\ a \\ a \end{bmatrix} \middle| a \in Z^o \right\} \subseteq V_3^2$$

and

$$W_4^2 = \{Z^o \times \{0\} \times Z^o \times \{0\} \times Z^o\} \subseteq V_4^2.$$

$W_2 = (W_1^2, W_2^2, W_3^2, W_4^2) \subseteq V_2$ is a special set vector subspace of $V_2$ over $Z^o$. Thus

$$\begin{aligned} W &= W_1 \cup W_2 \\ &= \{W_1^1, W_2^1, W_3^1, W_4^1\} \cup \{W_1^2, W_2^2, W_3^2, W_4^2\} \\ &\subseteq V_1 \cup V_2 \\ &= \{V_1^1, V_2^1, V_3^1, V_4^1\} \cup \{V_1^2, V_2^2, V_3^2, V_4^2\} \end{aligned}$$

is a special set vector subbispace of $V = V_1 \cup V_2$.

Now we proceed onto define the notion of special set linear subbialgebra of $V = V_1 \cup V_2$.

**DEFINITION 2.2.5**: *Let*
$$V = V_1 \cup V_2 = (V_1^1, V_2^1, \ldots, V_{n_1}^1) \cup (V_1^2, V_2^2, \ldots, V_{n_2}^2)$$
*be a special set linear bialgebra over the set S. Suppose $W = W_1 \cup W_2 = (W_1^1, W_2^1, \ldots, W_{n_1}^1) \cup (W_1^2, W_2^2, \ldots, W_{n_2}^2) \subseteq V_1 \cup V_2$ is such that $W_i$ is a special set linear algebra of $V_i$ over the set S for $i = 1, 2$ then we call $W = W_1 \cup W_2$ to be a special set linear subbialgebra of $V = V_1 \cup V_2$ over the set S.*

We illustrate this by some examples.

***Example 2.2.6:*** Let $V = (V_1^1, V_2^1, V_3^1, V_4^1) \cup (V_1^2, V_2^2, V_3^2, V_4^2) = V_1 \cup V_2$ be a special set linear bialgebra over the set $S = \{0, 1\}$. Here $V_1 = (V_1^1, V_2^1, V_3^1, V_4^1)$ is a special set linear algebra over the set $S = \{0, 1\}$ where $V_1^1 = \{Z_2 \times Z_2 \times Z_2 \times Z_2 \times Z_2\}$,

$$V_2^1 = \left\{ \begin{pmatrix} a & b \\ c & d \end{pmatrix} \middle| a,b,c,d \in Z_5 \right\},$$

$$V_3^1 = \{(a\ a\ a\ a\ a),\ (a\ a\ a\ a) \mid a \in Z_7\}$$

and

$V_4^1 = \{Z_2[x]$ all polynomials of degree less than or equal to 10$\}$.

$V_2 = (V_1^2, V_2^2, V_3^2, V_4^2)$ is a special set linear algebra over the set $S = \{0, 1\}$ where

$$V_1^2 = \{Z_{10} \times Z_{10} \times Z_{10} \times Z_{10} \times Z_{10}\},$$

$$V_2^2 = \left\{ \begin{pmatrix} a & b & c & d \\ e & f & g & h \end{pmatrix} \middle| a,b,c,d,e,f,g,h \in Z_9 \right\},$$

$$V_3^2 = \left\{ \begin{bmatrix} a & a \\ a & a \\ a & a \\ a & a \end{bmatrix}, \begin{bmatrix} a & a & a & a & a \\ a & a & a & a & a \end{bmatrix} \middle| a \in Z_{12} \right\}$$

and

$$V_4^2 = \{(a\ a\ a\ a\ a\ a\ a\ a) \mid a \in Z^+ \cup \{0\}\}$$

is also a special set linear algebra over the set $S = \{0, 1\}$. Now take

$$\begin{aligned} W &= W_1 \cup W_2 \\ &= \{\{W_1^1, W_2^1, W_3^1, W_4^1\} \cup \{W_1^2, W_2^2, W_3^2, W_4^2\} \\ &\subseteq V_1 \cup V_2 \end{aligned}$$

where

$$W_1^1 \subseteq V_1^1, W_2^1 \subseteq V_2^1, W_3^1 \subseteq V_3^1 \text{ and } W_4^1 \subseteq V_4^1;$$

similarly

$$W_i^2 \subseteq V_i^2;\ 1 \le i \le 4.$$

Here

$$W_1^1 = \{Z_2 \times Z_2 \times Z_2 \times \{0\} \times \{0\}\} \subseteq V_1^1$$

$$W_2^1 = \left\{ \begin{pmatrix} a & a \\ a & a \end{pmatrix} \middle| a \in Z_5 \right\} \subseteq V_2^1,$$

$$W_3^1 = \{(a\ a\ a\ a\ a) \mid a \in Z_7\}$$

and

$W_4^1 = \{$all polynomials of degree less than or equal to 5 with coefficients from $S = \{0, 1\}\}$.

Clearly
$$W_1 = (W_1^1, W_2^1, W_3^1, W_4^1) \subseteq (V_1^1, V_2^1, V_3^1, V_4^1) = V_1$$
is a special set linear subalgebra of $V_1$.
$$W_2 = (W_1^2, W_2^2, W_3^2, W_4^2)$$
is such that
$$W_1^2 = \{Z_{10} \times Z_{10} \times \{0\} \times \{0\} \times Z_{10}\} \subseteq V_1^2,$$

$$W_2^2 = \left\{ \begin{pmatrix} a & a & a & a \\ a & a & a & a \end{pmatrix} \middle| a \in Z_9 \right\} \subseteq V_2^2,$$

$$W_3^2 = \left\{ \begin{bmatrix} a & a & a & a & a \\ a & a & a & a & a \\ a & a & a & a & a \end{bmatrix} \middle| a \in Z_{12} \right\} \subseteq V_3^2$$

and
$$W_4^2 = \{(a\ a\ a\ a\ a\ a\ a\ a) \mid a \in 5Z^+ \cup \{0\}\} \subseteq V_4^2.$$

Clearly $W_2 = (W_1^2, W_2^2, W_3^2, W_4^2) \subseteq V_2 = (V_1^2, V_2^2, V_3^2, V_4^2)$ is a special set linear subalgebra of $V_2$. Thus

$$\begin{aligned} W_1 \cup W_2 &= (W_1^1, W_2^1, W_3^1, W_4^1) \cup (W_1^2, W_2^2, W_3^2, W_4^2) \\ &\subseteq V_1 \cup V_2 \\ &= V. \end{aligned}$$

V is a special set linear bisubalgebra of $V = V_1 \cup V_2$.

We give another example of special set linear bisubalgebra.

***Example 2.2.7***: Let $V = (V_1 \cup V_2) = \{V_1^1, V_2^1, V_3^1, V_4^1, V_5^1\} \cup \{V_1^2, V_2^2, V_3^2\}$ where

$$V_1^1 = \left\{ \begin{pmatrix} a & b \\ c & d \end{pmatrix} \middle| a, b, c, d \in Z^o = Z^+ \cup \{0\} \right\},$$

$V_2^1 = \{(a\ a\ a\ a\ a), (a\ a) \mid a \in Z^o\}$, $V_3^1 = \{Z^o[x]$ of all polynomials in the variable x of degree less than or equal to 8 with coefficients from $Z^o\}$, $V_4^1 = \{Z^o \times Z^o \times Z^o \times Z^o\}$ and $V_5^1 = \{$all n × n upper triangular matrices with entries from $Z^o\}$. Thus $V_1 = (V_1^1, V_2^1, V_3^1, V_4^1, V_5^1)$ is a special set linear algebra over the set $Z^o = Z^+ \cup \{0\}$. Now $V_2 = (V_1^2, V_2^2, V_3^2)$ where $V_1^2 = \{(a\ a\ a) \mid a \in Z^o\}$, $V_2^2 = \{$all 10 × 10 matrices with entries from $Z^o\}$ and $V_3^2 = \{Z^o \times 2Z^o \times 5Z^o \times 7Z^o \times 8Z^o \times 9Z^o\}$ is a special set linear algebra over the set $Z^o$. Thus $V = V_1 \cup V_2$ is a special set linear bialgebra over the set $Z^o$. Now consider

$$\begin{aligned}
W &= W_1 \cup W_2 \\
&= \{W_1^1, W_2^1, W_3^1, W_4^1, W_5^1\} \cup \{W_1^2, W_2^2, W_3^2\} \\
&\subseteq V_1 \cup V_2 \\
&= (V_1^1, V_2^1, V_3^1, V_4^1, V_5^1) \cup (V_1^2, V_2^2, V_3^2)
\end{aligned}$$

where

$$W_1^1 = \left\{ \begin{pmatrix} a & a \\ a & a \end{pmatrix} \middle| a \in Z^o \cup \{0\} \right\} \subseteq V_1^1,$$

$W_2^1 = \{(a\ a\ a\ a\ a) \mid a \in 5Z^o\} \subseteq V_2^1$, $W_3^1 = \{Z_4^o[x]$ all polynomials of degree less than or equal to 4 with coefficients from $Z^o\} \subseteq V_3^1$, $W_4^1 = \{Z^o \times \{0\} \times Z^o\{0\}\} \subseteq V_4^1$ and $W_5^1 = \{$all n × n upper triangular matrices with entries from $5Z^o\} \subseteq V_5^1$.

Clearly
$W_1 = (W_1^1, W_2^1, W_3^1, W_4^1, W_5^1) \subseteq V = (V_1^1, V_2^1, V_3^1, V_4^1, V_5^1)$

is a special set linear subalgebra of $V_1$. Now $W_1^2 = \{(a\ a\ a) \mid a \in 7Z^o\} \subseteq V_1^2$, $W_2^2 = \{$all $10 \times 10$ matrices with entries from $5Z^o\}$ and $W_3^2 = \{Z^o \times \{0\} \times 5Z^o \times 7Z^o \times \{0\} \times \{0\}\} \subseteq V_3^2$. It is easily verified that $W_2 = (W_1^2, W_2^2, W_3^2) \subseteq V_2 = (V_1^2, V_2^2, V_3^2)$ is a special set linear subalgebra of $V_2$. Thus

$$W_1 \cup W_2 = \{W_1^1, W_2^1, W_3^1, W_4^1, W_5^1)\} \cup (W_1^2, W_2^2, W_3^2)\}$$
$$\subseteq V_1 \cup V_2$$

is a special set linear subbialgebra of V over the set $Z^o$.

Now we proceed onto define yet another new substructure.

**DEFINITION 2.2.6**: *Let*
$$V = (V_1^1, V_2^1, \ldots, V_{n_1}^1) \cup (V_1^2, V_2^2, \ldots, V_{n_2}^2) = V_1 \cup V_2$$
*be a special set linear bialgebra defined over the set S. Suppose*
$$\begin{aligned} W &= W_1 \cup W_2 \\ &= (W_1^1, W_2^1, \ldots, W_{n_1}^1) \cup (W_1^2, W_2^2, \ldots, W_{n_2}^2) \\ &\subseteq V_1 \cup V_2 \\ &= (V_1^1, V_2^1, \ldots, V_{n_1}^1) \cup (V_1^2, V_2^2, \ldots, V_{n_2}^2) \end{aligned}$$
*is such that $W_1 = (W_1^1, W_2^1, \ldots, W_{n_1}^1) \subseteq V_1$ is also a pseudo special set vector subspace of $V_1$ over S and $W_2 = (W_1^2, W_2^2, \ldots, W_{n_2}^2) \subseteq V_2$ is only a pseudo special set vector subspace of $V_2$ over S. Then we call $W = W_1 \cup W_2 \subseteq V_1 \cup V_2$ to be a pseudo special set vector bisubspace of V over S. Now if in this definition one of $W_1$ or $W_2$ is taken to be a special set linear subalgebra of $V_i$ (i = 1 or 2) then we call $W = W_1 \cup W_2 \subseteq V_1 \cup V_2$ to be a pseudo special set linear subbialgebra of V over the set S.*

We illustrate this set up by one example.

*Example 2.2.8*: Let
$$V = V_1 \cup V_2 = \{V_1^1, V_2^1, V_3^1, V_4^1\} \cup \{V_1^2, V_2^2, V_3^2\}$$
be a special set linear bialgebra over the set S = {0, 1}, where

$$V_1^1 = \{Z_2 \times Z_2 \times Z_2 \times Z_2\},$$

$$V_2^1 = \left\{ \begin{pmatrix} a & b \\ c & d \end{pmatrix} \middle| a, b, c, d \in Z_{20} \right\},$$

$V_3^1 = \{$All 5×5 matrices with entries from the set $\{-1, 0\ 1\}\}$

and

$$V_4^1 = \left\{ \begin{bmatrix} a_1 \\ a_2 \\ a_3 \\ a_4 \end{bmatrix} \middle| a_i \in Z_{10}, 1 \leq i \leq 4 \right\}.$$

$V_1 = (V_1^1, V_2^1, V_3^1, V_4^1)$ is a special set linear algebra over the set $S = \{0, 1\}$.

Now $V_1^2 = \{Z_5 \times Z_5 \times Z_5\}$, $V_2^2 = \{$all $4 \times 4$ matrices with entries from $Z_2 = \{0, 1\}$; under matrix addition using the fact $1 + 1 \equiv 0 \mod 2$, is a semigroup$\}$, $V_3^2 = \{Z_2[x]$, all polynomials of degree less than or equal to 15 with entries from the set $\{0, 1\}\}$. Clearly $V_2 = (V_1^2, V_2^2, V_3^2)$ is a special set linear algebra over the set $S = \{0, 1\}$. Thus $V_1 \cup V_2$ is the given special set linear bialgebra over the set $S = \{0, 1\}$. Take

$$\begin{aligned} W &= W_1 \cup W_2 \\ &= (W_1^1, W_2^1, W_3^1, W_4^1) \cup (W_1^2, W_2^2, W_3^2) \\ &\subseteq V_1 \cup V_2 \end{aligned}$$

where $W_1^1 = \{(1\ 1\ 0\ 0), (0\ 1\ 1\ 1), (1\ 1\ 1\ 1), (0\ 0\ 0\ 0)\} \subseteq V_1^1$; $W_1^1$ is just a pseudo special set vector subspace of $V_1^1$,

$$W_2^1 = \left\{ \begin{pmatrix} a & 0 \\ b & 0 \end{pmatrix} \begin{pmatrix} a_1 & a_2 \\ 0 & 0 \end{pmatrix} \middle| a_1, a_2, a, b \in Z_{20} \right\} \subseteq V_2^1,$$

clearly $W_2^1$ is also a pseudo special set vector subspace of $V_2^1$. Now

$$W_3^1 = \left\{ \begin{pmatrix} 0 & 0 & 0 & 0 & 0 \\ 0 & 0 & 0 & 0 & 0 \\ 0 & 0 & 0 & 0 & 0 \\ 0 & 0 & 0 & 0 & 0 \\ 0 & 0 & 0 & 0 & 0 \end{pmatrix}, \begin{pmatrix} 1 & -1 & 0 & 1 & 1 \\ 1 & 1 & 0 & 0 & -1 \\ 1 & 1 & 1 & 1 & -1 \\ 1 & -1 & 0 & 1 & 0 \\ 0 & 0 & -1 & 0 & 0 \end{pmatrix}, \right.$$

$$\left. \begin{pmatrix} 1 & 1 & 1 & 1 & 1 \\ 1 & 1 & 1 & 1 & 1 \\ 1 & 1 & 1 & 1 & 1 \\ 1 & 1 & 1 & 1 & 1 \\ 1 & 1 & 1 & 1 & 1 \end{pmatrix}, \begin{pmatrix} 0 & 1 & 1 & 0 & -1 \\ 1 & 0 & 0 & -1 & 0 \\ 0 & 0 & 1 & 0 & 1 \\ 1 & 0 & 1 & 0 & 1 \\ 0 & 1 & -1 & 0 & 0 \end{pmatrix} \right\} \subseteq V_3^1$$

$W_3^1$ is only just a pseudo special set vector subspace of $V_3^1$.
Now

$$W_4^1 = \left\{ \begin{bmatrix} 0 \\ 0 \\ 0 \\ 0 \end{bmatrix}, \begin{bmatrix} a_1 \\ 0 \\ a_2 \\ 0 \end{bmatrix}, \begin{bmatrix} 0 \\ a_1 \\ 0 \\ a_2 \end{bmatrix} \middle| a_1, a_2 \in Z_{10} \right\} \subseteq V_4^1$$

is only a pseudo special set vector subspace of $V_4^1$. Thus $W_1 = (W_1^1, W_2^1, W_3^1, W_4^1) \subseteq V_1$ is only a pseudo special set vector subspace of $V_1$ over the set $S = \{0, 1\}$. We choose $W_1^2 = \{(1\ 2\ 3), (4\ 1\ 0), (1\ 1\ 2), (0\ 0\ 0), (4\ 1\ 1)\} \subseteq V_1^2$ is only a pseudo set vector subspace of $V_1^2$ over the set $S = \{0, 1\}$,

$$W_2^2 = \left\{ \begin{bmatrix} 1 & 1 & 1 & 0 \\ 0 & 0 & 0 & 0 \\ 1 & 1 & 1 & 1 \\ 0 & 0 & 0 & 0 \end{bmatrix}, \begin{bmatrix} 0 & 0 & 0 & 0 \\ 0 & 0 & 0 & 0 \\ 0 & 0 & 0 & 0 \\ 0 & 0 & 0 & 0 \end{bmatrix}, \right.$$

$$\left\{ \begin{bmatrix} 1 & 1 & 1 & 1 \\ 0 & 1 & 1 & 1 \\ 0 & 0 & 1 & 1 \\ 0 & 0 & 0 & 1 \end{bmatrix}, \begin{bmatrix} 1 & 1 & 1 & 1 \\ 1 & 0 & 0 & 0 \\ 1 & 0 & 0 & 0 \\ 1 & 0 & 0 & 0 \end{bmatrix}, \begin{bmatrix} 1 & 0 & 0 & 0 \\ 1 & 1 & 0 & 0 \\ 1 & 0 & 0 & 0 \\ 1 & 0 & 0 & 1 \end{bmatrix} \right\} \subseteq V_2^2$$

is only a pseudo set vector subspace of $V_2^2$ over the set $S = \{0, 1\}$. Finally $W_3^2 = \{x^3 + 1, 0, x^7 + 1 + x, x^2 + x + x^3 + x^4 + 1, 1 + x^{10} + x^9, x^{10} + x^2 + x^6 + x^8 + x^4 + 1\} \subseteq V_3^2$, $W_3^2$ is only a pseudo set vector subspace of $V_3^2$. Thus $W_2 = (W_1^2, W_2^2, W_3^2) \subseteq V_2$ is only a pseudo special set vector subspace of $V_2$. So

$$\begin{aligned} W &= W_1 \cup W_2 \\ &= (W_1^1, W_2^1, W_3^1, W_4^1) \cup (W_1^2, W_2^2, W_3^2) \\ &\subseteq V_1 \cup V_2 \end{aligned}$$

is a pseudo special set vector bi subspace $V = V_1 \cup V_2$ over the set $S = \{0, 1\}$.

Now we show that $V = V_1 \cup V_2$ has also pseudo special set linear bi subalgebras. To this end we define

$$\begin{aligned} W &= W_1 \cup W_2 \\ &= \{W_1^1, W_2^1, W_3^1, W_4^1\} \cup \{W_1^2, W_2^2, W_3^2\} \\ &\subseteq V_1 \cup V_2 \end{aligned}$$

as follows: $W_1^1 = \{Z_2 \times Z_2 \times \{0\} \times \{0\}\} \subseteq V_1^1$ is a set linear subalgebra of $V_1^1$,

$$W_2^1 = \left\{ \begin{pmatrix} a & a \\ a & a \end{pmatrix} \middle| a \in Z_{20} \right\} \subseteq V_2^1$$

is again a set linear subalgebra of $V_2^1$, $W_3^1 = \{5 \times 5$ lower triangular matrices with entries from $\{0, 1\}\} \subseteq V_3^1$ is just a set vector subspace of $V_3^1$.

$$W_4^1 = \left\{ \begin{bmatrix} a \\ a \\ a \\ a \\ a \end{bmatrix} \middle| a \in Z_{10} \right\} \subseteq V_4^1$$

is again a set linear subalgebra of $V_4^1$. Thus

$$W_1 = (W_1^1, W_2^1, W_3^1, W_4^1) \subseteq V_1$$

is a special set linear subalgebra of $V_1$. Take

$$W_2 = (W_1^2, W_2^2, W_3^2)$$

as above so that $W_2$ is only a pseudo special set vector subspace of $V_2$. Thus $W = W_1 \cup W_2 \subseteq V_1 \cup V_2$ is a pseudo special set linear subbialgebra of $V = V_1 \cup V_2$.

Now we define the special bigenerator set of the special set vector bispace.

**DEFINITION 2.2.7**: *Let*

$$V = V_1 \cup V_2 = \{V_1^1, V_2^1, \ldots, V_{n_1}^1\} \cup \{V_1^2, V_2^2, \ldots, V_{n_2}^2\}$$

*be a special set vector bispace over the set S. If $X = X_1 \cup X_2 = (X_1^1, X_2^1, \ldots, X_{n_1}^1) \cup (X_1^2, X_2^2, \ldots, X_{n_2}^2) \subseteq V_1 \cup V_2$, i.e., $X_{p_i}^i \subseteq V_{p_i}^i$ $1 \le i \le 2$ with $i = 1, 2$; $1 \le p_1 \le n_1$ and $1 \le p_2 \le n_2$ is such that each $X_i^1$ set generates the set vector space $V_i^1$ over S, $1 \le i \le n_1$ and $X_i^2$ set generates the set vector space $V_i^2$ over S, $1 \le i \le n_2$; then we define $X = X_1 \cup X_2 = (X_1^1, X_2^1, \ldots, X_{n_1}^1) \cup (X_1^2, X_2^2, \ldots, X_{n_2}^2) \subseteq V_1 \cup V_2$ to be the special bigenerator subset of V. If each $X_{t_i}^i$ is finite, $1 \le i \le 2$ with $1 \le t_1 \le n_1$ and $1 \le t_2 \le n_2$ then we say $V = V_1 \cup V_2$ is bigenerated finitely over S. Even if one of the sets $X_t^i$ is of infinite cardinality then we say V is bigenerated infinitely. The bicardinality of X is given by*

$$|X| = |X_1| \cup |X_2|$$
$$= (|X_1^1|, |X_2^1|, \ldots, |X_{n_1}^1|) \cup (|X_1^2|, |X_2^2|, \ldots, |X_{n_2}^2|).$$

*In case $V = V_1 \cup V_2$ is a special set linear bialgebra then the special bigenerating biset $X = X_1 \cup X_2$ is such that*

$(X_1^1, X_2^1, \ldots, X_{n_1}^1)$ is a special generator of the special set linear algebra $V_1$ over $S$ and $(X_1^2, X_2^2, \ldots, X_{n_2}^2)$ is the special set bigenerator of $V_2$ over $S$ as a special set linear algebra then we call $X = X_1 \cup X_2 \subseteq V_1 \cup V_2$ as the special bigenerator set of the special set linear bialgebra over $S$.

We illustrate this by some examples.

***Example 2.2.9***: Let $V = V_1 \cup V_2 = \{V_1^1, V_2^1, V_3^1\} \cup \{V_1^2, V_2^2\}$ where

$$V_1^1 = \left\{ \begin{pmatrix} a & a \\ a & a \end{pmatrix} \middle| a \in Z_2 \right\},$$

$V_2^1 = Z_2 \times Z_2$ and $V_3^1 = \{$all polynomials of the form $1 + x$, $1 + x^2$, $0$, $x^2$, $x + x^2$, $1 + x + x^2$, $x\}$. Clearly $V_1 = \{V_1^1, V_2^1, V_3^1\}$ is a special set vector space over the set $S = \{0, 1\}$. Now take

$$V_1^2 = \left\{ \begin{pmatrix} a & a \\ a & a \\ a & a \end{pmatrix} \middle| a \in Z_2 = \{0, 1\} \right\}$$

and $V_2^2 = \{(1\ 1\ 1), (0\ 0\ 0), (1\ 1\ 0), (0\ 1\ 1)\}$. $V_2 = (V_1^2, V_2^2)$ is also a special set vector space over the set $S = \{0, 1\}$. Thus $V = V_1 \cup V_2$ is a special set vector bispace over the set $S = \{0, 1\}$.
Take
$X = X_1 \cup X_2$

$= \left\{ \begin{pmatrix} 1 & 1 \\ 1 & 1 \end{pmatrix},\ \{(1\ 1), (0\ 1), (1\ 0)\}, \right.$

$(x, x^2, x + x^2, 1 + x, 1+x^2\ 1 + x + x^2)\} \cup$

$\left\{ \begin{bmatrix} 1 & 1 \\ 1 & 1 \\ 1 & 1 \end{bmatrix}, ((1\ 1\ 0), (1\ 1\ 1), (0\ 1\ 1)) \right\}$

$\subseteq\ V_1 \cup V_2.$

It is easily verified that $X = X_1 \cup X_2$ special bigenerates $V_1 \cup V_2$. Now the special bidimension of $V = V_1 \cup V_2$ is $|X| = |X_1| \cup |X_2| = (1, 6, 1) \cup (1, 3)$.

We give yet another example.

*Example 2.2.10*: Let
$$V = V_1 \cup V_2 = \{V_1^1, V_2^1, V_3^1, V_4^1\} \cup \{V_1^2, V_2^2, V_3^2, V_4^2, V_5^2\}$$
be a special set vector bispace over the set $Z^+ \cup \{0\} = S$. Here $V_1^1 = \{S \times S\}$, $V_2^1 = \{(a\ a\ a) \mid a \in S\}$,

$$V_3^1 = \left\{ \begin{pmatrix} a & a \\ a & a \end{pmatrix} \middle| a \in S \right\}$$

and

$$V_4^1 = \left\{ \begin{bmatrix} a \\ a \\ a \\ a \end{bmatrix}, \begin{bmatrix} a & a & a & a & a \\ a & a & a & a & a \end{bmatrix} \middle| a \in S^o \right\}.$$

Now
$$X_1 = (X_1^1, X_2^1, X_3^1, X_4^1) \subseteq V_1 = (V_1^1, V_2^1, V_3^1, V_4^1),$$
i.e., $V_i^1$ contains $X_i^1$ as a subset for $i = 1, 2, 3, 4$. $X_1$ is a special set generator of $V_1$ over $S$ where $X_1^1 = \{$an infinite set of the form $(x, y)\}$;

$$X_2^1 = \{(1\ 1\ 1)\} \subseteq V_2^1,\ X_3^1 = \left\{ \begin{pmatrix} 1 & 1 \\ 1 & 1 \end{pmatrix} \right\} \subseteq V_3^1$$

and

$$X_4^1 = \left\{ \begin{bmatrix} 1 \\ 1 \\ 1 \\ 1 \end{bmatrix}, \begin{bmatrix} 1 & 1 & 1 & 1 & 1 \\ 1 & 1 & 1 & 1 & 1 \end{bmatrix} \right\} \subseteq V_4^1.$$

Now $X_1 = (X_1^1, X_2^1, X_3^1, X_4^1)$ is a special set generator of $V_1$. Let $X_2 = (X_1^2, X_2^2, X_3^2, X_4^2, X_5^2) \subseteq V_2 = (V_1^2, V_2^2, \ldots, V_5^2)$; $X_i^2 \subseteq V_i^2$; $1 \le i \le 5$.

$$V_1^2 = \left\{ \begin{pmatrix} a & a & a \\ a & a & a \\ a & a & a \end{pmatrix} \middle| a \in S^o = Z^+ \cup \{0\} \right\},$$

$$V_2^2 = \{(a\ a\ a\ a\ a) \mid a \in S^o\},\ V_3^2 = \left\{ \begin{pmatrix} a & b \\ c & d \end{pmatrix} \middle| a, b, c, d \in S^o \right\},$$

$$V_4^2 = \left\{ \begin{bmatrix} a \\ a \\ a \\ a \\ a \end{bmatrix}, \begin{bmatrix} a \\ a \\ a \\ a \\ a \\ a \end{bmatrix} \middle| a \in S^o \right\}$$

and
$V_5^2 = \{S^o (x) = (Z^+ \cup \{0\})(x)$, all polynomials in x with coefficients from $S^o\}$. $V_2 = (V_1^2, V_2^2, V_3^2, V_4^2, V_5^2)$ is a special set vector space over $S^o$. Now
$$X_2 = (X_1^2, X_2^2, X_3^2, X_4^2, X_5^2) \subseteq (V_1^2, V_2^2, \ldots, V_5^2)$$
be the special set generator of $V_2$ over $S^o$, where

$$X_1^2 = \left\{ \begin{pmatrix} 1 & 1 & 1 \\ 1 & 1 & 1 \\ 1 & 1 & 1 \end{pmatrix} \right\} \subseteq V_1^2,$$

$$X_2^2 = \{(1\ 1\ 1\ 1)\} \subseteq V_2^2,$$

$X_3^2 = \{$infinite set of matrices with elements from $S^o\} \subseteq V_3^2$,

$$X_4^2 = \left\{ \begin{bmatrix} 1 \\ 1 \\ 1 \\ 1 \end{bmatrix}, \begin{bmatrix} 1 \\ 1 \\ 1 \\ 1 \\ 1 \end{bmatrix} \right\} \subseteq V_4^2$$

and

$$X_5^2 = \{\text{infinite set of polynomials}\}.$$

Clearly $X_2 = (X_1^2, X_2^2, \ldots, X_5^2)$ special generates $V_2$. Now we see $|X_1 \cup X_2| = |X_1| \cup |X_2| = \{(\infty, 1, 1, 2) \cup (1, 1, \infty, 2, \infty)\}$. We see $V = V_1 \cup V_2$ is a bigenerated by the special biset $X = X_1 \cup X_2$.

We have just seen examples of special set bivector spaces which are finite dimensional as well as examples of infinite dimensional.

Now we wish to record if $V = V_1 \cup V_2$ is a special set vector bispace and if it is made into a special set linear bialgebra then we see the special bidimension of $V$ as a special set vector bispace is always greater than or equal to the special bidimension of $V$ as a special set linear bialgebra.

We illustrate this by an example as every special set linear bialgebra is a special set vector bispace.

*Example 2.2.11*: Let
$$V = V_1 \cup V_2 = (V_1^1, V_2^1, V_3^1, V_4^1) \cup (V_1^2, V_2^2, V_3^2, V_4^2, V_5^2)$$
be a special set linear bialgebra over the set $S = Z^+ \cup \{0\}$. Here $V_1 = (V_1^1, V_2^1, V_3^1, V_4^1)$ where

$$V_1^1 = \left\{ \begin{pmatrix} a & b \\ 0 & d \end{pmatrix} \middle| a, b, c, d \in Z^+ \cup \{0\} \right\},$$

$V_2^1 = \{S \times S \times S \times S\}$, $V_3^1 = \{Z^+ \cup \{0\} \times Z^+ \cup \{0\}\}$

and

$$V_4^1 = \left\{ \begin{bmatrix} a_1 \\ a_2 \\ a_3 \end{bmatrix}, \begin{bmatrix} a_1 & a_2 & a_3 \\ a_4 & a_5 & a_6 \end{bmatrix} \middle| a_i \in S; 1 \le i \le 6 \right\}$$

is a special set linear algebra over $S = Z^+ \cup \{0\}$.
Let
$$X_1 = (X_1^1, X_2^1, X_3^1, X_4^1) \subseteq (V_1^1, V_2^1, V_3^1, V_4^1)$$
where
$$X_1^1 = \left\{ \begin{pmatrix} 1 & 0 \\ 0 & 0 \end{pmatrix}, \begin{pmatrix} 0 & 1 \\ 0 & 0 \end{pmatrix}, \begin{pmatrix} 0 & 0 \\ 0 & 1 \end{pmatrix} \right\} \subseteq V_1^1$$

set generates $V_1^1$ as a linear algebra over $Z^+ \cup \{0\}$, $X_2^1 = \{(1\ 0\ 0\ 0), (0\ 1\ 0\ 0), (0\ 0\ 1\ 0), (0\ 0\ 0\ 1)\} \subseteq V_2^1$ set generates $V_2^1$ as a set linear algebra over $S = Z^+ \cup \{0\}$, $X_3^1 = \{(1\ 0), (0\ 1)\} \subseteq V_3^1$ set generates $V_3^1$ as a set linear algebra over S. Finally $V_4^1$ is generated infinitely by $X_4^1 = \{$this is an infinite set$\}$.

So $X_1 = (X_1^1, X_2^1, X_3^1, X_4^1)$ is a special generator of $V_1$ and $|X_1| = (|X_1^1|, |X_2^1|, |X_3^1|, |X_4^1|) = \{3, 4, 2, \infty\}$ as a special set linear algebra. Now let $V_2 = (V_1^2, V_2^2, V_3^2, V_4^2, V_5^2)$ be given by $V_1^2 = \{Z^o \times Z^o \times Z^o\}$ where $Z^o = Z^+ \cup \{0\} = S$,

$$V_2^2 = \left\{ \begin{bmatrix} a \\ a \\ a \\ a \\ a \end{bmatrix} \text{ such that } a \in Z^o \right\},$$

$$V_3^2 = \left\{ \begin{bmatrix} a & a & a \\ a & a & a \\ a & a & a \end{bmatrix}, \begin{bmatrix} a & a & a & a & a \\ a & a & a & a & a \end{bmatrix} \middle| a \in Z^o \right\},$$

$V_4^2 = \{Z^o[x]$ all polynomials of degree less than or equal to three$\}$

and

$$V_5^2 = \left\{ \begin{bmatrix} a & 0 & 0 & 0 \\ a & a & 0 & 0 \\ a & a & a & 0 \\ a & a & a & a \end{bmatrix} \middle| a \in Z^o \cup \{0\} \right\}.$$

Thus $V_2 = (V_1^2, V_2^2, V_3^2, V_4^2, V_5^2)$ is a special set linear algebra over the set $S = Z^o = Z^+ \cup \{0\}$. Consider
$$X_2 = (X_1^2, X_2^2, X_3^2, X_4^2, X_5^2) \subseteq (V_1^2, V_2^2, V_3^2, V_4^2, V_5^2)$$
where
$$X_1^2 = \{(1\ 0\ 0), (0\ 1\ 0), (0\ 0\ 1)\} \subseteq V_1^2$$
generates $V_1^2$ as a set linear algebra over $S = Z^o$,

$$X_2^2 = \left\{ \begin{bmatrix} 1 \\ 1 \\ 1 \\ 1 \\ 1 \end{bmatrix} \right\} \subseteq V_2^2$$

is the set generator of the set linear algebra $V_2^2$.

$$X_3^2 = \left\{ \begin{bmatrix} 1 & 1 & 1 \\ 1 & 1 & 1 \end{bmatrix}, \begin{bmatrix} 1 & 1 & 1 & 1 & 1 \\ 1 & 1 & 1 & 1 & 1 \end{bmatrix} \right\} \subseteq V_3^2$$

is the set generator of the set vector space $V_3^2$,

$$X_4^2 = \{1, x, x^2, x^3\} \subseteq V_4^2$$

is a set generator of the set linear algebra $V_4^2$ over $Z^o$.

$$X_5^2 = \left\{ \begin{bmatrix} 1 & 0 & 0 & 0 \\ 1 & 1 & 0 & 0 \\ 1 & 1 & 1 & 0 \\ 1 & 1 & 1 & 1 \end{bmatrix} \right\} \subseteq V_5^2$$

is the set generator of the set linear algebra $V_5^2$ over $Z^o$.

Thus $X_2 = (X_1^2, X_2^2, X_3^2, X_4^2, X_5^2)$ is the special set generator of the special set linear algebra over $Z^o$. Now $|X_2| = \{3, 1, 2, 4, 1\}$ is the set dimension of $V_2$ as a special set linear algebra over $Z^o$. We see $V = V_1 \cup V_2$ as a special set linear algebra is special set bigenerated by

$$\begin{align} X &= X_1 \cup X_2 \\ &= \{X_1^1, X_2^1, X_3^1, X_4^1\} \cup \{X_1^2, X_2^2, X_3^2, X_4^2, X_5^2\} \\ &\subseteq V_1 \cup V_2 \end{align}$$

and the special set bicardinality of $X = |X_1| \cup |X_2| = (3, 4, 2, \infty) \cup (3, 1, 2, 4, 1)\}$.

Now we will find the special set bigenerator of the special set vector bispace $V = V_1 \cup V_2$, i.e.; we are treating the same $V = V_1 \cup V_2$ now only as a special set vector bispace over $Z^o = S$. Now the biset which is the special bigenerator is given by $X = X_1 \cup X_2$ where $X_1 = (X_1^1, X_2^1, X_3^1, X_4^1, X_5^1)$ with $X_1^1 = \{$An infinite set of upper singular $2 \times 2$ matrices with entries from $Z^o\} \subseteq V_1^1$; $X_2^1 = \{$An infinite set of 4-tuples with entries from $Z^o\} \subseteq V_2^1$, $X_3^1 = \{$an infinite set of pairs with entries from $Z^o\} \subseteq V_3^1$; $X_4^1 = \{$An infinite set of $3 \times 1$ matrix and an infinite set of $2 \times 3$ matrix with entries from $Z^o\} \subseteq V_4^1$. $X_1 = (X_1^1, X_2^1, X_3^1, X_4^1)$ set generates $V_1 = ((V_1^1, V_2^1, V_3^1, V_4^1)$ infinitely as a special set vector space and the special dimension of $V_1$ is $|X_1| = (|X_1^1|, |X_2^1|, |X_3^1|, |X_4^1|) = (\infty, \infty, \infty, \infty)$. Thus we see the difference between the special set dimensions in case of the special set vector space $V_1$ and the same special set linear algebra $V_1$ defined over the same set special dimension of $V_1$ as a special set vector space over $S = Z^+ \cup \{0\}$, is $|X_1| = (\infty, \infty, \infty, \infty)$ and the special dimension of $V_1$ as a special set linear algebra over the set $S = Z^+ \cup \{0\}$ is $(3, 4, 2, \infty)$.

Now we obtain the special set dimension of the special set vector space $V_2 = (V_1^2, V_2^2, V_3^2, V_4^2, V_5^2)$ over the set $Z^o$. Let $X_2 = (X_1^2, X_2^2, X_3^2, X_4^2, X_5^2) \subseteq (V_1^2, V_2^2, V_3^2, V_4^2, V_5^2)$ where

$X_1^2 = \{$is an infinite set of triples with entries from $Z^o\} \subseteq V_1^2$;

$$X_2^2 = \left\{ \begin{bmatrix} 1 \\ 1 \\ 1 \\ 1 \\ 1 \end{bmatrix} \right\} \subseteq V_2^2,$$

$$X_3^2 = \left\{ \begin{bmatrix} 1 & 1 & 1 \\ 1 & 1 & 1 \\ 1 & 1 & 1 \end{bmatrix}, \begin{bmatrix} 1 & 1 & 1 & 1 & 1 \\ 1 & 1 & 1 & 1 & 1 \end{bmatrix} \right\} \subseteq V_3^2,$$

$X_4^2 = \{$an infinite set of polynomial in x with coefficients from $Z^o$ of degree less than or equal to three$\} \subseteq V_4^2$

and

$$X_5^2 = \left\{ \begin{bmatrix} 1 & 0 & 0 & 0 \\ 1 & 1 & 0 & 0 \\ 1 & 1 & 1 & 0 \\ 1 & 1 & 1 & 1 \end{bmatrix} \right\} \subseteq V_5^2.$$

Thus $X_2 = (X_1^2, X_2^2, X_3^2, X_4^2, X_5^2)$ is a special set generator of the special set vector space $V_2 = (V_1^2, V_2^2, V_3^2, V_4^2, V_5^2)$. Now the special dimension of $V_2$ as a special set vector space over $Z^o$ is $|X_2| = (\infty, 1, 2, \infty, 1) > (3, 1, 2, 4, 1)$ which is the special set dimension $V_2$ as a special set linear algebra over $Z^o$.

Thus we see as a special set linear algebra $V_2$ over $Z^o$ is finite dimensional whereas $V_2$ as a special set vector space over $Z^o$ is infinite dimension. Thus $V = V_1 \cup V_2$ which we utilized as a special set linear bialgebra over $Z^o$ is of special set bidimension $|X| = |X_1| \cup |X_2| = (3, 4, 2, \infty) \cup (3, 1, 2, 4, 1)$ but the same $V = V_1 \cup V_2$ treated as a special set vector bispace over the set $Z^o$ is of special set bidimension $|X| = |X_1| \cup |X_2| = (\infty, \infty, \infty, \infty) \cup (\infty, 1, 2, \infty, 1)$.

Now having defined the special set bidimension we proceed to define special set linear bitransformation and special set linear bioperator.

**DEFINITION 2.2.8**: *Let*
$$V = V_1 \cup V_2 = (V_1^1, V_2^1, \ldots, V_{n_1}^1) \cup (V_1^2, V_2^2, \ldots, V_{n_2}^2)$$
*be a special set vector bispace over the set S. Let $W = W_1 \cup W_2$ $= (W_1^1, W_2^1, \ldots, W_{n_1}^1) \cup (W_1^2, W_2^2, \ldots W_{n_2}^2)$ be a special set vector bispace over the set S. Let*
$$T = T_1 \cup T_2 = ((T_1^1, T_2^1, \ldots, T_{n_1}^1) \cup (T_1^2, T_2^2, \ldots, T_{n_2}^2)$$
*be a special set bimap from $V_1 \cup V_2 \to W_1 \cup W_2$, i.e., $(T_1^1, T_2^1, \ldots, T_{n_1}^1): V_1 \to W_1$ and $(T_1^2, T_2^2, \ldots, T_{n_2}^2): V_2 \to W_2$ such that $T_i^1 : V_i^1 - W_i^1$ with $T_i^1(\alpha v_i^1) = \alpha T_i^1(v_i^1)$ for all $\alpha \in S$ and $v_i^1 \in V_i^1$ true for $i = 1, 2, \ldots, n_1$, and $T_i^2 : V_i^2 \to W_i^2$ such that $T_i^2(\alpha v_i^2) = \alpha T_i^2(v_i^2)$ true for $i = 1, 2, \ldots, n_2$ for $v_i^2 \in V_i^2$, and for all $\alpha \in S$. We call*
$$T = T_1 \cup T_2 = (T_1^1, T_2^1, \ldots, T_{n_1}^1) \cup (T_1^2, T_2^2, \ldots, T_{n_2}^2) :$$
$$V = (V_1^1, \ldots, V_{n_1}^1) \cup (V_1^2, \ldots, V_{n_2}^2)$$
$$\to (W_1^1, W_2^1, \ldots, W_{n_1}^1) \cup (W_1^2, W_2^2, \ldots, W_{n_2}^2) = W$$
*a special set linear bitransformation of the special set bivector space $V = V_1 \cup V_2$ into the special set bivector space $W = W_1 \cup W_2$ both V and W defined over the same set S.*

*The following observations are important.*

1. *Both $V = V_1 \cup V_2$ and $W = W_1 \cup W_2$ must be defined over the same set S.*
2. *Clearly both $V_1$ and $W_1$ must have same number of set vector spaces; likewise $V_2$ and $W_2$ must also have the same number of set vector spaces.*

3. If $SHom_S(V, W) = SHom_S(V_1, W_1) \cup SHom_S(V_2, W_2) =$
   $\{Hom_S(V_1^1, W_1^1), Hom_S(V_2^1, W_2^1), \ldots, Hom_S(V_{n_1}^1, W_{n_1}^1)\} \cup$
   $\{Hom_S(V_1^2, W_1^2), Hom_S(V_2^2, W_2^2), \ldots, Hom_S(V_{n_2}^2, W_{n_2}^2)\}$
   then, $SHom_S(V, W)$ is also a special set vector bispace over the set S.

We illustrate this by the following examples.

*Example 2.2.12*: Let
$$V = V_1 \cup V_2 = \{V_1^1, V_2^1, V_3^1\} \cup \{V_1^2, V_2^2, V_3^2, V_4^2\}$$
and
$$W = W_1 \cup W_2 = \{W_1^1, W_2^1, W_3^1\} \cup \{W_1^2, W_2^2, W_3^2, W_4^2\}$$
be two special set vector bispaces over the set $Z^o = Z^+ \cup \{0\}$ where $V_1 = (V_1^1, V_2^1, V_3^1)$ is such that $V_1^1 = \{Z^o \times Z^o \times Z^o\}$,

$$V_2^1 = \left\{ \begin{bmatrix} a & b \\ c & d \end{bmatrix} \middle| a, b, c, d \in Z^o \right\}$$

and

$$V_3^1 = \left\{ \begin{bmatrix} a & a & a & a \\ a & a & a & a \end{bmatrix} \middle| a \in Z^o \right\}$$

is a special set vector space over $Z^o$. Now $V_2 = (V_1^2, V_2^2, V_3^2, V_4^2)$ is such that $V_1^2 = Z^o \times Z^o \times Z^o \times Z^o$,

$$V_2^2 = \left\{ \begin{pmatrix} a_1 & a_2 & a_3 \\ a_4 & a_5 & a_6 \end{pmatrix} \middle| a_i \in Z^o, 1 \leq i \leq 6 \right\},$$

$$V_3^2 = \left\{ \begin{bmatrix} a_1 & a_2 \\ a_3 & a_4 \\ a_5 & a_6 \\ a_7 & a_8 \end{bmatrix} \middle| a_i \in Z^o; 1 \leq i \leq 8 \right\}$$

and $V_4^2$ = {all 4 × 4 upper triangular matrices with entries from $Z^o$} is a special set vector space over the set $Z^o$. Define $W = W_1 \cup W_2$ where $W_1 = (W_1^1, W_2^1, W_3^1)$ is such that

$$W_1^1 = \left\{ \begin{pmatrix} a & b \\ c & 0 \end{pmatrix} \middle| a, b, c \in Z^o \right\},$$

$$W_2^1 = \{Z^o \times Z^o \times Z^o \times Z^o\}$$

and

$$W_3^1 = \left\{ \begin{bmatrix} a & a \\ a & a \\ a & a \\ a & a \end{bmatrix} \middle| a \in Z^o \right\}$$

is also a special set vector space over the set $Z^o$. Now $W_2 = (W_1^2, W_2^2, W_3^2, W_4^2)$ is such that

$$W_1^2 = \left\{ \begin{pmatrix} a & b \\ c & d \end{pmatrix} \middle| a, b, c, d \in Z^o \right\},$$

$$W_2^2 = \left\{ \begin{bmatrix} a_1 & a_2 \\ a_3 & a_4 \\ a_5 & a_6 \end{bmatrix} \middle| a_i \in Z^o, 1 \le i \le 6 \right\},$$

$$W_3^2 = \left\{ \begin{bmatrix} a_1 & a_2 & a_3 & a_4 \\ 0 & 0 & a_5 & a_6 \\ 0 & 0 & a_8 & a_7 \end{bmatrix} \text{ such that } a_i \in Z^o; 1 \le i \le 8 \right\}$$

and $W_4^2$ = {all 4 × 4 matrices with entries from the set $Z^o$}. Clearly $W_2 = (W_1^2, W_2^2, W_3^2, W_4^2)$ is a special set vector space over the set $Z^o$. Thus $V = V_1 \cup V_2$ and $W = W_1 \cup W_2$ are special set vector bispaces over the set $Z^o$. Now define the special set bimap

$$T = T_1 \cup T_2 = (T_1^1, T_2^1, T_3^1) \cup (T_1^2, T_2^2, T_3^2, T_4^2):$$
$$V = (V_1^1, V_2^1, V_3^1) \cup (V_1^2, V_2^2, V_3^2, V_4^2) \to$$
$$W = (W_1^1, W_2^1, W_3^1) \cup (W_1^2, W_2^2, W_3^2, W_4^2)$$

such that
$$T_1^1 : V_1^1 \to W_1^1$$
$$T_2^1 : V_2^1 \to W_2^1$$
$$T_3^1 : V_3^1 \to W_3^1$$

and
$$T_1^2 : V_1^2 \to W_1^2$$
$$T_2^2 : V_2^2 \to W_2^2$$
$$T_3^2 : V_3^2 \to W_3^2$$

and
$$T_4^2 : V_4^2 \to W_4^2$$

given by

$$T_1^1 (a\ b\ c) = \begin{pmatrix} a & b \\ c & 0 \end{pmatrix},$$

$$T_2^1 \begin{pmatrix} a & b \\ c & d \end{pmatrix} = (a\ b\ c\ d)$$

and

$$T_3^1 \begin{pmatrix} a & a & a & a \\ a & a & a & a \end{pmatrix} = \begin{pmatrix} a & a \\ a & a \\ a & a \\ a & a \end{pmatrix}.$$

Now

$$T_1^2 (a\ b\ c\ d) = \begin{pmatrix} a & b \\ c & d \end{pmatrix},$$

$$T_2^2 \begin{pmatrix} a_1 & a_2 & a_3 \\ a_4 & a_5 & a_6 \end{pmatrix} = \begin{pmatrix} a_1 & a_2 \\ a_3 & a_4 \\ a_5 & a_6 \end{pmatrix},$$

$$T_3^2 \begin{pmatrix} a_1 & a_2 \\ a_3 & a_4 \\ a_5 & a_6 \\ a_7 & a_8 \end{pmatrix} = \begin{pmatrix} a_1 & a_2 & a_3 & a_4 \\ 0 & 0 & a_5 & a_6 \\ 0 & 0 & a_7 & a_8 \end{pmatrix}$$

and

$$T_4^2 \begin{pmatrix} a & b & c & d \\ 0 & e & f & g \\ 0 & 0 & i & j \\ 0 & 0 & 0 & k \end{pmatrix} = \begin{pmatrix} a & b & c & d \\ b & e & f & g \\ c & f & c & j \\ d & g & j & k \end{pmatrix}.$$

Clearly $T = T_1 \cup T_2$ is the special set linear bitransformation of $V = V_1 \cup V_2$ into $W = W_1 \cup W_2$.

*Example 2.2.13*: Let $V = V_1 \cup V_2$ and $W = W_1 \cup W_2$ be two special set vector bispaces over the set $S = \{0, 1\}$. Here

$$V = V_1 \cup V_2 = \{V_1^1, V_2^1\} \cup \{V_1^2, V_2^2, V_3^2\}$$

and $W = W_1 \cup W_2 = \{W_1^1, W_2^1\} \cup \{W_1^2, W_2^2, W_3^2\}$

with

$$V_1^1 = \{Z_2 \times Z_2 \times Z_2 \times Z_2\},$$

$$V_2^1 = \left\{ \begin{pmatrix} a & b & c \\ d & e & f \end{pmatrix} \middle| a, b, c, d, e, f \in Z_2 \right\},$$

$$V_1^2 = \left\{ [a_{ij}]_{4 \times 4}, a_{ij} \in Z_4 \right\}, \quad V_2^2 = \{Z_7\} \text{ and } V_3^2 \setminus \{Z_8\}.$$

$W_1 = (W_1^1, W_2^1)$ where

$$W_1^1 = \left\{ \begin{pmatrix} a & b \\ c & d \end{pmatrix} \middle| a, b, c, d \in Z_2 \right\}$$

and

$$W_2^1 = \left\{ \begin{bmatrix} a & b \\ c & d \\ e & f \end{bmatrix} \middle| a, b, c, d, e, f \in Z_2 \right\}.$$

$W_2 = \{W_1^2, W_2^2, W_3^2\}$ where

$$W_1^2 = \left\{ \begin{pmatrix} a_1 & a_2 & a_3 & a_4 \\ a_5 & a_6 & a_7 & a_8 \end{pmatrix} \middle| a_i \in Z_4; 1 \le i \le 8 \right\},$$

$$W_2^2 = \{(Z_7)\} \text{ and } W_3^2 = \{0, 1, 2, \ldots, 23 = Z_{24}\}.$$

Define
$$T = T_1 \cup T_2 = (T_1^1, T_2^1) \cup (T_1^2, T_2^2, T_3^2)$$

as $T_1^1: V_1^1 \to W_1^1$ where

$$T_1^1(a\ b\ c\ d) = \begin{pmatrix} a & b \\ c & d \end{pmatrix},$$

$T_2^1: V_2^1 \to W_2^1$ given by

$$T_2^1 \begin{pmatrix} a & b & c \\ d & e & f \end{pmatrix} = \begin{bmatrix} a & b \\ c & d \\ e & f \end{bmatrix}.$$

Now $T_1^2: V_1^2 \to W_1^2$;

$$T_1^2 \begin{pmatrix} a_1 & a_2 & a_3 & a_4 \\ a_5 & a_6 & a_7 & a_8 \\ a_9 & a_{10} & a_{11} & a_{12} \\ a_{13} & a_{14} & a_{15} & a_{16} \end{pmatrix} = \begin{pmatrix} a_1 & a_2 & a_3 & a_4 \\ a_5 & a_6 & a_7 & a_8 \end{pmatrix},$$

$T_2^2: V_2^2 \to W_2^2$; $T_2^2(a) = a$,
$T_3^2: V_3^2 \to W_3^2$; $T_3^2(a) = a \in Z_{24}$

i.e., $\qquad T_3^2(0) = 0,\ T_3^2(1) = (3),$
$$T_3^2(2) = 6,\ T_3^2(3) = 9,$$
$$T_3^2(4) = 12,\quad T_3^2(5) = 15,$$
$$T_3^2(6) = 18 \text{ and } T_3^2(7) = 21.$$

Thus
$$T = T_1 \cup T_2 : V \to W = V_1 \cup V_2 \to W_1 \cup W_2 \text{ is}$$
$$(T_1^1, T_2^1) \cup (T_1^2, T_2^2, T_3^2):$$
$$(V_1^1, V_2^1) \cup (V_1^2, V_2^2, V_3^2) \to (W_1^1, W_2^1) \cup (W_1^2, W_2^2, W_3^2)$$

is a special set linear bitransformation of V into W.

Now are proceed onto define the notion of special set linear bioperators on V.

**DEFINITION 2.2.9**: *Let*

$$V = V_1 \cup V_2 = (V_1^1, V_2^1, \ldots, V_{n_1}^1) \cup (V_1^2, V_2^2, \ldots, V_{n_2}^2)$$

*is a special set vector bispace over the set S. Let $T = T_1 \cup T_2 = (T_1^1, T_2^1, \ldots, T_{n_1}^1) \cup (T_1^2, T_2^2, \ldots, T_{n_2}^2) : V_1 \cup V_2 = (V_1^1, V_2^1, \ldots, V_{n_1}^1) \cup (V_1^2, V_2^2, \ldots, V_{n_2}^2) \to V_1 \cup V_2$ such that $T_i^1 : V_i^1 \to V_i^1$, $1 \leq i \leq n_1$ and $T_i^2 : V_i^2 \to V_i^2$, $1 \leq i \leq n_2$. If $T_1 \cup T_2 = T$ is a special set linear bitransformation of V into V then we call $T = T_1 \cup T_2$ to be the special set linear bioperator on V; i.e., if in the definition of a special set linear bitransformation the domain space coincides with the range space then we call T to be a special set linear bioperator on V.*

We illustrate this by a simple example.

***Example 2.2.14***: Let $V = V_1 \cup V_2$ where $V_1 = (V_1^1, V_2^1, V_3^1, V_4^1)$ and $V_2 = (V_1^2, V_2^2, V_3^2, V_4^2, V_5^2)$ be a special set vector bispace over the set $Z^o = Z^+ \cup \{0\}$. Here $V_1^1 = Z^o \times Z^o \times Z^o \times Z^o$,

$$V_2^1 = \left\{ \begin{pmatrix} a & b \\ c & d \end{pmatrix} \middle| a, b, c, d \in Z^o \right\},$$

$$V_3^1 = \left\{ \begin{bmatrix} a_1 & a_2 & a_3 & a_4 & a_5 \\ a_6 & a_7 & a_8 & a_9 & a_{10} \end{bmatrix}, \begin{bmatrix} a_1 \\ a_2 \end{bmatrix} \middle| a_i \in Z^o; 1 \leq i \leq 10 \right\},$$

and $V_4^1 = $ {all polynomials in the variable x with coefficients from $Z^o$ of degree less than or equal to 4}. Now

$$V_1^2 = \left\{ \begin{pmatrix} a & b & c \\ 0 & d & e \\ 0 & 0 & f \end{pmatrix} \middle| a, b, c, d, e, f \in Z^o \right\},$$

$$V_2^2 = \{Z^o \times Z^o \times Z^o \times Z^o \times Z^o\},$$

$V_3^2 = \{$all polynomial in x of degree less than or equal to 6 with coefficients from $Z^o\}$,

$$V_4^2 = \left\{ (a_1\ a_2\ a_3\ a_4\ a_5\ a_6),\ \begin{bmatrix} a_1 & a_2 \\ a_3 & a_4 \\ a_5 & a_6 \end{bmatrix} \middle| a_i \in Z^o,\ 1 \le i \le 6 \right\}$$

and
$$V_5^2 = \{4 \times 4 \text{ matrices with entries from } Z^o\}.$$

$V = V_1 \cup V_2$ is a special set vector bispace over the set $Z^o$. Define $T = T_1 \cup T_2 : V \to V$ by

$$T_1 \cup T_2 = (T_1^1, T_2^1, T_3^1, T_4^1) \cup (T_1^2, T_2^2, T_3^2, T_4^2, T_5^2) :$$
$$V = V_1 \cup V_2$$
$$= (V_1^1, V_2^1, V_3^1, V_4^1) \cup (V_1^2, V_2^2, V_3^2, V_4^2, V_5^2)$$
$$\to (V_1^1, V_2^1, V_3^1, V_4^1) \cup (V_1^2, V_2^2, V_3^2, V_4^2, V_5^2)$$

as $\quad T_i^1 : V_i^1 \to V_i^1,\ 1 \le i \le 4$
and $\quad T_i^2 : V_i^2 \to V_i^2,\ 1 \le i \le 5.$
$T_1^1 : V_1^1 \to V_1^1$
$$T_1^1 (a\ b\ c\ d) = (b\ c\ d\ a),$$
$T_2^1 : V_2^1 \to V_2^1$
defined by
$$T_2^1 \begin{pmatrix} a & b \\ c & d \end{pmatrix} = \begin{pmatrix} a & c \\ b & d \end{pmatrix},$$
$T_3^1 : V_3^1 \to V_3^1$

defined by

$$T_3^1 \begin{pmatrix} a_1 & a_2 & a_3 & a_4 & a_5 \\ a_6 & a_7 & a_8 & a_9 & a_{10} \end{pmatrix} = \begin{pmatrix} a_6 & a_7 & a_8 & a_9 & a_{10} \\ a_1 & a_2 & a_3 & a_4 & a_5 \end{pmatrix}$$

and

$$T_3^1 \begin{bmatrix} a_1 \\ a_2 \end{bmatrix} = \begin{bmatrix} a_2 \\ a_1 \end{bmatrix}.$$

$T_4^1 : V_4^1 \to V_4^1$
defined by

$$T_4^1 (a_0 + a_1x + a_2x^2 + a_3x^3 + a_4x^4) = (a_1x + a_3 x^3).$$

$T_1 = (T_1^1, T_2^1, T_3^1, T_4^1)$ is a special set linear operator on $V_1 = (V_1^1, V_2^1, V_3^1, V_4^1)$. Define $T_2 = (T_1^2, T_2^2, T_3^2, T_4^2, T_5^2)$ from $V_2 \to V_2 = (V_1^2, V_2^2, V_3^2, V_4^2, V_5^2)$ by

$T_1^2 : V_1^2 \to V_1^2$ ;

$$T_2^1 \begin{pmatrix} a & b & c \\ 0 & d & e \\ 0 & 0 & f \end{pmatrix} = \begin{pmatrix} a & d & e \\ 0 & b & c \\ 0 & 0 & f \end{pmatrix},$$

$T_2^2 : V_2^2 \to V_2^2$ by

$$T_2^2 (a\ b\ c\ d\ e) = (a\ c\ b\ e\ d),$$

$T_3^2 : V_3^2 \to V_3^2$ is defined by

$$T_3^2 (a_0 + a_1x + a_2x^2 + a_3x^3 + a_4 x^4 + a_5 x^5 + a_6x^6)$$
$$= a_1x^6 + a_1x^5 + a_2x^4 + a_3x^3 + a_4x^2 + a_5x + a_6.$$

$T_4^2 : V_4^2 \to V_4^2$ is given by

$$T_4^2 (a_1\ a_2\ a_3\ a_4\ a_5\ a_6) = \begin{bmatrix} a_1 & a_2 \\ a_3 & a_4 \\ a_5 & a_6 \end{bmatrix}$$

and

$$T_4^2 \begin{bmatrix} a_1 & a_2 \\ a_3 & a_4 \\ a_5 & a_6 \end{bmatrix} = (a_1\ a_2\ a_3\ a_4\ a_5\ a_6).$$

Finally $T_5^2 : V_5^2 \to V_5^2$ is defined as

$$T_5^2 \begin{pmatrix} a_{11} & a_{12} & a_{13} & a_{14} \\ a_{21} & a_{22} & a_{23} & a_{24} \\ a_{31} & a_{32} & a_{33} & a_{34} \\ a_{41} & a_{42} & a_{43} & a_{44} \end{pmatrix} = \begin{pmatrix} a_{11} & a_{21} & a_{31} & a_{41} \\ a_{12} & a_{22} & a_{32} & a_{42} \\ a_{13} & a_{23} & a_{33} & a_{43} \\ a_{14} & a_{24} & a_{34} & a_{44} \end{pmatrix}.$$

Clearly $T_2 = (T_1^2, T_2^2, T_3^2, T_4^2, T_5^2)$ is a special set linear operator on $V_2$. Thus $T_1 \cup T_2 = (T_1^1, T_2^1, T_3^1, T_4^1) \cup (T_1^2, T_2^2, T_3^2, T_4^2, T_5^2)$ is a special set linear bioperator on $V = V_1 \cup V_2$.

Now we can also define as in case of special set vector spaces the notion of special set pseudo linear bioperator on $V = V_1 \cup V_2$.

**DEFINITION 2.2.10**: *Let*

$$V = V_1 \cup V_2 = (V_1^1, V_2^1, \ldots, V_{n_1}^1) \cup (V_1^2, V_2^2, \ldots, V_{n_2}^2)$$

*be a special set vector bispace over a set S. Let*

$$T = T_1 \cup T_2 = (T_1^1, T_2^1, \ldots, T_{n_1}^1) \cup (T_1^2, T_2^2, \ldots, T_{n_2}^2)$$

*be a special set vector transformation such $T_i^1 : V_i^1 \to V_j^1$, $(i \neq j)$; $1 \leq i, j \leq n_1$ for some i and j.*

*Similarly $T_i^2 : V_i^2 \to V_j^2$, $i \neq j$ for some i and j; $1 \leq i, j \leq n_2$ we call this $T = T_1 \cup T_2$ to be a special set pseudo linear bioperator on V. Clearly*

$$SHom_S(V,V) = \{Hom_S(V_{i_1}^1, V_{j_1}^1) \ldots Hom_S(V_{in_1}^1, V_{jn_1}^1)\}$$

$$\cup \{Hom_S(V_{i_2}^2, V_{j_2}^2) \cup \ldots \cup Hom_S(V_{in_2}^2, V_{jn_2}^2)\}$$

*is only a special set vector bispace over the set S.*

We illustrate this by an example.

*Example 2.2.15:* Let
$$V = V_1 \cup V_2 = \{V_1^1, V_2^1, V_3^1, V_4^1\} \cup \{V_1^2, V_2^2, V_3^2\}$$
where
$$V_1^1 = \left\{ \begin{pmatrix} a & b \\ c & d \end{pmatrix} \middle| a,b,c,d \in Z^+ \cup \{0\} \right\},$$
$$V_2^1 = \{(Z^+ \cup \{0\}) \times (Z^+ \cup \{0\}) \times Z^+ \cup \{0\} \times (Z^+ \cup \{0\})\},$$
$$V_3^1 = \{(Z^+ \cup \{0\})[x], \text{ all polynomials of degree less than or equal to 5}\}$$
and
$$V_4^1 = \left\{ \begin{pmatrix} a_1 & a_2 & a_3 \\ a_4 & a_5 & a_6 \end{pmatrix} \right\} \text{ such that } a_i \in Z^+ \cup \{0\}; 1 \le i \le 6\}.$$

$$V_1^2 = \left\{ \begin{pmatrix} a & b & c & d \\ 0 & e & f & g \\ 0 & 0 & h & i \\ 0 & 0 & 0 & k \end{pmatrix} \middle| a,b,c,d,e,f,g,h,i,k \in Z^+ \cup \{0\} \right\},$$

$$V_2^2 = \{(Z^+ \cup \{0\})[x]; \text{ all polynomials of degree less than or equal to 9}\}$$
and
$$V_3^2 = \left\{ \begin{pmatrix} a_1 & a_2 & a_3 \\ a_4 & a_5 & a_6 \\ a_7 & a_8 & a_9 \\ a_{10} & 0 & 0 \end{pmatrix} \middle| a_i \in Z^+ \cup \{0\}; 1 \le i \le 9 \right\};$$

thus $V = V_1 \cup V_2 = (V_1^1, V_2^1, V_3^1, V_4^1) \cup (V_1^2, V_2^2, V_3^2)$ is a special set vector bispace over the set $S = Z^+ \cup \{0\}$. Now define a special set linear bitransformation T from V to V by
$$T = T_1 \cup T_2 = \{T_1^1, T_2^1, T_3^1, T_4^1\} \cup \{T_1^2, T_2^2, T_3^2\};$$
$$V = V_1 \cup V_2 = (V_1^1, V_2^1, V_3^1, V_4^1) \cup (V_1^2, V_2^2, V_3^2) \to$$
$$V = V_1 \cup V_2 = (V_1^1, V_2^1, V_3^1, V_4^1) \cup (V_1^2, V_2^2, V_3^2)$$

as follows :

$T_1^1 : V_1^1 \to V_2^1$ defined by

$$T_1^1 \begin{pmatrix} a & b \\ c & d \end{pmatrix} = (a, b, c, d),$$

$T_2^1 : V_2^1 \to V_1^1$ by

$$T_2^1 (a\,b\,c\,d) = \begin{pmatrix} a & b \\ c & d \end{pmatrix},$$

$T_3^1 : V_3^1 \to V_4^1$ by

$$T_3^1 (a_0 + a_1 x + a_2 x^2 + a_3 x^3 + a_4 x^4 + a_5 x^5) = \begin{bmatrix} a_0 & a_1 & a_2 \\ a_3 & a_4 & a_5 \end{bmatrix}$$

and $T_4^1 : V_4^1 \to V_3^1$ by

$$T_4^1 \begin{bmatrix} a_1 & a_2 & a_3 \\ a_4 & a_5 & a_6 \end{bmatrix} = a_1 x^5 + a_2 x^4 + a_3 x^3 + a_4 x^2 + a_5 x + a_6.$$

Thus $T_1 = (T_1^1, T_2^1, T_3^1, T_4^1) : V_1 \to V_1$ is only a pseudo set linear operator on V. Now define

$T_1^2 : V_1^2 \to V_2^2$ by

$$T_1^2 \begin{pmatrix} a & b & c & d \\ 0 & e & f & g \\ 0 & 0 & h & i \\ 0 & 0 & 0 & j \end{pmatrix} = (a + bx + cx^2 + dx^3 \\ + ex^4 + fx^5 + gx^6 + hx^7 + ix^8 + jx^9),$$

$T_2^2 : V_2^2 \to V_3^2$ by

$$T_2^2 = (a_0 + a_1 x + a_2 x^2 + \ldots + a_9 x^9) = \begin{pmatrix} a_0 & a_1 & a_2 \\ a_3 & a_4 & a_5 \\ a_6 & a_7 & a_8 \\ a_9 & 0 & 0 \end{pmatrix}$$

and

$T_3^2 : V_3^2 \to V_2^2$ by

$$T_3^2 \begin{pmatrix} a_1 & a_2 & a_3 \\ a_4 & a_5 & a_6 \\ a_7 & a_8 & a_9 \\ a_{10} & 0 & 0 \end{pmatrix} = a_1 x^9 + a_2 x^8 + a_3 x^7 + x_4 x^6$$
$$+ a_5 x^6 + a_5 x^5 + a_6 x^4 + a_7 x^3 + a_8 x^2 + a_9 x + a_{10}.$$

Clearly
$T_2 = (T_1^2, T_2^2, T_3^2) : V_2 = (V_1^2, V_2^2, V_3^2) \to V_2 = (V_1^2, V_2^2, V_3^2)$

is a pseudo set linear operator on $V_2$.

Thus $T = T_1 \cup T_2 = (T_1^1, T_2^1, T_3^1, T_4^1) \cup (T_1^2, T_2^2, T_3^2)$ is a pseudo special set linear bioperator from $V = V_1 \cup V_2$ into $V = V_1 \cup V_2$.

## 2.3 Special Set Vector n-spaces

Now we proceed onto generalize the notion of special set vector bispaces to special set vector n-spaces; $n \geq 3$.

**DEFINITION 2.3.1**: *Let $V = (V_1 \cup V_2 \cup ... \cup V_n)$ where each $V_i = (V_1^i, V_2^i, ..., V_{n_i}^i)$ is a special set vector space over the same set S; i = 1, 2, ..., n and $V_i \neq V_j$; if $i \neq j$, i.e.; $V_i$'s are distinct, $V_i \not\subseteq V_j$ and $V_j \not\subseteq V_i$, $1 \leq i, j \leq n$. Then we call $V = (V_1 \cup V_2 \cup ... \cup V_n) = (V_1^1, V_2^1, ..., V_{n_1}^1) \cup (V_1^2, V_2^2, ..., V_{n_2}^2) \cup ... \cup (V_1^n, V_2^n, ..., V_{n_n}^n)$ to be a special set vector n-space over the set S. If n = 2 then $V = V_1 \cup V_2$ is the special set vector bispace that is why we assume $n \geq 3$; when n = 3 we get the special set vector trispace.*

We illustrate this by some simple examples.

*Example 2.3.1.*: Let
$$\begin{aligned} V &= V_1 \cup V_2 \cup V_3 \cup V_4 \\ &= (V_1^1, V_2^1, V_3^1) \cup (V_1^2, V_2^2, V_3^2, V_4^2) \cup \\ & \quad (V_1^3, V_2^3, V_3^3, V_4^3) \cup (V_1^4, V_2^4, V_3^4) \end{aligned}$$

be a special set vector 4-space over set $S = Z^+ \cup \{0\}$. Here $V_1 = (V_1^1, V_2^1, V_3^1)$ where

$$V_1^1 = \left\{ \begin{pmatrix} a & b \\ c & d \end{pmatrix} \middle| a,b,c,d \in S \right\},$$

$$V_2^1 = \{S \times S \times S\}$$

and

$$V_3^1 = \{(a\ a\ a\ a\ a\ a) \mid a \in S\}$$

is a special set vector space over the set S. Now

$$V_2 = (V_1^2, V_2^2, V_3^2, V_4^2)$$

where

$$V_1^2 = (S \times S \times S \times S \times S)$$

$$V_2^2 = \left\{ \begin{pmatrix} a_1 & a_2 & a_3 \\ a_4 & a_5 & a_6 \\ a_7 & a_8 & a_9 \end{pmatrix} \middle| a_i \in S, 1 \leq i \leq 9 \right\},$$

$$V_3^2 = \{(a\ a\ a\ a) \mid a \in S\}$$

and

$$V_4^2 = \left\{ \begin{bmatrix} a_1 & a_2 \\ a_3 & a_4 \\ a_5 & a_6 \\ a_7 & a_8 \end{bmatrix} \middle| a_i \in S; 1 \leq i \leq 8 \right\};$$

$V_2 = (V_1^2, V_2^2, V_3^2, V_4^2)$ is a special set vector space over the set $S = Z^+ \cup \{0\}$. Now $V_3 = (V_1^3, V_2^3, V_3^3, V_4^3)$ where $V_1^3 = \{4 \times 4$ upper triangular matrices with entries from $S\}$, $V_2^3 = \{(Z^o \cup \{0\})[x]$ all polynomials is the variable x with coefficients from S of degree less than or equal to 4$\}$, $V_3^3 = \{7 \times 7$ diagonal matrices with entries from S$\}$ and

$$V_4^3 = \left\{ \begin{pmatrix} a & b & c \\ d & e & f \\ g & h & i \\ j & k & l \end{pmatrix} \middle| a,b,c,d,e,f,g,h,i,j,k,l \in S \right\}.$$

Thus $V_3$ is a special set vector space over the set S. Now consider $V_4 = (V_1^4, V_2^4, V_3^4)$

$$V_1^4 = \left\{ \begin{pmatrix} a & b \\ c & d \end{pmatrix}, \begin{pmatrix} a & a & a \\ a & a & a \end{pmatrix} \middle| a,b,c,d \in Z^o = Z^+ \cup \{0\} \right\},$$

$$V_2^4 = \left\{ (a\ a\ a\ a), \begin{bmatrix} a \\ a \\ a \\ a \\ a \end{bmatrix} \middle| a \in Z^+ \cup \{0\} = S \right\}$$

and

$V_3^4 = \{$all $5 \times 5$ lower triangular matrices with entries from $S = Z^o = Z^+ \cup \{0\}\}$.

$V_4 = (V_1^4, V_2^4, V_3^4)$ is a special set vector space over the set S. Thus $V = V_1 \cup V_2 \cup V_3 \cup V_4$ is a special set vector 4-space over the set S.

We give another example of a special set n-vector space over a set S $n \geq 3$.

***Example 2.3.2:*** Let $V = V_1 \cup V_2 \cup V_3 = (V_1^1, V_2^1, V_3^1, V_4^1, V_5^1) \cup (V_1^2, V_2^2, V_3^2, V_4^2) \cup (V_1^3, V_2^3, V_3^3)$ be a special vector trispace over the set $S = \{0\ 1\}$. $V_1 = (V_1^1, V_2^1, V_3^1, V_4^1, V_5^1)$ where $V_1^1 = \{$all $3 \times 3$ matrices with entries from S$\}$, $V_2^1 = \{S \times S \times S \times S \times S \times S\}$, $V_3^1 = \{$all polynomials of degree less than or equal to 6$\}$,

$$V_4^1 = \left\{ \begin{bmatrix} a & a & a \\ a & a & a \end{bmatrix}, \begin{bmatrix} a & a \\ a & a \\ a & a \\ a & a \end{bmatrix} \middle| a \in 3Z^+ \cup \{0\} \right\}$$

and $V_5^1$ = {all 7 × 7 upper triangular matrices with entries from S}. $V_1 = (V_1^1, V_2^1, V_3^1, V_4^1, V_5^1)$ is a special set vector space over the set S. $V_2 = (V_1^2, V_2^2, V_3^2, V_4^2)$ where $V_1^2 = \{S \times S \times S \times S\}$, $V_2^2$ = {all 8 × 8 matrices with entries from S}, $V_3^2$ = {(a a a a a), (a a), (a a a a a a a) | a ∈ $Z^+ \cup \{0\}$} and

$$V_4^2 = \left\{ \begin{bmatrix} a & a \\ a & a \\ a & a \\ a & a \end{bmatrix}, \begin{bmatrix} a & a & a & a \\ a & a & a & a \end{bmatrix} \middle| a \in Z^+ \cup \{0\} \right\}.$$

$V_2 = (V_1^2, V_2^2, V_3^2, V_4^2)$ is again a special set vector space over S. Finally $V_3 = (V_1^3, V_2^3, V_3^3)$ is such that $V_1^3 = \{S \times S\}$,

$$V_2^3 = \left\{ \begin{pmatrix} a & b \\ c & d \end{pmatrix} \middle| a, b, c, d \in S \right\}$$

and $V_3^3$ = {all 5 × 5 lower triangular matrices}. We see $V_3$ is again a special set vector space over the set S. Thus $V = V_1 \cup V_2 \cup V_3$ is a special set trivector space over the set S = {0, 1}.

Now we proceed onto define the notion of special set n-vector subspace.

**DEFINITION 2.3.2.** *Let*
$V = (V_1 \cup V_2 \cup V_3 \cup ... \cup V_n)$
$= (V_1^1, V_2^1, ..., V_{n_1}^1) \cup (V_1^2, V_2^2, ..., V_{n_2}^2) \cup ... \cup (V_1^n, V_2^n, ..., V_{n_n}^n)$
*be a special set n-vector space over the set S. Suppose*
$W = W_1 \cup W_2 \cup ... \cup W_n$
$= (W_1^1, W_2^1, ..., W_{n_1}^1) \cup (W_1^2, W_2^2, ..., W_{n_2}^2) \cup ... \cup (W_1^n, W_2^n, ..., W_{n_n}^n)$
$\subseteq V_1 \cup V_2 \cup ... \cup V_n$
*be such that* $(W_1^i, W_2^i, ..., W_{n_i}^i) \subseteq (V_1^i, V_2^i, ..., V_{n_i}^i)$,

*i.e.;* $W_t^i \subseteq V_t^i$ *is a set vector subspace of* $V_t^i$, $1 \leq t \leq n_i$ *so that* $W_i = (W_1^i, W_2^i, ..., W_{n_i}^i)$ *is a special set vector subspace of* $V_i = (V_1^i, ..., V_{n_i}^i)$; *true for i = 1, 2, ..., n then we define* $W = W_1 \cup ... \cup W_n$ *to be a special set vector n-subspace of V over the set S.*

We will illustrate this by some simple examples.

*Example 2.3.3.* Let
$$V = V_1 \cup V_2 \cup V_3 \cup V_4$$
$$= (V_1^1, V_2^1, V_3^1) \cup (V_1^2, V_2^2, V_3^2) \cup$$
$$(V_1^3, V_2^3, V_3^3, V_4^3) \cup (V_1^4, V_2^4, V_3^4, V_4^4)$$

be a special set 4-vector space over the set $S = \{0, 1\}$. Here $V_1 = (V_1^1, V_2^1, V_3^1)$ is such that $V_1^1 = (S \times S \times S \times S)$,

$$V_2^1 = \left\{ \begin{bmatrix} a & a \\ a & a \\ a & a \\ a & a \end{bmatrix}, \begin{bmatrix} a & a & a & a & a \\ a & a & a & a & a \end{bmatrix} \middle| a \in Z^+ \cup \{0\} \right\}$$

and $V_1^3 = \{$all $5 \times 5$ upper triangular matrices with entries from $Z^+ \cup \{0\}\}$; $V_1$ is a special set vector space over the set $S = \{0, 1\}$. $V_2 = (V_1^2, V_2^2, V_3^2)$ where

$$V_1^2 = \left\{ \begin{pmatrix} a & b & c \\ d & e & f \\ g & h & i \end{pmatrix} \middle| a, b, ..., i \in Z_2 = \{0, 1\} \right\},$$

$$V_2^2 = \{Z^o \times Z^o \times Z^o \mid Z^o = Z^+ \cup \{0\}\}$$

and

$$V_3^2 = \left\{ \begin{bmatrix} a_1 & a_2 & a_3 \\ a_4 & a_5 & a_6 \end{bmatrix}, \begin{bmatrix} a_1 & a_2 \\ a_3 & a_4 \\ a_5 & a_6 \\ a_7 & a_8 \end{bmatrix} \middle| a_i \in Z_2 = \{0, 1\} \right\}$$

is also a special set vector space over the set $S = \{0, 1\}$.
$V_3 = (V_1^3, V_2^3, V_3^3, V_4^3)$ where

$$V_1^3 = \{Z^o \times Z^o \times Z^o \times Z^o \mid Z^o = Z^+ \cup \{0\}\},$$

$$V_2^3 = \left\{ \begin{pmatrix} a_1 & a_2 & a_3 \\ a_4 & a_5 & a_6 \\ a_7 & a_8 & a_9 \end{pmatrix} \middle| a_i \in Z^+ \cup \{0\} \right\},$$

$$V_3^3 = \left\{ \begin{bmatrix} a & a & a & a \\ a & a & a & a \end{bmatrix}, (a\ a\ a\ a\ a\ a) \middle| a \in Z^o = Z^+ \cup \{0\} \right\}$$

and

$$V_4^3 = \left\{ \begin{bmatrix} a \\ a \\ a \\ a \\ a \end{bmatrix}, \begin{bmatrix} a & a & a & a & a \\ a & a & a & a & a \\ a & a & a & a & a \end{bmatrix} \middle| a \in Z^+ \cup \{0\} \right\}$$

is again a special set vector space over the set $S = \{0, 1\}$.
Finally $V_4 = \{V_1^4, V_2^4, V_3^4, V_4^4\}$ is such that

$$V_1^4 = \{S \times S \times S \times S \times S\};$$

$$V_2^4 = \left\{ \begin{pmatrix} a & b & c \\ 0 & d & e \\ 0 & 0 & f \end{pmatrix} \middle| a, b, c, d, e, f \in Z^+ \cup \{0\} \right\},$$

$$V_3^4 = \left\{ \begin{bmatrix} a & a & a & a \\ a & a & a & a \end{bmatrix}, \begin{bmatrix} a & a \\ a & a \\ a & a \\ a & a \\ a & a \\ a & a \end{bmatrix} \middle| a \in Z^+ \cup \{0\} \right\}$$

and

$$V_4^4 = \{(a\ a\ a\ a\ a\ a), \begin{bmatrix} a \\ a \\ a \end{bmatrix} \mid a \in Z^+ \cup \{0\}\}$$

is a special set vector space over the set $S = \{0,1\}$. Thus $V = V_1 \cup V_2 \cup V_3 \cup V_4$ is a special set vector 4-space over the set $S = \{0,1\}$. Take

$$\begin{aligned} W &= W_1 \cup W_2 \cup W_3 \cup W_4 \\ &= (W_1^1, W_2^1, W_3^1) \cup (W_1^2, W_2^2, W_3^2) \cup \\ &\quad (W_1^3, W_2^3, W_3^3, W_4^3) \cup (W_1^4, W_2^4, W_3^4, W_4^4) \\ &\subseteq (V_1^1, V_2^1, V_3^1) \cup (V_1^2, V_2^2, V_3^2) \cup \\ &\quad (V_1^3, V_2^3, V_3^3, V_4^3) \cup (V_1^4, V_2^4, V_3^4, V_4^4) \end{aligned}$$

where $W_1^1 \subseteq V_1^1$ and $W_1^1 = S \times S \times \{0\} \times \{0\} \subseteq V_1^1$,

$$W_2^1 = \left\{ \begin{bmatrix} a & a \\ a & a \\ a & a \\ a & a \end{bmatrix} \mid a \in Z^+ \cup \{0\} \right\} \subseteq V_2^1$$

and

$$W_3^1 = \{5 \times 5 \text{ upper triangular matrices with entries from } 3Z^+ \cup \{0\}\} \subseteq V_3^1.$$

Thus $(W_1^1, W_2^1, W_3^1) \subseteq (V_1^1, V_2^1, V_3^1)$ is a special set vector subspace of $V_1$. Now

$$W_1^2 = \left\{ \begin{pmatrix} a & a & a \\ a & a & a \\ a & a & a \end{pmatrix} \mid a \in \{0, 1\} \subseteq V_1^2 \right.,$$

$$W_2^2 = \{Z^o \times Z^o \times \{0\}\} \subseteq V_2^2$$

and

$$W_3^2 = \left\{ \begin{bmatrix} a_1 & a_2 \\ a_3 & a_4 \\ a_5 & a_6 \\ a_7 & a_8 \end{bmatrix} \middle| a_i \in Z_2 = \{0,1\} \ 1 \le i \le 8 \right\} \subseteq V_3^2.$$

Thus $W_2 = (W_1^2, W_2^2, W_3^2) \subseteq (V_1^2, V_2^2, V_3^2) = V_2$ is a special set vector subspace of $V_2$ over $5 = \{0, 1\}$. Take

$$W_1^3 = \{Z^o \times \{0\} \times Z^o \{0\}\} \subseteq V_1^3,$$

$$W_2^3 = \left\{ \begin{pmatrix} a & a & a \\ a & a & a \\ a & a & a \end{pmatrix} \middle| a \in Z^+ \cup \{0\} \right\} \subseteq V_2^3,$$

$$W_3^3 = \{(a\ a\ a\ a\ a\ a) \mid a \in Z^+ \cup \{0\}\} \subseteq V_3^3$$

and

$$W_4^3 = \left\{ \begin{bmatrix} a \\ a \\ a \\ a \\ a \\ a \end{bmatrix} \middle| a \in Z^+ \cup \{0\} \right\} \subseteq V_4^3$$

so $W_3 = (W_1^3, W_2^3, W_3^3, W_4^3) \subseteq (V_1^3, V_2^3, V_3^3, V_4^3) = V_3$ is a special set vector subspace of $V_3$ over the set $S = \{0, 1\}$. Finally $W_1^4 = \{S \times S \times S \times S \times \{0\}\} \subseteq V_1^4$ is a set vector subspace of $V_1^4$ over $S = \{0, 1\}$.

$$W_2^4 = \left\{ \begin{pmatrix} a & a & a \\ 0 & a & a \\ 0 & 0 & a \end{pmatrix} \middle| a \in Z^+ \cup \{0\} \right\} \subseteq V_2^4$$

over the set $S = \{0,1\}$.

$$W_3^4 = \left\{ \begin{bmatrix} a & a & a & a \\ a & a & a & a \end{bmatrix} \middle| \ a \in Z^+ \cup \{0\} \right\} \subseteq V_3^4$$

is a set vector subspace of $V_3^4$ over the set $S = \{0, 1\}$ and

$$W_4^4 = \left\{ \begin{bmatrix} a \\ a \\ a \end{bmatrix} \middle| \ a \in Z^+ \cup \{0\} \right\} \subseteq V_4^4$$

is a set vector subspace of $V_4^4$ over the set $S = \{0,1\}$. Thus
$$W_4 = (W_1^4, W_2^4, W_3^4, W_4^4) \subseteq (V_1^4, V_2^4, V_3^4, V_4^4) = V_4$$
is a special set vector subspace of $V_4$ over the set $S = \{0,1\}$.

Hence
$$\begin{aligned}
W &= W_1 \cup W_2 \cup W_3 \cup W_4 \\
&= (W_1^1, W_2^1, W_3^1) \cup (W_1^2, W_2^2, W_3^2) \cup (W_1^3, W_2^3, W_3^3, W_4^3) \\
&\quad \cup (W_1^4, W_2^4, W_3^4, W_4^4) \\
&\subseteq (V_1^1, V_2^1, V_3^1) \cup (V_1^2, V_2^2, V_3^2) \cup (V_1^3, V_2^3, V_3^3, V_4^3) \cup \\
&\quad (V_1^4, V_2^4, V_3^4, V_4^4) \\
&= V_1 \cup V_2 \cup V_3 \cup V_4 = V
\end{aligned}$$

is a special set vector 4-subspace of $V$ over set $S = \{0,1\}$.

***Example 2.3.4:*** Let
$$\begin{aligned}
V &= V_1 \cup V_2 \cup V_3 \\
&= \{V_1^1, V_2^1, V_3^1, V_4^1\} \cup \{V_1^2, V_2^2, V_3^2, V_4^2, V_5^2\} \cup \\
&\quad \{V_1^3, V_2^3, V_3^3, V_4^3\}
\end{aligned}$$
be a special set vector 3-space over the set $S = Z^+ \cup \{0\}$. Here $V_1 = \{V_1^1, V_2^1, V_3^1, V_4^1\}$ is such that $V_1^1 = \{S \times S \times S\}$, $V_2^1 = \{2 \times 2$ matrices with entries from $S\}$,

$$V_3^1 = \left\{ \begin{bmatrix} a & a & a & a & a \\ a & a & a & a & a \end{bmatrix}, \begin{bmatrix} a \\ a \\ a \\ a \\ a \\ a \end{bmatrix} \middle| \ a \in S \right\}$$

and

$$V_4^1 = \left\{ \begin{bmatrix} a & a \\ a & a \\ a & a \\ a & a \end{bmatrix}, [a\,a\,a\,a\,a\,a] \mid a \in S \right\}.$$

$V_1$ is a special set vector space over the set S.

Take $W_1 = (W_1^1, W_2^1, W_3^1, W_4^1) \subseteq (V_1^1, V_2^1, V_3^1, V_4^1)$ where

$$W_1^1 = \{S \times S \times \{0\}\} \subseteq V_1^1,$$

$$W_2^1 = \left\{ \begin{pmatrix} a & a \\ a & a \end{pmatrix} \bigg| a \in S \right\} \subseteq V_2^1,$$

$$W_3^1 = \left\{ \begin{bmatrix} a & a & a & a & a \\ a & a & a & a & a \end{bmatrix} \bigg| a \in S \right\} \subseteq V_3^1, \text{ and}$$

$$W_4^1 = \{[a\,a\,a\,a\,a\,a] \mid a \in S\} \subseteq V_4^1.$$

Thus $W_1$ is a special set vector subspace of $V_1$ over the set S.
Now $V_2 = (V_1^2, V_2^2, V_3^2, V_4^2, V_5^2)$ where $V_1^2 = \{S \times S \times S \times S \times S\}$,

$$V_2^2 = \left\{ \begin{pmatrix} a & a & a & a \\ a & a & a & a \end{pmatrix}, (a\,a\,a\,a\,a) \mid a \in S \right\},$$

$V_3^2 = \{3 \times 3 \text{ matrices with entries from S}\}$, $V_4^2 = \{s[x]$ polynomials in x with coefficients from S of degree less than or equal to 5$\}$ and $V_5^2 = \{4 \times 4 \text{ upper triangular matrices with entries from S}\}$. Take

$$W_2 = (W_1^2, W_2^2, W_3^2, W_4^2, W_5^2)\} \subseteq (V_1^2, V_2^2, V_3^2, V_4^2, V_5^2)$$

where

$$W_1^2 = S \times S \times S \times \{0\} \times \{0\} \subseteq V_1^2,$$

$$W_2^2 = \left\{ \begin{pmatrix} a & a & a & a \\ a & a & a & a \end{pmatrix} \bigg| a \in Z^+ \cup \{0\} \right\} \subseteq V_2^2,$$

$W_3^2 = \{3 \times 3 \text{ matrices with entries from } 5Z^+ \cup \{0\}\} \subseteq V_3^2$,

$W_4^2$ = {All polynomials of degree less than or equal to 3 with coefficients from S} $\subseteq V_4^2$

and

$W_5^2$ = {4 × 4 upper triangular matrices with entries from $3Z^+ \cup \{0\}$} $\subseteq V_5^2$.

Thus $W_2 = (W_1^2, W_2^2, W_3^2, W_4^2, W_5^2) \subseteq (V_1^2, V_2^2, V_3^2, V_4^2, V_5^2) = V_2$ is a special set vector subspace of $V_2$ over S.

Let $V_3 = (V_1^3, V_2^3, V_3^3, V_4^3)$ where

$V_1^3$ = {all lower triangular 5 × 5 matrices with entries from $S = Z^+ \cup \{0\}$},

$V_2^3$ = {[a a a a a], [a a], (a a a) | a ∈ S},

$$V_3^3 = \left\{ \begin{bmatrix} a & a \\ a & a \\ a & a \\ a & a \end{bmatrix}, \begin{bmatrix} a & a & a & a & a \\ a & a & a & a & a \end{bmatrix} \middle| a \in S \right\} \text{ and}$$

$V_4^3$ = {All polynomials in the variable x with coefficients from S of degree less than or equal to 7}.

Clearly $V_3 = (V_1^3, V_2^3, V_3^3, V_4^3)$ is a special set vector space over the set $S = Z^+ \cup \{0\}$. Take $W_3 = (W_1^3, W_2^3, W_3^3, W_4^3)$ where

$W_1^3$ = {all lower triangular 5 × 5 matrices with entries from $3Z^+ \cup \{0\}$} $\subseteq V_1^3$,

$W_2^3$ = {[a a a a], [a a] | a ∈ S} $\subseteq V_2^3$,

$$W_3^3 = \left\{ \begin{bmatrix} a & a \\ a & a \\ a & a \\ a & a \end{bmatrix} \;\middle|\; a \in 3Z^+ \cup \{0\} \right\} \subseteq V_3^3$$

and

$W_4^3 = \{$all polynomials in the variable x with coefficients from x of degree less than or equal to 5$\} \subseteq V_4^3$.

$W_3 = (W_1^3, W_2^3, W_3^3, W_4^3) \subseteq V_3 = (V_1^3, V_2^3, V_3^3, V_4^3)$ is a special set vector subspace of $V_3$ over the set S.

Thus $W = W_1 \cup W_2 \cup W_3 \subseteq V_1 \cup V_2 \cup V_3$ is a special set 3-vector subspace of V over the set S.

Now we proceed on to define the notion of special set n-linear algebra or special set linear n-algebra over a set S.

**DEFINITION 2.3.3**: *Let*
$$V = (V_1 \cup V_2 \cup ... \cup V_n)$$
$$= (V_1^1, V_2^1, ..., V_{n_1}^1) \cup (V_1^2, V_2^2, ..., V_{n_2}^2) \cup ... \cup (V_1^n, V_2^n, ..., V_{n_n}^n)$$
*be a special set n vector space over the set S.*
*If each $V_i = (V_1^i, V_2^i, ..., V_{n_i}^i)$ is a special set linear algebra over the set S for each i, i = 1, 2, ..., n then we call $V = (V_1 \cup V_2 \cup ... \cup V_n)$ to be a special set linear n-algebra over the set S.*

We illustrate this by the following example.

*Example 2.3.5*: Let
$$\begin{aligned} V &= (V_1 \cup V_2 \cup V_3 \cup V_4) \\ &= \{V_1^1, V_2^1, V_3^1, V_4^1\} \cup \{V_1^2, V_2^2, V_3^2\} \cup \\ &\quad \{V_1^3, V_2^3, V_3^3, V_4^3, V_5^3\} \cup \{V_1^4, V_2^4, V_3^4\} \end{aligned}$$
be a special set linear 4-algebra over the set $S = Z^+ \cup \{0\}$. Here each $V_i$ is a special set linear algebra over S, $1 \leq i \leq 4$. Now $V_1 = (V_1^1, V_2^1, V_3^1, V_4^1)$ is such that

$$V_1^1 = \{\text{All } 3 \times 3 \text{ matrices with entries in } S\},$$
$$V_2^1 = \{(a\ a\ a\ a\ a) \mid a \in S\},$$
$$V_3^1 = \left\{ \begin{bmatrix} a & a \\ a & a \\ a & a \end{bmatrix}, \begin{bmatrix} a & a & a \\ a & a & a \end{bmatrix} \middle| a \in S \right\}$$

and
$$V_4^1 = \{S[x] \text{ all polynomials in the variable}$$
x with coefficients form S of degree less than or equal to 4\}.

Clearly $V_1 = (V_1^1, V_2^1, V_3^1, V_4^1)$ is special set linear algebra over the set S. Now $V_2 = (V_1^2, V_2^2, V_3^2)$ where $V_1^2 = \{\text{all } 4 \times 4 \text{ upper triangular matrices with entries from the set } S\}$,

$$V_2^2 = \left\{ \begin{bmatrix} a & a & a & a \\ a & a & a & a \end{bmatrix}, \begin{bmatrix} a & a \\ a & a \\ a & a \end{bmatrix} \middle| a \in S \right\}$$

and
$$V_3^2 = \left\{ \begin{bmatrix} a & b \\ c & d \end{bmatrix} \middle| a, b, c, d \in S \right\}.$$

$V_2$ is a special set linear algebra over the set S.

$V_3 = (V_1^3, V_2^3, V_3^3, V_4^3, V_5^3)$ where $V_1^3 = S \times S \times S \times S$,

$$V_2^3 = \left\{ (a\ a\ a\ a\ a\ a), \begin{pmatrix} a \\ a \\ a \\ a \\ a \end{pmatrix} \middle| a \in S \right\},$$

$V_3^3 = \{\text{all } 4 \times 4 \text{ matrices with entries from } S\}$,
$V_4^3 = \{\text{all diagonal } 5 \times 5 \text{ matrices with entries from } S\}$
and

$$V_5^3 = \left\{ \begin{pmatrix} a & a \\ a & a \end{pmatrix}, \begin{pmatrix} a & a & a \\ a & a & a \\ a & a & a \end{pmatrix} \middle| a \in S \right\}.$$

$V_3 = (V_1^3, V_2^3, V_3^3, V_4^3, V_5^3)$ is a special set linear algebra over the set S. Finally $V_4 = (V_1^4, V_2^4, V_3^4)$ is such that $V_1^4 = S \times S \times S$,

$$V_2^4 = \left\{ (a\ a\ a\ a\ a\ a), \begin{bmatrix} a \\ a \\ a \end{bmatrix} \middle| a \in S \right\}$$

and

$$V_3^4 = \{\text{all } 7 \times 7 \text{ matrices with entries from } S\}.$$

$V_4 = (V_1^4, V_2^4, V_3^4)$ is a special set linear algebra over the set S. Thus $V = V_1 \cup V_2 \cup V_3 \cup V_4$ is a special set linear 4-algebra over the set S.

**Example 2.3.6:** Let $V = V_1 \cup V_2 \cup V_3$ be a special set linear trialgebra over the set $S = \{0, 1\}$. Here
$$V_1 = (V_1^1, V_2^1, V_3^1, V_4^1, V_5^1),$$
$$V_2 = (V_1^2, V_2^2, V_3^2) \text{ and } V_3 = (V_1^3, V_2^3)$$
where $V_1$, $V_2$ and $V_3$ are special set linear algebras over the set $S = \{0, 1\}$. $V_1 = (V_1^1, V_2^1, V_3^1, V_4^1, V_5^1)$ is given by

$$V_1^1 = S \times S \times S \times S.$$

$$V_2^1 = \left\{ \begin{pmatrix} a & b \\ c & a \end{pmatrix} \middle| a, b, c, d \in Z_2 \right\},$$

$$V_3^1 = \left\{ \begin{bmatrix} a & a & a & a & a \\ a & a & a & a & a \end{bmatrix}, \begin{bmatrix} a & a & a \\ a & a & a \end{bmatrix}, \begin{bmatrix} a & a \\ a & a \\ a & a \\ a & a \end{bmatrix} \middle| a \in Z^+ \cup \{0\} \right\},$$

$V_4^1 = $ {All polynomials with coefficients from S of degree less than or equal to 4} and $V_5^1 = $ {All $3 \times 5$ matrices with elements from S}. $V_1 = (V_1^1, V_2^1, V_3^1, V_4^1, V_5^1)$ is a special set linear algebra over the set S. We take $V_2 = (V_1^2, V_2^2, V_3^2)$ where

$$V_1^2 = (S \times S \times S \times S),$$

$$V_2^2 = \left\{ \begin{pmatrix} a & a \\ a & a \end{pmatrix} \middle| a \in Z^+ \cup \{0\} \right\}$$

and

$$V_3^2 = \{(a\ a\ a\ a\ a), (a\ a), (a\ a\ a) \mid a \in Z^+ \cup \{0\}\};$$

$V_2 = \left(V_1^2, V_2^2, V_3^2\right)$ is a special set linear algebra over the set S = $\{0, 1\}$. Now $V_3 = \left(V_1^3, V_2^3\right)$ where $V_1^3 = S \times S \times S \times S \times S$ and $V_2^3 = \{(a\ a\ a\ a\ a), (a\ a\ a), (a\ a\ a\ a\ a\ a) \mid a \in Z^+ \cup \{0\}\}$. Thus $V_3 = \left(V_1^3, V_2^3\right)$ is a special set linear algebra over S = $\{0, 1\}$. Hence $V = (V_1 \cup V_2 \cup V_3)$ is a special set linear trialgebra over the set S = $\{0, 1\}$.

Now we proceed onto define the notion of special set linear n-subalgebra.

**DEFINITION 2.3.4:** *Let*
$V = (V_1 \cup V_2 \cup ... \cup V_n)$
$= \left(V_1^1, V_2^1, ..., V_{n_1}^1\right) \cup \left(V_1^2, V_2^2, ..., V_{n_2}^2\right) \cup ... \cup \left(V_1^n, V_2^n, ..., V_{n_n}^n\right)$
*be a special set linear n-algebra over the set S. Suppose*
$W = (W_1 \cup W_2 \cup ... \cup W_n)$
$= \left(W_1^1, W_2^1, ..., W_{n_1}^1\right) \cup \left(W_1^2, W_2^2, ..., W_{n_2}^2\right) \cup ... \cup \left(W_1^n, W_2^n, ..., W_{n_n}^n\right)$
$\subseteq (V_1 \cup V_2 \cup ... \cup V_n)$
$= V$ *and if each* $W_i = \left(W_1^i, W_2^i, ..., W_{n_i}^i\right) \subseteq V_i$ *is a special set linear subalgebra of $V_i$ over the set S, for i = 1, 2, ..., n, then we call $W = W_1 \cup W_2 \cup ... \cup W_n$ to be a special set linear n-subalgebra of V over the set S.*

*Now in case if in $W = W_1 \cup ... \cup W_n$ some of the $W_i$'s are not special set linear subalgebra of $V_i$'s but only a special set vector subspace of $V_i$'s, $1 \le i \le n$ then we call $W = (W_1 \cup W_2 \cup ... \cup W_n) \subseteq V = (V_1 \cup V_2 \cup ... \cup V_n)$ to be a pseudo special set linear n-subalgebra of V over the set S. If in $W = (W_1 \cup W_2 \cup ... \cup W_n)$ all the $W_i$'s are only special set vector subspaces of $V_i$'s for each $i = 1, 2, ..., n$ then we call $W = W_1 \cup W_2 \cup ... \cup W_n$ to be a pseudo special set vector subspace of V over the set S.*

Now we illustrate these concepts by some simple examples.

***Example 2.3.7:*** Let
$$V = (V_1 \cup V_2 \cup V_3 \cup V_4)$$
$$= (V_1^1, V_2^1, V_3^1) \cup (V_1^2, V_2^2, V_3^2) \cup (V_1^3, V_2^3) \cup (V_1^4, V_2^4, V_3^4, V_4^4)$$

be a special set linear 4-algebra over the set $S = Z^+ \cup \{0\}$. Here $V_1 = (V_1^1, V_2^1, V_3^1)$ where

$$V_1^1 = \left\{ \begin{pmatrix} a & b \\ c & d \end{pmatrix} \middle| a,b,c,d \in S \right\},$$
$$V_2^1 = S \times S \times S \times S \times S,$$
$$V_3^1 = \{(a\ a\ a\ a\ a\ a), (a\ a\ a) \mid a \in S\}.$$

Take $W_1 = (W_1^1, W_2^1, W_3^1) \subseteq (V_1^1, V_2^1, V_3^1)$ where

$$W_1^1 = \left\{ \begin{pmatrix} a & a \\ a & a \end{pmatrix} \middle| a \in S \right\} \subseteq V_1^1,$$
$$W_2^1 = \{(a\ a\ a\ a\ a) \mid a \in S\} \subseteq V_2^1$$

and
$$W_3^1 = \{(a\ a\ a\ a\ a\ a), (a\ a\ a) \mid a \in 5Z^+ \cup \{0\}\} \subseteq V_3^1.$$

Clearly $W_1 = (W_1^1, W_2^1, W_3^1) \subseteq (V_1^1, V_2^1, V_3^1)$ is a special set linear subalgebra of $V_1$ over the set $S = Z^+ \cup \{0\}$.

Consider $V_2 = (V_1^2, V_2^2, V_3^2)$ where $V_1^2 = $ {all $3 \times 3$ upper triangular matrices with entries from S},

$$V_2^2 = \left\{ \begin{pmatrix} a & a & a \\ a & a & a \end{pmatrix}, \begin{pmatrix} a & a \\ a & a \\ a & a \\ a & a \end{pmatrix} \middle| a \in S \right\}$$

and $V_3^2 = $ {S[x] all polynomials of degree less than or equal to 4 with coefficients from S}. Take in $V_2$ a special subset $W_2 = (W_1^2, W_2^2, W_3^2) \subseteq (V_1^2, V_2^2, V_3^2) = V_2$ where $W_1^2 = $ {all upper triangular $3 \times 3$ matrices with entries from $3Z^+ \cup \{0\}\} \subseteq V_1^2$,

$$W_2^2 = \left\{ \begin{pmatrix} a & a & a \\ a & a & a \end{pmatrix}, \begin{pmatrix} a & a \\ a & a \\ a & a \\ a & a \end{pmatrix} \middle| a \in 5S^+ \cup \{0\} \right\} \subseteq V_2^2$$

and $W_3^2 = $ {S[x] all polynomials of degree less than or equal to two with coefficients from S}. $W_2 = (W_1^2, W_2^2, W_3^2) \subseteq V_2 = (V_1^2, V_2^2, V_3^2)$ is a special set sublinear algebra of $V_2$. Set $V_3 = (V_1^3, V_2^3)$ as $V_1^3 = $ {set of all $5 \times 5$ matrices with entries from S} and $V_2^3 = \{(a\ a\ a), (a\ a)| a \in S\}$. In $V_3$ take $W_3 = (W_1^3, W_2^3)$ where $W_1^3 = $ {set of all $5 \times 5$ matrices with entries from the set $5Z^+ \cup \{0\} \subseteq V_1^3$ and $W_2^3 = \{(a\ a\ a),(a\ a)| a \in 7Z^+ \cup \{0\}\} \subseteq V_2^3$. Clearly $W_3 = (W_1^3, W_2^3)$ is a special set linear subalgebra of $V_3$ over the set S. Finally assume $V_4 = (V_1^4, V_2^4, V_3^4, V_4^4)$ as $V_1^4 = \{S \times S \times S \times S \times S \times S\}$,

$$V_2^4 = \left\{ \begin{pmatrix} a \\ a \\ a \\ a \\ a \end{pmatrix}, \begin{pmatrix} a \\ a \\ a \\ a \end{pmatrix}, \begin{pmatrix} a \\ a \\ a \end{pmatrix} \middle| a \in S \right\},$$

$V_3^4$ = {All upper triangular $2 \times 2$ matrices with entries S} and
$V_4^4$ = {(a a a a a), (a a a) | a ∈ S} where

$$W_4 = \left( W_1^4, W_2^4, W_3^4, W_4^4 \right) \subseteq \left( V_1^4, V_2^4, V_3^4, V_4^4 \right) = V_4$$

as a special set linear subalgebra of $V_4$. Thus

$$\begin{aligned} W &= (W_1 \cup W_2 \cup W_3 \cup W_4) \\ &= \left( W_1^1, W_2^1, W_3^1 \right) \cup \left( W_1^2, W_2^2, W_3^2 \right) \cup \\ &\quad \left( W_1^3, W_2^3 \right) \cup \left( W_1^4, W_2^4, W_3^4, W_4^4 \right) \\ &\subseteq (V_1 \cup V_2 \cup V_3 \cup V_4) \end{aligned}$$

is a special set 4-linear subalgebra of V over S.

*Example 2.3.8:* Let $V = (V_1 \cup V_2 \cup V_3)$ be a special set linear trialgebra over the set $S = \{0,1\}$ where $V_1 = (V_1^1, V_2^1, V_3^1, V_4^1)$, $V_2 = (V_1^2, V_2^2)$ and $V_3 = (V_1^3, V_2^3, V_3^3)$ with

$$V_1^1 = \left\{ \begin{pmatrix} a & b \\ c & d \end{pmatrix} \middle| a,b,c,d \in Z^+ \cup \{0\} \right\},$$

$V_2^1$ = {(a a a a), (a a a)| a ∈ $Z^+ \cup \{0\}$},

$$V_3^1 = \left\{ \begin{pmatrix} a & a & a & a & a \\ a & a & a & a & a \end{pmatrix}, \begin{pmatrix} a & a \\ a & a \\ a & a \\ a & a \end{pmatrix} \middle| a \in Z^+ \cup \{0\} \right\}$$

and

$V_4^1 = \{$all polynomials in the variable x with coefficients from $Z^+ \cup \{0\}$ of degree less than or equal to 5$\}$.

Take
$$W_1 = \left(W_1^1, W_2^1, W_3^1, W_4^1\right)$$
where
$$W_1^1 = \left\{ \begin{pmatrix} a & a \\ a & a \end{pmatrix} \middle| a \in Z^+ \cup \{0\} \right\} \subseteq V_1^1$$

is a set linear subalgebra of $V_1^1$, $W_2^1 = \{(a\ a\ a)|\ a \in Z^+ \cup \{0\}\} \subseteq V_2^1$ is a set linear subvector space of $V_2^1$.

$$W_3^1 = \left\{ \begin{pmatrix} a & a & a & a & a \\ a & a & a & a & a \end{pmatrix} \middle| a \in Z^+ \cup \{0\} \right\} \subseteq V_3^1$$

is a set linear subvector space of $V_3^1$ and $W_4^1 = \{$all polynomials in the variable x of degree less than or equal to 3$\} \subseteq V_4^1$ is a set linear subalgebra of $V_4^1$. Thus $W_1 = \left(W_1^1, W_2^1, W_3^1, W_4^1\right) \subseteq V_1$ is a special set linear subalgebra of $V_1$ over the set $S = \{0, 1\}$. Now consider $V_2 = \left(V_1^2, V_2^2\right)$ where

$$V_1^2 = \left\{ \begin{pmatrix} a & a & a \\ a & a & a \\ a & a & a \end{pmatrix} \middle| a \in Z^+ \cup \{0\} \right\}$$

and

$$V_2^2 = \left\{ \begin{pmatrix} a \\ a \\ a \\ a \\ a \end{pmatrix}, \begin{pmatrix} a & a \\ a & a \\ a & a \\ a & a \\ a & a \end{pmatrix} \middle| a \in Z^+ \cup \{0\} \right\}.$$

Take a subset $W_2 = \left(W_1^2, W_2^2\right)$ of $V_2$ as

$$W_1^2 = \left\{ \begin{pmatrix} a & a & a \\ a & a & a \\ a & a & a \end{pmatrix} \middle| a \in 3Z^+ \cup \{0\} \right\} \subseteq V_1^2$$

is set linear subalgebra of $V_1^2$ and

$$W_2^2 = \left\{ \begin{pmatrix} a \\ a \\ a \\ a \\ a \end{pmatrix} \middle| a \in 5Z^+ \cup \{0\} \right\}$$

is a set vector subspace of $V_2^2$. Thus $W_2 = \left( W_1^2, W_2^2 \right) \subseteq V_2$ is a special set linear subalgebra of $V_2$ over the set $S = \{0, 1\}$.

Define $V_3 = \left( V_1^3, V_2^3, V_3^3 \right)$ where $V_1^3 = $ {all 4×4 upper triangular matrices with entries from $Z^+ \cup \{0\}$}, $V_2^3 = Z^o \times Z^o \times Z^o \mid Z^o = Z^+ \cup \{0\}$} and

$$V_3^3 = \left\{ \begin{pmatrix} a & a \\ a & a \\ a & a \\ a & a \end{pmatrix}, \begin{pmatrix} a & a & a & a & a \\ a & a & a & a & a \end{pmatrix} \middle| a \in Z^o \right\}.$$

Take $W_3 = \left( W_1^3, W_2^3, W_3^3 \right) \subseteq V_3 = \left( V_1^3, V_2^3, V_3^3 \right)$ such that

$$W_1^3 = \left\{ \begin{pmatrix} a & a & a & a \\ 0 & a & a & a \\ 0 & 0 & a & a \\ 0 & 0 & 0 & a \end{pmatrix} \middle| a \in Z^+ \cup \{0\} \right\} \subseteq V_1^3$$

is a set linear subalgebra of $V_1^3$, $W_2^3 = \{Z^o \times Z^o \times \{0\}\} \subseteq V_2^3$ is again a set linear subalgebra of $V_2^3$ and

$$W_3^3 = \left\{ \begin{pmatrix} a & a \\ a & a \\ a & a \end{pmatrix} \middle| a \in Z^o \right\}$$

is a set vector subspace of $V_3^3$. Hence $W_3 = \left( W_1^3, W_2^3, W_3^3 \right)$ $\subseteq V_3$ is a special set linear subalgebra of $V_3$.

Thus $W = (W_1 \cup W_2 \cup W_3) \subseteq (V_1 \cup V_2 \cup V_3) = V$ is a special set linear sub trialgebra of $V$ over the set $S = \{0, 1\}$.

Now we proceed onto illustrate a pseudo special set vector n subspace of V over the set S.

***Example 2.3.9:*** Let $V = (V_1 \cup V_2 \cup V_3 \cup V_4)$ where
$$V_1 = ( V_1^1, V_2^1 ),$$
$$V_2 = \left( V_1^2, V_2^2, V_3^2 \right),$$
$$V_3 = \left( V_1^3, V_2^3 \right) \text{ and } V_4 = \left( V_1^4, V_2^4, V_3^4, V_4^4 \right)$$

be a special set linear 4-algebra over the set $S = Z^+ \cup \{0\}$.
Take
$$V_1^1 = \left\{ \begin{pmatrix} a & b \\ c & d \end{pmatrix} \middle| a,b,c,d \in Z^+ \cup \{0\} \right\}$$

and $V_2^1 = \{(a\ a\ a\ a), (a\ a), (a\ a\ a\ a)$ such that $a \in S\}$.
Take
$$W_1^1 = \left\{ \begin{pmatrix} 2a & 2a \\ 2a & 2a \end{pmatrix}, \begin{pmatrix} 3a & 3a \\ 3a & 3a \end{pmatrix}, \begin{pmatrix} 7a & 7a \\ 7a & 7a \end{pmatrix} \middle| a \in Z^+ \cup \{0\} \right\} \subseteq V_1^1,$$

$W_1^1$ is only a pseudo set vector subspace of $V_1^1$. $W_2^1 = \{(a\ a\ a\ a), (a\ a)|\ a \in S\} \subseteq V_2^1$ is again only a pseudo set vector subspace of $V_2^1$. Hence $W_1 = \left(W_1^1, W_2^1\right)$ is a pseudo special set vector subspace of V over the set S.

Consider $V_2 = \left(V_1^2, V_2^2, V_3^2\right)$ where $V_1^2 = \{$all $4 \times 5$ matrices with entries in S$\}$, $V_2^2 = \{2 \times 7$ matrices with entries in S$\}$ and $V_3^2 = \{$(a a a a a), (a a) | a $\in$ S$\}$. Now take $W_1^2 = \{$all $4 \times 5$ matrices with entries from $2Z^+$ and all $4 \times 5$ matrices with entries from $3Z^+ \cup \{0\}\} \subseteq V_1^2$ is a pseudo set vector subspace of $V_1^2$. $W_2^2 = \{2 \times 7$ matrices with entries from $5Z^+ \cup \{0\}$ and $2 \times 7$ matrices with entries from $7Z^+ \cup \{0\}\} \subseteq V_2^2$ is only a pseudo set vector subspace of $V_2^2$ over $S = Z^+ \cup \{0\}$. $W_3^2 = \{$(a a a a a), (a a)| a $\in 3Z^+ \cup \{0\}\} \subseteq V_3^2$ is only a set vector subspace of $V_3^2$. Thus $W_2 = \left(W_1^2, W_2^2, W_3^2\right)$ is a pseudo special set vector subspace of $V_2$ over the set $S = Z^+ \cup \{0\}$. We take $V_3 = \left(V_1^3, V_2^3\right)$ where $V_1^3 = \{$set of all $3 \times 5$ matrices with entries from $Z^+ \cup \{0\}\}$ and $V_2^3 = \{$(a a a a a), (a a a), (a a), (a a a a)| a $\in Z^+ \cup \{0\}\}$. Consider $W_1^3 = \{$set of all $3 \times 5$ matrices with entries from the set $5Z^+ \cup \{0\}$ and $3 \times 5$ matrices with entries form the set $3Z^+ \cup \{0\}\} \subseteq V_1^3$. Take $W_2^3 = \{$(a a), (a a a a)| a $\in Z^+ \cup \{0\}\} \subseteq V_2^3$; Clearly $W_2^3$ is only a set vector subspace of $W_2^3$. Take $W_3 = \left(W_1^3, W_2^3\right)$ is only a pseudo special set vector subspace of $V_3$.

Finally consider $V_4 = \left(V_1^4, V_2^4, V_3^4, V_4^4\right)$ where $V_1^4 = \{2 \times 4$ matrices with entries from $Z^+ \cup \{0\}\}$,

$$V_2^4 = \left\{ \begin{pmatrix} a & a \\ a & a \\ a & a \end{pmatrix}, \begin{pmatrix} a & a & a & a & a & a \\ a & a & a & a & a & a \\ a & a & a & a & a & a \end{pmatrix} \middle| a \in Z^+ \cup \{0\} \right\},$$

$V_3^4$ = {set of all polynomial in the variable x with coefficients from $Z^+ \cup \{0\}$ of degree less than or equal to 9} and $V_4^4$ = {6×6 upper triangular matrices with entries from $Z^+ \cup \{0\}$}. Take $W_4 = \left(W_1^4, W_2^4, W_3^4, W_4^4\right) \subseteq \left(V_1^4, V_2^4, V_3^4, V_4^4\right)$ where $W_1^4$ = {2 × 4 matrices with entries from $7Z^+ \cup \{0\}$ and 2 × 4 matrices with entries from $11Z^+ \cup \{0\}\} \subseteq V_1^4$ is only a pseudo set vector subspace of $V_1^4$ over the set $Z^+ \cup \{0\}$,

$$W_2^4 = \left\{ \begin{pmatrix} a & a \\ a & a \\ a & a \end{pmatrix}, \begin{pmatrix} a & a & a & a & a & a \\ a & a & a & a & a & a \\ a & a & a & a & a & a \end{pmatrix} \middle| a \in 3Z^+ \cup \{0\} \right\} \subseteq V_2^4,$$

$W_2^4$ is only a subspace of $V_2^4$. Take $W_3^4$ = {set of all polynomials n $(x^6 + x^5 + 2x^4 + 3x^3 + 5x + 1)$; n ∈ $3Z^+ \cup \{0\}$ and $m(3x^9 + 5x^8 + 4x^7 + 3x^2 - 7)$ such that m ∈ $7Z^+ \cup \{0\}\} \subseteq V_3^4$; clearly $W_3^4$ is a pseudo set vector subspace of $V_3^4$ over the set $Z^+ \cup \{0\}$. Hence $W_4 = \left(W_1^4, W_2^4, W_3^4, W_4^4\right)$ is a pseudo special set vector subspace of $V_4$ over the set $Z^+ \cup \{0\}$. Now W = ($W_1 \cup W_2 \cup W_3 \cup W_4$) ⊆ ($V_1 \cup V_2 \cup V_3 \cup V_4$) = V is a pseudo special set vector 4-subspace of V defined over the set $Z^+ \cup \{0\}$.

Now we proceed onto define the new notion of pseudo special set linear n subalgebra of V over the set S.

*Example 2.3.10:* Let V = ($V_1 \cup V_2 \cup V_3 \cup V_4$) where
$$V_1 = \left(V_1^1, V_2^1, V_3^1\right),$$
$$V_2 = \left(V_1^2, V_2^2\right),$$
$$V_3 = \left(V_1^3, V_2^3, V_3^3\right)$$
and $$V_4 = \left(V_1^4, V_2^4, V_3^4, V_4^4\right)$$

be a special set linear 4-algebra over the set $S = \{0, 1\}$. Choose in $V_1 = \left(V_1^1, V_2^1, V_3^1\right)$; $V_1^1 = \{2 \times 2$ matrices with entries from $Z^o = Z^+ \cup \{0\}\}$, $V_2^1 = \{Z^o \times Z^o \times Z^o \times Z^o\}$, $V_3^1 = \{(a\ a\ a\ a\ a\ a), (a\ a\ a) / a \in Z^o\}$ so that $V_1$ is a special set linear algebra over the set $S = \{0,1\}$. Let $V_2 = \left(V_1^2, V_2^2\right)$ such that $V_1^2 = \{$all $3 \times 3$ matrices with entries from $Z^o\}$ and $V_2^2 = \{$set of all $2 \times 5$ matrices with entries from $Z^o$ and the set of all $4 \times 3$ matrices with entries from $Z^o\}$ so that $V_2$ is also a special set linear algebra over $S = \{0, 1\}$. $V_3 = (V_1^3, V_2^3, V_3^3)$ is such that

$V_1^3 = \{4 \times 4$ upper triangular matrices with entries from $Z^o\}$,
$V_2^3 = \{(a\ a\ a\ a), (a\ a), (a\ a\ a\ a\ a) \mid a \in Z^o\}$

and

$V_3^3 = \{Z^o[x]$ all polynomials in the variable x of degree less than or equal to 6 with coefficients from $Z^o\}$.

Thus $V_3 = \left(V_1^3, V_2^3, V_3^3\right)$ is a special set linear algebra over the set $S = \{0, 1\}$. Now take $V_4 = \left(V_1^4, V_2^4, V_3^4, V_4^4\right)$ where $V_1^4 = \{Z^o \times Z^o\}$, $V_2^4 = \{$upper triangular $2 \times 2$ matrices with entries from $Z^o\}$, $V_3^4 = \{2 \times 3$ matrices with entries from $3Z^+ \cup \{0\}$ and $2 \times 3$ matrices with entries from $5Z^+ \cup \{0\}\}$ and

$$V_4^4 = \left\{(a\ a\ a\ a\ a), \begin{bmatrix} a \\ a \\ a \end{bmatrix}, \begin{bmatrix} a & a \\ a & a \\ a & a \\ a & a \\ a & a \end{bmatrix} \middle| a \in Z \right\}.$$

$V_4 = \left(V_1^4, V_2^4, V_3^4, V_4^4\right)$ is a special set linear algebra over the set $S = \{0, 1\}$.

Take $W = (W_1 \cup W_2 \cup W_3 \cup W_4)$ with $W_1 = \left(W_1^1, W_2^1, W_3^1\right) \subseteq \left(V_1^1, V_2^1, V_3^1\right)$ such that

$$W_1^1 = \left\{ \begin{pmatrix} a & a \\ a & a \end{pmatrix} \middle| a \in Z^o \right\} \subseteq V_1^1,$$

$W_2^1 = \{Z^o \times Z^o \times \{0\} \times \{0\}\} \subseteq V_2^1$, $W_3^1 = \{(a\ a\ a\ a\ a\ a), (a\ a\ a) | a \in 5Z^o\} \subseteq V_3^1$. $W_1$ is a special set linear subalgebra of $V_1$ over $S = \{0, 1\}$. $W_2 = \left(W_1^2, W_2^2\right) \subseteq V_2$ where $W_1^2 = \{3 \times 3$ matrices with entries from $3Z^o$; $3 \times 3$ matrices with entries from $5Z^o\} \subseteq V_1^2$, $W_2^2 = \{$set of all $2 \times 5$ matrices with entries from $3Z^o$ and the set of all $4 \times 3$ matrices with entries from $5Z^o\} \subseteq V_2^2$. Thus $W_2 = \left(W_1^2, W_2^2\right)$ is only a pseudo special vector subspace of $V_2$ over the set $S = \{0, 1\}$. $W_3 = \left(W_1^3, W_2^3, W_3^3\right) \subseteq V_3$ where $W_1^3 = \{4 \times 4$ upper triangular matrices with entries from $3Z^o$ and $4 \times 4$ upper triangular matrices with entries from $5Z^o\} \subseteq V_1^3$, $W_2^3 = \{(a\ a), (a\ a\ a\ a) | a \in Z^o\} \subseteq V_2^3$ and $W_3^3 = \{Z^o[x]$ all polynomials of degree less than or equal to 3 and all polynomials of the form $n(1 + x + x^2 + x^3 + x^4 + x^5 + x^6)$ such that $n \in Z^o\}\} \subseteq V_3^3$. Clearly $W_3 = \left(W_1^3, W_2^3, W_3^3\right) \subseteq V_3$ is only a pseudo special set vector subspace of $V_3$ over the set $S = \{0, 1\}$. $W_4 = \left(W_1^4, W_2^4, W_3^4, W_4^4\right) \subseteq V_4$ is such that $W_1^4 = \{3Z^o \times 3Z^o\} \subseteq V_1^4$,

$$W_2^4 = \left\{ \begin{pmatrix} a & a \\ a & a \end{pmatrix} \middle| a \in Z^o \right\} \subseteq V_2^4,$$

$W_3^4 = \{2 \times 3$ matrices with entries from $6Z^o \cup \{0\}$ and $2 \times 3$ matrices with entries from $10Z^o \cup \{0\}\} \subseteq V_3^4$ and

$$W_4^4 = \left\{ (a \ \ a \ \ a \ \ a \ \ a), \begin{bmatrix} a \\ a \\ a \end{bmatrix} \middle| a \in Z^o \right\} \subseteq V_4^4.$$

Thus $W_4 = \left( W_1^4, W_2^4, W_3^4, W_4^4 \right) \subseteq V_4 = \left( V_1^4, V_2^4, V_3^4, V_4^4 \right)$ is a special set linear subalgebra of $V_4$ over $S = \{0, 1\}$. Thus $W = (W_1 \cup W_2 \cup W_3 \cup W_4)$ is only a pseudo special set linear 4-subalgebra of V over the set $\{0, 1\}$.

Now we proceed on to define the new notion of special set n-dimension generator of a special set n-vector space over the set S.

**DEFINITION 2.3.5:** *Let*
$$V = (V_1 \cup V_2 \cup \ldots \cup V_n)$$
$$= \left( V_1^1, V_2^1, \ldots, V_{n_1}^1 \right) \cup \left( V_1^2, V_2^2, \ldots, V_{n_2}^2 \right) \cup \ldots \cup \left( V_1^n, V_2^n, \ldots, V_{n_n}^n \right)$$
*be a special set vector n-space over the set S. Suppose*
$$X = (X_1 \cup X_2 \cup \ldots \cup X_n)$$
$$= \left( X_1^1, X_2^1, \ldots, X_{n_1}^1 \right) \cup \left( X_1^2, X_2^2, \ldots, X_{n_2}^2 \right) \cup \ldots \cup \left( X_1^n, X_2^n, \ldots, X_{n_n}^n \right)$$
$$\subseteq \left( V_1^1, V_2^1, \ldots, V_{n_1}^1 \right) \cup \left( V_1^2, V_2^2, \ldots, V_{n_2}^2 \right) \cup \ldots \cup \left( V_1^n, V_2^n, \ldots, V_{n_n}^n \right)$$

*such that $X_t^i \subseteq V_t^i$; $1 \leq i \leq n$ and $1 \leq t \leq n_1, n_2, \ldots, n_n$, is such that each $X_t^i$ generates $V_t^i$ over the set S, $i = 1, 2, \ldots, n$, then we call X to be the special set n-dimension generator of V over the set S. If each $X_t^i$ is of finite cardinality, $1 \leq t \leq n_1, n_2, \ldots, n_n$, $i = 1, 2, \ldots, n$, then we say the special set n-vector space over S is finite n-dimensional. If even one of the $X_t^i$ is of infinite cardinality then we say the special set n-vector space over S is infinite n-dimensional.*

We shall describe this by the following example:

***Example 2.3.11:*** Let $V = (V_1 \cup V_2 \cup V_3 \cup V_4)$ be a special set 4 vector space over the set $S = Z^+ \cup \{0\} = Z^o$. Here
$$V_1 = \left(V_1^1, V_2^1, V_3^1, V_4^1\right),$$
$$V_2 = \left(V_1^2, V_2^2\right),$$
$$V_3 = \left(V_1^3, V_2^3, V_3^3\right)$$
and $\qquad V_4 = \left(V_1^4, V_2^4\right)$

where

$$V_1^1 = \left\{ \begin{pmatrix} a & b \\ c & d \end{pmatrix} \middle| a, b, c, d \in Z^o = S \right\},$$

$V_2^1 = \{(a\ a\ a\ a), (a\ a\ a) \mid a \in Z^o\}$; $V_3^1 = \{Z^o[x]$; all polynomials of degree less than or equal to 5$\}$ and

$$V_4^1 = \left\{ \begin{pmatrix} a & a & a & a \\ a & a & a & a \end{pmatrix}, \begin{bmatrix} a & a \\ a & a \\ a & a \end{bmatrix} \middle| a \in Z^o \right\}.$$

$V_1 = \left(V_1^1, V_2^1, V_3^1, V_4^1\right)$ has a subset

$$X_1 = \left(X_1^1, X_2^1, X_3^1, X_4^1\right)$$
$$= \left\{ \begin{pmatrix} 1 & 0 \\ 0 & 0 \end{pmatrix}, \begin{pmatrix} 0 & 1 \\ 0 & 0 \end{pmatrix}, \begin{pmatrix} 0 & 0 \\ 1 & 0 \end{pmatrix}, \begin{pmatrix} 0 & 0 \\ 0 & 1 \end{pmatrix} \right\}, \{(1\ 1\ 1\ 1), (1\ 1\ 1)\},$$

{an infinite set of elements which are polynomials},

$$\left\{ \begin{pmatrix} 1 & 1 & 1 & 1 \\ 1 & 1 & 1 & 1 \end{pmatrix}, \begin{pmatrix} 1 & 1 \\ 1 & 1 \\ 1 & 1 \end{pmatrix} \right\} \subseteq \left(V_1^1, V_2^1, V_3^1, V_4^1\right)$$

is a special set generator of $V_1$ and $|X_1| = (4, 2, \infty, 2)$ … I.

Now $V_2 = \left(V_1^2, V_2^2\right)$ where

$$V_1^2 = \left\{ \begin{pmatrix} a & a & a \\ a & a & a \\ a & a & a \end{pmatrix} \text{ such that } a \in Z^o \right\},$$

$$V_2^2 = \left\{ \begin{bmatrix} a & a & a \\ a & a & a \end{bmatrix}, \begin{bmatrix} a & a \\ a & a \\ a & a \\ a & a \end{bmatrix} \middle| a \in Z^o \right\}$$

with a special set generator

$$X_2 = \left( X_1^2, X_2^2 \right) = \left( \left\{ \begin{pmatrix} 1 & 1 & 1 \\ 1 & 1 & 1 \end{pmatrix} \right\}, \left\{ \begin{pmatrix} 1 & 1 & 1 \\ 1 & 1 & 1 \end{pmatrix}, \begin{bmatrix} 1 & 1 \\ 1 & 1 \\ 1 & 1 \\ 1 & 1 \end{bmatrix} \right\} \right)$$

and $|X_2| = (1, 2) \ldots$ II.

$V_3 = \left( V_1^3, V_2^3, V_3^3 \right)$ where

$$V_1^3 = \{ n(1 + x + x^2 + x^3 + x^4 + x^5) | n \in Z^o \},$$

$$V_2^3 = \left\{ \begin{pmatrix} a & b & c \\ 0 & d & e \\ 0 & 0 & g \end{pmatrix} \middle| a,b,c,d,e,f,g \in Z^o \right\},$$

$V_3^3 = \{(a\ a\ a\ a\ a\ a), (a\ a), (a\ a\ a) \mid a \in Z^o\}$. Let $X_3 = \left( X_1^3, X_2^3, X_3^3 \right)$ be the special set generating $V_3$, then $X_3 = (\{(1 + x + x^2 + x^3 + x^4 + x^5)\}$, {an infinite collection of $3 \times 3$ upper triangular matrices}, $\{(1\ 1\ 1\ 1\ 1\ 1), (1\ 1), (1\ 1\ 1)\} \subseteq V_3 = \left( V_1^3, V_2^3, V_3^3 \right)$.

Thus $|X_3| = (1, \infty, 3) \ldots$ III.

Finally $V_4 = \left( V_1^4, V_2^4 \right)$ where

$$V_1^4 = \left\{ \begin{pmatrix} a & a \\ a & a \end{pmatrix} \begin{pmatrix} a & a & a \\ 0 & a & a \\ 0 & 0 & a \end{pmatrix} \begin{pmatrix} a & a & a \\ a & a & a \end{pmatrix} \middle| a \in Z^\circ \right\}$$

and $V_2^4 = \{n(1+ x + x^2), n(x^3 + 3x^2 + x + 3), n(x^7 + 7) \mid n \in Z^\circ\}$. The special set generator subset $X_4 = (X_1^4, X_2^4) \subseteq (V_1^4, V_2^4)$ is given by

$$X_4 = \left\{ \left\{ \begin{pmatrix} 1 & 1 \\ 1 & 1 \end{pmatrix}, \begin{pmatrix} 1 & 1 & 1 \\ 0 & 1 & 1 \\ 0 & 0 & 1 \end{pmatrix}, \begin{pmatrix} 1 & 1 & 1 \\ 1 & 1 & 1 \end{pmatrix} \right\}, \right.$$

$\{(1 + x + x^2, x^7 + 7, x^3 + 3x^2 + x + 3)\}\} \subseteq (V_1^4, V_2^4)$ and $|X_4| = (3, 3) \ldots$ IV.

$X_4$ is the special set generator of $V_4$. Thus we see $X = (X_1 \cup X_2 \cup X_3 \cup X_4) \subseteq (V_1 \cup V_2 \cup V_3 \cup V_4) = V$ is a special set 4-generator of $V$ over $S$ and the 4-dimension of $V$ is $|X| = (|X_1|, |X_2|, |X_3|, |X_4|) = \{(4, 2, \infty, 2), (1, 2), (1, \infty, 3), (3, 3)\}$. Thus $V$ is an infinite dimensional, special set 4-vector space over $S$.

Now we proceed on to define the notion of special set linear n-transformation of special set n-vector spaces defined over the same set S.

**DEFINITION 2.3.6:** *Let $V = (V_1 \cup V_2 \cup \ldots \cup V_n)$ and $W = (W_1 \cup W_2 \cup \ldots \cup W_n)$ be two special set n-vector spaces defined over the same set S; such that in both V and W each $V_i$ in V and $W_i$ in W have same number of sets which are set vector spaces, i.e., $V_i = (V_1^i, \ldots, V_{n_i}^i)$ and $W_i = (W_1^i, \ldots, W_{n_i}^i)$ true for $i = 1, 2, \ldots, n$.*

*That is*
$$V = (V_1 \cup V_2 \cup \ldots \cup V_n)$$
$$= (V_1^1, V_2^1, \ldots, V_{n_1}^1) \cup (V_1^2, V_2^2, \ldots, V_{n_2}^2) \cup \ldots \cup (V_1^n, V_2^n, \ldots, V_{n_n}^n)$$
*and*
$$W = (W_1 \cup W_2 \cup \ldots \cup W_n)$$

$$= \left(W_1^1, W_2^1, \ldots, W_{n_1}^1\right) \cup \left(W_1^2, W_2^2, \ldots, W_{n_2}^2\right) \cup \ldots \cup \left(W_1^n, W_2^n, \ldots, W_{n_n}^n\right)$$

be such special set n-vector spaces over the set S.
Let $T : V \to W$ such that
$$T = (T_1 \cup T_2 \cup \ldots \cup T_n)$$
$$= \left(T_1^1, T_2^1, \ldots, T_{n_1}^1\right) \cup \left(T_1^2, T_2^2, \ldots, T_{n_2}^2\right) \cup \ldots \cup \left(T_1^n, T_2^n, \ldots, T_{n_1}^n\right):$$
$$(V_1 \cup V_2 \cup \ldots \cup V_n) \to (W_1 \cup W_2 \cup \ldots \cup W_n)$$
where each $T_i : V_i \to W_i$ is
$$\left(T_1^i, T_2^i, \ldots, T_{n_i}^i\right) : V_i = \left(V_1^i, V_2^i, \ldots, V_{n_i}^i\right) \to W_i = \left(W_1^i, W_2^i, \ldots, W_{n_i}^i\right)$$
with $T_j^i : V_j^i \to W_j^i$; $1 \le j \le n_i$; $i = 1, 2, \ldots, n$; i.e., each $T_i: V_i \to W_i$ is a special set linear transformation of the special set vector space $V_i$ to special set vector space $W_i$ over S, for $i = 1, 2, \ldots, n$. Then we call $T = (T_1 \cup T_2 \cup \ldots \cup T_n): V \to W$ as the special set linear n-transformation.

$SHom_S (V, W) = \{SHom_S (V_1, W_1) \cup SHom_S (V_2, W_2) \cup \ldots \cup SHom_S (V_n, W_n)\}$ where $SHom_S (V_i, W_i) = \{ Hom_S \left(V_1^i, W_1^i\right), Hom_S \left(V_2^i, W_2^i\right), \ldots, Hom_S \left(V_{n_i}^i, W_{n_i}^i\right)\}$ true for $i = 1, 2, \ldots, n$. Thus $SHom_S = (V, W)$ is again a special set vector n space over the set S.

We illustrate this by a simple example.

*Example 2.3.12:* Let
$$V = (V_1 \cup V_2 \cup V_3 \cup V_4)$$
$$= \left(V_1^1, V_2^1\right) \cup \left(V_1^2, V_2^2, V_3^2\right) \cup \left(V_1^3, V_2^3\right) \cup \left(V_1^4, V_2^4, V_3^4, V_4^4\right)$$
be a special set 4-vector space over the set $S = Z^o \cup \{0\}$.
$$W = (W_1 \cup W_2 \cup W_3 \cup W_4)$$
$$= \left(W_1^1, W_2^1\right) \cup \left(W_1^2, W_2^2, W_3^2\right) \cup \left(W_1^3, W_2^3\right) \cup \left(W_1^4, W_2^4, W_3^4, W_4^4\right)$$
be a special set 4-vector space over the same set $S = Z^+ \cup \{0\}$.
Here $V_1 = \left(V_1^1, V_2^1\right)$ with
$$V_1^1 = \left\{ \begin{pmatrix} a & b \\ c & d \end{pmatrix} \middle| a, b, c, d \in Z^o \cup \{0\} \right\}$$

and $V_2^1 = S \times S$, such that $V_1$ is a special set vector space over $Z^o \cup \{0\}$. $V_2 = \left(V_1^2, V_2^2, V_3^2\right)$ where $V_1^2 = S \times S \times S$,

$$V_2^2 = \left\{ \begin{bmatrix} a_1 & a_2 & a_3 \\ a_4 & a_5 & a_6 \end{bmatrix} \middle| a_i \in Z^o \cup \{0\}, 1 \le i \le 6 \right\}$$

and

$$V_3^2 = \left\{ \begin{bmatrix} a_1 \\ a_2 \\ a_3 \end{bmatrix} \middle| a_i \in S, 1 \le i \le 3 \right\}$$

is a special set vector space over S. $V_3 = \left(V_1^3, V_2^3\right)$ where $V_1^3 =$ {All polynomials of degree less than or equal to 3 with coefficients from S} and

$$V_2^3 = \left\{ \begin{pmatrix} a & b & c \\ 0 & d & e \\ 0 & 0 & f \end{pmatrix} \middle| a, b, c, d, e, f \in S \right\};$$

$V_3$ is again a special set vector space over S.

Finally $V_4 = \left(V_1^4, V_2^4, V_3^4, V_4^4\right)$ is such that $V_1^4 = S \times S$

$$V_2^4 = \left\{ \begin{pmatrix} a & b \\ o & d \end{pmatrix} \middle| a, b, c, d \in Z^o \cup \{0\} \right\},$$

$V_3^4 =$ {all polynomials in x of degree less than or equal to 4 with coefficients from S} and

$$V_4^4 = \left\{ \begin{pmatrix} a & b & c & d \\ e & f & g & h \end{pmatrix} \middle| a, b, c, d, e, f, g, h \in S = Z^o \cup \{0\} \right\}.$$

$V_4$ is again a special set vector space over S.

Now $W_1 = \left(W_1^1, W_2^1\right)$ where $W_1^1 = S \times S \times S \times S$ and

$$W_2^1 = \left\{ \begin{pmatrix} a \\ b \end{pmatrix} \middle| a, b \in S \right\};$$

$W_1$ is special set vector space over $Z^+ \cup \{0\}$. $W_2 = \left(W_1^2, W_2^2, W_3^2\right)$ where

$$W_1^2 = \left\{ \begin{pmatrix} a & b \\ 0 & c \end{pmatrix} \middle| a, b, c \in Z^+ \cup \{0\} \right\},$$

$W_2^2 = S \times S \times S \times S \times S \times S$ and $W_3^2 = S \times S \times S$; $W$ is a special set vector space over $Z^+ \cup \{0\}$; $W_3 = \left(W_1^3, W_2^3\right)$ where $W_1^3 = S \times S \times S \times S$ and $W_2^3 =$ {all polynomials of degree less than or equal to 5 with coefficients from S}, so that $W_3$ is again a special set vector space over S. $W_4 = \left(W_1^4, W_2^4, W_3^4, W_4^4\right)$ where

$$W_1^4 = \left\{ \begin{pmatrix} a \\ b \end{pmatrix} \middle| a, b \in S \right\},$$

$$W_2^4 = \left\{ \begin{pmatrix} a & 0 \\ b & d \end{pmatrix} \middle| a, b, c, d \in S \right\},$$

$$W_3^4 = S \times S \times S \times S \times S,$$

$$W_4^4 = \left\{ \begin{bmatrix} a_1 & a_2 \\ a_3 & a_4 \\ a_5 & a_6 \\ a_7 & a_8 \end{bmatrix} \middle| a_i \in S, 1 \leq i \leq 8 \right\}.$$

$W_4$ is a special set vector space over S.
   Now define
$$T = (T_1 \cup T_2 \cup T_3 \cup T_4)$$

$$= \left(T_1^1, T_2^1\right) \cup \left(T_1^2, T_2^2, T_3^2\right) \cup \left(T_1^3, T_2^3\right) \cup \left(T_1^4, T_2^4, T_3^4, T_4^4\right):$$
$$V = (V_1 \cup V_2 \cup V_3 \cup V_4)$$
$$= \left(V_1^1, V_2^1\right) \cup \left(V_1^2, V_2^2, V_3^2\right) \cup \left(V_1^3, V_2^3\right) \cup \left(V_1^4, V_2^4, V_3^4, V_4^4\right)$$
$$\to W = (W_1 \cup W_2 \cup W_3 \cup W_4)$$
$$= \left(W_1^1, W_2^1\right) \cup \left(W_1^2, W_2^2, W_3^2\right) \cup \left(W_1^3, W_2^3\right) \cup \left(W_1^4, W_2^4, W_3^4, W_4^4\right)$$

as follows:
$$T_1^1 : V_1^1 \to W_1^1 \text{ defined by } T_1^1 \begin{pmatrix} a & b \\ c & d \end{pmatrix} = (a\ b\ c\ d).$$

$$T_2^1 : V_2^1 \to W_2^1 \text{ is given by } T_2^1(a\ \ b) = \begin{bmatrix} a \\ b \end{bmatrix};$$

thus $T_1 = \left(T_1^1, T_2^1\right) : V_1 \to W_1$ is a special set linear transformation over S.

Now define $T_2 = \left(T_1^2, T_2^2, T_3^2\right) : V_2 \to W_2$ by

$$T_1^2(a\ \ b\ \ c) = \begin{pmatrix} a & b \\ 0 & c \end{pmatrix}.$$

$$T_2^2 \begin{pmatrix} a_1 & a_2 & a_3 \\ a_4 & a_5 & a_6 \end{pmatrix} = (a_1, a_2, a_3, a_4, a_5, a_6),$$

$$T_3^2 \begin{pmatrix} a_1 \\ a_2 \\ a_3 \end{pmatrix} = (a_1, a_2, a_3).$$

Thus $T_2 : V_2 \to W_2$ is a special set linear transformation over S.

$$T_3 : V_3 \to W_3;$$
$$T_3 = \left(T_1^3, T_2^3\right) : \left(V_1^3, V_2^3\right) \to \left(W_1^3, W_2^3\right)$$

given by
$$T_1^3 : V_1^3 \to W_1^3$$
$$T_1^3 (a_0 + a_1 x + a_2 x^2 + a_3 x^3) = (a_0, a_1, a_2, a_3)$$

$T_2^3 : V_2^3 \to W_2^3$ is defined as

$$T_2^3 \begin{pmatrix} a_1 \\ a_2 \\ a_3 \end{pmatrix} = (a_1, a_2, a_3).$$

Thus $T_3: V_3 \to W_3$ is a special set linear transformation on S.

$T_4: V_4 \to W_4$ i.e.,
$$T_4 = \left(T_1^4, T_2^4, T_3^4, T_4^4\right) : \left(V_1^4, V_2^4, V_3^4, V_4^4\right) \to \left(W_1^4, W_2^4, W_3^4, W_4^4\right)$$
is given by

$$T_1^4 (a \quad b) = \begin{bmatrix} a \\ b \end{bmatrix};$$

$T_1^4 : V_1^4 \to W_1^4 \quad T_2^4 : V_2^4 \to W_2^4$

is given by

$$T_2^4 \begin{pmatrix} a & b \\ 0 & d \end{pmatrix} = \begin{pmatrix} a & 0 \\ b & d \end{pmatrix}.$$

$$T_3^4 : V_3^4 \to W_3^4$$

is such that
$$T_3^4 (a_0 + a_1x + a_2x^2 + a_3x^3 + a_4x^4) = (a_0, a_1, a_2, a_3, a_4)$$

and
$$T_4^4 : V_4^4 \to W_4^4$$

is such that

$$T_4^4 \begin{pmatrix} a & b & c & d \\ e & f & g & h \end{pmatrix} = \begin{bmatrix} a & e \\ b & f \\ c & g \\ d & h \end{bmatrix}.$$

Thus $T_4 = \left(T_1^4, T_2^4, T_3^4, T_4^4\right)$ is a special set linear transformation of $V_4$ to $W_4$ over S. Thus $T = (T_1, T_2, T_3, T_4): V = (V_1 \cup V_2 \cup V_3 \cup V_4) \to (W_1 \cup W_2 \cup W_3 \cup W_4) = W$ is a special set linear 4-transformation over the set S.

Now we proceed on to define the notion of special set linear operator of a special set vector space V over a set S.

**DEFINITION 2.3.7:** *Let*
$$V = (V_1 \cup \ldots \cup V_n)$$
$$= \left(V_1^1, V_2^1, \ldots, V_{n_1}^1\right) \cup \left(V_1^2, V_2^2, \ldots, V_{n_2}^2\right) \cup \ldots \cup \left(V_1^n, V_2^n, \ldots, V_{n_n}^n\right)$$
*be a special set n-vector space over the set S. Suppose $T = (T_1 \cup T_2 \cup \ldots \cup T_n)$ be a special set linear operator on V; i.e.,*
$$T = (T_1 \cup T_2 \cup \ldots \cup T_n)$$
$$= \left(T_1^1, T_2^1, \ldots, T_{n_1}^1\right) \cup \left(T_1^2, T_2^2, \ldots, T_{n_2}^2\right) \cup \ldots \cup \left(T_1^n, T_2^n, \ldots, T_{n_n}^n\right):$$
$$V = \{V_1 \cup V_2 \cup \ldots \cup V_n\}$$
$$= \left(V_1^1, V_2^1, \ldots, V_{n_1}^1\right) \cup \ldots \cup \left(V_1^n, V_2^n, \ldots, V_{n_n}^n\right) \rightarrow$$

$V = \{V_1 \cup V_2 \cup \ldots \cup V_n\}$ *is a special set linear n-operator on V if T is a special set linear transformation from V to V. Let $SHom_S(V,V) = \{SHom_S(V_1,V_1) \cup \ldots \cup SHom_S(V_n,V_n)\}$ where $SHom_S(V_i, V_i) = \{Hom_S\left(V_1^i, V_1^i\right), \ldots, Hom_S\left(V_{n_i}^i, V_{n_i}^i\right)\}$ for $i = 1, 2, \ldots, n$.*

*We can verify each $Hom_S\left(V_{n_i}^i, V_{n_i}^i\right)$ is a set linear algebra over S. Thus $SHom_S(V_i, V_i)$ is a special set linear algebra over the set S under the composition of mappings. Thus $SHom_S(V,V) = (SHom_S(V_1, V_1) \cup SHom_S(V_2, V_2) \cup \ldots \cup SHom_S(V_n, V_n)\}$ is a special set linear n-algebra over the set S.*

Now we illustrate this by a simple example.

*Example 2.3.13:* Let
$$V = (V_1 \cup V_2 \cup V_3)$$
$$= \left(V_1^1, V_2^1, V_3^1, V_4^1\right) \cup \left(V_1^2, V_2^2\right) \cup \left(V_1^3, V_2^3, V_3^3\right)$$
where
$$V_1^1 = \left\{\begin{pmatrix} a & b \\ c & d \end{pmatrix} \middle| a,b,c,d \in Z^+ \cup \{0\}\right\},$$

$V_2^1 = \{S \times S \times S \times S \times S \text{ where } S = Z^+ \cup \{0\}\}$,
$V_3^1 = \{S[x] \text{ all polynomials of degree less than or equal to four}\}$

and

$$V_4^1 = \left\{ \begin{pmatrix} a & a & a & a \\ a & a & a & a \end{pmatrix}, \begin{pmatrix} a & a \\ a & a \\ a & a \\ a & a \\ a & a \end{pmatrix} \middle| a \in Z^+ \cup \{0\} \right\},$$

$V_1 = \left( V_1^1, V_2^1, V_3^1, V_4^1 \right)$ is a special set vector space over the set S. Consider $V_2 = \left( V_1^2, V_2^2 \right)$ where

$$V_1^2 = \left\{ \begin{pmatrix} a & b & c \\ d & e & f \\ g & h & i \end{pmatrix} \middle| a,b,c,d,e,f,g,h,i \in S \right\},$$

$$V_2^2 = \left\{ [a\ a\ a\ a\ a\ a], \begin{bmatrix} a \\ a \\ a \\ a \\ a \\ a \\ a \end{bmatrix}, \begin{bmatrix} a & a \\ a & a \\ a & a \end{bmatrix}, \begin{pmatrix} a & a & a \\ a & a & a \end{pmatrix} \middle| a \in Z^+ \cup \{0\} = S \right\},$$

thus $V_2 = \left( V_1^2, V_2^2 \right)$ is a special set vector space over the set $S = Z^+ \cup \{0\}$. $V_3 = \left( V_1^3, V_2^3, V_3^3 \right)$ where

$$V_1^3 = \left\{ \begin{bmatrix} a_1 & a_2 & a_3 & a_4 \\ a_5 & a_6 & a_7 & a_8 \end{bmatrix} \middle| a_i \in S; 1 \leq i \leq 8 \right\}$$

$V_2^3 = \{$All polynomials in x with coefficients from S of degree less than or equal to 5$\}$ and $V_3^3 = \{$All upper triangular $4 \times 4$ matrices with entries from S and all lower triangular matrices with entries from S}. Hence $V_3 = \left(V_1^3, V_2^3, V_3^3\right)$ is a special set vector space over S. Thus $V = (V_1 \cup V_2 \cup V_3)$ is a special set vector 3-space over S.

Define
$$T = (T_1 \cup T_2 \cup T_3)$$
$$= \left(T_1^1, T_2^1, T_3^1, T_4^1\right) \cup \left(T_1^2, T_2^2\right) \cup \left(T_1^3, T_2^3, T_3^3\right)$$
from $V = (V_1 \cup V_2 \cup V_3)$ into $V = (V_1 \cup V_2 \cup V_3)$ as follows:
$T_1: V_1 = V_1$ such that
$$T_1 = \left(T_1^1, T_2^1, T_3^1, T_4^1\right) : \left(V_1^1, V_2^1, V_3^1, V_4^1\right) \text{ into } \left(V_1^1, V_2^1, V_3^1, V_4^1\right)$$
defined by
$T_1^1 : V_1^1 \to V_1^1$ where
$$T_1^1 \begin{pmatrix} a & b \\ c & d \end{pmatrix} = \begin{pmatrix} c & d \\ a & b \end{pmatrix}.$$

$T_2^1 : V_2^1 \to V_2^1$ such that
$$T_2^1 \, (a\ b\ c\ d\ e) = (d\ e\ c\ a\ b)$$

$T_3^1 : V_3^1 \to V_3^1$ by
$$T_3^1 \, (a_0 + a_1x + a_2x^2 + a_3x^3 + a_4x^4) = (a_0 + a_2x^2 + a_4x^4)$$
and $T_4^1 : V_4^1 \to V_4^1$ defined as
$$T_4^1 \begin{pmatrix} a & a & a & a \\ a & a & a & a \end{pmatrix} = \begin{bmatrix} a & a \\ a & a \\ a & a \\ a & a \end{bmatrix}$$

and
$$T_4^1 \begin{bmatrix} a & a \\ a & a \\ a & a \\ a & a \end{bmatrix} = \begin{pmatrix} a & a & a & a \\ a & a & a & a \end{pmatrix}.$$

Thus $T_1 = \left(T_1^1, T_2^1, T_3^1, T_4^1\right) : V_1 \to V_1$ is a special set linear operator on $V_1$. Now
$$T_2 = \left(T_1^2, T_2^2\right): V_2 = \left(V_1^2, V_2^2\right) \text{ into } V_2 = \left(V_1^2, V_2^2\right).$$

$T_2 : V_2 \to V_2$ is $T_1^2 : V_1^2 \to V_1^2$ is such that defined as

$$T_1^2 \begin{pmatrix} a & b & c \\ d & e & f \\ g & h & i \end{pmatrix} = \begin{pmatrix} a & b & c \\ 0 & e & f \\ 0 & 0 & i \end{pmatrix}$$

and $T_2^2 : V_2^2 \to V_2^2$ is given by

$$T_2^2 \begin{pmatrix} a & a & a & a & a & a \end{pmatrix} = \begin{pmatrix} a \\ a \\ a \\ a \\ a \\ a \end{pmatrix}$$

and

$$T_2^2 \begin{pmatrix} a \\ a \\ a \\ a \\ a \\ a \end{pmatrix} = \begin{pmatrix} a & a & a & a & a & a \end{pmatrix}.$$

$$T_2^2 \begin{pmatrix} a & a \\ a & a \\ a & a \end{pmatrix} = \begin{pmatrix} a & a & a \\ a & a & a \end{pmatrix}$$

and

$$T_1^2 \begin{pmatrix} a & a & a \\ a & a & a \end{pmatrix} = \begin{pmatrix} a & a \\ a & a \\ a & a \end{pmatrix}.$$

Thus $T_2 = \left( T_1^2, T_2^2 \right)$ is a special set linear operator on $V_2$.

Now consider $T_3 = \left( T_1^3, T_2^3, T_3^3 \right) : V_3 = \left( V_1^3, V_2^3, V_3^3 \right)$ into $V_3 = \left( V_1^3, V_2^3, V_3^3 \right)$ defined as follows.

$T_1^3 : V_1^3 \to V_1^3$

$$T_1^3 \begin{pmatrix} a_1 & a_2 & a_3 & a_4 \\ a_5 & a_6 & a_7 & a_8 \end{pmatrix} = \begin{pmatrix} a_5 & a_6 & a_7 & a_8 \\ a_1 & a_2 & a_3 & a_4 \end{pmatrix},$$

$T_2^3 : V_2^3 \to V_2^3$ by
$T_2^3 (a_0 + a_1x + a_2x^2 + a_3x^3 + a_4x^4 + a_5x^5) = (a_0 + a_1x + a_3x^3 + a_5x^5)$
and
$T_3^3 : V_3^3 \to V_3^3$ is given by

$$T_3^3 \begin{pmatrix} a & b & c & d \\ 0 & e & f & g \\ 0 & 0 & h & i \\ 0 & 0 & 0 & j \end{pmatrix} = \begin{pmatrix} a & 0 & 0 & 0 \\ b & e & 0 & 0 \\ c & f & h & 0 \\ d & g & i & j \end{pmatrix}$$

and

$$T_3^3 \begin{pmatrix} a & 0 & 0 & 0 \\ b & e & 0 & 0 \\ c & f & h & 0 \\ d & g & i & j \end{pmatrix} = \begin{pmatrix} a & b & c & d \\ 0 & e & f & g \\ 0 & 0 & h & i \\ 0 & 0 & 0 & j \end{pmatrix}$$

Thus $T_3 = \left( T_1^3, T_2^3, T_3^3 \right) : V_3 \to V_3$ is also a special set linear operator on $V_3$. Hence
$T = (T_1, T_2, T_3) = \left( T_1^1, T_2^1, T_3^1, T_4^1 \right) \cup \left( T_1^2, T_2^2 \right) \cup \left( T_1^3, T_2^3, T_3^3 \right)$
is a special set linear 3-operator on $V = (V_1 \cup V_2 \cup V_3)$.

Now we proceed to define the new notion of pseudo special set linear n-operator on V.

**DEFINITION 2.3.8:** *Let*
$$V = (V_1 \cup V_2 \cup ... \cup V_n)$$
$$= \left( V_1^1, V_2^1, ..., V_{n_1}^1 \right) \cup \left( V_1^2, V_2^2, ..., V_{n_2}^2 \right) \cup ... \cup \left( V_1^n, V_2^n, ..., V_{n_n}^n \right)$$
*be a special set n-vector space over the set S. Take*
$$T = (T_1 \cup T_2 \cup ... \cup T_n)$$
$$= \left( T_1^1, T_2^1, ..., T_{n_1}^1 \right) \cup \left( T_1^2, T_2^2, ..., T_{n_2}^2 \right) \cup ... \cup \left( T_1^n, T_2^n, ..., T_{n_n}^n \right):$$

$V = (V_1 \cup V_2 \cup ... \cup V_n) \rightarrow (V_1 \cup V_2 \cup ... \cup V_n)$
*such that T is a special set linear n-transformation on V; i.e., $T_i$: $V_i \rightarrow V_j$ where $i \neq j$ and $n_i \leq n_j$, j, i = 1, 2, ...,n.*

*Then we call $T = (T_1 \cup T_2 \cup ... \cup T_n)$ to be a pseudo special set linear n-operator on V. Further $SHom_S (V, V) = \{SHom_S (V_1, V_{j1}) \cup ... \cup SHom_S (V_n, V_{jn})\}$ here if $i \neq j_i$ then $n_i \leq n_{ji}, 1 \leq j_i \leq n$, i = 1, 2,..., n; where $SHom_S(V_t, V_{jt}) = \{$Special set linear transformation from $V_t$ into $V_{jt}\}$.*

$\{Hom_S(V_1^i, V_{t1}^i), Hom_S(V_2^i, V_{t2}^i), ..., Hom_S(V_{ni}^i, V_{tn_i}^i)\} = SHom_S(V_i, V_i)$ *where $1 \leq t_1, t_2, ..., t_{ni} \leq n_i$, this is true for each i, i = 1, 2, ..., n. So that $SHom_S(V, V)$ is a special set n-vector space over the set S.*

We now illustrate this definition by some examples.

**Example 2.3.14:** Let $V = (V_1 \cup V_2 \cup V_3 \cup V_4)$ be a special set vector 4-space over the set $S = \{0, 1\}$. Here

$$V_1 = \left(V_1^1, V_2^1, V_3^1\right),$$

$$V_2 = \left(V_1^2, V_2^2, V_3^2, V_4^2\right),$$

$$V_3 = \left(V_1^3, V_2^3, V_3^3\right)$$

and $V_4 = \left(V_1^4, V_2^4, V_3^4, V_4^4\right)$

are described in the following.

$$V_1^1 = \left\{ \begin{pmatrix} a & b & e \\ c & d & f \end{pmatrix} \middle| a,b,c,d,e,f \in S \right\}$$

$$V_2^1 = S \times S \times S \times S,$$

$V_3^1 = \{$All polynomials of degree less than or equal to 5 with entries from S$\}$.

$V_1 = \left(V_1^1, V_2^1, V_3^1\right)$ is a special set vector space over S.
Now

$$V_1^2 = \left\{ \begin{pmatrix} a & b \\ c & d \end{pmatrix} \middle| a,b,c,d \in Z^+ \cup \{0\} \right\},$$

$V_2^2 = \{(a\ a\ a\ a\ a), (a\ a\ a) \mid a \in Z^+ \cup \{0\}\}$; $V_3^2 = \{(Z^+ \cup \{0\})[x]$, polynomials of degree less than or equal to four$\}$ and $V_4^2 = \{$All $4 \times 4$ upper triangular matrices with entries from S$\}$; Clearly $V_2 = \left( V_1^2, V_2^2, V_3^2, V_4^2 \right)$ is a special set vector space over the set S. In $V_3$, take

$$V_1^3 = \left\{ \begin{pmatrix} a & b \\ c & d \end{pmatrix} \middle| a,b,c,d \in S \right\},$$

$$V_2^3 = \{S \times S \times S \times S \times S \times S\},$$

$$V_3^3 = \left\{ \begin{pmatrix} a & b \\ c & d \\ e & f \end{pmatrix} \middle| a,b,c,d,e,f \in S \right\}$$

so that $V_3$ is again a special set vector space over S. Consider $V_4 = \left( V_1^4, V_2^4, V_3^4, V_4^4 \right)$ here $V_1^4 = \{(Z^+ \cup \{0\})[x]$ all polynomials of degree less than or equal to 8$\}$,

$$V_2^4 = \left\{ \begin{bmatrix} a \\ a \\ a \\ a \end{bmatrix}, \begin{bmatrix} a \\ a \\ a \\ a \\ a \\ a \end{bmatrix} \middle| a \in Z^+ \cup \{0\} \right\},$$

$V_3^4 = \{Z^o \times Z^o \times Z^o \times Z^o$ where $Z^o = Z^+ \cup \{0\}\}$ and $V_4^4 = \{$All $4 \times 4$ lower triangular matrices with entries from S$\}$; $V_4 = \left( V_1^4, V_2^4, V_3^4, V_4^4 \right)$ is a special set vector space over the set S.

Now $V = (V_1 \cup V_2 \cup V_3 \cup V_4)$ is a special set vector 4-space over S. Define $T : V \to V$; i.e., $T = (T_1 \cup T_2 \cup T_3 \cup T_4) : V = (V_1 \cup V_2 \cup V_3 \cup V_4) \to (V_1 \cup V_2 \cup V_3 \cup V_4)$ by $T_1 : V_1 \to V_3$ that is

$$T_1 = \left(T_1^1, T_2^1, T_3^1\right) : \left(V_1^1, V_2^1, V_3^1\right) = V_1 \to \left(V_1^3, V_2^3, V_3^3\right)$$

defined by

$$T_1^1 : V_1^1 \to V_3^3$$
$$T_2^1 : V_2^1 \to V_1^3$$
$$T_3^1 : V_3^1 \to V_2^3$$

as follows

$$T_1^1 \begin{pmatrix} a & b & e \\ c & d & f \end{pmatrix} = \begin{bmatrix} a & b \\ c & d \\ e & f \end{bmatrix},$$

$$T_2^1 \begin{pmatrix} a & b & c & d \end{pmatrix} = \begin{bmatrix} a & b \\ c & d \end{bmatrix}$$

and

$$T_3^1 \left(a_0 + a_1 x + a_2 x^2 + a_3 x^3 + a_4 x^4 + a_5 x^5\right) = (a_0, a_1, a_2, a_3, a_4, a_5).$$

Thus $T_1 = \left(T_1^1, T_2^1, T_3^1\right) : V_1 \to V_3$ is a pseudo special set linear transformation of $V_1$ into $V_3$.

Define $T_2 : V_2 \to V_4$ where

$$T_2 = \left(T_1^2, T_2^2, T_3^2, T_4^2\right) :$$
$$V_2 = \left(V_1^2, V_2^2, V_3^2, V_4^2\right) \to \left(V_1^4, V_2^4, V_3^4, V_4^4\right) = V_4$$

as follows:

$T_1^2 : V_1^2 \to V_3^4$ as

$$T_1^2 \begin{pmatrix} a & b \\ c & d \end{pmatrix} = \begin{pmatrix} a & b & c & d \end{pmatrix}$$

and $T_2^2 : V_2^2 \to V_2^4$ by

$$T_2^2(a \ a \ a \ a \ a) = \begin{pmatrix} a \\ a \\ a \\ a \\ a \end{pmatrix}$$

and

$$T_2^2 (a \ a \ a) = \begin{bmatrix} a \\ a \\ a \end{bmatrix},$$

$T_3^2 : V_3^2 \to V_1^4$ by

$T_3^2 (a_0 + a_1x + a_2x^2 + a_3x^3 + a_4x^4) = (a_0 + a_1x^2 + a_2x^4 + a_3x^6 + a_4x^8)$

and $T_4^2 : V_4^2 \to V_4^4$ by

$$T_4^2 \begin{pmatrix} a & b & c & d \\ 0 & e & f & g \\ 0 & 0 & h & i \\ 0 & 0 & 0 & j \end{pmatrix} = \begin{pmatrix} a & 0 & 0 & 0 \\ b & e & 0 & 0 \\ c & f & h & 0 \\ d & g & i & j \end{pmatrix}.$$

Thus $T_2 = \left(T_1^2, T_2^2, T_3^2, T_4^2\right)$ is a pseudo special set linear transformation of $V_2$ into $V_4$.

Define $T_3 = \left(T_1^3, T_2^3, T_4^3\right) : V_3 \to V_1$ as follows.

$T_1^3 : V_1^3 \to V_2^1$ such that

$$T_1^3 \begin{pmatrix} a & b \\ c & d \end{pmatrix} = (a \ b \ c \ d),$$

$T_2^3 : V_2^3 \to V_1^1$ defined as

$$T_2^3 (a \ b \ c \ d \ e \ f) = \begin{pmatrix} a & b & e \\ c & d & f \end{pmatrix}$$

and $T_3^3 : V_3^3 \to V_3^1$ is defined by

$$T_3^3 \ [a_0 + a_1x + a_2x^2 + a_3x^3 + a_4x^4 + a_5x^5] = \begin{bmatrix} a_0 & a_1 \\ a_2 & a_3 \\ a_4 & a_5 \end{bmatrix}.$$

Thus $T_3 : V_3 \to V_1$ is again a pseudo special set linear transformation. Finally $T_4 : V_4 \to V_2$ is defined as follows.
$$T_4 = \left(T_1^4, T_2^4, T_3^4, T_4^4\right) :$$
$$V_4 = \left(V_1^4, V_2^4, V_3^4, V_4^4\right) \to V_2 = \left(V_1^2, V_2^2, V_3^2, V_4^2\right)$$
such that
$T_1^4 : V_1^4 \to V_3^2$ is defined by
$$T_1^4 \ (a_0 + a_1x + a_2x^2 + a_3x^3 + a_4x^4 + a_5x^5 + a_6x^6 + a_7x^7 + a_8x^8)$$
$$= (a_0 + a_1x + a_2x^2 + a_3x^3 + a_4x^4),$$

$T_2^4 : V_2^4 \to V_2^2$ is defined as
$$T_2^4 \begin{bmatrix} a \\ a \\ a \\ a \\ a \\ a \end{bmatrix} = \begin{pmatrix} a & a & a & a & a \end{pmatrix}$$
and
$$T_2^4 \begin{bmatrix} a \\ a \\ a \end{bmatrix} = (a\ a\ a),$$

$T_3^4 : V_3^4 \to V_1^2$ is defined by
$$T_3^4 \begin{pmatrix} a & b & c & d \end{pmatrix} = \begin{pmatrix} a & b \\ c & d \end{pmatrix}$$
and $T_4^4 : V_4^4 \to V_4^2$ is defined as

$$T_4^4 \begin{pmatrix} a & 0 & 0 & 0 \\ b & e & 0 & 0 \\ c & f & h & 0 \\ d & g & i & j \end{pmatrix} = \begin{pmatrix} a & b & c & d \\ 0 & e & f & g \\ 0 & 0 & h & i \\ 0 & 0 & 0 & j \end{pmatrix}.$$

Thus $T_4 = \left(T_1^4, T_2^4, T_3^4, T_4^4\right)$ is a pseudo special linear transformation of $V_4^4$ into $V_4^2$.

Hence
$$T = (T_1 \cup T_2 \cup T_3 \cup T_4)$$
$$= \left(T_1^1, T_2^1, T_3^1\right) \cup \left(T_1^2, T_2^2, T_3^2, T_4^2\right) \cup \left(T_1^3, T_2^3, T_4^3\right)$$
$$\cup \left(T_1^4, T_2^4, T_3^4, T_4^4\right)$$
$$: V = (V_1 \cup V_2 \cup V_3 \cup V_4)$$
$$= \left(V_1^1, V_2^1, V_3^1\right) \cup \left(V_1^2, V_2^2, V_3^2, V_4^2\right) \cup \left(V_1^3, V_2^3, V_4^3\right)$$
$$\cup \left(V_1^4, V_2^4, V_3^4, V_4^4\right)$$
$$\to (V_1 \cup V_2 \cup V_3 \cup V_4)$$

is a pseudo special linear 4-operator from V into V.

Now as in case of special set linear bioperators in case of special set linear n-operators also we can define n-projection and idempotent special set linear n-operators.

Now we proceed on to define the notion of special set direct sum and of special set vector n-subspaces of V.

**DEFINITION 2.3.9:** *Let*
$$V = (V_1 \cup V_2 \cup ... \cup V_n)$$
$$= \left(V_1^1, V_2^1, ..., V_{n_1}^1\right) \cup \left(V_1^2, V_2^2, ..., V_{n_2}^2\right) \cup ... \cup \left(V_1^n, V_2^n, ..., V_{n_n}^n\right)$$
*be a special set vector n space over the set S. Suppose each $V_{t_i}^i$ is a direct summand of subspaces $W_{t_1}^i \oplus ... \oplus W_{t_r}^i$ for each i = 1, 2, ..., n and $1 \le t_i \le n_i$ then we call $V = (V_1 \cup V_2 \cup ... \cup V_n)$ be the special set direct n union of subspaces.*

We illustrate this by a simple example.

*Example 2.3.15:* Let
$$V = (V_1 \cup V_2 \cup V_3)$$
$$= \left(V_1^1, V_2^1, V_3^1\right) \cup \left(V_1^2, V_2^2\right) \cup \left(V_1^3, V_2^3, V_3^3\right)$$
where

$$V_1^1 = \left\{ \begin{pmatrix} a & b \\ c & d \end{pmatrix} \middle| a,b,c,d \in Z^+ \cup \{0\} = S \right\},$$

$$V_2^1 = \left\{ (a\ a\ a), \begin{bmatrix} a \\ a \\ a \end{bmatrix}, [a\ a] \middle| a \in S = Z^+ \cup \{0\} \right\}$$

$$V_3^1 = \{S \times S \times S\},$$

$$V_1^2 = \left\{ (a\ a\ a\ a), \begin{bmatrix} a \\ a \end{bmatrix}, [a\ a\ a] \middle| a \in S \right\}$$

$$V_2^2 = \left\{ \begin{pmatrix} a_1 & a_2 & a_3 \\ a_4 & a_5 & a_6 \\ a_7 & a_8 & a_9 \end{pmatrix} \text{ such that } a_i \in S;\ 1 \le i \le 9 \right\},$$

$$V_1^3 = \{S \times S \times S \times S\},$$

$V_2^3 = \{S[x];$ polynomials of degree less than or equal to 2$\}$

and

$$V_3^3 = \left\{ \begin{pmatrix} a_1 \\ a_2 \\ a_3 \\ a_4 \end{pmatrix}, [a\ a\ a\ a\ a\ a], \begin{bmatrix} a & a \\ a & a \\ a & a \\ a & a \\ a & a \end{bmatrix} \middle| a, a_i \in S;\ 1 \le i \le 4 \right\}.$$

Now we express each $V_{t_i}^i$ as a direct union of vector subspaces.

$$V_1^1 = \left\{ \begin{pmatrix} a & b \\ c & d \end{pmatrix} \middle| a,b,c,d \in 3Z^+ \cup \{0\} \right\} \oplus$$

$$\left\{ \begin{pmatrix} a & b \\ c & d \end{pmatrix} \middle| a,b,c,d \in 5Z^+ \cup \{0\} \right\} \oplus$$

$$\left\{ \begin{pmatrix} a & b \\ c & d \end{pmatrix} \middle| a,b,c,d \in S \setminus \{5Z^+ \cup 3Z^+\} \right\} = \left( W_{11}^1 \oplus W_{21}^1 \oplus W_{31}^1 \right).$$

$$V_2^1 = \{(a \quad a \quad a)\} \oplus \left\{ \begin{bmatrix} a \\ a \\ a \end{bmatrix} \right\} \oplus \{[a \quad a]\}$$

$$= \left( W_{12}^1 \oplus W_{22}^1 \oplus W_{32}^1 \right) \ (a \in S).$$

$$V_3^1 = (\{S \times S \times \{0\} \cup \{0\}) \oplus (\{0\} \times \{0\} \times S \times S)$$
$$= \left( W_{11}^1 \oplus W_{21}^1 \oplus W_{31}^1 \right).$$

$$V_1^2 = \left\{ \{(a \quad a \quad a \quad a)\} \oplus \left\{ \begin{pmatrix} a \\ a \end{pmatrix} \right\} \oplus \{(a \quad a \quad a)\} \mid a \in S \right\}$$
$$= \left( W_{12}^1 \oplus W_{22}^1 \oplus W_{32}^1 \right)$$

$$V_2^2 = \left\{ \begin{pmatrix} a_1 & a_2 & a_3 \\ a_4 & a_5 & a_6 \\ a_7 & a_8 & a_9 \end{pmatrix} \middle| a_i \in 3Z^+ \cup \{0\}; \ 1 \le i \le 9 \right\} \oplus$$

$$\left\{ \begin{pmatrix} a_1 & a_2 & a_3 \\ a_4 & a_5 & a_6 \\ a_7 & a_8 & a_9 \end{pmatrix} \middle| a_i \in 7S \setminus 3Z^+ \cup \{0\}; \ 1 \le i \le 9 \right\}$$

$$= \left(W_{12}^2 \oplus W_{22}^2\right).$$

We represent
$$V_1^3 = \{S \times S \times \{0\} \times \{0\}\} \oplus \{\{0\} \times \{0\} \times S \times S\}$$
$$= \left(W_{13}^1 \oplus W_{23}^1\right)$$

$$V_2^3 = \{ax + b \mid a, b \in 3Z^+ \cup \{0\}\}$$
$$\oplus \{ax + b \mid a, b \in 5Z^+ \cup \{0\}\}$$
$$\oplus \{ax + b \mid a, b \in S - \{3Z^+ \cup 5Z^+\}\}$$
$$= \left(W_{13}^2 \oplus W_{23}^2 \oplus W_{33}^2\right).$$

$$V_3^3 = \left\{ \begin{pmatrix} a_1 \\ a_2 \\ a_3 \\ a_4 \end{pmatrix} \middle| a_i \in 9Z^+ \cup \{0\} \quad 1 \le i \le 4 \right\} \oplus$$

$$\left\{ \begin{pmatrix} a_1 \\ a_2 \\ a_3 \\ a_4 \end{pmatrix} \middle| a_i \in S \setminus 9Z^+ \cup \{0\} \quad 1 \le i \le 4 \right\} \oplus$$

$$\left\{ [a\ a\ a\ a\ a\ a] \middle| a \in S \right\} \oplus \left\{ \begin{bmatrix} a & a \\ a & a \\ a & a \\ a & a \end{bmatrix} \middle| a \in S \right\}$$

$$= \left(W_{13}^3 \oplus W_{23}^3 \oplus W_{33}^3 \oplus W_{43}^3\right).$$

Thus
$$V = V_1 \cup V_2 \cup V_3 \cup V_4$$
$$= \{\left(W_{11}^1 \oplus W_{21}^1 \oplus W_{31}^1, W_{11}^2 \oplus W_{21}^2 \oplus W_{31}^2, W_{11}^3 \oplus W_{21}^3 \oplus W_{31}^3\right) \cup$$
$$\left(W_{12}^1 \oplus W_{22}^1 \oplus W_{32}^1, W_{12}^2 \oplus W_{22}^2\right)$$

$$= \left(V_1^1, V_2^1, V_3^1\right) \cup \left(V_1^2, V_2^2\right) \cup \left(V_1^3, V_2^3, V_3^3\right).$$

Now it is pertinent to mention here that where $V = (V_1 \cup V_2 \cup \ldots \cup V_n)$ happen to be a special set linear n-algebra we can replace in the definition the special set direct n union by special set direct n-sum.

We shall represent this by one example when V is a special set linear n-algebra.

***Example 2.3.16***: Let $V = (V_1 \cup V_2 \cup V_3 \cup V_4)$ be a special set linear 4 algebra over the set $Z^+ \cup \{0\}$. Take
$$V_1 = \left(V_1^1, V_2^1, V_3^1, V_4^1\right),$$
$$V_2 = \left(V_1^2, V_2^2\right),$$
$$V_3 = \left(V_1^3, V_2^3, V_3^3\right)$$
and
$$V_4 = \left(V_1^4, V_2^4\right).$$

In $V_1 = \left(V_1^1, V_2^1, V_3^1, V_4^1\right)$,

$$V_1^1 = \left\{ \begin{pmatrix} a & b \\ c & d \end{pmatrix} \middle| a, b, c, d \in Z^+ \cup \{0\} \right\},$$

$$V_2^1 = \{ S \times S \times S \mid S = Z^+ \cup \{0\} \},$$

$$V_3^1 = \left\{ \begin{pmatrix} a_1 & a_2 \\ a_3 & a_4 \\ a_5 & a_6 \end{pmatrix} \middle| a_i \in S; 1 \le i \le 6 \right\}$$

and

$$V_4^1 = \left\{ \begin{pmatrix} a_1 & a_2 & a_3 \\ 0 & a_4 & a_5 \\ 0 & 0 & a_6 \end{pmatrix} \middle| a_i \in S; 1 \le i \le 6 \right\}.$$

Now we see $V_1$ is a special set linear algebra over the set S. Here we write $V_1$ as a special set direct summand as follows.

$$V_1^1 = \left\{\left\langle \begin{pmatrix} 1 & 0 \\ 0 & 0 \end{pmatrix} \begin{pmatrix} 0 & 0 \\ 1 & 0 \end{pmatrix} \right\rangle\right\} \oplus \left\{\left\langle \begin{pmatrix} 0 & 1 \\ 0 & 0 \end{pmatrix} \right\rangle\right\} \oplus \left\{\begin{pmatrix} 0 & 0 \\ 0 & 1 \end{pmatrix}\right\}$$

$$= \left(W_{11}^1 \oplus W_{21}^1 \oplus W_{31}^1\right)$$

is a special set direct summand of $V_1^1$.

$$V_1^2 = \{\langle(1, 00)\rangle\} \oplus \{\langle(0, 11)\rangle\} = \left(W_{12}^2 \oplus W_{22}^2\right)$$

is a special set direct summand of $V_2^1$.

$$V_3^1 = \left\{\left\langle \begin{pmatrix} 1 & 0 \\ 0 & 0 \\ 0 & 0 \end{pmatrix} \right\rangle\right\} \oplus \left\{\left\langle \begin{pmatrix} 0 & 1 \\ 0 & 0 \\ 0 & 0 \end{pmatrix} \right\rangle\right\} \oplus \left\{\left\langle \begin{pmatrix} 0 & 0 \\ 1 & 0 \\ 0 & 0 \end{pmatrix} \right\rangle\right\} \oplus$$

$$\left\{\left\langle \begin{pmatrix} 0 & 0 \\ 0 & 1 \\ 0 & 0 \end{pmatrix} \right\rangle\right\} \oplus \left\{\left\langle \begin{pmatrix} 0 & 0 \\ 0 & 0 \\ 1 & 0 \end{pmatrix} \right\rangle\right\} \oplus \left\{\left\langle \begin{pmatrix} 0 & 0 \\ 0 & 0 \\ 0 & 1 \end{pmatrix} \right\rangle\right\}$$

$$= \left(W_{13}^3 \oplus W_{23}^3 \oplus W_{33}^3 \oplus W_{43}^3 \oplus W_{53}^3 \oplus W_{63}^3\right)$$

is a special set direct summand of $V_3^1$.
Finally

$$V_4^1 = \left\{\left\langle \begin{bmatrix} 1 & 0 & 0 \\ 0 & 0 & 0 \\ 0 & 0 & 0 \end{bmatrix} \right\rangle\right\} \oplus \left\{\left\langle \begin{bmatrix} 0 & 1 & 0 \\ 0 & 0 & 0 \\ 0 & 0 & 0 \end{bmatrix} \right\rangle\right\} \oplus$$

$$\left\{\left\langle \begin{bmatrix} 0 & 0 & 0 \\ 0 & 0 & 0 \\ 0 & 0 & 1 \end{bmatrix} \right\rangle\right\} \oplus \left\{\left\langle \begin{bmatrix} 0 & 0 & 0 \\ 0 & 1 & 0 \\ 0 & 0 & 0 \end{bmatrix} \right\rangle\right\} \oplus$$

$$\left\{\left\langle\begin{bmatrix} 0 & 0 & 0 \\ 0 & 1 & 0 \\ 0 & 0 & 0 \end{bmatrix}\right\rangle\right\} \oplus \left\{\left\langle\begin{bmatrix} 0 & 0 & 0 \\ 0 & 0 & 1 \\ 0 & 0 & 0 \end{bmatrix}\right\rangle\right\}$$

$$= \left(W_{14}^4 \oplus W_{24}^4 \oplus W_{34}^4 \oplus W_{44}^4 \oplus W_{54}^4 \oplus W_{64}^4\right)$$

is a special set direct summand of $V_4^1$.
Thus
$$V_1 = \left(V_1^1, V_2^1, V_3^1, V_4^1\right)$$
$$= \left(\left(W_{11}^1 \oplus W_{21}^1 \oplus W_{31}^1\right), \left(W_{12}^2 \oplus W_{22}^2\right),\right.$$
$$\left(W_{13}^3 \oplus W_{23}^3 \oplus W_{33}^3 \oplus W_{43}^3 \oplus W_{53}^3 \oplus W_{63}^3\right),$$
$$\left.\left(W_{14}^4 \oplus W_{24}^4 \oplus W_{34}^4 \oplus W_{44}^4 \oplus W_{54}^4 \oplus W_{64}^4\right)\right) \quad \ldots \text{ I.}$$

Consider $V_2 = \left(V_1^2, V_2^2\right)$ where

$$V_1^2 = \left\{\begin{pmatrix} a & a_1 \\ a & a_1 \end{pmatrix} \middle| a_1, a \in S = Z^+ \cup \{0\}\right\}$$

and

$$V_2^2 = \left\{\begin{pmatrix} a & a & a_2 & a_1 \\ a & a & a_2 & a_1 \\ a & a & a_2 & a_1 \end{pmatrix} \middle| a_1, a, a_2 \in Z^+ \cup \{0\} = S\right\}$$

be the special set vector space over S.
Now
$$V_1^2 = \left\{\left\langle\begin{pmatrix} 1 & 0 \\ 1 & 0 \end{pmatrix}\right\rangle\right\} \oplus \left\{\left\langle\begin{pmatrix} 0 & 1 \\ 0 & 1 \end{pmatrix}\right\rangle\right\}$$
$$= \left(W_{12}^1 \oplus W_{22}^1\right)$$

and

$$V_2^2 = \left\{ \left\langle \begin{bmatrix} 1 & 1 & 0 & 0 \\ 1 & 1 & 0 & 0 \\ 1 & 1 & 0 & 0 \end{bmatrix} \right\rangle \right\} \oplus \left\{ \left\langle \begin{bmatrix} 0 & 0 & 1 & 0 \\ 0 & 0 & 1 & 0 \\ 0 & 0 & 1 & 0 \end{bmatrix} \right\rangle \right\}$$

$$\oplus \left\{ \left\langle \begin{bmatrix} 0 & 0 & 0 & 1 \\ 0 & 0 & 0 & 1 \\ 0 & 0 & 0 & 1 \end{bmatrix} \right\rangle \right\}$$

$$= \left( W_{13}^2 \oplus W_{23}^2 \oplus W_{33}^2 \right).$$

Thus

$$V_2 = \left( V_1^2, V_2^2 \right) = \left\{ \left( W_{12}^1 \oplus W_{22}^1 \right), \left( W_{13}^2 \oplus W_{23}^2 \oplus W_{33}^2 \right) \right\} \quad \ldots \quad \text{II}$$

is a special set direct summand of the special set linear algebra $V_2$. Now $V_3 = \left( V_1^3, V_2^3, V_3^3 \right)$ where

$$V_1^3 = \left\{ \begin{pmatrix} a & a & a & a_1 \\ a & a & a & a_1 \end{pmatrix} \middle| a, a_1 \in S \right\},$$

$$V_2^3 = \left\{ \begin{bmatrix} a & a & a_1 \\ a & a & a_1 \\ a & a & a_1 \\ a & a & a_1 \end{bmatrix} \middle| a, a_1 \in S \right\}$$

and $V_3^3 = \{S \times S \times S \times S\}$ be the special set vector space over S.

$$V_1^3 = \left\{ \left\langle \begin{bmatrix} 1 & 1 & 1 & 0 \\ 1 & 1 & 1 & 0 \end{bmatrix} \right\rangle \right\} \oplus \left\{ \left\langle \begin{bmatrix} 0 & 0 & 0 & 1 \\ 0 & 0 & 0 & 1 \end{bmatrix} \right\rangle \right\}$$

$$= \left( W_{13}^1 \oplus W_{23}^1 \right)$$

$$V_2^3 = \left\{ \left\langle \begin{bmatrix} 1 & 1 & 0 \\ 1 & 1 & 0 \\ 1 & 1 & 0 \\ 1 & 1 & 0 \end{bmatrix} \right\rangle \right\} \oplus \left\{ \left\langle \begin{bmatrix} 0 & 0 & 1 \\ 0 & 0 & 1 \\ 0 & 0 & 1 \\ 0 & 0 & 1 \end{bmatrix} \right\rangle \right\}$$

$$= \left( W_{23}^2 \oplus W_{33}^2 \right)$$

and

$$V_3^3 = \{ S \times S \times \{0\} \times \{0\} \} \oplus \{ \{0\} \times \{0\} \times S \times S \}$$
$$= \left( W_{13}^3 \oplus W_{23}^3 \right).$$

Thus $V_3 = \left( V_1^3, V_2^3, V_3^3 \right)$

$$= \left\{ \left( W_{13}^1 \oplus W_{23}^1 \right), \left( W_{23}^2 \oplus W_{33}^2 \right), \left( W_{13}^3 \oplus W_{23}^3 \right) \right\} \quad \ldots \text{III}$$

is a special set direct summand of $V_3$. Hence using I, II and III we get

$$V = (V_1 \cup V_2 \cup V_3)$$
$$= \{ \left( W_{11}^1 \oplus W_{21}^1 \oplus W_{31}^1 \right), \left( W_{12}^2 \oplus W_{22}^2 \right),$$
$$\left( W_{13}^3 \oplus W_{23}^3 \oplus \ldots \oplus W_{63}^3 \right), \left( W_{14}^4 \oplus W_{24}^4 \oplus W_{34}^4 \oplus \ldots \oplus W_{64}^4 \right) \}$$
$$\cup \{ \left( W_{12}^1 \oplus W_{22}^1 \right), \left( W_{13}^2 \oplus W_{23}^2 \oplus W_{33}^2 \right) \} \cup$$
$$\{ \left( W_{13}^1 \oplus W_{23}^1 \right), \left( W_{23}^2 \oplus W_{33}^2 \right), \left( W_{13}^3 \oplus W_{23}^3 \right) \}$$

to be special set 3-direct summand of V over S.

Now it is pertinent here to mention the following:

1. When we use special set vector n-space the problem of representation in terms of generating sets or a direct union of subspaces are no doubt tedious when compared with special set linear n-algebras. Special set linear n-algebras yield easy direct summand and a small or manageable n-generating subsets.
2. But yet we cannot rule out the use of special set vector n-spaces or completely replace special set n-spaces for we need them when the component subspaces of $V = (V_1, \ldots,$ $V_n$) takes the form as $V_{t_i}^i = \{(1\ 1\ 1\ 1\ 1), (1\ 1\ 1), (0\ 0\ 0), (0\ 0\ 0\ 0), (1\ 0\ 1), (1\ 0), (0\ 1), (0\ 0)\}$ or of the form only a few polynomials; i.e., a finite set in such cases we see special set vector n-spaces are preferable to special set linear n-algebras.

## 2.4 Special Set Fuzzy Vector Spaces

Now in the following section. We proceed onto define the notion of fuzzy analogue of these new types of special set vector spaces.

**DEFINITION 2.4.1:** *Let $V = \{S_1, S_2, ..., S_n\}$ be a special set vector space over the set p. If $\eta = (\eta_1, ..., \eta_n)$ is a n-map from V into [0, 1] where each $\eta_i: S_i \to [0,1]$, $i = 1, 2, ..., n$; $\eta_i(ra_i) \geq \eta_i(a_i)$ for all $a_i \in S_i$ and $r \in p$ then we call $V_\eta = V_{(\eta_1, \eta_2, ..., \eta_n)} = (V_1, V_2, ..., V_n)_\eta = (V_1\eta_1, V_2\eta_2, ..., V_n\eta_n)$ to be a special set fuzzy vector space.*

We illustrate this by an example.

*Example 2.4.1:* Let $V = (V_1, V_2, V_3, V_4)$ where

$$V_1 = \left\{ \begin{pmatrix} a & b \\ c & d \end{pmatrix} \middle| a, b, c, d \in Z_2 \right\},$$

$$V_2 = \{(Z^o \times Z^o \times Z^o); Z^o = Z^+ \cup \{0\}\},$$

$$V_3 = \left\{ \begin{pmatrix} a & a & a & a \\ a & a & a & a \end{pmatrix} \middle| a \in Z^+ \cup \{0\} \right\}$$

and $V_4 = \{$all polynomials of finite degree over $Z_2\}$ is a special set vector space over the set $S = \{0, 1\}$.

Define $\eta = (\eta_1, \eta_2, \eta_3, \eta_4) : V = (V_1, V_2, V_3, V_4) \to [0, 1]$ such that each $\eta_i: V_i \to [0, 1]$; $i = 1, 2, 3, 4$ with

$\eta_1: V_1 \to [0, 1]$ is defined by

$$\eta_1 \begin{pmatrix} a & b \\ c & d \end{pmatrix} = \begin{cases} \dfrac{1}{(ad-bc)} & \text{if } ad-bc \neq 0 \\ 1 & \text{if } ad-bc = 0 \end{cases}$$

$\eta_2: V_2 \to [0,1]$

$$\eta_2(a, b, c) = \begin{cases} \dfrac{1}{a+b+c} & \text{if } a+b+c \neq 0 \\ 1 & \text{if } a+b+c = 0 \end{cases}$$

$\eta_3: V_3 \to [0,1]$

$$\eta_3 \begin{bmatrix} a & a & a & a \\ a & a & a & a \end{bmatrix} = \begin{cases} \dfrac{1}{a} & \text{if } a \neq 0 \\ 1 & \text{if } a = 0 \end{cases}$$

$\eta_4: V_4 \to [0,1]$

$$\eta_4(p(x)) = \begin{cases} \dfrac{1}{\deg p(x)} & \text{if } \deg p(x) \neq 0 \\ 1 & \text{if } p(x) = 0 \text{ or a constant polynomial} \end{cases}$$

Clearly

$$V_\eta = (V_1, V_2, V_3, V_4)_{(\eta_1, \eta_2, \eta_3, \eta_4)} = (V_1\eta_1, V_2\eta_2, V_3\eta_3, V_4\eta_4)$$

is a special set fuzzy vector space.

We define the notion of special set fuzzy vector subspace.

**DEFINITION 2.4.2** *Let $V = \{V_1, ..., V_n\}$ be a special set vector space over the set $S$. Let $W = (W_1, W_2, ..., W_n) \subseteq (V_1, V_2, ..., V_n) = V$, that is each $W_i \subseteq V_i$; $i = 1, 2, ..., n$ be a special set vector subspace of $V$. Define $\eta: W \to [0, 1]$, i.e., $\eta = (\eta_1, \eta_2, ..., \eta_n)$: $(W_1, W_2, ..., W_n) \to [0,1]$; that is $\eta_i: W_i \to [0,1]$, for $i = 1, 2, ..., n$ such that $\eta_i(sa_i) \geq \eta_i(a_i)$ for all $a_i \in W_i$ for all $s \in S$; $1 \leq i \leq n$. Then $W_\eta = (W_1, W_2, ..., W_n)_{(\eta_1, \eta_2, ...\eta_n)} = (W_1\eta_1, W_2\eta_2, ..., W_n\eta_n)$ is a special set fuzzy vector subspace.*

We shall illustrate this by a simple example.

**Example 2.4.2:** Let $V = (V_1, V_2, \ldots, V_5)$ where

$$V_1 = \left\{ \begin{pmatrix} a & b & c \\ 0 & d & e \\ 0 & 0 & f \end{pmatrix} \middle| a,b,c,d,e,f \in Z^+ \cup \{0\} \right\},$$

$V_2 = \{(a\ a\ a\ a\ a) \mid a \in Z^+ \cup \{0\}\}$, $V_3 = \{S[x]$, all polynomials of degree less than or equal to if with coefficients from $S\}$,

$$V_4 = \left\{ \begin{bmatrix} a_1 \\ a_2 \\ a_3 \\ a_4 \\ a_5 \end{bmatrix} \begin{bmatrix} a_1 & a_2 & a_3 & a_4 \end{bmatrix} \middle| a_i \in Z^+ \cup \{0\}; 1 \leq i \leq 5 \right\}$$

and $V_5 = \{S \times S \times S \times S \times S \times S\}$ is a special set vector space over the set $P = \{0, 1\}$.

Let $W = (W_1, W_2, W_3, W_4, W_5) \subseteq (V_1, V_2, V_3, V_4, V_5) = V$ be a proper subset of V where

$$W_1 = \left\{ \begin{pmatrix} a & a & a \\ 0 & a & a \\ 0 & 0 & a \end{pmatrix} \middle| a \in S \right\} \subseteq V_1,$$

$W_2 = \{(a\ a\ a\ a\ a) \mid a \in 7Z^+ \cup \{0\}\} \subseteq V_2$,
$W_3 = \{$All polynomials of degree less than or equal to 26 with coefficients from $9Z^+ \cup \{0\}\} \subseteq V_3$,

$$W_4 = \left\{ \begin{bmatrix} a_1 \\ a_2 \\ a_3 \\ a_4 \\ a_5 \end{bmatrix} \middle| a_i \in S; \ 1 \leq i \leq 5 \right\} \subseteq V_4$$

and $W_5 = \{S \times S \times \{0\} \times S \times S \times \{0\}\} \subseteq V_5$. W is clearly a special set vector subspace of V over the set P.

Define $\eta : W \to [0,1]$ as follows.

$$\eta = (\eta_1, \eta_2, \eta_3, \eta_4, \eta_5) : W = (W_1, W_2, W_3, W_4, W_5) \to [0,1]$$

such that each $\eta_i : W_i \to [0,1]$ ; $i = 1, 2, 3, 4, 5$.

$\eta_1 : W_1 \to [0,1]$ such that

$$\eta_1 \begin{pmatrix} a & a & a \\ 0 & a & a \\ 0 & 0 & a \end{pmatrix} = \begin{cases} \dfrac{1}{a} & \text{if } a \neq 0 \\ 1 & \text{if } a = 0 \end{cases}$$

$\eta_2 : W_2 \to [0,1]$ such that

$$\eta_2 [a\ a\ a\ a\ a] = \begin{cases} \dfrac{1}{5a} & \text{if } a \neq 0 \\ 1 & \text{if } a = 0 \end{cases}$$

$\eta_3 : W_3 \to [0, 1]$

$$\eta_3(p(x)) = \begin{cases} \dfrac{1}{\deg p(x)} & \text{if } \deg p(x) \neq 0 \\ 1 & \text{if } \deg p(x) = 0 \text{ or a constant polynomial} \end{cases}$$

$\eta_4 : W_4 \to [0,1]$

$$\eta_4 \left( \begin{bmatrix} a_1 \\ a_2 \\ a_3 \\ a_4 \\ a_5 \end{bmatrix} \right) = \begin{cases} \dfrac{1}{a_1 + a_2 + a_3 + a_4 + a_5} & \text{if } \sum_{i=1}^{5} a_i \neq 0 \\ 1 & \text{if } \sum_{i=1}^{5} a_i = 0 \end{cases}$$

$\eta_5: W_5 \to [0, 1]$

$$\eta_5(a\ b\ 0\ c\ d\ 0) = \begin{cases} \dfrac{1}{a+b} & \text{if } b+a \neq 0 \\ 1 & \text{if } b+a = 0 \end{cases}$$

Clearly
$$W_\eta = \left(W_1, W_2, W_3, W_4, W_5\right)_{(\eta_1, \eta_2, \eta_3, \eta_4, \eta_5)}$$
$$= (W_1\eta_1, W_2\eta_2, W_3\eta_3, W_4\eta_4, W_5\eta_5)$$

is a special set fuzzy vector subspace.

Now we proceed on to define the new notion of special set fuzzy linear algebra.

**DEFINITION 2.4.3:** *Let $V = (V_1, V_2, ..., V_n)$ be a special set linear algebra over the set S. Take a n-map $\eta = (\eta_1, \eta_2, ..., \eta_n)$ : $V = (V_1, V_2, ..., V_n) \to [0,1]$ such that $\eta_i : V_i \to [0,1]$ satisfying the condition $\eta_i(a v_i) \geq \eta_i(v_i)$ for all $a \in S$ and $v_i \in V$; $1 \leq i \leq n$. Then we call $(V_1\eta_1, V_2\eta_2, ..., V_n\eta_n)$ to be a special set fuzzy linear algebra.*

It is important to note that both special set fuzzy vector space and special set fuzzy linear algebra are the same. Thus the concept of fuzziness makes them identical or no difference exist. We call such algebraic structures which are distinct but become identical due to fuzzyfying them are said to be fuzzy equivalent. Thus we see special set vector spaces and special set linear algebras are fuzzy equivalent.

We shall exhibit the above definition by an example.

*Example 2.4.3:* Let $V = (V_1, V_2, V_3, V_4)$ be a special set linear algebra over the set $\{0, 1\} = S$, where

$$V_1 = \left\{ \begin{pmatrix} a_1 & a_2 & a_3 & a_4 \\ a_5 & a_6 & a_7 & a_8 \\ a_9 & a_{10} & a_{11} & a_{12} \end{pmatrix} \middle| a_i \in \{0,1\}; 1 \leq i \leq 12 \right\},$$

$$V_2 = (S \times S \times S \times S \times S);$$

$$V_3 = \left\{ \begin{bmatrix} a_1 & a_2 \\ a_3 & a_4 \\ a_5 & a_6 \\ a_7 & a_8 \end{bmatrix} \middle| a_i \in \{0, 1\}; 1 \leq i \leq 8 \right\}$$

and $V_4 = \{$all $4 \times 4$ matrices with entries from $\{0,1\}\}$.

Define $\eta = (\eta_1, \eta_2, \eta_3, \eta_4) : V = (V_1, V_2, V_3, V_4) \to [0,1]$ such that each $\eta_i : V_i \to [0,1]$ ; $i = 1, 2, 3, 4$ as follows.

$\eta_1 : V_1 \to [0,1]$

$$\eta_1\left(\begin{bmatrix} a_1 & a_2 & a_3 & a_4 \\ a_5 & a_6 & a_7 & a_8 \\ a_9 & a_{10} & a_{11} & a_{12} \end{bmatrix}\right) = \begin{cases} \dfrac{1}{\sum_{i=1}^{12} a_i} & \text{if } \sum_{i=1}^{12} a_i \neq 0 \\ 1 & \text{if } \sum_{i=1}^{12} a_i = 0 \end{cases}$$

$\eta_2 : V_2 \to [0, 1]$ defined by

$$\eta_2(a, b, c, d, e) = \begin{cases} \dfrac{1}{a+e} & \text{if } e+a \neq 0 \\ 1 & \text{if } e+a = 0 \end{cases}$$

$\eta_3 : V_3 \to [0, 1]$ is such that

$$\eta_3 \begin{bmatrix} a_1 & a_2 \\ a_3 & a_4 \\ a_5 & a_6 \\ a_7 & a_8 \end{bmatrix} = \begin{cases} \dfrac{1}{\sum_{i=1}^{8} a_i} & \text{if } \sum_i a_i \neq 0 \\ 1 & \text{if } \sum_i a_i = 0 \end{cases}$$

$\eta_4 : V_4 \to [0, 1]$

$$\eta_4(A) = \begin{cases} \dfrac{1}{|A|} & \text{if } |A| \neq 0 \\ 1 & \text{if } |A| = 0 \end{cases}$$

Here A is a 4×4 matrix. Thus
$$V_\eta = \left(V_1, V_2, V_3, V_4\right)_{(\eta_1, \eta_2, \eta_3, \eta_4)} = (V_1\eta_1, V_2\eta_2, V_3\eta_3, V_4\eta_4)$$
is a special set fuzzy linear algebra.

This is the same as a special set fuzzy vector space.

Now we proceed onto define the notion of special set fuzzy linear subalgebra.

**DEFINITION 2.4.4:** *Let $V = (V_1, V_2, ..., V_n)$ be a special set linear algebra over the set S. Let $W = (W_1, W_2, ..., W_n) \subseteq (V_1, V_2, ..., V_n)$ where W is a special set linear subalgebra of V over the set S. Define $\eta : W = (W_1, W_2, ..., W_n) \to [0,1]$ as follows:*
$$\eta = (\eta_1, \eta_2, ..., \eta_n) : (W_1, W_2, ..., W_n) \to [0,1]$$
*such that $\eta_i : W_i \to [0,1]$, $i=1,2,...,n$. $W_\eta = \left(W_1, W_2, ..., W_n\right)_\eta =$*
*$(W_1\eta_1, W_2\eta_2, ..., W_n\eta_n)$ is a special set fuzzy linear subalgebra.*

Let us now proceed onto define the notion of special set fuzzy bispaces.

**DEFINITION 2.4.5:** *Let $V = \left(V_1^1, V_2^1, ..., V_{n_1}^1\right) \cup \left(V_1^2, V_2^2, ..., V_{n_2}^2\right) = V_1 \cup V_2$ be a special set bivector space over the set S. Let the bifuzzy map $\eta = \eta_1 \cup \eta_2 = \left(\eta_1^1, \eta_2^1, ..., \eta_{n_1}^1\right) \cup \left(\eta_1^2, \eta_2^2, ..., \eta_{n_2}^2\right)$: $V_1 \cup V_2 \to [0,1]$ be defined as $\eta_i^1 : V_i^1 \to [0,1]$; $1 \leq i \leq n_1$ and $\eta_j^2 : V_j^2 \to [0, 1]$, $j = 1, 2, ..., n_2$ satisfy the condition $\eta_i^1(ra_i^1) \geq \eta_i^1(a_i^1)$ and $\eta_j^2(ra_j^2) \geq \eta_j^2(a_j^2)$ for all $r \in S$ and $a_i^1 \in V_i^1$, $a_j^2 \in V_j^2$, $i = 1, 2, ..., n_1$ and $j = 1, 2, ..., n_2$. Thus*
$$V_\eta = (V_1 \cup V_2)_{\eta_1 \cup \eta_2}$$
$$= (V_1^1, V_2^1, ..., V_{n_1}^1)_{\left(\eta_1^1, \eta_2^1, ..., \eta_{n_1}^1\right)} \cup (V_1^2, V_2^2, ..., V_{n_2}^2)_{\left(\eta_1^2, \eta_2^2, ..., \eta_{n_2}^2\right)}$$

$$= (V^1_{1\eta^1_1}, V^1_{2\eta^1_2},...,V^1_{n_1\eta^1_{n_1}}) \cup (V^2_{1\eta^2_1},...,V^2_{n_2\eta^2_{n_2}})$$

*is a special set fuzzy vector bispace.*

We illustrate this by the following example.

***Example 2.4.4:*** Let
$$V = (V_1 \cup V_2) = (V^1_1, V^1_2, V^1_3) \cup (V^2_1, V^2_2, V^2_3, V^2_4)$$
be a special set vector bispace over the set $S = \{0, 1\}$. Here $V_1 = (V^1_1, V^1_2, V^1_3)$ is given by

$$V^1_1 = \left\{ \begin{pmatrix} a & b \\ c & d \end{pmatrix} \middle| a,b,c,d \in Z^+ \right\},$$

$$V^1_2 = \{[a\ a\ a\ a\ a] \mid a \in Z^o\}$$

and

$$V^1_3 = \left\{ \begin{pmatrix} a_1 & a_5 \\ a_2 & a_6 \\ a_3 & a_7 \\ a_4 & a_8 \end{pmatrix} \middle| a_i \in Z^o; 1 \le i \le 8 \right\},$$

$V_2 = (V^2_1, V^2_2, V^2_3, V^2_4)$ where

$V^2_1 = \{3 \times 3$ upper triangular matrices with entries from $Z^o\}$,
$$V^2_2 = \{Z^o \times Z^o \times Z^o\},$$
$$V^2_3 = \{[a\ a\ a\ a\ a\ a] \mid a \in Z^o\}$$

and

$$V^2_4 = \left\{ \begin{pmatrix} a_1 \\ a_2 \\ a_3 \\ a_4 \\ a_5 \\ a_6 \end{pmatrix} \middle| a_i \in Z^o; 1 \le i \le 6 \right\}.$$

Define the bifuzzy set
$$\begin{aligned}\eta &= \eta_1 \cup \eta_2 \\ &= \left(\eta_1^1, \eta_2^1, \eta_3^1\right) \cup \left(\eta_1^2, \eta_2^2, \eta_3^2, \eta_4^2\right) \\ &= (V_1 \cup V_2) \\ &= \left(V_1^1, V_2^1, V_3^1\right) \cup \left(V_1^2, V_2^2, V_3^2, V_4^2\right)\end{aligned}$$

as follows :
$\eta_1^1 : V_1^1 \to [0, 1]$ defined by

$$\eta_1^1 \begin{pmatrix} a & b \\ c & d \end{pmatrix} = \begin{cases} \dfrac{1}{ad} & \text{if } ad \neq 0 \\ 1 & \text{if } ad = 0 \end{cases}$$

$\eta_2^1 : V_2^1 \to [0, 1]$ is defined by

$$\eta_2^1 [a\ a\ a\ a\ a] = \begin{cases} \dfrac{1}{5a} & \text{if } a \neq 0 \\ 1 & \text{if } a = 0 \end{cases}$$

$\eta_3^1 : V_3^1 \to [0,1]$ by

$$\eta_3^1 \begin{pmatrix} a_1 & a_5 \\ a_2 & a_6 \\ a_3 & a_7 \\ a_4 & a_8 \end{pmatrix} =$$

$$\begin{cases} \dfrac{1}{a_1 a_5 + a_2 a_6 + a_3 a_7 + a_4 a_8} & \text{if } a_1 a_5 + a_2 a_6 + a_3 a_7 + a_4 a_8 \neq 0 \\ 1 & \text{if } a_1 a_5 + a_2 a_6 + a_3 a_7 + a_4 a_8 = 0 \end{cases}$$

$\eta_1^2 : V_1^2 \to [0, 1]$ is given by

$$\eta_1^2 (A) = \begin{cases} \dfrac{1}{|A|} & \text{if } |A| \neq 0 \\ 1 & \text{if } |A| = 0 \end{cases}$$

where A is the 3×3 upper triangular matrix with entries from $Z^o$.

$\eta_2^2 : V_2^2 \to [0, 1]$ is defined by

$$\eta_2^2 (a\ b\ c) = \begin{cases} \dfrac{1}{abc} & \text{if } abc \neq 0 \\ 1 & \text{if } abc = 0 \end{cases}.$$

$\eta_3^2 : V_3^2 \to [0, 1]$ is such that

$$\eta_3^2 [a\ a\ a\ a\ a\ a] = \begin{cases} \dfrac{1}{6a} & \text{if } a \neq 0 \\ 1 & \text{if } a = 0 \end{cases}$$

$\eta_4^2 : V_4^2 \to [0,1]$ is given by

$$\eta_4^2 \begin{pmatrix} a_1 \\ a_2 \\ a_3 \\ a_4 \\ a_5 \\ a_6 \end{pmatrix} = \begin{cases} \dfrac{1}{a_1 a_2 a_3 + a_4 a_5 a_6} & \text{if atleast one of } a_1 a_2 a_3 \neq 0 \\ 1 & \text{if both } a_1 a_2 a_3 = 0,\ a_4 a_5 a_6 = 0 \end{cases}$$

Thus the bifuzzy set
$$\eta = \eta_1 \cup \eta_2 = \left(\eta_1^1, \eta_2^1, \eta_3^1\right) \cup \left(\eta_1^2, \eta_2^2, \eta_3^2, \eta_4^2\right) :$$
$$(V_1 \cup V_2) = \left(V_1^1, V_2^1, V_3^1\right) \cup \left(V_1^2, V_2^2, V_3^2, V_4^2\right) \to [0, 1]$$

is a special set fuzzy vector bispace and is denoted by

$$\begin{aligned} V &= (V_1 \cup V_2)_{(\eta_1 \cup \eta_2)} \\ &= (V_1 \eta_1 \cup V_2 \eta_2) \\ &= \left(V_{1\eta_1^1}^1, V_{2\eta_2^1}^1, V_{3\eta_3^1}^1\right) \cup \left(V_{1\eta_1^2}^2, V_{2\eta_2^2}^2, V_{3\eta_3^2}^2, V_{4\eta_4^2}^2\right). \end{aligned}$$

Now having defined the notion of special set fuzzy vector bispace now we proceed on to define the notion of special set fuzzy vector subbispace.

**DEFINITION 2.4.6:** *Let*
$$V = (V_1 \cup V_2) = \left(V_1^1, V_2^1, ..., V_{n_1}^1\right) \cup \left(V_1^2, V_2^2, ..., V_{n_2}^2\right)$$
*be a special set vector bispace over a set S. Take*
$$W = (W_1 \cup W_2) = \left(W_1^1, W_2^1, ..., W_{n_1}^1\right) \cup \left(W_1^2, W_2^2, ..., W_{n_2}^2\right)$$
$\subseteq (V_1 \cup V_2) = V$ *to be a special set vector subbispace of V over the set S. Let* $\eta = \eta_1 \cup \eta_2$ *be a bimap from* $W = (W_1 \cup W_2)$ *into the set [0, 1]; i.e.,* $\eta_1 \cup \eta_2 = \left(\eta_1^1, \eta_2^1, ..., \eta_{n_1}^1\right) \cup \left(\eta_1^2, \eta_2^2, ..., \eta_{n_2}^2\right)$:
$W = (W_1 \cup W_2) = \left(W_1^1, W_2^1, ..., W_{n_1}^1\right) \cup \left(W_1^2, W_2^2, ..., W_{n_2}^2\right) \to [0, 1]$
*such that* $\eta_i^1 : W_i^1 \to [0,1]$ *for each* $i = 1, 2, ..., n$, $\eta_i^2 : W_i^2 \to [0, 1]$ *for each* $i = 1, 2, ..., n$ *such that* $W_{i\eta_i^1}^1$ *and* $W_{i\eta_i^2}^2$ *are fuzzy set vector bispace of* $V_i^1$ *and* $V_i^2$ *respectively. Then we call*

$$\begin{aligned}
W\eta &= (W_1 \cup W_2)_{\eta_1 \eta_2} \\
&= \left(W_1^1, W_2^1, ..., W_{n_1}^1\right)_{\eta_1} \cup \left(W_1^2, W_2^2, ..., W_{n_2}^2\right)_{\eta_2} \\
&= \left(W_1^1, W_2^1, ..., W_{n_1}^1\right)_{(\eta_1^1, \eta_2^1 ... \eta_{n_1}^1)} \cup \left(W_1^2, W_2^2, ..., W_{n_2}^2\right)_{(\eta_1^2, \eta_2^2 ... \eta_{n_2}^2)} \\
&= \left(W_{1\eta_1^1}^1, W_{2\eta_2^1}^1, ..., W_{n_1\eta_{n_1}^1}^1\right) \cup \left(W_{1\eta_1^2}^2, W_{2\eta_2^2}^2, ..., W_{n_2\eta_{n_2}^2}^2\right)
\end{aligned}$$

*the special set fuzzy vector subbispace or bisubspace.*

Now having defined a special set fuzzy bisubspace we proceed onto illustrate it by some example.

*Example 2.4.5:* Let
$$V = (V_1 \cup V_2) = \left(V_1^1, V_2^1, V_3^1\right) \cup \left(V_1^2, V_2^2, V_3^2, V_4^2\right)$$
be a special set vector bispace over the set $S = Z^o = Z^+ \cup \{0\}$ where $V_1^1 = S \times S \times S$,

$$V_2^1 = \left\{ \begin{pmatrix} a & b \\ c & d \end{pmatrix} \middle| a,b,c,d \in Z^o \right\},$$

$$V_3^1 = \left\{ \begin{pmatrix} a_1 & a_2 & a_3 \\ a_4 & a_5 & a_6 \end{pmatrix} \middle| a_i \in Z^o \right\}.$$

$V_2 = \left( V_1^2, V_2^2, V_3^2, V_4^2 \right)$ where

$V_1^2 = \{3 \times 3$ upper triangular matrices with entries from $Z^o\}$,
$$V_2^2 = \{S \times S \times S \times S \times S\},$$

$$V_3^2 = \left\{ \begin{pmatrix} a_1 & a_5 \\ a_2 & a_6 \\ a_3 & a_7 \\ a_4 & a_8 \end{pmatrix} \middle| a_i \in Z^o; 1 \le i \le 8 \right\}$$

and

$V_4^2 = \{$all $4 \times 4$ lower triangular matrices with entries from $Z^o\}$.

Now consider the set vector subspaces of these spaces; take

$$W_1^1 = \{S \times S \times \{0\}\} \subseteq V_1^1,$$

$$W_2^1 = \left\{ \begin{pmatrix} a & b \\ c & d \end{pmatrix} \middle| a,b,c,d \in 5Z^o \right\} \subseteq V_2^1,$$

$$W_3^1 = \left\{ \begin{pmatrix} a & a & a \\ a & a & a \end{pmatrix} \middle| a \in Z^o \right\} \subseteq V_3^1,$$

$W_1^2 = \{$all $3 \times 3$ upper triangular matrices with entries from $Z^o\}$
$\subseteq V_1^2$, $W_2^2 = \{S \times S \times S \times \{0\} \times S\} \subseteq V_2^2$,

$$W_3^2 = \left\{ \begin{pmatrix} a & a \\ a & a \\ a & a \\ a & a \end{pmatrix} \middle| a_i \in Z^o \right\} \subseteq V_3^2$$

and $W_4^2 = \{$all $4 \times 4$ lower triangular matrices with entries from $19Z^o\} \subseteq V_4^2$. Clearly

$$W = (W_1 \cup W_2)$$
$$= \left(W_1^1, W_2^1, W_3^1\right) \cup \left(W_1^2, W_2^2, W_3^2, W_4^2\right) \subseteq (V_1 \cup V_2)$$

is a special set vector subbispace of V over the set $S = Z^o$.

Define a fuzzy bimap
$$\eta = \eta_1 \cup \eta_2$$
$$= \left(\eta_1^1, \eta_2^1, \eta_3^1\right) \cup \left(\eta_1^2, \eta_2^2, \eta_3^2, \eta_4^2\right):$$
$$(W_1 \cup W_2) = \left(W_1^1, W_2^1, W_3^1\right) \cup \left(W_1^2, W_2^2, W_3^2\right) \to [0, 1] \text{ as}$$

$\eta_1^1 : W_1^1 \to [0, 1]$ is defined by

$$\eta_1^1 (a\ b\ 0) = \begin{cases} \dfrac{1}{a+b} & \text{if } a+b \neq 0 \\ 1 & \text{if } a+b = 0 \end{cases}$$

$\eta_2^1 : W_2^1 \to [0,1]$ is given by

$$\eta_2^1 \begin{pmatrix} a & b \\ c & d \end{pmatrix} = \begin{cases} \dfrac{1}{a+b+c+d} & \text{if } a+b+c+d \neq 0 \\ 1 & \text{if } a+b+c+d = 0 \end{cases}$$

$\eta_3^1 : W_3^1 \to [0,1]$ is given by

$$\eta_3^1 \begin{bmatrix} a & a & a \\ a & a & a \end{bmatrix} = \begin{cases} \dfrac{1}{6a} & \text{if } a \neq 0 \\ 1 & \text{if } a = 0 \end{cases}$$

$\eta_1^2 : W_1^2 \to [0, 1]$ is given by

$$\eta_1^2 \begin{bmatrix} a & a & a \\ 0 & a & a \\ 0 & 0 & a \end{bmatrix} = \begin{cases} \dfrac{1}{6a} & \text{if } a \neq 0 \\ 1 & \text{if } a = 0 \end{cases}$$

$\eta_2^2 : W_2^2 \to [0, 1]$ is defined by

$$\eta_2^2 \,(a\ b\ c\ 0\ d) = \begin{cases} \dfrac{1}{abcd} & \text{if } abcd \neq 0 \\ 1 & \text{if } abcd = 0 \end{cases}$$

$\eta_3^2 : W_3^2 \to [0,1]$ is such that

$$\eta_3^2 \begin{bmatrix} a & a \\ a & a \\ a & a \\ a & a \end{bmatrix} = \begin{cases} \dfrac{1}{8a} & \text{if } a \neq 0 \\ 1 & \text{if } a = 0 \end{cases}$$

$\eta_4^2 : W_4^2 \to [0,1]$ is such that

$$\eta_4^2 \begin{bmatrix} a_1 & 0 & 0 & 0 \\ a_2 & a_3 & 0 & 0 \\ a_4 & a_5 & a_6 & 0 \\ a_7 & a_8 & a_9 & a_{10} \end{bmatrix} = \begin{cases} \dfrac{1}{a_1 a_3 a_6 a_{10}} & \text{if } a_1 a_3 a_6 a_{10} \neq 0 \\ 1 & \text{if } a_1 a_3 a_6 a_{10} = 0 \end{cases}$$

Thus
$W\eta = (W_1 \cup W_2)_{\eta_1 \cup \eta_2}$
$\quad = \left(W_1^1, W_2^1, W_3^1\right)_{(\eta_1^1, \eta_2^1, \eta_3^1)} \cup \left(W_1^2, W_2^2, W_3^2, W_4^2\right)_{(\eta_1^2, \eta_2^2, \eta_3^2, \eta_4^2)}$
$\quad = \left(W_{1\eta_1^1}^1, W_{2\eta_2^1}^1, W_{3\eta_3^1}^1\right) \cup \left(W_{1\eta_1^2}^2, W_{2\eta_2^2}^2, W_{3\eta_3^2}^2, W_{4\eta_4^2}^2\right)$
is the special set fuzzy vector subbispace.

Now we proceed onto define the notion of special set fuzzy linear bialgebra.

**DEFINITION 2.4.7:** *Let*
$$V = (V_1 \cup V_2) = \left(V_1^1, V_2^1, ..., V_{n_1}^1\right) \cup \left(V_1^2, V_2^2, ..., V_{n_2}^2\right)$$
*be a special set linear bialgebra over the set S. Let $\eta = \eta_1 \cup \eta_2$ = be a bimap from $(V_1 \cup V_2)$ into [0,1], i.e., $\eta = \eta_1 \cup \eta_2 =$ $\left(\eta_1^1, \eta_2^1, ..., \eta_{n_1}^1\right) \cup \left(\eta_1^2, \eta_2^2, ..., \eta_{n_2}^2\right)$: $V = (V_1 \cup V_2) = \left(V_1^1, V_2^1, ..., V_{n_1}^1\right)$ $\cup \left(V_1^2, V_2^2, ..., V_{n_2}^2\right) \to [0,1]$. Such that $\eta_i^1 : V_i^1 \to [0,1]$ such that $V_{i\eta_i^1}^1$ is a set fuzzy linear algebra for $i = 1, 2, ..., n$. and $\eta_i^2 : V_i^2 \to [0,1]$ such that $V_{2\eta_i^2}$ is a set fuzzy linear algebra for $i = 1, 2, ..., n_2$; then we call*
$$\begin{aligned} V\eta &= (V_1 \cup V_2)_{\eta_1 \eta_2} \\ &= V_1 \eta_1 \cup V_2 \eta_2 \\ &= \left(V_1^1, V_2^1, ..., V_{n_1}^1\right)_{(\eta_1^1, \eta_2^1 ... \eta_{n_1}^1)} \cup \left(V_1^2, V_2^2, ..., V_{n_2}^2\right)_{(\eta_1^2, \eta_2^2 ... \eta_{n_2}^2)} \\ &= \left(V_{1\eta_1^1}^1, V_{2\eta_2^1}^1, ..., V_{n_1 \eta_{n_1}^1}^1\right) \cup \left(V_{1\eta_1^2}^2, V_{2\eta_2^2}^2, ..., V_{n_2 \eta_{n_2}^2}^2\right) \end{aligned}$$
*the special set fuzzy linear bialgebra.*

We illustrate this by some simple examples.

*Example 2.4.6:* Let
$$V = (V_1 \cup V_2) = \left(V_1^1, V_2^1, V_3^1, V_4^1\right) \cup \left(V_1^2, V_2^2, V_3^2, V_4^2\right)$$
be a special set vector bispace over the set $S = Z^o = Z^+ \cup \{0\}$ where

$$V_1^1 = \left\{ \begin{pmatrix} a_1 & a_2 & a_3 \\ a_4 & a_5 & a_6 \end{pmatrix} \middle| a_i \in Z^o = S = \{Z^+ \cup \{0\}\} \right\},$$

$V_2^1 = S \times S \times S \times S$, $V_3^1 = \{(a\ a\ a\ a\ a\ a) \mid a \in S\}$, $V_4^1 = \{$all $5 \times 5$ upper triangular matrices with entries from $S\}$ and

$$V_1^2 = \left\{ \begin{pmatrix} a & b \\ c & d \end{pmatrix} \middle| a, b, c, d \in Z^o \right\},$$

$$V_2^2 = \left\{ \begin{pmatrix} a \\ a \\ a \\ a \\ a \\ a \end{pmatrix} \middle| a \in Z^o \right\},$$

$V_3^2 = S \times S \times S \times S \times S \times S$ and $V_4^2 = \{4 \times 4$ lower triangular matrices with entries from $S\}$. Clearly $V = (V_1 \cup V_2)$ is a special set linear bialgebra over the set S.

Let us define the fuzzy bimap

$$\begin{aligned} \eta &= \eta_1 \cup \eta_2 \\ &= \left(\eta_1^1, \eta_2^1, \eta_3^1, \eta_4^1\right) \cup \left(\eta_1^2, \eta_2^2, \eta_3^2, \eta_4^2\right): \\ V &= (V_1 \cup V_2) \\ &= \left(V_1^1, V_2^1, V_3^1, V_4^1\right) \cup \left(V_1^2, V_2^2, V_3^2, V_4^2\right) \to [0, 1] \end{aligned}$$

by
$\eta_1^1 : V_1^1 \to [0, 1]$ such that

$$\eta_1^1 \begin{pmatrix} a_1 & a_2 & a_3 \\ a_4 & a_5 & a_6 \end{pmatrix} = \begin{cases} \dfrac{1}{a_1 a_2 a_3} & \text{if } a_1 a_2 a_3 \neq 0 \\ 1 & \text{if } a_1 a_2 a_3 = 0 \end{cases}$$

$\eta_2^1 : V_2^1 \to [0,1]$ is defined by

$$\eta_2^1 [a\ b\ c\ d] = \begin{cases} \dfrac{1}{abcd} & \text{if } abcd \neq 0 \\ 1 & \text{if } abcd = 0 \end{cases}$$

$\eta_3^1 : V_3^1 \to [0, 1]$ is defined by

$$\eta_3^1 (a\ a\ a\ a\ a\ a) = \begin{cases} \dfrac{1}{6a} & \text{if } a \neq 0 \\ 1 & \text{if } a = 0 \end{cases}$$

$\eta_4^1 : V_4^1 \to [0, 1]$ is defined by

$$\eta_4^1 \begin{pmatrix} a_1 & a_2 & a_3 & a_4 & a_5 \\ 0 & a_6 & a_7 & a_8 & a_9 \\ 0 & 0 & a_{10} & a_{11} & a_{12} \\ 0 & 0 & 0 & a_{13} & a_{14} \\ 0 & 0 & 0 & 0 & a_{15} \end{pmatrix} = \begin{cases} \dfrac{1}{a_1 a_6 a_{10} a_{13} a_{15}} & \text{if } a_1 a_6 a_{10} a_{13} a_{15} \neq 0 \\ 1 & \text{if } a_1 a_6 a_{10} a_{13} a_{15} = 0 \end{cases}$$

$\eta_1^2 : V_1^2 \to [0,1]$ is defined by

$$\eta_1^2 \begin{pmatrix} a & b \\ c & d \end{pmatrix} = \begin{cases} \dfrac{1}{bc} & \text{if } bc \neq 0 \\ 1 & \text{if } bc = 0 \end{cases}$$

$\eta_2^2 : V_2^2 \to [0, 1]$ is such that

$$\eta_2^2 \begin{pmatrix} a \\ a \\ a \\ a \\ a \end{pmatrix} = \begin{cases} \dfrac{1}{5a} & \text{if } a \neq 0 \\ 1 & \text{if } a = 0 \end{cases}$$

$\eta_3^2 : V_3^2 \to [0, 1]$ is defined by

$$\eta_3^2 [a\ b\ c\ d\ e\ f] = \begin{cases} \dfrac{1}{ace} & \text{if } ace \neq 0 \\ 1 & \text{if } ace = 0 \end{cases}$$

$\eta_4^2 : V_4^2 \to [0, 1]$ is given by

$$\eta_4^2 \begin{pmatrix} a & 0 & 0 & 0 \\ b & c & 0 & 0 \\ d & e & f & 0 \\ g & h & i & j \end{pmatrix} = \begin{cases} \dfrac{1}{acbj} & \text{if } acbj \neq 0 \\ 1 & \text{if } acbj = 0 \end{cases}$$

Thus $\eta = \eta_1 \cup \eta_2 : (V_1 \cup V_2) \to [0, 1]$ is a bimap such that

$$V\eta = (V_1 \cup V_2)_{\eta_1 \cup \eta_2}$$
$$= \left(V_1^1, V_2^1, ..., V_{n_1}^1\right)_{(\eta_1^1, \eta_2^1 ... \eta_{n_1}^1)} \cup \left(V_1^2, V_2^2, ..., V_{n_2}^2\right)_{(\eta_1^2, \eta_2^2 ... \eta_{n_2}^2)}$$
$$= \left(V_{1\eta_1^1}^1, V_{2\eta_2^1}^1, ..., V_{n_1\eta_{n_1}^1}^1\right) \cup \left(V_{1\eta_1^2}^2, V_{2\eta_2^2}^2, ..., V_{n_2\eta_{n_2}^2}^2\right)$$

is a special set fuzzy linear bialgebra, here $n_1 = n_2 = 4$.

Now we proceed on to define the notion of special set fuzzy linear subbialgebra.

**DEFINITION 2.4.8:** *Let*
$$V = (V_1 \cup V_2) = \left(V_1^1, V_2^1, ..., V_{n_1}^1\right) \cup \left(V_1^2, V_2^2, ..., V_{n_2}^2\right)$$
*be a special set linear bialgebra over the set S. Take $W = (W_1 \cup W_2) = \left(W_1^1, W_2^1, ..., W_{n_1}^1\right) \cup \left(W_1^2, W_2^2, ..., W_{n_2}^2\right) \subseteq (V_1 \cup V_2)$ be a special set linear subbialgebra of V over the set S. If*
$$\eta = \eta_1 \cup \eta_2$$
$$= \left(\eta_1^1, \eta_2^1, ..., \eta_{n_1}^1\right) \cup \left(\eta_1^2, \eta_2^2, ..., \eta_{n_2}^2\right):$$
$$(W_1 \cup W_2) = \left(W_1^1, W_2^1, ..., W_{n_1}^1\right) \cup \left(W_1^2, W_2^2, ..., W_{n_2}^2\right) \to [0,1]$$
*be such that $\eta_i^1 : W_i^1 \to [0,1]$ such that $W_{i\eta_i^1}^1$ is a set fuzzy linear algebra for every i = 1, 2, ..., n, and $\eta_j^2 : W_j^2 \to [0,1]$ is such that $W_{j\eta_j^2}^2$ is a set fuzzy linear algebra for each j = 1, 2, ..., $n_2$ then we call*

$$W\eta = (W_1 \cup W_2)_{\eta_1 \cup \eta_2}$$
$$= \left(W^1_{1\eta^1_1}, W^1_{2\eta^1_2}, \ldots, W^1_{n_1\eta^1_{n_1}}\right) \cup \left(W^2_{1\eta^2_1}, W^2_{2\eta^2_2}, \ldots, W^2_{n_2\eta^2_{n_2}}\right)$$

*is a special set fuzzy linear subbialgebra.*

We illustrate this by a simple example.

***Example 2.4.7:*** Let $V = \left(V^1_1, V^1_2, V^1_3\right) \cup \left(V^2_1, V^2_2, V^2_3, V^2_4, V^2_5\right)$ be a special set linear bialgebra of V over the set $Z^o = S = Z^+ \cup \{0\}$. Here $V^1_1 = S \times S \times S$,

$$V^1_2 = \left\{ \begin{pmatrix} a & b \\ 0 & c \end{pmatrix} \middle| a,b,c, \in S \right\}$$

and

$$V^1_3 = \{(a\ a\ a\ a\ a\ a\ a) | a \in Z^o\}$$

$$V^2_1 = \left\{ \begin{pmatrix} a & b \\ c & d \end{pmatrix} \middle| a,b,c,d \in Z^o \right\},$$

$$V^2_2 = \left\{ \begin{pmatrix} a & a \\ a & a \\ a & a \\ a & a \end{pmatrix} \middle| a \in Z^o \right\},$$

$V^2_3 = \{\text{all } 4\times4 \text{ lower triangular matrices}\}$,

$$V^2_4 = \left\{ \begin{pmatrix} a & a & a & a & a & a \\ a & a & a & a & a & a \end{pmatrix} \middle| a \in S \right\}$$

and $V^2_5 = \{\text{all } 5 \times 5 \text{ upper triangular matrices}\}$. Consider
$$W = (W_1 \cup W_2) = \left(W^1_1, W^1_2, W^1_3\right) \cup \left(W^2_1, W^2_2, W^2_3, W^2_4, W^2_5\right)$$
$$\subseteq (V_1 \cup V_2)$$

be a subbialgebra of V over the set S, where

$$W_1^1 = \{S \times S \times \{0\}\} \subseteq V_1^1,$$

$$W_2^1 = \left\{ \begin{pmatrix} a & b \\ 0 & c \end{pmatrix} \middle| a,b,c \in 7Z^o \right\} \subseteq V_2^1$$

and

$$W_3^1 = \{(a\ a\ a\ a\ a\ a\ a) \mid a \in 5Z^o\} \subseteq V_3^1,$$

$$W_1^2 = \left\{ \begin{pmatrix} a & b \\ c & d \end{pmatrix} \middle| a,b,c,d \in Z^o \right\} \subseteq V_1^2,$$

$$W_2^2 = \left\{ \begin{pmatrix} a & a \\ a & a \\ a & a \\ a & a \end{pmatrix} \middle| a \in 10Z^o \right\} \subseteq V_2^2$$

$$W_3^2 = \{\text{all } 4 \times 4 \text{ lower triangular matrices}\} \subseteq V_3^2,$$

$$W_4^2 = \left\{ \begin{pmatrix} a & a & a & a & a & a \\ a & a & a & a & a & a \end{pmatrix} \middle| a \in 2Z^o \right\} \subseteq V_4^2$$

and

$$W_5^2 = \{\text{all } 5 \times 5 \text{ upper triangular matrices}\} \subseteq V_5^2.$$

Define $\eta = \eta_1 \cup \eta_2: (W_1 \cup W_2) \to [0, 1]$ as follows.

$\eta_1^1: W_1^1 \to [0,1]$ is such that

$$\eta_1^1 (a\ b\ 0) = \begin{cases} \dfrac{1}{ab} & \text{if } ab \neq 0 \\ 1 & \text{if } ab = 0 \end{cases}$$

$\eta_2^1: W_2^1 \to [0,1]$ is defined by

$$\eta_2^1 \begin{pmatrix} a & b \\ 0 & c \end{pmatrix} = \begin{cases} \dfrac{1}{abc} & \text{if } abc \neq 0 \\ 1 & \text{if } abc = 0 \end{cases}$$

$\eta_3^1 : W_3^1 \to [0, 1]$ is given by

$$\eta_3^1 (a\ a\ a\ a\ a\ a\ a) = \begin{cases} \dfrac{1}{7a} & \text{if } a \neq 0 \\ 1 & \text{if } a = 0 \end{cases}$$

$\eta_1^2 : W_1^2 \to [0, 1]$ is given by

$$\eta_1^2 \begin{pmatrix} a & b \\ c & d \end{pmatrix} = \begin{cases} \dfrac{1}{\begin{vmatrix} a & b \\ c & d \end{vmatrix}} & \text{if } ab - bc \neq 0 \\ 1 & \text{if } ab - bc = 0 \end{cases}$$

$\eta_2^2 : W_2^2 \to [0,1]$ is given by

$$\eta_2^2 \begin{pmatrix} a & a \\ a & a \\ a & a \\ a & a \end{pmatrix} = \begin{cases} \dfrac{1}{8a} & \text{if } a \neq 0 \\ 1 & \text{if } a = 0 \end{cases}$$

$\eta_3^2 : W_3^2 \to [0,1]$ is such that

$$\eta_3^2 \begin{bmatrix} a & 0 & 0 & 0 \\ b & c & 0 & 0 \\ d & e & f & 0 \\ g & h & i & j \end{bmatrix} = \begin{cases} \dfrac{1}{|A|} & \text{if } |A| \neq 0 \\ 1 & \text{if } |A| = 0 \end{cases}$$

where

$$A = \begin{bmatrix} a & 0 & 0 & 0 \\ b & c & 0 & 0 \\ d & e & f & 0 \\ g & h & i & j \end{bmatrix}.$$

$\eta_4^2 : W_4^2 \to [0, 1]$ is given by

$$\eta_4^2 \begin{bmatrix} a & a & a & a & a & a \\ a & a & a & a & a & a \end{bmatrix} = \begin{cases} \dfrac{1}{12a} & \text{if } a \neq 0 \\ 1 & \text{if } a = 0 \end{cases}$$

$\eta_5^2 : W_5^2 \to [0, 1]$ is such that

$$\eta_5^2 \begin{bmatrix} a & b & c & d & e \\ 0 & f & g & h & i \\ 0 & 0 & j & k & l \\ 0 & 0 & 0 & m & n \\ 0 & 0 & 0 & 0 & p \end{bmatrix} = \begin{cases} \dfrac{1}{\text{afjmp}} & \text{if afjmp} \neq 0 \\ 1 & \text{if afjmp} = 0 \end{cases}$$

Thus
$$W\eta = (W_1 \cup W_2)_{\eta_1 \cup \eta_2}$$
$$= \left( W_{1\eta_1^1}^1, W_{2\eta_1^2}^1, ..., W_{n_1\eta_1^{n_1}}^1 \right) \cup \left( W_{1\eta_1^2}^2, W_{2\eta_2^2}^2, ..., W_{n_2\eta_{n_2}^2}^2 \right)$$

is the special set fuzzy linear subbialgebra.

Now we are bound to make the following observations:

1. The notion of special set fuzzy vector bispace and special set fuzzy linear bialgebra are fuzzy equivalent.
2. Likewise the notion of special set fuzzy vector subbispace and special set fuzzy linear subbialgebra are also fuzzy equivalent.
3. Now using the special set linear bialgebra and special set vector bispace we may define a infinite number of special set fuzzy linear bialgebra and special set fuzzy vector bispace.

The same holds good for special set fuzzy linear subbialgebra and special set fuzzy vector subbispaces.

We now proceed onto define the notion of special set fuzzy n vector spaces and special set fuzzy linear n- algebras; though we know both the concepts are fuzzy equivalent.

**DEFINITION 2.4.9:** *Let*
$$V = (V_1 \cup V_2 \cup \ldots \cup V_n)$$
$$= \left(V_1^1, V_2^1, \ldots, V_{n_1}^1\right) \cup \left(V_1^2, V_2^2, \ldots, V_{n_2}^2\right) \cup \ldots \cup \left(V_1^n, V_2^n, \ldots, V_{n_n}^n\right)$$
*be a special set vector n – space over the set S. Let*
$$\eta = \eta_1 \cup \eta_2 \cup \ldots \cup \eta_n$$
$$= \left(\eta_1^1, \eta_2^1, \ldots, \eta_{n_1}^1\right) \cup \left(\eta_1^2, \eta_2^2, \ldots, \eta_{n_2}^2\right) \cup \ldots \cup \left(\eta_1^n, \eta_2^n, \ldots, \eta_{n_n}^n\right):$$
$$V = (V_1 \cup V_2 \cup \ldots \cup V_n) \to [0, 1]$$
*such that for each $i_1$, $\eta_{i_1}^1 : V_{i_1}^1 \to [0, 1]$ so that $V_{i_1 \eta_{i_1}^1}^1$ is a set fuzzy vector space and this is true for $1 \le i_1 \le n_1$; $\eta_{i_2}^2 : V_{i_2}^2 \to [0,1]$ such that $V_{i_2 \eta_{i_2}^2}^2$ is a set fuzzy vector space for each $i_2$; $1 \le i_2 \le n_2$; and so on $\eta_{i_n}^n : V_{i_n}^n \to [0,1]$ such that $V_{i_n \eta_{i_n}^n}^n$ is a set fuzzy vector space for each $i_n$; $1 \le i_n \le n_n$.*
*Thus*
$$V\eta = (V_1 \cup V_2 \cup \ldots \cup V_n)_{\eta_1 \cup \ldots \cup \eta_n}$$
$$= V_1 \eta_1 \cup \ldots \cup V_n \eta_n$$
$$= \left(V_{1\eta_1^1}^1, V_{2\eta_2^1}^1, \ldots, V_{n_1 \eta_{n_1}^1}^1\right) \cup \left(V_{1\eta_1^2}^2, V_{2\eta_2^2}^2, \ldots, V_{n_2 \eta_{n_2}^2}^2\right) \cup \ldots \cup \left(V_{n\eta_1^n}^n, \ldots, V_{n_n \eta_{n_n}^n}^n\right)$$

*is a special set fuzzy n-vector space.*

We illustrate this by the following example:

***Example 2.4.8:*** Let
$$V = V_1 \cup V_2 \cup V_3 \cup V_4$$
$$= \left(V_1^1, V_2^1, V_3^1\right) \cup \left(V_1^2, V_2^2\right) \cup \left(V_1^3, V_2^3, V_3^3, V_4^3\right) \cup \left(V_1^4, V_2^4, V_3^4\right)$$
be a special set 4 vector space over the set $S = Z^+ \cup \{0\}$ where

$$V_1^1 = S \times S \times S \times S,$$

$$V_2^1 = \left\{ \begin{pmatrix} a & b \\ c & d \end{pmatrix} \middle| a,b,c,d \in S \right\},$$

$$V_3^1 = \left\{ \begin{bmatrix} a & a & a & a \\ a & a & a & a \end{bmatrix} \begin{bmatrix} a \\ a \\ a \\ a \\ a \end{bmatrix} \middle| a \in S \right\},$$

$V_1^2 = \{3 \times 3$ matrices with entries from $S\}$, $V_2^2 = S \times S \times S$, $V_1^3 = S \times S \times S \times S \times S \times S$, $V_2^3 = \{$set of all $4 \times 4$ matrices with entries from the set $S\}$,

$$V_3^3 = \left\{ \begin{bmatrix} a \\ a \\ a \\ a \\ a \end{bmatrix} [a \quad a \quad a] \middle| a \in S \right\}$$

$$V_4^3 = \left\{ \begin{bmatrix} a & a \\ a & a \end{bmatrix}, \begin{bmatrix} a & a & a & a \\ a & a & a & a \end{bmatrix} \middle| a \in S \right\},$$

$V_1^4 = \{$set of $3 \times 3$ upper triangular matrices with entries from $S\}$, $V_2^4 = \{S \times S \times S \times S\}$, $V_3^4 = \{(a\,a\,a\,a\,a\,a\,a) \mid a \in S\}$.

Define $\eta = \eta_1 \cup \eta_2 \cup ... \cup \eta_4 : V = V_1 \cup V_2 \cup V_3 \cup V_4 \to [0, 1]$ such that $\eta_i: V_i \to [0, 1]$ for $i = 1, 2, 3, 4$.
Now
$$\eta_1 = \left(\eta_1^1, \eta_2^1, \eta_3^1\right) : \left(V_1^1, V_2^1, V_3^1\right) \to [0,1]$$
is such that
$\eta_1^1 : V_1^1 \to [0, 1]$ is defined as

$$\eta_1^1(a\ b\ c\ d) = \begin{cases} \dfrac{1}{ad} & \text{if } ad \neq 0 \\ 1 & \text{if } ad = 0 \end{cases}$$

$\eta_2^1 : V_2^1 \to [0, 1]$ given by

$$\eta_2^1 \begin{pmatrix} a & b \\ c & d \end{pmatrix} = \begin{cases} \dfrac{1}{ad - bc} & \text{if } ad - bc \neq 0 \\ 1 & \text{if } ad - bc = 0 \end{cases}$$

$\eta_3^1 : V_3^1 \to [0,1]$ is defined by

$$\eta_3^1 \begin{pmatrix} a & a & a & a \\ a & a & a & a \end{pmatrix} = \begin{cases} \dfrac{1}{8a} & \text{if } a \neq 0 \\ 1 & \text{if } a = 0 \end{cases}$$

Thus $V_1 \eta_1 = \left( V_{1\eta_1^1}^1, V_{2\eta_2^1}^1, V_{3\eta_3^1}^1 \right)$ is a special set fuzzy vector space.

Now $\eta_2 : V_2 \to [0,1]$, i.e., $\eta_2 = \left( \eta_1^2, \eta_2^2 \right) : \left( V_1^2, V_2^2 \right) \to [0,1]$

where
$\eta_1^2 : V_1^2 \to [0,1]$ is such that

$$\eta_1^2 \begin{pmatrix} a & b & c \\ d & e & f \\ g & h & i \end{pmatrix} = \begin{cases} \dfrac{1}{aei} & \text{if } aei \neq 0 \\ 1 & \text{if } aei = 0 \end{cases}$$

$\eta_2^2 : V_2^2 \to [0, 1]$ is given by

$$\eta_2^2 (a\ b\ c) = \begin{cases} \dfrac{1}{abc} & \text{if } abc \neq 0 \\ 1 & \text{if } abc = 0 \end{cases}$$

Thus $V_2 \eta_2 = \left( V_{1\eta_1^2}^2, V_{2\eta_2^2}^2 \right)$ is again special set fuzzy vector space.

$\eta_3: V_3 \to [0,1]$ is given by
$$\eta_{3=}\left(\eta_1^3, \eta_2^3, \eta_3^3, \eta_4^3\right): V_3 = \left(V_1^3, V_2^3, V_3^3, V_4^3\right) \to [0,1]$$
is such that

$\eta_1^3 : V_1^3 \to [0,1]$ is given by

$\eta_1^3([a\ b\ c\ d\ e\ f]) =$

$$\begin{cases} \dfrac{1}{a+b+c+d+e+f} & \text{if } a+b+c+d+e+f \neq 0 \\ 1 & \text{if } a+b+c+d+e+f = 0 \end{cases}$$

$\eta_2^3 : V_2^3 \to [0,1]$ is given by

$$\eta_2^3 \begin{pmatrix} a_{11} & a_{12} & a_{13} & a_{14} \\ a_{21} & a_{22} & a_{23} & a_{24} \\ a_{31} & a_{32} & a_{33} & a_{34} \\ a_{41} & a_{42} & a_{43} & a_{44} \end{pmatrix} = \begin{cases} \dfrac{1}{a_{11}a_{22}a_{33}a_{44}} & \text{if } a_{11}a_{22}a_{33}a_{44} \neq 0 \\ 1 & \text{if } a_{11}a_{22}a_{33}a_{44} = 0 \end{cases}$$

$\eta_3^3 : V_3^3 \to [0,1]$ is defined by

$$\eta_3^3 \begin{pmatrix} a & a \\ a & a \end{pmatrix} = \begin{cases} \dfrac{1}{4a} & \text{if } a \neq 0 \\ 1 & \text{if } a = 0 \end{cases}$$

$\eta_4^3 : V_4^3 \to [0,1]$ is such that

$$\eta_4^3 \begin{bmatrix} a & a & a & a & a \\ a & a & a & a & a \end{bmatrix} = \begin{cases} \dfrac{1}{10a} & \text{if } a \neq 0 \\ 1 & \text{if } a = 0 \end{cases}$$

Thus $V_3\eta_3 = \left(V_{1\eta_1^3}^3, V_{2\eta_2^3}^3, V_{3\eta_3^3}^3, V_{4\eta_4^3}^3\right)$ is a special set fuzzy vector space.

Finally $\eta_4: V_4 \to [0, 1]$, i.e.,
$$\eta_4 = \left(\eta_1^4, \eta_2^4, \eta_3^4\right) : \left(V_1^4, V_2^4, V_3^4\right) \to [0,1]$$
is such that

$\eta_1^4 : V_1^4 \to [0,1]$ such that

$$\eta_1^4 \begin{pmatrix} a & b \\ d & c \end{pmatrix} = \begin{cases} \dfrac{1}{ac} & \text{if } ac \neq 0 \\ 1 & \text{if } ac = 0 \end{cases}$$

$\eta_2^4 : V_2^4 \to [0,1]$ is defined such that

$$\eta_2^4 \, (a \ b \ c \ d) = \begin{cases} \dfrac{1}{ab+cd} & \text{if } ab+cd \neq 0 \\ 1 & \text{if } ab+cd = 0 \end{cases}$$

$\eta_3^4 : V_3^4 \to [0,1]$ is given by

$$\eta_3^4 \, (a \ a \ a \ a \ a \ a \ a) = \begin{cases} \dfrac{1}{7a} & \text{if } a \neq 0 \\ 1 & \text{if } a = 0 \end{cases}$$

Thus $V_4 \eta_4 = \left(V^4_{1\eta_1^4}, V^4_{2\eta_2^4}, V^4_{3\eta_3^4}\right)$ is a special set fuzzy vector space.

Thus

$$\begin{aligned} V &= (V_1 \cup V_2 \cup V_3 \cup V_4)_\eta \\ &= (V_1 \eta_1 \cup V_2 \eta_2 \cup V_3 \eta_3 \cup V_4 \eta_4) \\ &= \left(V^1_{1\eta_1^1}, V^1_{2\eta_2^1}, V^1_{3\eta_3^1}\right) \cup \left(V^2_{1\eta_1^2}, V^2_{2\eta_2^2}\right) \cup \\ &\quad \left(V^3_{1\eta_1^3}, V^3_{2\eta_2^3}, V^3_{3\eta_3^3}, V^3_{4\eta_4^3}\right) \cup \left(V^4_{1\eta_1^4}, V^4_{2\eta_2^4}, V^4_{3\eta_3^4}\right) \end{aligned}$$

is a special set fuzzy vector 4-space.

Now we proceed onto define the notion of special set fuzzy vector n-subspace.

**DEFINITION 2.4.10:** *Let*
$$V = (V_1 \cup V_2 \cup ... \cup V_n)$$
$$= \left(V_1^1, V_2^1, ..., V_{n_1}^1\right) \cup \left(V_1^2, V_2^2, ..., V_{n_2}^2\right) \cup ... \cup \left(V_1^n, V_2^n, ..., V_{n_n}^n\right)$$
*be a special set vector n-space over the set S. Suppose*
$$W = (W_1, W_2, ..., W_n)$$
$$= \left(W_1^1, W_2^1, ..., W_{n_1}^1\right) \cup \left(W_1^2, W_2^2, ..., W_{n_2}^2\right) \cup ... \cup \left(W_1^n, W_2^n, ..., W_{n_n}^n\right)$$
$$\subseteq (V_1 \cup V_2 \cup ... \cup V_n)$$
*i.e.,* $W_{t_i}^i \subseteq V_{t_i}^i$, $1 \le t_i \le n_i$ *and* $i = 1, 2, ..., n$, *is a special set vector n-subspace of V over the set S.*

*Now define a n-map*
$$\eta = (\eta_1, \cup \eta_2, \cup ... \cup \eta_n)$$
$$= \left(\eta_1^1, \eta_2^1, ..., \eta_{n_1}^1\right) \cup \left(\eta_1^2, \eta_2^2, ..., \eta_{n_2}^2\right) \cup ... \cup \left(\eta_1^n, \eta_2^n, ..., \eta_{n_n}^n\right):$$
$$(W_1, W_2, ..., W_n) \to [0, 1]$$
*such that* $\eta_i: W_i \to [0, 1]$; $i = 1, 2, ..., n$ *and* $\eta_{t_i}^i : W_{t_i}^i \to [0,1]$ *is so defined that* $W_{t_i\eta_{t_i}^i}^i$ *is a set fuzzy vector subspace true for* $t_i = 1, 2, ..., n_i$, $1 \le i \le n$.

$$W\eta = (W_1 \cup W_2 \cup ... \cup W_n)_\eta$$
$$W\eta = (W_1 \cup W_2 \cup ... \cup W_n)_{\eta_1 \cup \eta_2 \cup ... \cup \eta_n}$$
$$= \left(W_1^1, W_2^1, ..., W_{n_1}^1\right)_{(\eta_1^1, \eta_2^1 ... \eta_{n_1}^1)} \cup \left(W_1^2, W_2^2, ..., W_{n_2}^2\right)_{(\eta_1^2, \eta_2^2 ... \eta_{n_2}^2)}$$
$$\cup ... \cup \left(W_1^n, W_2^n, ..., W_{n_n}^n\right)_{(\eta_1^n, \eta_2^n ... \eta_{n_n}^n)}$$
$$= \left(W_{1\eta_1^1}^1, W_{2\eta_2^1}^1, ..., W_{n_1\eta_{n_1}^1}^1\right) \cup \left(W_{1\eta_1^2}^2, W_{2\eta_2^2}^2, ..., W_{n_2\eta_{n_2}^2}^2\right)$$
$$\cup ... \cup \left(W_{1\eta_1^n}^n, W_{2\eta_2^n}^n, ..., W_{n_n\eta_{n_n}^n}^n\right)$$
*is a special set fuzzy vector n-subspace.*

We illustrate this by the following example:

**Example 2.4.9:** Let $V = (V_1 \cup V_2 \cup V_3 \cup V_4 \cup V_5)$ be a special set 5-vector space over the set $S = Z^o = Z^+ \cup \{0\}$. Here
$$V_1 = \left(V_1^1, V_2^1, V_3^1\right), V_2 = \left(V_1^2, V_2^2\right),$$

$V_3 = \left(V_1^3, V_2^3, V_3^3\right), V_4 = \left(V_1^4, V_2^4\right)$ and $V_5 = \left(V_1^5, V_2^5, V_3^5, V_4^5\right)$

where $\quad V_1^1 = S \times S \times S \times S,$

$$V_2^1 = \left\{ \begin{pmatrix} a & b \\ c & d \end{pmatrix} \middle| a, b, c, d \in S \right\},$$

$$V_3^1 = \left\{ \begin{pmatrix} a & a & a & a & a & a \\ b & b & b & b & b & b \end{pmatrix} \middle| a, b \in S \right\},$$

$V_1^2 = \{$all upper triangular $4 \times 4$ matrices of the form

$$\left\{ \begin{pmatrix} a & a & a & a \\ 0 & b & b & b \\ 0 & 0 & a & a \\ 0 & 0 & 0 & b \end{pmatrix} \middle| a, b \in S \right\},$$

$$V_2^2 = S \times S \times S,$$

$$V_1^3 = \left\{ \begin{pmatrix} a & a & a \\ b & b & b \\ c & c & c \end{pmatrix} \middle| a, b, c \in S \right\},$$

$$V_2^3 = S \times S \times S \times S \times S,$$

$$V_3^3 = \left\{ \begin{pmatrix} a & a & a & a & a \\ b & b & b & b & b \\ c & c & c & c & c \end{pmatrix} \middle| a, b, c \in S \right\},$$

$V_1^4 = S \times S \times S \times S \times S$, $V_2^4 = \{$all upper triangular $5 \times 5$ matrices with entries from $S\}$,

$V_1^5 = $ {lower triangular 3×3 matrices of the form

$$\left\{ \begin{pmatrix} a & 0 & 0 \\ a & a & 0 \\ a & a & a \end{pmatrix} \middle| a \in S \right\},$$

$$V_2^5 = \{S \times S \times S \times S \times S\},$$

$$V_3^5 = \left\{ \begin{pmatrix} a & b & c & d \\ a & b & c & d \\ a & b & c & d \end{pmatrix} \middle| a,b,c,d \in S \right\}$$

and

$$V_4^5 = \left\{ \begin{pmatrix} a & a & a & a \\ b & b & b & b \\ c & c & c & c \\ d & d & d & d \end{pmatrix} \middle| a,b,c,d \in S \right\}$$

where $V_{j_i}^i$ is a special set vector space; $1 \leq i \leq 5$; $1 \leq j_i \leq n_i$. Now we define the 5-map

$$\begin{aligned} \eta &= (\eta_1 \cup \eta_2 \cup \eta_3 \cup \eta_4 \cup \eta_5) \\ &= \left(\eta_1^1, \eta_2^1, \eta_3^1\right) \cup \left(\eta_1^2, \eta_2^2\right) \cup \left(\eta_1^3, \eta_2^3, \eta_3^3\right) \cup \left(\eta_1^4, \eta_2^4\right) \\ &\quad \cup \left(\eta_1^5, \eta_2^5, \eta_3^5, \eta_4^5\right): \end{aligned}$$

$V = (V_1 \cup V_2 \cup V_3 \cup V_4 \cup V_5) \to [0,1]$

where $\eta_{j_i}^i : V_{j_i}^i \to [0, 1]$;

$i = 1, 2, 3, 4, 5$; $1 \leq j_i \leq n_i$; $i = 1, 2, 3, 4, 5$.

$\eta_1^1 : V_1^1 \to [0, 1]$ is defined by

$$\eta_1^1(a\ b\ c\ d) = \begin{cases} \dfrac{1}{abcd} & \text{if } abcd \neq 0 \\ 1 & \text{if } abcd = 0 \end{cases}$$

$\eta_2^1 : V_2^1 \to [0, 1]$ is such that

$$\eta_2^1 \begin{pmatrix} a & b \\ c & d \end{pmatrix} = \begin{cases} \dfrac{1}{ad-bc} & \text{if } ad \neq bc \\ 1 & \text{if } ad = bc \end{cases}$$

$\eta_3^1 : V_3^1 \to [0, 1]$ is given by

$$\eta_3^1 \begin{pmatrix} a & a & a & a & a & a \\ b & b & b & b & b & b \end{pmatrix} = \begin{cases} \dfrac{1}{6ab} & \text{if } ab \neq 0 \\ 1 & \text{if } ab = 0 \end{cases}$$

$\eta_1^2 : V_1^2 \to [0, 1]$ is given by

$$\eta_1^2 \begin{pmatrix} a & a & a & a \\ 0 & b & b & b \\ 0 & 0 & a & a \\ 0 & 0 & a & b \end{pmatrix} = \begin{cases} \dfrac{1}{ab} & \text{if } ab \neq 0 \\ 1 & \text{if } ab = 0 \end{cases}$$

$\eta_2^2 : V_2^2 \to [0, 1]$ is such that

$$\eta_2^2 (a\ b\ c) = \begin{cases} \dfrac{1}{abc} & \text{if } abc \neq 0 \\ 1 & \text{if } abc = 0 \end{cases}$$

$\eta_1^3 : V_1^3 \to [0, 1]$ is such that

$$\eta_1^3 \begin{pmatrix} a & a & a \\ b & b & b \\ c & c & c \end{pmatrix} = \begin{cases} \dfrac{1}{a+b+c} & \text{if } a+b+c \neq 0 \\ 1 & \text{if } a+b+c = 0 \end{cases}$$

$\eta_2^3 : V_2^3 \to [0, 1]$ is given by

$$\eta_2^3 (a\ b\ c\ d\ e) = \begin{cases} \dfrac{1}{abc + de} & \text{if } abc + de \neq 0 \\ 1 & \text{if } abc + de = 0 \end{cases}$$

$\eta_3^3 : V_3^3 \to [0, 1]$ is define by

$$\eta_3^3 \begin{pmatrix} a & a & a & a & a \\ b & b & b & b & b \\ c & c & c & c & c \end{pmatrix} = \begin{cases} \dfrac{1}{abc} & \text{if } abc \neq 0 \\ 1 & \text{if } abc = 0 \end{cases}$$

$\eta_1^4 : V_1^4 \to [0,1]$ such that

$$\eta_1^4 (a\ b\ c\ d\ e) = \begin{cases} \dfrac{1}{a+b+c+d+e} & \text{if } a+b+c+d+e \neq 0 \\ 1 & \text{if } a+b+c+d+e = 0 \end{cases}$$

$\eta_2^4 : V_2^4 \to [0, 1]$ is defined by

$$\eta_2^4 \begin{pmatrix} a & b & c & d & e \\ 0 & f & g & h & i \\ 0 & 0 & j & k & l \\ 0 & 0 & 0 & m & n \\ 0 & 0 & 0 & 0 & p \end{pmatrix} = \begin{cases} \dfrac{1}{afjmp} & \text{if } afjmp \neq 0 \\ 1 & \text{if } afjmp = 0 \end{cases}$$

$\eta_1^5 : V_1^5 \to [0,1]$ is defined by

$$\eta_1^5 \begin{pmatrix} a & 0 & 0 \\ a & a & 0 \\ a & a & a \end{pmatrix} = \begin{cases} \dfrac{1}{6a} & \text{if } a \neq 0 \\ 1 & \text{if } a = 0 \end{cases}$$

$\eta_2^5 : V_2^5 \to [0,1]$ is such that

$$\eta_2^5 \,(a\ b\ c\ d\ e) = \begin{cases} \dfrac{1}{a} & \text{if } a \neq 0 \\ 1 & \text{if } a = 0 \end{cases}$$

$\eta_3^5 : V_3^5 \to [0,1]$ is defined by

$$\eta_3^5 \begin{bmatrix} a & b & c & d \\ a & b & c & d \\ a & b & c & d \end{bmatrix} = \begin{cases} \dfrac{1}{abcd} & \text{if } abcd \neq 0 \\ 1 & \text{if } abcd = 0 \end{cases}$$

$\eta_4^5 : V_4^5 \to [0,1]$ is given by

$$\eta_4^5 \begin{pmatrix} a & a & a & a \\ b & b & b & b \\ c & c & c & c \\ d & d & d & d \end{pmatrix} = \begin{cases} \dfrac{1}{abcd} & \text{if } abcd \neq 0 \\ 1 & \text{if } abcd = 0 \end{cases}$$

Thus
$$\begin{aligned}
V\eta &= V(\eta_1 \cup \eta_2 \cup \eta_3 \cup \eta_4 \cup \eta_5) \\
&= (V_1 \cup V_2 \cup V_3 \cup V_4 \cup V_5)(\eta_1 \cup \eta_2 \cup \eta_3 \cup \eta_4 \cup \eta_5) \\
&= (V_1\eta_1 \cup V_2\eta_2 \cup V_3\eta_3 \cup V_4\eta_4 \cup V_5\eta_5) \\
&= \left(V^1_{1\eta^1_1}, V^1_{2\eta^1_2}, V^1_{3\eta^1_3}\right) \cup \left(V^2_{1\eta^2_1}, V^2_{2\eta^2_2}\right) \cup \left(V^3_{1\eta^3_1}, V^3_{2\eta^3_2}, V^3_{3\eta^3_3}\right) \\
&\quad \cup \left(V^4_{1\eta^4_1}, V^4_{2\eta^4_2}\right) \cup \left(V^5_{1\eta^5_1}, V^5_{2\eta^5_2}, V^5_{3\eta^5_3}, V^5_{4\eta^5_4}\right)
\end{aligned}$$

is a special set fuzzy 5-vector space.

Now for this we give a special set fuzzy 5-vector subspace in the following.

Now take
$$W_1^1 = (S \times S \times \{0\} \times \{0\}) \subseteq V_1^1,$$

$$W_2^1 = \left\{ \begin{pmatrix} a & b \\ c & d \end{pmatrix} \middle| a,b,c,d \in 7Z^\circ \right\} \subseteq V_2^1,$$

$$W_3^1 = \left\{ \begin{pmatrix} a & a & a & a & a & a \\ a & a & a & a & a & a \end{pmatrix} \middle| a \in S \right\} \subseteq V_3^1,$$

$W_1^2 = $ {all 4×4 upper triangular matrices of the form

$$\left\{ \begin{pmatrix} a & a & a & a \\ 0 & b & b & b \\ 0 & 0 & a & a \\ 0 & 0 & 0 & b \end{pmatrix} \middle| a,b \in 3Z^o \right\} \subseteq V_1^2,$$

$$W_2^2 = S \times \{0\} \times \{0\} \subseteq V_2^2,$$

$$W_1^3 = \left\{ \begin{pmatrix} a & a & a \\ b & b & b \\ c & c & c \end{pmatrix} \middle| a,b,c \in 2Z^o \right\} \subseteq V_1^3,$$

$$W_2^3 = (S \times \{0\} \times S \times \{0\} \times S) \subseteq V_2^3,$$

$$W_3^3 = \left\{ \begin{pmatrix} a & a & a & a & a \\ b & b & b & b & b \\ c & c & c & c & c \end{pmatrix} \middle| a,b,c \in 11Z^o \right\} \subseteq V_3^3,$$

$W_1^4 = (S \times S \times \{0\} \times \{0\} \times S) \subseteq V_1^4$, $W_2^4 = $ {all 5×5 upper triangular matrices with entries from $13Z^o$} $\subseteq V_2^4$,

$W_1^5 = $ {low triangular matrices of the form

$$\left\{ \begin{pmatrix} a & 0 & 0 \\ a & a & 0 \\ a & a & a \end{pmatrix} \middle| a \in 12Z^o \right\} \subseteq V_1^5,$$

$$W_2^5 = \{S \times \{0\} \times \{0\} \times \{0\} \times S\} \subseteq V_2^5,$$

$$W_3^5 = \left\{ \begin{pmatrix} a & b & a & b \\ a & b & a & b \\ a & b & a & b \end{pmatrix} \middle| a,b \in S \right\} \subseteq V_3^5,$$

$$W_4^5 = \left\{ \begin{pmatrix} a & a & a & a \\ b & b & b & b \\ c & c & c & c \\ d & d & d & d \end{pmatrix} \middle| a,b,c,d \in 23Z^o \right\} \subseteq V_4^5,$$

here $W_{j_i}^i \subseteq V_{j_i}^i$ is a special set vector subspace of $V_{j_i}^i$, $1 \le j_i \le n_i$, $1 \le i \le 5$. Now define a 5-map $\eta = (\eta_1 \cup \eta_2 \cup \eta_3 \cup \eta_4 \cup \eta_5)$: W = $(W_1 \cup W_2 \cup \ldots \cup W_5) \to [0, 1]$ by $\eta_i: W_i \to [0,1]$, i = 1, 2, 3, 4, 5 such that

$$\eta_1 = (\eta_1^1, \eta_2^1, \eta_3^1) : (W_1^1, W_2^1, W_3^1) \to [0, 1]$$

by $\eta_1^1 : W_1^1 \to [0, 1]$ is given by

$$\eta_1^1 (a\ b\ 0\ 0) = \begin{cases} \dfrac{1}{ab} & \text{if } ab \ne 0 \\ 1 & \text{if } ab = 0 \end{cases}$$

$\eta_2^1 : W_2^1 \to [0, 1]$ is defined by

$$\eta_2^1 \begin{bmatrix} a & b \\ c & d \end{bmatrix} = \begin{cases} \dfrac{1}{a+d} & \text{if } a+d \ne 0 \\ 1 & \text{if } a+d = 0 \end{cases}$$

$\eta_3^1 : W_3^1 \to [0, 1]$ is such that

$$\eta_3^1 \begin{bmatrix} a & a & a & a & a & a \\ a & a & a & a & a & a \end{bmatrix} = \begin{cases} \dfrac{1}{12a} & \text{if } a \ne 0 \\ 1 & \text{if } a = 0 \end{cases}$$

$W_1\eta_1 = \left(W^1_{1\eta^1_1}, W^1_{2\eta^1_2}, W^1_{3\eta^1_3}\right)$ is a special fuzzy vector subspace.

$$\eta_2 = \left(\eta^2_1, \eta^2_2\right) : W_2 = \left(W^2_1, W^2_2\right) \to [0,1]$$

is given by

$\eta^2_1 : W^2_1 \to [0,1]$ is such that

$$\eta^2_1 \begin{bmatrix} a & a & a & a \\ 0 & b & b & b \\ 0 & 0 & a & a \\ 0 & 0 & 0 & b \end{bmatrix} = \begin{cases} \dfrac{1}{ab} & \text{if } ab \neq 0 \\ 1 & \text{if } ab = 0 \end{cases}$$

$\eta^2_2 : W^2_2 \to [0,1]$ is defined by

$$\eta^2_2 (a\ 0\ 0) = \begin{cases} \dfrac{1}{a} & \text{if } a \neq 0 \\ 1 & \text{if } a = 0 \end{cases}$$

$$\eta_3 = \left(\eta^3_1, \eta^3_2, \eta^3_3\right) : W_3 \to [0,1];$$

$\eta^3_1 : W^3_1 \to [0,1]$ is such that

$$\eta^3_1 \begin{pmatrix} a & a & a \\ b & b & b \\ c & c & c \end{pmatrix} = \begin{cases} \dfrac{1}{a+b+c} & \text{if } a+b+c \neq 0 \\ 1 & \text{if } a+b+c = 0 \end{cases}$$

$\eta^3_2 : W^3_2 \to [0,1]$ is given by

$$\eta^3_2 (a\ 0\ a\ 0\ a) = \begin{cases} \dfrac{1}{3a} & \text{if } a \neq 0 \\ 1 & \text{if } a = 0 \end{cases}$$

$\eta^3_3 : W^3_3 \to [0,1]$ is such that

$$\eta_3^3 \begin{pmatrix} a & a & a & a & a \\ b & b & b & b & b \\ c & c & c & c & c \end{pmatrix} = \begin{cases} \dfrac{1}{a+b+c} & \text{if } a+b+c \neq 0 \\ 1 & \text{if } a+b+c = 0 \end{cases}$$

Thus $W_3 \eta_3 = \left( W^3_{1\eta_1^3}, W^3_{2\eta_2^3}, W^3_{3\eta_3^3} \right)$ is a special set fuzzy vector subspace.

$$\eta_4 = \left( \eta_1^4, \eta_2^4 \right) : W_4 = \left( W_1^4, W_2^4 \right) \to [0,1]$$

such that

$\eta_1^4 : W_1^4 \to [0,1]$ is given by

$$\eta_1^4 \, (a\ b\ 0\ 0\ e) = \begin{cases} \dfrac{1}{a+b+e} & \text{if } a+b+e \neq 0 \\ 1 & \text{if } a+b+e = 0 \end{cases}$$

$\eta_2^4 : W_2^4 \to [0,1]$ is such that

$$\eta_2^4 \begin{pmatrix} a & b & c & d & e \\ 0 & f & g & h & i \\ 0 & 0 & j & k & l \\ 0 & 0 & 0 & m & n \\ 0 & 0 & 0 & 0 & p \end{pmatrix} = \begin{cases} \dfrac{1}{a+p} & \text{if } a+p \neq 0 \\ 1 & \text{if } a+p = 0 \end{cases}$$

$$\eta_5 : W_5 = \left( W_1^5, W_2^5, W_3^5, W_4^5 \right) \to [0,1]$$

$\eta_1^5 : W_1^5 \to [0,1]$ is given by

$$\eta_1^5 \begin{bmatrix} a & 0 & 0 \\ a & a & 0 \\ a & a & a \end{bmatrix} = \begin{cases} \dfrac{1}{a} & \text{if } a \neq 0 \\ 1 & \text{if } a = 0 \end{cases}$$

$\eta_2^5 : W_2^5 \to [0, 1]$ is such that

$$\eta_2^5 (a\ 0\ 0\ 0\ b) = \begin{cases} \dfrac{1}{a+b} & \text{if } a+b \neq 0 \\ 1 & \text{if } a+b = 0 \end{cases}$$

$\eta_3^5 : W_3^5 \to [0, 1]$ so that

$$\eta_3^5 \begin{bmatrix} a & b & a & b \\ a & b & a & b \\ a & b & a & b \end{bmatrix} = \begin{cases} \dfrac{1}{ab} & \text{if } ab \neq 0 \\ 1 & \text{if } ab = 0 \end{cases}$$

$\eta_4^5 : W_4^5 \to [0, 1]$ is such that

$$\eta_4^5 \begin{bmatrix} a & a & a & a \\ b & b & b & b \\ c & c & c & c \\ d & d & d & d \end{bmatrix} = \begin{cases} \dfrac{1}{a+b+c+d} & \text{if } a+b+c+d \neq 0 \\ 1 & \text{if } a+b+c+d = 0 \end{cases}$$

Thus $W_5\eta_5 = \left( W_{1\eta_1^5}^5, W_{2\eta_2^5}^5, W_{3\eta_3^5}^5, W_{4\eta_4^5}^5 \right)$ is a special set fuzzy vector subspace.

Thus
$W\eta = (W_1 \cup W_2 \cup W_3 \cup W_4 \cup W_5)_\eta$
$= \left( W_{1\eta_1^1}^1, W_{2\eta_2^1}^1, W_{3\eta_3^1}^1 \right) \cup \left( W_{1\eta_1^2}^2, W_{2\eta_2^2}^2 \right) \cup \left( W_{1\eta_1^3}^3, W_{2\eta_2^3}^3, W_{3\eta_3^3}^3 \right)$
$\cup \left( W_{1\eta_1^4}^4, W_{2\eta_2^4}^4 \right) \cup \left( W_{1\eta_1^5}^5, W_{2\eta_2^5}^5, W_{3\eta_3^5}^5, W_{4\eta_4^5}^5 \right)$

is a special set fuzzy vector 4-subspace.

We now proceed onto define the notion of special set fuzzy linear n-algebra. It is pertinent to mention here that infact the special set fuzzy vector n-space and special set linear n-algebra are fuzzy equivalent.

**DEFINITION 2.4.11:** *Let*
$$V = (V_1 \cup V_2 \cup \ldots \cup V_n)$$
$$= \left(V_1^1, V_2^1, \ldots, V_{n_1}^1\right) \cup \left(V_1^2, V_2^2, \ldots, V_{n_2}^2\right) \cup \ldots \cup \left(V_1^n, V_2^n, \ldots, V_{n_n}^n\right)$$
*be a special set n-linear algebra over the set S. Take $\eta = (\eta_1 \cup \eta_2 \cup \ldots \cup \eta_n)$ a n-map from $V = (V_1 \cup V_2 \cup \ldots \cup V_n)$ into [0,1], i.e.,*
$$(\eta_1 \cup \eta_2 \cup \ldots \cup \eta_n)$$
$$= \left(\eta_1^1, \eta_2^1, \ldots, \eta_{n_1}^1\right) \cup \left(\eta_1^2, \eta_2^2, \ldots, \eta_{n_2}^2\right) \cup \ldots \cup \left(\eta_1^n, \eta_2^n, \ldots, \eta_{n_n}^n\right):$$
$$V = (V_1 \cup V_2 \cup \ldots \cup V_n)$$
$$= \left(V_1^1, V_2^1, \ldots, V_{n_1}^1\right) \cup \left(V_1^2, V_2^2, \ldots, V_{n_2}^2\right) \cup \ldots \cup \left(V_1^n, V_2^n, \ldots, V_{n_n}^n\right) \to$$
*[0, 1]; $\eta_i = \left(\eta_1^i, \ldots, \eta_{n_i}^i\right) : V_i \to [0,1]$ is defined such that $\eta_{t_i}^i : W_{t_i}^i \to [0,1]$ so that $V_{t_i \eta_{t_i}^i}^i$ is a set fuzzy linear algebra, $1 \le i \le n$, $1 \le t_i \le n_i$; $i = 1, 2, \ldots, n$.*

*Thus $V_i \eta_i$ is a special set fuzzy linear algebra.*

$$V\eta = (V_1 \cup V_2 \cup \ldots \cup V_n)_{\eta_1 \cup \ldots \cup \eta_n}$$
$$= \left(V_{1\eta_1^1}^1, V_{2\eta_2^1}^1, \ldots, V_{n_1 \eta_{n_1}^1}^1\right) \cup \left(V_{1\eta_1^2}^2, V_{2\eta_2^2}^2, \ldots, V_{n_2 \eta_{n_2}^2}^2\right) \cup \ldots \cup$$
$$\left(V_{1\eta_1^n}^n, V_{2\eta_2^n}^n, \ldots, V_{n_n \eta_{n_n}^n}^n\right)$$

*is a special set fuzzy n linear algebra.*

We shall illustrate this by the following example.

***Example 2.4.10:*** Let
$$V = V_1 \cup V_2 \cup V_3 \cup V_4$$
$$= \left(V_1^1, V_2^1\right) \cup \left(V_1^2, V_2^2, V_3^2\right) \cup \left(V_1^3, V_2^3, V_3^3, V_4^3\right) \cup \left(V_1^4, V_2^4, V_3^4\right)$$
be a special set 4 linear algebra over the set $S = \{0, 1\}$.
Here

$$V_1^1 = \left\{ \begin{pmatrix} a & b \\ c & d \end{pmatrix} \middle| a, b, c, d \in Z^o \right\}$$

$$V_2^1 = \left\{ \begin{pmatrix} a & a & a & a \\ a & a & a & a \end{pmatrix} \middle| a \in Z^o \right\},$$

$$V_1^2 = \left\{ \begin{pmatrix} a & b & c \\ d & e & f \\ g & h & i \end{pmatrix} \middle| a,b,c,d,e,f,g,h,i \in Z^o \right\},$$

$$V_2^2 = \left\{ \begin{pmatrix} a & a \\ a & a \\ a & a \\ a & a \end{pmatrix} \middle| a \in Z^o \right\},$$

$$V_3^2 = Z^o \times Z^o \times Z^o \times Z^o, V_1^3 = Z^o \times Z^o \times Z^o,$$

$$V_2^3 = \left\{ \begin{pmatrix} a & o \\ b & d \end{pmatrix} \middle| a,b,d \in Z^o \right\},$$

$$V_3^3 = \left\{ \begin{pmatrix} a & a & a & a & a \\ b & b & b & b & b \end{pmatrix} \middle| a,b \in Z^o \right\},$$

$$V_4^3 = \left\{ \begin{pmatrix} a & a & a \\ b & b & b \\ c & c & c \\ d & d & d \\ e & e & e \end{pmatrix} \middle| a,b,c,d,e \in Z^o \right\},$$

$$V_1^4 = \left\{ \begin{pmatrix} a & a & a & a & a \\ b & b & b & b & b \\ c & c & c & c & c \end{pmatrix} \middle| a,b,c \in Z^o \right\},$$

$$V_2^4 = \left\{ \begin{pmatrix} a & 0 & 0 & 0 \\ a & a & 0 & 0 \\ b & b & b & 0 \\ b & b & b & b \end{pmatrix} \middle| \text{ s.t } a, b \in Z^o \right\}$$

and

$$V_3^4 = \{Z^o \times Z^o \times Z^o \times Z^o \times Z^o \times Z^o\}.$$

$$\eta = \eta_1 \cup \eta_2 \cup \eta_3 \cup \eta_4:$$
$$\left(\eta_1^1, \eta_2^1\right) \cup \left(\eta_1^2, \eta_2^2, \eta_3^2\right) \cup \left(\eta_1^3, \eta_2^3, \eta_3^3, \eta_4^3\right) \cup \left(\eta_1^4, \eta_2^4, \eta_3^4\right):$$
$$V = (V_1 \cup V_2 \cup V_3 \cup V_4) \to [0,1]$$

where $\eta_{j_i}^i : V_{j_i}^i \to [0, 1]$; $i = 1, 2, 3, 4$; $1 \le t_i \le n_i$ ; $1 \le i \le 4$ such that each $V_{t_i \eta_{t_i}^i}^i$ is a set fuzzy vector space, so that $V_i \eta_i$ is a special set fuzzy vector space; we define for each $1 \le t_i \le n_i$, $1 \le i \le 4$ as follows:

$\eta_1^1 : V_1^1 \to [0,1]$ is such that

$$\eta_1^1 \begin{pmatrix} a & b \\ c & d \end{pmatrix} = \begin{cases} \dfrac{1}{ad - bc} & \text{if } ad \ne bc \\ 1 & \text{if } ad = bc \end{cases}$$

$\eta_2^1 : V_2^1 \to [0, 1]$ is given by

$$\eta_2^1 \begin{pmatrix} a & a & a & a \\ a & a & a & a \end{pmatrix} = \begin{cases} \dfrac{1}{8a} & \text{if } a \ne 0 \\ 1 & \text{if } a = 0 \end{cases}$$

$\eta_1^2 : V_1^2 \to [0, 1]$ is defined by

$$\eta_1^2 \begin{pmatrix} a & b & c \\ d & e & f \\ g & h & i \end{pmatrix} = \begin{cases} \dfrac{1}{abc} & \text{if } abc \ne 0 \\ 1 & \text{if } abc = 0 \end{cases}$$

$\eta_2^2 : V_2^2 \to [0, 1]$ is given by

$$\eta_2^2 \begin{bmatrix} a & a \\ a & a \\ a & a \\ a & a \end{bmatrix} = \begin{cases} \dfrac{1}{2a} & \text{if } a \neq 0 \\ 1 & \text{if } a = 0 \end{cases}$$

$\eta_3^2 : V_3^2 \to [0,1]$ is defined by

$$\eta_3^2 (a\ b\ c\ d) = \begin{cases} \dfrac{1}{abc + d} & \text{if } abc + d \neq 0 \\ 1 & \text{if } abc + d = 0 \end{cases}$$

$\eta_1^3 : V_1^3 \to [0,1]$ is such that

$$\eta_1^3 (a\ b\ c) = \begin{cases} \dfrac{1}{a + b + c} & \text{if } a + b + c \neq 0 \\ 1 & \text{if } a + b + c = 0 \end{cases}$$

$\eta_2^3 : V_2^3 \to [0,1]$ is defined to be

$$\eta_2^3 \begin{pmatrix} a & 0 \\ b & d \end{pmatrix} = \begin{cases} \dfrac{1}{ad} & \text{if } ad \neq 0 \\ 1 & \text{if } ad = 0 \end{cases}$$

$\eta_3^3 : V_3^3 \to [0,1]$ is given by

$$\eta_3^3 \begin{pmatrix} a & a & a & a & a \\ b & b & b & b & b \end{pmatrix} = \begin{cases} \dfrac{1}{3a + 2b} & \text{if } 3a + 2b \neq 0 \\ 1 & \text{if } 3a + 2b = 0 \end{cases}$$

$\eta_4^3 : V_4^3 \to [0,1]$ is given by

$$\eta_4^3 \begin{bmatrix} a & a & a \\ b & b & b \\ c & c & c \\ d & d & d \\ e & e & e \end{bmatrix} = \begin{cases} \dfrac{1}{a+b+c+d+e} & \text{if } a+b+c+d+e \neq 0 \\ 1 & \text{if } a+b+c+d+e = 0 \end{cases}$$

$\eta_1^4 : V_1^4 \to [0, 1]$ is such that

$$\eta_1^4 \begin{pmatrix} a & a & a & a & a \\ b & b & b & b & b \\ c & c & c & c & c \end{pmatrix} = \begin{cases} \dfrac{1}{5a+5b+5c} & \text{if } 5a+5b+5c \neq 0 \\ 1 & \text{if } 5a+5b+5c = 0 \end{cases}$$

$\eta_2^4 : V_2^4 \to [0, 1]$ is given by

$$\eta_2^4 \begin{bmatrix} a & 0 & 0 & 0 \\ a & a & 0 & 0 \\ b & b & b & 0 \\ b & b & b & b \end{bmatrix} = \begin{cases} \dfrac{1}{3a+7b} & \text{if } 3a+7b \neq 0 \\ 1 & \text{if } 3a+7b = 0 \end{cases}$$

$\eta_3^4 : V_3^4 \to [0,1]$ is such that

$$\eta_3^4 (a\ b\ c\ d\ e\ f) = \begin{cases} \dfrac{1}{a+f} & \text{if } a+f \neq 0 \\ 1 & \text{if } a+f = 0 \end{cases}$$

Thus
$$V\eta = (V_1 \cup V_2 \cup V_3 \cup V_4)_{\eta1 \cup \eta2 \cup \eta3 \cup \eta4}$$
$$= \left( V^1_{1\eta_1^1}, V^1_{2\eta_2^1} \right) \cup \left( V^2_{1\eta_1^2}, V^2_{2\eta_2^2}, V^2_{3\eta_3^2} \right) \cup \left( V^3_{1\eta_1^3}, V^3_{2\eta_2^3}, V^3_{3\eta_3^3}, V^3_{4\eta_4^3} \right)$$
$$\cup \left( V^4_{1\eta_1^4}, V^4_{2\eta_2^4}, V^4_{3\eta_3^4} \right)$$
is a special set fuzzy linear 4-algebra.

Next we proceed onto define the notion of special set fuzzy linear n-subalgebra.

**DEFINITION 2.4.12:** *Let*
$$V = (V_1 \cup V_2 \cup ... \cup V_n)$$
$$= \left(V_1^1, V_2^1,...,V_{n_1}^1\right) \cup \left(V_1^2, V_2^2,...,V_{n_2}^2\right) \cup ... \cup \left(V_1^n, V_2^n,...,V_{n_n}^n\right)$$
*be a special set linear n – algebra defined over the set S. Suppose*
$$W = (W_1, W_2, ... , W_n)$$
$$= \left(W_1^1, W_2^1,...,W_{n_1}^1\right) \cup \left(W_1^2, W_2^2,...,W_{n_2}^2\right) \cup ... \cup \left(W_1^n, W_2^n,...,W_{n_n}^n\right)$$
$$\subseteq (V_1 \cup V_2 \cup ... \cup V_n)$$
*be a proper special set linear n-subalgebra of V over S. Define* $\eta = (\eta_1, \cup \eta_2, \cup ... \cup \eta_n) : : (W_1, W_2, ..., W_n) \to [0,1]$ *such that*

$$W\eta = (W_1 \cup W_2 \cup ... \cup W_n)_{\eta_1 \cup \eta_2 \cup ... \cup \eta_n}$$
$$= \left(W_{1\eta_1^1}^1, W_{2\eta_1^2}^1,...,W_{n_1\eta_{n_1}^1}^1\right) \cup \left(W_{1\eta_1^2}^2, W_{2\eta_2^2}^2,...,W_{n_2\eta_{n_2}^2}^2\right)$$
$$\cup ... \cup \left(W_{1\eta_1^n}^n, W_{2\eta_1^n}^n,...,W_{n_n\eta_{n_n}^n}^n\right)$$

*be a special set fuzzy linear n-algebra then we call W$\eta$ to be the special set fuzzy linear n-subalgebra where*
$$\eta = (\eta_1, \cup \eta_2, \cup ... \cup \eta_n)$$
$$= \left(\eta_1^1, \eta_2^1,...,\eta_{n_1}^1\right) \cup \left(\eta_1^2, \eta_2^2,...,\eta_{n_2}^2\right) \cup ... \cup \left(\eta_1^n, \eta_2^n,...,\eta_{n_n}^n\right)$$
*is a map such that $W\eta = (W_1\eta_1, W_2\eta_2, ..., W_n\eta_n)$ where $W_i\eta_i$ is a special set fuzzy linear subalgebra, i = 1, 2, ..., n.*

Now we illustrate this by a simple example.

***Example 2.4.11:*** Let $V = (V_1 \cup V_2 \cup V_3 \cup V_4)$ where $V_1 = \left(V_1^1, V_2^1, V_3^1, V_4^1\right)$, $V_2 = \left(V_1^2, V_2^2, V_3^2, V_4^2\right)$, $V_3 = \left(V_1^3, V_2^3\right)$ and $V_4 = \left(V_1^4, V_2^4, V_3^4\right)$ be a special set linear 4-algebra over the set $S = Z^o$, where $V_1^1 = S \times S \times S$, $V_2^1 = $ {all polynomials of degree less than or equal to four with coefficients from $Z^o$}, $V_3^1 = $ {all 4×4 upper triangular matrices with entries from $Z^o$} and

$$V_4^1 = \left\{ \begin{pmatrix} a & a & a & a & a \\ b & b & b & b & b \end{pmatrix} \middle| a, b \in Z^o \right\}.$$

$V_1^2 = \{3 \times 3$ matrices with entries from $Z^o\}$,

$$V_2^2 = \left\{ \begin{pmatrix} a \\ a \\ a \\ a \\ a \end{pmatrix} \middle| a \in Z^o \right\}$$

$V_3^2 = S \times S \times S \times S$ and $V_4^2 = \{5 \times 5$ lower triangular matrices with entries from $S\}$. $V_1^3 = S \times S \times S \times S \times S$ and

$$V_2^3 = \left\{ \begin{pmatrix} a & a & a \\ b & b & b \end{pmatrix} \middle| a, b \in Z^o \right\}.$$

$V_1^4 = S \times S \times S \times S$, $V_2^4 = \{4 \times 4$ lower triangular matrices with entries from $Z^o\}$ and

$$V_3^4 = \left\{ \begin{pmatrix} a & b \\ a & b \\ a & b \\ a & b \\ a & b \end{pmatrix} \middle| a, b \in Z^o \right\}.$$

Consider $W = (W_1 \cup W_2 \cup W_3 \cup W_4) \subseteq V$ where $W_1 = \left( W_1^1, W_2^1, W_3^1, W_4^1 \right) \subseteq V_1$ is such that $W_1^1 = S \times \{0\} \times S \subseteq V_1^1$, $W_2^1 = \{$all polynomials of degree less than or equal to four with coefficients from $3Z^o$} $\subseteq V_2^1$, $W_3^1 = $ {all $4 \times 4$ upper triangular matrices with entries from $5Z^o$} $\subseteq V_3^1$, and

$$W_4^1 = \left\{ \begin{pmatrix} a & a & a & a & a \\ b & b & b & b & b \end{pmatrix} \middle| a, b \in 7Z^o \right\} \subseteq V_4^1,$$

$W_1^2 = $ {$3 \times 3$ matrices with entries from $11Z^o$} $\subseteq V_1^2$,

$$V_2^2 = \left\{ \begin{pmatrix} a \\ a \\ a \\ a \\ a \end{pmatrix} \middle| a \in 2Z^o \right\} \subseteq V_2^2,$$

$W_3^2 = \{S \times S \times \{0\} \times \{0\}\} \subseteq V_3^2$ and $W_4^2 = $ {$5 \times 5$ lower triangular matrices with entries from S} $\subseteq V_4^2$, $W_1^3 = S \times S \times \{0\} \times \{0\} \times \{0\} \subseteq V_1^3$,

$$W_2^3 = \left\{ \begin{pmatrix} a & a & a \\ a & a & a \end{pmatrix} \middle| a \in Z^o \right\} \subseteq V_2^3,$$

$W_1^4 = \{S \times \{0\} \times \{0\} \times \{0\}\} \subseteq V_1^4$, $W_2^4 = $ {$4 \times 4$ lower triangular matrices with entries from $Z^o$} $\subseteq V_2^4$ and

$$W_3^4 = \left\{ \begin{pmatrix} a & a \\ a & a \\ a & a \\ a & a \\ a & a \end{pmatrix} \middle| a \in Z^o \right\} \subseteq V_3^4.$$

Define

$$\eta = (\eta_1 \cup \eta_2 \cup \eta_3 \cup \eta_4):$$
$$W = (W_1 \cup W_2 \cup W_3 \cup W_4) \to [0, 1]$$

by the following rule,

$\eta_1 : W_1 \to [0,1]$ is such that
$$\eta_1 = \left(\eta_1^1, \eta_2^1, \eta_3^1, \eta_4^1\right) : \left(W_1^1, W_2^1, W_3^1, W_4^1\right) \to [0, 1]$$
given by

$\eta_1^1 : W_1^1 \to [0, 1]$ is defined by

$$\eta_1^1 \, (a \; 0 \; b) = \begin{cases} \dfrac{1}{ab} & \text{if } ab \neq 0 \\ 1 & \text{if } ab = 0 \end{cases}$$

$\eta_2^1 : W_2^1 \to [0, 1]$ is defined by

$$\eta_2^1 \, (p\,(x)) = \begin{cases} \dfrac{1}{\deg p(x)} & \text{if } p(x) \text{ is not a constant} \\ 1 & \text{if } p(x) \text{ is a constant} \end{cases}$$

$\eta_3^1 : W_3^1 \to [0, 1]$ is given by

$$\eta_3^1 \begin{bmatrix} a & b & c & d \\ 0 & e & f & g \\ 0 & 0 & h & i \\ 0 & 0 & 0 & j \end{bmatrix} = \begin{cases} \dfrac{1}{abcd + efg} & \text{if } abcd + efg \neq 0 \\ 1 & \text{if } abcd + efg = 0 \end{cases}$$

$\eta_4^1 : W_4^1 \to [0,1]$ is defined by

$$\eta_4^1 \begin{bmatrix} a & a & a & a & a \\ b & b & b & b & b \end{bmatrix} = \begin{cases} \dfrac{1}{ab} & \text{if } ab \neq 0 \\ 1 & \text{if } ab = 0 \end{cases}$$

Thus
$$W_1\eta_1 = \left(W^1_{1\eta^1_1}, W^1_{2\eta^1_2}, W^1_{3\eta^1_3}, W^1_{4\eta^1_4}\right)$$
is a special set fuzzy linear subalgebra.

$$\eta_2 = \left(\eta^2_1, \eta^2_2, \eta^2_3, \eta^2_4\right) : W_2 = \left(W^2_1, W^2_2, W^2_3, W^2_4\right) \to [0, 1]$$
is such that

$\eta^2_1 : W^2_1 \to [0, 1]$ is defined by

$$\eta^2_1 \begin{bmatrix} a & b & c \\ d & e & f \\ g & h & i \end{bmatrix} = \begin{cases} \dfrac{1}{abcdg} & \text{if } abcdg \neq 0 \\ 1 & \text{if } abcdg = 0 \end{cases}$$

$\eta^2_2 : W^2_2 \to [0,1]$ such that

$$\eta^2_2 \begin{bmatrix} a \\ a \\ a \\ a \\ a \end{bmatrix} = \begin{cases} \dfrac{1}{5a} & \text{if } a \neq 0 \\ 1 & \text{if } a = 0 \end{cases}$$

$\eta^2_3 : W^2_3 \to [0,1]$ is such that

$$\eta^2_3 \,(a\ b\ 0\ 0) = \begin{cases} \dfrac{1}{a+b} & \text{if } a+b \neq 0 \\ 1 & \text{if } a+b = 0 \end{cases}$$

$\eta^2_4 : W^2_4 \to [0,1]$ is given by

$$\eta_4^2 \begin{pmatrix} a & 0 & 0 & 0 & 0 \\ d & e & 0 & 0 & 0 \\ g & h & i & 0 & 0 \\ j & k & l & m & 0 \\ o & p & q & r & s \end{pmatrix} =$$

$$\begin{cases} \dfrac{1}{a+e+i+m+s} & \text{if } a+e+i+m+s \neq 0 \\ 1 & \text{if } a+e+i+m+s = 0 \end{cases}$$

Thus $W_2\eta_2 = \left(W^2_{1\eta_1^2}, W^2_{2\eta_2^2}, W^2_{3\eta_3^2}, W^2_{4\eta_4^2}\right)$ is a special set fuzzy linear subalgebra.

$\eta_3: W_3 \to [0, 1]; \eta_3 = \left(\eta_1^3, \eta_2^3\right) : \left(W_1^3, W_2^3\right) \to [0, 1]$ is given by;

$\eta_1^3 : W_1^3 \to [0, 1]$ is such that

$$\eta_1^3 (a\ b\ 0\ 0\ 0) = \begin{cases} \dfrac{1}{ab} & \text{if } ab \neq 0 \\ 1 & \text{if } ab = 0 \end{cases}$$

$\eta_2^3 : W_2^3 \to [0,1]$ is defined by

$$\eta_2^3 \begin{bmatrix} a & a & a \\ a & a & a \end{bmatrix} = \begin{cases} \dfrac{1}{6a} & \text{if } a \neq 0 \\ 1 & \text{if } a = 0 \end{cases}$$

Thus $W_3\eta_3 = \left(W^3_{1\eta_1^3}, W^3_{2\eta_2^3}\right)$ is a special set fuzzy linear sub algebra.

Define $\eta_4 : W_4 \to [0,1]$ as follows.

$$\eta_4 = \left(\eta_1^4, \eta_2^4, \eta_3^4\right) : \left(W_1^4, W_2^4, W_3^4\right) \to [0, 1];$$

$\eta_1^4 : W_1^4 \to [0, 1]$ is defined by

$$\eta_1^4 \,(a\ 0\ 0\ 0) = \begin{cases} \dfrac{1}{a} & \text{if } a \neq 0 \\ 1 & \text{if } a = 0 \end{cases}$$

$\eta_2^4 : W_2^4 \to [0, 1]$ is given by

$$\eta_2^4 \begin{pmatrix} a & 0 & 0 & 0 \\ b & c & 0 & 0 \\ d & e & f & 0 \\ g & h & i & j \end{pmatrix} = \begin{cases} \dfrac{1}{acfj} & \text{if } acfj \neq 0 \\ 1 & \text{if } acfj = 0 \end{cases}$$

$\eta_3^4 : W_3^4 \to [0, 1]$

$$\eta_3^4 \begin{pmatrix} a & a \\ a & a \\ a & a \\ a & a \\ a & a \end{pmatrix} = \begin{cases} \dfrac{1}{10a} & \text{if } a \neq 0 \\ 1 & \text{if } a = 0 \end{cases}$$

Thus $W_4\eta_4 = \left(W_{1\eta_1^4}^4, W_{2\eta_2^4}^4, W_{3\eta_3^4}^4\right)$ is a special set fuzzy linear subalgebra of $V_4$.

Now

$$\begin{aligned}
W\eta &= (W_1 \cup W_2 \cup W_3 \cup W_4 \cup W_5)_\eta \\
&= (W_1\eta_1 \cup W_2\eta_2 \cup W_3\eta_3 \cup W_4\eta_4) \\
&= \left(W_{1\eta_1^1}^1, W_{2\eta_2^1}^1, W_{3\eta_3^1}^1, W_{4\eta_4^1}^1\right) \cup \left(W_{1\eta_1^2}^2, W_{2\eta_2^2}^2, W_{3\eta_3^2}^2, W_{4\eta_4^2}^2\right) \\
&\quad \cup \left(W_{1\eta_1^3}^3, W_{2\eta_2^3}^3\right) \cup \left(W_{1\eta_1^4}^4, W_{2\eta_2^4}^4, W_{3\eta_3^4}^4\right)
\end{aligned}$$

is a special set fuzzy linear 4-subalgebra.

It is important to mention that the notion of special set fuzzy linear n-subalgebra and special set vector n-subspace are equivalent. However we can associate with a given special set linear n-subalgebra infinite number of special set fuzzy linear n-subalgebras. Secondly the notion of special set fuzzy linear n-subalgebras and special set fuzzy vector spaces are fuzzy equivalent.

    So even if the models are different in reality by fuzzifying they can be made identical or equivalent.

**Chapter Three**

# SPECIAL SEMIGROUP SET VECTOR SPACES AND THEIR GENERALIZATION

This chapter has two sections. First section just recalls some of the basic definitions from [60]. Section two defines the notion of special semigroup set vector semigroup set vector spaces and generalizes them.

## 3.1 Introduction to Semigroup Vector Spaces

In this section we just recall the basic definitions essential to make this chapter a self contained one.

**DEFINITION 3.1.1:** *Let V be a set, S any additive semigroup with 0. We call V to be a semigroup vector space over S if the following conditions hold good.*
1. *$sv \in V$ for all $s \in S$ and $v \in V$.*
2. *$0 \cdot v = 0 \in V$ for all $v \in V$ and $0 \in S$; 0 is the zero vector.*
3. *$(s_1 + s_2) v = s_1 v + s_2 v$*

*for all $s_1, s_2 \in S$ and $v \in V$.*

We illustrate this by the following examples.

***Examples 3.1.1:*** Let V = ($Z^+ \cup \{0\}$) × $2Z^+ \cup \{0\}$ × ($3Z^+ \cup \{0\}$) be a set and S = $Z^+ \cup \{0\}$ be a semigroup under addition. V is a semigroup vector space over S.

***Example 3.1.2:*** Let

$$V = \left\{ \begin{pmatrix} a_1 & a_2 & a_3 & a_4 \\ a_5 & a_6 & a_7 & a_8 \end{pmatrix} \middle| a_i \in Z^+ \cup \{0\}; 1 \le i \le 8 \right\}$$

be a set. Suppose S = $2Z^+ \cup \{0\}$ be a semigroup under addition. V is a semigroup vector space over S.

***Example 3.1.3:*** Let V = $3Z^+ \cup \{0\}$ be a set and S = $Z^+ \cup \{0\}$ be a semigroup under addition. V is a semigroup vector space over S.

**DEFINITION 3.1.2:** *Let V be semigroup vector space over the semigroup S. A set of vectors $\{v_1, …, v_n\}$ in V is said to be a semigroup linearly independent set if*

*(i) $v_i \ne sv_j$*

*for any $s \in S$ for $i \ne j$; $1 \le i, j \le n$.*

**DEFINITION 3.1.3:** *Let V be a semigroup vector space over the semigroup S under addition. Let T = $\{v_1, …, v_n\} \subseteq V$ be a subset of V we say T generates the semigroup vector space V over S if every element v of V can be got as $v = sv_i$, $v_i \in T$; $s \in S$.*

**DEFINITION 3.1.4:** *Let V be a semigroup vector space over the semigroup S. Suppose P is a proper subset of V and P is also a semigroup vector space over the semigroup S, then we call P to be semigroup subvector space of V.*

***Example 3.1.4:*** Let V = $\{Q^+ \cup \{0\}\}$ be a semigroup vector space over the semigroup S = $Z^+ \cup \{0\}$. Take W = $2Z^+ \cup \{0\}$, W is a semigroup subvector space of V over S.

***Example 3.1.5:*** Take V = {(1 1 0 1 0), (0 0 0 0 0), (1 1 1 1 1), (1 0 1 0 1), (0 1 1 1 0), (0 0 0), (0 1 1), (0 1 0), (1 0 1)} be a set. V is a semigroup vector space over the semigroup S = {0, 1} under addition 1 + 1 = 1.

Take W = {(1 0 1), (0 1 1), (0 0 0), (0 1 0)} ⊆ V; W is a semigroup subvector space over the semigroup S = {0, 1}. In fact every subset of V with the two elements (0 0 0) and or (0 0 0 0 0) is a semigroup subvector space over the semigroup S = {0, 1}.

**DEFINITION 3.1.5:** *Let V be a semigroup vector space over the semigroup S. If V is itself a semigroup under '+' then we call V to be a semigroup linear algebra over the semigroup S, if $s(v_1 + v_2) = sv_1 + sv_2$ ; $v_1, v_2 \in V$ and $s \in S$.*

**DEFINITION 3.1.6:** *Let V be a semigroup linear algebra over the semigroup S. Suppose P is a proper subset of V and P is a subsemigroup of V. Further if P is a semigroup linear algebra over the same semigroup S then we call P a semigroup linear subalgebra of V over S.*

**DEFINITION 3.1.7:** *Let V be a semigroup vector space over the semigroup S. Let P ⊂ V be a proper subset of V and T a subsemigroup of S. If P is a semigroup vector space over T then we call P to be a subsemigroup subvector space over T.*

**DEFINITION 3.1.8:** *Let V be a semigroup linear algebra over the semigroup S. Let P ⊂ V be a proper subset of V which is a subsemigroup under '+'. Let T be a subsemigroup of S. If P is a semigroup linear algebra over the semigroup T then we call P to be a subsemigroup linear subalgebra over the subsemigroup T.*

**DEFINITION 3.1.9:** *Let V be a semigroup linear algebra over the semigroup S. If V has no subsemigroup linear subalgebras over any subsemigroup of S then we call V to be a pseudo simple semigroup linear algebra.*

**DEFINITION 3.1.10:** *Let V be a semigroup linear algebra over a semigroup S. Suppose V has a proper subset P which is only a semigroup vector space over the semigroup S and not a semigroup linear algebra then we call P to be the pseudo semigroup subvector space over S.*

**DEFINITION 3.1.11:** *Let V be a semigroup linear algebra over the semigroup S. Let P be a proper subset of V and P is not a semigroup under the operations of V. Suppose $T \subseteq S$, a proper subset of S and T is also a semigroup under the same operations of S; i.e., T a subsemigroup of S, then we call P to be a pseudo subsemigroup subvector space over T if P is a semigroup vector space over T.*

**DEFINITION 3.1.12:** *Let V be a semigroup linear algebra over the semigroup S. If V has no subsemigroup linear algebras over subsemigroups of S then we call V to be a simple semigroup linear algebra.*

**DEFINITION 3.1.13:** *Let V be a semigroup under addition and S a semigroup such that V is a semigroup linear algebra over the semigroup S. If V has no proper subset P ($\subseteq V$) such that V is a pseudo subsemigroup vector subspace over a subsemigroup, T of S then we call V to be a pseudo simple semigroup linear algebra.*

**DEFINITION 3.1.14:** *Let V and W be any two semigroup linear algebras defined over the same semigroup, S we say T from V to W is a semigroup linear transformation if $T(c\alpha + \beta) = cT(\alpha) + T(\beta)$ for all $c \in S$ and $\alpha, \beta \in V$.*

**DEFINITIONS 3.1.15:** *Let V be a semigroup linear algebra over the semigroup S. A map T from V to V is said to be a semigroup linear operator on V if $T(cu + v) = cT(u) + T(v)$ for every $c \in S$ and $u, v \in V$.*

**DEFINITION 3.1.16:** *Let V be a semigroup vector space over the semigroup S. Let $W_1, …, W_n$ be a semigroup subvector spaces of*

V over the semigroup S. If $V = \bigcup_{i=1}^{n} W_i$ but $W_i \cap W_j \neq \phi$ or $\{0\}$ if $i \neq j$ then we call V to be the pseudo direct union of semigroup vector spaces of V over the semigroup S.

**DEFINITION 3.1.17:** *Let V be a semigroup linear algebra over the semigroup S. We say V is a direct sum of semigroup linear subalgebras $W_1, \ldots, W_n$ of V if*
1. $V = W_1 + \ldots + W_n$
2. $W_i \cap W_j = \{0\}$ or $\phi$ if $i \neq j$ $(1 \leq i, j \leq n)$.

**DEFINITION 3.1.18:** *Let V be a set with zero, which is non empty. Let G be a group under addition. We call V to be a group vector space over G if the following condition are true.*

1. *For every $v \in V$ and $g \in G$ $gv$ and $vg \in V$.*
2. *$0.v = 0$ for every $v \in V$, 0 the additive identify of G.*

We illustrate this by the following examples.

*Example 3.1.6:* Let $V = \{0, 1, 2, \ldots, 15\}$ integers modulo 15. $G = \{0, 5, 10\}$ group under addition modulo 15. Clearly V is a group vector space over G, for $gv \equiv v_1$ (mod 15), for $g \in G$ and $v, v_1 \in V$.

*Example 3.1.7:* Let $V = \{0, 2, 4, \ldots, 10\}$ integers 12. Take $G = \{0, 6\}$, G is a group under addition modulo 12. V is a group vector space over G, for $gv \equiv v_1$ (mod 12) for $g \in G$ and $v, v_1 \in V$.

*Example 3.1.8:* Let

$$M_{2 \times 3} = \left\{ \begin{pmatrix} a_1 & a_2 & a_3 \\ a_4 & a_5 & a_6 \end{pmatrix} \middle| a_i \in \{-\infty, \ldots, -4, -2, 0, 2, 4, \ldots, \infty\} \right\}.$$

Take $G = Z$ be the group under addition. $M_{2 \times 3}$ is a group vector space over $G = Z$.

**DEFINITION 3.1.19:** *Let V be the set which is a group vector space over the group G. Let P $\subseteq$ V be a proper subset of V. We say P is a group vector subspace of V if P is itself a group vector space over G.*

**DEFINITION 3.1.20:** *Let V be a group vector space over the group G. We say a proper subset P of V to be a linearly dependent subset of V if for any $p_1, p_2 \in P$, $(p_1 \neq p_2)$ $p_1 = ap_2$ or $p_2 = a'p_1$ for some $a, a' \in G$.*

*If for no distinct pair of elements $p_1, p_2 \in P$ we have $a, a_1 \in G$ such that $p_1 = ap_2$ or $p_2 = a_1 p_1$ then we say the set P is a linearly independent set.*

**DEFINITION 3.1.21:** *Let V be a group vector space over the group G.*

*Suppose T is a subset of V which is linearly independent and if T generates V i.e., using $t \in T$ and $g \in V$ we get every $v \in V$ as $v = gt$ for some $g \in G$ then we call T to be the generating subset of V over G.*

*The number of elements in V gives the dimension of V. If T is of finite cardinality V is said to be finite dimensional otherwise V is said to be of infinite dimension.*

**DEFINITION 3.1.22:** *Let V be a group vector space over the group G.*

*Let W $\subseteq$ V be a proper subset of V. H $\subset$ G be a proper subgroup of G. If W is a group vector space over H and not over G then we call W to be a subgroup vector subspace of V.*

**DEFINITION 3.1.23:** *Let V be a group vector space over the group G. Suppose W $\subset$ V is a subset of V. Let S be a subset of G. If W is a set vector space over S then we call W to be a pseudo set vector subspace of the group vector space.*

**DEFINITION 3.1.24:** *Let V be a group linear algebra over the group G. Suppose $W_1, W_2, …, W_n$ are distinct group linear subalgebras of V. We say V is a pseudo direct sum if*

1. $W_1 + ... + W_n = V$
2. $W_i \cap W_j \neq \{0\}$, even if $i \neq j$
3. We need $W_i$'s to be distinct i.e., $W_i \cap W_j \neq W_i$ or $W_j$ if $i \neq j$.

For more please refer [60].

## 3.2 Special Semigroup Set Vector Spaces and Special Group Set Vector Spaces

Now in this section we proceed on to define yet another new type of vector spaces called special semigroup set vector space over sets.

**DEFINITION 3.2.1:** *Let $V = (S_1, S_2, ..., S_n)$ be a set of collection of semigroups. Suppose P is any nonempty set such that for every $p \in P$ and $s_i \in S_i$, $ps_i \in S_i$ true for each $i = 1, 2, ..., n$. Then we call V to be a special semigroup set vector space over the set P.*

We now illustrate this by the following example.

*Example 3.2.1:* Let $V = \{S_1, S_2, S_3, S_4\}$, where

$$S_1 = \left\{ \begin{pmatrix} a & b \\ c & d \end{pmatrix} \middle| a,b,c,d \in Z_2 = \{0,1\} \right\}$$

is a semigroup under addition modulo 2, $S_2 = Z_2 \times Z_2 \times Z_2 \times Z_2$ is a semigroup under component wise addition, $S_3 = \{Z_2[x];$ all polynomials of degree less than or equal to 3$\}$ is a semigroup under polynomial addition and

$$S_4 = \left\{ \begin{pmatrix} a & b \\ c & d \\ e & f \\ g & h \end{pmatrix} \middle| a,b,c,d,e,f,g,h \in Z_2 = \{0,1\} \right\}$$

is a semigroups under matrix addition modulo 2.

Then $V = \{S_1, S_2, S_3, S_4\}$ is a special semigroup set vector space over the set $P = \{0, 1\}$.

***Example 3.2.2:*** Let $V = \{S_1, S_2, S_3\}$ where

$$S_1 = Z^+ \times Z^+ \times Z^+, \; S_2 = \left\{ \begin{pmatrix} a_1 & a_2 & a_3 \\ a_4 & a_5 & a_6 \\ a_7 & a_8 & a_9 \end{pmatrix} \middle| a_i \in Z^+; 1 \leq i \leq 9 \right\}$$

and

$S_3 = \{Z^+[x]$; all polynomials of degree less than or equal to 5 with coefficients from $Z^+\}$;

clearly $S_1$, $S_2$ and $S_3$ are semigroups under addition. Thus $V = \{S_1, S_2, S_3\}$ is a special semigroup set vector space over the set $P \subseteq Z^+$.

***Example 3.2.3:*** Let $V = \{S_1, S_2, S_3, S_4, S_5\}$ where $S_1 = \{Z_6 \times Z_6\}$ is a semigroup under addition modulo 6,

$$S_2 = \left\{ \begin{pmatrix} a & b \\ c & d \end{pmatrix} \middle| a, b, c, d \in Z_6 \right\}$$

and

$$S_3 = \left\{ \begin{pmatrix} x & x \\ y & y \\ z & z \\ w & w \end{pmatrix} \middle| x, y, z, w \in Z_6 \right\}$$

are semigroups under matrix addition modulo 6.

$$S_4 = \left\{ \begin{pmatrix} a & a & a & a \\ a & a & a & a \end{pmatrix} \middle| a \in Z_6 \right\}$$

and

$$S_5 = \left\{ \begin{pmatrix} a & a & a & a & a \\ 0 & a & a & a & a \\ 0 & 0 & a & a & a \\ 0 & 0 & 0 & a & a \\ 0 & 0 & 0 & 0 & a \end{pmatrix} \middle| a \in Z_6 \right\}$$

are semigroups under addition modulo 6. $V = \{S_1, S_2, S_3, S_4, S_5\}$ is a special semigroup set vector space over any subset of $Z_6$.

Now we proceed to define the notion of special semigroup set vector subspace of a special semigroup set vector space over a set P.

**DEFINITION 3.2.2:** *Let $V = \{S_1, S_2, …, S_n\}$ be a special semigroup set vector space over the set P, that is each $S_i$ is a semigroup and for each $s_i \in S_i$ and $p \in P$, $ps_i \in S_i$; $1 \leq i \leq n$. Let $W = \{T_1, …, T_n\}$ where each $T_i$ is a proper subsemigroup of $S_i$ for $i = 1, 2, …, n$. If W is a special semigroup set vector space over the set P, then we call W to be the special semigroup set vector subspace of V over P.*

Now we illustrate this by a few examples.

*Example 3.2.4:* Let $V = \{S_1, S_2, S_3, S_4, S_5\}$ where $S_1 = \{Z_5 \times Z_5 \times Z_5 \times Z_5\}$, a semigroup under component wise addition modulo 5,

$$S_2 = \left\{ \begin{pmatrix} a & a \\ a & b \end{pmatrix} \middle| a, b \in Z_5 \right\}$$

a semigroup under matrix addition modulo 5, $S_3 = \{Z_5[x];$ all polynomials of degree less than or equal to six$\}$, is a semigroup under polynomial addition,

$$S_4 = \left\{ \begin{pmatrix} a & a & a_1 & a_1 \\ a & a & a_1 & a_1 \\ a & a & a_1 & a_1 \end{pmatrix} \middle| a_1, a \in Z_5 \right\}$$

is a semigroup under matrix addition modulo 5 and

$$S_5 = \left\{ \begin{pmatrix} a & 0 & 0 & 0 \\ b & e & 0 & 0 \\ c & f & i & 0 \\ d & g & j & k \end{pmatrix} \middle| a,b,c,d,e,f,g,h,i,j,k \in Z_5 \right\}$$

is again a semigroup under matrix addition modulo 5. Thus V is a special semigroup set vector space over the set P = {0, 1, 2, 3, 4} = $Z_5$. Take W = ($W_1$, $W_2$, $W_3$, $W_4$, $W_5$) where $W_1$ = {(a a a a) / a ∈ $Z_5$} ⊆ $V_1$,

$$W_2 = \left\{ \begin{pmatrix} a & a \\ a & a \end{pmatrix} \middle| a \in Z_5 \right\} \subseteq V_2$$

a subsemigroup of $V_2$ ; $W_3$ = {$Z_5$[x]; all polynomials of even degree with coefficients from $Z_5$} is a subsemigroup of $V_3$,

$$W_4 = \left\{ \begin{pmatrix} a & a & a & a \\ a & a & a & a \end{pmatrix} \middle| a \in Z_5 \right\} \subseteq V_4;$$

a subsemigroup of $V_4$ and

$$W_5 = \left\{ \begin{pmatrix} a & 0 & 0 & 0 \\ a & a & 0 & 0 \\ a & a & a & 0 \\ a & a & a & a \end{pmatrix} \middle| a \in Z_5 \right\} \subseteq V_5;$$

a subsemigroup of $V_5$. Clearly W = ($W_1$, $W_2$, $W_3$, $W_4$, $W_5$) is a special semigroup set subvector space over the set P = {0, 1, 2, 3, 4} ⊆ $Z_5$.

Now we proceed on to give yet another example.

***Example 3.2.5:*** Let V = ($V_1$, $V_2$, $V_3$, $V_4$) where

$$V_1 = \left\{ \begin{pmatrix} a & b & c \\ d & e & f \\ g & h & i \end{pmatrix} \middle| a,b,c,d,e,f,g,h,i \in Z^+ \right\}$$

is a semigroup under matrix addition; $V_2 = \{Z^+ \times Z^+ \times Z^+\}$ a semigroup under component wise addition, $V_3 = \{Z^+[x]$; all polynomials in the variable x with coefficients from $Z^+\}$ is a semigroup under polynomial addition and

$$V_4 = \left\{ \begin{pmatrix} a_1 & 0 & 0 & 0 \\ 0 & a_2 & 0 & 0 \\ 0 & 0 & a_3 & 0 \\ 0 & 0 & 0 & a_4 \end{pmatrix} \middle| a_i \in Z^+; 1 \leq i \leq 4 \right\}$$

is again a semigroup under matrix addition. Clearly $V = (V_1, V_2, V_3, V_4)$ is a special semigroup set vector space over the set $Z^+$. Take $W = \{W_1, W_2, W_3, W_4\}$ where

$$W_1 = \left\{ \begin{pmatrix} a & a & a \\ a & a & a \\ a & a & a \end{pmatrix} \middle| a \in Z^+ \right\}$$

is a subsemigroup of $V_1$ under matrix addition, $W_2 = (3Z^+ \times 3Z^+ \times 5Z^+)$ is a subsemigroup of $Z^+ \times Z^+ \times Z^+$, $W_3 = \{$all polynomials of even degree in x with coefficients from $Z^+\}$ is a subsemigroup under polynomial addition of $V_3$ and

$$W_4 = \left\{ \begin{pmatrix} a & 0 & 0 & 0 \\ 0 & a & 0 & 0 \\ 0 & 0 & b & 0 \\ 0 & 0 & 0 & b \end{pmatrix} \middle| a,b \in 2Z^+ \right\}$$

is again a subsemigroup of $V_4$. Thus $W = \{W_1, W_2, W_3, W_4\} \subseteq V$ is a special semigroup set vector subspace of V.

Now we define the notion of special semigroup set linear algebra over a semigroup P.

**DEFINITION 3.2.3:** *Let $V = (V_1, V_2, ..., V_n)$ be a special semigroup set vector space over the set P. If P is an additive semigroup then we call V to be a special semigroup set linear algebra over the semigroup P.*

It is important to note that all special semigroup set linear algebras are special semigroup set vector spaces but special semigroup set vector spaces in general are not special semigroup set linear algebras.

We see by the very definitions of special semigroup set vector spaces and special semigroup set linear algebras all special semigroup set linear algebras are special semigroup set vector spaces; as every semigroup can be realized also as a set.

On the contrary a special semigroup set vector space in general is not a special semigroup set linear algebra. This is proved by an example.

Take $V = (V_1, V_2, V_3)$ where $V_1 = \{Z^+ \times Z^+\}$, $V_2 = \{Z^+[x]\}$ and

$$V_3 = \left\{ \begin{pmatrix} a & b \\ c & d \end{pmatrix} \bigg| a,b,c,d \in Z^+ \right\};$$

clearly $V_1$, $V_2$ and $V_3$ are semigroups under addition. Suppose $P = \{1, 2, 5, 7, 3, 8, 12, 15\}$ then V is a special semigroup set vector space over the set P. Clearly P is not a semigroup so V cannot be a special semigroup set linear algebra.

Now we give some examples of special semigroup set linear algebras.

***Example 3.2.6:*** Let $V = \{V_1, V_2, V_3, V_4\}$ be such that $V_1 = \{Z_{10} \times Z_{10}\}$, $V_2 = \{Z_{10}[x]$ where x an indeterminate, and $Z_{10}[x]$ contains all polynomials in x with coefficients from $Z_{10}\}$,

$$V_3 = \left\{ \begin{pmatrix} a & a \\ a & a \\ a & a \\ a & a \end{pmatrix} \middle| a \in Z_{10} \right\},$$

$$V_4 = \left\{ \begin{pmatrix} a & 0 & 0 & 0 \\ b & e & 0 & 0 \\ c & f & g & 0 \\ d & g & h & i \end{pmatrix} \middle| a,b,e,c,f,g,d,h,f,i \in Z_{10} \right\}$$

and

$$V_5 = \left\{ \begin{pmatrix} a & b & c & d & e \\ f & g & h & i & j \end{pmatrix} \middle| a,b,c,d,e,f,g,h,i,j \in Z_{10} \right\}$$

are semigroups under addition. Clearly if we take $P = Z_{10}$, P is also a semigroup under addition modulo 10. V is a special semigroup set linear algebra over the semigroup $Z_{10}$.

*Example 3.2.7:* Let $V = (V_1, V_2, V_3, V_4)$ where

$$V_1 = \{Z^+ \times Z^+ \times Z^+\},$$

$$V_2 = \left\{ \begin{pmatrix} a & b \\ 0 & d \end{pmatrix} \middle| a,c,d \in Z^+ \right\},$$

$$V_3 = \left\{ \begin{pmatrix} 0 & a & -b & c \\ -a & 0 & d & -e \\ b & -d & 0 & f \\ -c & e & -f & 0 \end{pmatrix} \middle| a,b,c,d,e,f \in Z^+ \right\}$$

and

$$V_4 = \left\{ \begin{pmatrix} a & b & c & d \\ b & e & f & g \\ c & f & h & i \\ d & g & i & j \end{pmatrix} \middle| a,b,c,d,e,f,g,h,i \text{ and } j \text{ are in } 2Z^+ \right\}.$$

Clearly $V_1$, $V_2$, $V_3$ and $V_4$ are semigroup under addition and V is a special semigroup set linear algebra over the semigroup $Z^+$.

***Example 3.2.8:*** Let $V = (V_1, V_2, V_3, V_4, V_5)$ where $V_1 = Z_5$, $V_2 = Z_7$, $V_3 = Z_3$, $V_4 = Z_{12}$ and $V_5 = Z_{13}$ are semigroups under addition modulo the appropriate n, n ∈ {5, 7, 3, 12, 13}. V is a special semigroup set vector space over the set $S = \{0, 1\}$.

***Remark:*** When we say $V = (V_1, \ldots, V_n)$ is a special semigroup set linear algebra over the semigroup T, we unlike in a linear algebra demand only two things

1. For every $p_i \in V_i$ and $t_i \in T$, $t_i p_i \in V_i$. $1 \le i \le n$
2. $V_i$ is a semigroup
3. T is a semigroup
4. We do not have $(a + b)p_i = ap_i + bp_i$ in general; $1 \le i \le n$;

that is this sort of distribution laws are never assumed in case of special semigroup set linear algebras over the semigroup T.

Now we proceed on to define the notion of special semigroup linear subalgebra and special subsemigroup linear subalgebras. Before we go for these two definitions we proceed onto show how this new structure varies from other structures.

***Example 3.2.9:*** Suppose $V = (S_1, S_2) = (Z_{12}, Z_7)$ be the semigroup special linear algebra over the semigroup $P = Z_6$. Then $P = \{0, 1, 2, \ldots, 5\}$, semigroup under addition modulo 6. For $7 \in Z_{12}$ and $5 \in Z_6$ $5.7 \equiv 11 \pmod{12}$. For $4 \in Z_6$ and $6 \in Z_7$ $4.6 = 24 \equiv 3 \pmod 7$ $4.6 \equiv 0 \pmod{12}$ if $6 \in Z_{12}$ and $4 \in Z_6$. So as an element of one semigroup its impact on the same element from $Z_6$ is distinctly different.

Thus when we have two different sets (semigroups) they act very differently on the same element from the semigroup over which they are defined, or in truth even the elements may represent the same property but its acting on the same element from the semigroup yields different values.

This is very evident from the example we give yet another example before we proceed to work with other properties.

***Example 3.2.10:*** Let $V = (V_1, V_2, V_3)$ where $V_1 = Z_{12}$, $V_2 = Z_9$ and $V_3 = Z_8$. Suppose V is a special semigroup set linear algebra over the semigroup $P = Z_6$. Take $6 \in Z_{12}$, $6 \in Z_9$ and $6 \in Z_8$. Take $4 \in P = Z_6$. Now $4.6 \equiv 0 \pmod{12}$ $4.6 \equiv 6 \pmod 9$ and $4.6 \equiv 0 \pmod 8$. So for these two semigroups $V_1$ and $V_3$, $4 \in P$ acts as an annulling element were as for $V_2$ it gives back the same element.

This special type of properties are satisfied by many real models, it may be useful in industries or experiments in such cases they can use these special algebraic structures.

**DEFINITION 3.2.4:** *Let $V = (V_1, V_2, ..., V_n)$ be a special semigroup set linear algebra over the semigroup S. Let $W = (W_1, W_2, ..., W_n) \subseteq (V_1, V_2, ..., V_n)$ ($W_i \subseteq V_i$; $W_i$ subsemigroup of $V_i$, $1 \le i \le n$). If W is itself a special semigroup set linear algebra over the semigroup S then we call W to be a special semigroup set linear subalgebra of V over the semigroup S.*

We illustrate this by some examples.

***Example 3.2.11:*** Let $V = (V_1, V_2, V_3, V_4)$ where $V_1 = Z_8$, $V_2 = Z_6$, $V_3 = Z_9$ and $V_3 = Z_{12}$. V is a semigroup set linear algebra over the semigroup $P = Z_4$. Take $W = (W_1, W_2, W_3, W_4)$ where $W_1 = \{0, 2, 4, 6\} \subseteq Z_8$, $W_2 = \{0, 2, 4\} \subseteq Z_6$, $W_3 = \{0, 3, 6\} \subseteq Z_9$ and $W_4 = \{0, 3, 6, 9\} \subseteq Z_{12}$. It is easily verified $W = (W_1, W_2, W_3, W_4)$ is a special semigroup set linear subalgebra over the semigroup $P = Z_4$.

*Example 3.2.12:* Let $V = (V_1, V_2, V_3, V_4, V_5)$ where $V_1 = Z_7$, $V_2 = Z_5$, $V_3 = Z_6$, $V_4 = Z_8$ and $V_5 = Z_4$. V is a special semigroup set linear algebra over the semigroup $P = Z_2 = \{0, 1\}$. We see V has no proper subset W which can be a special semigroup set linear algebra over the semigroup P. For $V_1$ and $V_2$ have no proper subsemigroups.

In view of this we define the notion of special semigroup set simple linear algebra.

**DEFINITION 3.2.5:** *Let $V = (V_1, V_2, \ldots, V_n)$ be a special semigroup set linear algebra over the semigroup P. If V has no proper subset $W = (W_1, W_2, \ldots, W_n)$ such that $W_i$ is a proper subsemigroup of $V_i$ for atleast for one i, $1 \leq i \leq n$ then we call V to be a special semigroup set simple linear algebra. In a special set simple linear algebra V we do not have proper special semigroup set linear subalgebras.*

We illustrate this situation by some examples.

*Example 3.2.13:* Let $V = (V_1, V_2, V_3, V_4, V_5)$ where $V_1 = Z_3$, $V_2 = Z_7$, $V_3 = Z_5$, $V_4 = Z_6$ and $V_5 = Z_{11}$ be a special semigroup set linear algebra over the semigroup $P = Z_8$. We see $V_1$, $V_2$, $V_3$ and $V_5$ has no proper subsemigroups. So V is only a special semigroup set simple linear algebra.

*Example 3.2.14:* Let $V = (V_1, V_2, V_3, V_4)$ where $V_1 = Z_7$, $V_2 = Z_{13}$, $V_3 = Z_{19}$ and $V_4 = Z_{11}$; V is a special semigroup set linear algebra over the semigroup $P = Z_6$. We see $V_1$, $V_2$, $V_3$ and $V_4$ do not have any proper subsemigroups, i.e., none of the semigroups $V_1$, $V_2$, $V_3$ and $V_4$ have subsemigroups. Thus $V = (V_1, V_2, V_3, V_4)$ is a special semigroup set simple linear algebra over $Z_6$.

In view of this we have the following definition.

**DEFINITION 3.2.6:** *Let $V = (V_1, V_2, \ldots, V_n)$ be a special semigroup set linear algebra over the semigroup P. In none of the semigroups $V_i$, $1 \leq i \leq n$ has proper subsemigroups then we* *call V to be a special semigroup set strong simple linear algebra.*

**Example 3.2.15:** *Let $V = (V_1, V_2, V_3, V_4)$ where $V_1 = Z_5$, $V_2 = Z_{23}$, $V_3 = Z_{11}$ and $V_4 = Z_{13}$ be the special semigroup set linear algebra over the semigroup $P = Z_4$. It is easily verified V is a special semigroup set strong simple linear algebra.*

We have the following interesting theorem.

**THEOREM 3.2.1:** *Let $V = (V_1, V_2, ..., V_n)$ be a special semigroup set strong simple linear algebra then V is a special semigroup set simple linear algebra. But a special semigroup set simple linear algebra need not in general be a special semigroup set strong simple linear algebra.*

*Proof:* Let $V = (V_1, V_2, ..., V_n)$ be a special semigroup set strong simple linear algebra over the semigroup P. This implies every semigroup $V_i$ of V has no proper subsemigroup, $1 \leq i \leq n$. So V is also a special semigroup set simple linear algebra.

We prove the converse by a simple example.
Let $V = (V_1, V_2, V_3, V_4, V_5)$ where $V_1 = Z_7$, $V_2 = Z_8$, $V_3 = Z_9$, $V_4 = Z_{11}$ and $V_5 = Z_{13}$. V is a special semigroup set linear algebra over the semigroup $Z_6$. Clearly the semigroups $V_1$, $V_4$ and $V_5$ have no subsemigroups but $V_2$ and $V_3$ have subsemigroups. Thus V is only a special semigroup set simple linear algebra over the semigroup $Z_6$ and V is not a special semigroup set strongly simple linear algebra over the semigroup $Z_6$. Hence the theorem.

Now in case of special semigroup set vector spaces also we have these concepts which is defined briefly.

**DEFINITION 3.2.7:** *Let $V = (V_1, V_2, ..., V_n)$ be a special semigroup set vector space over the set P, where $V_1, V_2, ..., V_n$ are semigroups under addition. If V does not contain a $W = (W_1, ..., W_n)$ such that $W_i$'s are proper subsemigroups of $V_i$ for atleast some i, $1 \leq i \leq n$, then we define V to be a special* *semigroup set simple vector space. If none of the semigroups $V_i$ in V has subsemigroups, $1 \leq i \leq n$ then we say the special semigroup set vector space V is a special set semigroup strong simple vector space over P.*

Now we illustrate this by some examples.

***Example 3.2.16:*** Let $V = (V_1, V_2, V_3, V_4, V_5)$ where $V_1 = Z_8$, $V_2 = Z_7$, $V_3 = Z_5$, $V_4 = Z_6$ and $V_5 = Z_{11}$ be a special semigroup set vector space over the set $P = \{0, 1, 2, 3, 4\}$. V is only a special semigroup set simple vector space over P, for $V_2$, $V_3$ and $V_5$ have no subsemigroups.

***Example 3.2.17:*** Let $V = (V_1, V_2, V_3, V_4)$ where $V_1 = Z_5$, $V_2 = Z_7$, $V_3 = Z_{13}$ and $V_4 = Z_{23}$ be a special semigroup set vector space over the set $P = \{0, 1, 2, 3, 4, 5, 6\}$. Clearly V is a special semigroup set strong simple vector space over the set $P = \{0, 1, 2, 3, 4, 5, 6\}$.

**DEFINITION 3.2.8:** *Let $V = (V_1, V_2, ..., V_n)$ be a special semigroup set linear algebra over the semigroup S. If $W = (W_1, W_2, ..., W_n)$ where $W_i \subseteq V_i$ is such that $W_i$ is a proper subsemigroup of the semigroup $V_i$, $1 \leq i \leq n$ and P be a proper subsemigroup of S. Suppose W is a special semigroup set linear algebra over P then we call W to be a special subsemigroup set linear subalgebra of V over P.*

We illustrate this by some simple examples.

***Example 3.2.18:*** Let $V = (V_1, V_2, V_3, V_4, V_5)$ where $V_1 = \{Z^+ \times Z^+ \times Z^+\}$, $V_2 = Z^+[x]$,

$$V_3 = \left\{ \begin{pmatrix} a & b \\ c & d \end{pmatrix} \middle| a,b,c,d \in Z^+ \right\},$$

$$V_4 = \left\{ \begin{pmatrix} a_1 & a_2 & a_3 \\ a_4 & a_5 & a_6 \end{pmatrix} \middle| a_1,a_2,a_3,a_4,a_5,a_6 \in Z^+ \right\}$$

and

$$V_5 = \left\{ \begin{pmatrix} a & b & c & d \\ b & t & e & f \\ c & e & g & h \\ d & f & h & k \end{pmatrix} \middle| a,b,c,d,e,f,g,h,k,t \in Z^+ \right\}$$

is a special semigroup linear algebra over the semigroup $S = Z^+$. Take $W = (W_1, W_2, W_3, W_4, W_5)$ where $W_1 = (2Z^+ \times 2Z^+ \times 2Z^+)$, $W_2 = \{$all polynomials of only degree with coefficients from $Z^+\}$,

$$W_3 = \left\{ \begin{pmatrix} a & a \\ a & a \end{pmatrix} \middle| a \in Z^+ \right\},$$

$$W_4 = \left\{ \begin{pmatrix} a & a & a \\ a & a & a \end{pmatrix} \middle| a \in Z^+ \right\}$$

and

$$W_5 = \left\{ \begin{pmatrix} a & a & a & a \\ a & a & a & a \\ a & a & a & a \\ a & a & a & a \end{pmatrix} \middle| a \in Z^+ \right\}.$$

Clearly $W = (W_1, W_2, W_3, W_4, W_5)$ is a special semigroup linear algebra over the semigroup $2Z^+$.

Thus $W = (W_1, W_2, W_3, W_4, W_5)$ is a special subsemigroup set linear subalgebra over the subsemigroup $P = 2Z^+ \subseteq S = Z^+$.

We give yet another example.

***Example 3.2.19:*** Let $V = (V_1, V_2, V_3, V_4)$ where $V_1 = Z_5$, $V_2 = Z_7$, $V_3 = Z_{11}$ and $V_4 = Z_{12}$, $V_1, V_2, V_3$ and $V_4$ are semigroups under addition modulo appropriate n, $n \in \{5, 7, 11, 12\}$. V is a special semigroup set linear algebra over the semigroup $S = Z_3$, semigroup under addition modulo 3. Now V has no proper subset $W = (W_1, W_2, W_3, W_4)$ such that each $W_i$ is a subsemigroup of $V_i$, $1 \leq i \leq 4$ and S has no proper subsemigroup.

So we see V has no special subsemigroup set linear subalgebra. This leads to the concept of a new algebra structure.

**DEFINITION 3.2.9:** *Let $V = (V_1, V_2, \ldots, V_n)$ be a special semigroup set linear algebra over the semigroup S. If V has no proper subset $W = (W_1, \ldots, W_n)$ such that each $W_i$ is a subsemigroup of the semigroup $V_i$; $i = 1, 2, \ldots, n$ and S has no proper subsemigroup then we call V to be a doubly simple special semigroup set linear algebra.*

We have the following interesting result.

**THEOREM 3.2.2:** *Let $V = (V_1, V_2, \ldots, V_n)$ be a special semigroup set linear algebra over the semigroup. If V is a doubly simple special semigroup set linear algebra then V is a special semigroup set strong simple linear algebra.*

*Proof:* Given $V = (V_1, V_2, \ldots, V_n)$ is a doubly simple special semigroup set linear algebra, so V does not contain a proper n-subset $W = (W_1, W_2, \ldots, W_n)$ such that each $W_i$ is a proper subsemigroup of the semigroup $V_i$; $(1 \leq i \leq n)$. So V has no special semigroup set linear subalgebra and each semigroup $V_i$ has no proper subsemigroup $W_i$, $1 \leq i \leq n$. Hence V is a special semigroup strong simple set linear algebra.

Now it is pertinent to mention here that we cannot say if V is a special semigroup strong simple set linear algebra then V is a doubly simple special semigroup linear algebra. For we may have a special semigroup strongly simple set linear algebra but it may not be a doubly simple special semigroup linear algebra.

For take $V = (V_1, V_2, V_3, V_4)$ where $V_1 = Z_5$, $V_2 = Z_7$, $V_3 = Z_{11}$ and $V_4 = Z_{19}$. V is clearly a special semigroup set strongly simple linear algebra over the semigroup $Z_4$, as V has no proper subset $W = (W_1, W_2, W_3, W_4)$ such that $W_i$ is a proper subsemigroup of $V_i$, $1 \leq i \leq 4$. For $Z_5$ has no subsemigroup, $Z_7$ has no subsemigroup, $Z_{11}$ cannot have a proper subsemigroup and $Z_{19}$ has no proper subsemigroup. But the semigroup V is defined over $Z_4$ and $Z_4$ has proper subset $P = \{0, 2\}$ which is a proper subsemigroup of $Z_4$ hence we cannot say $V = (V_1, V_2, V_3, V_4)$ is a doubly simple special semigroup set linear algebra.

Thus every doubly simple special semigroup set linear algebra is a special semigroup set linear algebra is a special semigroup set strong simple linear algebra.

Next we proceed on to define the generating n set of a special semigroup set vector space over a set S.

**DEFINITION 3.2.10** *Let $V = (V_1, ..., V_n)$ be a special semigroup set vector space over the set S. Let $X = (X_1, ..., X_n)$ be a proper subset of V i.e., each $X_i \subseteq V_i$, $1 \leq i \leq n$. If each $X_i$ generates $V_i$ over S then we say X is the generating n set of the special semigroup set vector space over the set S.*

We illustrate this by the following example.

*Example 3.2.20* Let $V = (V_1, V_2, V_3, V_4)$ where $V_1 = \{(a\ a\ a) \,/\, a \in Z^+\}$,

$$V_2 = \left\{ \begin{pmatrix} a & a \\ a & a \end{pmatrix} \middle| a \in Z^+ \right\},$$

$$V_3 = \left\{ \begin{pmatrix} a & a & a \\ a & a & a \end{pmatrix} \middle| a \in Z^+ \right\}$$

and $V_4 = \{a + ax + ax^2 + ax^3 \mid a \in Z^+\}$. Clearly $V_1$, $V_2$, $V_3$ and $V_4$ are semigroups under addition. $V = (V_1, V_2, V_3, V_4)$ is a special semigroup set vector space over the set $Z^+$. Take $X = \{X_1, X_2, X_3, X_4\}$ where $X_1 = (1\ 1\ 1)$,

$$X_2 = \begin{pmatrix} 1 & 1 \\ 1 & 1 \end{pmatrix}, X_3 = \begin{pmatrix} 1 & 1 & 1 \\ 1 & 1 & 1 \end{pmatrix},$$

$X_4 = \{1 + x + x^2 + x^3\} \subseteq V = (V_1, V_2, V_3, V_4)$, i.e., $X_i \subseteq V_i$, $i = 1, 2, 3, 4$. We see X is a 4-generating subset of V over the set $Z^+$.

*Example 3.2.21:* Let $V = (V_1, V_2, V_3, V_4, V_5)$ where $V_1 = \{Z_{12} \times Z_{12} \times Z_{12}\}$,

$$V_2 = \left\{ \begin{pmatrix} a & b & e \\ c & d & f \end{pmatrix} \middle| a,b,c,d,e,f \in Z_{12} \right\},$$

$$V_3 = \left\{ \begin{pmatrix} a & a \\ a & a \\ a & a \\ a & a \end{pmatrix} \middle| a \in Z_{12} \right\},$$

$V_4 = \{Z_{12}[x]$, all polynomials of degree less than or equal to 4$\}$
and

$$V_5 = \left\{ \begin{pmatrix} 0 & a & b & c \\ a & 0 & d & e \\ b & d & 0 & f \\ c & e & f & 0 \end{pmatrix} \middle| a = b = c = d = f = x \in Z_{12} \right\}$$

be a special semigroup set linear algebra over the semigroup $S = Z_{12}$. $X = (X_1, X_2, X_3, X_4, X_5) = \{(1\ 0\ 0\ 0), (0\ 1\ 0\ 0), (0\ 0\ 1\ 0), (0\ 0\ 0\ 1)\}$,

$$\left\{ \begin{pmatrix} 1 & 0 & 0 \\ 0 & 0 & 0 \end{pmatrix} \begin{pmatrix} 0 & 1 & 0 \\ 0 & 0 & 0 \end{pmatrix} \begin{pmatrix} 0 & 0 & 1 \\ 0 & 0 & 0 \end{pmatrix} \right.$$

$$\left. \begin{pmatrix} 0 & 0 & 0 \\ 1 & 0 & 0 \end{pmatrix} \begin{pmatrix} 0 & 0 & 0 \\ 0 & 1 & 0 \end{pmatrix} \begin{pmatrix} 0 & 0 & 0 \\ 0 & 0 & 1 \end{pmatrix} \right\}, \left\{ \begin{bmatrix} 1 & 1 \\ 1 & 1 \\ 1 & 1 \\ 1 & 1 \end{bmatrix} \right\}, \{1, x, x^2, x^3, x^4\},$$

$$\left\{ \begin{pmatrix} 0 & 1 & 1 & 1 \\ 1 & 0 & 1 & 1 \\ 1 & 1 & 0 & 1 \\ 1 & 1 & 1 & 0 \end{pmatrix} \right\} \subseteq (V_1, V_2, V_3, V_4, V_5)$$

is the 5-generating 5-subset of V over the semigroup $S = Z_{12}$. The dimension of V is (4, 6, 1, 5, 1).

Now we define special semigroup set linear transformation for special semigroup set vector spaces defined only on the same set S.

**DEFINITION 3.2.11:** *Let $V = (V_1, ..., V_n)$ and $W = (W_1, ..., W_n)$ be two special semigroup set vector spaces defined over the same set S. Let $T = (T_1, T_2, ..., T_n)$ be a n-map from V into W such that $T_i: V_i \to W_j$, $1 \leq i, j \leq n$. Clearly no two $V_i$'s can be mapped into the same $W_j$ that is each $V_i$ has a unique $W_j$ associated with it.*
  *If $T_i(\mu + v) = T_i(\mu) + T(v)$; $\mu, v \in V_i$, this is true for each i, $i = 1, 2$; $T = (T_1, T_2, T_3, T_4)$ then we call T to be special semigroup set linear transformation of V into W.*

We illustrate this by a simple example.

*Example 3.2.22:* Let $V = (V_1, V_2, V_3, V_4)$ and $W = (W_1, W_2, W_3, W_4)$ where

$$V_1 = \left\{ \begin{pmatrix} a & b \\ c & d \end{pmatrix} \middle| a, b, c, d \in Z^+ \right\},$$

$V_2 = \{Z^+ [x]$; all polynomials of degree less than or equal to 4$\}$,

$$V_3 = \left\{ \begin{pmatrix} a & a \\ a & a \\ a & a \\ a & a \end{pmatrix} \middle| a \in Z^+ \right\}$$

and $V_4 = Z^+ \times Z^+ \times Z^+ \times Z^+ \times Z^+ \times Z^+$ be a special semigroup set vector space over the set $S = Z^+$. $W_1 = \{Z^+ \times Z^+ \times Z^+ \times Z^+ \times Z^+\}$,

$$W_2 = \left\{ \begin{pmatrix} a & a \\ a & a \end{pmatrix} \middle| a \in Z^+ \right\},$$

$$W_3 = \left\{ \begin{pmatrix} a & b & f \\ c & d & g \end{pmatrix} \middle| a,b,c,d \in Z^+ \right\}$$

and

$$W_4 = \left\{ \begin{pmatrix} a & a & a & a \\ a & a & a & a \end{pmatrix} \middle| a \in Z^+ \right\}$$

be a special semigroup vector space over the same set $S = Z^+$.
Define $T = (T_1, T_2, T_3, T_4): V \to W$ as

$T_1: V_1 \to W_2$ defined by

$$T_1 \begin{pmatrix} a & b \\ c & d \end{pmatrix} = \begin{pmatrix} a & a \\ a & a \end{pmatrix}$$

$T_2 : V_2 \to W_1$ defined by

$$T_2(a_0 + a_1x + a_2x^2 + a_3x^3 + a_4x^4) = (a_0, a_1, a_2, a_3, a_4).$$

$T_3 : V_3 \to W_4$ defined by

$$T_3 \begin{pmatrix} a & a \\ a & a \\ a & a \\ a & a \end{pmatrix} = \begin{pmatrix} a & a & a & a \\ a & a & a & a \end{pmatrix}$$

and $T_4 : V_4 \to W_3$ defined by

$$T_4 (a_1\ a_2\ a_3\ a_4\ a_5\ a_6) = \begin{pmatrix} a_1 & a_2 & a_3 \\ a_4 & a_5 & a_6 \end{pmatrix}.$$

Clearly $T = (T_1, T_2, T_3, T_4)$ is a special semigroup set linear transformation from V to W.

Thus if one wishes to study using a special semigroup set vector space another special semigroup set vector space provided both of them are defined on the same set one by using the special semigroup set linear transformations study them.

Now we proceed on to define the notion of special semigroup set linear operator on a special semigroup set vector space V over a set S.

**DEFINITION 3.2.12:** *Let $V = (V_1, V_2, \ldots, V_n)$ denote the special semigroup set vector space V over the set S. Let $T = (T_1, T_2, \ldots, T_n)$ be a special semigroup set linear transformation from V to V, then we call T to be a special semigroup set linear operator on V. Thus if in a special semigroup set linear transformation W is replaced by V itself then we call that linear transformation as special semigroup set linear operator on V.*

We now illustrate it by some examples.

*Example 3.2.23:* Let $V = (V_1, \ldots, V_5)$ be a special semigroup set vector space over the set $S = Z_4$. Here $V_1 = Z_4 \times Z_4 \times Z_4 \times Z_4$, $V_2 = \{Z_4[x]$ all polynomials of degree less than or equal to three$\}$.

$$V_3 = \left\{ \begin{pmatrix} a & b \\ c & d \end{pmatrix} \middle| a,b,c,d \in Z_4 \right\},$$

$$V_4 = \left\{ \begin{pmatrix} a & 0 & 0 \\ b & 0 & 0 \\ c & 0 & d \end{pmatrix} \middle| a,b,c,d \in Z_4 \right\}$$

and

$$V_5 = \left\{ \begin{pmatrix} a & 0 & b & 0 & c \\ 0 & d & 0 & 0 & 0 \end{pmatrix} \middle| a,b,c,d \in Z_4 \right\}$$

are semigroups under appropriate addition modulo 4.

Let $T = (T_1, T_2, T_3, T_4, T_5)$ be a 5-map from V into V defined by

$$T_1 : V_1 \to V_2$$
$$T_2 : V_2 \to V_4$$
$$T_3 : V_3 \to V_1$$
$$T_4 : V_4 \to V_5$$

$$T_5 : V_5 \to V_3$$

where $T_1(a\ b\ c\ d) = a + bx + cx^2 + dx^3$, $T_2 : V_2 \to V_4$ is defined by

$$T_2(a+bx+cx^2+dx^3) = \begin{pmatrix} a & 0 & 0 \\ b & 0 & 0 \\ c & 0 & d \end{pmatrix},$$

$$T_3 \begin{pmatrix} a & b \\ c & d \end{pmatrix} = (a\ b\ c\ d),$$

$$T_4 \begin{pmatrix} a & 0 & 0 \\ b & 0 & 0 \\ c & 0 & d \end{pmatrix} = \begin{pmatrix} a & 0 & b & 0 & c \\ 0 & d & 0 & 0 & 0 \end{pmatrix}$$

and

$$T_5 \begin{pmatrix} a & 0 & b & 0 & c \\ 0 & d & 0 & 0 & 0 \end{pmatrix} = \begin{pmatrix} a & b \\ c & d \end{pmatrix}.$$

Thus T is a special semigroup set linear operator on V. How many such T's can be defined on V? For we have one set $SH_{Z_4}(V, V)$ given by ($\text{Hom}_{Z_4}(V_1, V_2)$, $\text{Hom}_{Z_4}(V_2, V_4)$, $\text{Hom}_{Z_4}(V_3, V_1)$, $\text{Hom}_{Z_4}(V_4, V_5)$, $\text{Hom}_{Z_4}(V_5, V_3)$) for $T = (T_1, T_2, T_3, T_4, T_5)$.

Thus we have many such T's yielding different sets of $SH_{Z_4}(V, V)$.

**DEFINITION 3.2.13:** *Let $V = (V_1, V_2, ..., V_n)$ be a special semigroup set linear algebra over the semigroup S. $W = (W_1, W_2, ..., W_n)$ be another special semigroup set linear algebra over the same semigroup S. If $T = (T_1, T_2, ..., T_n) : V \to W$ where $T_i : V_i \to W_j$ such that no two $V_i$'s are mapped into the same $W_j$ and $T_i(\mu + v) = T_i(\mu) + T_i(v)$ for all $\alpha \in S$ and $\mu, v \in V_i$, true for $i = 1, 2, ..., n$ and $1 \le j \le n$, then we call T to be a special semigroup set linear transformation of V to W.*

*Note:* If $\text{SHom}_s(V, W) = \{\text{Hom}_s(V_1, W_{i_1}), \text{Hom}_s(V_2, W_{i_2}), \ldots, \text{Hom}_s(V_n, W_{i_n})\}$, what is the algebraic structure of $\text{SHom}_s(V, W)\}$? Is $\text{SHom}_s(V, W)$ again a special semigroup set linear algebra over the same semigroup S?

If in the definition we replace W by V then we call $T = (T_1, T_2, \ldots, T_n)$ as the special semigroup linear operator of V over the semigroup S. If $\text{SHom}_s(V, V) = \{\text{Hom}_s(V_1, V_{i_1}), \text{Hom}_s(V_2, V_{i_2}) \ldots \text{Hom}_s(V_n, V_{i_n})\}$; is $\text{SHom}_s(V,V)$ a special semigroup set linear algebra over S?

This work is left as an exercise for the reader.

Now we illustrate these definitions by some examples.

*Example 3.2.24:* Let $V = (V_1, V_2, V_3, V_4)$ where $V_1 = \{Z_{10} \times Z_{10} \times Z_{10} \times Z_{10}\}$, $V_2 = \{Z_{10}[x]$, all polynomials of degree less than or equal to 3$\}$,

$$V_3 = \left\{ \begin{pmatrix} a & b & c \\ d & e & f \end{pmatrix} \middle| a,b,c,d,e,f \in Z_{10} \right\}$$

and

$$V_4 = \left\{ \begin{pmatrix} a & a & a \\ a & a & a \\ a & a & a \\ a & a & a \\ a & a & a \end{pmatrix} \middle| a \in Z_{10} \right\}.$$

Clearly $V_1, V_2, V_3$ and $V_4$ are semigroups under addition. V is a special semigroup set linear algebra over the semigroup $Z_{10}$.

Let $W = (W_1, W_2, W_3, W_4)$ where

$$W_1 = \left\{ \begin{pmatrix} a & a & a \\ a & a & a \end{pmatrix} \middle| a \in Z_{10} \right\},$$

$$W_2 = \left\{ \begin{pmatrix} a & a & a & a & a \\ a & a & a & a & a \\ a & a & a & a & a \end{pmatrix} \middle| a \in Z_{10} \right\}$$

$$W_3 = \left\{ \begin{pmatrix} a & b \\ c & d \end{pmatrix} \middle| a,b,c,d \in Z_{10} \right\}$$

and

$$W_4 = \left\{ \begin{pmatrix} a & 0 & b \\ 0 & c & 0 \\ 0 & 0 & d \end{pmatrix} \middle| a,b,c,d \in Z_{10} \right\}.$$

Clearly $W_1$, $W_2$, $W_3$ and $W_4$ are semigroups under addition. Thus $W = (W_1, W_2, W_3, W_4)$ is a special semigroup set linear algebra over the semigroup $Z_{10}$. Now let us define $T = (T_1, T_2, T_3, T_4)$ from V into W by

$T_1 : V_1 \to W_3$ given by

$$T_1 (a\ b\ c\ d) = \begin{pmatrix} a & b \\ c & d \end{pmatrix}.$$

$T_2 : V_2 \to W_4$ defined by

$$T_2(a_0 + a_1x + a_2x^2 + a_3x^3) = \begin{pmatrix} a_0 & 0 & a_1 \\ 0 & a_2 & 0 \\ 0 & 0 & a_3 \end{pmatrix}.$$

$T_3 : V_3 \to W_1$ defined by

$$T_3 \begin{pmatrix} a & b & c \\ d & e & f \end{pmatrix} = \begin{pmatrix} a & a & a \\ a & a & a \end{pmatrix}.$$

$T_4 : V_4 \to W_2$ defined by

$$T_4 \begin{pmatrix} a & a & a \\ a & a & a \\ a & a & a \\ a & a & a \\ a & a & a \end{pmatrix} = \begin{bmatrix} a & a & a & a & a \\ a & a & a & a & a \\ a & a & a & a & a \end{bmatrix}.$$

Clearly $T = (T_1, T_2, T_3, T_4)$ is a special semigroup linear transformation V to W. Suppose

$SHom_{Z_{10}}(V,W) = \{ Hom_{Z_{10}}(V_1,W_3), Hom_{Z_{10}}(V_2,W_4),$

$Hom_{Z_{10}}(V_3,W_1), Hom_{Z_{10}}(V_4,W_2)\}.$

Show $Hom_{Z_{10}}(V,W)$ is a special semigroup set vector space over $Z_{10}$. We give an example of a special semigroup set linear operator on V.

***Example 3.2.25:*** Let $V = (V_1, V_2, V_3, V_4, V_5)$ where $V_1 = Z_{20}$ a semigroup under addition modulo 20, $V_2 = Z_{12}$ a semigroup under addition modulo 12, $V_3 = \{Z_{10}$ a semigroup under addition modulo 10$\}$, $V_4 = Z_{23}$ a semigroup under addition modulo 23 and $V_5 = Z_8$ a semigroup under addition modulo 8. Now $V = (V_1, V_2, \ldots, V_5)$ is a special semigroup set linear algebra over the semigroup $S = Z_2$ under addition modulo 2.

Now find $SHom_{Z_2}(V,V) = \{ Hom_{Z_2}(V_1, V_3), Hom_{Z_2}(V_2, V_4), Hom_{Z_2}(V_3, V_1), Hom_{Z_2}(V_4, V_5), Hom_{Z_2}(V_5, V_2)\}$.

Now we proceed onto define yet another new notion called special group set vector space over the set S.

**DEFINITION 3.2.14:** *Let $V = (V_1, \ldots, V_n)$ where $V_1, \ldots, V_n$ are distinct additive abelian groups. Let S be any non empty set. If for every $x \in V_i$ and $a \in S$, $a x \in V_i$, $1 \leq i \leq n$ then we call V to be a special group set vector space over the set S.*

We illustrate this definition by some examples.

***Example 3.2.26:*** Let $V = (V_1, V_2, V_3, V_4, V_5)$ where $V_1 = Z_{10}$ group under addition modulo 10, $V_2 = Z_{26}$ group under addition modulo 26, $V_3 = Z_{15}$ group under a addition modulo 15, $V_4 = Z_{12}$ group under addition modulo 12 and $V_5 = Z_{14}$ group under addition modulo 14.

Take $S = \{0, 1, 2, 3, 4\}$ it is easily verified V is a special group set vector space over the set S. For we see if $4 \in S$ and $10 \in Z_{15} = V_3$, $4.10 = 40 \equiv 10 \pmod{15}$ so $4.10 \in V_3$; likewise for $14 \in Z_{26}$ and $3 \in S$; $3.14 = 42 \equiv 16 \pmod{26}$ so $3.14 \in Z_{26} = V_2$ and so on.

Now we can give yet another example of a special group set vector space over a set.

***Example 3.2.27:*** Let $V = \{V_1, V_2, V_3, V_4, V_5, V_6\}$ where $V_1 = 2Z$, $V_2 = 5Z$, $V_3 = 7Z$, $V_4 = 3Z$, $V_5 = 11Z$ and $V_6 = 11Z$.

Clearly $V_1, V_2, \ldots, V_6$ are groups under addition. Take $S = \{0\ 1\ 2\ 3\ 4\ 5\ 6\ 7\ 8\ 9\ 10\}$. Clearly $V = (V_1, \ldots, V_6)$ is a special group set vector space over the set S.

Now we can define their substructures.

**DEFINITION 3.2.15:** *Let $V = (V_1, V_2, V_3, \ldots, V_n)$ be n distinct abelian groups under addition. Suppose V be a special group set vector over the set S. Let $W = (W_1, W_2, \ldots, W_n)$ where each $W_i$ is a subgroup of the group $V_i$, $1 \le i \le n$. Suppose W is a special group set vector space over the same set S then we call W to be a special group set vector subspace of V over the set S.*

We now illustrate this situation by some simple examples.

***Example 3.2.28:*** Let $V = (V_1, V_2, V_3, V_4, V_5)$ where $V_1 = Z_{10}$, $V_2 = Z_{12}$, $V_3 = Z_{25}$, $V_4 = Z_9$ and $V_5 = Z_{21}$ all groups under modulo addition for appropriate n. ($n \in \{10, 12, 25, 9$ and $21\}$). V is a special group set vector space over the set $S = \{0, 1, 2, 3, 4\}$. Let $W = (W_1, W_2, \ldots, W_5)$ where $W_1 = \{0, 2, 4, 6, 8\} \subseteq Z_{10}$, $W_2 = \{0, 6\} \subseteq Z_{12}$, $W_3 = \{0, 5, 10, 15, 20\} \subseteq Z_5$, $W_4 = \{0, 3, 6\} \subseteq Z_9$ and $W_5 = \{0, 7, 14\} \subseteq Z_{21}$. $W = (W_1, W_2, W_3, W_4, W_5)$ is a special group set vector subspace of V over the same set S.

Suppose $T = (T_1, T_2, \ldots, T_5)$ where $T_1 = \{0, 5\} \subseteq Z_{10} = V_1$, $T_2 = \{0, 2, 4, 6, 8, 10\} \subseteq Z_{12} = V_2$, $T_3 = \{0, 5, 10, 15, 20\} \subseteq Z_{25}$

$= V_3$, $T_4 = \{0, 3, 6\} \subseteq Z_9 = V_4$ and $T_5 = \{0, 3, 6, 9, 12, 15, 18\}$ $\subseteq Z_{21} = V_5$. Clearly $T_1, T_2, \ldots, T_5$ are subgroups of $V_1, \ldots, V_5$ respectively. Thus $T = (T_1, T_2, \ldots, T_5) \subseteq (V_1, V_2, \ldots, V_5)$ is a special group set vector subspace of V over the same set S.

***Example 3.2.29:*** Let $V = (V_1, V_2, V_3, V_4)$ where

$$V_1 = \left\{ \begin{pmatrix} a & b \\ c & d \end{pmatrix} \middle| a, b, c, d \in Z \right\},$$

$$V_2 = \left\{ \begin{pmatrix} a & a & a \\ a & a & a \end{pmatrix} \middle| a \in Z \right\},$$

$V_3 = (Z \times Z \times Z \times Z)$ and $V_4 = (3Z \times 5Z \times 7Z \times 11Z \times 15Z)$. Clearly $V = (V_1, V_2, V_3, V_4)$ is a special group set vector space over the set $S = \{0, 1, 2, \ldots, 20\}$. Now take $W = (W_1, W_2, W_3, W_4)$ where

$$W_1 = \left\{ \begin{pmatrix} a & a \\ a & a \end{pmatrix} \middle| a \in Z \right\}$$

a proper additive subgroup of $V_1$,

$$W_2 = \left\{ \begin{pmatrix} a & a & a \\ a & a & a \end{pmatrix} \middle| a \in 3Z \right\} \subseteq V_2$$

is a proper subgroup of $V_2$, $W_3 = \{Z \times \{0\} \times Z \times \{0\}\} \subseteq V_3$ is also a subgroup of $V_3$ and $W_4 = (3Z \times 5Z \times \{0\} \times \{0\} \times 15Z) \subseteq V_4$ is a subgroup $V_4$. Thus $W = (W_1, W_2, W_3, W_4) \subseteq (V_1, V_2, V_3, V_4) = V$ is a special group set vector subspace of V over the set S.

***Example 3.2.30:*** Let $V = (V_1, V_2, V_3, V_4, V_5)$ where $V_1 = Z_5$, $V_2 = Z_7$, $V_3 = Z_{11}$, $V_4 = Z_2$ and $V_5 = Z_{13}$ all of them are additive groups modulo n where $n \in \{5, 7, 11, 2, 13\}$. Clearly V is a special group set vector space over the set $S = \{0, 1\}$.

Now we see V does not contain a proper subset $W = (W_1, W_2, W_3, W_4, W_5)$ such that W is a special group vector subspace of V that is each of the groups $V_i$ in V have no proper subgroup, the only subgroups of each of these $V_i$ is just $\{0\}$ and $V_i$, this is true for i =1, 2, 3, 4, 5. Thus W does not exists for this V.

In view of this example we are interested to define yet a new type of special group set vector spaces over the set S.

**DEFINITION 3.2.16:** *Let $V = (V_1, V_2, \ldots, V_n)$ be a special group set vector space over the set S. If V does not contain a proper n-subset $W = (W_1, W_2, \ldots, W_n) \subseteq (V_1, \ldots, V_n) = V$ such that $V_i$'s cannot have $W_i$'s which are different from $V_i$'s and $\{0\}$ i.e., if V has no proper special group set vector subspace over S then, we call V to be a special group set simple vector space over the set S.*

We have just now given an example of a special group set simple vector space. Now we give yet another example.

*Example 3.2.31:* Let $V = (V_1, V_2, V_3, V_4)$ where $V_1 = Z_{11}$, $V_2 = Z_{19}$, $V_3 = Z_{23}$ and $V_4 = Z_{29}$. V is a special group set vector space over the set $\{0, 1, 2, 3, 4, 5, 6\}$. Clearly each $V_i$ is a simple group for they do not have proper subgroups, hence $V = (V_1, V_2, V_3, V_4)$ is a special group set simple vector space over the set S.

Now we give yet another example.

*Example 3.2.32:* Let $V = (V_1, V_2, V_3, V_4, V_5)$ where $V_1 = Z_5$, $V_2 = Z_6$, $V_3 = Z_7$, $V_4 = Z_9$ and $V_5 = Z_{11}$. V is a special group set vector space over the set $\{0, 1, 3\}$. Take $W = (W_1, W_2, W_3, W_4, W_5)$ where $W_1 = Z_5$, $W_2 = \{0, 2, 4\} \subseteq Z_6 = V_2$, $W_3 = V_3$, $W_4 = \{0, 3, 6\} \subseteq Z_9 = V_4$ and $W_5 = \{0\}$. Now $W = (W_1, W_2, W_3, W_4, W_5)$ is a special group set vector space over the set S. We call this W by a different name for only some groups $V_i$ in V are simple and other groups $V_j$ in V are not simple.

In view of this we give yet another new definition.

**DEFINITION 3.2.17:** *Let $V = (V_1, V_2, ..., V_n)$ be a special group set vector space such that some $V_i$'s do not have proper subgroups and other $V_j$'s have proper subgroups $i \neq j$ ($1 \leq i, j \leq n$). Let $W = (W_1, W_2, ..., W_n) \subseteq (V_1, V_2, ..., V_n) = V$ such that some of the subgroups $W_i$ in $W$ are trivial subgroups of $V_i$ in $V$ and some of the subgroups $W_j$ in $W$ are nontrivial subgroups of $V_j$ ($i \neq j$); $1 \leq i, j \leq n$. We call this $W$ to be a special group set semisimple vector subspace over $S$ provided $W$ is a special group set vector space over the set $S$.*

We now illustrate this by some simple examples.

*Example 3.2.33:* Let $V = (V_1, V_2, V_3, V_4)$ where $V_1 = Z \times Z$, $V_2 = Z_7$, $V_3 = Z_5$ and

$$V_4 = \left\{ \begin{pmatrix} a & b \\ c & d \end{pmatrix} \middle| a, b, c, d \in Z \right\}.$$

$V$ is clearly a special group set vector space over the set $S = \{0, 1\}$. Take $W = \{W_1, W_2, W_3, W_4\}$ where $W_1 = 2Z \times 3Z \subseteq V_1$, $W_2 = Z_7 = V_2$, $W_3 = Z_5 = V_3$ and

$$W_4 = \left\{ \begin{pmatrix} a & b \\ c & d \end{pmatrix} \middle| a, b, c, d \in 3Z \right\} \subseteq V_4.$$

Clearly $W = (W_1, W_2, W_3, W_4)$ is a special group set vector space over the set $S = \{0, 1\}$. Thus $W = (W_1, W_2, W_3, W_4)$ is a special group set semisimple vector subspace of $V$ over the set $S = \{0, 1\}$.

Now we define yet another substructure of these special group set vector spaces over a set $S$.

**DEFINITION 3.2.18:** *Let $V = (V_1, ..., V_n)$ be a special group set vector space over the set $S$. Let $W = (W_1, W_2, ..., W_n) \subseteq (V_1, V_2, ..., V_n)$ such that each $W_i \subseteq V_i$ is a subsemigroup of $V_i$ under the same operations of $W_i$ and if $W = (W_1, ..., W_n)$ happens to be a special semigroup set vector space over the set $S$ then we call $W$ to be a pseudo special semigroup set vector subspace of $V$.*

The following observations are pertinent at this juncture.

1. All special group set vector subspaces of V are trivially pseudo special semigroup set vector subspaces of V.
2. Pseudo special semigroup set vector subspaces of V is never a special group set vector subspace of V for atleast one of the $W_i$'s in W will not be a group only semigroup.

We illustrate this situation by some simple examples.

***Example 3.2.34:*** Let $V = (V_1, V_2, V_3, V_4, V_5)$ where $V_1 = Z \times Z \times Z$,

$$V_2 = \left\{ \begin{pmatrix} a & b \\ c & d \end{pmatrix} \middle| a,b,c,d \in Z \right\},$$

$$V_3 = \left\{ \begin{pmatrix} a & a & a \\ a & a & a \end{pmatrix} \middle| a \in Z \right\},$$

$V_4 = \{Z[x]$ / all polynomials of degree less than or equal to 5 in the variable x with coefficients from $Z\}$

and

$$V_5 = \left\{ \begin{pmatrix} a & 0 & 0 & 0 \\ b & d & 0 & 0 \\ c & e & f & 0 \\ g & h & i & j \end{pmatrix} \middle| a,b,c,d,e,f,g,h,i,j \in Z \right\},$$

is a special group set vector space over the set $S = \{0, 1, 3, 5, 4, 9\}$. Take $W = \{(W_1, W_2, W_3, W_4, W_5\} \subseteq (V_1, V_2, V_3, V_4, V_5) = V$ where $W_1 = Z^+ \cup \{0\} \times Z^+ \cup \{0\} \times Z^+ \cup \{0\} \subseteq Z \times Z \times Z = V_1$,

$$W_2 = \left\{ \begin{pmatrix} a & b \\ c & d \end{pmatrix} \middle| a,b,c,d \in Z^+ \cup \{0\} \right\} \subseteq V_2,$$

$$W_3 = \left\{ \begin{pmatrix} a & a & a \\ a & a & a \end{pmatrix} \middle| a \in Z^+ \cup \{0\} \right\} \subseteq V_3,$$

$W_4 = \{Z^+ \cup \{0\}[x]$ / all polynomials of degree less than or equal to 5 in the variable x with coefficients from $Z^+ \cup \{0\}\} \subseteq V_4$ and

$$W_5 = \left\{ \begin{pmatrix} a & 0 & 0 & 0 \\ a & a & 0 & 0 \\ a & a & a & 0 \\ a & a & a & a \end{pmatrix} \middle| a \in Z^+ \cup \{0\} \right\} \subseteq V_5.$$

Clearly $W = (W_1, \ldots, W_5) \subseteq (V_1, V_2, \ldots, V_5)$ and W is only a pseudo special semigroup set vector subspace of V over S.

We see in this example all $W_i$ in W are only semigroups and not subgroups of $V_i$ in V.

*Example 3.2.35:* Let $V = (V_1, V_2, V_3, V_4)$ where $V_1 = Z \times Z$,

$$V_2 = \left\{ \begin{pmatrix} a & b \\ c & d \end{pmatrix} \middle| a, b, c, d \in 2Z \right\},$$

$V_3 = Z_{12}$ and $V_4 = \left\{ \begin{pmatrix} a & a & a & a \\ a & a & a & a \end{pmatrix} \middle| a \in 3Z \right\}$,

V is a special group set vector space over the set {0, 1, 2, 3, 4, 5, 6, 7, 8, 9, 10}. Take $W = \{W_1, W_2, W_3, W_4\} \subset (V_1, V_2, V_3, V_4) = V$ where $W_1 = Z^+\cup\{0\}\times Z^+\cup\{0\}$, $W_3 = \{0, 2, 4, 6, 8, 10\} \subseteq V_3$,

$$W_2 = \left\{ \begin{pmatrix} a & b \\ c & d \end{pmatrix} \middle| a, b, c, d \in 2Z^+ \cup \{0\} \right\} \subseteq V_2$$

and

$$W_4 = \left\{ \begin{pmatrix} a & a & a & a \\ a & a & a & a \end{pmatrix} \middle| a \in 3Z^+ \cup \{0\} \right\} \subseteq V_4.$$

$W = (W_1, W_2, W_3, W_4) \subseteq V$ is a pseudo special semigroup set vector subspace of V but $W_3$ is only a subgroup.

In view of this we define yet another substructure for V.

**DEFINITION 3.2.19:** *Let $V = (V_1, \ldots, V_n)$ be a special group set vector space over the set S. Suppose $W = (W_1, W_2, \ldots, W_n) \subseteq (V_1, \ldots, V_n) = V$ such that only some of the $W_i \subseteq V_i$ are semigroups and other $W_j \subseteq V_j$ are only groups (so trivially a semigroup) $i \neq j$, $(1 \leq i, j \leq n)$, then we call $W = (W_1, W_2, \ldots, W_n) \subseteq V$ to be a pseudo special quasi semigroup set vector subspace of V over the set S.*

We illustrate this by some simple examples.

***Example 3.2.36:*** Let $V = (V_1, V_2, V_3, V_4, V_5)$ where $V_1 = Z \times Z \times Z$, $V_2 = Z_{29}$, $V_3 = Z_{13}$,

$$V_4 = \left\{ \begin{pmatrix} a & b \\ c & d \end{pmatrix} \middle| a, b, c, d \in Z \right\}$$

and

$$V_5 = \left\{ \begin{pmatrix} a & a \\ a & a \\ a & a \\ a & a \\ a & a \end{pmatrix} \middle| a \in Z \right\};$$

V is a special group set vector space over the set $S = \{0, 1, 2, 3, 4, 5, 6\}$. Take $W = (W_1, W_2, W_3, W_4, W_5)$ where $W_1 = Z^o \times Z^o \times Z^o$ (where $Z^o = Z^+ \cup \{0\}) \subseteq V_1$, $W_2 = V_2 = Z_{29}$, $W_3 = Z_{13} = V_3$,

$$W_4 = \left\{ \begin{pmatrix} a & b \\ c & d \end{pmatrix} \middle| a, b, c, d \in Z^o = Z^+ \cup \{0\} \right\} \subseteq V_4$$

and

$$W_5 = \left\{ \begin{pmatrix} a & a \\ a & a \\ a & a \\ a & a \\ a & a \end{pmatrix} \middle| a \in Z^o = Z^+ \cup \{0\} \right\} \subseteq V_5.$$

Clearly $W = (W_1, W_2, W_3, W_4, W_5) \subseteq (V_1, V_2, V_3, V_4, V_5) = V$ is only a pseudo special quasi semigroup set vector subspace of V over S.

It is important to note that all pseudo special semigroup set vector subspaces of V are pseudo special quasi semigroup vector subspaces of V, however every pseudo quasi semigroup set vector subspaces of V need not in general be a pseudo special semigroup set vector subspaces of V over S in general.

The above example is an illustration of the above statement.

Now we proceed into define the new notion of special group set linear algebras over the group G.

**DEFINITION 3.2.20:** *Let $V = (V_1, V_2, \ldots, V_n)$ be a special group set vector space over the set S. If the set S is closed under addition and is an additive abelian group then we call V to be a special group set linear algebra over the group S, we only demand for every $v \in V_i$ and $s \in S$, $sv \in V_i$; $1 \leq i \leq n$.*

We illustrate this by some simple examples.

*Example 3.2.37:* $V = (V_1, V_2, V_3, V_4)$ is a special group set linear algebra over the group $G = \{Z\}$ where $V_1 = Z \times Z \times Z$,

$$V_2 = \left\{ \begin{pmatrix} a & b \\ c & d \end{pmatrix} \middle| a, b, c, d \in Z \right\},$$

$V_3 = \{Z[x]$ all polynomials of degree less than or equal to 4 in the variable x with coefficients from Z$\}$ and

$$V_4 = \left\{ \begin{pmatrix} a & a & a & a & a \\ a & a & a & a & a \end{pmatrix} \middle| a \in Z \right\}.$$

*Example 3.2.38:* Let $V = (V_1, V_2, V_3, V_4, V_5)$ where

$$V_1 = Z_{10} \times Z_{10} \times Z_{10},$$

$$V_2 = \left\{ \begin{pmatrix} a & b & c \\ d & e & f \end{pmatrix} \middle| a,b,c,d,e,f \in Z_{10} \right\},$$

$$V_3 = \left\{ \begin{pmatrix} x & 0 & 0 \\ y & z & 0 \\ g & h & i \end{pmatrix} \middle| x,y,z,g,h,i \in Z_{10} \right\},$$

$$V_4 = \left\{ \begin{pmatrix} a & b \\ c & d \\ e & f \end{pmatrix} \middle| a,b,c,d,e,f \in Z_{10} \right\}$$

and $V_5 = \{Z_{10}[x]$; all polynomials of degree less than or equal to three. Then $V = (V_1, V_2, \ldots, V_5)$ is a special group set linear algebra over the group $Z_{10}$ under addition modulo 10.

Now we proceed onto illustrate various substructure of a special group set linear algebra over a group G.

**DEFINITION 3.2.21:** *Let $V = (V_1, V_2, \ldots, V_n)$ be a special group set linear algebra over the group G. If $W = (W_1, \ldots, W_n) \subseteq (V_1, V_2, \ldots, V_n) = V$ be a proper subset of V such that for each i, $W_i \subseteq V_i$; is a subgroup of $V_i$, $1 \leq i \leq n$ then we say $W = (W_1, \ldots, W_n)$ is a special n-subgroup of V. If $W = (W_1, \ldots, W_n) \subseteq V$ is a special group set linear algebra over the same group G, then we call W to be a special group set linear subalgebra of V over G.*

We illustrate this by some simple examples.

***Example 3.2.39:*** Let $V = (V_1, V_2, \ldots, V_5)$ where

$$V_1 = \left\{ \begin{pmatrix} a & b \\ c & d \end{pmatrix} \middle| a,b,c,d \in Z \right\},$$

$V_2 = \{Z[x];$ all polynomials of degree less than or equal to five in the variable x with coefficients from $Z\}$,

$$V_3 = \left\{ \begin{pmatrix} a & a \\ a & a \\ a & a \\ a & a \\ a & a \end{pmatrix} \middle| a \in Z \right\}, \quad V_4 = \{Z \times Z \times Z \times Z \times Z\}$$

and

$$V_5 = \left\{ \begin{pmatrix} a & a & a & a \\ a & a & a & a \end{pmatrix} \middle| a \in Z \right\}.$$

$V = (V_1, V_2, V_3, V_4, V_5)$ is a special group set linear algebra over the group Z. Let $W = (W_1, W_2, \ldots, W_5) \subseteq (V_1, V_2, \ldots, V_5)$ where

$$W_1 = \left\{ \begin{pmatrix} a & b \\ c & d \end{pmatrix} \middle| a,b,c,d \in 2Z \right\} \subseteq V_1,$$

is a proper subgroup of $V_1$, $W_2 = \{$All polynomials in the variable x of degree less than or equal to five with coefficients from $3Z\}$, $W_3 \subseteq V_3$, where

$$W_3 = \left\{ \begin{pmatrix} a & a \\ a & a \\ a & a \\ a & a \\ a & a \end{pmatrix} \middle| a \in 5Z \right\},$$

$$W_4 = \{3Z \times 3Z \times 3Z \times 3Z \times 3Z\} \subseteq V_4$$

and

$$W_5 = \left\{ \begin{pmatrix} a & a & a & a \\ a & a & a & a \end{pmatrix} \middle| a \in 5Z \right\} \subseteq V_5,$$

$W = (W_1, W_2, W_3, W_4, W_5) \subseteq (V_1, V_2, \ldots, V_5)$ is a special group set linear subalgebra over the group Z.

***Example 3.2.40:*** Let $V = (V_1, V_2, V_3, V_4)$ where $V_1 = \{Z_{12} \times Z_{12} \times Z_{12}\}$, $V_2 = \{Z_{12}[x] /$ all polynomials of degree less than or equal to 5 in the variable x with coefficients from $Z_{12}\}$,

$$V_3 = \left\{ \begin{pmatrix} a & b \\ c & d \end{pmatrix} \middle| a,b,c,d \in Z_{12} \right\}$$

and

$$V_4 = \left\{ \begin{pmatrix} a & a & a & a & a \\ a & a & a & a & a \end{pmatrix} \middle| a \in Z_{12} \right\}.$$

$V = (V_1, V_2, V_3, V_4)$ is a special group set linear algebra over the group $Z_{12}$. Let $W = (W_1, W_2, W_3, W_4)$ where $W_1 = Z_{12} \times \{0\} \times Z_{12} \subseteq V_1$, $W_2 = \{$All polynomials in the variable x of degree less than or equal to 5 with coefficients from $\{0, 2, 4, 6, 8, 10\}\} \subseteq V_2$,

$$W_3 = \left\{ \begin{pmatrix} a & a \\ a & a \end{pmatrix} \middle| a \in Z_{12} \right\} \subseteq V_3$$

and

$$W_4 = \left\{ \begin{pmatrix} a & a & a & a & a \\ a & a & a & a & a \end{pmatrix} \middle| a \in \{0,2,4,6,8,10\} \right\} \subseteq V_4,$$

$W = (W_1, W_2, W_3, W_4) \subseteq V$ is a special group set linear subalgebra of V over $Z_{12}$.

We define yet another new substructure.

**DEFINITION 3.2.22:** *Let $V = (V_1, V_2, ..., V_n)$ be a special group set linear algebra over the group G. Let $W = (W_1, ..., W_n) \subseteq (V_1, V_2, ..., V_n)$ be a n-subgroup of V, i.e., each $W_i$ is a subgroup of $V_i$, $1 \leq i \leq n$. If for some proper subgroup H of G we have W to be a special group set linear algebra over H then we call $W = (W_1, W_2, ..., W_n)$ to be a special subgroup set linear subalgebra of V over the subgroup H of G.*

We illustrate this by some examples.

*Example 3.2.41:* Let $V = (V_1, V_2, V_3, V_4)$ be a special group set linear algebra over the group Z, where $V_1 = Z \times Z \times Z \times Z$,

$$V_2 = \left\{ \begin{pmatrix} a & b \\ c & d \end{pmatrix} \middle| a,b,c,d \in Z \right\},$$

$$V_3 = \left\{ \begin{pmatrix} a & a & a & a \\ a & a & a & a \end{pmatrix} \middle| a \in Z \right\}$$

and

$$V_4 = \left\{ \begin{pmatrix} a & b \\ c & d \\ a & b \\ c & d \end{pmatrix} \middle| a,b,c,d \in Z \right\}.$$

Take $W = (W_1, W_2, W_3, W_4) \subseteq (V_1, V_2, V_3, V_4) = V$ where $W_1 = 2Z \times 2Z \times 2Z \times 2Z \subseteq V_1$,

$$W_2 = \left\{ \begin{pmatrix} a & b \\ c & d \end{pmatrix} \middle| a,b,c,d \in 2Z \right\} \subseteq V_2,$$

$$W_3 = \left\{ \begin{pmatrix} a & a & a & a \\ a & a & a & a \end{pmatrix} \middle| a \in 2Z \right\} \subseteq V_3$$

and

$$W_4 = \left\{ \begin{pmatrix} a & b \\ c & d \\ a & b \\ c & d \end{pmatrix} \middle| a,b,c,d \in 2Z \right\} \subseteq V_4.$$

Clearly $W = (W_1, W_2, W_3, W_4) \subseteq (V_1, V_2, V_3, V_4) = V$ is a special subgroup linear subalgebra of V over the subgroup $2Z \subseteq Z$.

*Example 3.2.42:* Let $V = (V_1, V_2, \ldots, V_5)$ where $V_1 = \{Z_6 \times Z_6 \times Z_6 \times Z_6 \times Z_6\}$,

$$V_2 = \left\{ \begin{pmatrix} a & a \\ a & a \end{pmatrix} \middle| a \in Z_6 \right\},$$

$$V_3 = \left\{ \begin{pmatrix} a & a & a & a & a \\ a & a & a & a & a \end{pmatrix} \middle| a \in Z_6 \right\},$$

$$V_4 = \left\{ \begin{pmatrix} a & a \\ a & a \\ a & a \\ a & a \end{pmatrix} \middle| a \in Z_6 \right\}$$

and $V_5 = \{Z_6[x]$ all polynomials in the variable x with coefficients from $Z_6$ of degree less than or equal to 6$\}$. $V = (V_1, V_2, V_3, V_4, V_5)$ is a special group set linear algebra over the group $Z_6$. Take $W = (W_1, W_2, \ldots, W_5) \subseteq (V_1, V_2, \ldots, V_5) = V$ where $W_1 = \{S \times S \times S \times S \times S$ where $S = \{0, 2, 4\}\} \subseteq V_1$,

$$W_2 = \left\{ \begin{pmatrix} a & a \\ a & a \end{pmatrix} \middle| a \in \{0, 2, 4\} \right\} \subseteq V_2,$$

$$W_3 = \left\{ \begin{pmatrix} a & a & a & a & a \\ a & a & a & a & a \end{pmatrix} \middle| a \in \{0, 2, 4\} \right\} \subseteq V_3,$$

$$W_4 = \left\{ \begin{pmatrix} a & a \\ a & a \\ a & a \\ a & a \end{pmatrix} \;\middle|\; a \in \{0, 2, 4\} \right\}$$

$\subseteq V_4$ and $W_5 = \{$all polynomials in the variable x with coefficients from the set $S = \{0, 2, 4\}$ of degree less than or equal to 6$\} \subseteq V_5$. $W = (W_1, W_2, W_3, W_4, W_5) \subseteq (V_1, V_2, \ldots, V_5) = V$ is a special subgroup set linear subalgebra of V over the subgroup $H = \{0, 2, 4\} \subseteq Z_6$.

***Example 3.2.43:*** Let $V = (V_1, V_2, V_3, V_4)$ where $V_1 = Z_7$, $V_2 = Z_5$, $V_3 = Z_{11}$, and $V_4 = Z_{13}$ be group such that V is a special group set linear algebra over the group $Z_2 = \{0, 1\} = G$ addition modulo 2.

We see G has no proper subgroups. Further each $V_i$ is such that, they do not contain proper subgroups; $1 \leq i \leq 4$. Thus V has no special group set linear subalgebra as well as V does not contain any special subgroup linear subalgebra.

In view of this example we define a special type of linear algebra.

**DEFINITION 3.2.23:** *Let $V = (V_1, V_2, \ldots, V_n)$ be a special group set linear algebra over the group G. Let $W = (W_1, W_2, \ldots, W_n) \subseteq (V_1, V_2, \ldots, V_n) = V$ such that for each i, $W_i \subseteq V_i$ is either $W_i = \{0\}$ or $W_i = V_i$, i.e., none of the $V_i$'s have proper subgroup; $1 \leq i \leq n$. Suppose G has no proper subgroup $H \subseteq G$, i.e., either $H = \{0\}$ or $H = G$. Then we call $V = (V_1, V_2, \ldots, V_n)$ to be a doubly simple special group linear algebra over G.*

***Example 3.2.44:*** Let $V = (V_1, V_2, V_3, V_4, V_5, V_6)$ where $V_1 = Z_3$, $V_2 = Z_5$, $V_3 = Z_7$, $V_4 = Z_{11}$, $V_5 = Z_{13}$ and $V_6 = Z_{17}$. V is a special group set linear algebra over the group $G = Z_2 = \{0, 1\}$ addition modulo 2. Clearly V is a doubly simple special group set linear algebra over G.

Now we proceed onto define yet another new algebraic substructure of V.

**DEFINITION 3.2.24:** *Let $V = (V_1, V_2, ..., V_n)$ be a special group linear algebra over the group G. Let $W = (W_1, W_2, ..., W_n) \subseteq (V_1, V_2, ..., V_n) = V$ be a proper n-subset of V such that each $W_i \subseteq V_i$ is a subsemigroup of $V_i$, under the same operations of $V_i$, i=1,2,..., n. If $W = (W_1, ..., W_n)$ is a special set semigroup linear algebra over the semigroup H of G, then we call W to be a pseudo special set subsemigroup linear subalgebra of V over $H \subseteq G$.*

We illustrate this situation by an example.

***Example 3.2.45:*** Let $V = (V_1, V_2, V_3, V_4)$ where $V_1 = Z \times Z \times Z \times Z$,

$$V_2 = \left\{ \begin{pmatrix} a & b \\ c & d \end{pmatrix} \middle| a,b,c,d \in Z \right\},$$

$V_3 = \{Z[x]$, all polynomials in the variable x with coefficients from Z of degree less than or equal to 5$\}$ and

$$V_4 = \left\{ \begin{pmatrix} a & a \\ a & a \\ a & a \\ a & a \\ a & a \\ a & a \end{pmatrix} \middle| a \in Z \right\}.$$

$V = (V_1, V_2, V_3, V_4)$ is a special group set linear algebra over the group Z. Let $W = (W_1, W_2, W_3, W_4) \subseteq (V_1, V_2, V_3, V_4)$ where $W_1 = Z_o \times Z_o \times Z_o$ where $Z_o = Z^+ \cup \{0\}\} \subseteq V_1$,

$$W_2 = \left\{ \begin{pmatrix} a & b \\ c & d \end{pmatrix} \middle| a,b,c,d \in Z^o = Z^+ \cup \{0\} \right\} \subseteq V_2,$$

$W_3 = \{Z^o[x]$ all polynomials in x of degree less than or equal to 5 with coefficients from $Z^o\} \subseteq V_3$ and

$$W_4 = \left\{ \begin{pmatrix} a & a \\ a & a \\ a & a \\ a & a \\ a & a \\ a & a \end{pmatrix} \middle| a \in Z^o = Z^+ \cup \{0\} \right\} \subseteq V_4.$$

$W = (W_1, W_2, W_3, W_4)$ is a special semigroup set linear algebra over the semigroup $Z^o \subseteq Z$. Thus W is a pseudo special subsemigroup set linear subalgebra of V over the subsemigroup $Z^o \subseteq Z$.

*Example 3.2.46:* Let $V = (V_1, V_2, V_3, V_4, V_5)$ where $V_1 = Z_{16} \times Z_{16} \times Z_{16}$,

$$V_2 = \left\{ \begin{pmatrix} a & b \\ c & d \end{pmatrix} \middle| a, b, c, d \in Z_{16} \right\},$$

$V_3 = \{Z_{16}[x]$ all polynomials of degree less than or equal to 2$\}$,

$$V_4 = \left\{ \begin{pmatrix} a & a & a \\ a & a & a \end{pmatrix} \middle| a \in Z_{16} \right\}$$

and

$$V_5 = \left\{ \begin{pmatrix} a & a & a \\ a & a & a \\ a & a & a \\ a & a & a \\ a & a & a \\ a & a & a \end{pmatrix} \middle| a \in Z_{16} \right\}$$

be a special group set linear algebra over the group $Z_{16}$. $W = (W_1, W_2, W_3, W_4, W_5)$ where $W_1 = \{S \times S \times S \mid S = \{0, 2, 4, 6, 8, 10, 12, 14\}\} \subseteq V_1$,

$$W_2 = \left\{ \begin{pmatrix} a & a \\ a & a \end{pmatrix} \middle| a \in Z_{16} \right\} \subseteq V_2,$$

$W_3 = $ {all polynomials in x of degree less than to 2 with coefficients from $S = \{0, 2, 4, \ldots, 14\}$},

$$W_4 = \left\{ \begin{pmatrix} a & a & a \\ a & a & a \end{pmatrix} \middle| a \in S = \{0,2,4,6,8,10,12,14\} \right\} \subseteq V_4$$

and

$$W_5 = \left\{ \begin{pmatrix} a & a & a \\ a & a & a \\ a & a & a \\ a & a & a \\ a & a & a \\ a & a & a \end{pmatrix} \middle| a \in S = \{0,2,4,6,8,10,12,14\} \right\} \subseteq V_5.$$

$W = (W_1, W_2, W_3, W_4, W_5) \subseteq (V_1, V_2, \ldots, V_5)$ is a special subgroup set linear subalgebra over the subgroup $S = \{0, 2, 4, \ldots, 14\} \subseteq Z_{16}$.

Now we proceed on to define another new substructure in $V = (V_1, V_2, \ldots, V_n)$.

**DEFINITION 3.2.25:** *Let $V = (V_1, V_2, \ldots, V_n)$ be a special group set linear algebra over the group G. Suppose $H = (H_1, \ldots, H_n) \subseteq (V_1, V_2, \ldots, V_n)$ such that each $H_i \subseteq V_i$, $1 \leq i \leq n$ is a subsemigroup of the group $V_i$ then we all $H = (H_1, \ldots, H_n)$ the pseudo special set subsemigroup linear subalgebra over the subsemigroup P of the group G if $H = (H_1, \ldots, H_n)$ is a special semigroup set linear algebra over the semigroup P of the group G.*

We illustrate this situation by a simple example.

***Example 3.2.47:*** Let $V = (V_1, V_2, V_3, V_4, V_5)$ where $V_1 = Z \times Z$,

$$V_2 = \left\{ \begin{pmatrix} a & b \\ c & d \end{pmatrix} \middle| a,b,c,d \in Z \right\},$$

$$V_3 = \left\{ \begin{pmatrix} a & a & a & a \\ a & a & a & a \end{pmatrix} \middle| a \in Z \right\},$$

$$V_4 = \left\{ \begin{pmatrix} a & a \\ a & a \\ a & a \\ a & a \\ a & a \end{pmatrix} \middle| a \in Z \right\}$$

and $V_5 = \{Z[x]$ all polynomials of degree less than or equal to 3 with coefficients from $Z\}$. $V$ is a special group set linear algebra over the group $Z$. Let $H = (H_1, H_2, \ldots, H_5) \subseteq (V_1, \ldots, V_5)$ where $H_1 = \{Z^+ \cup \{0\} \times Z^+ \cup \{0\}\} \subseteq V_1$,

$$H_2 = \left\{ \begin{pmatrix} a & a \\ a & a \end{pmatrix} \middle| a \in Z^+ \cup \{0\} \right\} \subseteq V_2,$$

$$H_3 = \left\{ \begin{pmatrix} a & a & a & a \\ a & a & a & a \end{pmatrix} \middle| a \in Z^+ \cup \{0\} \right\} \subseteq V_3,$$

$$H_4 = \left\{ \begin{pmatrix} a & a \\ a & a \\ a & a \\ a & a \\ a & a \end{pmatrix} \middle| a \in Z^+ \cup \{0\} \right\}$$

$\subseteq V_4$ and $H_5 = \{(Z^+ \cup \{0\})[x]$ / this consists of all polynomials of degree less than or equal to 3 with coefficients from $Z^+ \cup \{0\}\} \subseteq V_5$. $H = (H_1, H_2, \ldots, H_5) \subseteq (V_1, V_2, \ldots, V_5) = V$ is a pseudo special subsemigroup linear subalgebra over the subsemigroup $Z^+ \cup \{0\} \subseteq Z$.

*Example 3.2.48:* Let $V = (V_1, V_2, V_3, V_4)$ where $V_1 = Q \times Q \times Q$, $V_2 = \{Q[x]$ / all polynomials of degree less than or equal to 5 with coefficients from Q in the variable x$\}$,

$$V_3 = \left\{ \begin{pmatrix} a & b \\ c & d \end{pmatrix} \middle| a, b, c, d \in Q \right\}$$

and

$$V_4 = \left\{ \begin{pmatrix} a & a \\ a & a \\ a & a \\ a & a \end{pmatrix} \middle| a \in Q \right\}.$$

V is a special group set linear algebra over the group Z. Take $W = (W_1, W_2, W_3, W_4)$ where $W_1 = \{Z \times Z \times Z^+ \cup \{0\}\} \subseteq V_1$, $W_2 = \{S[x]$, where $S = Z^+ \cup \{0\}$, all polynomials in the variable x with coefficients from S of degree less than or equal to 5$\} \subseteq V_2$,

$$W_3 = \left\{ \begin{pmatrix} a & b \\ c & d \end{pmatrix} \middle| a, b, c, d \in S = Z^+ \cup \{0\} \right\} \subseteq V_3$$

and

$$W_4 = \left\{ \begin{pmatrix} a & a \\ a & a \\ a & a \\ a & a \end{pmatrix} \middle| a \in Z^+ \cup \{0\} = S \right\} \subseteq V_4.$$

Clearly each $W_i$ is a only semigroup under addition; $1 \leq i \leq 4$; and $W = (W_1, W_2, W_3, W_4) \subseteq (V_1, V_2, V_3, V_4)$ is a pseudo special subsemigroup set linear subalgebra of V over the subsemigroup $S = Z^+ \cup \{0\}$ contained in the group Q.

Now we having seen some substructures in V we now proceed onto define the notion of linear transformations and linear operators.

**DEFINITION 3.2.26:** *Let $V = (V_1, ..., V_n)$ and $W = (W_1, ..., W_n)$ be two special group set vector spaces over the same group G. Let $T = (T_1, ..., T_n): V \to W$ be defined by $T_i : V_i \to W_i$, $i = 1, 2, ..., n$, such that $T_i(x\mu + \upsilon) = xT_i(\mu) + T_i(\upsilon)$ for $x \in G$ and $\mu, \upsilon$ in $V_i$; $i = 1, 2, ..., n$. We call T to be a special group set linear transformation on V. Suppose $SHom_G(V, W) = \{Hom_G(V_1, W_1), Hom_G(V_2, W_2), ..., Hom_G(V_n, W_n)\}$ where for each i, $Hom_G(V_i, W_i) = \{$all maps $T_i$ from $V_i$ to $W_i$ such that $T_i$'s are group special linear transformations$\}$. Clearly $\theta : V_i \to W_i$; $\theta(\upsilon_i) = 0$ for all $v_i \in V_i$, for every $T_i$, $P_i$ from $V_i$ to $W_i$ we have*

$$
\begin{aligned}
(T_i + P_i)(x\mu+\upsilon) &= T_i(x\mu+\upsilon) + P_i(x\mu+\upsilon) \\
&= xT_i(\mu) + T_i(\upsilon)) + x(P_i(\mu) + P_i(\upsilon)) \\
&= x[(T_i + P_i)(\mu)] + (T_i + P_i)(\upsilon)
\end{aligned}
$$

*for all $\mu, \upsilon \in V_i$ and $x \in G$. Now for every $T_i \in Hom_G(V_i, W_i)$ there exists, $-T_i \in Hom_G(V_i, W_i)$ by*

$$T_i(x\mu+\upsilon) = x T_i(\mu) + T_i(\upsilon)$$

*then*

$$-T_i(x\mu+\upsilon) = -xT_i(\mu) + (-T_i(\upsilon))$$

*so that*

$$
\begin{aligned}
T_i + (-T_i)(x\mu+\upsilon) &= (T_i - T_i)(x\mu) + (T_i - T_i)(\upsilon) \\
&= \theta(x\mu) + \theta(\upsilon) \\
&= 0 + 0 \\
&= 0(x\mu+\upsilon) \\
&= 0 + 0.
\end{aligned}
$$

*Thus if for $a \in G$ and $T_i \in Hom_G(V_i, W_i)$, $aT_i : V_i \to W_i$. Hence each $Hom_G(V_i, W_i)$ is a special group vector space over the group G.*

   *Thus $SHom(V,W) = \{Hom_G(V_1, W_1), Hom_G(V_2, W_2), ..., Hom_G(V_n, W_n)\}$ is also a special group set vector space over the group G. Suppose $T = (T_1, T_2, ..., T_n): V \to W$ be defined by $T_i :$*

$V_i \to W_j$ such that no two $V_i$'s are mapped onto same $W_j$ then also T is a special group set linear transformation from V to W and it is also denoted by $SHom_G(V, W) = \{Hom_G(V_1, W_{i_1}), Hom_G(V_2, W_{i_2}), \ldots, Hom_G(V_n, W_{i_n})\}$, i.e., $Hom_G(V_m, W_{i_m}) = $ {Set of all group linear transformation from $V_m$ to $W_{i_m}$}; clearly $Hom_G(V_m, W_{i_m})$ is an additive group i.e., each $Hom_G(V_m, V_{i_m})$ is a special group vector space over the group G. Thus $SHom_G(V, W)$ is also a special group set vector space over the group G.

We first illustrate this situation by some simple examples.

*Example 3.2.49:* Let $V = (V_1, V_2, V_3, V_4)$ and $W = (W_1, W_2, W_3, W_4)$ where $V_1 = Z_{10} \times Z_{10} \times Z_{10}$,

$$V_2 = \left\{ \begin{pmatrix} a & b \\ c & d \end{pmatrix} \middle| a,b,c,d \in Z_{10} \right\},$$

$$V_3 = \left\{ \begin{pmatrix} a & a & a & a \\ a & a & a & a \end{pmatrix} \middle| a \in Z_{10} \right\}$$

and $V_4 = \{Z_{10}[x]$ all polynomials of degree less than or equal to $5\}$ and $W_1 = Z_{10} \times Z_{10} \times Z_{10} \times Z_{10} \times Z_{10} \times Z_{10}$, $W_2 = \{Z_{10}[x]$ all polynomials of degree less than or equal to $2\}$,

$$W_3 = \left\{ \begin{pmatrix} a & a \\ a & a \\ a & a \\ a & a \end{pmatrix} \middle| a \in Z_{10} \right\}$$

and

$$W_4 = \left\{ \begin{pmatrix} a & b & 0 \\ 0 & 0 & 0 \\ c & d & 0 \end{pmatrix} \middle| a,b,c,d \in Z_{10} \right\}$$

are special group set vector space over the group $Z_{10}$. Define T = $(T_1, T_2, T_3, T_4)$: V → W defined by

$T_1 : V_1 → W_2$ where
$$T_1(a\ b\ c) = (a + bx + cx^2),$$

$T_2 : V_2 → W_4$ defined by
$$T_2\begin{pmatrix} a & b \\ c & d \end{pmatrix} = \begin{pmatrix} a & b & 0 \\ 0 & 0 & 0 \\ c & d & 0 \end{pmatrix},$$

$T_3 : V_3 → W_3$ defined by
$$T_3\begin{pmatrix} a & a & a & a \\ a & a & a & a \end{pmatrix} = \begin{pmatrix} a & a \\ a & a \\ a & a \\ a & a \end{pmatrix}$$

and $T_4 : V_4 → W_1$ defined by

$$T_4(a_0 + a_1x + a_2x^2 + a_3x^3 + a_4x^4 + a_5x^5) = (a_0, a_1, a_2, a_3, a_4, a_5).$$

Clearly T = $(T_1, T_2, T_3, T_4)$ is a special group set linear transformation of V into W. It is left as an exercise for the reader to find $SH_{Z_{10}}(V,W) = \{Hom_{Z_{10}}(V_1,W_2), Hom_{Z_{10}}(V_2,W_4), Hom_{Z_{10}}(V_3,W_3), Hom_{Z_{10}}(V_4, W_1)\}$. Is $SH_{Z_{10}}(V, W)$ is a special group set vector space over $Z_{10}$?

***Example 3.2.50:*** Let V = $(V_1, V_2, V_3, V_4)$ and W = $(W_1, W_2, W_3, W_4)$ where

$$V_1 = \left\{ \begin{pmatrix} a & b \\ c & d \end{pmatrix} \middle| a,b,c,d \in Z \right\},$$

$$V_2 = \left\{ \begin{pmatrix} a & a & a & a \\ a & a & a & a \end{pmatrix} \middle| a \in Z \right\},$$

$V_3 = \{Z[x]$ all polynomials of degree less than or equal to 5$\}$

and

$$V_4 = \left\{ \begin{pmatrix} a & 0 & 0 \\ b & d & 0 \\ c & e & f \end{pmatrix} \middle| a,b,c,d,e,f \in Z \right\}$$

is a special group set vector space over Z. $W_1 = \{Z \times Z \times Z \times Z\}$,

$$W_2 = \left\{ \begin{pmatrix} a & a \\ a & a \\ a & a \\ a & a \end{pmatrix} \middle| a \in Z \right\},$$

$$W_3 = \left\{ \begin{pmatrix} a & b & c \\ 0 & d & e \\ 0 & 0 & f \end{pmatrix} \middle| a,b,c,d,e,f \in Z \right\}$$

and

$$W_4 = \left\{ \begin{pmatrix} a & a & a \\ a & a & a \\ a & a & a \end{pmatrix} \middle| a \in Z \right\}$$

be a special group set vector space over Z. Define $T = (T_1, T_2, T_3, T_4)$ where
$T_1: V_1 \to W_1$ is defined by

$$T_1 \begin{pmatrix} a & b \\ c & d \end{pmatrix} = (a\ b\ c\ d),$$

$T_2 : V_2 \to W_2$ is given by

$$T_2 \begin{pmatrix} a & a & a & a \\ a & a & a & a \end{pmatrix} = \begin{pmatrix} a & a \\ a & a \\ a & a \\ a & a \end{pmatrix}$$

$T_3: V_3 \to W_3$ is such that

$$T_3(a_0 + a_1x + a_2x^2 + a_3x^3 + a_4x^4 + a_5x^5) = \begin{pmatrix} a_0 & a_1 & a_2 \\ 0 & a_3 & a_4 \\ 0 & 0 & a_5 \end{pmatrix}$$

and $T_4 : V_4 \to W_4$ is given by

$$T_4 \begin{pmatrix} a & 0 & 0 \\ b & d & 0 \\ c & e & f \end{pmatrix} = \begin{pmatrix} a & a & a \\ a & a & a \\ a & a & a \end{pmatrix}.$$

$T = (T_1, T_2, T_3, T_4)$ is a special group set linear transformation from V to W.

Now we proceed on to define special group set linear operator on the special group set vector space over a group G.

**DEFINITION 3.2.27:** *Let $V = (V_1, V_2, ..., V_n)$ be a special group set vector space over a group G. Let $T = (T_1,..., T_n)$ such that $T_i: V_i \to V_i$, $i = 1, 2, ..., n$; if each $T_i$ is a special group vector space over G, then we call T to be a special group set linear operator of V over G.*

We illustrate this by some examples.

***Example 3.2.51:*** Let $V = (V_1, V_2, V_3, V_4, V_5)$ be a special group set vector space over the group Z where $V_1 = Z \times Z \times Z \times Z$,

$$V_2 = \left\{ \begin{pmatrix} a & b \\ c & d \end{pmatrix} \middle| a,b,c,d \in Z \right\},$$

$V_3 = \{Z[x]$ all polynomials of degree less than or equal to 6$\}$
and

$$V_4 = \left\{ \begin{pmatrix} a & a \\ a & a \\ a & a \\ a & a \end{pmatrix} \middle| a \in Z \right\}.$$

Define $T = (T_1, T_2, T_3, T_4)$ from V to V by $T_1: V_1 \to V_1$ defined by
$$T_1(x\ y\ z\ \omega) = (x + y, y + z, z + \omega, \omega + x),$$

$T_2 : V_2 \to V_2$ is given by
$$T_2\begin{pmatrix} a & b \\ c & d \end{pmatrix} = \begin{pmatrix} a & a \\ a & a \end{pmatrix},$$

$T_3 : V_3 \to V_3$ is defined by
$$T_3(a_0 + a_1x + \ldots + a_6x^6) = a_0 + a_2x^2 + a_4x^4 + a_6x^6$$

and $T_4 : V_4 \to V_4$ is given by

$$T_4 \begin{pmatrix} a & a \\ a & a \\ a & a \\ a & a \end{pmatrix} = \begin{pmatrix} 2a & 2a \\ 2a & 2a \\ 2a & 2a \\ 2a & 2a \end{pmatrix}.$$

It is easily verified that $T = (T_1, T_2, T_3, T_4)$ is a special group set linear operator on V.

***Example 3.2.52:*** Let $V = (V_1, V_2, V_3, V_4, V_5)$ where

$$V_1 = \left\{ \begin{pmatrix} a & a \\ a & a \\ a & a \\ a & a \end{pmatrix} \middle| a \in Z_5 \right\},$$

$$V_2 = Z_5 \times Z_5 \times Z_5 \times Z_5,$$

$$V_3 = \left\{ \begin{pmatrix} a & b & c \\ d & e & f \\ g & h & i \end{pmatrix} \middle| a,b,c,d,e,f,g,h,i \in Z_5 \right\},$$

$$V_4 = \left\{ \begin{pmatrix} a & a & a & a \\ a & a & a & a \end{pmatrix} \middle| a \in Z_5 \right\}$$

and

$$V_5 = \left\{ \begin{pmatrix} a & 0 & 0 & 0 \\ b & c & 0 & 0 \\ d & e & f & 0 \\ g & h & i & j \end{pmatrix} \middle| a,b,c,d,e,f,g,h,i,j \in Z_5 \right\}.$$

V is a special group set vector space over the group $Z_5$. Let $T = (T_1, T_2, T_3, T_4, T_5)$ where

$T_1: V_1 \to V_1$ is defined by

$$T_1 \begin{pmatrix} a & a \\ a & a \\ a & a \\ a & a \end{pmatrix} = \begin{pmatrix} 2a & 2a \\ 2a & 2a \\ 2a & 2a \\ 2a & 2a \end{pmatrix},$$

$T_2 : V_2 \to V_2$ is given by
$$T_2(x\ y\ z\ w) = (x + \omega, y, z, x - \omega),$$

$T_3 : V_3 \to V_3$ is such that

$$T_3 \begin{pmatrix} a & b & c \\ d & e & f \\ g & h & i \end{pmatrix} = \begin{bmatrix} a & a & a \\ a & a & a \\ a & a & a \end{bmatrix},$$

$T_4: V_4 \to V_4$ is given by

$$T_4 \begin{pmatrix} a & a & a & a \\ a & a & a & a \end{pmatrix} = \begin{pmatrix} 2a & 2a & 2a & 2a \\ 2a & 2a & 2a & 2a \end{pmatrix}$$

and $T_5 : V_5 \to V_5$ is defined by

$$T_5 \begin{pmatrix} a & 0 & 0 & 0 \\ b & c & 0 & 0 \\ d & e & f & 0 \\ g & h & i & j \end{pmatrix} = \begin{pmatrix} a & 0 & 0 & 0 \\ a & a & 0 & 0 \\ a & a & a & 0 \\ a & a & a & a \end{pmatrix}.$$

$T = (T_1, T_2, T_3, T_4, T_5)$ is a special group set linear operator on V.

Suppose for $V = (V_1, \ldots, V_n)$ we define $T = (T_1 \ldots T_n)$ by $T_i: V_i \to V_j$, $(i \neq j)$; $1 \leq i, j \leq n$ i.e., no two $V_i$'s are mapped on to the same $V_j$ and $V_i$ is mapped onto $V_j$ with $j \neq i$ then we define such T to be a quasi special group set linear operator on V.

We illustrate this by a few examples.

***Example 3.2.53:*** Let $V = (V_1, V_2, V_3, V_4)$ where

$$V_1 = \left\{ \begin{pmatrix} a & b \\ c & d \end{pmatrix} \middle| a, b, c, d \in Z \right\},$$

$V_2 = \{Z[x]$ all polynomials of degree less than or equal to 3$\}$,

$$V_3 = \left\{ \begin{pmatrix} a & a \\ a & a \\ a & a \\ a & a \\ a & a \end{pmatrix} \middle| a \in Z \right\}$$

and

$$V_4 = \left\{ \begin{pmatrix} a & a & a & a & a \\ a & a & a & a & a \end{pmatrix} \middle| a \in Z \right\}$$

be the special group set linear operator on the group Z. Take $T = (T_1, T_2, T_3, T_4)$ where $T_1: V_1 \to V_2$, $T_2: V_2 \to V_1$, $T_3: V_3 \to V_4$ and $T_4: V_4 \to V_3$ defined by

$$T_1 \begin{pmatrix} a & b \\ c & d \end{pmatrix} = \{a + bx + cx^2 + dx^3\},$$

$$T_2 (a + bx + cx^2 + dx^3) = \begin{pmatrix} a & b \\ c & d \end{pmatrix},$$

$$T_3 \begin{pmatrix} a & a \\ a & a \\ a & a \\ a & a \\ a & a \end{pmatrix} = \begin{pmatrix} a & a & a & a & a \\ a & a & a & a & a \end{pmatrix}$$

and

$$T_4 \begin{pmatrix} a & a & a & a & a \\ a & a & a & a & a \end{pmatrix} = \begin{pmatrix} a & a \\ a & a \\ a & a \\ a & a \\ a & a \end{pmatrix}.$$

It is easily verified that T is a quasi special group set linear operator on V.

**Example 3.2.54:** Let $V = (V_1, V_2, V_3, V_4, V_5)$ where

$$V_1 = \left\{ \begin{pmatrix} a & b \\ c & d \end{pmatrix} \middle| a,b,c,d \in Z_7 \right\},$$

$V_2 = \{Z_7[x]$ all polynomials of degree less than or equal to 5 with coefficients from $Z_7$ in the indeterminate $x\}$;

$$V_3 = \left\{ \begin{pmatrix} a_0 & a_1 & a_2 \\ a_3 & a_4 & a_5 \end{pmatrix} \middle| a_i \in Z_7; 0 \le i \le 5 \right\},$$

$$V_4 = \left\{ \begin{pmatrix} a & b & g \\ b & p & d \\ g & d & c \end{pmatrix} \middle| p, a, b, c, d, g \in Z_7 \right\}$$

and $V_5 = Z_7 \times Z_7 \times Z_7 \times Z_7$ be a special group set vector space over the group $Z_7$.

Define $T = (T_1, T_2, T_3, T_4, T_5) : V \to V$ by

$$T_1 : V_1 \to V_5,$$
$$T_2 : V_2 \to V_3,$$
$$T_3 : V_3 \to V_4,$$
$$T_4 : V_4 \to V_2$$

and $\quad T_5 : V_5 \to V_1$

where

$$T_1 \begin{pmatrix} a & b \\ c & d \end{pmatrix} = (a\ b\ c\ d).$$

$$T_2 (a_0 + a_1 x + a_2 x^2 + a_3 x^3 + a_4 x^4 + a_5 x^5) = \begin{pmatrix} a_0 & a_1 & a_2 \\ a_3 & a_4 & a_5 \end{pmatrix}.$$

$$T_3 \begin{pmatrix} a_0 & a_1 & a_2 \\ a_3 & a_4 & a_5 \end{pmatrix} = \begin{pmatrix} a_0 & a_1 & a_2 \\ a_1 & a_3 & a_4 \\ a_2 & a_4 & a_5 \end{pmatrix}.$$

$$T_4 \begin{pmatrix} a` & b & g \\ b & p & d \\ g & d & c \end{pmatrix} = (a + bx + gx^2 + px^3 + dx^4 + cx^5)$$

and

$$T_5(a\ b\ c\ d) = \begin{pmatrix} a & b \\ c & d \end{pmatrix}.$$

Clearly $T = (T_1, T_2, T_3, T_4, T_5)$ is a quasi special group linear operator on V.

Now we proceed onto define the notion of inverse of a special group linear operator on V.

**DEFINITION 3.2.28:** *Let $V = (V_1, V_2, ..., V_n)$ be a special group set vector space over the group G. Let $T = (T_1, ..., T_n)$ be a special group set linear operator on V for each i, i.e., $T_i: V_i \to V_i$ if their exists for each $T_i$, a group linear operator $T_i^{-1}: V_i \to V_i$ such that $T_i^{-1} o\, T_i = T_i\, o\, T_i^{-1} = I_i$, $I_i$ is the group identity map on $V_i$. This is true for i=1, 2, ..., n. If $T^{-1} = \left(T_1^{-1}, T_2^{-1}, ..., T_n^{-1}\right)$ then we call $T^{-1}$ the special group set inverse linear operator on V; and $T\, o\, T^{-1} = T^{-1}\, o\, T = (I_1, I_2, ..., I_n)$.*

We illustrate this by a simple example.

*Example 3.2.55:* Let $V = (V_1, V_2, V_3, V_4)$ where

$$V_1 = \left\{ \begin{pmatrix} a & b \\ c & d \end{pmatrix} \middle| a,b,c,d \in Z \right\},$$

$$V_2 = Z \times Z \times Z \times Z \times Z,$$

$$V_3 = \left\{ \begin{pmatrix} a & a & a & a & a \\ b & b & b & b & b \end{pmatrix} \middle| a,b \in Z \right\}$$

and

$V_4 = $ {all polynomials in the variable x with coefficients from Z of degree less than or equal to 6}.

V is a special group set vector space over Z. Take $T = (T_1, T_2, T_3, T_4)$ defined by

$T_1 : V_1 \to V_1$ is defined by

$$T_1 \begin{pmatrix} a & b \\ c & d \end{pmatrix} = \begin{pmatrix} d & c \\ b & a \end{pmatrix}.$$

$T_2 : V_2 \to V_2$ is given by
$$T_2 \,(a\ b\ c\ d\ e) = (e\ d\ c\ b\ a).$$

$T_3 : V_3 \to V_3$ is defined by

$$T_3 \begin{pmatrix} a & a & a & a & a \\ b & b & b & b & b \end{pmatrix} = \begin{pmatrix} b & b & b & b & b \\ a & a & a & a & a \end{pmatrix}$$

and $T_4 : V_4 \to V_4$ is defined by

$$T_4\{a_0 + a_1x + a_2x^2 + a_3x^3 + a_4x^4 + a_5x^5 + a_6x^6\}$$
$$= \{a_0x^6 + a_1x^5 + a_2x^3 + a_3x^3 + a_4x^2 + a_5x + a_6\}.$$

It is easily verified that $T = (T_1, T_2, T_3, T_4)$ is a special group linear operator on V. Define $T^{-1} = \left(T_1^{-1}, T_2^{-1}, T_3^{-1}, T_4^{-1}\right)$ on V by

$T_1^{-1} : V_1 \to V_1$ defined by
$$T_1^{-1}\begin{pmatrix} a & b \\ c & d \end{pmatrix} = \begin{pmatrix} d & c \\ b & a \end{pmatrix}.$$

$T_2^{-1} : V_2 \to V_2$ is given by
$$T_2^{-1}(a\ b\ c\ d\ e) = (e\ d\ c\ b\ a).$$

$T_3^{-1} : V_3 \to V_3$ is defined by

$$T_3^{-1}\begin{pmatrix} a & a & a & a & a \\ b & b & b & b & b \end{pmatrix} = \begin{pmatrix} b & b & b & b & b \\ a & a & a & a & a \end{pmatrix}$$

and $T_4^{-1} : V_4 \to V_4$ is given by

$$T_4^{-1}(a_0 + a_1x + a_2x^2 + a_3x^3 + a_4x^4 + a_5x^5 + a_6x^6)$$
$$= (a_0x^6 + a_1x^5 + a_2x^3 + a_3x^3 + a_4x^2 + a_5x + a_6).$$

Now

$$T_1 \circ T_1^{-1} \begin{pmatrix} a & b \\ c & d \end{pmatrix} = T_1 \begin{pmatrix} d & c \\ b & a \end{pmatrix} = \begin{pmatrix} a & b \\ c & d \end{pmatrix}.$$

$$T_1^{-1} \circ T_1 \begin{pmatrix} a & b \\ c & d \end{pmatrix} = T_1^{-1} \begin{pmatrix} d & c \\ b & a \end{pmatrix} = \begin{pmatrix} a & b \\ c & d \end{pmatrix}.$$

Thus $T_1 \circ T_1^{-1} = T_1^{-1} \circ T_1 = I_1$ (identity on $V_1$)

$$T_2 \circ T_2^{-1} (a\ b\ c\ d\ e) = T_2(e\ d\ c\ b\ a) = (a\ b\ c\ d\ e).$$
$$T_2^{-1} \circ T_2 (a\ b\ c\ d\ e) = T_2^{-1}(e\ d\ c\ b\ a) = (a\ b\ c\ d\ e).$$

Thus $T_2^{-1} \circ T_2 = T_2 \circ T_2^{-1} = I_2$ (identity on $V_2$).

$$T_3^{-1} \circ T_3 \begin{pmatrix} a & a & a & a & a \\ b & b & b & b & b \end{pmatrix} = T_3^{-1} \begin{pmatrix} b & b & b & b & b \\ a & a & a & a & a \end{pmatrix}$$
$$= \begin{pmatrix} a & a & a & a & a \\ b & b & b & b & b \end{pmatrix}.$$

$$T_3 \circ T_3^{-1} \begin{pmatrix} a & a & a & a & a \\ b & b & b & b & b \end{pmatrix} = T_3 \begin{pmatrix} b & b & b & b & b \\ a & a & a & a & a \end{pmatrix}$$
$$= \begin{pmatrix} a & a & a & a & a \\ b & b & b & b & b \end{pmatrix}.$$

Thus $T_3 \circ T_3^{-1} = T_3^{-1} \circ T_3 = I_3$ (identity on $V_3$).

$$T_4 \circ T_4^{-1} (a_0 + a_1x + a_2x^2 + a_3x^3 + a_4x^4 + a_5x^5 + a_6x^6)$$
$$= T_4 (a_0x^6 + a_1x^5 + a_2x^3 + a_3x^3 + a_4x^2 + a_5x + a_6)$$
$$= (a_0 + a_1x + a_2x^2 + a_3x^3 + a_4x^4 + a_5x^5 + a_6x^6).$$

Now

$$T_4^{-1} \circ T_4 (a_0 + a_1x + a_2x^2 + a_3x^3 + a_4x^4 + a_5x^5 + a_6x^6)$$
$$= T_4^{-1} (a_0x^6 + a_1x^5 + a_2x^3 + a_3x^3 + a_4x^2 + a_5x + a_6)$$
$$= (a_0 + a_1x + a_2x^2 + a_3x^3 + a_4x^4 + a_5x^5 + a_6x^6).$$

Thus $T_4^{-1} \circ T_4 = T_4 \circ T_4^{-1} = I_4$ (identity on $V_4$).

Thus $T^{-1} = \left(T_1^{-1}, T_2^{-1}, T_3^{-1}, T_4^{-1}\right)$ is the special group set inverse linear operator on V.

Now we proceed onto define the quasi special group set inverse linear operator on V.

**DEFINITION 3.2.29:** *Let $V = (V_1, V_2, ..., V_n)$ be a special group set vector space over a set G. Let $T = (T_1, ..., T_n)$ be a quasi special set inverse linear operator on V. We call $T^{-1} = \left(T_1^{-1}, T_2^{-1}, ..., T_n^{-1}\right)$ to be a quasi special set inverse linear operator of T on V if $T \circ T^{-1} = T^{-1} \circ T = (I_1, ..., I_n)$; i.e., each $T_i^{-1} \circ T_i = T_i \circ T_i^{-1} = I_i$ for $i = 1, 2, ..., n$.*

We illustrate this situation by an example.

***Example 3.2.56:*** Let $V = (V_1, V_2, V_3, V_4)$ where

$$V_1 = Z \times Z \times Z \times Z,$$

$$V_2 = \left\{ \begin{pmatrix} a & b \\ c & d \end{pmatrix} \middle| a, b, c, d \in Z \right\},$$

$$V_3 = \left\{ \begin{pmatrix} a_1 & a_2 & a_3 \\ a_4 & a_5 & a_6 \\ a_7 & a_8 & a_9 \end{pmatrix} \middle| a_i \in Z \right\}$$

and $V_4 = \{Z[x] \mid$ all polynomials of degree less than or equal to 8 with coefficients from Z in the variable x$\}$. Clearly V is a special group set vector space over the set Z. Now define $T = (T_1, T_2, T_3, T_4)$ on V by

$T_1: V_1 \to V_2$, is such that

$$T_1 (a\ b\ c\ d) = \begin{pmatrix} a & b \\ c & d \end{pmatrix}.$$

$T_2 : V_2 \to V_1$ is defined by

$$T \begin{pmatrix} a & b \\ c & d \end{pmatrix} = (a\ b\ c\ d).$$

$T_3 : V_3 \to V_4$ is defined by

$$T_3 \begin{pmatrix} a_1 & a_2 & a_3 \\ a_4 & a_5 & a_6 \\ a_7 & a_8 & a_9 \end{pmatrix} = (a_1 + a_2 x + a_3 x^2 + a_4 x^3 + a_5 x^4 + a_6 x^5 + a_7 x^6 + a_8 x^7 + a_9 x^8)$$

and $T_4 : V_4 \to V_3$ is defined by

$$T_4 (a_1 + a_2 x + \ldots + a_9 x^8) = \begin{pmatrix} a_1 & a_2 & a_3 \\ a_4 & a_5 & a_6 \\ a_7 & a_8 & a_9 \end{pmatrix}.$$

It is easily verified that T is a special group set linear operator on V. Define $T^{-1} = (T_1^{-1}, T_2^{-1}, T_3^{-1}, T_4^{-1})$ from V to V as follows:

$T_1^{-1} : V_2 \to V_1$ is defined by

$$T_1^{-1} \begin{pmatrix} a & b \\ c & d \end{pmatrix} = (a\ b\ c\ d).$$

$T_2^{-1} : V_1 \to V_2$ defined by

$$T_2^{-1} (a\ b\ c\ d) = \begin{pmatrix} a & b \\ c & d \end{pmatrix}.$$

Now

$$T_1 \circ T_1^{-1} \begin{pmatrix} a & b \\ c & d \end{pmatrix} = T_1 \left( T^{-1} \begin{pmatrix} a & b \\ c & d \end{pmatrix} \right).$$

$$= T_1 \,(a\ b\ c\ d) = \begin{pmatrix} a & b \\ c & d \end{pmatrix}.$$

That is $T_1 \circ T_1^{-1} = I_2 : V_2 \to V_2$

$$T_1^{-1} \circ T_1 = T_1^{-1} \,(T_1 \,(a\ b\ c\ d)) = T_1^{-1} \begin{pmatrix} a & b \\ c & d \end{pmatrix} = (a\ b\ c\ d).$$

So $T_1^{-1} \circ T_1 = I_1 : V_1 \to V_1$. Consider

$$T_2 \circ T_2^{-1} \,(a\ b\ c\ d) = T_2 \begin{pmatrix} a & b \\ c & d \end{pmatrix} = (a\ b\ c\ d)$$

That is $T_2 \circ T_2^{-1} = I_1 : V_1 \to V_1$. Now

$$T_2^{-1} \circ T_2 \begin{pmatrix} a & b \\ c & d \end{pmatrix} = T_2^{-1} \left( T_2 \begin{pmatrix} a & b \\ c & d \end{pmatrix} \right)$$

$$T_2^{-1} \,(a\ b\ c\ d) = \begin{pmatrix} a & b \\ c & d \end{pmatrix}$$

So $T_2^{-1} \circ T_2 = I_2 : V_2 \to V_2$. We see

$T_3^{-1} : V_4 \to V_3$ defined by

$$T_3^{-1} \,(a_1 + a_2 x + \ldots + a_8 x^9) = \begin{pmatrix} a_1 & a_2 & a_3 \\ a_4 & a_5 & a_6 \\ a_7 & a_8 & a_9 \end{pmatrix}$$

and $T_4^{-1} : V_3 \to V_4$ is defined by

$$T_4^{-1} \begin{pmatrix} a_1 & a_2 & a_3 \\ a_4 & a_5 & a_6 \\ a_7 & a_8 & a_9 \end{pmatrix} = (a_1 + a_2 x + a_3 x^2 + \ldots + a_8 x^7 + a_9 x^8).$$

Clearly

$$T_3 \circ T_3^{-1} (a_1 + a_2x + \ldots + a_8x^9)$$

$$= T_3 \begin{pmatrix} a_1 & a_2 & a_3 \\ a_4 & a_5 & a_6 \\ a_7 & a_8 & a_9 \end{pmatrix} = a_1 + a_2x + \ldots + a_8x^9.$$

Thus $T_3 \circ T_3^{-1} : V_4 \to V_4$ and $T_3 \circ T_3^{-1} = I_4$ on $V_4$

$$T_3^{-1} \circ T_3 \begin{pmatrix} a_1 & a_2 & a_3 \\ a_4 & a_5 & a_6 \\ a_7 & a_8 & a_9 \end{pmatrix} = T_3^{-1} \left( T_3 \begin{pmatrix} a_1 & a_2 & a_3 \\ a_4 & a_5 & a_6 \\ a_7 & a_8 & a_9 \end{pmatrix} \right)$$

$$= T_3^{-1} (a_1 + a_2x + \ldots + a_8x^7 + a_9x^8) \begin{pmatrix} a_1 & a_2 & a_3 \\ a_4 & a_5 & a_6 \\ a_7 & a_8 & a_9 \end{pmatrix}.$$

Thus $T_3^{-1} \circ T_3: V_3 \to V_3$ and $T_3^{-1} \circ T_3 = I_3$. Finally

$$T_4^{-1} \circ T_4 (a_1 + a_2x + \ldots + a_8x^7 + a_9x^8)$$

$$T_4^{-1} \begin{pmatrix} a_1 & a_2 & a_3 \\ a_4 & a_5 & a_6 \\ a_7 & a_8 & a_9 \end{pmatrix} = a_1 + a_2x + \ldots + a_8x^7 + a_9x^8.$$

Thus $T_4^{-1} \circ T_4 : V_4 \to V_4$ i.e., $T_4^{-1} \circ T_4 = I_4 : V_4 \to V_4$. But

$$T_4 \circ T_4^{-1} \begin{pmatrix} a_1 & a_2 & a_3 \\ a_4 & a_5 & a_6 \\ a_7 & a_8 & a_9 \end{pmatrix} = T_4 \left( T_4^{-1} \begin{pmatrix} a_1 & a_2 & a_3 \\ a_4 & a_5 & a_6 \\ a_7 & a_8 & a_9 \end{pmatrix} \right)$$

$$= T_4 (a_1 + a_2x + \ldots + a_8x^7 + a_9x^8) = \begin{pmatrix} a_1 & a_2 & a_3 \\ a_4 & a_5 & a_6 \\ a_7 & a_8 & a_9 \end{pmatrix}.$$

Thus $T_4 \circ T_4^{-1} : V_3 \to V_3$ is such that $T_4 \circ T_4^{-1} = I_3 : V_3 \to V_3$. Thus we have verified $T^{-1} = \left( T_1^{-1}, T_2^{-1}, T_3^{-1}, T_4^{-1} \right)$ is the special group set quasi inverse linear operator of T. Thus
$$T \circ T^{-1} = T^{-1} \circ T = (I_1, I_2, I_3, I_4) = I.$$

On similar lines one can always define the notion of inverse for any special group set linear operator on a special group set vector space as well as inverse of special group set linear transformation on a special group set vector spaces V and W defined on the same set S.

We define the notion of direct sum in case of special group set vector spaces defined over a set S.

**DEFINITION 3.2.30:** *Let $V = (V_1, \ldots, V_n)$ be a special group set vector space over the set S. We say V is a direct sum if each $V_i$ can be represented as a direct sum of subspaces $i = 1, 2, \ldots, n$. that is*
$$V = (W_1^1 \oplus \ldots \oplus W_{t_1}^1, W_1^2 \oplus \ldots \oplus W_{t_2}^2, \ldots, W_1^n \oplus \ldots \oplus W_{t_n}^n)$$
*where each $V_i = \left( W_1^i \oplus W_2^i \oplus \ldots \oplus W_{t_i}^i \right)$ is a direct sum, i.e., each $v_i \in V_i$ can be represented uniquely as a sum of elements from $\left( W_1^i, W_2^i, \ldots, W_{t_i}^i \right)$ and $W_p^i \cap W_q^i = (0)$ if $p \neq q$. This is true for each i, $i = 1, 2, \ldots, n$.*

We represent this by some examples.

***Example 3.2.57:*** Let $V = (V_1, V_2, V_3, V_4)$ where

$$V_1 = \left\{ \begin{pmatrix} a & b \\ c & d \end{pmatrix} \middle| a, b, c, d \in Z \right\},$$

$V_2 = \{Z[x]$; all polynomials of degree less than or equal to $4\}$, $V_3 = Z \times Z \times Z$

and

$$V_4 = \left\{ \begin{pmatrix} a_1 & a_2 \\ a_3 & a_4 \\ a_5 & a_6 \\ a_7 & a_8 \\ a_9 & a_{10} \end{pmatrix} \middle| a_i \in Z : 1 \leq i \leq 10 \right\}$$

be a special group set vector space over the group Z.
Let

$$W_1^1 = \left\{ \begin{pmatrix} a & 0 \\ 0 & d \end{pmatrix} \middle| a, d \in Z \right\}$$

and

$$W_2^1 = \left\{ \begin{pmatrix} 0 & e \\ g & 0 \end{pmatrix} \middle| e, g \in Z \right\}$$

be a special group vector subspaces of $V_1$. We see

$$V_1 = W_1^1 \oplus W_2^1 \text{ and } W_1^1 \cap W_2^1 = \begin{pmatrix} 0 & 0 \\ 0 & 0 \end{pmatrix}.$$

Consider $W_1^2 = \{$all polynomials of degree one or three with coefficients from Z i.e., $ax + bx^3$; $a, b \in Z\} \subseteq V_2$. $W_1^2$ is a special group vector subspace of $V_2$. $W_2^2 = \{$all polynomials of the form $a_0 + a_1 x^2 + a_2 x^4 \mid a_0, a_1, a_2 \in Z\} \subseteq V_2$ is also a special group vector subspace of $V_2$. We see $V_2 = W_1^2 \oplus W_2^2$ and $W_1^2 \cap W_2^2 = (0)$. $W_1^3 = \{Z \times \{0\} \times Z\} \subseteq V_3$ is a special group vector subspace of $V_3$ and
$$W_2^3 = \{\{0\} \times Z \times \{0\}\} \subseteq V_3$$
is a special group vector subspace of $V_3$.

Further $V_3 = W_1^3 \oplus W_2^3$ and $W_1^3 \cap W_2^3 = (0)$.

Take

$$W_1^4 = \left\{ \begin{pmatrix} a_1 & a_2 \\ a_3 & a_4 \\ 0 & 0 \\ 0 & 0 \\ 0 & 0 \end{pmatrix} \middle| a_i \in Z; 1 \leq i \leq 4 \right\} \subseteq V_4$$

is a special group set vector space over Z. Take

$$W_2^4 = \left\{ \begin{pmatrix} 0 & 0 \\ 0 & 0 \\ a_1 & a_2 \\ a_3 & a_4 \\ a_5 & a_6 \end{pmatrix} \middle| a_i \in Z; 1 \leq i \leq 6 \right\} \subseteq V_4$$

is a special group set vector space over Z. For $V_4 = W_1^4 \oplus W_2^4$ and

$$W_1^4 \cap W_2^4 = \begin{bmatrix} 0 & 0 \\ 0 & 0 \\ 0 & 0 \\ 0 & 0 \\ 0 & 0 \end{bmatrix}.$$

Now we show this representation is not unique. For take

$$W_1^1 = \begin{pmatrix} a & 0 \\ 0 & 0 \end{pmatrix} \subseteq V_1,$$

$$W_2^1 = \begin{pmatrix} 0 & b \\ 0 & 0 \end{pmatrix} \subseteq V_1,$$

$$W_3^1 = \begin{pmatrix} 0 & 0 \\ c & 0 \end{pmatrix} \subseteq V_1$$

and
$$W_4^1 = \begin{pmatrix} 0 & 0 \\ 0 & d \end{pmatrix} \subseteq V_1.$$

$W_1^1, W_2^1, W_3^1$ and $W_4^1$ are special group subspaces of $V_1$ and
$$V_1 = W_1^1 \oplus W_2^1 \oplus W_3^1 \oplus W_4^1$$
with
$$W_j^1 \cap W_i^1 = \begin{pmatrix} 0 & 0 \\ 0 & 0 \end{pmatrix} \text{ if } i \neq j.$$

Take the space $V_2$, let
$$W_1^2 = \{a + bx \mid a, b \in Z\} \subseteq V_2,$$
$$W_2^2 = \{ax^2 + cx^3 \mid a, c \in Z\} \subseteq V_2$$
and
$$W_3^2 = \{dx_4 \mid d \in Z\} \subseteq V_2.$$

It is easily verified $W_1^2$, $W_2^2$ and $W_3^2$ are special group vector subspaces of $V_2$.
Further
$$V_2 = W_1^2 \oplus W_2^2 \oplus W_3^2$$
and $W_i^2 \cap W_j^2 = \{0\}$ if $i \neq j$. $1 \leq j, i \leq 3$. So $V_2$ is a direct sum of $W_1^2$, $W_2^2$ and $W_3^2$.

Now consider $V_3$,
$$W_1^3 = Z \times \{0\} \times \{0\} \subseteq V_3,$$
$$W_2^3 = \{0\} \times Z \times \{0\} \subseteq V_3 \text{ and}$$
$$W_3^3 = \{0\} \times \{0\} \times Z \subseteq V_3$$

are special group vector subspaces of $V_3$ such that
$$V_3 = W_1^3 \oplus W_2^3 \oplus W_3^3$$
with
$$W_i^3 \cap W_j^3 = \{0\} \times \{0\} \times \{0\},$$
if $i \neq j$; $1 \leq i, j \leq 3$.

Consider $V_4$, take

$$W_1^4 = \left\{ \begin{pmatrix} a_1 & a_2 \\ 0 & 0 \\ a_3 & a_4 \\ 0 & 0 \\ 0 & 0 \end{pmatrix} \middle| \begin{array}{l} a_i \in Z \\ 1 \le i \le 4 \end{array} \right\} \subseteq V_4,$$

$$W_2^4 = \left\{ \begin{pmatrix} 0 & 0 \\ a_1 & a_2 \\ 0 & 0 \\ 0 & 0 \\ 0 & 0 \end{pmatrix} \middle| a_1, a_2 \in Z \right\} \subseteq V_4$$

and $W_3^4 = \left\{ \begin{pmatrix} 0 & 0 \\ 0 & 0 \\ 0 & 0 \\ a_1 & a_2 \\ a_3 & a_4 \end{pmatrix} \middle| \begin{array}{l} a_i \in Z \\ 1 \le i \le 4 \end{array} \right\}.$

We see $V_4 = W_1^4 \oplus W_2^4 \oplus W_3^4$ with

$$W_i^4 \cap W_j^4 = \begin{bmatrix} 0 & 0 \\ 0 & 0 \\ 0 & 0 \\ 0 & 0 \\ 0 & 0 \end{bmatrix} \text{ if } i \ne j; \ 1 \le \underline{i}, j \le 3.$$

Thus we see $V = (V_1, V_2, V_3, V_4)$ is the direct sum of

$$V = (W_1^1 \oplus W_2^1 \oplus W_3^1 \oplus W_4^1, W_1^2 \oplus W_2^2 \oplus W_3^2,$$

$$W_1^3 \oplus W_2^3 \oplus W_3^3 , \; W_1^4 \oplus W_2^4 \oplus W_3^4 ) \qquad \ldots \quad I$$

Thus the representation of the special group set vector spaces as a direct sum of special group set vector subspaces is not unique. Further one of the uses of finding the direct sum, is we can get many sets of special group set vector subspaces of V. For instance in the representation of V given in I. We have 108 proper special group set vector subspaces of V.

$$W_1 = (W_1^1, W_2^1, W_3^1, W_4^1),$$
$$W_2 = (W_1^2, W_2^2, W_3^2)$$
$$W_3 = (W_1^3, W_2^3, W_3^3) \text{ and}$$
$$W_4 = (W_1^4, W_2^4, W_3^4).$$

Now having seen the direct sum and special group set vector subspaces we can define special group set projection operators and special group set idempotent operator.

Let $V = (V_1, V_2, \ldots, V_n)$ be a special group set vector space over a set S. If $T = (T_1, \ldots, T_n)$ is a special group set linear operator on V we call T to be a special group set linear idempotent operator on V if $T \circ T = T$, i.e., $T_i \circ T_i = T_i$ for every $i = 1, 2, \ldots, n$.

Now for us to define the notion of special group set projection linear operator on V we in the first place need, some proper special group set vector subspace of V. Let $W = (W_1, W_2, \ldots, W_n)$ be a special group set vector subspace of V, i.e., each $W_i \subseteq V_i$ is a special group vector subspace of V; $i = 1, 2, \ldots, n$; i.e., $W = (W_1, W_2, \ldots, W_n) \subseteq (V_1, V_2, \ldots, V_n)$. Let $P = (P_1, P_2, \ldots, P_n)$ be a special group set linear operator on V such that $P: V \to V$, that is $P_i: V_i \to V_i$ is such that $P_i: V_i \to W_i \subseteq V_i$; $i = 1, 2, \ldots, n$; that is $P_i$ is a special group projection operator on $V_i$ for each i and $P_i \circ P_i = P_i$ for each i, $1 \leq i \leq n$. We call $P = (P_1, P_2, \ldots, P_n)$ to be the special group set projection operator of V. Clearly $P \circ P = P$.

Now we illustrate this by an example.

**Example 3.2.58:** Let $V = (V_1, V_2, V_3, V_4, V_5)$ be a special group set vector space over the set Z where $V_1 = Z \times Z$,

$$V_2 = \left\{ \begin{pmatrix} a \\ b \end{pmatrix} \middle| a, b \in Z \right\},$$

$$V_3 = \left\{ \begin{pmatrix} a & b \\ 0 & d \end{pmatrix} \middle| a, b, d \in Z \right\},$$

$$V_4 = \left\{ \begin{pmatrix} a & b & c \\ d & e & f \end{pmatrix} \middle| a, b, c, d, e, f \in Z \right\}$$

and $V_5 = \{Z[x]$ all polynomials of degree less than or equal to three with coefficients from $Z\}$. Take $W = (W_1, W_2, W_3, W_4, W_5) \subseteq (V_1, V_2, V_3, V_4, V_5)$ where

$$W_1 = \{(a, 0) \mid a \in Z\} \subseteq V_1,$$

$$W_2 = \left\{ \begin{pmatrix} 0 \\ b \end{pmatrix} \middle| b \in Z \right\} \subseteq V_2,$$

$$W_3 = \left\{ \begin{pmatrix} a & b \\ 0 & 0 \end{pmatrix} \middle| a, b \in Z \right\} \subseteq V_3,$$

$$W_4 = \left\{ \begin{pmatrix} a & b & c \\ 0 & 0 & 0 \end{pmatrix} \middle| a, b, c \in Z \right\} \subseteq V_4$$

and $W_5 = \{a_0 + a_1 x^3 \mid a_0, a_1 \in Z\} \subseteq V_5$. Clearly $W = (W_1, W_2, W_3, W_4, W_5)$ is a special group set vector subspace of V. Define $P = (P_1, P_2, \ldots, P_5)$ from V to V by $P_1: V_1 \to V_1$ where $P_1(a, b) = (a, 0)$. $P_2 : V_2 \to V_2$ is defined by

$$P_2 \left\{ \begin{bmatrix} a \\ b \end{bmatrix} \right\} = \begin{bmatrix} 0 \\ b \end{bmatrix}$$

$P_3: V_3 \to V_3$ is given by
$$P_3 \left\{ \begin{bmatrix} a & b \\ 0 & d \end{bmatrix} \right\} = \begin{bmatrix} a & b \\ 0 & 0 \end{bmatrix}$$

$P_4 : V_4 \to V_4$ is such that
$$P_4 \left\{ \begin{bmatrix} a & b & c \\ d & e & f \end{bmatrix} \right\} = \begin{bmatrix} a & b & c \\ 0 & 0 & 0 \end{bmatrix}$$

and $P_5: V_5 \to V_5$ is defined by
$$P_5(a_0 + a_1x + a_2x^2 + a_3x^3) = a_0 + a_3x^3.$$

Now it is easily verified that $P = (P_1, P_2, P_3, P_4, P_5)$ is a special group set linear operator on V which is a projection of V onto W. It is further evident $P \circ P = P$. For $P_1 \circ P_1(a, b) = P_1(a, 0) = (a, 0)$; i.e., $P_1 \circ P_1 = P_1$.

$$P_2 \circ P_2 \left\{ \begin{bmatrix} a \\ b \end{bmatrix} \right\} = P_2 \left\{ \begin{bmatrix} 0 \\ b \end{bmatrix} \right\} = \begin{bmatrix} 0 \\ b \end{bmatrix}$$

Hence $P_2 \circ P_2 = P_2$.

$$P_3 \circ P_3 \left\{ \begin{bmatrix} a & b \\ 0 & d \end{bmatrix} \right\} = P_3 \begin{bmatrix} a & b \\ 0 & 0 \end{bmatrix} = \begin{bmatrix} a & b \\ 0 & 0 \end{bmatrix}$$

So $P_3 \circ P_3 = P_3$.

$$P_4 \circ P_4 \left\{ \begin{bmatrix} a & b & c \\ d & e & f \end{bmatrix} \right\} = P_4 \begin{bmatrix} a & b & c \\ 0 & 0 & 0 \end{bmatrix} = \begin{bmatrix} a & b & c \\ 0 & 0 & 0 \end{bmatrix},$$

thus $P_4 \circ P_4 = P_4$ and $P_5 \circ P_5(a_0 + a_1x + a_2x^2 + a_3x^3) = P_5(a_0 + a_3x^3) = (a_0 + a_3x^3)$; hence $P_5 \circ P_5 = P_5$.

Thus P o P = P. Hence P is an idempotent special group set linear operator on V.

It may so happen that we have $V = (V_1, V_2, \ldots, V_n)$ to be a special group set vector space over a set S. Suppose

$$V = \left( W_1^1 \oplus \ldots \oplus W_{t_1}^1, W_1^2 \oplus \ldots \oplus W_{t_2}^2, \ldots, W_1^n \oplus \ldots \oplus W_{t_n}^n \right)$$

is the direct sum of special group set vector subspaces of V over the set S. Suppose in addition to these conditions we assume $t_1 = t_2 = \ldots = t_n = t$. Let

$$W_1 = \left( W_1^1, W_1^2, \ldots, W_1^n \right) \subseteq V, \; W_2 = \left( W_2^1, W_2^2, \ldots, W_2^n \right) \subseteq V,$$

and so on;

$$W_t = \left( W_t^1, W_t^2, \ldots, W_t^n \right) \subseteq V.$$

If $P_1, P_2, \ldots, P_n$ be projection operators of V to $W_1$, V to $W_2$, …, V to $W_t$ respectively. We have $P_i$ o $P_i = P_i$ and $P_i$ o $P_j = 0$ if $i \neq j$.

Further if

$$P_1 = \left( P_1^1, P_1^2, \ldots, P_1^n \right),$$
$$P_2 = \left( P_2^1, P_2^2, \ldots, P_2^n \right)$$

and so on …

$$P_t = \left( P_t^1, P_t^2, \ldots, P_t^n \right)$$

then $P_1 + \ldots + P_t$
$= \left( P_1^1 + P_2^1 + \ldots + P_t^1, P_1^2 + P_2^2 + \ldots + P_t^2, \ldots, P_1^n + P_2^n + \ldots + P_t^n \right)$
$= (I_1, I_2, \ldots, I_n).$

We see
$$P_i^r o P_j^r = 0 \; ; \text{ if } i \neq j, \; 1 \leq r \leq n.$$
and
$$P_i^r o P_j^r = P_i^r \text{ if } i = j, \; 1 \leq r \leq n.$$

As in case of vector spaces in the case of special group set vector spaces we would not be in a position to define precisely the notion of eigen values eigen vectors etc. We to over come this difficulty make use of subspaces and their direct sum concept. We also wish to state that we need to write each and every special group vector space as a sum of the same number of subspace that is each $V_i$ is expressed as

$$V_i = \left(W_1^i \oplus \ldots \oplus W_t^i\right)$$

for i = 1, 2, …, n then also we can have the above mentioned results to be true.

If they ($V_i$'s) are of different sums we cannot get the above mentioned results.

We shall illustrate this by some examples.

***Example 3.2.59:*** Let $V = (V_1, V_2, V_3, V_4)$ where

$$V_1 = Z \times Z \times Z,$$

$$V_2 = \left\{ \begin{pmatrix} a \\ b \end{pmatrix} \middle| a, b \in Z \right\},$$

$$V_3 = \left\{ \begin{pmatrix} a & b \\ c & d \end{pmatrix} \middle| a, b, c, d \in Z \right\}$$

and $V_4 = \{Z[x]$; all polynomials of degree less than or equal to three in the variable x with coefficients from Z$\}$ be a special group set vector space over Z. Take the substructures in $V_1$ as

$$W_1^1 = Z \times Z \times \{0\} \subseteq V_1$$

and

$$W_2^1 = \{0\} \times \{0\} \times Z \subseteq V_1,$$
$$V_1 = W_1^1 \oplus W_2^1,$$

where $W_1^1$ and $W_2^1$ are special group set vector subspaces of $V_1$. Consider $V_2$, take

$$W_1^2 = \left\{ \begin{pmatrix} a \\ 0 \end{pmatrix} \middle| a \in Z \right\} \subseteq V_2$$

and

$$W_2^2 = \left\{ \begin{pmatrix} 0 \\ b \end{pmatrix} \middle| b \in Z \right\}$$

as special group set vector subspaces of $V_2$.

$$V_2 = W_1^2 \oplus W_2^2$$

i.e.,

$$W_1^2 \cap W_2^2 = \begin{pmatrix} 0 \\ 0 \end{pmatrix}.$$

In $V_3$ define

$$W_1^3 = \left\{ \begin{pmatrix} a & b \\ 0 & 0 \end{pmatrix} \middle| a, b \in Z \right\} \subseteq V_3$$

and

$$W_2^3 = \left\{ \begin{pmatrix} 0 & 0 \\ c & d \end{pmatrix} \middle| c, d \in Z \right\} \subseteq V_3.$$

$W_1^3$ and $W_2^3$ are special group set vector subspaces of $V_3$ and $V_3 = W_1^3 \oplus W_2^3$. In the special group set vector space $V_4$ take

$$W_1^4 = \{a_0 + a_2 x^2 \mid a_0, a_2 \in Z\} \subseteq V_4$$

and

$$W_2^4 = \{a_1 x + a_3 x^3 \mid a_1, a_3 \in Z\} \subseteq V_4.$$

Clearly $V_4 = W_1^4 \oplus W_2^4$ and $W_1^4 \cap W_2^4 = \{0\}$.

Thus
$$V = (V_1, V_2, V_3, V_4)$$
$$= (W_1^1 \oplus W_2^1, W_1^2 \oplus W_2^2, W_1^3 \oplus W_2^3, W_1^4 \oplus W_2^4).$$

Take $W_1 = (W_1^1, W_1^2, W_1^3, W_1^4)$ and $W_2 = (W_2^1, W_2^2, W_2^3, W_2^4)$.

Let $P_1 : V$ into $W_1$, a projection; i.e., a special group set linear operator on V.
$P = \left( P_1^1, P_1^2, P_1^3, P_1^4 \right)$ defined by

$P_1^1 : V_1 \to W_1$ is such that

$$P_1^1 \ (a\ b\ c) = (a, b, 0)$$
$$P_1^1 \ o \ P_1^1 \ (a\ b\ c)$$
$$P_1^1 \ (a\ b\ 0) = (a\ b\ 0)$$

that is $P_1^1 \ o \ P_1^1 = P_1^1$.

Consider $P_1^2 : V_2 \to W_2$ is defined by

$$P_1^2 \begin{bmatrix} a \\ b \end{bmatrix} = \begin{bmatrix} a \\ 0 \end{bmatrix}$$

we see

$$P_1^2 \text{ o } P_1^2 \begin{bmatrix} a \\ b \end{bmatrix} = P_1^2 \begin{bmatrix} a \\ 0 \end{bmatrix} = \begin{bmatrix} a \\ 0 \end{bmatrix}$$

$$P_1^2 \text{ o } P_1^2 = P_1^2.$$

Now $P_1^3 : V_3 \to W_3$ is defined by

$$P_1^3 \begin{bmatrix} a & b \\ c & d \end{bmatrix} = \begin{bmatrix} a & b \\ 0 & 0 \end{bmatrix}$$

also

$$P_1^3 \text{ o } P_1^3 \begin{bmatrix} a & b \\ c & d \end{bmatrix} = P_1^3 \begin{bmatrix} a & b \\ 0 & 0 \end{bmatrix} = \begin{bmatrix} a & b \\ 0 & 0 \end{bmatrix}$$

$$P_1^3 \text{ o } P_1^3 = P_1^3.$$

Finally $P_1^4 : V_4 \to W_4$ is defined by

$$P_1^4 (a_0 + a_1x + a_2x^2 + a_3x^3) = (a_0 + a_2x^2).$$

Thus

$$P_1^4 \text{ o } P_1^4 (a_0 + a_1x + a_2x^2 + a_3x^3) = P_1^4 (a_0 + a_2x^2)$$
$$= a_0 + a_2x^2.$$

Thus $P_1^4 \text{ o } P_1^4 = P_1^4$. So we see $P_1 \text{ o } P_1 = P_1$,

$$P_1 : V \to (W_1^1, W_1^2, W_1^3, W_1^4).$$

Now take $P_2 : V = (W_2^1, W_2^2, W_2^3, W_2^4)$ where

$$P_2 = (P_2^1, P_2^2, P_2^3, P_2^4) : (V_1, V_2, V_3, V_4) \to (W_2^1, W_2^2, W_2^3, W_2^4);$$

defined by
$$P_2^1 : V_1 \to W_2^1; \; P_2^1 (a \; b \; c) = (0 \; 0 \; c)$$
so
$$P_2^1 \circ P_2^1 (a \; b \; c) = P_2^1 (0 \; 0 \; c) = (0 \; 0 \; c).$$
Thus
$$P_2^1 \circ P_2^1 = P_2^1.$$
Consider
$$P_1^1 \circ P_2^1 (a \; b \; c) = P_1^1 (0 \; 0 \; c) = (0 \; 0 \; 0).$$
Also
$$P_2^1 \circ P_1^1 (a \; b \; c) = P_2^1 (a \; b \; 0) = (0 \; 0 \; 0).$$
Thus we get
$$P_2^1 \circ P_1^1 = P_1^1 \circ P_2^1 = \theta_1 : V_1 \to W_1.$$

Now consider $P_2^2 : V_2 \to W_2^2$; defined by

$$P_2^2 \begin{bmatrix} a \\ b \end{bmatrix} = \begin{bmatrix} 0 \\ b \end{bmatrix}$$

so

$$P_2^2 \circ P_2^2 \begin{bmatrix} a \\ b \end{bmatrix} = P_2^2 \begin{bmatrix} 0 \\ b \end{bmatrix} = \begin{bmatrix} 0 \\ b \end{bmatrix}$$

which implies
$$P_2^2 \circ P_2^2 = P_2^2.$$
Now
$$P_2^2 \circ P_1^2 \begin{bmatrix} a \\ b \end{bmatrix} = P_2^2 \begin{bmatrix} a \\ 0 \end{bmatrix} = \begin{bmatrix} 0 \\ 0 \end{bmatrix}$$
and
$$P_1^2 \circ P_2^2 \begin{bmatrix} a \\ b \end{bmatrix} = P_1^2 \begin{bmatrix} 0 \\ b \end{bmatrix} = \begin{bmatrix} 0 \\ 0 \end{bmatrix}.$$
Thus
$$P_1^2 \circ P_2^2 = P_2^2 \circ P_1^2 = \theta_2 : V_2 \to V_2$$
that is

$$\theta_2 \begin{bmatrix} a \\ b \end{bmatrix} = \begin{bmatrix} 0 \\ 0 \end{bmatrix}.$$

$P_2^3 : V_3 \to W_2^3$ is defined by

$$P_2^3 \begin{bmatrix} a & b \\ c & d \end{bmatrix} = \begin{bmatrix} 0 & 0 \\ c & d \end{bmatrix};$$

$$P_2^3 \text{ o } P_2^3 \begin{bmatrix} a & b \\ c & d \end{bmatrix} = P_2^3 \begin{bmatrix} 0 & 0 \\ c & d \end{bmatrix} = \begin{bmatrix} 0 & 0 \\ c & d \end{bmatrix}.$$

Thus
$$P_2^3 \text{ o } P_2^3 = P_2^3.$$

Now
$$P_2^3 \text{ o } P_1^3 \begin{bmatrix} a & b \\ c & d \end{bmatrix} = P_2^3 \begin{bmatrix} a & b \\ 0 & 0 \end{bmatrix} = \begin{bmatrix} 0 & 0 \\ 0 & 0 \end{bmatrix}$$

$$P_1^3 \text{ o } P_2^3 \begin{bmatrix} a & b \\ c & d \end{bmatrix} = P_1^3 \begin{bmatrix} 0 & 0 \\ c & d \end{bmatrix} = \begin{bmatrix} 0 & 0 \\ 0 & 0 \end{bmatrix}.$$

$$P_1^3 \text{ o } P_2^3 = P_2^3 \text{ o } P_1^3 = \theta_3 : V_3 \to V_3;$$

such that
$$\theta_3 \begin{bmatrix} a & b \\ c & d \end{bmatrix} = \begin{bmatrix} 0 & 0 \\ 0 & 0 \end{bmatrix}.$$

Finally $P_2^4 : V_4 \to W_2^4$ is defined by

$$P_2^4 (a_0 + a_1x + a_2x^2 + a_3x^3) = a_1x + a_3x^3.$$
$$P_1^4 \text{ o } P_2^4 (a_0 + a_1x + a_2x^2 + a_3x^3) = P_1^4 (a_1x + a_3x^3) = 0$$
$$P_2^4 \text{ o } P_1^4 (a_0 + a_1x + a_2x^2 + a_3x^3) = P_2^4 (a_0 + a_2x^2) = 0.$$

Thus
$$P_2^4 \text{ o } P_1^4 = P_1^4 \text{ o } P_2^4 = \theta_4 : V_4 \to V_4$$

defined by
$$\theta_4 (a_0 + a_1x + a_2x^2 + a_3x^3) = 0.$$

Thus we see
$$P_1 + P_2 = \left(P_1^1 + P_2^1, P_1^2 + P_2^2, P_1^3 + P_2^3, P_1^4 + P_2^4\right) : V \to V$$

$$\begin{aligned} P_1^1 + P_2^1 \;(a\;b\;c) &= P_1^1\;(a\;b\;c) + P_2^1\;(a\;b\;c) \\ &= (a\;b\;0) + (0\;0\;c) \\ &= (a\;b\;c). \end{aligned}$$

Thus $P_1^1 + P_2^1 = I_1 : V \to V$.

$$\left(P_1^2 + P_2^2\right)\left(\begin{bmatrix} a \\ b \end{bmatrix}\right)$$

$$P_1^2 \begin{bmatrix} a \\ b \end{bmatrix} + P_2^2 \begin{bmatrix} a \\ b \end{bmatrix} = \begin{bmatrix} a \\ 0 \end{bmatrix} + \begin{bmatrix} 0 \\ b \end{bmatrix} = \begin{bmatrix} a \\ b \end{bmatrix}.$$

Thus $P_1^2 + P_2^2 = I_2 : V_2 \to V_2$ (identity special group operator on $V_2$).

Now

$$\left(P_1^3 + P_2^3\right)\left(\begin{bmatrix} a & b \\ c & d \end{bmatrix}\right) = P_1^3 \begin{bmatrix} a & b \\ c & d \end{bmatrix} + P_2^3 \begin{bmatrix} a & b \\ c & d \end{bmatrix}$$

$$= \begin{bmatrix} a & b \\ 0 & 0 \end{bmatrix} + \begin{bmatrix} 0 & 0 \\ c & d \end{bmatrix} = \begin{bmatrix} a & b \\ c & d \end{bmatrix}.$$

Hence $\left(P_1^3 + P_2^3\right) = I_3$. $I_3 : V_3 \to V_3$ (identity special group operator on $V_3$).

Finally
$$P_1^4 + P_2^4\;(a_0 + a_1x + a_2x^2 + a_3x^3)$$
$$= P_1^4\;(a_0 + a_1x + a_2x^2 + a_3x^3) + P_2^4\;(a_0 + a_1x + a_2x^2 + a_3x^3)$$
$$= a_0 + a_2x^2 + a_1x + a_3x^3 = a_0 + a_1x + a_2x^2 + a_3x^3.$$

that is $P_1^4 + P_2^4 = I_4 : V_4 \to V_4$. Thus $P_1 + P_2 = (I_1, I_2, I_3, I_4)$.

Now we give yet another example.

***Example 3.2.60:*** Let $V = (V_1, V_2, V_3)$ where $V_1 = Z_{12} \times Z_{12} \times Z_{12} \times Z_{12}$, $V_2 = \{Z_{12}[x]$ all polynomials in the variable x with coefficients from $Z_{12}$ of degree less than or equal to 6$\}$ and

$$V_3 = \left\{ \begin{pmatrix} a & b & c \\ d & e & f \\ g & h & i \end{pmatrix} \middle| a,b,c,d,e,f,g,h,i \in Z_{12} \right\}.$$

$V = (V_1, V_2, V_3)$ is a special group set vector space over the set $Z_{10}$. Now let

$$W_1 = \left(W_1^1, W_1^2, W_1^3\right) \subseteq V$$
$$W_2 = \left(W_2^1, W_2^2, W_2^3\right) \subseteq V$$
and $$W_3 = \left(W_3^1, W_3^2, W_3^3\right) \subseteq V$$

Such that
$$V = (W_1 + W_2 + W_3)$$
$$= \left(W_1^1 + W_2^1 + W_3^1, W_1^2 + W_2^2 + W_3^2, W_1^3 + W_2^3 + W_3^3\right)$$
$$= (V_1, V_2, V_3)$$

that is $V_i = \left(W_1^i + W_2^i + W_3^i\right)$, $1 \le i \le 3$. For the special group vector space $V_1$,

$$W_1^1 = Z_{12} \times Z_{12} \times \{0\} \times \{0\},$$
$$W_2^1 = \{0\} \times Z_{12} \times \{0\}$$
and $$W_3^1 = \{0\} \times \{0\} \times \{0\} \times Z_{12};$$

clearly $W_1^1$, $W_2^1$ and $W_3^1$ are special group vector subspaces of $V_1$. Here $V_1 = W_1^1 \oplus W_2^1 \oplus W_3^1$, with $W_i^1 \cap W_j^1 = (0\ 0\ 0\ 0)$ if $i \ne j$, $1 \le i, j \le 3$. Now consider the special group vector space $V_2$, take

$$W_1^2 = \{a_0 + a_1 x \mid a_0, a_1 \in Z_{12}\} \subseteq V_2,$$
$$W_2^2 = \{a_2 x^2 + a_3 x^3 \mid a_2, a_3 \in Z_{12}\} \subseteq V_2 \text{ and}$$
$$W_3^2 = \{a_4 x^4 + a_5 x^5 + a_6 x^6 \mid a_4, a_5, a_6 \in Z_{12}\} \subseteq V_2.$$

Clearly $W_1^2$, $W_2^2$ and $W_3^2$ are special group vector subspaces of $V_2$ with $V_2 = W_1^2 \oplus W_2^2 \oplus W_3^2$ where $W_i^2 \cap W_j^2 = 0$ if $i \neq j$; $1 \leq i, j \leq 3$.

Finally take

$$V_3 = \left\{ \begin{bmatrix} a & b & c \\ d & e & f \\ g & h & i \end{bmatrix} \middle| a,b,c,d,e,f,g,h,i \in Z_{12} \right\}$$

with

$$W_1^3 = \left\{ \begin{bmatrix} 0 & b & c \\ 0 & 0 & f \\ 0 & 0 & 0 \end{bmatrix} \middle| b,c,f \in Z_{12} \right\},$$

$$W_2^3 = \left\{ \begin{bmatrix} a & 0 & 0 \\ 0 & e & 0 \\ 0 & 0 & i \end{bmatrix} \middle| a,e,i \in Z_{12} \right\}$$

and

$$W_3^3 = \left\{ \begin{bmatrix} a & 0 & 0 \\ d & 0 & 0 \\ g & h & 0 \end{bmatrix} \middle| d,g,h \in Z_{12} \right\}$$

as special group vector subspaces of $V_3$.

Clearly $V_3 = W_1^3 \oplus W_2^3 \oplus W_3^3$ with

$$W_i^3 \cap W_j^3 = \begin{bmatrix} 0 & 0 & 0 \\ 0 & 0 & 0 \\ 0 & 0 & 0 \end{bmatrix}; i \neq j.\ 1 \leq i, j \leq 3.$$

$$\left\{ \begin{bmatrix} a & 0 & 0 \\ d & 0 & 0 \\ g & h & 0 \end{bmatrix} \middle| d,g,h \in Z_{12} \right\}$$

$$V = \left(W_1^1 \oplus W_2^1 \oplus W_3^1, W_1^2 \oplus W_2^2 \oplus W_3^2, W_1^3 \oplus W_2^3 \oplus W_3^3\right).$$

$$V_1 = W_1^1 \oplus W_2^1 \oplus W_3^1$$
$$V_2 = W_1^2 \oplus W_2^2 \oplus W_3^2 \text{ and } V_3 = W_1^3 \oplus W_2^3 \oplus W_3^3.$$

Define
$$P_1 = \left(P_1^1, P_1^2, P_1^3\right), \ P_2 = \left(P_2^1, P_2^2, P_2^3\right) \text{ and } P_3 = \left(P_3^1, P_3^2, P_3^3\right)$$
on V by
$$P_1 : V \to \left(W_1^1, W_1^2, W_1^3\right)$$
where
$$P_1^1 : V_1 \to W_1^1,$$
$$P_1^2 : V_2 \to W_1^2 \text{ and}$$
$$P_1^3 : V_3 \to W_1^3$$
defined by
$$P_1^1 \ [(x \ y \ z \ \omega)] = (x \ y \ 0 \ 0),$$
$$P_1^2 \ [a_0 + a_1 x + \ldots + a_6 x^6] = (a_0 + a_1 x)$$
and
$$P_1^3 \left(\begin{bmatrix} a & b & c \\ d & e & f \\ g & h & i \end{bmatrix}\right) = \begin{bmatrix} 0 & b & c \\ 0 & 0 & f \\ 0 & 0 & 0 \end{bmatrix}.$$

Thus $P_1 = \left(P_1^1, P_1^2, P_1^3\right)$ is a special group set projection operator from V into $\left(W_1^1, W_1^2, W_1^3\right)$.

Now consider $P_2 = \left(P_2^1, P_2^2, P_2^3\right)$ where
$P_2 : V \to \left(W_2^1, W_2^2, W_2^3\right)$ defined by
$P_2^1 : V_1 \to W_2^1$ as
$$P_2^1 (x \ y \ z \ \omega) = (0 \ 0 \ z \ 0).$$
$P_2^2 : V_2 \to W_2^2$ defined by
$$P_2^2 (a_0 + a_1 x + \ldots + a_6 x^6) = a_2 x^2 + a_3 x^3$$
and

$P_2^3 : V_3 \to W_2^3$ defined by

$$P_2^3 \begin{bmatrix} a & b & c \\ d & e & f \\ g & h & i \end{bmatrix} = \begin{bmatrix} a & 0 & 0 \\ 0 & e & 0 \\ 0 & 0 & i \end{bmatrix}.$$

Thus $P_2 = \left(P_2^1, P_2^2, P_2^3\right)$ is a special group set projection operator of $V = (V_1, V_2, V_3)$ into the special group set vector subspace $W_2 = \left(W_2^1, W_2^2, W_2^3\right)$.

Now $P_3 = \left(P_1^3, P_2^3, P_3^3\right) : V \to W_3 = \left(W_3^1, W_3^2, W_3^3\right)$ where
$P_1^3 : V_1 \to W_3^1$ is defined by
$$P_1^3 (x\ y\ z\ \omega) = (0\ 0\ 0\ \omega)$$
$P_2^3 : V_2 \to W_3^2$ is defined by
$$P_2^3 (a_0 + a_1 x + \ldots + a_6 x^6) = a_4 x^4 + a_5 x^5 + a_6 x^6)$$
and
$P_3^3 : V_3 \to W_3^3$ by

$$P_3^3 \begin{bmatrix} a & b & c \\ d & e & f \\ g & h & i \end{bmatrix} = \begin{bmatrix} 0 & 0 & 0 \\ d & 0 & 0 \\ g & h & 0 \end{bmatrix}.$$

It is easily verified
$$P_3 = \left(P_1^3, P_2^3, P_3^3\right) : V \to W_3 = \left(W_3^1, W_3^2, W_3^3\right)$$
is a special group set projection operator of $V$ into $W_3$.

Now consider

$P_1 + P_2 + P_3 = \left(P_1^1 + P_2^1 + P_3^1, P_1^2 + P_2^2 + P_3^2, P_1^3 + P_2^3 + P_3^3\right) : V \to V$.

$\left(P_1^1 + P_2^1 + P_3^1\right) (x\ y\ z\ \omega)$

$\quad = \quad P_1^1 (x\ y\ z\ \omega) + P_2^1 (x\ y\ z\ \omega) + P_3^1 (x\ y\ z\ \omega)$
$\quad = \quad (x\ y\ 0\ 0) + (0\ 0\ z\ 0) + (0\ 0\ 0\ \omega) = (x\ y\ z\ \omega);$

Thus $\left(P_1^1 + P_2^1 + P_3^1\right) = I_1$, that is identity map on $V_1$, i.e., $I_1: V_1 \to V_1$ is the special group set linear operator which is an identity operator on $V_1$.

Now consider $P_1^2 + P_2^2 + P_3^2 : V_2 \to V_2$ defined by

$$\left(P_1^2 + P_2^2 + P_3^2\right)(a_0 + a_1x + \ldots + a_6x^6)$$
$$= P_1^2(a_0 + a_1x + \ldots + a_6x^6) + P_2^2(a_0 + a_1x + \ldots + a_6x^6) +$$
$$P_3^2(a_0 + a_1x + \ldots + a_6x^6)$$
$$= (a_0 + a_1x) + (a_2x^2 + a_3x^3) + (a_4x^4 + a_5x^5 + a_6x^6)$$
$$= (a_0 + a_1x + \ldots + a_6x^6),$$

that is
$$\left(P_1^2 + P_2^2 + P_3^2\right) = I_2 : V_2 \to V_2.$$

Now consider
$$P_1^3 + P_2^3 + P_3^3 : V_3 \to V_3;$$

$$\left(P_1^3 + P_2^3 + P_3^3\right)\begin{pmatrix} a & b & c \\ d & e & f \\ g & h & i \end{pmatrix}$$

$$= P_1^3\begin{pmatrix} a & b & c \\ d & e & f \\ g & h & i \end{pmatrix} + P_2^3\begin{pmatrix} a & b & c \\ d & e & f \\ g & h & i \end{pmatrix} + P_3^3\begin{pmatrix} a & b & c \\ d & e & f \\ g & h & i \end{pmatrix}$$

$$= \begin{pmatrix} 0 & b & c \\ 0 & 0 & f \\ 0 & 0 & 0 \end{pmatrix} + \begin{pmatrix} a & 0 & 0 \\ 0 & e & 0 \\ 0 & 0 & i \end{pmatrix} + \begin{pmatrix} 0 & 0 & 0 \\ d & 0 & 0 \\ g & h & 0 \end{pmatrix}$$

$$= \begin{pmatrix} a & b & c \\ d & e & f \\ g & h & i \end{pmatrix};$$

thus
$$\left(P_1^3 + P_2^3 + P_3^3\right) = I_3 : V_3 \to V_3.$$

Thus
$$\begin{aligned}P_1 + P_2 + P_3 &= \left(P_1^1 + P_2^1 + P_3^1, P_1^2 + P_2^2 + P_3^2, P_1^3 + P_2^3 + P_3^3\right)\\ &= (I_1, I_2, I_3): V \to V.\end{aligned}$$

Thus we see the sum of the projection operators from a direct sum of subspaces is a special group set linear identity operator on $V = (V_1, V_2, V_3)$.

Now consider the composition of
$$P_j^i o P_t^i, j \neq t; \quad 1 \leq i \leq 3.$$

$$P_2^1 \text{ o } P_1^1 \text{ (a b c d)} = P_2^1 \text{ (a b 0 0)} = (0\ 0\ 0\ 0)$$
and
$$P_1^1 \text{ o } P_2^1 \text{ (a b c d)} = P_1^1 \text{ (0 0 c d)} = (0\ 0\ 0\ 0).$$
Thus
$$P_1^1 \text{ o } P_2^1 = P_2^1 \text{ o } P_1^1 = \theta_1$$
Also
$$P_2^1 \text{ o } P_3^1 \text{ (a b c d)} = P_2^1 (\ 0\ 0\ 0\ d) = (0\ 0\ 0\ 0)$$
and
$$P_3^1 \text{ o } P_2^1 \text{ (a b c d)} = P_3^1 \text{ (0 0 c 0)} = (0\ 0\ 0\ 0).$$
So
$$P_3^1 \text{ o } P_2^1 = P_2^1 \text{ o } P_3^1 = \theta_1.$$
Likewise
$$P_1^2 \text{ o } P_3^2 \ (a_0 + a_1 x + \ldots + a_6 x^6) = P_1^2 \ (a_4 x^4 + a_5 x^5 + a_6 x^6) = \theta_2.$$
Now
$$P_3^2 \text{ o } P_1^2 \ (a_0 + a_1 x + \ldots + a_6 x^6) = P_3^2 \ (a_0 + a_1 x) = 0.$$
Hence
$$P_3^2 \text{ o } P_1^2 = P_1^2 \text{ o } P_3^2 = \theta_2.$$

$$P_1^3 \text{ o } P_3^3 \left[\begin{pmatrix} a & b & c \\ d & e & f \\ g & h & i \end{pmatrix}\right] = P_1^3 \begin{pmatrix} 0 & 0 & 0 \\ d & 0 & 0 \\ g & h & 0 \end{pmatrix}$$

$$= \begin{pmatrix} 0 & 0 & 0 \\ 0 & 0 & 0 \\ 0 & 0 & 0 \end{pmatrix}.$$

Also

$$P_3^3 \text{ o } P_1^3 \begin{pmatrix} a & b & c \\ d & e & f \\ g & h & i \end{pmatrix}$$

$$P_3^3 \begin{pmatrix} a & b & c \\ 0 & 0 & f \\ 0 & 0 & 0 \end{pmatrix} = \begin{pmatrix} 0 & 0 & 0 \\ 0 & 0 & 0 \\ 0 & 0 & 0 \end{pmatrix}.$$

Thus we have

$$P_3^3 \text{ o } P_1^3 = P_1^3 \text{ o } P_3^3 = \theta_3. \quad \theta_3: V_3 \to V_3$$

such that $\theta_3(x) = 0$ for all $x \in V_3$ so

$$P_j \text{ o } P_k = (\theta_1, \theta_2, \theta_3).$$

where $j \neq k$, $1 \leq k, j \leq 3$.

These concepts can be extended to special group set linear algebras defined over the group G. As mentioned earlier when we use special group set linear algebra in the place of special group set vector spaces we see the cardinality of the generating set becomes considerably small in case of special group set linear algebras apart from that all properties defined for special group set vector spaces can be easily extended to special group set linear algebras.

Infact while using the set vector spaces or group vector spaces or semigroup vector spaces or special set vector spaces and so on we do not have the means to transform these algebraic structures to matrix equivalent. But when we want to use them in problems in which the notion of matrix is not used we find these structures can be easily replaced by the conventional ones.

Now we proceed onto describe the fuzzy analogue of these notions.

**Chapter Four**

# SPECIAL FUZZY SEMIGROUP SET VECTOR SPACE AND THEIR GENERALIZATIONS

This chapter has two sections. Section one defines the notion of special fuzzy semigroup set vector spaces and gives some important properties related to them. Section two defines the new concept of special semigroup fuzzy set n-vector spaces and special group fuzzy set n-vector spaces and describe some properties about them.

## 4.1 Special Fuzzy Semigroup Set Vector Spaces and their Properties

In this section we introduce the fuzzy structures defined in the earlier chapters.
   More specifically we define the new notion of special fuzzy semigroup set vector spaces and describe some of its properties.

**DEFINITION 4.1.1:** *Let $V = (S_1\ S_2,\ldots, S_n)$ be a special semigroup set vector space over the set $P.\mu: V \rightarrow [0, 1]$ is said to be a* *special fuzzy semigroup set vector space (special semigroup set fuzzy vector space) if the following conditions hold good.*

$$\mu = (\mu_1, \mu_2, ..., \mu_n): V = (S_1, S_2, ..., S_n) \to [0, 1]$$

*is such that for each i,*

$$\mu_i: S_i \to [0, 1];$$

*$1 \leq i \leq n$ satisfies the following conditions:*

1. *$\mu_i(x_i + y_i) \geq \min \{\mu_i(x_i), \mu_i(y_i)\}$*
2. *$\mu_i(s\, x_i) \geq \mu_i(x_i)$*

*for all $s \in P$ and for all $x_i$, $y_i$ in $S_i$ for every i, $1 \leq i \leq n$. We denote the special semigroup set fuzzy vector space by*

$$V_\mu = (S_1, ..., S_n)_\mu = \left(S_1, S_2, ..., S_n\right)_{(\mu_1, ..., \mu_n)}$$
$$= \left(S_{1\mu_1}, S_{2\mu_2}, ..., S_{n\mu_n}\right).$$

We illustrate this by a simple example.

***Example 4.1.1:*** Let $V = \{S_1, S_2, S_3, S_4\}$ be a special semigroup set vector space over the set $S = \{0, 1\}$ where

$$S_1 = \left\{ \begin{pmatrix} a & b \\ c & d \end{pmatrix} \,\bigg|\, a, b, c, d \in Z_2 = \{0,1\} \right\},$$

$$S_2 = \{Z_2 \times Z_2 \times Z_2 \times Z_2\},$$

$S_3 = \{Z_2[x]$ all polynomials of degree less than or equal to 5$\}$

and

$$S_4 = \left\{ \begin{bmatrix} a & b \\ c & d \\ e & f \\ g & h \end{bmatrix} \,\bigg|\, a, b, c, d, e, f, g, h \in Z_2 = \{0, 1\} \right\}.$$

Define

$$\eta = (\eta_1, \eta_2, \eta_3, \eta_4): V = (S_1\ S_2\ S_3\ S_4) \to [0,1]$$

by

$$\eta_1: S_1 \to [0,1]$$
$$\eta_2: S_2 \to [0,1]$$
$$\eta_3: S_3 \to [0,1]$$

and $\eta_4 : S_4 \to [0,1]$.

where

$$\eta_1 \begin{bmatrix} a & b \\ c & d \end{bmatrix} = \begin{cases} \dfrac{1}{2} & \text{if } a \neq 0 \\ 1 & \text{if } a = 0 \end{cases}$$

$$\eta_2 [a\ b\ c\ d] = \begin{cases} \dfrac{1}{5} & \text{if } a+b+c+d \neq 0 \\ 1 & \text{if } a+b+c+d = 0 \end{cases}$$

$$\eta_3[p(x)] = \begin{cases} \dfrac{1}{\deg p(x)} & \text{if } p(x) \text{ is not a constant} \\ 1 & \text{if } p(x) = \text{constant} \end{cases}$$

$$\eta_4 \begin{bmatrix} a & b \\ c & d \\ e & f \\ g & h \end{bmatrix} = \begin{cases} \dfrac{1}{6} & \text{if } a+b \neq 0 \\ 1 & \text{if } a+b = 0 \end{cases}$$

Thus $V\eta = \left( S_{1\eta_1}, S_{2\eta_2}, S_{3\eta_3}, S_{4\eta_4} \right)$ is a special semigroup set fuzzy vector space.

Now we proceed onto define the notion of special semigroup set fuzzy vector subspace.

**DEFINITION 4.1.2:** *Let $V = (S_1, S_2, \ldots, S_n)$ be a special semigroup set vector space over the set P. Take $W = (W_1, W_2, W_3, \ldots, W_n) \subseteq (V_1, \ldots, V_n)$ be a special semigroup set vector subspace of V over the same set P.*

*Define $\eta = (\eta_1, \eta_2, \ldots, \eta_n): W \to [0,1]$ as $\eta_i: W_i \to [0,1]$ for every i such that $\left( W_{1\eta_1}, W_{2\eta_2}, W_{3\eta_3}, \ldots, W_{n\eta_n} \right)$ is a special semigroup set fuzzy vector space then we call*
$$W_\eta = \left( W_{1\eta_1}, W_{2\eta_2}, W_{3\eta_3}, \ldots, W_{n\eta_n} \right)$$

*to be a special semigroup set fuzzy vector subspace.*

We shall illustrate this by an example.

***Example 4.1.2:*** Let $V = (S_1, S_2, S_3, S_4, S_5)$ where $S_1 = \{S \times S \times S \mid S = Z^o = Z^+ \cup \{0\}\}$, $S_2 = \{[a\ a\ a\ a\ a\ a] \mid a \in S\}$,

$$S_3 = \left\{ \begin{pmatrix} a & b \\ c & d \end{pmatrix} \middle| a, b, c, d \in 2Z^o \right\},$$

$S_4 = \{$all polynomials of degree less than or equal to 5 with coefficients from $Z^o\}$ and

$$S_5 = \left\{ \begin{bmatrix} a & a \\ a & a \\ a & a \\ a & a \end{bmatrix} \middle| a \in Z^o \right\}$$

be a special semigroup set vector space over the set S. Choose $W = (W_1, W_2, \ldots, W_5)$ such that $W_1 = (S \times S \times \{0\}) \subseteq S_1$, $W_2 = \{[a\ a\ a\ a\ a\ a] \mid a \in 3Z^o\} \subseteq S_2$,

$$W_3 = \left\{ \begin{pmatrix} a & a \\ a & a \end{pmatrix} \middle| a \in 6Z^o \right\} \subseteq S_3,$$

$W_4 = \{$all polynomials of degree less than or equal to 5 with coefficients from $3Z^o\} \subseteq S_4$ and

$$W_5 = \left\{ \begin{bmatrix} a & a \\ a & a \\ a & a \\ a & a \end{bmatrix} \middle| a \in 7Z^o \right\}$$

be a special semigroup set vector subspace of V over the set

$S = Z^o$. Define
$$\eta = (\eta_1, \eta_2, \ldots, \eta_5) : W = (W_1, W_2, \ldots, W_5) \to [0,1]$$
such that $\eta_i : W_i \to [0,1]$ for $i = 1, 2, \ldots, 5$ is defined as follows:

$\eta_1 : W_1 \to [0,1]$ is defined by

$$\eta_1(a\ b\ 0) = \begin{cases} \dfrac{1}{2} & \text{if } a+b \text{ is even} \\ \dfrac{1}{3} & \text{if } a+b \text{ is odd} \\ 1 & \text{if } a+b = 0 \end{cases}$$

$\eta_2 : W_2 \to [0,1]$ is such that

$$\eta_2(a\ a\ a\ a\ a\ a) = \begin{cases} \dfrac{1}{6} & \text{if } a \text{ is odd}, a \neq 0 \\ 1 & \text{if } a \text{ is even or } 0 \end{cases}$$

$\eta_3: W_3 \to [0, 1]$ is given by

$$\eta_3 \begin{pmatrix} a & a \\ a & a \end{pmatrix} = \begin{cases} \dfrac{1}{5} & \text{if } a \neq 0 \\ 1 & \text{if } a = 0 \end{cases}$$

$\eta_4: W_4 \to [0, 1]$ is such that

$$\eta_4(p(x)) = \begin{cases} \dfrac{1}{\deg p(x)} & \text{if } p(x) \neq \text{a constant} \\ 1 & \text{if } p(x) = \text{constant polynomial} \end{cases}$$

$\eta_5 : W_5 \to [0,1]$ is defined by

$$\eta_5 \begin{pmatrix} a & a \\ a & a \\ a & a \\ a & a \end{pmatrix} = \begin{cases} \dfrac{1}{3} & \text{if } a \text{ is even } a \neq 0 \\ \dfrac{1}{6} & \text{if } a \text{ is odd } a \neq 0 \\ 1 & \text{if } a = 0 \end{cases}$$

$$W\eta = \left(W_1, W_2, W_3, W_4, W_5\right)_{(\eta_1, \eta_2, \ldots, \eta_5)} = \left(W_{1\eta_1}, W_{2\eta_2}, \ldots, W_{5\eta_5}\right)$$

is a special semigroup set fuzzy vector subspace.

**DEFINITION 4.1.3:** *Let $V = \{V_1, V_2, \ldots, V_n\}$ be a special semigroup set linear algebra over the semigroup $S$, $S$ a semigroup under addition. If $\eta = (\eta_1, \eta_2, \ldots, \eta_n): V = (V_1, \ldots, V_n) \to [0, 1]$ is such that $\eta_i: V_i \to [0, 1]$; $\eta_i(a_i + b_i) \geq \min(\eta_i(a_i), \eta_i(b_i))$ and $\eta_i(sa_i) \geq \eta_i(a_i)$ for every $a_i, b_i, \in V_i$ and for all $s \in S$ and this is true for every $i = 1, 2, \ldots, n$; then we call*

$$V_\eta = \left(V_1, V_2, \ldots, V_n\right)_{(\eta_1, \ldots, \eta_n)} = \left(V_{1\eta_1}, V_{2\eta_2}, \ldots, V_{n\eta_n}\right)$$

*be a special semigroup set fuzzy linear algebra. It is pertinent to mention here that if $V_\eta$ is a special semigroup set fuzzy linear algebra is the same as the special semigroup set fuzzy vector space. Thus the notion of special semigroup set fuzzy vector space and special semigroup set fuzzy linear algebra are fuzzy equivalent.*

However we illustrate this by a simple example.

***Example 4.1.3:*** Let $V = (V_1, V_2, V_3, V_4)$ where $V_1 = \{S \times S \times S \times S \text{ such that } S = Z^o = Z^+ \cup \{0\}\}$, $V_2 = \{$all polynomials of degree less than or equal to 5 with coefficients from $S\}$,

$$V_3 = \left\{ \begin{bmatrix} a & a & a & a & a \\ a & a & a & a & a \end{bmatrix} \middle| a \in S \right\}$$

and $V_4 = \{$All $3 \times 3$ matrices with entries from $S\}$ be a special semigroup set linear algebra over the semigroup $S = Z^o = Z^+$

∪ {0}. Define η = (η₁, η₂, η₃, η₄): V = (V₁, V₂, V₃, V₄) → [0,1] such that $\eta_i: V_i \to [0,1]$ for every i = 1, 2, 3, 4. $\eta_i$'s are defined in the following way

$\eta_1: V_1 \to [0,1]$ is such that

$$\eta_1(a\ b\ c\ d) = \begin{cases} \dfrac{1}{2} & \text{if } ab+cd \text{ is even} \\ 1 & \text{if } ab+cd \text{ is odd or zero} \end{cases}$$

$\eta_2: V_2 \to [0,1]$ is given by

$$\eta_2(p(x)) = \begin{cases} \dfrac{1}{5} & \text{if deg } p(x) \text{ is even} \\ \dfrac{1}{2} & \text{if deg } p(x) \text{ is odd} \\ 1 & \text{if } p(x) \text{ is constant} \end{cases}$$

$\eta_3: V_3 \to [0,1]$ is defined by

$$\eta_3 \begin{pmatrix} a & a & a & a & a \\ a & a & a & a & a \end{pmatrix} = \begin{cases} \dfrac{1}{5} & \text{if } a \text{ is even} \\ \dfrac{1}{2} & \text{if } a \text{ is odd} \\ 1 & \text{if } a = 0 \end{cases}$$

$\eta_4: V_4 \to [0,1]$ is such that

$$\eta_4 \begin{pmatrix} a & b & c \\ d & e & f \\ g & h & i \end{pmatrix} = \begin{cases} \dfrac{1}{2} & \text{if trace sum is even} \\ \dfrac{1}{3} & \text{if trace sum is odd} \\ 1 & \text{if } a+e+i = 0 \end{cases}$$

Thus

$$V_\eta = (V_1, V_2, V_3, V_4)_{(\eta_1, \eta_2, \eta_3, \eta_4)} = \left(V_{1\eta_1}, V_{2\eta_2}, V_{3\eta_3}, V_{4\eta_4}\right)$$

is a special semigroup set fuzzy linear algebra. Clearly $V_\eta$ is also a special semigroup set fuzzy vector space.

Next we proceed on to define the notion of special semigroup set fuzzy linear subalgebra.

**DEFINITION 4.1.4:** *Let $V = (V_1, V_2, ..., V_n)$ be a special semigroup set linear algebra over the additive semigroup S.*

*Let $W = (W_1, W_2, ..., W_n) \subseteq (V_1, V_2, ..., V_n)$ be such that W is a special semigroup set linear subalgebra of V over S. Now define $\eta = (\eta_1, \eta_2, ..., \eta_n) : W = (W_1, W_2, ..., W_n) \to [0,1]$ by $\eta_i: W_i \to [0,1]$ for each i such that $W_\eta = \left(W_{1\eta_1}, ..., W_{n\eta_n}\right)$ is also a special semigroup set fuzzy linear algebra then we call $W_\eta$ to be a special semigroup set fuzzy linear subalgebra.*

We shall illustrate this by an example.

***Example 4.1.4:*** Let $V = (V_1, V_2, V_3, V_4, V_5)$ where

$$V_1 = \left\{ \begin{bmatrix} a & a & a & a \\ a & a & a & a \end{bmatrix} \middle| a \in Z_6 \right\},$$

$$V_2 = \left\{ \begin{bmatrix} a & a & a \\ a & a & a \\ a & a & a \\ a & a & a \\ a & a & a \end{bmatrix} \middle| a \in Z_6 \right\},$$

$V_3 = \{$all polynomials in the variable x with coefficients from $Z_6$ of degree less than or equal to 4$\}$, $V_4 = \{3 \times 3$ matrices with entries from $Z_6\}$ and $V_5 = \{Z_6 \times Z_6 \times Z_6 \times Z_6\}$ be a special semigroup set linear algebra over the semigroup $S = Z_6$. Consider $W = (W_1, W_2, W_3, W_4, W_5) \subseteq (V_1, V_2, V_3, V_4, V_5) = V$ where

$$W_1 = \left\{ \begin{bmatrix} a & a & a & a \\ a & a & a & a \end{bmatrix} \middle| \; a \in \{0, 2, 4\}; \right.$$

addition of matrices modulo $6\} \subseteq V_1$,

$$W_2 = \left\{ \begin{bmatrix} a & a & a \\ a & a & a \\ a & a & a \\ a & a & a \\ a & a & a \end{bmatrix} \middle| \; a \in \{0, 2, 4\}; \right.$$

addition of matrices modulo $6\} \subseteq V_2$,

$W_3 = \{$all polynomials in x of degree less than or equal to three with coefficients from $Z_6\} \subseteq V_3$,

$W_4 = \{$upper triangular $3 \times 3$ matrices with entries from $Z_6\}$
$$\subseteq V_4$$
and
$$W_5 = \{Z_6 \times Z_6 \times \{0\} \times \{0\}\} \subseteq V_5;$$

$W = (W_1, W_2, W_3, W_4, W_5)$ is a special semigroup set linear subalgebra of V over the semigroup $Z_6$.

Take
$$\eta = (\eta_1, \eta_2, \eta_3, \eta_4, \eta_5): W = (W_1, W_2, \ldots, W_5) \to [0,1]$$

where $\eta_i : W_i \to [0,1]$ for $i = 1, 2, \ldots, 5$ such that

$\eta_1: W_1 \to [0,1]$ is defined by

$$\eta_1 \begin{bmatrix} a & a & a & a \\ a & a & a & a \end{bmatrix} = \begin{cases} \dfrac{1}{2} & \text{if } a \neq 0 \\ 1 & \text{if } a = 0 \end{cases}$$

$\eta_2$: $W_2 \to [0,1]$ is such that

$$\eta_2 \begin{bmatrix} a & a & a \\ a & a & a \\ a & a & a \\ a & a & a \\ a & a & a \end{bmatrix} = \begin{cases} \dfrac{1}{9} & \text{if } a \neq 0 \\ 1 & \text{if } a = 0 \end{cases}$$

$\eta_3$: $W_3 \to [0,1]$ is given by

$$\eta_3(p(x)) = \begin{cases} \dfrac{1}{\deg p(x)} & \text{if } p(x) \neq \text{a constant} \\ 1 & \text{if } p(x) = \text{constant polynomial} \end{cases}$$

$\eta_4$: $W_4 \to [0,1]$ is defined by

$$\eta_4 \begin{bmatrix} a & a & a \\ 0 & d & e \\ 0 & 0 & f \end{bmatrix} = \begin{cases} \dfrac{1}{2} & \text{if } a+d+e \neq 0 \\ 1 & \text{if } a+d+e = 0 \end{cases}$$

$\eta_5$: $W_5 \to [0,1]$ is given by

$$\eta_5\,(a\ b\ 0\ 0) = \begin{cases} \dfrac{1}{2} & \text{if } a+b \neq 0 \\ 1 & \text{if } a+b = 0 \end{cases}$$

Thus
$$W_\eta = \left(W_1, W_2, W_3, W_4, W_5\right)_{(\eta_1, \eta_2, \ldots, \eta_5)} = \left(W_{1\eta_1}, W_{2\eta_2}, \ldots, W_{5\eta_5}\right)$$
is a special semigroup set fuzzy linear subalgebra.

However it is pertinent to mention here that the notion of special semigroup set fuzzy vector subspace and the special semigroup set fuzzy linear subalgebra are fuzzy equivalent.

Now we proceed on to define the notion of special group set fuzzy vector space.

**DEFINITION 4.1.5:** *Let $V = (V_1, V_2, \ldots, V_n)$ be a special group set vector space over the set S. Let $\eta = (\eta_1, \eta_2, \ldots, \eta_n) : V = (V_1, \ldots V_n) \to [0,1]$ be such that $\eta_i : V_i \to [0,1]$ for each i, $1 \le i \le n$ satisfying the following conditions*

1. $\eta_i(a_i + b_i) \ge min\{\eta_i(a_i), \eta_i(b_i)\}$
2. $\eta_i(-a_i) = \eta(a_i)$
3. $\eta_i(0) = 1$
4. $\eta_i(ra_i) \ge \eta_i(a_i)$ *for all $r \in S$ and for all $a_i, b_i \in V_i$, true for each $i = 1, 2, \ldots, n$.*

*Thus $\left(V_{i\eta_i}\right)$ is a group set fuzzy vector space, $1 \le i \le n$. Now*
$$V_\eta = \left(V_1, V_2, \ldots, V_n\right)_{(\eta_1, \eta_2, \ldots, \eta_n)} = \left(V_{1\eta_1}, V_{2\eta_2}, \ldots, V_{n\eta_n}\right)$$
*is a special set group fuzzy vector space.*

We illustrate this by a simple example.

***Example 4.1.5:*** Let $V = (V_1, V_2, V_3, V_4, V_5)$ where $V_1 = Z_{10}$ additive abelian group modulo 10. $V_2 = Z_{15} \times Z_{15} \times Z_{15}$ is again a group under component wise addition modulo 15, $V_3 = Z_{12}$, $V_4 = Z_{14} \times Z_{14}$ and $V_5 = Z_7 \times Z_7 \times Z_7 \times Z_7$ is a special set group vector space over the set $S = \{0, 1\}$. Define
$\eta = (\eta_1, \eta_2, \eta_3, \eta_4, \eta_5): V = (V_1, V_2, V_3, V_4, V_5) \to [0, 1]$
by $\eta_i : V_i \to [0,1]$, $1 \le i \le 5$ as

$\eta_1: V_1 \to [0,1]$ is defined by

$$\eta_1(a) = \begin{cases} \dfrac{1}{2} & \text{if } a \in \{2, 4, 6, 8\} \\ 1 & \text{if } a = 0 \\ \dfrac{1}{3} & \text{if } a \in \{1, 3, 5, 7, 9\} \end{cases}$$

$\eta_2: V_2 \rightarrow [0,1]$ is given by

$$\eta_2 (a\ b\ c) = \begin{cases} \dfrac{1}{6} & \text{if } a+b+c \in \{1,\ 2,\ 8\} \\ \dfrac{1}{5} & \text{if } a+b+c \in \{4,\ 6,\ 10,\ 12\} \\ \dfrac{1}{4} & \text{if } a+b+c \in \{3,\ 7,\ 11\} \\ \dfrac{1}{3} & \text{if } a+b+c \in \{5,\ 9,\ 13,\ 14\} \\ 1 & \text{if } a = 0 \text{ or } a+b+c = 0 \end{cases}$$

$\eta_3: V_3 \rightarrow [0,\ 1]$ is such that

$$\eta_3(a) = \begin{cases} \dfrac{1}{9} & \text{if } a \in \{2, 4,\ 6,\ 8,\ 10\} \\ 1 & \text{if } a = 0 \\ \dfrac{1}{8} & \text{if } a \in \{1,\ 3,\ 5,\ 7,\ 9,\ 11\} \end{cases}$$

$\eta_4: V_4 \rightarrow [0,1]$ is given by

$$\eta_4 (a,\ b) = \begin{cases} \dfrac{1}{3} & \text{if } a+b \in \{2, 4,\ 6,\ 8,\ 10,\ 12\} \\ 1 & \text{if } a+b = 0 \\ \dfrac{1}{5} & \text{if } a+b \in \{1,\ 3,\ 5,\ 7,\ 9,\ 11,\ 13\} \end{cases}$$

$\eta_5: V_5 \rightarrow [0,1]$

$$\eta_5 (a\ b\ c\ d) = \begin{cases} \dfrac{1}{3} & \text{if } a+b+c+d \in \{2, 4,\ 6\} \\ 1 & \text{if } a+b = 0 \\ \dfrac{1}{2} & \text{if } a+b+c+d \in \{1,\ 3,\ 5\} \end{cases}$$

Thus
$$V\eta = (V_1, V_2, V_3, V_4, V_5)_{(\eta_1,\eta_2,\eta_3,\eta_4,\eta_5)} = (V_{1\eta_1}, V_{2\eta_2}, ..., V_{n\eta_n})$$
is a special set group fuzzy vector space.

Now we proceed on to define the notion of special group set fuzzy vector subspace.

**DEFINITION 4.1.6:** *Let $V = (V_1, V_2, ..., V_n)$ be a special group set vector over the set S. Let $W = (W_1, W_2, ..., W_n) \subseteq (V_1, V_2, ..., V_n) = V$; $(W_i \subseteq V_i, 1 \leq i \leq n)$ be a special set group vector subspace of V over the set S. Let $\eta = (\eta_1, \eta_2, ..., \eta_n): W = (W_1, ..., W_n) \to [0,1]$ be such that $\eta_i : W_i \to [0,1]$ gives $W_{i\eta_i}$ to be a group set vector subspace for each $i = 1, 2, ..., n$. Then*
$$W_\eta = (W_{1\eta_1}, W_{2\eta_2}, ..., W_{n\eta_n})$$
*is defined as the special set group fuzzy vector subspace.*

We illustrate this by an example.

*Example 4.1.6:* Let $V = (V_1, V_2, V_3, V_4, V_5)$ be a special set group vector space over the set $S = \{0, 1, 2, ..., 10\}$ where $V_1 = 2Z$, $V_2 = 3Z$, $V_3 = 4Z$, $V_4 = 5Z$ and $V_5 = 7Z$. Now take $W = (W_1, W_2, W_3, W_4, W_5) \subseteq (V_1, V_2, V_3, V_4, V_5)$, a special group set vector subspace of V where $W_1 = 14Z \subseteq V_1$, $W_2 = 6Z \subseteq V_2$, $W_3 = 12Z \subseteq V_3$, $W_4 = 10Z \subseteq V_4$ and $W_5 = 21Z\ V_5$. Define $\eta = (\eta_1, \eta_2, \eta_3, \eta_4, \eta_5): W = (W_1, W_2, W_3, W_4, W_5) \to [0, 1]$ by

$\eta_1 : W_1 \to [0,1]$ such that
$$\eta_1(a) = \begin{cases} \frac{1}{5} & \text{if } a \neq 0 \\ 1 & \text{if } a = 0 \end{cases}$$

$\eta_2 : W_2 \to [0,1]$ defined by
$$\eta_2(a) = \begin{cases} \frac{1}{6} & \text{if } a \neq 0 \\ 1 & \text{if } a = 0 \end{cases}$$

$\eta_3 : W_3 \to [0,1]$ is such that

$$\eta_3(a) = \begin{cases} \dfrac{1}{2} & \text{if } a \neq 0 \\ 1 & \text{if } a = 0 \end{cases}$$

$\eta_4 : W_4 \to [0,1]$ is such that

$$\eta_4(a) = \begin{cases} \dfrac{1}{9} & \text{if } a \neq 0 \\ 1 & \text{if } a = 0 \end{cases}$$

and $\eta_5 : W_5 \to [0,1]$ is defined by

$$\eta_5(a) = \begin{cases} \dfrac{1}{7} & \text{if } a \neq 0 \\ 1 & \text{if } a = 0 \end{cases}$$

Thus $W_\eta = \left( W_{1\eta_1}, W_{2\eta_2}, ..., W_{5\eta_5} \right)$ is a special group set fuzzy vector subspace.

We will now proceed on to define special group set fuzzy linear algebra.

**DEFINITION 4.1.7:** *Let $V = (V_1, V_2, ..., V_n)$ be a special group set linear algebra defined over the set S. Let $\eta = (\eta_1, ..., \eta_n)$ be a map from V into [0,1], that is $\eta = (\eta_1, \eta_2, ... \eta_n): V = (V_1, ..., V_n) \to [0, 1]$; such that $\eta_i: V_i \to [0,1]$ for each i; $V_{i\eta_i}$ is a group set fuzzy linear algebra, $1 \leq i \leq n$. Then $V_\eta = \left( V_{1\eta_1}, V_{2\eta_2}, ..., V_{n\eta_n} \right)$ is a special group set fuzzy linear algebra.*

We illustrate this by an example, we also observe that the notion of special group set fuzzy linear algebra and special group set fuzzy vector space are fuzzy equivalent.

**Example 4.1.7:** Let $V = (V_1, V_2, V_3, V_4)$ be a special group set linear algebra over the set $Z^+ \cup \{0\} = S$. Here $V_1 = \{$set of all $2 \times 2$ matrices with entries from $Z\}$, $V_2 = Z \times Z \times Z \times Z$,

$$V_3 = \left\{ \begin{bmatrix} a \\ a \\ a \\ a \\ a \\ a \end{bmatrix} \middle| \ a \in Z \right\}$$

and $V_4 = \{Z^7[x]$ all polynomials of degree less than or equal to 7 with coefficients from $Z\}$. Now define a map
$$\eta = (\eta_1, \eta_2, \eta_3, \eta_4) : V = (V_1, V_2, V_3, V_4) \to [0,1]$$
by $\eta_i : V_i \to [0, 1]$ such that $V_{i\eta_i}$ is a group fuzzy linear algebra for $i = 1, 2, 3, 4$ as follows:

$\eta_1: V_1 \to [0, 1]$ is defined by

$$\eta_1 \begin{pmatrix} a & b \\ c & d \end{pmatrix} = \begin{cases} \dfrac{1}{2} & \text{if } a \text{ is even} \\ \dfrac{1}{3} & \text{if } a \text{ is odd} \\ 1 & \text{if } a = 0 \end{cases}$$

$\eta_2: V_2 \to [0, 1]$ is such that

$$\eta_2(a\ b\ c\ d) = \begin{cases} \dfrac{1}{5} & \text{if } a+b+c+d \text{ is even} \\ \dfrac{1}{7} & \text{if } a+b+c+d \text{ is odd} \\ 1 & \text{if } a+b+c+d = 0 \end{cases}$$

$\eta_3: V_3 \to [0, 1]$ is defined by

$$\eta_3 \begin{bmatrix} a \\ a \\ a \\ a \\ a \end{bmatrix} = \frac{1}{12} \begin{cases} \dfrac{1}{12} & \text{if a is even} \\ \dfrac{1}{13} & \text{if a is odd} \\ 1 & \text{if } a = 0 \end{cases}$$

Finally $\eta_4 : V_4 \to [0,1]$ is given by

$$\eta_4(p(x)) = \begin{cases} \dfrac{1}{23} & \text{if deg } p(x) \text{ is even} \\ \dfrac{1}{2} & \text{if deg } p(x) \text{ is odd} \\ 1 & \text{if } p(x) \text{ is constant} \end{cases}$$

Thus $V_\eta = \left(V_{1\eta_1}, V_{2\eta_2}, V_{3\eta_3}, V_{4\eta_4}\right)$ is a special group set fuzzy linear algebra.

Now we proceed on to define the notion of special group set fuzzy linear subalgebra.

**DEFINITION 4.1.8:** *Let $V = (V_1, V_2, \ldots, V_n)$ be a special group set linear algebra over the set S. Let $W = (W_1, W_2, \ldots, W_n) \subseteq (V_1, V_2, \ldots, V_n) = V$ be a special set linear subalgebra of V over the set S. Let $\eta = (\eta_1, \ldots, \eta_n) : W \to [0,1]$ be such that $\eta_i : W_i \to [0,1]$ for each i; $i = 1, 2, \ldots, n$. $W_{i\eta_i}$ is a group set fuzzy linear subalgebra for every $1 \leq i \leq n$. So that $W_\eta = \left(W_{1\eta_1}, W_{2\eta_2}, \ldots, W_{n\eta_n}\right)$ is a special group set fuzzy linear subalgebra.*

It is interesting to note that the notion of special group set fuzzy vector subspaces and special group set fuzzy linear subalgebras are fuzzy equivalent.

We how ever illustrate special group set fuzzy linear subalgebra.

***Example 4.1.8:*** Let $V = (V_1, V_2, V_3, V_4, V_5)$ be a special group set linear algebra over the set $S = 2Z^+ \cup \{0\}$ where $V_1 = Z \times Z \times Z$, $V_2 = \{$All $5 \times 5$ matrices with entries from $Z\}$, $V_3 = \{[a\ a\ a\ a\ a\ a\ a] \mid a \in Z\}$,

$$V_4 = \left\{ \begin{bmatrix} a \\ a \\ a \\ a \\ a \end{bmatrix} \middle| a \in Z \right\}$$

and $V_5 = \{Z[x]$ polynomials of degree less than or equal to 9 with coefficients from $Z\}$. Take $W = (W_1, W_2, W_3, W_4, W_5) \subseteq (V_1, V_2, \ldots, V_5)$ that is $W_i \subseteq V_i$, $1 \le i \le 5$, where $W_1 = Z \times \{0\} \times \{0\} \subseteq V_1$, $W_2 = \{$all $5 \times 5$ matrices with entries from $3Z\}$ contain as a group set linear subalgebra of $V_2$, $W_3 = \{[a\ a\ a\ a\ a\ a\ a] \mid a \in 5Z\} \subseteq V_3$,

$$W_4 = \left\{ \begin{bmatrix} a \\ a \\ a \\ a \\ a \end{bmatrix} \middle| a \in 7Z \right\} \subseteq V_4$$

and $W_5 = \{$All polynomials in $x$ of degree less than or equal to 5 with coefficients from $Z\} \subseteq V_5$. Thus $W = (W_1, W_2, \ldots, W_5) \subseteq V$ is a special group set linear subalgebra of $V$ over the set $2Z^+ \cup \{0\}$. Now define $\eta = (\eta_1, \eta_2, \eta_3, \eta_4, \eta_5): W = (W_1, W_2, W_3, W_4, W_5) \to [0, 1]$ by $\eta_i: W_i \to [0, 1]$, $1 \le i \le 5$ such that $W_{i\eta_i}$ is the group set fuzzy linear subalgebra for $i = 1, 2, \ldots, 5$.

Define $\eta_1: W_1 \to [0, 1]$ by

$$\eta_1(a\ 0\ 0) = \begin{cases} \dfrac{1}{7} & \text{if } a \text{ is even} \\ \dfrac{1}{5} & \text{if } a \text{ is odd} \\ 1 & \text{if } a = 0 \end{cases}$$

$\eta_2 : W_2 \to [0,1]$ by

$$\eta_2 \begin{pmatrix} a_1 & a_2 & a_3 & a_4 & a_5 \\ a_6 & a_7 & a_8 & a_9 & a_{10} \\ a_{11} & a_{12} & a_{13} & a_{14} & a_{15} \\ a_{16} & a_{17} & a_{18} & a_{19} & a_{20} \\ a_{21} & a_{22} & a_{23} & a_{24} & a_{25} \end{pmatrix} =$$

$$\begin{cases} \dfrac{1}{11} & \text{if } a_1 + a_7 + a_{13} + a_{19} + a_{25} \text{ is even} \\ \dfrac{1}{12} & \text{if } a_1 + a_7 + a_{13} + a_{19} + a_{25} \text{ is odd} \\ 1 & \text{if } a_1 + a_7 + a_{13} + a_{19} + a_{25} = 0 \end{cases}$$

Define $\eta_3: W_3 \to [0,1]$ by

$$\eta_3 \,(a\,a\,a\,a\,a\,a\,a) = \begin{cases} \dfrac{1}{19} & \text{if } a \text{ is even} \\ \dfrac{1}{5} & \text{if } a \text{ is odd} \\ 1 & \text{if } a = 0 \end{cases}$$

$\eta_4: W_4 \to [0, 1]$ is defined by

$$\eta_4 \begin{pmatrix} a \\ a \\ a \\ a \\ a \end{pmatrix} = \begin{cases} 1 & \text{if } a \text{ is even or } 0 \\ \dfrac{1}{4} & \text{if } a \text{ is odd} \end{cases}$$

$\eta_5 : W_5 \to [0,1]$ is given by

$$\eta_5(p(x)) = \begin{cases} \dfrac{1}{5} & \text{if } p(x) \text{ is odd degree} \\ \dfrac{1}{7} & \text{if } p(x) \text{ is even degree} \\ 1 & \text{if } p(x) \text{ is constant} \end{cases}$$

Thus
$$W\eta = \left(W_1, W_2, W_3, W_4, W_5\right)_{(\eta_1, \eta_2, \eta_3, \eta_4, \eta_5)} = \left(W_{1\eta_1}, W_{2\eta_2}, \ldots, W_{5\eta_5}\right)$$
is a special group set fuzzy linear subalgebra.

It is important mention at this juncture when we define various special group set substructures most of them happen to be fuzzy equivalent.

### 4.2 Special Semigroup Set n-vector Spaces

In this section we proceed on to define the new notion of special semigroup set n-vector spaces, special group set n-vector spaces and their fuzzy analogue.

Now we proceed on to define special semigroup set n-vector spaces, special group set n-vector spaces and their fuzzy analogue.

**DEFINITION 4.2.1:** *Let*
$$V = (V_1 \cup V_2) = \left(V_1^1, V_2^1, \ldots, V_{n_1}^1\right) \cup \left(V_1^2, V_2^2, \ldots, V_{n_2}^2\right)$$
*be a special semigroup set vector bispace over the set S if both $V_1$ and $V_2$ are two distinct special semigroup set vector spaces over the set S, $V_1 \neq V_2$, $V_1 \not\subseteq V_2$ or $V_2 \not\subseteq V_1$.*

We illustrate this by a simple example.

*Example 4.2.1:* Let
$$V = (V_1 \cup V_2) = \left(V_1^1, V_2^1, V_3^1, V_4^1\right) \cup \left(V_1^2, V_2^2, V_3^2, V_4^2, V_5^2\right)$$

where

$$V_1^1 = \left\{ \begin{pmatrix} a & b \\ c & d \end{pmatrix} \middle| a,b,c,d \in Z^+ \cup \{0\} \right\},$$

$$V_2^1 = S \times S \times S, \quad V_3^1 = \{[a\ a\ a\ a\ a] \mid a \in S\},$$

$$V_4^1 = \left\{ \begin{bmatrix} a \\ a \\ a \end{bmatrix} \middle| a \in S \right\}. \quad V_1^2 = S \times S \times S \times S,$$

$$V_2^2 = \left\{ \begin{bmatrix} a \\ a \\ a \\ a \\ a \end{bmatrix} \middle| a \in S \right\},$$

$V_3^2 = \{$set of all $3 \times 3$ matrices with entries from $S\}$, $V_4^2 = \{$All polynomials of degree less than or equal to 3 with coefficients from $S\}$ and $V_5^2 = \{$all $4 \times 4$ lower triangular matrices with entries from $S\}$. Thus $V = (V_1 \cup V_2)$ is a special semigroup set vector bispace over the set $S = Z^+ \cup \{0\}$.

Now we proceed to define the notion of special semigroup set vector subbispace of V over the set S.

**DEFINITION 4.2.2:** *Let*
$$V = (V_1 \cup V_2) = \left( V_1^1, V_2^1, ..., V_{n_1}^1 \right) \cup \left( V_1^2, V_2^2, ..., V_{n_2}^2 \right)$$
*be a special semigroup set vector bispace over the set S. Take*
$$W = (W_1 \cup W_2) = \left( W_1^1, W_2^1, ..., W_{n_1}^1 \right) \cup \left( W_1^2, W_2^2, ..., W_{n_2}^2 \right)$$
$$\subseteq (V_1 \cup V_2) = \left( V_1^1, V_2^1, ..., V_{n_1}^1 \right) \cup \left( V_1^2, V_2^2, ..., V_{n_2}^2 \right)$$

that is $W^i_{t_i} \subseteq V^i_{t_i}$; $i = 1,2$ and $1 \leq t_i \leq n_i$ and $1 \leq t_i \leq n_2$ is such that $W^i_{t_i}$ is a semigroup set vector subspace of $V^i_{t_i}$. Then we call $W = W_1 \cup W_2$ to be a special semigroup set vector subbispace of V over the set S.

We illustrate this situation by a simple example.

*Example 4.2.2:* Let
$$V = (V_1 \cup V_2) = \left(V^1_1, V^1_2, V^1_3\right) \cup \left(V^2_1, V^2_2, V^2_3, V^2_4, V^2_5\right)$$
where $V^1_1 = \{$All $5 \times 5$ matrices with entries from $S = Z^+ \cup \{0\}\}$,
$V^1_2 = \{[a\ a\ a\ a\ a\ a\ a] \mid a \in S\}$,

$$V^1_3 = \left\{ \begin{bmatrix} a & a \\ a & a \\ a & a \\ a & a \\ a & a \end{bmatrix} \middle| a \in S \right\}.$$

$$V^2_1 = \left\{ \begin{bmatrix} a & a & a & a & a \\ b & b & b & b & b \end{bmatrix} \middle| a,b \in S \right\},$$

$V^2_2 = \{$set of all $4 \times 4$ lower triangular matrices with entries from S$\}$, $V^2_3 = \{$set of all polynomials of degree less than or equal to 4 with coefficients form S$\}$, $V^2_4 = \{S \times S \times S \times S \times S\}$ and

$$V^2_5 = \left\{ \begin{bmatrix} a \\ a \\ a \\ a \\ a \end{bmatrix} \middle| a \in S \right\}$$

be a special semigroup set vector bispace over the set $S = Z^+ \cup \{0\} = Z^0$. Take

$$W = (W_1 \cup W_2) = \left(W_1^1, W_2^1, W_3^1\right) \cup \left(W_1^2, W_2^2, W_3^2, W_4^2, W_5^2\right)$$
$$\subseteq \left(V_1^1, V_2^1, V_3^1\right) \cup \left(V_1^2, V_2^2, V_3^2, V_4^2, V_5^2\right)$$

where $W_1^1 = \{$all 5×5 matrices with entries from $5S\} \subseteq V_1^1$, $W_2^1 = \{(a\ a\ a\ a\ a\ a\ a) \mid a \in 7S\} \subseteq V_2^1$,

$$W_3^1 = \left\{ \begin{bmatrix} a & a \\ a & a \\ a & a \\ a & a \\ a & a \end{bmatrix} \middle| a \in 7S \right\} \subseteq V_3^1;$$

$$W_1^2 = \left\{ \begin{bmatrix} a & a & a & a & a \\ b & b & b & b & b \end{bmatrix} \middle| a, b \in 11S \right\} \subseteq V_1^2,$$

$W_2^2 = \{$set of all 4×4 lower triangular matrices with entries from $3Z^0 \cup \{0\} \subseteq V_2^2$, $W_3^2 = \{$set of all polynomials in the variable x with coefficients from $12S\} \subseteq V_3^2$, $W_4^2 = \{S \times S \times S \times \{0\} \times \{0\}\} \subseteq V_4^2$ and

$$W_5^2 = \left\{ \begin{bmatrix} a \\ a \\ a \\ a \\ a \end{bmatrix} \middle| a \in 11S \right\} \subseteq V_5^2.$$

Clearly
$$W = (W_1 \cup W_2) = \left(W_1^1, W_2^1, W_3^1\right) \cup \left(W_1^2, W_2^2, W_3^2, W_4^2, W_5^2\right)$$
$$\subseteq (V_1 \cup V_2) \subseteq V$$

is a special semigroup set vector bisubspace of V over the set S.

Now having seen an example of a special semigroup set vector subbispace we now proceed onto define the notion of special semigroup set linear algebra.

**DEFINITION 4.2.3:** *Let*
$$V = (V_1 \cup V_2) = \left(V_1^1, V_2^1, ..., V_{n_1}^1\right) \cup \left(V_1^2, V_2^2, ..., V_{n_2}^2\right)$$
*be a special semigroup set vector bispace over the set S if each $V_i$ is a special semigroup set linear algebra over S then we call $V = (V_1 \cup V_2)$ to be a special set semigroup linear bialgebra.*

We now illustrate this definition by a simple example.

*Example 4.2.3:* Let
$$V = (V_1 \cup V_2) = \left(V_1^1, V_2^1, V_3^1, V_4^1\right) \cup \left(V_1^2, V_2^2, V_3^2, V_4^2, V_5^2\right)$$
where $V_1^1 = S \times S$ where $S = Z^+ \cup \{0\} = Z^o$, $V_2^1 = \{$set of all $3 \times 3$ matrices with entries from $Z^o\}$

$$V_3^1 = \left\{ \begin{bmatrix} a \\ b \\ c \\ d \end{bmatrix} \;\middle|\; a,b,c,d \in Z^o \right\},$$

$$V_4^1 = \left\{ \begin{bmatrix} a & a & a & a & a \\ a & a & a & a & a \end{bmatrix} \;\middle|\; a \in Z^o \right\}$$

is such that $V_1 = \left(V_1^1, V_2^1, V_3^1, V_4^1\right)$ is a special semigroup set linear algebra over $Z^o$. Now let $V_2 = (V_1^2, V_2^2, V_3^2, V_4^2, V_5^2)$ where $V_1^2 = S \times S \times S \times S$, $V_2^2 = \{$all $4 \times 4$ lower triangular matrices with entries from $Z^o\}$, $V_3^2 = \{$set of all polynomials of degree less than or equal to 7 with coefficients from $S\}$,

$$V_4^2 = \left\{ \begin{bmatrix} a & b \\ a & b \\ a & b \\ a & b \end{bmatrix} \middle| a, b \in S \right\}$$

and

$$V_5^2 = \left\{ \begin{bmatrix} a & a & a & a & a \\ b & b & b & b & b \\ c & c & c & c & c \end{bmatrix} \middle| a, b, c \in S \right\}.$$

Clearly $V = V_1 \cup V_2$ is a special semigroup set linear bialgebra over S.

Now we proceed onto define the notion of special semigroup set linear subbialgebra of a special semigroup set linear bialgebra.

**DEFINITION 4.2.4:** *Let*
$$V = (V_1 \cup V_2) = \left( V_1^1, V_2^1, \ldots, V_{n_1}^1 \right) \cup \left( V_1^2, V_2^2, \ldots, V_{n_2}^2 \right)$$
*be a special semigroup set linear bialgebra defined over the additive semigroup S. Suppose*
$$W = (W_1 \cup W_2) = \left( W_1^1, W_2^1, \ldots, W_{n_1}^1 \right) \cup \left( W_1^2, W_2^2, \ldots, W_{n_2}^2 \right)$$
$\subseteq (V_1 \cup V_2)$ *be such that $W_i$ is a special semigroup set linear algebra over the semigroup S of the special semigroup set linear algebra $V_i$, i=1, 2. Then we call $W = (W_1 \cup W_2)$ to be a special semigroup set linear subbialgebra of V over the semigroup S.*

We illustrate this by an example.

*Example 4.2.4:* Let
$$V = (V_1 \cup V_2) = \left( V_1^1, V_2^1, V_3^1 \right) \cup \left( V_1^2, V_2^2, V_3^2, V_4^2 \right)$$
where $V_1^1$ = {set of all 2 × 2 matrices with entries from $Z_5$}, $V_2^1$ = $Z_5 \times Z_5 \times Z_5$,

$$V_3^1 = \left\{ \begin{bmatrix} a_1 & a_2 & a_3 \\ a_4 & a_5 & a_6 \end{bmatrix} \middle| a_i \in Z_5, 1 \leq i \leq 6 \right\},$$

$$V_1^2 = Z_5 \times Z_5 \times Z_5 \times Z_5 \times Z_5,$$

$$V_2^2 = \left\{ \begin{bmatrix} a_1 & a_2 & a_3 \\ a_4 & a_5 & a_6 \\ a_7 & a_8 & a_9 \\ a_{10} & a_{11} & a_{12} \end{bmatrix} \middle| a_i \in Z_5, i = 1,2,3,...,12 \right\},$$

$V_3^2 = \{$All polynomials in the variable x with coefficients from $Z_5$ of degree less than or equal to 6$\}$

and

$$V_4^2 = \left\{ \begin{bmatrix} a & a & a & a & a \\ b & b & b & b & b \\ c & c & c & c & c \end{bmatrix} \middle| a,b,c \in Z_5 \right\}$$

be a special semigroup set linear bialgebra over the semigroup $Z_5$.

Choose

$$W = (W_1 \cup W_2) = \left(W_1^1, W_2^1, W_3^1\right) \cup \left(W_1^2, W_2^2, W_3^2, W_4^2, W_5^2\right)$$
$$\subseteq (V_1 \cup V_2)$$

where $W_1^1 = \{$set of all $2 \times 2$ upper triangular matrices with entries from $Z_5\} \subseteq V_1^1$,

$$W_2^1 = \{Z_5 \times Z_5 \times \{0\}\} \subseteq V_2^1,$$

$$W_3^1 = \left\{ \begin{bmatrix} a & a & a \\ a & a & a \end{bmatrix} \middle| a \in Z_5 \right\} \subseteq V_3^1;$$

$$W_1^2 = \{Z_5 \times \{0\} \times Z_5 \times \{0\} \times Z_5\} \subseteq V_1^2,$$

$$W_2^2 = \left\{ \begin{bmatrix} a & a & a \\ a & a & a \\ a & a & a \\ a & a & a \end{bmatrix} \middle| a \in Z_5 \right\} \subseteq V_2^2,$$

$W_3^2 = $ {all polynomials of degree less than or equal to four in the variable x with coefficients from $Z_5$} $\subseteq V_3^2$,

$$W_4^2 = \left\{ \begin{bmatrix} a & a & a & a & a \\ a & a & a & a & a \\ a & a & a & a & a \end{bmatrix} \middle| a \in Z_5 \right\}.$$

$W = (W_1 \cup W_2) \subseteq (V_1 \cup V_2) \subseteq V$ is a special semigroup set linear subbialgebra of V over the semigroup S.

Now we proceed on to define the notion of special semigroup set linear bitransformation of a special semigroup set vector bispace $V = V_1 \cup V_2$.

**DEFINITION 4.2.5:** *Let*
$$V = (V_1 \cup V_2) = \left(V_1^1, V_2^1, ..., V_{n_1}^1\right) \cup \left(V_1^2, V_2^2, ..., V_{n_2}^2\right)$$
*be a special semigroup set vector bispace over the set S. Suppose $W = (W_1 \cup W_2) = \left(W_1^1, W_2^1, ..., W_{n_1}^1\right) \cup \left(W_1^2, W_2^2, ..., W_{n_2}^2\right)$ be another special semigroup set vector bispace over the same set S. Let $T = T_1 \cup T_2: V = V_1 \cup V_2 \rightarrow W = (W_1 \cup W_2)$ be a bimap such that each $T_i: V_i \rightarrow W_i$; i = 1, 2 is a special semigroup set linear transformation of the special semigroup set vector space $V_i$ into the special semigroup set vector space $W_i$, i = 1, 2 then we define $T = T_1 \cup T_2$ to be a special semigroup set linear bitransformation from V into W. If in the above definition V = W then we call $T = T_1 \cup T_2 : V = V_1 \cup V_2 \rightarrow V = V_1 \cup V_2$ to be a special semigroup set linear bioperator on V.*

*As in case of special set vector spaces we can define pseudo special set linear bioperators on V. Also if the special semigroup set vector bispaces are replaced by special semigroup set linear bialgebras we get the special semigroup set linear bioperators or special semigroup set linear bioperators and so on.*

We illustrate this definition by some examples.

*Example 4.2.5:* Let
$$V = (V_1 \cup V_2) = \left(V_1^1, V_2^1, V_3^1\right) \cup \left(V_1^2, V_2^2, V_3^2, V_4^2\right)$$
and $W = (W_1 \cup W_2) = \left(W_1^1, W_2^1, W_3^1\right) \cup \left(W_1^2, W_2^2, W_3^2, W_4^2\right)$ be two special semigroup set vector bispaces over the same set $S = Z^o$. Let $T = T_1 \cup T_2 \to W_1 \cup W_2$ where $T_1 : V_1 \to W_1$ with

$$T = \left(T_1^1, T_2^1, T_3^1\right) : V_1 = \left(V_1^1, V_2^1, V_3^1\right) \to W_1 = \left(W_1^1, W_2^1, W_3^1\right)$$

and

$$T_2 = \left(T_1^2, T_2^2, T_3^2, T_4^2\right) : V_2 = \left(V_1^2, V_2^2, V_3^2, V_4^2\right)$$
$$\to W_2 = \left(W_1^2, W_2^2, W_3^2, W_4^2\right).$$

We have
$$V_1^1 = \left\{ \begin{pmatrix} a & b \\ c & d \end{pmatrix} \middle| a,b,c,d \in Z^o = Z^+ \cup \{0\} \right\},$$

$$V_2^1 = \{S \times S \times S \times S \times S \times S\},$$

$V_3^1 = \{$All polynomials in the variable x of degree less than or equal to 5 of with coefficients from $Z^o\}$.

$$W_1^1 = \{S \times S \times S \times S\},$$

$$W_2^1 = \left\{ \begin{pmatrix} a & b & c \\ 0 & d & e \\ 0 & 0 & f \end{pmatrix} \middle| a,b,c,d,e,f \in Z^o \right\},$$

$$W_3^1 = \left\{ \begin{pmatrix} a_1 & a_2 & a_3 \\ a_4 & a_5 & a_6 \end{pmatrix} \middle| a_i \in Z^o; 1 \le i \le 6 \right\}.$$

$$V_1^2 = \{S \times S \times S \times S \times S\}$$

$$V_2^2 = \left\{ \begin{bmatrix} a & b \\ c & d \\ e & f \\ g & h \end{bmatrix} \middle| a,b,c,d,e,f,g,h \in S \right\},$$

$V_3^2 = \{S(x)$ is the set all polynomials of degree less than or equal to 8 with coefficients from $S\}$

and

$$V_4^2 = \left\{ \begin{bmatrix} a_1 & a_2 & a_3 & a_4 \\ a_5 & a_6 & a_7 & a_8 \end{bmatrix} \middle| a_i \in Z^o \right\}.$$

Now in $W_2$, $W_1^2 = \{$all polynomials of degree less than or equal to 4 with coefficients from $S\}$,

$$W_2^2 = \left\{ \begin{bmatrix} a & b & c & d \\ e & f & g & h \end{bmatrix} \middle| a,b,c,d,e,f,g,h \in S \right\},$$

$$W_3^2 = \left\{ \begin{bmatrix} a_1 & a_2 & a_3 \\ a_4 & a_5 & a_6 \\ a_7 & a_8 & a_9 \end{bmatrix} \middle| a_i \in S; 1 \le i \le 9 \right\}$$

and

$$W_4^2 = \left\{ \begin{bmatrix} a_1 & a_5 \\ a_2 & a_6 \\ a_3 & a_7 \\ a_4 & a_8 \end{bmatrix} \middle| a_i \in S; 1 \leq i \leq 8 \right\}.$$

Now

$$T = T_1 \cup T_2 = \left(T_1^1, T_2^1, T_3^1\right) \cup \left(T_1^2, T_2^2, T_3^2, T_4^2\right) : V = V_1 \cup V_2$$

into $W = (W_1 \cup W_2)$ is defined as follows

$$T_1 = \left(T_1^1, T_2^1, T_3^1\right) : V_1 = \left(V_1^1, V_2^1, V_3^1\right) \rightarrow \left(W_1^1, W_2^1, W_3^1\right)$$

and

$$T_2 = \left(T_1^2, T_2^2, T_3^2, T_4^2\right):$$
$$V_2 = \left(V_1^2, V_2^2, V_3^2, V_4^2\right) \rightarrow \left(W_1^2, W_2^2, W_3^2, W_4^2\right) = W_2;$$

such that

$T_1^1 : V_1^1 \rightarrow W_1^1$ is defined as

$$T_1^1 \begin{pmatrix} a & b \\ c & d \end{pmatrix} = (a, b, c, d).$$

$T_2^1 : V_2^1 \rightarrow W_2^1$ is given by

$$T_2^1 (a, b, c, d, e, f) = \begin{pmatrix} a & b & c \\ 0 & d & e \\ 0 & 0 & f \end{pmatrix}$$

and $T_3^1 : V_3^1 \rightarrow W_3^1$ is such that

$$T_2^1 (a_0 + a_1 x + a_2 x^2 + a_3 x^3 + a_4 x^5) = \begin{pmatrix} a_0 & a_1 & a_2 \\ a_3 & a_4 & a_5 \end{pmatrix}.$$

$$T_2 = \left(T_1^2, T_2^2, T_3^2, T_4^2\right): V_2 = \left(V_1^2, V_2^2, V_3^2, V_4^2\right) \to$$
$$W_2 = \left(W_1^2, W_2^2, W_3^2, W_4^2\right)$$

is defined as
$T_1^2 : V_1^2 \to W_1^2$ is such that

$$T_1^2 \,(a\ b\ c\ d\ e) = (a + bx + cx^2 + dx^3 + ex^4)$$

$T_2^2 : V_2^2 \to W_2^2$ is given by

$$T_2^2 \begin{bmatrix} a & b \\ c & d \\ e & f \\ g & h \end{bmatrix} = \begin{pmatrix} a & c & e & g \\ b & d & f & h \end{pmatrix}.$$

$T_3^2 : V_3^2 \to W_3^2$ is defined by

$$T_3^2 \,(a_0 + a_1 x + \ldots + a_8 x^8) = \begin{bmatrix} a_0 & a_1 & a_2 \\ a_3 & a_4 & a_5 \\ a_6 & a_7 & a_8 \end{bmatrix}$$

and $T_4^2 : V_4^2 \to W_4^2$ is defined by

$$T_4^2 \begin{bmatrix} a_1 & a_2 & a_3 & a_4 \\ a_5 & a_6 & a_7 & a_8 \end{bmatrix} = \begin{bmatrix} a_1 & a_5 \\ a_2 & a_6 \\ a_3 & a_7 \\ a_4 & a_8 \end{bmatrix}.$$

Thus
$$T = T_1 \cup T_2 = \left(T_1^1, T_2^1, T_3^1\right) \cup \left(T_1^2, T_2^2, T_3^2, T_4^2\right)$$
$$V = V_1 \cup V_2 \to W = W_1 \cup W_2$$
is a special semigroup set linear bitransformation.
   Now if
$T = T_1 \cup T_2 : V = V_1 \cup V_2$

$$= \left(V_1^1, V_2^1, \ldots, V_{n_1}^1\right) \cup \left(V_1^2, V_2^2, \ldots, V_{n_2}^2\right) \rightarrow$$
$$W = (W_1 \cup W_2) = \left(W_1^1, W_2^1, \ldots, W_{n_1}^1\right) \cup \left(W_1^2, W_2^2, \ldots, W_{n_2}^2\right)$$

is such that
$$T_i : V_i \rightarrow W_i; \; T_2 : V_2 \rightarrow W_2;$$
$$T_i = \left(T_1^i, T_2^i, \ldots, T_{n_1}^i\right) : \left(V_1^i, V_2^i, \ldots, V_{n_1}^i\right) \rightarrow \left(W_1^i, W_2^i, \ldots, W_{n_1}^i\right);$$
$$i = 1, 2, \; T_{k_i}^i : V_{k_i}^i \rightarrow W_{k_i}^i \; ; i = 1, 2, \; 1 \leq k_i \leq n_i.$$

Here
$$\text{SHom}(V, W) = \text{SHom}(V_1, W_1) \cup \text{SHom}(V_2, W_2)$$
$$= \{\text{Hom}\left(V_1^1, W_1^1\right), \text{Hom}\left(V_2^1, W_2^1\right), \ldots, \left(V_{n_1}^1, W_{n_1}^1\right)\} \cup \{\text{Hom}\left(V_1^2, W_1^2\right), \left(V_2^2, W_2^2\right), \ldots, \text{Hom}\left(V_{n_2}^2, W_{n_2}^2\right)\}$$

is again a special semigroup set vector bispace over S.

Now we will illustrate by an example the special semigroup set linear bioperator on a special semigroup set vector bispace V.

*Example 4.2.6:* Let
$$V = (V_1 \cup V_2) = \left(V_1^1, V_2^1, V_3^1, V_4^1\right) \cup \left(V_1^2, V_2^2, V_3^2, V_4^2\right)$$
be a special semigroup set vector bispace over the set $S = Z^+ \cup \{0\} = Z^o$.
Here
$$V_1^1 = \left\{ \begin{pmatrix} a & b \\ c & d \end{pmatrix} \middle| a, b, c, d \in Z^o \right\},$$

$$V_2^1 = \{S \times S \times S \times S\},$$

$$V_3^1 = \left\{ \begin{pmatrix} a \\ b \\ c \\ d \\ e \end{pmatrix} \middle| a, b, c, d, e \in S \right\},$$

$V_4^1 = \{$Upper triangular $3\times 3$ matrices$\}$. $V_1^2 = \{S\times S\times S\times S\times S\times S\}$,

$$V_2^2 = \left\{ \begin{bmatrix} a_1 & a_2 & a_3 \\ a_4 & a_5 & a_6 \end{bmatrix} \middle| a_i \in S\, 1\leq i \leq 6 \right\},$$

$V_3^2 = \{$all polynomials of degree less than or equal to 4$\}$,

$$V_4^2 = \left\{ \begin{bmatrix} a & b \\ c & d \\ e & f \\ g & h \\ i & j \end{bmatrix} \middle| a,b,c,d,e,f,g,h,i,j \in Z^o \right\}$$

and
$V_5^2 = \{$All $4\times 4$ lower triangular matrices with entries from $Z^o\}$.

Define the special semigroup set linear bioperator on V by

$$T = T_1 \cup T_2$$
$$= \left( T_1^1, T_2^1, T_3^1, T_4^1 \right) \cup \left( T_1^2, T_2^2, T_3^2, T_4^2, T_5^2 \right):$$
$$V = V_1 \cup V_2 = \left( V_1^1, V_2^1, V_3^1, V_4^1 \right) \cup \left( V_1^2, V_2^2, V_3^2, V_4^2, V_5^2 \right) \rightarrow$$
$$V = V_1 \cup V_2 = \left( V_1^1, V_2^1, V_3^1, V_4^1 \right) \cup \left( V_1^2, V_2^2, V_3^2, V_4^2, V_5^2 \right)$$

as follows:

$T_1^1 : V_1^1 \rightarrow V_1^1$ is defined by
$$T_1^1 \begin{pmatrix} a & b \\ c & d \end{pmatrix} = \begin{pmatrix} b & a \\ d & c \end{pmatrix},$$

$T_2^1 : V_2^1 \rightarrow V_2^1$ is given by
$$T_2^1 \,(a,\, b,\, c,\, d) = (b\ c\ d\ a),$$

$T_3^1 : V_3^1 \rightarrow V_3^1$ is such that

$$T_3^1 \left( \begin{bmatrix} a \\ b \\ c \\ d \\ e \end{bmatrix} \right) = \begin{pmatrix} a \\ a \\ a \\ a \\ a \end{pmatrix}$$

and $T_4^1 : V_4^1 \to V_4^1$ is defined as

$$T_4^1 \begin{pmatrix} a & b & c \\ 0 & d & e \\ 0 & 0 & f \end{pmatrix} = \begin{pmatrix} a & d & f \\ 0 & b & e \\ 0 & 0 & c \end{pmatrix}$$

and $T_1^2 : V_1^2 \to V_1^2$ is such that

$$T_1^2 \,(a\ b\ c\ d\ e\ f) = (\,a\ a\ a\ d\ d\ f),$$

$T_2^2 : V_2^2 \to V_2^2$ is defined by

$$T_2^2 \begin{bmatrix} a_1 & a_2 & a_3 \\ a_4 & a_5 & a_6 \end{bmatrix} = \begin{pmatrix} a_4 & a_5 & a_6 \\ a_1 & a_2 & a_3 \end{pmatrix},$$

$T_3^2 : V_3^2 \to V_3^2$ is such that

$$T_3^2 \,(a_0 + a_1 x + a_2 x^2 + a_3 x^3 + a_4 x^4) = (a_0 + a_2 x^2 + a_4 x^4),$$

$T_4^2 : V_4^2 \to V_4^2$ is given by

$$T_4^2 \begin{bmatrix} a & b \\ c & d \\ e & f \\ g & h \\ i & j \end{bmatrix} = \begin{bmatrix} a & b \\ a & b \\ a & b \\ a & b \\ a & b \end{bmatrix}$$

and $T_5^2 : V_5^2 \to V_5^2$ is defined by

$$T_5^2 \begin{pmatrix} a & 0 & 0 & 0 \\ b & c & 0 & 0 \\ d & e & f & 0 \\ g & h & i & j \end{pmatrix} = \begin{pmatrix} a & 0 & 0 & 0 \\ a & b & 0 & 0 \\ a & b & e & 0 \\ a & b & e & f \end{pmatrix}.$$

Thus we see
$$T = T_1 \cup T_2 = \left(T_1^1, T_2^1, T_3^1, T_4^1\right) \cup \left(T_1^2, T_2^2, T_3^2, T_4^2, T_5^2\right)$$
$$V = V_1 \cup V_2 \to V_1 \cup V_2 = V$$
is a special semigroup set linear bioperator on V.

Further it is still important to note
$$\text{SHom}(V,V) = \{\text{SHom}(V_1,V_1) \cup \text{SHom}(V_2,V_2)\}$$
$$= \{\text{Hom}\left(V_1^1, V_1^1\right), \text{Hom}\left(V_2^1, V_2^1\right), ..., \text{Hom}\left(V_{n_1}^1, V_{n_1}^1\right)\}$$
$$\cup \{\text{Hom}\left(V_1^2, V_1^2\right), \text{Hom}\left(V_2^2, V_2^2\right), ..., \text{Hom}\left(V_{n_2}^2, V_{n_2}^2\right)\}$$
is only a special semigroup set vector bispace over S as even if each $\text{Hom}\left(V_{t_i}^i, V_{t_i}^i\right)$ is a semigroup under the composition of maps yet S over which they are defined is not a semigroup. Thus it is pertinent at this juncture to mention the following. Suppose $V = V_1 \cup V_2 = \left(V_1^1, V_2^1, ..., V_{n_1}^1\right) \cup \left(V_1^2, V_2^2, ..., V_{n_2}^2\right)$ be a special semigroup set linear bialgebra over the semigroup S and if $\text{SHom}_s(V,V) = \text{SHom}_s(V_1,V_1) \cup \text{SHom}_s(V_2,V_2)$ denotes the collection of all special semigroup set linear bioperators on $V = V_1 \cup V_2$ then
$$\text{SHom}_s(V,V) = \text{SHom}_s(V_1 V_1) \cup \text{SHom}_s(V_2,V_2)$$
$$= \{\text{Hom}_s\left(V_1^1, V_1^1\right), \text{Hom}_s\left(V_2^1, V_2^1\right), ..., \text{Hom}_s\left(V_{n_1}^1, V_{n_1}^1\right)\}$$
$$\cup \{\text{Hom}_s\left(V_1^2, V_1^2\right), \text{Hom}_s\left(V_2^2, V_2^2\right), ..., \text{Hom}_s\left(V_{n_2}^2, V_{n_2}^2\right)\}$$
is a special semigroup set linear bialgebra over the semigroup S.

From the notion of special semigroup set linear bialgebra we can have pseudo special semigroup linear bioperator which may not in general be a special semigroup set linear bialgebra. As in case of special semigroup set linear algebras we can define the new substructure in special semigroup set linear bialgebras also.

Further as in case of special semigroup set linear algebras (vector space) we can also in the case of special semigroup set linear bialgebras (vector bispaces) define the notion of special semigroup set projections and the notion of special semigroup set direct union / direct summand. These tasks are left as simple exercises for the reader.

Now we proceed onto generalize these notions to special semigroup set linear n-algebras (vector n-spaces), $n \geq 3$.

**DEFINITION 4.2.6:** *Let*
$$V = V_1 \cup V_2 \cup \ldots \cup V_n$$
$$= \left(V_1^1, V_2^1, \ldots, V_{n_1}^1\right) \cup \left(V_1^2, V_2^2, \ldots, V_{n_2}^2\right) \cup \ldots \cup \left(V_1^n, V_2^n, \ldots, V_{n_n}^n\right)$$
*be such that each $V_i$ is a special semigroup set vector space over the set S for $i = 1, 2, \ldots, n$. Then we call V to be the special semigroup set vector n-space over S, $n \geq 3$. For if $n = 2$ we get the special semigroup set vector bispace, when $n = 3$ we also call the special semigroup set 3-vector space as special semigroup set trivector space or special semigroup set vector trispace.*

We illustrate this by a simple example.

*Example 4.2.7:* Let
$$V = (V_1 \cup V_2 \cup V_3 \cup V_4)$$
$$= \left(V_1^1, V_2^1, V_3^1\right) \cup \left(V_1^2, V_2^2\right) \cup \left(V_1^3, V_2^3, V_3^3, V_4^3\right) \cup \left(V_1^4, V_2^4\right)$$
be a special semigroup set 4-vector space over the set $S = Z_6$ where $V_1^1 = Z_6 \times Z_6 \times Z_6$, $V_2^1 = \{Z_6[x]$ all polynomials of degree less than or equal to 5 with coefficients from $Z_0\}$,

$$V_3^1 = \left\{ \begin{pmatrix} a & b \\ c & d \end{pmatrix} \middle| a, b, c, d \in Z_6 \right\},$$

$$V_1^2 = \{Z_6 \times Z_6 \times Z_6 \times Z_6\},$$

$V_2^2 = $ {all 5×5 lower triangular matrices with entries from $Z_6$},
$V_1^3 = \{Z_6 \times Z_6\}$, $V_2^3 = $ {all 7×2 matrices with entries from $Z_6$},

$$V_3^3 = \left\{ \begin{pmatrix} a \\ b \\ c \end{pmatrix} \middle| a, b, c \in Z_6 \right\},$$

$V_4^3 = $ {all 4×4 upper triangular matrices with entries from $Z_6$},
$$V_1^4 = \{Z_6 \times Z_6 \times Z_6 \times Z_6 \times Z_6\}$$

and

$$V_2^4 = \left\{ \begin{bmatrix} a & a & a \\ a & a & a \\ a & a & a \end{bmatrix} \middle| a \in Z_6 \right\}.$$

It is easily verified that $V = (V_1 \cup V_2 \cup V_3 \cup V_4)$ is a special semigroup set vector 4-space over the set $Z_6$.

Now we proceed on to define the notion of special semigroup set vector n-subspace of V.

**DEFINITION 4.2.7:** *Let*
$$V = V_1 \cup V_2 \cup \ldots \cup V_n$$
$$= \left(V_1^1, V_2^1, \ldots, V_{n_1}^1\right) \cup \left(V_1^2, V_2^2, \ldots, V_{n_2}^2\right) \cup \ldots \cup \left(V_1^n, V_2^n, \ldots, V_{n_n}^n\right)$$
*be a special semigroup set vector n-space over the set S. Take*
$$W = W_1 \cup W_2 \cup \ldots \cup W_n$$
$$= \left(W_1^1, W_2^1, \ldots, W_{n_1}^1\right) \cup \left(W_1^2, W_2^2, \ldots, W_{n_2}^2\right) \cup \ldots \cup \left(W_1^n, W_2^n, \ldots, W_{n_n}^n\right)$$
$$\subseteq V = (V_1 \cup V_2 \cup \ldots \cup V_n)$$
*if each $W_i$ is a special semigroup set vector subspace of $V_i$ for i = 1, 2, …, n, then we call W to be the special semigroup set vector n-subspace of V over the set S.*

We shall illustrate this situation by an example.

*Example 4.2.8:* Let
$$V = (V_1 \cup V_2 \cup V_3 \cup V_4)$$
$$= \left(V_1^1, V_2^1\right) \cup \left(V_1^2, V_2^2, V_3^2, V_4^2\right) \cup \left(V_1^3, V_2^3, V_3^3, V_4^3\right) \cup \left(V_1^4, V_2^4\right)$$
be a special semigroup set vector space over the set $Z_{12}$; where $V_1^1 = Z_{12} \times Z_{12} \times Z_{12}$, $V_2^1 = \{$all $3 \times 3$ matrices with entries from $Z_{12}\}$, $V_1^2 = \{$all $4 \times 4$ lower triangular matrices with entries from $Z_{12}\}$, $V_2^2 = Z_{12} \times Z_{12} \times Z_{12} \times Z_{12}$, $V_3^2 = \{[a\ a\ a\ a\ a\ a\ a] \mid a \in Z_{12}\}$ and

$$V_4^2 = \left\{ \begin{bmatrix} a & a \\ b & b \\ c & c \\ d & d \end{bmatrix} \middle| a,b,c,d \in Z_{12} \right\}.$$

$$V_1^3 = \left\{ \begin{bmatrix} a_1 & a_2 & a_3 & a_4 \\ a_5 & a_6 & a_7 & a_8 \end{bmatrix} \middle| a_i \in Z_{12}; 1 \le i \le 8 \right\},$$

$$V_2^3 = \{Z_{12} \times Z_{12} \times Z_{12} \times Z_{12} \times Z_{12} \times Z_{12}\},$$

$$V_3^3 = \left\{ \begin{pmatrix} a \\ a \\ a \\ a \\ a \end{pmatrix} \middle| a \in Z_{12} \right\}$$

and $V_4^3 = \{$all $5 \times 5$ upper triangular matrices with entries from $Z_{12}\}$. $V_1^4 = \{Z_{12} \times Z_{12}\}$ and $V_2^4 = \{$all $7 \times 7$ upper triangular matrices with entries from $Z_{12}\}$. Take
$$W = (W_1 \cup W_2 \cup W_3 \cup W_4) \subseteq (V_1 \cup V_2 \cup V_3 \cup V_4)$$
where
$$W_1 = \left\{W_1^1, W_2^1\right\} \subset V_1,$$
$$\left\{W_1^2, W_2^2, W_3^2\right\} = W_2 \subseteq V_2,$$

$$W_3 = \left\{ W_1^3, W_2^3, W_3^3, W_4^3 \right\} \subseteq V_3$$

and

$$W_4 = \left\{ W_1^4, W_2^4 \right\} \subseteq V_4$$

with $W_1^1 = \{Z_{12} \times Z_{12} \times \{0\}\} \subseteq V_1^1$, $W_2^1 = \{3 \times 3$ matrices with entries from $\{0, 2, 4, 6, 8, 10\}\} \subseteq V_2^1$, $W_1^2 = \{$all $4 \times 4$ lower triangular matrices with entries from $\{0, 2, 4, 6, 8, 10\}\} \subseteq V_1^2$, $W_2^2 = \{Z_{12} \times Z_{12} \times Z_{12} \times \{0\}\} \subseteq V_2^2$, $W_3^2 = \{[a\ a\ a\ a\ a\ a\ a] \mid a \in \{0, 2, 4, 6, 8, 10\}\} \subseteq V_3^2$,

$$W_4^2 = \left\{ \begin{bmatrix} a & a \\ a & a \\ a & a \\ a & a \end{bmatrix} \middle| a \in Z_{12} \right\} \subseteq V_4^2,$$

$$W_1^3 = \left\{ \begin{bmatrix} a & a & a & a \\ a & a & a & a \end{bmatrix} \middle| a \in Z_{12} \right\} \subseteq V_1^3,$$

$W_2^3 = \{Z_{12} \times \{0\} \times Z_{12} \times \{0\} \times Z_{12} \times \{0\}\} \subseteq V_2^3$,

$$W_3^3 = \left\{ \begin{bmatrix} a \\ a \\ a \\ a \\ a \end{bmatrix} \middle| a \in \{0, 2, 4, 6, 8, 10\} \right\} \subseteq V_3^3$$

and $W_4^3 = \{$all $5 \times 5$ upper triangular matrices with entries from $\{0, 2, 4, 6, 8, 10\}\} \subseteq V_4^3$, $W_1^4 = \{S \times S \mid s = \{0, 3, 6, 9\}\} \subseteq V_1^4$ and $W_2^4 = \{$all $7 \times 7$ upper triangular matrices with entries from $S = \{0, 3, 6, 9\}\} \subseteq V_2^4$. Thus

$$W = (W_1 \cup W_2 \cup W_3 \cup W_4)$$
$$= \{W_1^1, W_2^1\} \cup \{W_1^2, W_2^2, W_3^2, W_4^2\} \cup \{W_1^3, W_2^3, W_3^3, W_4^3\} \cup$$
$$\{W_1^4, W_2^4\} \subseteq \{V_1^1, V_2^1\} \cup \{V_1^2, V_2^2, V_3^2, V_4^2\} \cup$$
$$\{V_1^3, V_2^3, V_3^3, V_4^3\} \cup \{V_1^4, V_2^4\}$$
$$= (V_1 \cup V_2 \cup V_3 \cup V_4) = V$$

is a special semigroup set vector 4-subspace of V over the set $Z_{12}$.

Now we proceed onto define the new notion of special semigroup set n-linear algebras or special semigroup set linear n-algebras $n \geq 3$. For when $n = 2$ we have the special semigroup set linear bialgebra.

**DEFINITION 4.2.8:** *Let $V = (V_1 \cup V_2 \cup ... \cup V_n)$ be a special semigroup set n-vector space over the semigroup S instead of a set S, n-algebra then we call V to be the special semigroup set linear n-algebra.*

*Thus if in the definition of a special semigroup set n-vector space over the set S, we replace the set by a semigroup we get the new notion of special semigroup set n-linear algebra over the semigroup S. It is still important to mention that all special semigroup set n-linear algebras are special semigroup set vector n-spaces, however in general all special semigroup set vector n-spaces are not special semigroup set linear n-algebras.*

We illustrate this by an example.

*Example 4.2.9:* Let
$$V = (V_1 \cup V_2 \cup V_3 \cup V_4 \cup V_5)$$
$$= \{V_1^1, V_2^1, V_3^1\} \cup \{V_1^2, V_2^2\} \cup \{V_1^3, V_2^3, V_3^3, V_4^3\}$$
$$\cup \{V_1^4, V_2^4\} \cup \{V_1^5, V_2^5, V_3^5\}$$

be a special semigroup set n-linear algebra over the semigroup S $= Z^o = Z^+ \cup \{0\}$, where

$$V_1^1 = \{S \times S \times S\},$$

$$V_2^1 = \left\{ \begin{pmatrix} a & b \\ c & d \end{pmatrix} \middle| a,b,c,d \in S \right\},$$

$$V_3^1 = \{(a\ a\ a\ a\ a\ a) \mid a \in S\},$$

$$V_1^2 = \left\{ \begin{pmatrix} a \\ a \\ a \\ a \\ a \end{pmatrix} \middle| a \in S \right\},$$

$$V_2^2 = S \times S \times S \times S,$$

$$V_1^3 = \left\{ \begin{pmatrix} a & a \\ a & a \\ a & a \\ a & a \\ a & a \end{pmatrix} \middle| a \in S \right\},$$

$$V_2^3 = \left\{ \begin{pmatrix} a & a & a & a & a & a \\ a & a & a & a & a & a \end{pmatrix} \middle| a \in S \right\},$$

$$V_3^3 = S \times S \times S \times S \times S,$$

$V_4^3 = \{\text{all } 4 \times 4 \text{ lower triangular matrices with entries from } S\}$, $V_1^4 = \{S \times S \times S\}$, $V_2^4 = \{\text{all } 7 \times 2 \text{ matrices with entries from } S\}$, $V_1^5 = \{\text{all } 3 \times 6 \text{ matrices with entries from } S\}$, $V_2^5 = S \times S \times S \times S$ and $V_3^5 = \{6 \times 3 \text{ matrices with entries from } S\}$.

Clearly V is a special semigroup set linear n-algebra over the semigroup S.

**DEFINITION 4.2.9:** *Let $V = (V_1 \cup V_2 \cup ... \cup V_n)$ be a special semigroup set linear n-algebra over the semigroup S. Take*

$$W = (W_1 \cup W_2 \cup ... \cup W_n)$$
$$= \{W_1^1, W_2^1, ..., W_{n_1}^1\} \cup \{W_1^2, W_2^2, ..., W_{n_2}^2\} \cup ... \cup \{W_1^n, W_2^n, ..., W_{n_n}^n\}$$
$$\subseteq \{V_1^1, V_2^1, ..., V_{n_1}^1\} \cup \{V_1^2, V_2^2, ..., V_{n_2}^2\} \cup ... \cup \{V_1^n, V_2^n, ..., V_{n_n}^n\}$$

*be a proper subset of V. If $W = (W_1 \cup W_2 \cup ... \cup W_n)$ is itself a special semigroup set n-linear algebra over the semigroup S then we call*
$$W = (W_1 \cup W_2 \cup ... \cup W_n)$$
$$= \{W_1^1, W_2^1, ..., W_{n_1}^1\} \cup \{W_1^2, W_2^2, ..., W_{n_2}^2\} \cup ... \cup \{W_1^n, W_2^n, ..., W_{n_n}^n\}$$
*to be the special semigroup set linear n-subalgebra of V over the semigroup S.*

We illustrate this by a simple example.

***Example 4.2.10:*** Let $V = (V_1 \cup V_2 \cup V_3 \cup V_4 \cup V_5)$ be a special semigroup set linear 5-algebra over the semigroup $Z_{15}$ where $V_1 = \{V_1^1, V_2^1, V_3^1\}$ with $V_1^1 = \{Z_{15} \times Z_{15}\}$, $V_2^1 = \{(a_1, a_2, a_3, a_4, a_5) \mid a_i \in Z_{15}\}, 1 \leq i \leq 5\}$,

$$V_3^1 = \left\{ \begin{bmatrix} a_1 \\ a_2 \\ a_3 \\ a_4 \\ a_5 \end{bmatrix} \;\middle|\; a \in Z_{15} \right\},$$

$V_2 = \{V_1^2, V_2^2\}$ where

$$V_1^2 = \left\{ \begin{pmatrix} a & b & c \\ d & e & f \\ g & h & i \end{pmatrix} \;\middle|\; a,b,c,d,e,f,g,h,i \in Z_{15} \right\}$$

and $V_2^2 = Z_{15} \times Z_{15} \times Z_{15} \times Z_{15}$, $V_3 = \{V_1^3, V_2^3, V_3^3, V_4^3\}$

where
$$V_1^3 = Z_{15} \times Z_{15} \times Z_{15},$$

$$V_2^3 = \left\{ \begin{pmatrix} a & a & a & a & a \\ b & b & b & b & b \\ c & c & c & c & c \end{pmatrix} \middle| a, b, c \in Z_{15} \right\},$$

$V_3^3 = \{4 \times 4 \text{ lower triangular matrices with entries for } Z_{15}\}$

and

$$V_4^3 = \left\{ \begin{bmatrix} a_1 & a_4 \\ a_2 & a_5 \\ a_3 & a_6 \end{bmatrix} \middle| a_i \in Z_{15}; 1 \leq i \leq 6 \right\}.$$

$V_4 = \{V_1^4, V_2^4\}$ such that

$$V_1^4 = \left\{ \begin{bmatrix} a_1 \\ a_2 \\ a_3 \\ a_4 \\ a_5 \end{bmatrix} \middle| a_i \in Z_{15}; 1 \leq i \leq 15 \right\}$$

and
$$V_2^4 = Z_{15} \times Z_{15} \times Z_{15} \times Z_{15} \times Z_{15}.$$

$V_5 = \{V_1^5, V_2^5, V_3^5\}$ where
$$V_1^5 = \{Z_{15} \times Z_{15} \times Z_{15} \times Z_{15}\},$$

$$V_2^5 = \left\{ \begin{pmatrix} a_1 & a_2 & a_3 & a_7 \\ a_4 & a_5 & a_6 & a_8 \\ a_9 & a_{10} & a_{11} & a_{12} \end{pmatrix} \middle| a_i \in Z_{15}; 1 \leq i \leq 12 \right\}$$

and
$V_3^5 = \{5 \times 5 \text{ lower triangular matrices with entries from } Z_{15}\}$.
Take
$$W_1 = \{W_1^1, W_2^1, W_3^1\} \subseteq V_1$$

where $W_1^1 = Z_{15} \times \{0\} \subseteq V_1^1$, $W_2^1 = \{(a\ a\ a\ a\ a) \mid a \in Z_{15}\} \subseteq V_2^1$
and

$$W_3^1 = \left\{ \begin{bmatrix} a \\ a \\ a \\ a \\ a \\ a \end{bmatrix} \middle| a \in Z_{15} \right\} \subseteq V_3^1.$$

$W_1$ a special semigroup set linear subalgebra of $V_1$ over $Z_{15}$.

$$W_2 = \{W_1^2, W_2^2\}$$

where

$$W_1^2 = \left\{ \begin{pmatrix} a & a & a \\ a & a & a \\ a & a & a \end{pmatrix} \middle| a \in Z_{15} \right\} \subseteq V_1^2$$

and $W_2^2 = \{Z_{15} \times Z_{15} \times \{0\} \times \{0\}\} \subseteq V_2^2$. Thus $W_2$ is also a special semigroup set linear subalgebra of $V_2$ over $Z_{15}$.

$$W_3 = \{W_1^3, W_2^3, W_3^3, W_4^3\} \subseteq V_3$$

with

$$W_1^3 = Z_{15} \times Z_{15} \times \{0\} \subseteq V_1^3,$$

$$W_2^3 = \left\{ \begin{pmatrix} a & a & a & a & a \\ a & a & a & a & a \\ a & a & a & a & a \end{pmatrix} \middle| a \in Z_{15} \right\} \subseteq V_2^3,$$

$$W_3^3 = \left\{ \begin{pmatrix} a & 0 & 0 & 0 \\ a & a & 0 & 0 \\ a & a & a & 0 \\ a & a & a & a \end{pmatrix} \middle| a \in Z_{15} \right\} \subseteq V_3^3$$

and

$$W_4^3 = \left\{ \begin{pmatrix} a & a \\ a & a \\ a & a \end{pmatrix} \middle| a \in Z_{15} \right\} \subseteq V_4^3.$$

$W_3$ is again a special semigroup set linear subalgebra of $V_3$ over the semigroup $Z_{15}$. $W_4 = \left\{ W_1^4, W_2^4 \right\}$ with

$$W_1^4 = \left\{ \begin{pmatrix} a \\ a \\ a \\ a \\ a \end{pmatrix} \middle| a \in Z_{15} \right\} \subseteq V_1^4$$

and
$$W_2^4 = \{Z_{15} \times \{0\} \times \{0\} \times \{0\} \times Z_{15}\} \subseteq V_2^4$$

so that $W_4 = \left\{ W_1^4, W_2^4 \right\}$ is a special semigroup set linear subalgebra of $V_4$ over the semigroup $Z_{15}$.

Finally
$$W_5 = \left\{ W_1^5, W_2^5, W_3^5 \right\}$$

where

$W_1^5 = \{0\} \times S \times S$ (where $S = \{0, 3, 6, 9, 12\} \subseteq Z_{15}) \subseteq V_1^5$,

$$W_2^5 = \left\{ \begin{pmatrix} a & a & a & a \\ a & a & a & a \\ a & a & a & a \end{pmatrix} \middle| a \in Z_{15} \right\} \subseteq V_2^5$$

and

$$W_3^5 = \left\{ \begin{pmatrix} a & 0 & 0 & 0 & 0 \\ a & b & 0 & 0 & 0 \\ a & b & c & 0 & 0 \\ a & b & c & d & 0 \\ a & b & c & d & a \end{pmatrix} \middle| a, b, c, d \in Z_{15} \right\} \subseteq V_3^5$$

and $W_5$ is a special semigroup set linear subalgebra of $V_5$. Hence $W = W_1 \cup W_2 \cup W_3 \cup W_4 \cup W_5$ is a special semigroup set linear 5-subalgebra of V over the semigroup $S = Z_{15}$.

Now as in case of special semigroup set linear algebra we can also define in case of special semigroup set linear n-algebra the notion of pseudo special semigroup set linear sub n-algebra. We shall define special semigroup set linear n-transformations both in case of special semigroup set vector n-spaces and special semigroup set linear n-algebra.

**DEFINITION 4.2.10:** *Let*
$$V = (V_1 \cup V_2 \cup ... \cup V_n)$$
$$= \{V_1^1, V_2^1, ..., V_{n_1}^1\} \cup \{V_1^2, V_2^2, ..., V_{n_2}^2\} \cup ... \cup \{V_1^n, V_2^n, ..., V_{n_n}^n\}$$
*be a special semigroup set vector n space over the set S and*
$$W = (W_1 \cup W_2 \cup ... \cup W_n)$$
$$= \{W_1^1, W_2^1, ..., W_{n_1}^1\} \cup \{W_1^2, W_2^2, ..., W_{n_2}^2\} \cup ... \cup \{W_1^n, W_2^n, ..., W_{n_n}^n\}$$
*be another special semigroup set vector n-space over the set S. A n-map*
$$T = (T_1 \cup T_2 \cup ... \cup T_n)$$
$$= \{T_1^1, T_2^1, ..., T_{n_1}^1\} \cup \{T_1^2, T_2^2, ..., T_{n_2}^2\} \cup ... \cup \{T_1^n, T_2^n, ..., T_{n_n}^n\}$$
*is said to be a special semigroup set linear n-transformation of V to W if each $T_i = \left(T_1^i, T_2^i, ..., T_{n_i}^i\right)$ is a special semigroup set linear transformation of $V_i$ into $W_i$ for i= 1, 2, ..., n; i.e., $T_{ji}^i : V_{ji}^i \to W_{ji}^i$; $1 \le i \le n$; $1 \le j_i \le n_i$ is a set linear transformation of the vector space $V_{ji}^i$ into $W_{ji}^i$, $1 \le j_i \le n_i$; i =1, 2, ..., n. Now if we replace the special semigroup set n-vector spaces V and W by special semigroup set n-linear algebras then we also we have the special semigroup set linear n-transformation of special semigroup set n-linear algebras. We denote the collection of all special semigroup set linear n-transformations from V into W by SHom (V, W) = {SHom ($V_1$, $W_1$) $\cup$ SHom ($V_2$, $W_2$) $\cup$ ... $\cup$ SHom ($V_n$, $W_n$)}.*

$$= \{(Hom(V_1^1, W_1^1), Hom(V_2^1, W_2^1), ..., Hom(V_{n_1}^1, W_{n_1}^1)\} \cup$$
$$\{Hom(V_1^2, W_1^2), Hom(V_2^2, W_2^2), ..., Hom(V_{n_2}^2, W_{n_2}^2)\} \cup ... \cup$$
$$\{Hom(V_1^n, W_1^n), Hom(V_2^n, W_2^n), ..., Hom(V_{n_n}^n \, W_{n_n}^n)\}.$$

We will illustrate this situation by a simple example.

*Example 4.2.11:* Let
$$V = (V_1 \cup V_2 \cup V_3 \cup V_4)$$
$$= \{V_1^1, V_2^1\} \cup \{V_1^2, V_2^2\} \cup \{V_1^3, V_2^3, V_3^3\} \cup \{V_1^4, V_2^4\}$$
and
$$W = (W_1 \cup W_2 \cup W_3 \cup W_4)$$
$$= \{W_1^1, W_2^1\} \cup \{W_1^2, W_2^2\} \cup \{W_1^3, W_2^3, W_3^3\} \cup \{W_1^4, W_2^4\}$$

be special semigroup set vector 4-space over the set $S = Z^+ \cup \{0\}$.
$$V_1^1 = \{S \times S \times S\},$$

$$V_2^1 = \left\{ \begin{pmatrix} a & b \\ c & d \end{pmatrix} \middle| a,b,c,d \in S \right\},$$

$$V_1^2 = \left\{ \begin{pmatrix} a_1 & a_2 & a_3 & a_4 \\ a_5 & a_6 & a_7 & a_8 \end{pmatrix} \middle| a_i \in S; \ 1 \le i \le 8 \right\}$$

and $V_2^2 = S \times S \times S \times S$. $V_1^3 = \{3 \times 3$ matrices with entries from $S\}$, $V_2^3 = S \times S \times S \times S \times S$, $V_3^3 = \{4 \times 4$ diagonal matrices with entries from $S\}$, $V_1^4 = S \times S \times S \times S \times S \times S$ and $V_2^4 = \{S[x]$ all polynomials of degree less than or equal to 5$\}$.

$$W_1^1 = \left\{ \begin{pmatrix} a & b \\ 0 & c \end{pmatrix} \middle| a,b,c \in S \right\}, \ W_2^1 = S \times S \times S \times S,$$

$W_1^2 = \{$all polynomials of degree less than or equal to 7 with coefficients from $S\}$,

$$W_2^2 = \left\{ \begin{pmatrix} a_1 \\ a_2 \\ a_3 \\ a_4 \end{pmatrix} \middle| \; a_i \in S; 1 \leq i \leq 4 \right\},$$

$W_1^3 = \{3 \times 3$ upper triangular matrices with entries from $S\}$, $W_2^3 = \{$all polynomials of degree less than or equal to 4 with coefficients from $S\}$, $W_3^3 = \{4 \times 4$ lower triangular matrices with entries from $S\}$.

$$W_1^4 = \left\{ \begin{pmatrix} a_1 & a_2 & a_3 \\ a_4 & a_5 & a_6 \end{pmatrix} \middle| \; a_i \in S; \; 1 \leq i \leq 6 \right\}$$

and $W_2^4 = \{6 \times 6$ diagonal matrices with entries from the set $S\}$.
Let
$$T = T_1 \cup T_2 \cup T_3 \cup T_4$$
$$= \{T_1^1, T_2^1\} \cup \{T_1^2, T_2^2\} \cup \{T_1^3, T_2^3, T_3^3\} \cup \{T_1^4, T_2^4\} :$$
$$V = V_1 \cup V_2 \cup V_3 \cup V_4 \to W_1 \cup W_2 \cup W_3 \cup W_4.$$

$T_1^1 : V_1^1 \to W_1^1$ is such that
$$T_1^1 (a\ b\ c) = \begin{pmatrix} a & b \\ 0 & c \end{pmatrix}.$$

$T_2^1 : V_2^1 \to W_2^1$ is given by
$$T_2^1 \begin{pmatrix} a & b \\ c & d \end{pmatrix} = (a\ b\ c\ d).$$

$T_1^2 : V_1^2 \to W_1^2$ is given by

$$T_1^2 \begin{pmatrix} a_1 & a_2 & a_3 & a_4 \\ a_5 & a_6 & a_7 & a_8 \end{pmatrix} =$$
$$(a_1 + a_2 x + a_3 x^2 + a_4 x^3 + a_5 x^4 + a_6 x^5 + a_7 x^6 + a_8 x^7)$$

and
$T_2^2 : V_2^2 \to W_2^2$ is defined by

$$T_2^2 (a\ b\ c\ d) = \begin{bmatrix} a \\ b \\ c \\ d \end{bmatrix}.$$

$T_1^3 : V_1^3 \to W_1^3$ is such that

$$T_1^3 \begin{bmatrix} a & b & c \\ d & e & f \\ g & h & i \end{bmatrix} = \begin{bmatrix} a & b & c \\ 0 & e & f \\ 0 & 0 & i \end{bmatrix},$$

$T_2^3 : V_2^3 \to W_2^3$ is given by

$$T_2^3 (a\ b\ c\ d\ e) = a + bx + cx^2 + dx^3 + ex^4$$

and
$T_3^3 : V_3^3 \to W_3^3$ is such that

$$T_3^3 \begin{bmatrix} a & 0 & 0 & 0 \\ 0 & b & 0 & 0 \\ 0 & 0 & c & 0 \\ 0 & 0 & 0 & d \end{bmatrix} = \begin{bmatrix} a & 0 & 0 & 0 \\ a & b & 0 & 0 \\ a & b & c & 0 \\ a & b & c & d \end{bmatrix}.$$

$T_1^4 : V_1^4 \to W_1^4$ is such that

$$T_1^4 (a\ b\ c\ d\ e\ f) = \begin{bmatrix} a & b & c \\ d & e & f \end{bmatrix}$$

and
$T_2^4 : V_2^4 \to W_2^4$ is defined by
$$T_2^4 (p_0 + p_1 x + p_2 x^2 + p_3 x^3 + p_4 x^4 + p_5 x^5) =$$

$$\begin{bmatrix} p_0 & 0 & 0 & 0 & 0 & 0 \\ 0 & p_1 & 0 & 0 & 0 & 0 \\ 0 & 0 & p_2 & 0 & 0 & 0 \\ 0 & 0 & 0 & p_3 & 0 & 0 \\ 0 & 0 & 0 & 0 & p_4 & 0 \\ 0 & 0 & 0 & 0 & 0 & p_5 \end{bmatrix}.$$

Thus $T = (T_1 \cup T_2 \cup T_3 \cup T_4): V \to W$ is a special semigroup set linear 4 – transformation of V into W. Note if in the special semigroup set linear n-transformation; V into W we replace the range space W by V itself then that specific special semigroup set linear n-transformation from V into V will be known as the special semigroup set linear n-operator on V. We shall denote the special semigroup set linear operator of a special semigroup set linear n-algebra V into itself by

$$\text{SHom}(V, V) = \{ \text{S Hom}(V_1, V_1) \cup \text{S Hom}(V_2, V_2) \cup \ldots \cup \text{S Hom}(V_n, V_n) \}$$
$$= \{(\text{Hom}(V_1^1, V_1^1), \text{Hom}(V_2^1, V_2^1), \ldots, \text{Hom}(V_{n_1}^1, V_{n_1}^1)) \cup$$
$$(\text{Hom}(V_1^2, V_1^2), \text{Hom}(V_2^2, V_2^2), \ldots, \text{Hom}(V_{n_2}^2, V_{n_2}^2)) \cup \ldots \cup$$
$$(\text{Hom}(V_1^n, V_1^n), \ldots, \text{Hom}(V_{n_1}^n, V_{n_1}^n))\} .$$

Clearly SHom(V, V) is again a special semigroup set linear n-algebra over the semigroup. However if we replace the special semigroup set linear n-algebra V by a special semigroup set vector n-space then we see SHom(V,V) is not a special semigroup set linear n-algebra it is only a special semigroup set vector n-space.

We shall give one example of a special semigroup set linear n-operator of V.

*Example 4.2.12:* Let $V = (V_1 \cup V_2 \cup V_3)$ where
$$V_1 = \{V_1^1, V_2^1, V_3^1, V_4^1\},$$

$$V_2 = \{V_1^2, V_2^2, V_3^2\}$$

and

$$V_3 = \{V_1^3, V_2^3\}$$

defined by $V_1^1 = \{S \times S \times S \times S \mid S = Z^o = Z^+ \cup \{0\}\}$, $V_2^1 = \{$set of all $3 \times 3$ matrices with entries from S$\}$,

$$V_3^1 = \left\{ \begin{bmatrix} a_1 & a_2 \\ a_3 & a_4 \\ a_5 & a_6 \\ a_7 & a_8 \\ a_9 & a_{10} \end{bmatrix} \middle| a_i \in Z^o; 1 \leq i \leq 10 \right\},$$

$V_4^1 = \{$set of all polynomials in the variable x with coefficients from S of degree less than or equal to 9$\}$, $V_1^2 = S \times S \times S$, $V_2^2 = \{$all $4 \times 4$ low triangular matrices with entries from S$\}$ and

$$V_3^2 = \left\{ \begin{bmatrix} a_1 & a_2 & a_3 & a_4 & a_5 \\ a_6 & a_7 & a_8 & a_9 & a_{10} \end{bmatrix} \middle| a_i \in S; 1 \leq i \leq 10 \right\}.$$

$V_1^3 = S \times S \times S \times S \times S \times S$ and $V_2^3 = \{$all $5 \times 5$ upper triangular matrices with entries from S$\}$, be a special semigroup set vector 3-space over the set S.

Define $T = T_1 \cup T_2 \cup T_3$
$$= \{T_1^1, T_2^1, T_3^1, T_4^1\} \cup \{T_1^2, T_2^2, T_3^2\} \cup \{T_1^3, T_2^3\}:$$
$$V = V_1 \cup V_2 \cup V_3$$
$$= \{V_1^1, V_2^1, V_3^1, V_4^1\} \cup \{V_1^2, V_2^2, V_3^2\} \cup \{V_1^3, V_2^3\} \rightarrow$$
$$V = V_1 \cup V_2 \cup V_3 = \{V_1^1, V_2^1, V_3^1, V_4^1\} \cup \{V_1^2, V_2^2, V_3^2\} \cup \{V_1^3, V_2^3\}$$

as follows:
$$T_i : V_i \rightarrow V_i$$

where

$T^i_j : V^i_j \to V^i_j$ ; $1 \leq j \leq n_i$, $i = 1, 2, 3$

$T^1_1 : V^1_1 \to V^1_1$ is such that

$$T^1_1 \text{ (a b c d)} = \text{(a a a a)},$$

$T^1_2 : V^1_2 \to V^1_2$ is defined by

$$T^1_2 \begin{pmatrix} a & b & c \\ d & e & f \\ g & h & i \end{pmatrix} = \begin{pmatrix} a & a & a \\ a & a & a \\ b & b & b \end{pmatrix},$$

$T^1_3 : V^1_3 \to V^1_3$ is given by

$$T^1_3 \begin{pmatrix} a_1 & a_2 \\ a_3 & a_4 \\ a_5 & a_6 \\ a_7 & a_8 \\ a_9 & a_{10} \end{pmatrix} = \begin{pmatrix} a & a \\ a & a \\ a & a \\ a & a \\ a & a \end{pmatrix}$$

and
$T^1_4 : V^1_4 \to V^1_4$ is defined by

$$T^1_4 (a_0 + a_1 x + \ldots a_9 x^9) = a^0 + a_3 x^3 + a_6 x^6 + a_9 x^9.$$

$T^2_1 : V^2_1 \to V^2_1$ is given by

$$T^2_1 \text{(a b c)} = \text{(b c a)}$$

and
$T^2_2 : V^2_2 \to V^2_2$ is defined by

$$T^2_2 \begin{pmatrix} a & 0 & 0 & 0 \\ b & c & 0 & 0 \\ d & e & f & 0 \\ g & h & i & j \end{pmatrix} = \begin{pmatrix} a & 0 & 0 & 0 \\ a & a & 0 & 0 \\ a & a & a & 0 \\ a & a & a & a \end{pmatrix}.$$

$T_3^2 : V_3^2 \to V_3^2$ is defined by

$$T_3^2 \begin{bmatrix} a_1 & a_2 & a_3 & a_4 & a_5 \\ a_6 & a_7 & a_8 & a_9 & a_{10} \end{bmatrix} = \begin{bmatrix} a_6 & a_7 & a_8 & a_9 & a_{10} \\ a_1 & a_2 & a_3 & a_4 & a_5 \end{bmatrix}.$$

$T_1^3 : V_1^3 \to V_1^3$ is given by

$$T_1^3 (a\ b\ c\ d\ e\ f) = (a\ a\ a\ a\ a\ a).$$

$T_2^3 : V_2^3 \to V_2^3$ is defined by

$$T_2^3 \begin{bmatrix} a & b & c & d & e \\ 0 & f & g & h & i \\ 0 & 0 & j & k & l \\ 0 & 0 & 0 & m & n \\ 0 & 0 & 0 & 0 & p \end{bmatrix} = \begin{bmatrix} a & a & a & a & a \\ 0 & b & b & b & b \\ 0 & 0 & c & c & c \\ 0 & 0 & 0 & d & d \\ 0 & 0 & 0 & 0 & e \end{bmatrix}.$$

Clearly $T = T_1 \cup T_2 \cup T_3: V \to V$ is a special semigroup set linear 3-operator of V.

We can also define special semigroup set pseudo linear n-operators on V. Now in case of special semigroup set n-vector spaces we can define the notion of direct sum of special semigroup set vector n-spaces. Also we can define the notion of special semigroup set linear n-idempotent operator and special semigroup set linear n-projections on V.

Now we proceed onto show by an example how a special semigroup set vector n-space is represented as a direct sum of special semigroup set vector n-subspaces.

*Example 4.2.13:* Let
$$V = (V_1 \cup V_2 \cup V_3 \cup V_4)$$
$$= \{V_1^1, V_2^1, V_3^1\} \cup \{V_1^2, V_2^2\} \cup \{V_1^3, V_2^3\} \cup \{V_1^4, V_2^4, V_3^4\}$$

be a special semigroup set linear 4-algebra over the set $Z^o = S = Z^+ \cup \{0\}$. Here $V_1^1 = S \times S \times S \times S$,

$$V_2^1 = \left\{ \begin{pmatrix} a & b \\ c & d \end{pmatrix} \middle| a,b,c,d \in S \right\},$$

$V_3^1 = $ {all polynomials in x of degree less than or equal 5 with coefficients from S}, $V_1^2 = S \times S \times S$ and

$$V_2^2 = \left\{ \begin{bmatrix} a_1 & a_2 & a_3 \\ a_4 & a_5 & a_6 \end{bmatrix} \middle| a_i \in S; 1 \leq i \leq 6 \right\},$$

$$V_1^3 = \left\{ \begin{bmatrix} a_1 & a_2 \\ a_3 & a_4 \\ a_5 & a_6 \end{bmatrix} \middle| a_i \in S; 1 \leq i \leq 6 \right\}$$

and
$$V_2^3 = S \times S \times S \times S \times S.$$

$$V_1^4 = \left\{ \begin{pmatrix} a & b \\ c & d \end{pmatrix} \middle| a,b,c,d \in S \right\},$$

$$V_2^4 = S \times S \times S \times S$$

and

$$V_3^4 = \left\{ \begin{bmatrix} a_1 \\ a_2 \\ a_3 \end{bmatrix} \middle| a_i \in S; 1 \leq i \leq 3 \right\}.$$

We now represent each semigroup set linear algebra as a direct sum of semigroup set linear subalgebras.

$$V_1^1 = S \times \{0\} \times S \times \{0\} \oplus \{0\} \times S \times \{0\} \times \{0\}$$
$$\oplus \{0\} \times \{0\} \times S \times \{0\}$$

$$= W_{11}^1 \oplus W_{21}^1 \oplus W_{31}^1.$$

$$V_2^1 = \left\{ \begin{pmatrix} a & 0 \\ c & 0 \end{pmatrix} \middle| a,c \in S \right\} \oplus \left\{ \begin{pmatrix} 0 & 0 \\ 0 & d \end{pmatrix} \middle| d \in S \right\} \oplus \left\{ \begin{pmatrix} 0 & b \\ 0 & 0 \end{pmatrix} \middle| b \in S \right\}$$

$$= W_{12}^1 \oplus W_{22}^1 \oplus W_{32}^1.$$

$V_1^3 = $ {all polynomials of degree less than or equal to 3 with coefficients from S} $\oplus$ {all polynomials of degree strictly greater than 3 but is of degree less than or equal to 5}

$$= W_{13}^1 \oplus W_{23}^1.$$

Now

$$V_2^2 = \left\{ \begin{pmatrix} a_1 & a_2 & a_3 \\ 0 & 0 & 0 \end{pmatrix} \middle| a_1, a_2, a_3 \in S \right\} \oplus$$

$$\left\{ \begin{pmatrix} 0 & 0 & 0 \\ a_4 & a_5 & a_6 \end{pmatrix} \middle| a_4, a_5, a_6 \in S \right\} = W_{12}^2 \oplus W_{22}^2.$$

$$V_1^3 = \left\{ \begin{pmatrix} a_1 & a_2 \\ 0 & 0 \\ 0 & 0 \end{pmatrix} \middle| a_1, a_2 \in S \right\} \oplus \left\{ \begin{pmatrix} 0 & 0 \\ a_3 & a_4 \\ 0 & 0 \end{pmatrix} \middle| a_3, a_4 \in S \right\} \oplus$$

$$\left\{ \begin{pmatrix} 0 & 0 \\ 0 & 0 \\ a_5 & a_6 \end{pmatrix} \middle| a_5, a_6 \in S \right\} = W_{11}^3 \oplus W_{21}^3 \oplus W_{31}^3.$$

$$V_2^3 = \{S \times \{0\} \times \{0\} \times \{0\} \times S\} \oplus$$
$$\{\{0\} \times S \times \{0\} \times \{S\} \times \{0\}\} \oplus$$
$$\{0\} \times \{0\} \times S \times \{0\} \times \{0\}$$
$$= W_{12}^3 \oplus W_{22}^3 \oplus W_{32}^3.$$

$$V_1^4 = \left\{ \begin{pmatrix} a & 0 \\ 0 & 0 \end{pmatrix} \middle| a \in S \right\} \oplus \left\{ \begin{pmatrix} 0 & b \\ 0 & 0 \end{pmatrix} \middle| b \in S \right\} \oplus$$

$$\left\{ \begin{pmatrix} 0 & 0 \\ c & 0 \end{pmatrix} \middle| c \in S \right\} \oplus \left\{ \begin{pmatrix} 0 & 0 \\ 0 & d \end{pmatrix} \middle| d \in S \right\}$$

$$= W_{11}^4 \oplus W_{21}^4 \oplus W_{31}^4 \oplus W_{41}^4.$$

$$V_2^4 = S \times \{0\} \times \{0\} \times S \oplus$$
$$\{0\} \times S \times \{0\} \times \{0\} \oplus$$
$$\{0\} \times \{0\} \times S \times \{0\} = W_{12}^4 \oplus W_{22}^4 \oplus W_{32}^4.$$

$$V_3^4 = \left\{ \begin{bmatrix} a_1 \\ a_2 \\ 0 \end{bmatrix} \middle| a_1, a_2 \in S \right\} \oplus \left\{ \begin{bmatrix} 0 \\ 0 \\ a_3 \end{bmatrix} \middle| a_3 \in S \right\}$$

$$= W_{13}^4 \oplus W_{23}^4.$$

Thus
$$V = (V_1 \cup V_2 \cup V_3 \cup V_4)$$
$$= ( W_{11}^1 \oplus W_{21}^1 \oplus W_{31}^1, \ W_{12}^1 \oplus W_{22}^1 \oplus W_{32}^1, \ W_{13}^1 \oplus W_{23}^1 )$$
$$\cup ( W_{11}^2 \oplus W_{21}^2 \oplus W_{12}^2 \oplus W_{22}^2 \cup W_{11}^3 \oplus W_{21}^3 \oplus W_{31}^3,$$
$$W_{12}^3 \oplus W_{22}^3 \oplus W_{32}^3 ) \cup ( W_{11}^4 \oplus W_{21}^4 \oplus W_{31}^4 \oplus W_{41}^4,$$
$$W_{12}^4 \oplus W_{22}^4 \oplus W_{32}^4, \ W_{13}^4 \oplus W_{23}^4 )$$

is a special semigroup set direct sum representation of the special set semigroup linear 4-algebra over the semigroup $S = Z^o = Z^+ \cup \{0\}$.

However we can represent in the same way the special semigroup set vector n-space as direct sum. When we get a representation of special semigroup set vector n-space as direct sum we can define special semigroup set linear operator on V which are projections P and these projection are defined as follows:

$$P: V \to (W_1 \cup W_2 \cup \ldots \cup W_n)$$

are such that P o P = P.

The interested reader is requested to construct examples.

Now we proceed onto introduce special group set vector n-spaces and special group set linear n-algebras.

**DEFINITION 4.2.11:** *Let*
$$V = (V_1 \cup V_2 \cup \ldots \cup V_n)$$
$$= \{V_1^1, V_2^1, \ldots, V_{n_1}^1\} \cup \{V_1^2, V_2^2, \ldots, V_{n_2}^2\} \cup \ldots \cup \{V_1^n, V_2^n, \ldots, V_{n_n}^n\}$$
*be such that each $V_i$ is a special group set vector space over a set S and each $V_{t_i}^i$ is a group set vector space over the same set S, $i \leq t_i \leq n_i$; $i = 1, 2, \ldots, n$.*

*Then we call V to be a special group set n-vector space (special group set vector n-space) over the set S.*

We shall illustrate this by an example.

*Example 4.2.14:* Let
$$V = (V_1 \cup V_2 \cup V_3 \cup V_4)$$
$$= \{V_1^1, V_2^1\} \cup \{V_1^2, V_2^2, V_3^2\} \cup \{V_1^3, V_2^3\} \cup \{V_1^4, V_2^4, V_3^4, V_4^4\}$$
where each $V_{j_i}^i$ is defined as

$$V_1^1 = \left\{ \begin{pmatrix} a & b \\ c & d \end{pmatrix} \middle| a,b,c,d \in S = Z \right\},$$

$$V_2^1 = S \times S \times S, \quad V_1^2 = S \times S \times S \times S,$$

$$V_2^2 = \left\{ \begin{bmatrix} a_1 & a_2 & a_3 \\ a_4 & a_5 & a_6 \end{bmatrix} \middle| a_i \in S; \quad 1 \leq i \leq 6 \right\},$$

$V_3^2 = $ {all $3 \times 3$ upper triangular matrices with entries from Z},
$$V_1^3 = S \times S \times S,$$

$$V_2^3 = \left\{ \begin{bmatrix} a \\ b \\ c \\ d \end{bmatrix} \middle| a,b,c,d \in S \right\},$$

$$V_1^4 = S \times S \times S \times S \times S,$$

$V_2^4 = \{\text{all } 4 \times 4 \text{ lower triangular matrices with entries from Z}\}$,

$$V_3^4 = \left\{ \begin{bmatrix} a & a & a & a & a \\ a & a & a & a & a \end{bmatrix} \middle| a \in Z \right\}$$

and

$$V_4^4 = \left\{ \begin{bmatrix} a_1 & a_2 \\ a_3 & a_4 \\ a_5 & a_6 \end{bmatrix} \middle| a_i \in Z;\ 1 \le i \le 6 \right\}.$$

Thus $V = (V_1 \cup V_2 \cup V_3 \cup V_4)$ is a special group set vector 4 space over the set S.

It is interesting to observe in the definition when $n = 2$ we get the special group set vector bispace and when $n = 3$ we get the special group set vector trispace.

Now we proceed onto define the notion of special group set vector n-subspace W of a special group set vector n-space V defined over the set S.

**DEFINITION 4.2.12:** *Let*
$$V = (V_1 \cup V_2 \cup \ldots \cup V_n)$$
$$= \{V_1^1, V_2^1, \ldots, V_{n_1}^1\} \cup \{V_1^2, V_2^2, \ldots, V_{n_2}^2\} \cup \ldots \cup \{V_1^n, V_2^n, \ldots, V_{n_n}^n\}$$
*be a special group set vector n-space defined over the set S. Let*
$$W = (W_1 \cup W_2 \cup \ldots \cup W_n)$$
$$= \{W_1^1, W_2^1, \ldots, W_{n_1}^1\} \cup \{W_1^2, W_2^2, \ldots, W_{n_2}^2\} \cup \ldots \cup \{W_1^n, W_2^n, \ldots, W_{n_n}^n\}$$
$$\subset (V_1 \cup \ldots \cup V_n)$$

i.e., $W_i \subseteq V_i$ for every i, $1 \leq i \leq n$ and $W_{t_i}^i \subseteq V_{t_i}^i$ is such that $W_{t_i}^i$ is a proper subgroup of $V_{t_i}^i$, $1 \leq t_i \leq n_i$, i =1, 2, ..., n. Then we call W to be a special group set vector n-subspace of V over S if $W = (W_1 \cup ... \cup W_n)$ is itself a special group set vector n-space over the set S.

We shall represent this by a simple example.

*Example 4.2.15:* Let
$$V = (V_1 \cup V_2 \cup V_3 \cup V_4)$$
$$= \{V_1^1, V_2^1, V_3^1\} \cup \{V_1^2, V_2^2, V_3^2, V_4^2\} \cup \{V_1^3, V_2^3\} \cup \{V_1^4, V_2^4, V_3^4\}$$
be a special group set vector 4-space over the set $Z_{12}$ where

$$V_1^1 = \left\{ \begin{pmatrix} a & b \\ c & d \end{pmatrix} \middle| a,b,c,d \in Z_{12} \right\},$$

$V_2^1 = Z_{12} \times Z_{12} \times Z_{12}$, $V_3^1 = $ {all 7×2 matrices with entries from $Z_{12}$}, $V_1^2 = Z_{12} \times Z_{12} \times Z_{12} \times Z_{12}$,

$$V_2^2 = \left\{ \begin{pmatrix} a & b \\ 0 & d \end{pmatrix} \middle| a,b,d \in Z_{12} \right\},$$

$V_3^2 = $ {all 2 × 9 matrices with entries from $Z_{12}$}, $V_4^2 = $ {(a a a a a) | a ∈ $Z_{12}$}, $V_1^3 = Z_{12} \times Z_{12}$, $V_2^3 = $ {all 3×3 matrices with entries from $Z_{12}$}, $V_1^4 = $ {all 4 × 4 lower triangular matrices with entries from $Z_{12}$}, $V_2^4 = Z_{12} \times Z_{12} \times Z_{12} \times Z_{12}$
and
$V_3^4 = $ {6 × 3 matrices with entries from $Z_{12}$}.

Consider a proper subset

$$W = (W_1 \cup W_2 \cup W_3 \cup W_4)$$

$$= \{W_1^1, W_2^1, W_3^1\} \cup \{W_1^2, W_2^2, W_3^2, W_4^2\} \cup$$
$$\{W_1^3, W_2^3\} \cup \{W_1^4, W_2^4, W_3^4\}$$
$$\subseteq V_1 \cup V_2 \cup V_3 \cup V_4$$

where

$$W_1^1 = \left\{ \begin{pmatrix} a & b \\ c & d \end{pmatrix} \middle| a,b,c,d \in \{0,2,4,6,8,10\} = T \subseteq Z_{12} \right\} \subseteq V_1^1$$

is a subgroup of $V_1^1$, $W_2^1 = Z_{12} \times Z_{12} \times \{0\} \subseteq V_2^1$ is a subgroup of $V_2^1$, $W_3^1 = \{$all $7 \times 2$ matrices with entries from $T \subseteq Z_{12}\} \subseteq V_3^1$ is again a subgroup of $V_3^1$, $W_1^2 = \{Z_{12} \times Z_{12} \times \{0\} \times \{0\}\} \subseteq V_1^2$ is a subgroup of $V_1^2$,

$$W_2^2 = \left\{ \begin{pmatrix} a & a \\ 0 & a \end{pmatrix} \middle| a \in Z_{12} \right\} \subseteq V_2^2$$

is a subgroup of $V_2^2$, $W_3^2 = \{2 \times 9$ matrices with entries from the set $P = \{0, 6, 9\} \subseteq Z_{12}\} \subseteq V_3^2$ is a subgroup of $V_3^2$, $W_4^2 \subseteq V_4^2$ is such that $W_4^2 = \{(a\ a\ a\ a\ a) \mid a \in \{0, 6, 9\} \subseteq Z_{12}\} \subseteq V_4^2$ is a subgroup of $V_4^2$, $W_1^3 = Z_{12} \times P \mid P = \{0, 3, 6, 9\} \subseteq Z_{12}\} \subseteq V_1^3$ is a subgroup of $V_1^3$, $W_2^3 = \{$all $3 \times 3$ matrices with entries from the set $B = \{0, 6\} \subseteq Z_{12}\} \subseteq V_2^3$, $W_1^4 = \{$all $4 \times 4$ lower triangular matrices with entries from $T = \{0, 2, 4, 6, 8, 10\} \subseteq Z_{12}\} \subseteq V_1^4$ is a subgroup of $V_1^4$, $W_2^4 = Z_{12} \times Z_{12} \times \{0\} \times \{0\} \subseteq V_2^4$ is a subgroup of $V_2^4$, $W_3^4 = \{$all $6 \times 3$ matrices with entries from $T = \{0, 2, 4, 6, 8, 10\} \subseteq Z_{12}\} \subseteq V_3^4$ is a subgroup of $V_3^4$. Thus $W = W_1 \cup W_2 \cup W_3 \cup W_4$ is a special group set vector 4-subspace of V over $Z_{12}$.

Now we proceed onto define the new notion of special group semi n-vector spaces which are simple.

**DEFINITION 4.2.13:** *Let*
$$V = (V_1 \cup V_2 \cup \ldots \cup V_n)$$
$$= \{V_1^1, V_2^1, \ldots, V_{n_1}^1\} \cup \{V_1^2, V_2^2, \ldots, V_{n_2}^2\} \cup \ldots \cup \{V_1^n, V_2^n, \ldots, V_{n_n}^n\}$$
*be a special group set vector n-space over the set S. If each of the groups $V_{j_i}^i$ are simple for $1 \le j_i \le n_i$, $i = 1, 2, \ldots, n$, then we call V to be a special group set simple vector n-space. If in the special group set vector n-space each $V_{j_i}^i$ ($1 \le j_i \le n_i$, $i = 1, 2, \ldots, n$) does not have proper subgroups then we call V to be a special group set strong simple vector n-space.*

We denote these concepts by an example.

***Example 4.2.16:*** Let
$$V = (V_1 \cup V_2 \cup V_3)$$
$$= \{V_1^1, V_2^1\} \cup \{V_1^2, V_2^2\} \cup \{V_1^3, V_2^3, V_3^3\}$$
where
$$V_1^1 = Z_{11}, \; V_2^1 = Z_{13},$$

$$V_1^2 = \left\{ \begin{pmatrix} a & 0 \\ 0 & 0 \end{pmatrix} \middle| a \in Z_5 \right\},$$

$$V_2^2 = \left\{ \begin{pmatrix} a & 0 & 0 & 0 \\ 0 & 0 & 0 & a \end{pmatrix} \middle| a \in Z_2 \right\},$$

$V_1^3 = Z_5$, $V_2^3 = Z_{17}$ and $V_3^3 = Z_{19}$ be a special group set vector n space over the set $S = \{0, 1\}$, $n = 3$. We see each $V_{j_i}^i$ has no subgroup so V is a special group set strong simple vector 3-space over the set $S = \{0, 1\}$.

Now we proceed onto define the notion of special group set linear n-algebra over a group G.

**DEFINITION 4.2.14:** *Let $V = V_1 \cup V_2 \cup ... \cup V_n$ be a special group set vector n-space over a group S under addition (instead of a set S) then we call V to be a special group set linear n-algebra over the group S. If n = 2 we call it the special group set linear bialgebra when n = 3 we call V as the special group set linear trialgebra.*

We shall illustrate this definition by an example.

*Example 4.2.17:* Let
$$V = (V_1 \cup V_2 \cup V_3 \cup V_4 \cup V_5)$$
$$= \{V_1^1, V_2^1\} \cup \{V_1^2, V_2^2, V_3^2\} \cup \{V_1^3, V_2^3, V_3^3\} \cup \{V_1^4, V_2^4\} \cup$$
$$\{V_1^5, V_2^5, V_3^5, V_4^5\}$$

be a special group set linear 5-algebra over the group $G = Z$ (the integers under additions); where $V_1^1 = \{Z \times Z \times Z\}$, $V_2^1 = \{$all $3 \times 3$ matrices with entries from Z$\}$, $V_1^2 = \{$all $2 \times 5$ matrices with entries from Z$\}$, $V_2^2 = \{$all $4 \times 4$ upper triangular matrices with entries from Z$\}$, $V_3^2 = Z \times Z \times Z \times Z \times Z$, $V_1^3 = Z \times Z$, $V_2^3 = \{7 \times 2$ matrices with entries from Z$\}$, $V_3^3 = \{$all polynomials in Z[x] of degree less than or equal to 5$\}$, $V_1^4 = \{Z \times Z \times Z \times Z\}$, $V_2^4 = \{7 \times 7$ lower triangular matrices with entries from Z$\}$, $V_1^5 = Z \times Z \times Z$, $V_2^5 = \{(a\ a\ a\ a\ a\ a\ a) \mid a \in Z\}$, $V_3^5 = \{$all $9 \times 9$ diagonal matrices with entries from Z$\}$ and $V_4^5 = \{$all $3 \times 5$ matrices with entries from Z$\}$. It can be easily verified that each $V_{j_i}^i$ is a group under addition $1 \le j_i \le n_i$, i = 1, 2, 3, 4, 5. Thus V $= (V_1 \cup V_2 \cup V_3 \cup V_4 \cup V_5)$ is a special group set linear 5-algebra over the group Z.

Now we proceed onto define their substructures.

**DEFINITION 4.2.15:** *Let*
$$V = (V_1 \cup V_2 \cup ... \cup V_n)$$
$$= \{V_1^1, V_2^1, ..., V_{n_1}^1\} \cup \{V_1^2, V_2^2, ..., V_{n_2}^2\} \cup ... \cup \{V_1^n, V_2^n, ..., V_{n_n}^n\}$$

*be a special group set linear n-algebra over a group G. Let*
$$W = (W_1 \cup W_2 \cup ... \cup W_n)$$
$$= \{W_1^1, W_2^1, ..., W_{n_1}^1\} \cup \{W_1^2, W_2^2, ..., W_{n_2}^2\} \cup ... \cup \{W_1^n, W_2^n, ..., W_{n_n}^n\}$$
$$\subseteq (V_1 \cup V_2 \cup ... \cup V_n)$$
*where $W_i \subseteq V_i$ for $i = 1, 2, ..., n$ such that each $W_i$ is a special set group linear subalgebra of $V_i$ i.e., each component $W_{j_i}^i \to V_{j_i}^i$ is a group set linear subalgebra of $V_{j_i}^i$; $1 \leq j_i \leq n_i$, $i = 1, 2, ..., n$. Then we call $W = (W_1 \cup W_2 \cup ... \cup W_n) \subseteq (V_1 \cup V_2 \cup ... \cup V_n) = V$ to be the special group set linear n-subalgebra of V over the group G.*

We will illustrate this by a simple example.

***Example 4.2.18:*** Let
$$V = (V_1 \cup V_2 \cup V_3 \cup V_4)$$
$$= \{V_1^1, V_2^1, V_3^1\} \cup \{V_1^2, V_2^2, V_3^2, V_4^2\} \cup \{V_1^3, V_2^3\} \cup$$
$$\{V_1^4, V_2^4, V_3^4, V_4^4, V_5^4\}$$

be a special group set linear 4-algebra over the group Z where $V_1^1 = Z \times Z \times Z \times Z$, $V_2^1 = \{$all $2 \times 2$ matrices with entries from $Z\}$, $V_3^1 = \{Z(x)$, all polynomials of degree less than or equal to $8\}$, $V_1^2 = Z \times Z \times Z$, $V_2^2 = \{3 \times 7$ matrices with entries from $Z\}$, $V_3^2 = \{$(a a a a a a) $| a \in Z\}$, $V_4^2 = \{4 \times 3$ matrices with entries from $Z\}$, $V_1^3 = Z \times Z$, $V_2^3 = \{4 \times 4$ upper triangular matrices with entries from $Z\}$, $V_1^4 = \{Z \times Z \times Z \times Z \times Z\}$, $V_2^4 = 3 \times 3$ lower triangular matrices with entries from $Z\}$, $V_3^4 = 3Z \times 5Z \times 11Z$,

$$V_4^4 = \left\{ \begin{bmatrix} a \\ a \\ a \\ a \\ a \\ a \\ a \end{bmatrix} \middle| a \in Z \right\}$$

and $V_5^4 = \{5 \times 2$ matrices with entries from $Z\}$. Take

$$W = (W_1 \cup W_2 \cup W_3 \cup W_4) = \{W_1^1, W_2^1, W_3^1\} \cup$$
$$\{W_1^2, W_2^2, W_3^2, W_4^2\} \cup \{W_1^3, W_2^3\} \cup \{W_1^4, W_2^4, W_3^4, W_4^4, W_5^4\}$$
$$\subseteq (V_1 \cup V_2 \cup V_3 \cup V_4)$$

be such that $W_1^1 = Z \times Z \times \{0\} \times Z \subseteq V_1^1$ is a subgroup of the group $V_1^1$, $W_2^1 = \{2 \times 2$ matrices with entries from $3Z\}$ is again a subgroup of the group $V_2^1$, $W_3^1 = \{$all polynomials in x of degree less than or equal to 4 with coefficients from $Z\}$ is again a subgroup of the group $V_3^1$. Thus $W_1 = \{W_1^1, W_2^1, W_3^1\} \subseteq V_1$ is a special group set linear subalgebra of $V_1 = \{V_1^1, V_2^1, V_3^1\}$. Now $W_1^2 = Z \times \{0\} \times \{0\}$ is a subgroup of $V_1^2$, $W_2^2 = \{4 \times 3$ matrices with entries from $7Z\}$ is a subgroup of $V_2^2$, $W_3^2 = \{(a\ a\ a\ a\ a\ a) \mid a \in 19Z\}$ is a subgroup of the group $V_3^2$, $W_4^2 = \{$set of all $4 \times 3$ matrices with entries from $11Z\} \subseteq V_4^2$ is a subgroup of the group $V_4^2$. Thus $W_2 = \{W_1^2, W_2^2, W_3^2, W_4^2\}$ is a special group set linear subalgebra of the special group set linear algebra $V_2$. $W_1^3 = 11Z \times 11Z$ is a subgroup of the group $Z \times Z$, $W_2^3 = \{4 \times 4$ upper triangular matrices with entries from $2Z\}$ is a subgroup of the group $V_2^3$. Thus $W_3 = \{W_1^3, W_2^3\} \subseteq V_3$ is again a special group linear subalgebra of $V_3$ over $Z$.

Consider
$$W_1^4 = \{7Z \times 2Z \times 3Z \times 5Z \times 7Z\} \subseteq V_1^4,$$
$W_1^4$ is a subgroup of the group $V_1^4$.

$W_2^4 = \{3 \times 3$ lower triangular matrices with entries from $11Z\}$ is a subgroup of the group $V_2^4$.
$$W_3^4 = 6Z \times 10Z \times 33Z \subseteq V_3^4$$
is a subgroup of $V_3^4$,

$$W_4^4 = \left\{ \begin{bmatrix} a \\ a \\ a \\ a \\ a \\ a \\ a \end{bmatrix} \middle| \; a \in 12Z \right\}$$

is a subgroup of $V_4^4$ and $W_5^4 = \{5 \times 2$ matrices with entries from $6Z\}$ is a subgroup of $V_5^4$. Thus

$$W_4 = \{W_1^4, W_2^4, W_3^4, W_4^4, W_5^4\} \subseteq V_4 = \{V_1^4, V_2^4, V_3^4, V_4^4, V_5^4\}$$

is a special group set linear subalgebra of $V_4$. Thus $W = (W_1 \cup W_2 \cup W_3 \cup W_4) \subseteq (V_1 \cup V_2 \cup V_3 \cup V_4) = V$ is a special group set linear 4-subalgebra of $V$ over the group $G = Z$.

Now we proceed onto define the notion of pseudo special group set vector n-subspace of $V$.

**DEFINITION 4.2.16:** *Let*

$$V = (V_1 \cup V_2 \cup \ldots \cup V_n)$$
$$= \{V_1^1, V_2^1, \ldots, V_{n_1}^1\} \cup \{V_1^2, V_2^2, \ldots, V_{n_2}^2\} \cup \ldots \cup \{V_1^n, V_2^n, \ldots, V_{n_n}^n\}$$

*be a special group set linear n-algebra over a group G. If*

$$W = (W_1 \cup W_2 \cup \ldots \cup W_n)$$
$$= \{W_1^1, W_2^1, \ldots, W_{n_1}^1\} \cup \{W_1^2, W_2^2, \ldots, W_{n_2}^2\} \cup \ldots \cup \{W_1^n, W_2^n, \ldots, W_{n_n}^n\}$$
$$\subseteq (V_1 \cup V_2 \cup \ldots \cup V_n)$$

*such that each $W_{j_i}^i \subseteq V_{j_i}^i$ is a subgroup of $V_{j_i}^i$, $1 \leq j_i \leq n_i$, $i = 1, 2, \ldots, n$ and if $S \subset G$ be only a proper subset $G$ and if $W = W_1 \cup W_2 \cup \ldots \cup W_n \subseteq V_1 \cup V_2 \cup \ldots \cup V_n$ is a special group set vector n-space over the set S then we call $W = W_1 \cup \ldots \cup W_n$ to be a pseudo special group set vector n-subspace of V over the subset $S \subseteq G$.*

We shall illustrate this by a simple example.

*Example 4.2.19:* Let
$$V = (V_1 \cup V_2 \cup V_3 \cup V_4 \cup V_5)$$
$$= \{V_1^1, V_2^1\} \cup \{V_1^2, V_2^2, V_3^2\} \cup \{V_1^3, V_2^3\} \cup \{V_1^4, V_2^4, V_3^4\} \cup$$
$$\{V_1^5, V_2^5, V_3^5, V_4^5\}$$
be a special group set linear 5-algebra over the group Z where
$$V_1^1 = \{Z \times Z \times Z\},$$
$$V_2^1 = \{\text{all } 3 \times 6 \text{ matrices with entries from Z}\},$$
$$V_1^2 = 3Z \times 5Z \times 7Z \times Z,$$
$$V_2^2 = \{Z[x] \text{ all polynomials of degree less than or equal to 5}\},$$
$$V_3^2 = \{6 \times 2 \text{ matrices with entries from Z}\},$$
$$V_1^3 = \{11Z \times 2Z \times 3Z \times 11Z \times 13Z\},$$
$$V_2^3 = \{\text{all } 3 \times 3 \text{ upper triangular matrices with entries from Z}\},$$
$$V_1^4 = \{Z \times Z \times Z \times Z\},$$
$$V_2^4 = \{\text{all } 4 \times 4 \text{ lower triangular matrices with entries from Z}\},$$
$$V_3^4 = \{\text{all polynomials in x of degree less than or equal to 7}$$
$$\text{with coefficients from Z}\},$$
$$V_1^5 = Z \times Z \times Z \times Z \times 11Z \times 13Z,$$
$$V_2^5 = \{7 \times 7 \text{ upper triangular matrices with entries from Z}\},$$
$$V_3^5 = \{2 \times 3 \text{ matrices with entries from Z}\} \text{ and}$$
$$V_4^5 = \{5 \times 2 \text{ matrices with entries from Z}\}.$$

Take
$$W = (W_1 \cup W_2 \cup W_3 \cup W_4 \cup W_5)$$
$$= \{W_1^1, W_2^1\} \cup \{W_1^2, W_2^2, W_3^2\} \cup \{W_1^3, W_2^3\} \cup \{W_1^4, W_2^4, W_3^4\}$$
$$\cup \{W_1^5, W_2^5, W_3^5, W_4^5\}$$
$$\subseteq (V_1 \cup V_2 \cup V_3 \cup V_4 \cup V_5) = V$$
with $W_1^1 = Z \times Z \times \{0\} \subseteq V_1^1$ is a subgroup of $V_1^1$, $W_2^1 = \{\text{all } 3\times 6 \text{ matrices with entries from } 2Z\} \subseteq V_2^1$ is a subgroup of $V_2^1$, $W_1^2 = \{30Z \times \{0\} \times \{0\} \times 10Z\}$ is a subgroup of $V_1^2$, $W_2^2 = \{\text{all}$ polynomials of degree less than or equal to two with coefficients from Z} is a subgroup of $V_2^2$, $W_3^2 = \{6 \times 2$ matrices with entries from 3Z} is a subgroup of $V_3^2$, $W_1^3 = \{\{0\} \times 4Z \times \{0\} \times 22Z\}$ is a subgroup of $V_1^3$, $W_2^3 = \{$all $3 \times 3$ upper triangular matrices with entries from 5Z} is a subgroup of $V_2^3$, $W_1^4 = \{Z \times \{0\} \times 3Z \times \{0\}\}$ is a subgroup of $V_1^4$, $W_2^4 = \{$all $4 \times 4$ lower triangular matrices with entries from 7Z} is a subgroup of $V_2^4$, $W_3^4 = \{$all polynomials in x of degree less than or equal 3 with coefficients from Z} is a subgroup of $V_3^4$, $W_1^5 = Z \times \{0\} \times \{0\} \times \{0\} \times \{0\} \times 13Z$ is a subgroup of $V_1^5$, $W_2^5 = \{7 \times 7$ upper triangular matrices with entries from 9Z} is a subgroup of $V_2^5$, $W_3^5 = \{2 \times 3$ matrices with entries from 13Z} is a subgroup of $V_3^5$ and $W_4^5 = \{5 \times 2$ matrices with entries from 14Z} is a subgroup of $V_4^5$. Thus $W = (W_1 \cup W_2 \cup W_3 \cup W_4 \cup W_5) \subseteq (V_1 \cup V_2 \cup V_3 \cup V_4 \cup V_5)$ is a pseudo special group set vector 5-subspace of V over the set $S = \{0, \pm 1\} \subseteq Z$.

Now we proceed onto define yet another new substructure.

**DEFINITION 4.2.17:** *Let $V = (V_1 \cup V_2 \cup ... \cup V_n)$ be a special group set linear n-algebra over the group G. Suppose $W = (W_1 \cup W_2 \cup ... \cup W_n) \subseteq (V_1 \cup V_2 \cup ... \cup V_n)$ where each $W_i \subseteq V_i$ is a special group set linear n-subalgebra over a proper subgroup H of G then we call $W = (W_1 \cup ... \cup W_n) \subseteq (V_1 \cup V_2 \cup ... \cup V_n)$ to be a pseudo special subgroup set linear n-subalgebra over of V over the subgroup $H \subseteq G$.*

We illustrate first this definition by an example.

*Example 4.2.20:* Let $V = (V_1 \cup V_2 \cup V_3 \cup V_4 \cup V_5)$ be a special group set linear algebra over the group Z, where
$V_1 = \{V_1^1, V_2^1, V_3^1\}$ with

$$V_1^1 = Z \times Z \times Z,$$

$V_2^1 = \{$all $3 \times 3$ matrices with entries from Z} and

$$V_3^1 = \{3 \times 5 \text{ matrices with entries from } Z\}.$$
$V_2 = \{V_1^2, V_2^2\}$ where
$$V_1^2 = 3Z \times 7Z \times Z \text{ and}$$
$$V_2^2 = \{2 \times 6 \text{ matrices with entries from } Z\}.$$
$V_3 = \{V_1^3, V_2^3, V_3^3, V_4^3, V_5^3\}$ with
$$V_1^3 = \{2 \times 5 \text{ matrices with entries from } Z\},$$
$$V_2^3 = Z \times Z \times Z \times 3Z,$$
$$V_3^3 = \{6 \times 2 \text{ matrices with entries from } Z\},$$
$V_4^3 = \{$all polynomials in the variable x with coefficients from Z of degree less than or equal to 5$\}$ and

$$V_5^3 = \left\{ \begin{bmatrix} a \\ a \\ a \\ a \\ a \\ a \end{bmatrix} \middle| a \in Z \right\}.$$

$V_4 = \{V_1^4, V_2^4\}$ with
$$V_1^4 = \{(a\ a\ a\ a\ a\ a) \mid a \in Z\} \text{ and}$$
$$V_2^4 = Z \times 3Z \times 2Z \times 5Z \times Z$$
and
$V_5 = \{V_1^5, V_2^5, V_3^5, V_4^5\}$ with
$$V_1^5 = Z \times Z \times Z \times Z \times 3Z,$$
$V_2^5 = \{$all $5 \times 5$ upper triangular matrices with entries from $Z\}$,
$$V_3^5 = \{(a\ a\ a\ a\ a) \mid a \in Z\} \text{ and}$$
$V_4^5 = \{2 \times 9$ matrices with entries from $Z\}$.

Take
$$W_1^1 = 3Z \times 3Z \times 3Z \subseteq V_1^1 \text{ is a subgroup of } V_1^1,$$
$$W_2^1 = \{3 \times 3 \text{ matrices with entries from } 5Z\} \subseteq V_2^1$$

is a subgroup of $V_2^1$ and

$W_3^1 = \{3 \times 5 \text{ matrices with entries from } 7Z\} \subseteq V_3^1$

is a subgroup of $V_3^1$.

$W_1^2 = \{6Z \times \{0\} \times 5Z\} \subseteq V_1^2$

is a subgroup of $V_1^2$ and

$W_2^2 = \{2 \times 6 \text{ matrices with entries from } 7Z\} \subseteq V_2^2$,

$W_1^3 = \{2 \times 5 \text{ matrices with entries from } 8Z\} \subseteq V_1^3$

is a subgroup of $V_1^3$,

$W_2^3 = 3Z \times \{0\} \times \{0\} \times 3Z \subseteq V_2^3$

is a subgroup of $V_2^3$,

$W_3^3 = \{6 \times 2 \text{ matrices with entries from } 10Z\} \subseteq V_3^3$

is a subgroup of $V_3^3$,

$W_4^3 = \{$all polynomials in the variable x with coefficients from Z of degree less than or equal to $3\} \subseteq V_4^3$

is a subgroup of $V_4^3$,

$$W_5^3 = \left\{ \begin{bmatrix} a \\ a \\ a \\ a \\ a \\ a \end{bmatrix} \,\middle|\, a \in 10Z \right\} \subseteq V_5^3$$

is a subgroup of $V_5^3$,

$W_1^4 = \{(a\ a\ a\ a\ a\ a) \mid a \in 2Z\} \subseteq V_1^4$

is a subgroup of $V_1^4$ and

$W_2^4 = \{3Z \times \{0\} \times \{0\} \times \{0\} \times 5Z\}$

is a subgroup of $V_2^4$.

$W_1^5 = \{0\} \times 3Z \times \{0\} \times \{0\} \times 3Z$

is a subgroup of $V_1^5$,

$W_2^5 = \{$all $5 \times 5$ upper triangular matrices with entries from $7Z\}$

is a subgroup of $V_2^5$,
$$W_3^5 = \{(a\ a\ a\ a\ a) / a \in 2Z\} \subseteq V_3^5$$
is a subgroup of $V_3^5$ and
$$W_4^5 = \{2 \times 9 \text{ matrices with entries from } 10Z\} \subseteq V_4^5$$
is a subgroup of $V_4^5$.

Thus
$$W = V = (W_1 \cup W_2 \cup W_3 \cup W_4 \cup W_5)$$
$$= \{W_1^1, W_2^1, W_3^1\} \cup \{W_1^2, W_2^2\} \cup \{W_1^3, W_2^3, W_3^3, W_4^3, W_5^3\}$$
$$\cup \{W_1^4, W_2^4\} \cup \{W_1^5, W_2^5, W_3^5, W_4^5\}$$
$$\subseteq V = (V_1 \cup V_2 \cup V_3 \cup V_4)$$

is a pseudo special subgroup set linear 5-subalgebra over the subgroup $2Z \subseteq Z$.

It is important to observe the following:
   If the group over which the special group set linear n algebra V defined does not contain any proper subgroups then we define those special group set linear n-algebras as pseudo simple special group set linear n-algebras.

We in this context prove the following result.

**THEOREM 4.2.1:** *Let*
$$V = (V_1 \cup V_2 \cup \ldots \cup V_n)$$
$$= \{V_1^1, V_2^1, \ldots, V_{n_1}^1\} \cup \{V_1^2, V_2^2, \ldots, V_{n_2}^2\} \cup \ldots \cup \{V_1^n, V_2^n, \ldots, V_{n_n}^n\}$$
*be a special group set linear n-algebra over the group $Z_p$. (The group of integers under addition modulo p, p a prime}. Then V is a pseudo simple special group set linear n-algebra over $Z_p$.*

*Proof:* Given $V = V = (V_1 \cup V_2 \cup \ldots \cup V_n)$ is a special group set linear n-algebra defined over the additive group $Z_p$. Thus $Z_p$ being simple $Z_p$ cannot have non trivial subgroup so even if W is proper subset of V such that if
$$W = (W_1 \cup W_2 \cup \ldots \cup W_n)$$

$$= \{W_1^1, W_2^1, ..., W_{n_1}^1\} \cup \{W_1^2, W_2^2, ..., W_{n_2}^2\} \cup ... \cup$$
$$\{W_1^n, W_2^n, ..., W_{n_n}^n\}$$

such that each $W_{j_i}^i$, $1 \leq j_i \leq n_i$; $i = 1, 2, ..., n$ happens to be a subgroup of $V_{j_i}^i$ still the absence of subgroups in $Z_p$ makes V a pseudo simple special group set linear n-algebra over $Z_p$. Thus from this theorem we have an infinite class of pseudo simple special group set linear n-algebras as the number of primes is infinite.

We shall illustrate this by a simple example.

***Example 4.2.21:*** Let
$$V = (V_1 \cup V_2 \cup V_3 \cup V_4)$$
$$= \{V_1^1, V_2^1\} \cup \{V_1^2, V_2^2, V_3^2\} \cup \{V_1^3, V_2^3\} \cup \{V_1^4, V_2^4\}$$
be a special group set linear 4-algebra over the group $Z_7$. Here $V_1^1 = Z_7 \times Z_7$,

$$V_2^1 = \left\{ \begin{bmatrix} a \\ b \\ c \end{bmatrix} \middle| a, b, c \in Z_7 \right\}$$

$V_1^2 = Z_7 \times Z_7 \times Z_7$, $V_2^2 = \{2 \times 2$ matrices with entries from $Z_7\}$, $V_3^2 = \{2 \times 5$ matrices with entries from $Z_7$, $V_1^3 = \{Z_7 \times Z_7 \times Z_7 \times Z_7\}$, $V_2^3 = \{4 \times 4$ upper triangular matrices with entries from $Z_7\}$, $V_1^4 = \{5 \times 5$ lower triangular matrices with entries from $Z_7\}$ and $V_2^4 = Z_7 \times Z_7 \times Z_7 \times Z_7 \times Z_7$.

Take
$$W = (W_1 \cup W_2 \cup W_3 \cup W_4)$$
$$= \{W_1^1, W_2^1\} \cup \{W_1^2, W_2^2, W_3^2\} \cup \{W_1^3, W_2^3\} \cup \{W_1^4, W_2^4\}$$
where $W_1^1 = Z_7 \times \{0\} \subseteq V_1^1$ is a subgroup of the group $V_1^1$,

$$W_2^1 = \left\{ \begin{bmatrix} a \\ a \\ a \end{bmatrix} \middle| a \in Z_7 \right\} \subseteq V_2^1$$

is a subgroup of the group $V_2^1$, $W_1^2 = Z_7 \times \{0\} \times Z_7 \subseteq V_1^2$ is a subgroup of the group $V_1^2$,

$$W_2^2 = \left\{ \begin{pmatrix} a & a \\ a & a \end{pmatrix} \middle| a \in Z_7 \right\} \subseteq V_2^2$$

is a subgroup of the group $V_2^2$,

$$W_3^2 = \left\{ \begin{pmatrix} a & a & a & a & a \\ a & a & a & a & a \end{pmatrix} \middle| a \in Z_7 \right\} \subseteq V_3^2$$

is a subgroup of the group $V_3^2$. $W_1^3 = \{Z_7 \times \{0\} \times Z_7 \times \{0\}\} \subseteq V_1^3$ is a subgroup of the group $V_1^3$,

$$W_2^3 = \left\{ \begin{pmatrix} a & a & a & a \\ 0 & a & a & a \\ 0 & 0 & a & a \\ 0 & 0 & 0 & a \end{pmatrix} \middle| a \in Z_7 \right\} \subseteq V_2^3$$

is a subgroup of the group $V_2^3$

$$W_1^4 = \left\{ \begin{pmatrix} a & 0 & 0 & 0 & 0 \\ a & a & 0 & 0 & 0 \\ a & a & a & 0 & 0 \\ a & a & a & a & 0 \\ a & a & a & a & a \end{pmatrix} \middle| a \in Z_7 \right\} \subseteq V_1^4$$

is a subgroup of the group $V_1^4$ and $W_2^4 = Z_7 \times \{0\} \times Z_7 \times \{0\} \times \{0\} \subseteq V_2^4$ is a subgroup of the group.

Though $W = (W_1 \cup W_2 \cup W_3 \cup W_4) \subseteq (V_1 \cup V_2 \cup V_3 \cup V_4) = V$ still as $Z_7$ has no subgroup $W$ cannot be a pseudo special subgroup set linear 4-subalgebra of $V$ so $V$ is itself a pseudo simple special group set linear 4-algebra.

Now we proceed onto define the notion of special group set linear n-transformation.

**DEFINITION 4.2.18:** *Let*
$$V = (V_1 \cup V_2 \cup \ldots \cup V_n)$$
$$= \{V_1^1, V_2^1, \ldots, V_{n_1}^1\} \cup \{V_1^2, V_2^2, \ldots, V_{n_2}^2\} \cup \ldots \cup \{V_1^n, V_2^n, \ldots, V_{n_n}^n\}$$
*and*
$$W = (W_1 \cup W_2 \cup \ldots \cup W_n) =$$
$$\{W_1^1, W_2^1, \ldots, W_{n_1}^1\} \cup \{W_1^2, W_2^2, \ldots, W_{n_2}^2\} \cup \ldots \cup \{W_1^n, W_2^n, \ldots, W_{n_n}^n\}$$
*be two special group set vector n-spaces over the set S.*
$$T: (T_1 \cup T_2 \cup \ldots \cup T_n)$$
$$= \{T_1^1, T_2^1, \ldots, T_{n_1}^1\} \cup \{T_1^2, T_2^2, \ldots, T_{n_2}^2\} \cup \ldots \cup \{T_1^n, T_2^n, \ldots, T_{n_n}^n\} :$$
$$V = (V_1 \cup V_2 \cup \ldots \cup V_n) \to W = (W_1 \cup W_2 \cup \ldots \cup W_n)$$
*where $T_i : V_i \to W_i$ i.e.,*
$$T_i : \left(V_1^i, V_2^i, \ldots, V_{n_i}^i\right) \to \left(W_1^i, W_2^i, \ldots, W_{n_i}^i\right)$$
*such that $T_i$ is a special group set linear transformation of $V_i$ into $W_i$ over the group S for every i, i = 1, 2, …, n.*

*That is*
$$T_{ji}^i : V_{ji}^i \to W_{ji}^i ;$$
*$1 \leq j_i \leq n_i$, i = 1, 2, …, n where $T_i = \left(T_1^i, T_2^i, \ldots, T_{n_i}^i\right)$; then we call $T = (T_1 \cup T_2 \cup \ldots \cup T_n)$ to be a special group set linear n-transformation of V into W.*

*We shall denote the collection of all such special group set linear n-transformation V into W by*
*SHom(V, W) = {SHom (V_1, W_1)} $\cup$ {SHom (V_2, W_2)} $\cup$ … $\cup$ SHom(V_n, W_n)}*

$$= \{Hom\left(V_1^1, W_1^1\right), Hom\left(V_2^1, W_2^1\right), ..., Hom\left(V_{n_1}^1, W_{n_1}^1\right)\}$$
$$\cup \{Hom\left(V_2^2, W_2^2\right), ..., Hom\left(V_{n_2}^2, W_{n_2}^2\right)\} \cup ... \cup$$
$$\{Hom\left(V_1^n, W_1^n\right), Hom\left(V_2^n, W_2^n\right), ..., Hom\left(V_{n_n}^n, W_{n_n}^n\right)\}.$$

*Clearly S Hom (V, W) need not in general be a special group set vector n-space over S. Infact SHom (V, W) can in general be a special set vector n space over S.*

***Example 4.2.22:*** Let $V = (V_1 \cup V_2 \cup V_3 \cup V_4)$ and $W = (W_1 \cup W_2 \cup W_3 \cup W_4)$ be two special group set vector spaces over the group $S = \{0, 1\}$. Here
$$V_1 = \{V_1^1, V_2^1\},\ V_2 = \{V_1^2, V_2^2, V_3^2\},$$
$$V_3 = \{V_1^3, V_2^3, V_3^3\},\ V_4 = \{V_1^4, V_2^4\}$$

with
$$V_1^1 = Z \times Z \times Z \times Z,$$
$V_2^1 = \{$all $2 \times 3$ matrices with entries from $Z\}$, $V_1^2 = Z \times Z \times Z$,
$V_2^2 = \{$all polynomials in the variable x with coefficients from Z of degree less than or equal to 4$\}$,
$V_3^2 = \{3 \times 3$ upper triangular matrices with entries from Z$\}$,
$V_1^3 = \{3 \times 3$ matrices with entries from Z$\}$,
$V_2^3 = Z \times Z$, $V_3^3 = \{9 \times 1$ column vector with entries from Z$\}$,
$$V_1^4 = Z \times Z \times Z \times Z \times Z \times Z,$$
$V_2^4 = \{3 \times 2$ matrices with entries from Z$\}$,

$$W_1^1 = \left\{ \begin{pmatrix} a & b \\ c & d \end{pmatrix} \middle| a,b,c,d \in Z \right\},$$

$$W_2^1 = \left\{ \begin{bmatrix} a_1 & a_2 \\ a_3 & a_4 \\ a_5 & a_6 \end{bmatrix} \middle| a_i \in Z;\ 1 \le i \le 6 \right\},$$

$$W_1^2 = \left\{ \begin{pmatrix} a & b \\ 0 & c \end{pmatrix} \middle| a,b,c \in Z \right\},$$

$$W_2^2 = \left\{ \begin{bmatrix} a_1 \\ a_2 \\ a_3 \\ a_4 \\ a_5 \end{bmatrix} \middle| a_i \in Z; \ 1 \le i \le 5 \right\},$$

$W_3^2 = \{\{3 \times 2$ matrices with entries from $Z\}$,
$W_1^3 = \{$all polynomials in x degree less than or equal to 8 with coefficients from $Z\}$,

$$W_2^3 = \left\{ \begin{pmatrix} a \\ b \end{pmatrix} \middle| a,b \in Z \right\},$$

$W_1^4 = \{3 \times 3$ lower triangular matrices with entries from $Z\}$
and

$$W_2^4 = \left\{ \begin{pmatrix} a & b & c \\ d & e & f \end{pmatrix} \middle| a,b,c,d,e,f \in Z \right\}.$$

Now define the map $T = T_1 \cup T_2 \cup T_3 \cup T_4 : V = V_1 \cup V_2 \cup V_3 \cup V_4 \to W = W_1 \cup W_2 \cup W_3 \cup W_4$ as follows;

$$T = (T_1 \cup T_2 \cup T_3 \cup T_4)$$
$$= \{T_1^1, T_2^1\} \cup \{T_1^2, T_2^2, T_3^2\} \cup \{T_1^3, T_2^3, T_3^3\} \cup \{T_1^4, T_2^4\} :$$
$$V \to W \text{ by}$$

$T_1^1 : V_1^1 \to W_1^1$ given by

$$T_1^1 \ (a \ b \ c \ d) = \begin{pmatrix} a & b \\ c & d \end{pmatrix},$$

$T_2^1 : V_2^1 \to W_2^1$ is defined by

$$T_2^1 \begin{pmatrix} a & b & c \\ d & e & f \end{pmatrix} = \begin{pmatrix} a & d \\ b & e \\ c & f \end{pmatrix},$$

$T_1^2 : V_1^2 \to W_1^2$ is given by

$$T_1^2 \, (a\ b\ c) = \begin{pmatrix} a & b \\ 0 & c \end{pmatrix},$$

$T_2^2 : V_2^2 \to W_2^2$ is such that

$$T_2^2 \, (a_0 + a_1 x + a_2 x^2 + a_3 x^3 + a_4 x^4) = \begin{pmatrix} a_0 \\ a_1 \\ a_2 \\ a_3 \\ a_4 \end{pmatrix}$$

and
$T_3^2 : V_3^2 \to W_3^2$ is defined by

$$T_3^2 \begin{pmatrix} a & b & c \\ 0 & d & e \\ 0 & 0 & f \end{pmatrix} = \begin{bmatrix} a & d \\ b & e \\ c & f \end{bmatrix}.$$

$T_1^3 : V_1^3 \to W_1^3$ is given by

$$T_1^3 \begin{pmatrix} a & b & c \\ d & e & f \\ g & h & i \end{pmatrix}$$
$$= \{a + bx + cx^2 + dx^3 + ex^4 + fx^5 + gx^6 + hx^7 + ix^8\},$$

$T_2^3 : V_2^3 \to W_2^3$ is defined by

$$T_2^3 \, (a, b) = \begin{bmatrix} a \\ b \end{bmatrix},$$

$T_3^3 : V_3^3 \to W_3^3$ is such that

$$T_3^3 \begin{bmatrix} a_1 \\ a_2 \\ a_3 \\ a_4 \\ \vdots \\ a_9 \end{bmatrix} = a_1 + a_2x + a_3x^2 + a_4x^3 + a_5x^4 + a_6x^5 + a_7x^6 + a_8x^7 + a_9x^8,$$

$T_1^4 : V_1^4 \to W_1^4$ is given by

$$T_1^4 \,(a\ b\ c\ d\ e\ f) = \begin{pmatrix} a & 0 & 0 \\ b & d & 0 \\ c & e & f \end{pmatrix}$$

and
$T_2^4 : V_2^4 \to W_2^4$ is such that

$$T_2^4 \begin{pmatrix} a_1 & a_4 \\ a_2 & a_5 \\ a_3 & a_6 \end{pmatrix} = \begin{pmatrix} a_1 & a_2 & a_3 \\ a_4 & a_5 & a_6 \end{pmatrix}.$$

Thus $T = T_1 \cup T_2 \cup T_3 \cup T_4$ is a special group set linear 4 transformation of V into V.

It is important at this stage to make a note that if W in the definition is replaced by V itself i.e., The domain and the range space are one and the same then we call the map T: V $\to$ V which is a special group set linear transformation from V into V is defined as the special group set linear n-operator on V.

Now if
SHom (V,V) = {S Hom (V$_1$, V$_1$) $\cup$ S Hom (V$_2$, V$_2$) $\cup$ … $\cup$ S Hom (V$_n$, V$_n$)}
= (Hom $\left(V_1^1, V_1^1\right)$, Hom $\left(V_2^1, V_2^1\right)$, …, Hom $\left(V_{n_1}^1, V_{n_1}^1\right)) \cup$
(Hom $\left(V_1^2, V_1^2\right)$, Hom $\left(V_2^2, V_2^2\right)$, …, Hom $\left(V_{n_2}^2, V_{n_2}^2\right)) \cup … \cup$

$$(\text{Hom}\left(V_1^n, V_1^n\right), \ldots, \text{Hom}\left(V_{n_1}^n, V_{n_1}^n\right))$$

then we see SHom (V, V) is a special group set vector n-space over the sets over which V is defined.

We illustrate this by the following example.

*Example 4.2.23:* Let
$$V = (V_1 \cup V_2 \cup V_3 \cup V_4)$$
$$= \{V_1^1, V_2^1, V_3^1\} \cup \{V_1^2, V_2^2\} \cup \{V_1^3, V_2^3, V_3^3, V_4^3\} \cup \{V_1^4, V_2^4, V_3^4\}$$
where
$$V_1^1 = Z \times Z \times Z \times Z,$$
$$V_2^1 = \left\{ \begin{bmatrix} a & b \\ c & d \end{bmatrix} \middle| a, b, c, d \in Z \right\},$$

$V_3^1 = \{$set of all polynomials in the variable x of degree less than or equal to 6$\}$, $V_1^2 = 3Z \times Z \times 8Z$,

$V_2^2 = \{3 \times 3$ upper triangular matrices with entries from $Z\}$,
$$V_1^3 = Z \times 2Z \times 3Z \times 4Z \times 5Z,$$
$V_2^3 = \{$all $4 \times 4$ upper triangular matrices with entries from $Z\}$,

$V_3^3 = \{$all polynomials in the variable x with coefficients from Z of degree less than or equal to 5$\}$,

$$V_4^3 = \left\{ \begin{pmatrix} a & a & a & a \\ b & b & b & b \\ c & c & c & c \end{pmatrix} \middle| a, b, c \in Z \right\},$$

$$V_1^4 = \left\{ \begin{pmatrix} a \\ b \\ c \\ d \end{pmatrix} \middle| a, b, c, d \in Z \right\},$$

$$V_2^4 = \{Z \times Z \times Z \times Z \times Z \times Z\}$$

and $V_3^4 = \{$all $6 \times 6$ lower triangular matrices with entries from Z$\}$, be the special group set 4-vector space over Z.

Now define
$$T = (T_1 \cup T_2 \cup T_3)$$
$$= \left(T_1^1, T_2^1, T_3^1\right) \cup \left(T_1^2, T_2^2\right) \cup \left(T_1^3, T_2^3, T_3^3, T_4^3\right) \cup \left(T_1^4, T_2^4, T_3^4\right):$$
$$V = (V_1 \cup V_2 \cup V_3 \cup V_4)$$
$$= \left(V_1^1, V_2^1, V_3^1\right) \cup \left(V_1^2, V_2^2\right) \cup \left(V_1^3, V_2^3, V_3^3, V_4^3\right) \cup \left(V_1^4, V_2^4, V_3^4\right)$$
$$\to V = (V_1 \cup V_2 \cup V_3 \cup V_4)$$

as follows:
$T_1^1 : V_1^1 \to V_1^1$ is defined by
$$T_1^1 \,(a\ b\ c\ d) = (a\ a\ d\ d);$$

$T_2^1 : V_2^1 \to V_2^1$ is defined by

$$T_2^1 \begin{pmatrix} a & b \\ c & d \end{pmatrix} = \begin{pmatrix} d & c \\ b & a \end{pmatrix},$$

$T_3^1 : V_3^1 \to V_3^1$ is given by

$$T_3^1 \,(a_o + a_1 x + a_2 x^2 + a_3 x^3 + a_4 x^4 + a_5 x^5 + a_6 x^6)$$
$$= a_o x^6 + a_1 x^5 + a_2 x^4 + a_3 x^3 + a_4 x^2 + a_5 x + a_6),$$

$T_1^2 : V_1^2 \to V_1^2$ is defined by
$$T_1^2 \,(x\ y\ z) = (x\ \ x + y + z\ \ z),$$

$T_2^2 : V_2^2 \to V_2^2$ is given by

$$T_2^2 \begin{pmatrix} a & b & c \\ 0 & d & e \\ 0 & 0 & f \end{pmatrix} = \begin{pmatrix} a & 0 & 0 \\ 0 & d & 0 \\ 0 & 0 & f \end{pmatrix},$$

$T_1^3: V_1^3 \to V_1^3$ is such that
$$T_1^3 (x\ y\ z\ \omega\ t) = (x + y + z,\ \omega,\ 3y,\ 2y,\ z),$$

$T_2^3: V_2^3 \to V_2^3$ is given by

$$T_2^3 \begin{pmatrix} a & b & c & d \\ 0 & e & f & g \\ 0 & 0 & h & i \\ 0 & 0 & 0 & j \end{pmatrix} = \begin{pmatrix} a & b & c & d \\ 0 & e & 0 & 0 \\ 0 & 0 & h & 0 \\ 0 & 0 & 0 & j \end{pmatrix},$$

$T_3^3: V_3^3 \to V_3^3$ is defined as

$$T_3^3 (a_o + a_1 x + a_2 x^2 + a_3 x^3 + a_4 x^4 + a_5 x^5) = a_0 + a_3 x^3 + a_5 x^5,$$

$T_4^3: V_4^3 \to V_4^3$ is given by

$$T_4^3 \begin{pmatrix} a & a & a & a \\ b & b & b & b \\ c & c & c & c \end{pmatrix} = \begin{pmatrix} c & c & c & c \\ a & a & a & a \\ b & b & b & b \end{pmatrix},$$

$T_1^4: V_1^4 \to V_1^4$ is given by

$$T_1^4 \begin{pmatrix} a \\ b \\ c \\ d \end{pmatrix} = \begin{pmatrix} a \\ d \\ b \\ c \end{pmatrix},$$

$T_2^4: V_2^4 \to V_2^4$ is such that
$$T_2^4 (a\ b\ c\ d\ e\ f) = (a\ a\ a\ f\ f\ f)$$
and
$T_3^4: V_3^4 \to V_3^4$ is given by

$$T_3^4 \begin{pmatrix} a & 0 & 0 & 0 & 0 & 0 \\ b & c & 0 & 0 & 0 & 0 \\ d & e & f & 0 & 0 & 0 \\ g & h & i & j & 0 & 0 \\ k & l & m & n & p & 0 \\ q & r & s & t & u & v \end{pmatrix} = \begin{pmatrix} a & 0 & 0 & 0 & 0 & 0 \\ b & b & 0 & 0 & 0 & 0 \\ d & d & d & 0 & 0 & 0 \\ g & g & g & g & 0 & 0 \\ k & k & k & k & k & 0 \\ q & q & q & q & q & q \end{pmatrix}.$$

Thus we see $T = (T_1 \cup T_2 \cup T_3 \cup T_4) : V \to V$ is a special group set linear 4 operator on V.

Now we can proceed onto define the notion of pseudo special group set linear n-operator on V.

**DEFINITION 4.2.19:** *Let*
$$V = (V_1 \cup V_2 \cup \ldots \cup V_n)$$
$$= \left(V_1^1, V_2^1, \ldots, V_{n_1}^1\right) \cup \left(V_1^2, V_2^2, \ldots, V_{n_2}^2\right) \cup \ldots \cup \left(V_1^n, V_2^n, \ldots, V_{n_n}^n\right)$$
*be a special group set vector n-space over the set S. Let*
$$T = (T_1 \cup T_2 \cup \ldots \cup T_n)$$
$$= \left(T_1^1, T_2^1, \ldots, T_{n_1}^1\right) \cup \ldots \cup \left(T_1^n, T_2^n, \ldots, T_{n_n}^n\right).$$
$$V = (V_1 \cup V_2 \cup \ldots \cup V_n)$$
$$= \left(V_1^1, V_2^1, \ldots, V_{n_1}^1\right) \cup \left(V_1^2, V_2^2, \ldots, V_{n_2}^2\right) \cup \ldots \cup \left(V_1^n, V_2^n, \ldots, V_{n_n}^n\right)$$
$$\to V = (V_1 \cup V_2 \cup \ldots \cup V_n)$$

*is such that $T_i : V_i \to V_i$ is given by $T_{j_i}^i : V_{j_i}^i \to V_k^i$; $j_i \neq k$ for atleast one $j_i$; $1 \leq j_i, k \leq n_i, i = 1, 2, \ldots, n$ be a pseudo special group set linear operator on $V_i$, true for each i. Then we define $T = (T_1 \cup T_2 \cup \ldots \cup T_n)$ to be the pseudo special group set linear n-operator on V.*

We illustrate this by a simple example

***Example 4.2.24:*** Let $V = (V_1 \cup V_2 \cup V_3 \cup V_4)$ be a special group set vector 4-space over the set $Z^+ \cup \{0\} = S$ where

$$V_1 = \left(V_1^1, V_2^1, V_3^1, V_4^1\right), V_2 = \left(V_1^2, V_2^2, V_3^2\right),$$
$$V_3 = \left(V_1^3, V_2^3\right) \text{ and } V_4 = \left(V_1^4, V_2^4, V_3^4, V_4^4, V_5^4\right)$$

with
$$V_1^1 = S \times S \times S \times S,$$

$$V_2^1 = \left\{ \begin{pmatrix} a & b \\ c & a \end{pmatrix} \middle| a,b,c,d \in S \right\},$$

$$V_3^1 = \left\{ \begin{pmatrix} a_1 \\ a_2 \\ a_3 \\ a_4 \end{pmatrix} \middle| a_i \in S; \ 1 \le i \le 4 \right\},$$

$V_4^1 = \{$all polynomials in the variable x with coefficients from S of degree less than or equal to 5$\}$, $V_1^2 = S \times S \times S \times S \times S \times S$, $V_2^2 = \{$all $3 \times 3$ upper triangular matrices with entries from S$\}$, $V_3^2 = \{$all $6 \times 2$ matrices with entries from S$\}$, $V_1^3 = S \times S \times S \times S \times S \times S$, $V_2^3 = \{$all $4 \times 4$ lower triangular matrices with entries from S$\}$, $V_1^4 = \{S \times S \times S \times S \times S \times S \times S\}$, $V_2^4 = \{5 \times 2$ matrices with entries from $Z^+ \cup \{0\}\}$ $V_3^4 = \{2 \times 3$ matrices with entries from $Z^+ \cup \{0\}\}$, $V_4^4 = \{7 \times 1$ column matrix with entries from S$\}$ and $V_5^4 = \{$set of all polynomials in the variable x with coefficients from S of degree less than or equal to 9$\}$. V is clearly a special group set vector 4-space over the set S.

Let
$$T = (T_1 \cup T_2 \cup T_3 \cup T_4)$$
$$= \left(T_1^1, T_2^1, T_3^1, T_4^1\right) \cup \left(T_1^2, T_2^2, T_3^2\right)$$
$$\cup \left(T_1^3, T_2^3\right) \cup \left(T_1^4, T_2^4, T_3^4, T_4^4, T_5^4\right)$$
$$= V = (V_1 \cup V_2 \cup V_3 \cup V_4)$$

$$= \left(V_1^1, V_2^1, V_3^1, V_4^1\right) \cup \left(V_1^2, V_2^2, V_3^2\right) \cup \left(V_1^3, V_2^3\right) \cup$$
$$\left(V_1^4, V_2^4, V_3^4, V_4^4, V_5^4\right) \to V = (V_1 \cup V_2 \cup V_3 \cup V_4)$$

be defined as follows.
$T_i : V_i \to V_i$, $1 \le i \le 4$; $T_1 : V_1 \to V_1$, i.e.,

$$T_1 = \left(T_1^1, T_2^1, T_3^1, T_4^1\right) : \left(V_1^1, V_2^1, V_3^1, V_4^1\right) \to \left(V_1^1, V_2^1, V_3^1, V_4^1\right)$$

such that

$$T_1^1 : V_1^1 \to V_2^1,$$
$$T_2^1 : V_2^1 \to V_3^1,$$
$$T_3^1 : V_3^1 \to V_4^1 \text{ and}$$
$$T_4^1 : V_4^1 \to V_1^1$$

are given by;

$$T_1^1 \ (a\ b\ c\ d) \ = \ \begin{pmatrix} a & b \\ c & d \end{pmatrix},$$

$$T_2^1 : \begin{pmatrix} a & b \\ c & d \end{pmatrix} = \begin{bmatrix} a \\ b \\ c \\ d \end{bmatrix},$$

$$T_3^1 : \begin{bmatrix} a \\ b \\ c \\ d \end{bmatrix} = a + bx + cx^3 + dx^5$$

and
$$T_4^1 : (a_o + a_1x + a_2x^2 + a_3x^3 + a_4\ x^4 + a_5\ x^5)$$
$$= (a_o, a_1 + a_2, a_3 + a_4, a_5).$$

Now $T_2 \left(T_1^2, T_2^2, T_3^2\right) : V_2 = \left(V_1^2, V_2^2, V_3^2\right) \to V_2 = \left(V_1^2, V_2^2, V_3^2\right)$
is defined to be from

$$T_1^2 : V_1^2 \to V_2^2,$$
$$T_2^2 : V_2^2 \to V_3^2 \text{ and}$$
$$T_3^1 : V_3^1 \to V_4^1 \text{ ; are given in the following:}$$
$$T_3^2 : V_3^2 \to V_1^2$$

$$T_1^2 (a\, b\, c\, d\, e\, f) = \begin{pmatrix} a & b & c \\ 0 & d & e \\ 0 & 0 & f \end{pmatrix},$$

$$T_2^2 \begin{pmatrix} a & b & c \\ 0 & d & e \\ 0 & 0 & f \end{pmatrix} = \begin{pmatrix} a & 0 & c & 0 & d & 0 \\ 0 & e & 0 & f & 0 & e \end{pmatrix}^T,$$

and

$$T_3^2 \begin{pmatrix} a & b \\ c & d \\ e & f \\ g & h \\ i & j \\ k & l \end{pmatrix} = (a+b,\, c+d,\, e+f,\, g+h,\, i+j,\, k+l)$$

$$T_3 = \left( T_1^3, T_2^3 \right) : V_3 = \left( V_1^3, V_2^3 \right) \to \left( V_1^3, V_2^3 \right)$$

are such that

$$T_1^3 : V_1^3 \to V_2^3$$

and $T_2^3 : V_2^3 \to V_1^3$ are defined by

$$T_1^3 (a\, b\, c\, d\, e\, f) = \begin{pmatrix} a & 0 & 0 & 0 \\ b & e & 0 & 0 \\ c & f & 0 & 0 \\ d & a & c & d \end{pmatrix}$$

and

$$T_2^3 \begin{pmatrix} a & 0 & 0 & 0 \\ b & c & 0 & 0 \\ d & e & f & 0 \\ g & h & i & j \end{pmatrix} = (a,\ b+c,\ d,\ e,\ f+g,\ h+i).$$

Now

$$T_4 = \left(T_1^4, T_2^4, T_3^4, T_4^4, T_5^4\right):$$

$$V_4 = \left(V_1^4, V_2^4, V_3^4, V_4^4, V_5^4\right) \to V_4 = \left(V_1^4, V_2^4, V_3^4, V_4^4, V_5^4\right)$$

is defined by;

$T_1^4 : V_1^4 \to V_4^4$ is such that

$$T_1^4\ (a\ b\ c\ d\ e\ f\ g) = \begin{bmatrix} a \\ b \\ c \\ d \\ e \\ f \\ g \end{bmatrix};$$

$T_2^4 : V_2^4 \to V_5^4$ is defined by

$$T_2^4 \begin{pmatrix} a & b \\ c & d \\ e & f \\ g & h \\ i & j \end{pmatrix} = a + bx + cx^2 + dx^3 + ex^4 + fx^5 + gx^6 + hx^7 + ix^8 + jx^9$$

$T_3^4 : V_3^4 \to V_1^4$ is such that

$$T_3^4 \begin{pmatrix} a & b & c \\ d & e & f \end{pmatrix} = (a,\ b,\ c,\ d,\ e,\ f,\ a+f);$$

$T_3^4 : V_3^4 \to V_5^4$ is given by

$$T_4^4 \begin{bmatrix} a_0 \\ a_1 \\ a_2 \\ a_3 \\ a_4 \\ a_5 \\ a_6 \end{bmatrix} = a_o + a_1x + a_2x^2 + a_3x^3 + a_4x^4 + a_5x^5 + a_6x^6$$

and $T_5^4 : V_5^4 \to V_2^4$ is such that

$$T_5^4 = a_o + a_1x + a_2x^2 + a_3x^3 + a_4x^4$$

$$+ a_5x^5 + a_6x^6 + a_7x^7 + a_8x^8 + a_9x^9 = \begin{bmatrix} a_0 & a_1 \\ a_2 & a_3 \\ a_4 & a_5 \\ a_6 & a_7 \\ a_8 & a_9 \end{bmatrix}.$$

Clearly
$$T = (T_1 \cup T_2 \cup T_3 \cup T_4)$$
$$= \left(T_1^1, T_2^1, T_3^1, T_4^1\right) \cup \left(T_1^2, T_2^2, T_3^2\right) \cup \left(T_1^3, T_2^3\right) \cup$$
$$\left(T_1^4, T_2^4, T_3^4, T_4^4, T_5^4\right):$$
$$V = (V_1 \cup V_2 \cup V_3 \cup V_4) \to V$$
is a pseudo special group set linear 4-operator on V.

The following observations are important: If
$$V = (V_1 \cup V_2 \cup ... \cup V_n)$$
$$= \left(V_1^1, V_2^1,..., V_{n_1}^1\right) \cup \left(V_1^2, V_2^2,..., V_{n_2}^2\right) \cup ... \cup \left(V_1^n, V_2^n,..., V_{n_n}^n\right)$$
is a special group set vector n-space over the set S. The collection of all pseudo operators from V into V denoted by SHom(V, V) = SHom($V_1$, $V_1$) $\cup$ SHom($V_2$, $V_2$) $\cup$ ... $\cup$ SHom($V_n$, $V_n$) = {Hom$\left(V_1^1, V_{i_1}^1\right)$ Hom$\left(V_1^1, V_{i_1}^1\right)$ ... Hom$\left(V_{n_1}^1, V_{i_{n_1}}^1\right)$}
$\cup$ {Hom$\left(V_1^2, V_{i_2}^2\right)$ Hom$\left(V_2^2, V_{j_2}^2\right)$ ... Hom$\left(V_{n_2}^2, V_{i_{n_2}}^2\right)$} $\cup$ ... $\cup$

{Hom$\left(V_1^n, V_{p_2}^n\right)$ Hom$\left(V_2^n, V_{p_2}^n\right)$ ... Hom$\left(V_{n_n}^n, V_{p_{n_n}}^n\right)$} is again a special group set vector n-space over the set S. Here $(i_1,...,i_{n_1})$ is a permutation of (1, 2, ..., n$_1$), $(j_1,...,j_{n_2})$ is a permutation of (1,2, ..., n$_2$) so on $(p_1,...,p_{n_n})$ is a permutation of (1, 2, ..., n$_n$).

As in case of special group set vector n-spaces we can in case of special group set linear n-algebras define special group set linear n-transformations, special group set linear n-operators and so on. Further we can define special projections of these n-spaces.

We now define the fuzzy analogue of these notions.

**DEFINITION 4.2.20:** *Let V = (S$_1$, S$_2$, ... , S$_n$) be a special semigroup set vector space defined over the set P. Let = ($\eta_1$, ..., $\eta_n$): V →[0, 1] be such that $\eta_i$: S$_i$ →[0, 1]; 1 ≤i ≤n,*
$$\eta_i(a_i + b_i) \geq min\{\eta_i(a_i), \eta_i(b_i)\}$$
$$\eta_i(ca_i) \geq \eta_i(a_i) \text{ for all } c \in P$$
*and for all a$_i$, b$_i$ ∈ S$_i$, true for i =1, 2, ..., n. Then we call*
$$V_\eta = (V_1, V_2, ..., V_n)_{(\eta_1,...,\eta_n)} = (V_{1\eta_1}, V_{2\eta_2}, ..., V_{n\eta_n})$$
*to be a special semigroup set fuzzy vector space.*

We illustrate this by some examples.

***Example 4.2.25:*** Let V = (V$_1$, V$_2$, V$_3$, V$_4$, V$_5$) be a special semigroup set vector space over the set {0, 1, 2, ..., ∞} = S where V$_1$ = {$Z^o \times Z^o \times Z^o$ | $Z^o$ S = $Z^+ \cup \{0\}$},

$$V_2 = \left\{ \begin{pmatrix} a \\ b \\ c \\ d \end{pmatrix} \middle| a,b,c,d \in S \right\},$$

V$_3$ = {all 3×3 matrices with entries from S},

$$V_4 = \left\{ \begin{pmatrix} a & a & a & a \\ b & b & b & b \end{pmatrix} \text{ such that } a, b \in S \right\}$$

and

$$V_5 = \left\{ \begin{pmatrix} a & b \\ c & d \\ a & b \\ c & d \\ a & b \end{pmatrix} \text{ such that } a, b, c, d \in S \right\}.$$

Define $\eta = (\eta_1, \eta_2, \eta_3, \eta_4, \eta_5): V = (V_1, V_2, V_3, V_4, V_5) \to [0, 1]$ by $\eta_i: V_i \to [0, 1]$; $1 \leq i \leq 5$ as follows:

$\eta_1 : V_1 \to [0,1]$ is given by

$$\eta_1(x, y, z) = \begin{cases} \dfrac{1}{x+y+z} & \text{if } x+y+z \neq 0 \\ 1 & \text{if } x+y+z = 0 \end{cases}$$

$\eta_2: V_2 \to [0, 1]$ is defined by

$$\eta_2 \begin{bmatrix} a \\ b \\ c \\ d \end{bmatrix} = \begin{cases} \dfrac{1}{a+b+c+d} & \text{if } a+b+c+d \neq 0 \\ 1 & \text{if } a+b+c+d = 0 \end{cases}$$

$\eta_3: V_3 \to [0, 1]$ is such that

$$\eta_3 \begin{pmatrix} a & b & c \\ d & e & f \\ g & h & i \end{pmatrix} = \begin{cases} \dfrac{1}{a+e+i} & \text{if } a+e+i \neq 0 \\ 1 & \text{if } a+e+i = 0 \end{cases}$$

$\eta_4 : V_4 \to [0,1]$ is defined by

$$\eta_4 \begin{bmatrix} a & a & a & a \\ b & b & b & b \end{bmatrix} = \begin{cases} \dfrac{1}{4b+2a} & \text{if } 4b+2a \neq 0 \\ 1 & \text{if } 4b+2a = 0 \end{cases}$$

$\eta_5 : V_5 \to [0,1]$ is given by

$$\eta_5 \begin{bmatrix} a & b \\ c & d \\ a & b \\ c & d \\ a & b \end{bmatrix} = \begin{cases} \dfrac{1}{3a+2c} & \text{if } 3a+2c \neq 0 \\ 1 & \text{if } 3a+2c = 0 \end{cases}$$

Clearly
$$V_\eta = (V_1, V_2, V_3, V_4, V_5)_\eta = \left( V_{1_{\eta_1}}, V_{2_{\eta_2}}, ..., V_{5_{\eta_5}} \right)$$
is a special semigroup set fuzzy vector space.

*Example 4.2.26:* Let $V = (V_1, V_2, V_3, V_4)$ where

$$V_1 = Z_5 \times Z_5 \times Z_5 \times Z_5 \times Z_5,$$

$$V_2 = \left\{ \begin{pmatrix} a \\ b \\ c \\ d \end{pmatrix} \middle| a,b,c,d \in Z_5 \right\},$$

$V_3 = \{Z_5[x]$ all polynomials of finite degree with coefficients from $Z_5\}$ and $V_4 = \{$all $4 \times 4$ upper triangular matrices with entries from $Z_5\}$ is a special semigroup set vector space over the set $S = \{0, 1\}$.

Define $\eta = (\eta_1, \eta_2, \eta_3, \eta_4): V = (V_1\ V_2\ V_3\ V_4) \to [0,1]$ by $\eta_i: V_i \to [0,1]$, $i = 1, 2, 3, 4$ such that

$\eta_1 : V_1 \to [0,1]$ defined by

$$\eta_1(a, b, c, d, e) = \begin{cases} \dfrac{1}{a} & \text{if } a \neq 0 \\ 1 & \text{if } a = 0 \end{cases}$$

$\eta_2: V_2 \to [0, 1]$ is given by

$$\eta_2 \begin{bmatrix} a \\ b \\ c \\ d \end{bmatrix} = \begin{cases} \dfrac{1}{a+d} & \text{if } a+d \neq 0 \\ 1 & \text{if } a+d = 0 \end{cases}$$

$\eta_3: V_3 \to [0, 1]$ is such that

$$\eta_3(p(x)) = \begin{cases} \dfrac{1}{\deg p(x)} & \text{if } p(x) \neq \text{constant} \\ 1 & \text{if } p(x) \text{ is a constant} \end{cases}$$

$\eta_4: V_4 \to [0,1]$ is defined by

$$\eta_4 \begin{bmatrix} a & b & c & d \\ 0 & e & f & g \\ 0 & 0 & h & i \\ 0 & 0 & 0 & j \end{bmatrix} = \begin{cases} \dfrac{1}{a+e+h+j} & \text{if } a+e+h+j \neq 0 \\ 1 & \text{if } a+e+h+j = 0 \end{cases}$$

Clearly
$$V_\eta = \left(V_1, V_2, V_3, V_4\right)_{(\eta_1, \eta_2, \eta_3, \eta_4)} = \left(V_{1\eta_1}, V_{2\eta_2}, \ldots, V_{4\eta_4}\right)$$
is a special semigroup set fuzzy vector space.

Now we proceed on to define the notion of special semigroup set fuzzy vector subspace.

**DEFINITION 4.2.21:** *Let $V = (V_1, V_2, \ldots, V_n)$ be a special semigroup set vector space over the set S. Suppose $W = (W_1, W_2, \ldots, W_n) \subseteq (V_1, V_2, \ldots, V_n)$; $W_i \subseteq V_i$; ($1 \leq i \leq n$) is a special* *semigroup set vector subspace of V over the set S. Now define $\eta = (\eta_1, \eta_2, \ldots, \eta_n) : W \to [0, 1]; W = (W_1, W_2, \ldots, W_n) \to [0, 1]$ such that $\eta_i : W_i \to [0, 1]$ where $W_\eta = \left(W_{\eta_1}, W_{\eta_2}, \ldots, W_{\eta_n}\right)$ is a special semigroup set fuzzy vector space then we call $W_\eta$ to be a special semigroup set fuzzy vector subspace of V.*

We illustrate this by a simple example.

***Example 4.2.27:*** Let $V = (V_1, V_2, V_3, V_4, V_5)$ be a special semigroup set vector space over the set $S = \{0, 1\}$ where $V_1 = \{Z^o \times Z^o \times Z^o \times Z^o \times Z^o \mid Z^o = Z^+ \cup \{0\}\}$, $V_2 = \{\text{all } 4 \times 5 \text{ matrices with entries from } Z^o\}$,

$$V_3 = \left\{ \begin{pmatrix} a \\ b \\ c \\ d \end{pmatrix} \,\middle|\, a, b, c, d, e \in Z^o \right\}$$

$V_4 = \{\text{all } 4 \times 4 \text{ lower triangular matrices with entries from } Z^o\}$ and $V_5 = \{Z^o[x] \text{ all polynomials in the variable x with coefficients from } Z^o\}$. Take $W = (W_1, W_2, W_3, W_4, W_5) \subseteq (V_1, V_2, V_3, V_4, V_5) = V$ where $W_1 = \{2Z^o \times 3Z^o \times 4Z^o \times 5Z^o \times 6Z^o\} \subseteq V_1$,

$$W_2 = \left\{ \begin{pmatrix} a_1 & a_2 & a_3 & a_4 & a_5 \\ b_1 & b_2 & b_3 & b_4 & b_5 \\ c_1 & c_2 & c_3 & c_4 & c_5 \\ d_1 & d_2 & d_3 & d_4 & d_5 \end{pmatrix} \,\middle|\, a_i, b_i, c_i, d_i \in 3Z^o, 1 \le i \le 5 \right\} \subseteq V_2,$$

$$W_3 = \left\{ \begin{pmatrix} a \\ a \\ a \\ a \\ a \end{pmatrix} \,\middle|\, a \in Z^o \right\} \subseteq V_3,$$

$$W_4 = \left\{ \begin{pmatrix} a & 0 & 0 & 0 \\ a & a & 0 & 0 \\ b & b & a & 0 \\ c & c & b & a \end{pmatrix} \middle| a, b, c \in Z^o \right\} \subseteq V_4$$

and $W_5 = \{$all polynomials in $Z^o[x]$ of degree less than or equal to 11 with coefficients from $Z^o\} \subseteq V_5$. Clearly $W = (W_1, W_2, W_3, W_4, W_5)$ is special semigroup set vector subspace of V over the set $\{0, 1\}$.

Define $\eta = (\eta_1, \eta_2, \eta_3, \eta_4, \eta_5) : W = (W_1, W_2, W_3, W_4, W_5) \to [0, 1]$ as follows:
$$\eta_i : W_i \to [0,1]; \, i = 1, 2, 3, 4, 5.$$

$\eta_1: W_1 \to [0, 1]$ is defined by

$$\eta_1(a\ b\ c\ d\ e) = \begin{cases} \dfrac{1}{a+b+c+d+e} & \text{if } a+b+c+d+e \neq 0 \\ 1 & \text{if } a+b+c+d+e = 0 \end{cases}$$

$\eta_2 : W_2 \to [0,1]$ is such that

$$\eta_2 \begin{bmatrix} a_1 & a_2 & a_3 & a_4 & a_5 \\ d_1 & d_2 & d_3 & d_4 & d_5 \\ b_1 & b_2 & b_3 & b_4 & b_5 \\ p_1 & p_2 & p_3 & p_4 & p_5 \end{bmatrix} =$$

$$\begin{cases} \dfrac{1}{a_1 + d_2 + b_3 + p_4} & \text{if } a_1 + d_2 + b_3 + p_4 \neq 0 \\ 1 & \text{if } a_1 + d_2 + b_3 + p_4 = 0 \end{cases}$$

$\eta_3 : W_3 \to [0,1]$ is given by

$$\eta_3 \begin{pmatrix} a \\ a \\ a \\ a \\ a \\ a \end{pmatrix} = \begin{cases} \dfrac{1}{a} & \text{if } a \neq 0 \\ 1 & \text{if } a = 0 \end{cases}$$

$\eta_4 : W_4 \to [0,1]$ is defined by

$$\eta_4 \begin{pmatrix} a & 0 & 0 & 0 \\ a & a & 0 & 0 \\ b & b & a & 0 \\ c & c & b & a \end{pmatrix} = \begin{cases} \dfrac{1}{a+b+c} & \text{if } a+b+c \neq 0 \\ 1 & \text{if } a+b+c = 0 \end{cases}$$

Clearly $W = \left( W_{1\eta_1}, W_{2\eta_2}, W_{3\eta_3}, W_{4\eta_4} \right)$ is special semigroup set fuzzy vector subspace.

Now we proceed on to define the notion of special semigroup set linear algebra.

**DEFINITION 4.2.22:** *Let $V = (V_1, V_2, ..., V_n)$ be a special semigroup set vector space defined over the set S. Let $\eta = (\eta_1, \eta_2, ..., \eta_n): V = (V_1, V_2, ..., V_n) \to [0,1]$ such that $\eta_i : V_i \to [0, 1]$; $i = 1, 2, ..., n$.*

*$\eta_i(a_i + b_i) \geq \min\{\eta_i(a_i), \eta_i(b_i)\}$*

*$\eta_i(ra_i) \geq \eta_i(a_i)$ for all $a_i, b_i \in V_i$ and for all $r \in S$; $1 \leq i \leq n$.*
*Then $V_\eta = \left( V_{1\eta_1}, V_{2\eta_2}, ..., V_{n\eta_n} \right)$ is a special semigroup set fuzzy linear algebra.*

We illustrate this by an example.

***Example 4.2.28:*** Let $V = (V_1, V_2, V_3, V_4, V_5)$ where $V_1 = \{Z_{12} \times Z_{12} \times Z_{12}\}$, $V_2 = \{$all $3 \times 3$ matrices with entries from $Z_{12}\}$, $V_3 = \{Z_{12}[x]$ be all polynomials of degree less than or equal to $7\}$,

$$V_4 = \left\{ \begin{pmatrix} a \\ b \\ c \\ d \\ e \end{pmatrix} \middle| a, b, c, d, e \in Z_{12} \right\}$$

and $V_5 = \{$all $2 \times 5$ matrices with entries from $Z_{12}\}$ be a special semigroup set linear algebra over $Z_{12}$.

$$\eta = (\eta_1, \eta_2, \eta_3, \eta_4, \eta_5): V = (V_1\ V_2\ V_3\ V_4\ V_5) \to [0, 1];$$

$$\eta_i: W_i \to [0, 1];\ 1 \leq i \leq 5.$$

$\eta_1: W_1 \to [0, 1]$ is such that

$$\eta_1(a\ b\ c) = \begin{cases} \dfrac{1}{a+b+c} & \text{if } a+b+c \neq 0 \\ 1 & \text{if } a+b+c = 0 \end{cases}$$

$\eta_2 : V_2 \to [0,1]$ is defined by

$$\eta_2 \begin{bmatrix} a & b & c \\ d & e & f \\ g & h & i \end{bmatrix} = \begin{cases} \dfrac{1}{2a+3b+4i} & \text{if } 2a+3b+4i \neq 0 \\ 1 & \text{if } 2a+3b+4i = 0 \end{cases}$$

$\eta_3 : V_3 \to [0,1]$ is given by

$$\eta_3(p(x)) = \begin{cases} \dfrac{1}{\deg p(x)} & \text{if } p(x) \neq \text{constant} \\ 1 & \text{if } p(x) \text{ is a constant} \end{cases}$$

$\eta_4 : V_4 \to [0,1]$ is given by

$$\eta_4 \begin{bmatrix} a \\ b \\ c \\ d \\ e \end{bmatrix} = \begin{cases} \dfrac{1}{a+e} & \text{if } a+e \neq 0 \\ 1 & \text{if } a+e = 0 \end{cases}$$

$\eta_5 : V_5 \to [0,1]$ is defined by

$$\eta_5 \begin{bmatrix} a_1 & a_2 & a_3 & a_4 & a_5 \\ b_1 & b_2 & b_3 & b_4 & b_5 \end{bmatrix} = \begin{cases} \dfrac{1}{2a_1 + b_5} & \text{if } 2a_1 + b_5 \neq 0 \\ 1 & \text{if } 2a_1 + b_5 = 0 \end{cases}$$

$V_\eta = \left( V_{1\eta_1}, V_{2\eta_2}, V_{3\eta_3}, V_{4\eta_4}, V_{5\eta_5} \right)$ is a special semigroup set fuzzy linear algebra.

It is important to state that like most of the structures defined in this book; we see special semigroup set fuzzy vector space and special semigroup set fuzzy linear algebra are fuzzy equivalent. Now we proceed on to define special semigroup set fuzzy linear sub algebras.

**DEFINITION 4.2.23:** *Let $V = (V_1, V_2, ..., V_n)$ be a special semigroup set linear algebra over a semigroup S. Let $W = (W_1, W_2, ..., W_n)$ be a special semigroup set linear subalgebra of V over the semigroup S.*

*Define $W = (W_1, W_2, ..., W_n) \to [0, 1]$; $\eta = (\eta_1, \eta_2, ..., \eta_n)$: $(W_1, W_2, ..., W_n) \to [0,1]$; where $\eta_i : V_i \to [0,1]$; $i = 1, 2, ..., n$. such that*
$$W_\eta = (W_1, W_2, ..., W_n)_\eta = \left( W_{1\eta_1}, W_{2\eta_2}, ..., W_{n\eta_n} \right)$$
*is a special semigroup set fuzzy linear algebra then we call $W_\eta$ is a special semigroup set fuzzy linear subalgebra.*

We now illustrate this by a simple example.

***Example 4.2.29:*** Let $V = (V_1, V_2, V_3, V_4, V_5)$ where $V_1 = \{Z^o \times Z^o \times Z^o \times Z^o \times Z^o\}$, $V_2 = \{Z^o[x]$ is a polynomial of degree less than or equal to 12 with coefficients from $Z^o\}$, $V_3 = \{$all $5 \times 5$ lower triangular matrices with entries from $Z^o\}$, $V_4 = \{3 \times 6$ matrices with entries from $Z^o\}$ and

$$V_5 = \left\{ \begin{pmatrix} a \\ a \\ a \\ a \\ a \\ a \\ a \end{pmatrix} \middle| a \in Z^o \right\}$$

be a special semigroup set linear algebra over $S = \{1, 2, \ldots, 10\}$. $W = (W_1, W_2, W_3, W_4, W_5) \subseteq (V_1, V_2, V_3, V_4, V_5) = V$ where

$$W_1 = \{Z^o \times 2Z^o \times \{0\} \times \{0\} \times 3Z^o\} \subseteq V_1,$$
$W_2 = \{$all polynomials in x with coefficients from $Z^o$ of degree less than or equal to 7$\} \subseteq V_2$,
$W_3 = \{$all $5 \times 5$ lower triangular matrices with entries $3Z^o\} \subseteq V_3$,
$W_4 = \{3 \times 6$ matrices with entries $5Z^o\} \subseteq V_4$

and

$$W_5 = \left\{ \begin{pmatrix} a \\ a \\ a \\ a \\ a \\ a \\ a \end{pmatrix} \middle| a \in 3Z^o \right\} \subseteq V_5,$$

W is a special semigroup set linear subalgebra of V over the set $S = \{0, 1, 2, \ldots, 10\}$.
Define $\eta: W \to [0, 1]$; by
$\eta = (\eta_1, \eta_2, \eta_3, \eta_4, \eta_5): (W_1, W_2, W_3, W_4, W_5) = W \to [0, 1]$, $1 \leq i \leq 5$.

$\eta_1: W_1 \to [0, 1]$ is defined by

$$\eta_1(a\ b\ 0\ 0\ c) = \begin{cases} \dfrac{1}{b} & \text{if } b \neq 0 \\ 1 & \text{if } b = 0 \end{cases}$$

$\eta_2: W_2 \to [0, 1]$ is given by

$$\eta_2(p(x)) = \begin{cases} \dfrac{1}{\deg p(x)} & \text{if } p(x) \text{ is not a constant} \\ 1 & \text{if } p(x) \text{ is a constant} \end{cases}$$

$\eta_3: W_3 \to [0,1]$ is defined by

$$\eta_3 \begin{bmatrix} a & 0 & 0 & 0 & 0 \\ b & b & 0 & 0 & 0 \\ c & c & c & 0 & 0 \\ d & d & d & d & 0 \\ e & e & e & e & e \end{bmatrix} = \begin{cases} \dfrac{1}{a+b+c+d+e} & \text{if } a+b+c+d+e \neq 0 \\ 1 & \text{if } a+b+c+d+e = 0 \end{cases}$$

$\eta_4: W_4 \to [0,1]$ is such that

$$\eta_4 \begin{bmatrix} a_1 & a_2 & a_3 & a_4 & a_5 & a_6 \\ b_1 & b_2 & b_3 & b_4 & b_5 & b_6 \\ c_1 & c_2 & c_3 & c_4 & c_5 & c_6 \end{bmatrix} = \begin{cases} \dfrac{1}{a_1 + b_3 + c_5} & \text{if } a_1 + b_3 + c_5 \neq 0 \\ 1 & \text{if } a_1 + b_3 + c_5 = 0 \end{cases}$$

$\eta_5: W_5 \to [0,1]$ is given by

$$\eta_5 \begin{bmatrix} a \\ a \\ a \\ a \\ a \\ a \\ a \end{bmatrix} = \begin{cases} \dfrac{1}{3a} & \text{if } a \neq 0 \\ 1 & \text{if } a = 0 \end{cases}$$

$$W_\eta = \left(W_1, W_2, W_3, W_4, W_5\right)_{(\eta_1,\eta_2,...,\eta_5)} =$$
$$\left(W_{1\eta_1}, W_{2\eta_2}, W_{3\eta_3}, W_{4\eta_4}, W_{5\eta_5}\right)$$

is a special semigroup set fuzzy linear subalgebra.

It is pertinent to mention here that the notion of special semigroup set fuzzy vector subspace and special semigroup set fuzzy linear subalgebra are fuzzy equivalent.

Now we proceed on to define the fuzzy notion of special group set vector spaces defined over the group G.

**DEFINITION 4.2.24:** *Let $V = (V_1, V_2, ..., V_n)$ be a special group set vector space over the set S. Let $\eta = (\eta_1, \eta_2,..., \eta_n): V \to (V_1, V_2, ..., V_n) \to [0, 1]$ such that $\eta_i: V_i \to [0, 1]; i = 1 \leq i \leq n$.*

$$\eta_i(a_i + b_i) \geq \min\{\eta_i(a_i), \eta_i(b_i)\}$$
$$\eta_i(-a_i) \geq \eta_i(a_i); \eta_i = 0$$
$$\eta_i(ra_i) \geq \eta_i(a_i)$$

*for all $a_i, b_i \in V_i$ and $r \in S$. If this is true for every $i$; $i = 1,2,...,n$.*

$$V_\eta = \left(V_1, V_2,..., V_n\right)_{(\eta_1,\eta_2,...,\eta_n)} = \left(V_{1\eta_1}, V_{2\eta_2},..., V_{n\eta_n}\right)$$

*is a special group set fuzzy vector space.*

We shall illustrate this by a simple example.

*Example 4.2.30:* Let $V = (V_1\ V_2\ V_3\ V_4\ V_5)$ to be special group set vector space over the set $Z^o$ where $V_1 = \{Z \times Z \times Z \times Z \times Z \times Z\}$, $V_2 = \{$all $3 \times 3$ matrices with entries from $Z\}$,

$$V_3 = \left\{ \begin{pmatrix} a \\ b \\ c \\ d \\ e \end{pmatrix} \;\middle|\; a,b,c,d,e \in Z \right\},$$

$V_4 = \{$all $3 \times 5$ matrices with entries from $Z\}$ and $V_5 = \{$all $5 \times 5$ lower triangular matrices with entries from $Z\}$.

Define $\eta = (\eta_1, \eta_2, \eta_3, \eta_4, \eta_5) : V = (V_1, V_2, V_3, V_4, V_5) \to [0,1]$; $1 \leq i \leq 5$.

$\eta_1: V_1 \to [0, 1]$ is such that

$$\eta_1(a\ b\ c\ d\ e) = \begin{cases} \dfrac{1}{a+c+f} & \text{if } a+c+f \neq 0 \\ 1 & \text{if } a+c+f = 0 \end{cases}$$

$\eta_2: V_2 \to [0, 1]$ is defined by

$$\eta_2 \begin{bmatrix} a & b & c \\ d & e & f \\ g & h & i \end{bmatrix} = \begin{cases} \dfrac{1}{a+b+d+i} & \text{if } a+b+d+i \neq 0 \\ 1 & \text{if } a+b+d+i = 0 \end{cases}$$

$\eta_3: V_3 \to [0,1]$ is given by

$$\eta_3 \begin{bmatrix} a \\ b \\ c \\ d \\ e \end{bmatrix} = \begin{cases} \dfrac{1}{5a+7d+e} & \text{if } 5a+7d+e \neq 0 \\ 1 & \text{if } 5a+7d+e = 0 \end{cases}$$

$\eta_4: V_4 \to [0, 1]$ is such that

$$\eta_4 \begin{bmatrix} a_1 & a_2 & a_3 & a_4 & a_5 \\ b_1 & b_2 & b_3 & b_4 & b_5 \\ c_1 & c_2 & c_3 & c_4 & c_5 \end{bmatrix} = \begin{cases} \dfrac{1}{a_5+b_5+c_5} & \text{if } a_5+b_5+c_5 \neq 0 \\ 1 & \text{if } a_5+b_5+c_5 = 0 \end{cases}$$

$\eta_5: V_5 \to [0, 1]$ is given by

$$\eta_5 \begin{bmatrix} a & 0 & 0 & 0 & 0 \\ b & c & 0 & 0 & 0 \\ d & e & f & 0 & 0 \\ g & h & i & j & 0 \\ k & l & m & n & p \end{bmatrix} = \begin{cases} \dfrac{1}{a+p} & \text{if } a+p \neq 0 \\ 1 & \text{if } a+p = 0 \end{cases}$$

Thus $V_\eta = \left(V_{1\eta_1}, V_{2\eta_2}, ..., V_{5\eta_5}\right)$ is a special group set fuzzy vector space.

Now we proceed on to define the notion of special group set fuzzy vector subspace.

**DEFINITION 4.2.25:** *Let $V = (V_1, V_2, ..., V_n)$ be a special group set vector space over the set S. Let $W = (W_1, W_2, ..., W_n) \subseteq (V_1, V_2, ..., V_n) \subseteq V$ be a proper subset of V which is a special group set vector subspace of V over S.*

*Let $\eta = (\eta_1, \eta_2, ..., \eta_n): W = (W_1, W_2, ..., W_n) \to [0,1]$ such that $\eta_i: W_i \to [0, 1]$ for every i and $W_\eta = \left(W_{1\eta_1}, W_{2\eta_2}, ..., W_{n\eta_n}\right)$ be a special group set fuzzy vector space, then we define $W_\eta$ to be a special group set fuzzy vector subspace of V.*

Now we will illustrate this situation by a simple example.

*Example 4.2.31:* Let $V = (V_1, V_2, V_3, V_4, V_5)$ be a space over the set $S = Z$ where
$$V_1 = S \times S \times S \times S \times S,$$
$V_2 = \{$all $3 \times 3$ matrices with entries from $S\}$,
$V_3 = \{$all $5 \times 2$ matrices with entries from $S = Z\}$
$V_4 = \{(a\ a\ a\ a\ a\ a\ a) \mid a \in S\}$
and

$$V_5 = \left\{ \begin{pmatrix} a_1 \\ a_2 \\ a_3 \\ a_4 \\ a_5 \end{pmatrix} \middle| a_i \in S = Z \quad 1 \leq i \leq 5 \right\}.$$

Take $W = (W_1, W_2, ..., W_5) \subseteq (V_1, V_2, V_3, V_4, V_5) = V$ where

$$W_1 = S \times S \times \{0\} \times \{0\} \times S \subseteq V_1,$$
$W_2 = \{$all $3 \times 3$ upper triangular matrices with entries from $S\} \subseteq V_2$,

$$W_3 = \left\{ \begin{pmatrix} a & b \\ a & b \\ a & b \\ a & b \\ a & b \\ a & b \end{pmatrix} \middle| a, b \in S \right\} \subseteq V_3,$$

$$W_4 = \{(a\ a\ a\ a\ a\ a) \mid a \in 5S\} \subseteq V_4$$

and

$$W_5 = \left\{ \begin{pmatrix} a \\ a \\ a \\ a \\ a \end{pmatrix} \middle| a \in S \right\} \subseteq V_5.$$

Clearly $W = (W_1, W_2, W_3, W_4, W_5) \subseteq (V_1, V_2, V_3, V_4, V_5) = V$ is a special group set vector subspace of V over S.

Define $\eta = (\eta_1, \eta_2, \eta_3, \eta_4, \eta_5) : W = (W_1, W_2, W_3, W_4, W_5) \to [0,1]$; as $\eta_i: W_i \to [0, 1]$ for every i such that

$\eta_1: W_1 \to [0,1]$ is defined by

$$\eta_1(a\ b\ 0\ 0\ c) = \begin{cases} \dfrac{1}{a+b+c} & \text{if } a+b+c \neq 0 \\ 1 & \text{if } a+b+c = 0 \end{cases}$$

$\eta_2 : W_2 \to [0,1]$ is given by

$$\eta_2 \begin{pmatrix} a & b & c \\ 0 & d & e \\ 0 & 0 & f \end{pmatrix} = \begin{cases} \dfrac{1}{a+d+f} & \text{if } a+d+f \neq 0 \\ 1 & \text{if } a+d+f = 0 \end{cases}$$

$\eta_3 : W_3 \to [0,1]$ is given by

$$\eta_3 \begin{bmatrix} a & b \\ a & b \\ a & b \\ a & b \\ a & b \end{bmatrix} = \begin{cases} \dfrac{1}{a+b} & \text{if } a+b \neq 0 \\ 1 & \text{if } a+b = 0 \end{cases}$$

$\eta_4 : W_4 \to [0,1]$ is given by

$$\eta_4(a\ a\ a\ a\ a\ a) = \begin{cases} \dfrac{1}{a} & \text{if } a \neq 0 \\ 1 & \text{if } a = 0 \end{cases}$$

$\eta_5 : W_5 \to [0,1]$ is such that

$$\eta_5 \begin{bmatrix} a \\ a \\ a \\ a \\ a \end{bmatrix} = \begin{cases} \dfrac{1}{5a} & \text{if } a \neq 0 \\ 1 & \text{if } a = 0 \end{cases}$$

Thus $W_\eta = \left( W_{1\eta_1}, W_{2\eta_2}, W_{3\eta_3}, W_{4\eta_4}, W_{5\eta_5} \right) = (W_1\eta_1, W_2\eta_2, W_3\eta_3, W_4\eta_4, W_5\eta_5)$ is a special group set fuzzy vector subspace.

Now we proceed on to define the notion of special group set linear algebra V defined over the group G.

**DEFINITION 4.2.26:** *Let $V = (V_1, V_2, ..., V_n)$ be a special group set linear algebra over the group G. Let $\eta : V = (V_1, V_2, ..., V_n) \to [0,1]$ be such that $\eta = (\eta_1, \eta_2, ..., \eta_n): (V_1, V_2, ..., V_n) \to [0, 1]$ such that $\eta_i : V_i \to [0,1]$ for each i, i = 1, 2, ..., n; satisfying the following conditions*

$$\eta_i(a_i + b_i) \geq \min\{\eta_i(a_i), \eta_i(b_i)\}$$
$$\eta_i(-a_i) \geq \eta_i(a_i); \ \eta_i(0) = 1$$
$$\eta_i(ra_i) \geq \eta_i(a_i)$$

*for atleast one pair of $a_i, b_i \in V_i$ and for some $r \in G$; true for each i, i = 1, 2, ..., n. We call $V_\eta = (V_{1\eta_1}, V_{2\eta_2}, ..., V_{n\eta_n})$ to be a pseudo special group set fuzzy linear algebra.*

We illustrate this by some simple examples.

***Example 4.2.32:*** Let $V = (V_1, V_2, V_3, V_4, V_5)$ be a special group set linear algebra over the group $Z_{18}$ where $V_1 = Z_{18} \times Z_{18} \times Z_{18}$, $V_2 = \{$all $4 \times 4$ matrices with entries from $Z_{18}\}$, $V_3 = \{$all $5 \times 5$ lower triangular matrices with entries from $Z_{18}\}$,

$$V_4 = \left\{ \begin{pmatrix} a \\ b \\ c \\ d \\ e \end{pmatrix} \middle| a,b,c,d,e \in Z_{18} \right\}$$

and

$$V_5 = \left\{ \begin{pmatrix} a & a \\ b & b \\ c & c \\ d & d \\ e & e \\ f & f \end{pmatrix} \middle| a,b,c,d,e,f \in Z_{18} \right\}.$$

Define $\eta = (\eta_1, \eta_2, \eta_3, \eta_4, \eta_5) : V = (V_1, V_2, V_3, V_4, V_5) \to [0,1]$; such that $\eta_i : V_i \to [0,1]$, i = 1, 2, ..., 5.

$\eta_1 : V_1 \to [0,1]$ is defined by

$$\eta_1 (a\ b\ c) = \begin{cases} \dfrac{1}{a+b+c} & \text{if } a+b+c \neq 0 \\ 1 & \text{if } a+b+c = 0 \end{cases}$$

$\eta_2: V_2 \to [0, 1]$ is such that

$$\eta_2 \begin{pmatrix} a & b & c & d \\ e & f & g & h \\ c & j & k & l \\ m & n & o & p \end{pmatrix} = \begin{cases} \dfrac{1}{d+f+j+m} & \text{if } d+f+j+m \neq 0 \\ 1 & \text{if } d+f+j+m = 0 \end{cases}$$

$\eta_3: V_3 \to [0, 1]$ is defined by

$$\eta_3 \begin{bmatrix} a & 0 & 0 & 0 & 0 \\ b & c & 0 & 0 & 0 \\ d & e & f & 0 & 0 \\ g & h & i & j & 0 \\ k & l & m & n & p \end{bmatrix} = \begin{cases} \dfrac{1}{a+d+b+g+k} & \text{if } a+d+b+g+k \neq 0 \\ 1 & \text{if } a+d+b+g+k = 0 \end{cases}$$

$\eta_4: W_4 \to [0, 1]$ is given by

$$\eta_4 \begin{bmatrix} a \\ b \\ c \\ d \\ e \end{bmatrix} = \begin{cases} \dfrac{1}{a+c+e} & \text{if } a+c+e \neq 0 \\ 1 & \text{if } a+c+e = 0 \end{cases}$$

$\eta_5 : W_5 \to [0, 1]$ is defined by

$$\eta_5 \begin{bmatrix} a & a \\ b & b \\ c & c \\ d & d \\ e & e \\ f & f \end{bmatrix} = \begin{cases} \dfrac{1}{a+b+c+d+e+f} & \text{if } a+b+c+d+e+f \neq 0 \\ 1 & \text{if } a+b+c+d+e+f = 0 \end{cases}$$

Clearly $V_\eta = \left(V_{1\eta_1}, V_{2\eta_2}, V_{3\eta_3}, V_{4\eta_4}, V_{5\eta_5}\right)$ is a special group set fuzzy linear algebra.

It is pertinent to record at this juncture the notion of special group set fuzzy vector space and special group set fuzzy linear algebra are fuzzy equivalent. However we wish to state they become the same under the fuzzification.

Now we proceed on to define the notion of special group set fuzzy linear subalgebra.

**DEFINITION 4.2.27:** *Let $V = (V_1, V_2, \ldots, V_n)$ be a special group set linear algebra defined over the group G. Let $W = (W_1, W_2, \ldots, W_n) \subseteq (V_1, V_2, \ldots, V_n) = V$ such that $W_i \subseteq V_i$, $i = 1, 2, \ldots, n$; be a special group set linear subalgebra of V over the same group G. Let $\eta : (\eta_1, \eta_2, \ldots, \eta_n) : W = (W_1, W_2, \ldots, W_n) \to [0,1]$ be defined such that $\eta_i : W_i \to [0,1]$ and*
$$W_\eta = \left(W_{1\eta_1}, W_{2\eta_2}, \ldots, W_{n\eta_n}\right)$$
*is a special group set fuzzy linear algebra, then we call W to be a special group set fuzzy linear set subalgebra.*

We illustrate this by some an example.

*Example 4.2.33:* Let $V = (V_1, V_2, V_3, V_4)$ be a special group set linear algebra over the group R. Here $V_1 = R \times R \times R \times R$; R the reals under addition, $V_2 = \{$all $4 \times 4$ matrices with entries from R$\}$,

$$V_3 = \left\{ \begin{pmatrix} a_1 \\ a_2 \\ a_3 \\ a_4 \\ a_5 \end{pmatrix} \,\middle|\, a_i \in R; 1 \leq i \leq 5 \right\}$$

and $V_4 = \{$all $2 \times 7$ matrices with entries from R$\}$. Choose $W = (W_1, W_2, W_3, W_4)$ where $W_1 = R \times \{0\} \times \{0\} \times R \subseteq V_1$; $W_2 = \{$all 4 matrices with entries from R of the form

$$\left.\begin{Bmatrix} \begin{pmatrix} a & a & a & a \\ b & b & b & b \\ c & c & c & c \\ d & d & d & d \end{pmatrix} \middle| a,b,c,d \in R \end{Bmatrix} \subseteq V_2,\right.$$

$$W_3 = \left\{ \begin{pmatrix} a \\ a \\ a \\ a \\ a \end{pmatrix} \middle| a \in R \right\} \subseteq V_3,$$

and

$$W_4 = \left\{ \begin{pmatrix} a & a & a & a & a & a & a \\ b & b & b & b & b & b & b \end{pmatrix} \middle| a,b \in R \right\} \subseteq V_4,$$

so that W is easily verified to be a special group set linear subalgebra of V over the real. Define

$\eta = (\eta_1, \eta_2, \eta_3, \eta_4) : W = (W_1, W_2, W_3, W_4) \to [0,1]$;
$\eta_i : V_i \to [0, 1]$; $i = 1, 2, 3, 4$ as follows.

$\eta_1 : W_1 \to [0,1]$ such that

$$\eta_1(a\ 0\ 0\ b) = \begin{cases} \dfrac{1}{a+b} & \text{if } a+b \geq 1 \\ a+b & \text{if } a+b < 1 \\ 1 & \text{if } a+b = 0 \end{cases}$$

$\eta_2 : W_2 \to [0,1]$ is defined by

$$\eta_2 \begin{pmatrix} a & a & a & a \\ b & b & b & b \\ c & c & c & c \\ d & d & d & d \end{pmatrix} = \begin{cases} \dfrac{1}{a+b+c+d} & \text{if } a+b+c+d \geq 1 \\ a+b+c+d & \text{if } a+b+c+d < 1 \\ 1 & \text{if } a+b+c+d = 0 \end{cases}$$

$\eta_3 : W_3 \to [0, 1]$ is given by

$$\eta_3 \begin{bmatrix} a \\ a \\ a \\ a \\ a \end{bmatrix} = \begin{cases} \dfrac{1}{a} & \text{if } a \geq 1 \\ a & \text{if } a < 1 \\ 1 & \text{if } a = 0 \end{cases}$$

$\eta_4 : W_4 \to [0,1]$ is given by

$$\eta_4 \begin{bmatrix} a & a & a & a & a & a & a \\ b & b & b & b & b & b & b \end{bmatrix} = \begin{cases} \dfrac{a}{b} & \text{if } \dfrac{a}{b} < 1 \\ \dfrac{b}{a} & \text{if } \dfrac{b}{a} \leq 1 \neq 0 \text{ i.e., } \dfrac{a}{b} \geq 1 \\ 1 & \text{if } a = 0 \text{ and } b = 0 \end{cases}$$

Thus $W_\eta = \left( W_{1\eta_1}, W_{2\eta_2}, W_{3\eta_3}, W_{4\eta_4} \right)$ is a special group set fuzzy linear subalgebra.

It is to be observed that the notion of special group set fuzzy vector subspace and special group set fuzzy linear subalgebra are also fuzzy equivalent.

Now we proceed on to define yet another new notion.

**DEFINITION 4.2.28:** *Let*
$$V = (V_1, V_2, \ldots, V_n)$$
$$= \left( V_1^1, V_2^1, \ldots, V_{n_1}^1 \right) \cup \left( V_1^2, V_2^2, \ldots, V_{n_2}^2 \right) \cup \ldots \cup \left( V_1^n, V_2^n, \ldots, V_{n_n}^n \right)$$
*be a special semigroup set n-vector space over the set a S. If*
$$\eta = (\eta_1 \cup \eta_2 \cup \ldots \cup \eta_n)$$
$$= \left( \eta_1^1, \eta_2^1, \ldots, \eta_{n_1}^1 \right) \cup \left( \eta_1^2, \eta_2^2, \ldots, \eta_{n_2}^2 \right) \cup \ldots \cup \left( \eta_1^n, \eta_2^n, \ldots, \eta_{n_n}^n \right) :$$
$$V = (V_1, V_2, \ldots, V_n)$$
$$= \left( V_1^1, V_2^1, \ldots, V_{n_1}^1 \right) \cup \left( V_1^2, V_2^2, \ldots, V_{n_2}^2 \right) \cup \ldots \cup \left( V_1^n, V_2^n, \ldots, V_{n_2}^n \right)$$
$$\to [0, 1]$$

such that for each i, $\eta^i_{j_i} : V^i_{j_i} \to [0, 1]$; $(1 \leq j_i \leq n_i; 1 \leq i \leq n)$ then $V_\eta = (V_{1\eta_1} \cup ... \cup V_{n\eta_n})$ is a special semigroup set fuzzy n-vector space; when n = 2 we get the special semigroup set fuzzy vector bispace. If n = 3 we get the special set fuzzy vector trispace.

If we replace the special semigroup set n-vector space V by a special semigroup set linear n-algebra V and define $\eta = (\eta_1 \cup ... \cup \eta_n): V = (V_1, V_2, ..., V_n) \to [0, 1]$; we get the special semigroup set fuzzy linear n-algebra. Infact both these concepts are fuzzy equivalent.

We illustrate this by an example.

*Example 4.2.34:* Let
$$V = (V_1 \cup V_2 \cup V_3 \cup V_4)$$
$$= (V^1_1, V^1_2, V^1_3) \cup (V^2_1, V^2_2) \cup (V^3_1, V^3_2, V^3_3, V^3_4) \cup (V^4_1, V^4_2, V^4_3)$$
be a special semigroup set vector 4 space over the semigroup S $= Z^o = Z^+ \cup \{0\}$, where
$$V^1_1 = Z^o \times Z^o \times Z^o,$$

$$V^1_2 = \left\{ \begin{pmatrix} a \\ b \\ c \\ d \end{pmatrix} \middle| a,b,c,d \in Z^o \right\},$$

$V^1_3 = \{$all $3 \times 3$ matrices with entries from $Z^o\}$,
$$V^2_1 = Z^o \times Z^o \times Z^o \times Z^o \times Z^o$$

and

$V^2_2 = \{5 \times 2$ matrices with entries from $Z^o\}$,
$V^3_1 = \{2 \times 6$ matrices with entries from $Z^o\}$,
$$V^3_2 = \{ Z^o \times Z^o \times Z^o \times Z^o \}$$

$V^3_3 = \{Z^o[x]$ all polynomials in the variable x with coefficients from $Z^o$ of degree less than or equal to four$\}$

and

$$V_4^3 = \left\{ \begin{pmatrix} a & a \\ b & b \\ c & c \\ d & d \\ e & e \end{pmatrix} \middle| \, a,b,c,d,e \in Z^o \right\}.$$

$$V_1^4 = \{6 \times 3 \text{ matrices with entries from } Z^o\},$$
$$V_2^4 = \{ Z^o \times Z^o\}$$

and

$$V_3^4 = \{2 \times 7 \text{ matrices with entries from } Z^o\}.$$

Define
$$\eta = (\eta_1 \cup \eta_2 \cup \eta_3 \cup \eta_4)$$
$$= \left(\eta_1^1, \eta_2^1, \eta_3^1\right) \cup \left(\eta_1^2, \eta_2^2\right) \cup \left(\eta_1^3, \eta_2^3, \eta_3^3, \eta_4^3\right) \cup \left(\eta_1^4, \eta_2^4, \eta_3^4\right):$$
$$V = (V_1 \cup V_2 \cup V_3 \cup V_4)$$
$$= \left(V_1^1, V_2^1, V_3^1\right) \cup \left(V_1^2, V_2^2\right) \cup \left(V_1^3, V_2^3, V_3^3, V_4^3\right) \cup$$
$$\left(V_1^4, V_2^4, V_3^4, V_4^4\right) \to [0, 1]$$

such that $\eta_{j_i}^i : V_{j_i}^i \to [0,1]$; $1 \leq j_i \leq n_i$; $i = 1, 2, 3, 4$. Now

$\eta_1^1 : V_1^1 \to [0,1]$ such that

$$\eta_1^1 \, (x \; y \; z) = \begin{cases} \dfrac{1}{x+y+z} & \text{if } x+y+z \neq 0 \\ 1 & \text{if } x+y+z = 0 \end{cases}$$

$\eta_2^1 : V_2^1 \to [0,1]$ defined by

$$\eta_2^1 \begin{bmatrix} a \\ b \\ c \\ d \end{bmatrix} = \begin{cases} \dfrac{1}{2a+b+c+3d} & \text{if } 2a+b+c+3d \neq 0 \\ 1 & \text{if } 2a+b+c+3d = 0 \end{cases}$$

$\eta_3^1 : V_3^1 \to [0,1]$ is given by

$$\eta_3^1 \begin{bmatrix} a & b & c \\ d & e & f \\ g & h & i \end{bmatrix} = \begin{cases} \dfrac{1}{c+e+g} & \text{if } c+e+g \neq 0 \\ 1 & \text{if } c+e+g = 0 \end{cases}$$

$\eta_2^2 : V_2^2 \to [0,1]$ defined by

$$\eta_2^2 \begin{bmatrix} a & b \\ c & d \\ e & f \\ g & h \\ i & j \end{bmatrix} = \begin{cases} \dfrac{1}{a+d+e+h+i} & \text{if } a+d+e+h+i \neq 0 \\ 1 & \text{if } a+d+e+h+i = 0 \end{cases}$$

$\eta_1^2 : V_1^2 \to [0,1]$ such that

$$\eta_1^2 (x\ y\ z\ \omega\ u) =$$
$$\begin{cases} \dfrac{1}{x+2y+3z+4\omega+5u} & \text{if } x+2y+3z+4\omega+5u \neq 0 \\ 1 & \text{if } x+2y+3z+4\omega+5u = 0 \end{cases}$$

$\eta_1^3 : V_1^3 \to [0,1]$ such that

$$\eta_1^3 \begin{bmatrix} a & b & c & d & e & f \\ g & h & i & j & k & l \end{bmatrix} = \begin{cases} \dfrac{1}{a+1} & \text{if } a+1 \neq 0 \\ 1 & \text{if } a+1 = 0 \end{cases}$$

$\eta_2^3 : V_2^3 \to [0, 1]$ defined by

$$\eta_2^3 (x\ y\ z\ \omega) = \begin{cases} \dfrac{1}{x} & \text{if } x \neq 0 \\ 1 & \text{if } x = 0 \end{cases}$$

$\eta_3^3 : V_3^3 \to [0, 1]$ is such that

$$\eta_3^3(p(x)) = \begin{cases} \dfrac{1}{\deg p(x)} & \text{if } p(x) \neq \text{constant} \\ 1 & \text{if } p(x) \text{ is a constant} \end{cases}$$

$\eta_4^3 : V_4^3 \to [0, 1]$ is such that

$$\eta_4^3 \begin{bmatrix} a & a \\ b & b \\ c & c \\ d & d \\ e & e \end{bmatrix} = \begin{cases} \dfrac{1}{2a + 2d} & \text{if } a + d \neq 0 \\ 1 & \text{if } a + d = 0 \end{cases}$$

$\eta_1^4 : V_1^4 \to [0, 1]$ is defined by

$$\eta_1^4 \begin{bmatrix} a & b & c \\ d & e & f \\ g & h & i \\ j & k & l \\ m & n & o \\ p & q & r \end{bmatrix} = \begin{cases} \dfrac{1}{a+d+g+h+g+r} & \text{if } a+d+g+h+g+r \neq 0 \\ 1 & \text{if } a+d+g+h+g+r = 0 \end{cases}$$

$\eta_2^4 : V_2^4 \to [0, 1]$ is given by

$$\eta_2^4(x\, y) = \begin{cases} \dfrac{1}{x+y} & \text{if } x+y \neq 0 \\ 1 & \text{if } x+y = 0 \end{cases}$$

$\eta_3^4 : V_3^4 \to [0,1]$ is given by

$$\eta_3^4 \begin{pmatrix} a & b & c & d & e & f & g \\ h & i & j & k & l & m & n \end{pmatrix} = \begin{cases} \dfrac{1}{a+n} & \text{if } a+n \neq 0 \\ 1 & \text{if } a+n = 0 \end{cases}$$

Clearly $\eta = (\eta_1 \cup \eta_2 \cup \eta_3 \cup \eta_4) : V \to [0, 1]$ is such that $V_\eta = \left(V_{1\eta_1} \cup V_{2\eta_2} \cup V_{3\eta_3} \cup V_{4\eta_4}\right)$ is a special semigroup set fuzzy vector 4 space.

Now we proceed on to define the notion of special semigroup set vector n-space over a set S.

**DEFINITION 4.2.29:** *Let*
$$V = (V_1, V_2, \ldots, V_n)$$
$$= \left(V_1^1, V_2^1, \ldots, V_{n_1}^1\right) \cup \left(V_1^2, V_2^2, \ldots, V_{n_2}^2\right) \cup \ldots \cup \left(V_1^n, V_2^n, \ldots, V_{n_n}^n\right)$$
*be a special semigroup set n-vector space over a set S. Let*
$$W = (W_1, W_2, \ldots, W_n)$$
$$= \left(W_1^1, W_2^1, \ldots, W_{n_1}^1\right) \cup \left(W_1^2, W_2^2, \ldots, W_{n_2}^2\right) \cup \ldots \cup \left(W_1^n, W_2^n, \ldots, W_{n_n}^n\right)$$
$$\subseteq (V_1 \cup V_2 \cup \ldots \cup V_n) = V$$
*be a special semigroup set vector n-subspace of V over the set S. Let*
$$\eta = (\eta_1 \cup \eta_2 \cup \ldots \cup \eta_n)$$
$$= \left(\eta_1^1, \eta_2^1, \ldots, \eta_{n_1}^1\right) \cup \left(\eta_1^2, \eta_2^2, \ldots, \eta_{n_2}^2\right) \cup \ldots \cup \left(\eta_1^n, \eta_2^n, \ldots, \eta_{n_n}^n\right) :$$
$$W = (W_1, W_2, \ldots, W_n)$$
$$= \left(W_1^1, W_2^1, \ldots, W_{n_1}^1\right) \cup \left(W_1^2, W_2^2, \ldots, W_{n_2}^2\right) \cup \ldots \cup \left(W_1^n, W_2^n, \ldots, W_{n_n}^n\right)$$
$$\to [0, 1]$$
*defined by* $\eta_{j_i}^i : V_{j_i}^i \to [0,1]$; $1 \leq j_i \leq n_i$; $i = 1, 2, \ldots, n$. *If*
$$W_\eta = \left(W_1 \cup W_2 \cup \ldots \cup W_n\right)_{(\eta_1 \cup \eta_2 \cup \ldots \cup \eta_n)} = \left(W_{1\eta_1} \cup W_{2\eta_2} \cup \ldots \cup W_{n\eta_n}\right)$$
*is a special semigroup set fuzzy vector space then we call $W\eta$ to be a special semigroup set fuzzy vector n-subspace.*

We illustrate this situation by a simple example.

***Example 4.2.35:*** Let $V = (V_1 \cup V_2 \cup V_3 \cup V_4 \cup V_5)$ be a special semigroup set vector 5 space over the set $Z_{20}$ where

$$V_1^1 = \left(V_1^1, V_2^1, V_3^1\right), V_2 = \left(V_1^2, V_2^2\right), V_3 = \left(V_1^3, V_2^3, V_3^3, V_4^3\right),$$
$$V_4 = \left(V_1^4, V_2^4\right), V_5 = \left(V_1^5, V_2^5, V_3^5, V_4^5, V_5^5\right)$$

is such that $V_1^1 = Z_{20} \times Z_{20}$,

$$V_2^1 = \left\{ \begin{pmatrix} a \\ b \\ c \end{pmatrix} \middle| a, b, c \in Z_{20} \right\},$$

$V_3^1 = \{$all $3 \times 3$ matrices with entries from $Z_{20}\}$, $V_1^2 = Z_{20} \times Z_{20} \times Z_{20} \times Z_{20}$, $V_2^2 = \{$all $4 \times 4$ matrices with entries from $Z_{20}\}$, $V_1^3 = Z_{20} \times Z_{20} \times Z_{20}$,

$$V_2^3 = \left\{ \begin{pmatrix} a \\ b \\ c \\ d \end{pmatrix} \middle| a, b, c, d \in Z_{20} \right\},$$

$$V_3^3 = \left\{ \begin{pmatrix} a & a & a & a & a \\ b & b & b & b & b \end{pmatrix} \middle| a, b \in Z_{20} \right\},$$

$V_4^3 = \{$all $2 \times 6$ matrices with entries from $Z_{20}\}$, $V_1^4 = Z_{20} \times Z_{20} \times Z_{20} \times Z_{20} \times Z_{20}$; $V_2^4 = \{$all $2 \times 2$ matrices with entries from $Z_{20}\}$, $V_1^5 = Z_{20} \times Z_{20} \times Z_{20}$, $V_2^5 = \{$all $5 \times 5$ upper triangular matrices with entries from $Z_{20}\}$,

$$V_3^5 = \left\{ \begin{pmatrix} a & b \\ a & b \\ a & b \\ a & b \\ a & b \\ a & b \end{pmatrix} \middle| a, b \in Z_{20} \right\},$$

$V_4^5 = \{Z_{20}[x]$ all polynomials in the variable x with coefficients from $Z_{20}\}$ and $V_5^5 = \{$all $4 \times 4$ symmetric matrices with entries from $Z_{20}\}$. Take

$$W = (W_1 \cup W_2 \cup W_3 \cup W_4 \cup W_5)$$
$$= \left(W_1^1, W_2^1, W_3^1\right) \cup \left(W_1^2, W_2^2\right) \cup \left(W_1^3, W_2^3, W_3^3, W_4^3\right) \cup$$
$$\left(W_1^4, W_2^4\right) \cup \left(W_1^5, W_2^5, W_3^5, W_4^5, W_5^5\right)$$
$$\subseteq (V_1 \cup V_2 \cup V_3 \cup V_4 \cup V_5)$$

is given below;

$$W_1^1 = Z_{20} \times \{0\} \subseteq V_1^1,$$

$$W_2^1 = \left\{ \begin{pmatrix} a \\ a \\ a \end{pmatrix} \middle| a \in Z_{20} \right\} \subseteq V_2^1,$$

$W_3^1 = \{$all $3 \times 3$ upper triangular matrices with entries from $Z_{20}\}$ $\subseteq V_3^1$, $W_1^2 = Z_{20} \times Z_{20} \times \{0\} \times \{0\} \subseteq V_1^2$, $W_2^2 = \{$all $4 \times 4$ upper triangular matrices with entries from $Z_{20}\} \subseteq V_2^2$,

$$W_1^3 = Z_{20} \times \{0\} \times \{0\} \subseteq V_1^3,$$

$$W_2^3 = \left\{ \begin{pmatrix} a \\ a \\ b \\ b \end{pmatrix} \middle| a, b \in Z_{20} \right\} \subseteq V_2^3,$$

$$W_3^3 = \left\{ \begin{pmatrix} a & a & a & a & a \\ b & b & b & b & b \end{pmatrix} \middle| a, b \in \{0, 2, 4, 6, 8, 10, 12, 14, 16, 18\} \right\}$$

$\subseteq V_3^3$, $W_4^3 = \{$all $2 \times 6$ matrices with entries from $\{0, 5, 10, 15\}$ $\subseteq Z_{20}\} \subseteq V_4^3$, $W_1^4 = Z_{20} \times Z_{20} \times Z_{20} \times \{0\} \times \{0\} \subseteq V_1^4$, $W_2^4 = \{$all $2 \times 2$ upper triangular matrices with entries from $Z_{20}\} \subseteq V_2^4$, $W_1^5 = Z_{20} \times Z_{20} \times \{0\} \subseteq V_1^5$, $W_2^5 = \{$all 5×5 diagonal matrices with entries from $Z_{20}\} \subseteq V_2^5$,

$$W_3^5 = \left\{ \begin{pmatrix} a & b \\ a & b \\ a & b \\ a & b \\ a & b \\ a & b \end{pmatrix} \middle| a, b \in \{0, 5, 10, 15\} \subseteq Z_{20} \right\} \subseteq V_3^5,$$

$W_4^5 = \{$all polynomials in x of degree less than or equal to 4 with coefficients from $Z_{20}\} \subseteq V_4^5$ and $W_5^5 = \{$all 4 ×4 symmetric matrices with entries from $Z_{20}\} \subseteq V_5^5$; is clearly a special semigroup set vector n-subspace of V over the set $Z_{20}$. Take $\eta$ : W → [0, 1] . Define

$$\eta = (\eta_1 \cup \eta_2 \cup \eta_3 \cup \eta_4 \cup \eta_5)$$
$$= \left(\eta_1^1, \eta_2^1, \eta_3^1\right) \cup \left(\eta_1^2, \eta_2^2\right) \cup \left(\eta_1^3, \eta_2^3, \eta_3^3, \eta_4^3\right)$$
$$\cup \left(\eta_1^4, \eta_2^4, \eta_3^4\right) \cup \left(\eta_1^5, \eta_2^5, \eta_3^5, \eta_4^5, \eta_5^5\right) :$$
$$W = (W_1 \cup W_2 \cup W_3 \cup W_4 \cup W_5)$$
$$= \left(W_1^1, W_2^1, W_3^1\right) \cup \left(W_1^2, W_2^2\right) \cup \left(W_1^3, W_2^3, W_3^3, W_4^3\right) \cup$$
$$\left(W_1^4, W_2^4\right) \cup \left(W_1^5, W_2^5, W_3^5, W_4^5, W_5^5\right) \to [0, 1]$$

as follows; $\eta_{j_i}^i : V_{j_i}^i \to [0, 1]$; $1 \le j_i \le n_i$; i = 1, 2, 3, 4, 5 is given by

$\eta_1^1 : W_1^1 \to [0,1]$ such that

$$\eta_1^1(a \ 0) = \begin{cases} \dfrac{1}{a} & \text{if } a \ne 0 \\ 1 & \text{if } a = 0 \end{cases}$$

$\eta_2^1 : W_2^1 \to [0, 1]$ is defined by

$$\eta_2^1 \begin{bmatrix} a \\ a \\ a \end{bmatrix} = \begin{cases} \dfrac{1}{3a} & \text{if } a \neq 0 \\ 1 & \text{if } a = 0 \end{cases}$$

$\eta_3^1 : W_3^1 \to [0,1]$ is given by

$$\eta_3^1 \begin{bmatrix} a & b & c \\ 0 & d & e \\ 0 & 0 & f \end{bmatrix} = \begin{cases} \dfrac{1}{c+e+f} & \text{if } c+e+f \neq 0 \\ 1 & \text{if } c+e+f = 0 \end{cases}$$

$\eta_1^2 : W_1^2 \to [0, 1]$ such that

$$\eta_1^2 (a\ b\ 0\ 0) = \begin{cases} \dfrac{1}{ab} & \text{if } ab \neq 0 \\ 1 & \text{if } ab = 0 \end{cases}$$

$\eta_2^2 : W_2^2 \to [0,1]$ is defined by

$$\eta_2^2 \begin{pmatrix} a & b & c & d \\ 0 & e & f & g \\ 0 & 0 & h & i \\ 0 & 0 & o & k \end{pmatrix} = \begin{cases} \dfrac{1}{ab+cf+hi} & \text{if } ab+cf+hi \neq 0 \\ 1 & \text{if } ab+cf+hi = 0 \end{cases}$$

$\eta_1^3 : W_1^3 \to [0,1]$ is such that

$$\eta_1^3 (a\ 0\ 0) = \begin{cases} \dfrac{1}{7a} & \text{if } a \neq 0 \\ 1 & \text{if } a = 0 \end{cases}$$

$\eta_2^3 : W_2^3 \to [0,1]$ is given by

$$\eta_2^3 \begin{pmatrix} a \\ a \\ b \\ b \end{pmatrix} = \begin{cases} \dfrac{1}{ab + a^2 + b^2} & \text{if } ab + a^2 + b^2 \neq 0 \\ 1 & \text{if } ab + a^2 + b^2 = 0 \end{cases}$$

$\eta_3^3 : W_3^3 \to [0,1]$ is such that

$$\eta_3^3 \begin{pmatrix} a & a & a & a & a \\ b & b & b & b & b \end{pmatrix} = \begin{cases} \dfrac{1}{5ab} & \text{if } ab \neq 0 \\ 1 & \text{if } ab = 0 \end{cases}$$

$\eta_4^3 : W_4^3 \to [0,1]$ is such that

$$\eta_4^3 \begin{bmatrix} a & b & c & d & e & f \\ g & h & i & j & k & l \end{bmatrix} = \begin{cases} \dfrac{1}{ag + fl} & \text{if } ag + fl \neq 0 \\ 1 & \text{if } ag + fl = 0 \end{cases}$$

$\eta_1^4 : W_1^4 \to [0,1]$ is defined by

$$\eta_1^4 (a\ b\ c\ 0\ 0) = \begin{cases} \dfrac{1}{abc} & \text{if } abc \neq 0 \\ 1 & \text{if } abc = 0 \end{cases}$$

$\eta_2^4 : W_2^4 \to [0,1]$ is given by

$$\eta_2^4 \begin{pmatrix} a & b \\ 0 & c \end{pmatrix} = \begin{cases} \dfrac{1}{ac + b} & \text{if } ac + b \neq 0 \\ 1 & \text{if } ac + b = 0 \end{cases}$$

$\eta_1^5 : W_1^5 \to [0, 1]$ is such that

$$\eta_1^5 (a\ b\ 0) = \begin{cases} \dfrac{1}{3a + 5b} & \text{if } 3a + 5b \neq 0 \\ 1 & \text{if } 3a + 5b = 0 \end{cases}$$

$\eta_2^5 : W_2^5 \to [0, 1]$ is defined by

$$\eta_2^5 \begin{bmatrix} a & 0 & 0 & 0 & 0 \\ 0 & b & 0 & 0 & 0 \\ 0 & 0 & c & 0 & 0 \\ 0 & 0 & 0 & d & 0 \\ 0 & 0 & 0 & 0 & e \end{bmatrix} = \begin{cases} \dfrac{1}{abc + de} & \text{if } abc + de \neq 0 \\ 1 & \text{if } abc + de = 0 \end{cases}$$

$\eta_3^5 : W_3^5 \to [0, 1]$ is given by

$$\eta_3^5 \begin{bmatrix} a & b \\ a & b \\ a & b \\ a & b \\ a & b \\ a & b \\ a & b \end{bmatrix} = \begin{cases} \dfrac{1}{6a + 3b} & \text{if } 6a + 3b \neq 0 \\ 1 & \text{if } 6a + 3b = 0 \end{cases}$$

$\eta_4^5 : W_4^5 \to [0,1]$ is given by

$$\eta_4^5 (p(x)) = \begin{cases} \dfrac{1}{\deg p(x)} & \text{if } p(x) \neq \text{constant} \\ 1 & \text{if } p(x) \text{ is a constant} \end{cases}$$

$\eta_5^5 : W_5^5 \to [0,1]$ is given by

$$\eta_5^5 \begin{bmatrix} a & 0 & 0 & 0 \\ 0 & b & 0 & 0 \\ 0 & 0 & c & 0 \\ 0 & 0 & 0 & d \end{bmatrix} = \begin{cases} \dfrac{1}{ad + bc} & \text{if } ad + bc \neq 0 \\ 1 & \text{if } ad + bc = 0 \end{cases}$$

Clearly $W_\eta = \left( W_{1\eta_1} \cup W_{2\eta_2} \cup W_{3\eta_3} \cup W_{4\eta_4} \cup W_{5\eta_5} \right)$ is a special semigroup set fuzzy vector 5-subspace.

Now we proceed on to define the notion of special semigroup fuzzy linear n-algebra.

**DEFINITION 4.2.30:** *Let*
$$V = (V_1, V_2, \ldots, V_n)$$
$$= \left(V_1^1, V_2^1, \ldots, V_{n_1}^1\right) \cup \left(V_1^2, V_2^2, \ldots, V_{n_2}^2\right) \cup \ldots \cup \left(V_1^n, V_2^n, \ldots, V_{n_n}^n\right)$$
*be a special semigroup set linear n-algebra over the semigroup S. Suppose*
$$\eta = (\eta_1 \cup \eta_2 \cup \ldots \cup \eta_n)$$
$$= \left(\eta_1^1, \eta_2^1, \ldots, \eta_{n_1}^1\right) \cup \left(\eta_1^2, \eta_2^2, \ldots, \eta_{n_2}^2\right) \cup \ldots \cup \left(\eta_1^n, \eta_2^n, \ldots, \eta_{n_n}^n\right):$$
$$W = (W_1, W_2, \ldots, W_n)$$
$$= \left(W_1^1, W_2^1, \ldots, W_{n_1}^1\right) \cup \left(W_1^2, W_2^2, \ldots, W_{n_2}^2\right) \cup \ldots \cup \left(W_1^n, W_2^n, \ldots, W_{n_n}^n\right)$$
$$\subseteq (V_1 \cup V_2 \cup \ldots \cup V_n) = V$$
*be a special semigroup set vector n-subspace of V over the set S. Let*
$$\eta = (\eta_1 \cup \eta_2 \cup \ldots \cup \eta_n)$$
$$= \left(\eta_1^1, \eta_2^1, \ldots, \eta_{n_1}^1\right) \cup \left(\eta_1^2, \eta_2^2, \ldots, \eta_{n_2}^2\right) \cup \ldots \cup \left(\eta_1^n, \eta_2^n, \ldots, \eta_{n_n}^n\right):$$
$$V = (V_1, V_2, \ldots, V_n)$$
$$= \left(V_1^1, V_2^1, \ldots, V_{n_1}^1\right) \cup \left(V_1^2, V_2^2, \ldots, V_{n_2}^2\right) \cup \ldots \cup \left(V_1^n, V_2^n, \ldots, V_{n_2}^n\right)$$
$$\to [0, 1]$$
*where $\eta_i: V_i \to [0,1]$ is a such that $V_{i\eta_i}$ is a special semigroup set fuzzy linear algebra true for $i = 1, 2, \ldots, n$. and*
$$\eta_{j_i}^i : V_{j_i}^i \to [0, 1];\ 1 \le j_i \le n_i;\ i = 1, 2, \ldots, n.$$
*We call $V_\eta = \left(V_{1\eta_1}, V_{2\eta_2}, \ldots, V_{n\eta_n}\right)$ to be a special semigroup set fuzzy linear n-algebra.*

It is important at this juncture to mention that the notion of special semigroup set fuzzy vector n-space and special semigroup set fuzzy linear n-algebras are fuzzy equivalent.

However we shall illustrate by an example a special semigroup set fuzzy linear n-algebra.

*Example 4.2.36:* Let
$$V = (V_1, V_2, V_3, V_4)$$
$$= \left(V_1^1, V_2^1, V_3^1\right) \cup \left(V_1^2, V_2^2, V_3^2, V_4^2\right) \cup \left(V_1^3, V_2^3\right) \cup$$
$$\left(V_1^4, V_2^4, V_3^4, V_4^4, V_5^4\right)$$

be a special semigroup set linear 4-algebra over the semigroup S $= Z^o = Z^+ \cup \{0\}$ where $V_1^1 = Z^o \times Z^o \times Z^o \times Z^o$, $V_2^1 = \{Z^o[x]$ all polynomials in the variable x of degree less than or equal to 5$\}$, $V_3^1 = \{$all $7 \times 2$ matrices with entries from $Z^o\}$, $V_1^2 = Z^o \times Z^o \times Z^o \times Z^o \times Z^o \times Z^o$, $V_2^2 = \{$all upper triangular $5 \times 5$ matrices with entries from $Z^o\}$,

$$V_3^2 = \left\{ \begin{pmatrix} a & a & a & a & a & a \\ b & b & b & b & b & b \\ c & c & c & c & c & c \end{pmatrix} \middle| a, b, c \in Z^o \right\},$$

$V_4^2 = \{$all $2 \times 6$ matrices with entries from $Z^o\}$,

$$V_1^3 = \left\{ \begin{pmatrix} a & a & a \\ b & b & b \\ c & c & c \\ d & d & d \\ e & e & e \\ f & f & f \end{pmatrix} \middle| a, b, c, d, e, f \in Z^o \right\},$$

$$V_2^3 = Z^o \times Z^o \times Z^o,$$

$V_1^4 = \{Z^o[x]$ all polynomials in the variable x with coefficient from $Z^o\}$,

$V_2^4 = \{$all $3 \times 3$ lower triangular matrices with entries from $Z^o\}$,

$V_3^4 = \{3 \times 7$ matrices with entries from $Z^o\}$,

$V_4^4 = \{7 \times 2$ matrices with entries from $Z^o\}$ and

$$V_5^4 = Z^o \times Z^o \times Z^o \times Z^o.$$

Define
$$\eta = (\eta_1, \eta_2, \eta_3, \eta_4, \eta_5)$$
$$= \left(\eta_1^1, \eta_2^1, \eta_3^1\right) \cup \left(\eta_1^2, \eta_2^2, \eta_3^2, \eta_4^2\right) \cup \left(\eta_1^3, \eta_2^3\right) \cup$$
$$\left(\eta_1^4, \eta_2^4, \eta_3^4, \eta_4^4, \eta_5^4\right) : V = (V_1, V_2, V_3, V_4)$$
$$= \left(V_1^1, V_2^1, V_3^1\right) \cup \left(V_1^2, V_2^2, V_3^2, V_4^2\right) \cup \left(V_1^3, V_2^3\right) \cup$$
$$\left(V_1^4, V_2^4, V_3^4, V_4^4, V_5^4\right) \to [0, 1]$$

such that, $\eta_i : V_i \to [0, 1]$, $1 \leq i \leq 4$;
$$\eta_{j_i}^i : V_{j_i}^i \to [0, 1]; \ 1 \leq j_i \leq n_i \ 1 \leq i \leq 4.$$

Now
$\eta_1^1 : V_1^1 \to [0, 1]$ such that

$$\eta_1^1 \ (a \ b \ c \ d) = \begin{cases} \dfrac{1}{abcd} & \text{if } abcd \neq 0 \\ 1 & \text{if } abcd = 0 \end{cases}$$

$\eta_2^1 : V_2^1 \to [0,1]$ defined by

$$\eta_2^1 \ (p(x)) = \begin{cases} \dfrac{1}{\deg p(x)} & \text{if } p(x) \neq \text{constant} \\ 1 & \text{if } p(x) \text{ is a constant} \end{cases}$$

$\eta_3^1 : V_3^1 \to [0,1]$ is given by

$$\eta_3^1 \begin{bmatrix} a & b \\ c & d \\ e & f \\ g & h \\ i & j \\ k & l \\ m & n \end{bmatrix} = \begin{cases} \dfrac{1}{abcdef + ifklmn} & \text{if } abcdef + ifklmn \neq 0 \\ 1 & \text{if } abcdef + ifklmn = 0 \end{cases}$$

$\eta_1^2 : V_1^2 \to [0, 1]$ such that

$$\eta_1^2 \, (a\ b\ c\ d\ e\ f) = \begin{cases} \dfrac{abc}{def} < 1 & \text{if } def \neq 0 \\ \dfrac{def}{abc} < 1 & \text{if } abc \neq 0 \\ 1 & \text{if } abc = 0 \text{ or } def = 0 \text{ or } abc + def = 0 \end{cases}$$

$\eta_2^2 : V_2^2 \to [0, 1]$ given by

$$\eta_2^2 \begin{bmatrix} a & b & c & d & e \\ 0 & e & f & g & h \\ 0 & 0 & h & i & j \\ 0 & 0 & 0 & j & k \\ 0 & 0 & 0 & 0 & 1 \end{bmatrix} = \begin{cases} \dfrac{1}{aehjl} & \text{if } aehjl \neq 0 \\ 1 & \text{if } aehjl = 0 \end{cases}$$

$\eta_3^2 : V_3^2 \to [0, 1]$ defined by

$$\eta_3^2 \begin{bmatrix} a & a & a & a & a & a \\ b & b & b & b & b & b \\ c & c & c & c & c & c \end{bmatrix} = \begin{cases} \dfrac{1}{abc} & \text{if } abc \neq 0 \\ 1 & \text{if } abc = 0 \end{cases}$$

$\eta_4^2 : V_4^2 \to [0,1]$ is given by

$$\eta_4^2 \begin{bmatrix} a & b & c & d & e & f \\ g & h & i & j & k & l \end{bmatrix} = \begin{cases} \dfrac{1}{agbh} & \text{if } agbh \neq 0 \\ 1 & \text{if } agbh = 0 \end{cases}$$

$\eta_1^3 : V_1^3 \to [0,1]$ defined that

$$\eta_1^3 \begin{pmatrix} a & a & a \\ b & b & b \\ c & c & c \\ d & d & d \\ e & e & e \\ f & f & f \end{pmatrix} = \begin{cases} \dfrac{1}{abcdef} & \text{if } abcdef \neq 0 \\ 1 & \text{if } abcdef = 0 \end{cases}$$

$\eta_2^3 : V_2^3 \to [0, 1]$ is defined by

$$\eta_2^3 (a\ b\ c) = \begin{cases} \dfrac{1}{abc} & \text{if } abc \neq 0 \\ 1 & \text{if } abc = 0 \end{cases}$$

$\eta_1^4 : V_1^4 \to [0, 1]$ is given by

$$\eta_1^4 (p(x)) = \begin{cases} \dfrac{1}{\deg p(x)} & \text{if } p(x) \neq \text{constant} \\ \dfrac{1}{k} & \text{if } k \neq 0 \text{ and } p(x) = k \\ 1 & \text{if } k = 0 \end{cases}$$

$\eta_2^4 : V_2^4 \to [0, 1]$ is such that

$$\eta_2^4 \begin{bmatrix} a & 0 & 0 \\ b & c & 0 \\ d & e & f \end{bmatrix} = \begin{cases} \dfrac{1}{abde} & \text{if } abde \neq 0 \\ 1 & \text{if } abde = 0 \end{cases}$$

$\eta_3^4 : V_3^4 \to [0, 1]$ is defined by

$$\eta_3^4 \begin{pmatrix} a & b & c & d & e & f & g \\ h & i & j & k & l & m & n \\ p & q & r & s & t & u & v \end{pmatrix} =$$

$$\begin{cases} \dfrac{1}{ajv} & \text{if } ajv \neq 0 \\ \dfrac{1}{bku} & \text{if } ajv = 0 \text{ and } bku \neq 0 \\ \dfrac{1}{gmt} & \text{if } ajv = 0,\ bku = 0,\ gmt \neq 0 \\ 1 & \text{if } ajv = 0,\ bku = 0,\ gmt = 0 \end{cases}$$

$\eta_4^4 : V_4^4 \to [0, 1]$ is given by

$$\eta_4^4 \begin{pmatrix} a & h \\ b & i \\ c & j \\ d & k \\ e & l \\ f & m \\ g & n \end{pmatrix} = \begin{cases} \dfrac{1}{ah + gn + dk} & \text{if } ah + gn + dk \neq 0 \\ 1 & \text{if } ah + gn + dk = 0 \end{cases}$$

$\eta_5^4 : V_5^4 \to [0,1]$ is defined by

$$\eta_5^4 (a\ b\ c\ d) = \begin{cases} \dfrac{1}{abcd} & \text{if } abcd \neq 0 \\ 1 & \text{if } abcd = 0 \end{cases}$$

Thus $V_\eta = \left( V_{1\eta_1}, V_{2\eta_2}, V_{3\eta_3}, V_{4\eta_4} \right)$ is a special semigroup set fuzzy 4-linear algebra. The notion of special semigroup set fuzzy vector n-space and special semigroup set fuzzy linear n-algebra are fuzzy equivalent. Like wise the notion of special semigroup set vector n-subspace and special semigroup set linear n-subalgebras are fuzzy equivalent. Now we leave it for the reader as an exercise to define these notions and construct examples.

Now we proceed on to define the notion of special group set vector space defined over a set S.

**DEFINITION 4.2.31:** *Let $V = (V_1 \cup V_2 \cup ... \cup V_n)$*
$= \left( V_1^1, V_2^1, ..., V_{n_1}^1 \right) \cup \left( V_1^2, V_2^2, ..., V_{n_2}^2 \right) \cup ... \cup \left( V_1^n, V_2^n, ..., V_{n_n}^n \right)$
*be a special semigroup set n-vector space over the set S. Let*
$\eta = (\eta_1 \cup \eta_2 \cup ... \cup \eta_n) = \left( \eta_1^1, \eta_2^1, ..., \eta_{n_1}^1 \right) \cup \left( \eta_1^2, \eta_2^2, ..., \eta_{n_2}^2 \right) \cup ...$
$\cup \left( \eta_1^n, \eta_2^n, ..., \eta_{n_n}^n \right): V = (V_1, V_2, ..., V_n) = \left( V_1^1, V_2^1, ..., V_{n_1}^1 \right) \cup \left( V_1^2, V_2^2, ..., V_{n_2}^2 \right) \cup ... \cup \left( V_1^n, V_2^n, ..., V_{n_2}^n \right) \to [0,1]$ *is such that*

$\eta_i : V_i \to [0, 1]$ where $V_{i\eta_i}$ is a special group set fuzzy vector space, $i = 1, 2, \ldots, n$; here $\eta^i_{j_i} : V^i_{j_i} \to [0,1]$ with $1 \le j_i \le n_i$; $i = 1, 2, \ldots n$. We call $V_\eta = (V_{1\eta_1}, V_{2\eta_2}, \ldots, V_{n\eta_n})$ to be a special group set fuzzy vector n-space.

We shall illustrate this by an example.

***Example 4.2.37:*** Let
$$V = (V_1 \cup V_2 \cup V_3 \cup V_4 \cup V_5)$$
$$= (V^1_1, V^1_2) \cup (V^2_1, V^2_2, V^2_3) \cup$$
$$(V^3_1, V^3_2, V^3_3, V^3_4) \cup (V^4_1, V^4_2) \cup (V^5_1, V^5_2)$$

be a special group set vector 5-space over the set $Z_{15}$. Here $V^1_1 = Z_{15} \times Z_{15}$, $V^1_2 = \{$all $4 \times 2$ matrices with entries from $Z_{15}\}$, $V^2_1 = Z_{15} \times Z_{15} \times Z_{15} \times Z_{15}$, $V^2_2 = \{Z_{15}[x]$ all polynomials in x of degree less than or equal to 15 with coefficients from $Z_{15}\}$, $V^2_3 = \{$all $4 \times 4$ matrices with entries from $Z_{15}\}$, $V^3_1 = \{$all $5 \times 5$ lower triangular matrices with entries from $Z_{15}\}$

$$V^3_2 = Z_{15} \times Z_{15} \times Z_{15} \times Z_{15} \times Z_{15},$$

$$V^3_3 = \left\{ \begin{pmatrix} a & b \\ a & b \\ a & b \\ a & b \\ a & b \\ a & b \\ a & b \end{pmatrix} \middle| a, b \in Z_{15} \right\},$$

$$V^3_4 = \left\{ \begin{pmatrix} a & b & c & d \\ e & f & g & h \\ i & j & k & l \end{pmatrix} \middle| a,b,c,d,e,f,g,h,i,j,k,l \in Z_{15} \right\},$$

$V_1^4 = Z_{15} \times Z_{15} \times Z_{15}$, $V_2^4 = $ {all $4 \times 4$ diagonal matrices with entries from $Z_{15}$} $V_1^5 = $ {all $3 \times 5$ matrices with entries from $Z_{15}$} and

$$V_2^5 = \left\{ \begin{pmatrix} a \\ a \\ a \\ a \\ b \\ b \end{pmatrix} \middle| a, b \in Z_{15} \right\}.$$

Define
$$\eta = (\eta_1 \cup \eta_2 \cup \eta_3 \cup \eta_4 \cup \eta_5)$$
$$= \left(\eta_1^1, \eta_2^1, \eta_3^1\right) \cup \left(\eta_1^2, \eta_2^2\right) \cup \left(\eta_1^3, \eta_2^3, \eta_3^3, \eta_4^3\right) \cup \left(\eta_1^4, \eta_2^4, \eta_3^4\right) \cup$$
$$\left(\eta_1^5, \eta_2^5, \eta_3^5, \eta_4^5, \eta_5^5\right) : V = (V_1 \cup V_2 \cup V_3 \cup V_4 \cup V_5)$$
$$= \left(V_1^1, V_2^1\right) \cup \left(V_1^2, V_2^2, V_3^2\right) \cup \left(V_1^3, V_2^3, V_3^3, V_4^3\right)$$
$$\cup \left(V_1^4, V_2^4\right) \cup \left(V_1^5, V_2^5\right) \to [0,1];$$
$$\eta_i: V_i \to [0,1]; 1 \leq i \leq 5,$$
$$\eta_{j_i}^i : V_{j_i}^i \to [0,1]; 1 \leq j_i \leq n_i, 1 \leq i \leq 5.$$

Now
$\eta_1^1 : V_1^1 \to [0,1]$ defined by

$$\eta_1^1(a\ b) = \begin{cases} \dfrac{1}{ab} & \text{if } ab \neq 0 \\ 1 & \text{if } ab = 0 \end{cases}$$

$\eta_2^1 : V_2^1 \to [0, 1]$ is given by

$$\eta_2^1 \begin{bmatrix} a & e \\ b & f \\ c & g \\ d & h \end{bmatrix} = \begin{cases} \dfrac{1}{abcd + efgh} & \text{if } abcd + efgh \neq 0 \\ 1 & \text{if } abcd + efgh = 0 \end{cases}$$

$\eta_1^2 : V_1^2 \to [0,1]$ is defined by

$$\eta_1^2 (a\ b\ c\ d) = \begin{cases} \dfrac{1}{abc} & \text{if } abc \neq 0 \\ \dfrac{1}{d} & \text{if } abc = 0,\ d \neq 0 \\ 1 & \text{if } d = 0,\ abc = 0 \end{cases}$$

$\eta_2^2 : V_2^2 \to [0,1]$ such that

$$\eta_2^2 (p(x)) = \begin{cases} \dfrac{1}{p(x)} & \text{if } p(x) \neq k,\ k \text{ a constant} \\ \dfrac{1}{k} & \text{if } k \neq 0 \text{ and } p(x) = k \\ 1 & \text{if } k = 0 \end{cases}$$

$\eta_3^2 : V_3^2 \to [0,1]$ is defined by

$$\eta_3^2 \begin{bmatrix} a & b & c & d \\ e & f & g & h \\ i & j & k & l \\ m & n & p & q \end{bmatrix} = \begin{cases} \dfrac{1}{afkg} & \text{if } afkq \neq 0,\ dgjm = 0 \\ \dfrac{1}{dgjm} & \text{if } afkq = 0,\ dgjm \neq 0 \\ 1 & \text{if } afkq = 0,\ dgjm = 0 \end{cases}$$

$\eta_1^3 : V_1^3 \to [0, 1]$ is given by

$$\eta_1^3 \begin{pmatrix} a & 0 & 0 & 0 & 0 \\ b & c & 0 & 0 & 0 \\ d & e & f & 0 & 0 \\ g & h & i & j & 0 \\ k & l & m & n & p \end{pmatrix} = \begin{cases} \dfrac{1}{abdgk} & \text{if } abdgk \neq 0 \\ 1 & \text{if } abdgk = 0 \end{cases}$$

$\eta_2^3 : V_2^3 \to [0, 1]$ is defined by

$$\eta_2^3 \, (a \; b \; c \; d \; e) = \begin{cases} \dfrac{1}{ae} & \text{if } ae \neq 0 \\ 1 & \text{if } ae = 0 \end{cases}$$

$\eta_3^3 : V_3^3 \to [0, 1]$ is given by

$$\eta_3^3 \begin{bmatrix} a & b \\ a & b \\ a & b \\ a & b \\ a & b \\ a & b \\ a & b \end{bmatrix} = \begin{cases} \dfrac{1}{a+b} & \text{if } a+b \neq 0 \\ 1 & \text{if } a+b = 0 \end{cases}$$

$\eta_4^3 : V_4^3 \to [0, 1]$ is such that

$$\eta_4^3 \begin{bmatrix} a & b & c & d \\ e & f & g & h \\ i & j & k & l \end{bmatrix} = \begin{cases} \dfrac{1}{a} & \text{if } a \neq 0 \\ \dfrac{1}{b} & \text{if } b \neq 0 \\ 1 & \text{if } a = 0 \text{ and } b = 0 \end{cases}$$

$\eta_1^4 : V_1^4 \to [0,1]$ is given by

$$\eta_1^4 \, (a \; b \; c) = \begin{cases} \dfrac{1}{ab + bc + ca} & \text{if } ab + bc + ca \neq 0 \\ 1 & \text{if } ab + bc + ca = 0 \end{cases}$$

$\eta_2^4 : V_2^4 \to [0,1]$ is such that

$$\eta_2^4 \begin{bmatrix} a & 0 & 0 & 0 \\ 0 & b & 0 & 0 \\ 0 & 0 & c & 0 \\ 0 & 0 & 0 & d \end{bmatrix} = \begin{cases} \dfrac{1}{ad + bc} & \text{if } ad + bc \neq 0 \\ 1 & \text{if } ad + bc = 0 \end{cases}$$

$\eta_1^5 : V_1^5 \to [0,1]$ is defined by

$$\eta_1^5 \begin{pmatrix} a & b & c & d & e \\ g & h & i & j & k \\ l & m & n & p & q \end{pmatrix} =$$

$$\begin{cases} \dfrac{1}{ab+hi+pq} & \text{if } ab+hi+pq \neq 0 \\ \dfrac{1}{cd+gh+mn} & \text{if } ab+hi+pq = 0 \text{ and } cd+gh+mn \neq 0 \\ 1 & \text{if } ab+hi+pq = 0 \text{ and } cd+gh+mn = 0 \end{cases}$$

$\eta_2^5 : V_2^5 \to [0,1]$ is such that

$$\eta_2^5 \begin{pmatrix} a \\ a \\ a \\ a \\ b \\ b \end{pmatrix} = \begin{cases} \dfrac{1}{4a+2b} & \text{if } 4a+2b \neq 0 \\ 1 & \text{if } 4a+2b = 0 \end{cases}$$

Thus $V_\eta = \left( V_{1\eta_1}, V_{2\eta_2}, V_{3\eta_3}, V_{4\eta_4}, V_{5\eta_5} \right)$ is a special group set fuzzy 5-vector space. One of the major advantages of fuzzifying these concepts is many a times these notions become fuzzy equivalent, there by overcoming some disadvantages in practical application. These new algebraic structures find their application in computer engineering, web testing / processing, in fuzzy models, industries, coding theory and cryptology.

Chapter Five

# SUGGESTED PROBLEMS

This chapter suggests 66 problems so that the reader becomes familiar with these new concepts

1. Let $V = (V_1, V_2, V_3, V_4)$ where $V_1 = Z^+[x]$, $V_2 = Z^+ \times Z^+ \times Z^+$

$$V_3 = \left\{ \begin{pmatrix} a & b \\ c & d \end{pmatrix} \middle| a,b,c,d \in Z^+ \right\}$$

and

$$V_4 = \left\{ \begin{pmatrix} a & a & a & a \\ a & a & a & a \end{pmatrix} \middle| a \in Z^+ \right\}$$

be a special semigroup set vector space over the set {2, 4, 6, …}. Find a generating 4-subset of V. Is V finite or infinite? Find a proper special semigroup set vector subspace of V.

2. Obtain some interesting properties about special semigroup set vector spaces.

3. Find the applications of the special semigroup set vector spaces to industries.

4. Let $V = (V_1, V_2, V_3, V_4)$ where $V_1 = Z_8$, $V_2 = Z_6$, $V_3 = Z_9$ and $V_4 = Z_{12}$, clearly $V_1$, $V_2$, $V_3$ and $V_4$ are semigroups under modulo addition. V is a special semigroup set vector space over the set $S = \{0, 1, 3, 5\}$.
   a. Find the n-generating subset of V.
   b. Find proper 4-subset of V which is a special semigroup set vector subspace of V.
   c. Is V a special semigroup set linear algebra over S? Justify your claim.

5. Let $V = (V_1, V_2, V_3, V_4, V_5) = (Z_5, Z_7, Z_{23}, Z_{19}, Z_{11})$ be a special semigroup set linear algebra over the additive semigroup $= \{0, 1\}$. Is V a doubly simple special semigroup set linear algebra or is V a special semigroup set strongly simple linear algebra? Give the generating 5 subset of V over S. Is V finite 5-dimensional?

6. Let $V = (V_1, V_2, V_3, V_4) = (Z_6, Z_7, Z_8, Z_9)$ be a special semigroup set vector space over the set $S = \{0, 2, 3, 4\}$, Find a 4-generating subset of V. What is 4-dimension of V as a special semigroup set vector space over the set? How many special semigroup set vector subspaces V has? Can V have more than one 4-generating subset? Justify your claim.

7. Let
$$V = \left( \left\{ \begin{pmatrix} a & b \\ c & a \end{pmatrix} \middle| a \in Z_6 \right\}, (Z_6 \times Z_6 \times Z_6), \right.$$

$$\left. \left\{ \begin{pmatrix} a & a & a & a & a \\ a & a & a & a & a \end{pmatrix} \middle| a \in Z_6 \right\}, \left\{ \begin{pmatrix} a & a \\ a & a \\ a & a \\ a & a \\ a & a \end{pmatrix} \middle| a \in Z_6 \right\}, \right.$$

{$Z_6$ [x]; all polynomials of degree less than or equal to 5})
($V_1$, $V_2$, $V_3$, $V_4$, $V_5$) be a special semigroup set linear algebra over the semigroup $S = Z_6$.

a. Does V have special semigroup set linear subalgebras over S?
b. Is V a special semigroup set simple linear algebra over $Z_6$?
c. Can V have special subsemigroup set linear subalgebra?
d. Find a n-generating subset of V?
e. Can V have more than one n-generating subset?
f. What is the n-dimension of V?
g. Is V a doubly strong special semigroup set linear algebra over $S = Z_6$?
h. Is V a special semigroup set strongly simple linear algebra?

8. Let V = ($V_1$, $V_2$, $V_3$, $V_4$, $V_5$) where

$$V_1 = \left\{ \begin{pmatrix} a & b \\ c & d \end{pmatrix} \middle| a,b,c,d \in Z_7 \right\}, V_2 = \{Z_7 \times Z_7 \times Z_7 \times Z_7\},$$

$V_3 = \{Z_7$ [x]; all polynomials of degree less than or equal to three with coefficients from $Z_7\}$,

$$V_4 = \left\{ \begin{pmatrix} a & a \\ a & a \\ a & a \\ a & a \end{pmatrix} \middle| a \in Z_7 \right\}$$

and

$$V_5 = \left\{ \begin{pmatrix} a & 0 & b \\ 0 & c & 0 \\ b & 0 & d \end{pmatrix} \middle| a,b,c,d \in Z_7 \right\}$$

are semigroups under addition. Thus V = ($V_1$, $V_2$, ... , $V_5$) is a special semigroup set linear algebra over the semigroup $Z_7$.

a. Find a 5-generating 5-subset of V.

b. How many 5-generating 5-subsets of V can be found in V?
   c. Is V a special semigroup set simple linear algebra?
   d. Is V a special semigroup set strongly simple linear algebra?
   e. Is V a doubly simple special semigroup set linear algebra?
   f. Can V have special semigroup set linear subalgebras?

9. Let $V = (V_1, V_2, V_3, V_4, V_5)$,

$$\left( \left\{ \begin{pmatrix} a & b \\ c & d \end{pmatrix} \middle| a,b,c,d \in Z^+ \right\}, \left\{ \begin{pmatrix} a & a & a & a & a \\ a & a & a & a & a \\ a & a & a & a & a \end{pmatrix} \middle| a \in Z^+ \right\}, \right.$$

$(Z^+ \times Z^+ \times Z^+ \times Z^+ \times Z^+)$, $\{Z^+[x]$, all polynomials of degree less than or equal to $3\}$,

$$\left\{ \begin{pmatrix} a & a & a & a & a \\ a & b & b & b & b \\ a & b & b & b & b \\ a & b & b & b & b \\ a & b & b & b & b \end{pmatrix} \middle| a,b \in Z^+ \right\}$$

be a special semigroup set vector space defined over the set $Z^+$. Let $W = (W_1, W_2, W_3, W_4, W_5)$ where

$$W_1 = \left\{ \begin{pmatrix} a & a \\ a & a \end{pmatrix} \middle| a \in Z^+ \right\}, W_2 = \left\{ \begin{pmatrix} a & a & a \\ a & a & a \\ a & a & a \\ a & a & a \\ a & a & a \end{pmatrix} \middle| a \in Z^+ \right\},$$

$W_3 = \{Z^+ [x]$ all polynomials of degree less than or equal to 4 with coefficients from $Z^+\}$, $W_4 = (Z^+ \times Z^+ \times Z^+ \times Z^+)$ and

$$W_5 = \left\{ \begin{pmatrix} a & b & c & d & e \\ b & f & g & h & i \\ c & g & j & k & l \\ d & h & k & m & n \\ c & i & l & n & p \end{pmatrix} \middle| \text{ the entries are from } Z^+ \right\}$$

be a special semigroup set vector space over the same set $S = Z^+$. How many different types of special semigroup set linear transformations $(T_1, T_2, T_3, T_4, T_5) = T$ can be got from V to W. If $U = (U_1, U_2, \ldots, U_5)$ is a set linear transformation from W to V can we find any relation between U and V?

10. Let $V = (V_1, V_2, V_3)$ and $W = (W_1, W_2, W_3)$ where $V_1 = \{Z_{15} \times Z_{15} \times Z_{15}\}$,

$$V_2 = \left\{ \begin{pmatrix} a & b & e \\ c & d & f \end{pmatrix} \middle| a,b,c,d \in Z_{15} \right\}$$

and $V_3 = \{Z_{15} [x]$ all polynomials of degree less than 5$\}$ is a special semigroup set vector space over the set $Z_{15}$. $W = (W_1, W_2, W_3)$ where

$$W_1 = \left\{ \begin{pmatrix} a & b \\ 0 & c \end{pmatrix} \middle| a,b,c \in Z_{15} \right\},$$

$$W_2 = \left\{ \begin{pmatrix} a & a \\ a & a \\ a & a \end{pmatrix} \middle| a \in Z_{15} \right\}$$

and $W_3 = \{Z_{15} \times Z_{15} \times Z_{15} \times Z_{15} \times Z_{15} \times Z_{15}\}$ is again a special semigroup set vector space over the set $Z_{15}$. Find special semigroup linear transformation from V to W. How many such special semigroup linear transformations be made from V to W? Let $SH_{Z_{12}}(V, W)$ denote the collection of all such special semigroup linear transformations from V to W. What is the algebraic structure enjoyed by $SH_{Z_{12}}(V, W)$?

11. Let $V = (V_1, V_2, V_3, V_4)$ where

$$V_1 = \left\{ \begin{pmatrix} a & b \\ c & d \end{pmatrix} \middle| a,b,c,d \in Z^+ \right\},$$

$$V_2 = \left\{ \begin{pmatrix} a & 0 & 0 & 0 \\ b & e & 0 & 0 \\ c & f & g & 0 \\ d & q & p & h \end{pmatrix} \middle| a,b,c,d,e,f,g,p,g,h \in Z^+ \right\},$$

$V_3 = \{Z^+ \times Z^+ \times Z^+\}$ and $V_4 = \{$All polynomial of degree less than or equal to 2 with coefficients from $Z^+\}$. $V_1$, $V_2$, $V_3$ and $V_4$ are semigroups under addition. Thus V is a special semigroup set vector space over the set $S = Z^+$.

   a. Find $T = (T_1, T_2, T_3, T_4)$, a special semigroup set linear operator on V.
   b. If $SHom_s(V, V) = \{Hom_s(V_1, V_2), Hom_s(V_2, V_4), Hom_s(V_3, V_1), Hom_s(V_4, V_3)\}$. What is the algebraic structure of $SHom_s(V, V)$?
   c. If $CS(Hom_s(V, V)) = \{SHom_s(V, V)\}$, what can we say about the algebraic structure of $CS(Hom_s(V, V))$?

12. Let $V = (V_1, V_2, V_3, V_4, V_5)$ where $V_1 = Z_{12}$; $V_2 = Z_{14}$; $V_3 = Z_{16}$; $V_4 = Z_7$ and $V_5 = Z_{10}$ are semigroups under modulo addition. V is a special semigroup set vector space over the set $S = \{0, 1\}$. Find $SHom_s(V, V) = \{Hom_s(V_1, V_2), Hom_s(V_2, V_3), Hom_s(V_4, V_5), Hom_s(V_5, V_1), Hom_s(V_3, V_4)\}$. Is $Hom_s(V, V)$ a special semigroup set vector space over S =

{0, 1}? What is the structure of CS (Hom$_s$ (V, V)) = {SHom$_s$ (V, V)}?

13. Let V = (V$_1$, V$_2$, V$_3$, V$_4$, V$_5$) where

$$V_1 = \left\{ \begin{pmatrix} a & b & c \\ d & e & f \\ g & h & i \end{pmatrix} \middle| a,b,c,d,e,f,g,h,i \in Z^+ \right\},$$

V$_2$ = {Z$^+$ × Z$^+$ × Z$^+$ × Z$^+$ × Z$^+$}, V$_3$ = {Z$^+$ [x]; all polynomials of degree less than or equal to 5},

$$V_4 = \left\{ \begin{pmatrix} a & a & a & a & a \\ b & b & b & b & b \\ c & c & c & c & c \end{pmatrix} \middle| a,b,c \in Z^+ \right\}$$

and

$$V_5 = \left\{ \begin{pmatrix} a & 0 & 0 & 0 \\ b & e & 0 & 0 \\ c & f & g & 0 \\ d & h & i & j \end{pmatrix} \middle| a,b,c,d,e,f,g,h,i,j \in Z^+ \right\}.$$

V$_1$, V$_2$, V$_3$, V$_4$ and V$_5$ are semigroups under addition. V is a special semigroup set linear algebra over the semigroup Z$^+$. Let W = (W$_1$, W$_2$, W$_3$, W$_4$, W$_5$) where

$$W_1 = \left\{ \begin{pmatrix} a & a & a \\ b & b & b \\ b & b & b \\ c & c & c \\ c & c & c \end{pmatrix} \middle| a,b,c \in Z^+ \right\},$$

$$W_2 = \left\{ \begin{pmatrix} a & a & a \\ b & b & b \\ c & c & c \end{pmatrix} \middle| a,b,c \in Z^+ \right\},$$

$W_3 = \{2Z^+ \times Z^+ \times 3Z^+ \times 5Z^+ \times 7Z^+\}$, $W_4 = \{Z^+[x];$ all polynomials of degree less than or equal to 5$\}$
and

$$W_5 = \left\{ \begin{pmatrix} 0 & a & a & a \\ b & 0 & c & c \\ b & c & 0 & d \\ b & c & d & 0 \end{pmatrix} \middle| a,b,c,d \in Z^+ \right\},$$

$W = (W_1, W_2, W_3, W_4, W_5)$ is a special semigroup set linear algebra over the same semigroup $Z^+$. Define $T = (T_1, T_2, T_3, T_4, T_5)$ from V to W such that T is a special semigroup set linear transformation of V to W.

If $SHom_{Z^+}(V, V) = \{Hom_{Z^+}(V_1, W_2), Hom_{Z^+}(V_2, W_3), Hom_{Z^+}(V_3, W_4), Hom_{Z^+}(V_4, W_1), Hom_{Z^+}(V_5, W_5)\}$, prove $SHom_{Z^+}(V, V)$ is atleast a special semigroup set vector space over the set $Z^+$.

14. Obtain some interesting results about $SHom_s(V, W)$ and $SHom_s(V, V)$.

15. Let $V = V = (V_1, V_2, V_3, V_4)$ where $V_1 = \{Z_{10} \times Z_{10} \times Z_{10}\}$,

$$V_2 = \left\{ \begin{pmatrix} a & b & c \\ d & e & f \\ g & h & i \end{pmatrix} \middle| a,b,c,d,e,f,g,h,i \in Z_{10} \right\},$$

$V_3 = \{Z_{10}[x];$ set of all polynomials of degree less than or equal to 7$\}$ and

$$V_4 = \left\{ \begin{pmatrix} a & b & c & d \\ 0 & 0 & e & f \\ 0 & 0 & g & h \\ 0 & 0 & p & q \end{pmatrix} \,\middle|\, p,q,g,b,a,bc,de,f \in Z_{10} \right\}.$$

$V_1, V_2, V_3, V_4$ are semigroups under addition.

Define $T = (T_1, T_2, T_3, T_4)$ from V to V.

Suppose $SHom_{Z_{10}}(V, V) = \{ Hom_{Z_{10}}(V_1, V_1), Hom_{Z_{10}}(V_2, V_2), Hom_{Z_{10}}(V_3, V_3), Hom_{Z_{10}}(V_4, V_4) \}$. Is $SHom_{Z_{10}}(V, V)$ a special semigroup set linear algebra over the semigroup $Z_{10}$? What is the algebraic structure of $\{ SHom_{Z_{10}}(V, V) \}$?

16. Let $V = (V_1, V_2, V_3, V_4)$ be a special semigroup set vector space over the set S. Suppose $W = (W_1, W_2, W_3, W_4)$, a proper 4-subset of V and $W_1 \subset V_1$, $W_2 \subset V_2$, $W_3 \subset V_3$ and $W_4 \subset V_4$ are subsemigroups and $W = (W_1, W_2, W_3, W_4)$ is a special semigroup set vector subspace of V. Define $P = (P_1, P_2, P_3, P_4)$ from V into W such that each $P_i : V_i \to W_i$ ; $1 \leq i \leq 4$ is a projection. Prove P. P = P.

17. Let $V = (V_1, V_2, \ldots, V_n)$ be a special semigroup set linear algebra over a semigroup S. If $W = (W_1, W_2, \ldots, W_n)$ be a proper subset of V which is a special semigroup set linear subalgebra of V over S. If $P = (P_1, P_2, \ldots, P_n)$ is defined such that $P_i : V_i \to W_i$ where each $P_i$ is a projection of $V_i$ to $W_i$, $1 \leq i \leq n$. Show P is a special semigroup set linear operator on V and P . P = P.

18. Obtain some interesting results about the special semigroup set linear operators which are projection on V.

19. Let $V = (V_1, \ldots, V_n)$. If each $V_i = W_1^i \oplus \cdots \oplus W_{t_i}^i$ ; $1 \leq i \leq n$ then we write $V = (W_1^1 \oplus \cdots \oplus W_{t_1}^1, W_1^2 \oplus \cdots \oplus W_{t_2}^2, \ldots, W_1^n \oplus \cdots \oplus W_{t_n}^n)$ and call this as a special semigroup set direct sum of the special semigroup set linear algebra.

Using this definition, write $V = (V_1, V_2, V_3, V_4)$ where $V_1 = \{Z^+ \times Z^+ \times Z^+ \times Z^+ \times Z^+\}$

$$V_2 = \left\{ \begin{pmatrix} a & b \\ c & d \end{pmatrix} \middle| a, b, c, d \in Z^+ \right\},$$

$V_3 = \{(a\ a\ a\ a) \mid a \in Z^+\}$ and $V_4 = \{3Z^+ \times 2Z^+ \times 5Z^+\}$ as a direct sum of special semigroup set linear subalgebras.

20. Let $V = (V_1, V_2, V_3)$ where

$$V_1 = \left\{ \begin{pmatrix} a & b & c \\ d & e & f \\ g & h & i \end{pmatrix} \middle| a, b, c, d, e, f, g, h, i \in Z^+ \cup \{0\} \right\},$$

$$V_2 = \{Z^+ \cup \{0\} \times Z^+ \cup \{0\} \times Z^+ \cup \{0\}\}$$

and

$$V_3 = \left\{ \begin{pmatrix} a & 0 & 0 & 0 \\ b & c & 0 & 0 \\ d & e & f & 0 \\ g & h & i & j \end{pmatrix} \middle| a, b, c, d, e, f, g, h, i, j \in Z^+ \cup \{0\} \right\}$$

are semigroups under addition. $V = (V_1, V_2, V_3)$ is a special semigroup set linear algebra over the semigroup $Z^+ \cup \{0\}$.

$$V = \left\{ \left\langle \begin{pmatrix} a & b & c \\ 0 & 0 & 0 \\ 0 & 0 & 0 \end{pmatrix} \right\rangle \right\} \oplus \left\{ \left\langle \begin{pmatrix} 0 & 0 & 0 \\ d & e & f \\ 0 & 0 & 0 \end{pmatrix} \right\rangle \right\} \oplus \left\{ \left\langle \begin{pmatrix} 0 & 0 & 0 \\ 0 & 0 & 0 \\ g & h & i \end{pmatrix} \right\rangle \right\}$$

$$= W_1^1 \oplus W_2^1 \oplus W_3^1.$$

$$V_2 = \{\langle Z^+ \cup \{0\} \times \{\phi\} \times \{\phi\} \rangle\} \oplus \{\langle \{\phi\} \times Z^+ \cup \{0\} \times Z^+ \cup \{0\} \rangle\}$$

$$= W_1^2 \oplus W_2^2$$

and

$$V_3 = \left\{ \left\langle \begin{pmatrix} a & 0 & 0 & 0 \\ b & 0 & 0 & 0 \\ d & 0 & 0 & 0 \\ g & 0 & 0 & 0 \end{pmatrix} \right\rangle \right\} \oplus \left\{ \left\langle \begin{pmatrix} 0 & 0 & 0 & 0 \\ 0 & c & 0 & 0 \\ 0 & e & 0 & 0 \\ 0 & h & 0 & 0 \end{pmatrix} \right\rangle \right\}$$

$$\oplus \left\{ \left\langle \begin{pmatrix} 0 & 0 & 0 & 0 \\ 0 & 0 & 0 & 0 \\ 0 & 0 & f & 0 \\ 0 & 0 & i & j \end{pmatrix} \right\rangle \right\} = W_1^3 \oplus W_2^3 \oplus W_3^3$$

Now $W = \left( W_1^1, W_2^2, W_3^3 \right)$ is a special semigroup set linear subalgebra of V over $Z^+ \cup \{0\}$.

a. Define special semigroup set linear projection $P_1$ from V to W show $P_1 \circ P_1 = P_1$.

b. Let $R = \left( W_2^1, W_1^2, W_2^3 \right)$ be a special semigroup set linear subalgebra of V over $Z^+ \cup \{0\}$. Define a special semigroup set linear projection $P_2$ from V into R.

21. Let $V = (V_1, V_2, V_3, V_4)$ where

$$V_1 = \left\{ \begin{pmatrix} a & b & c \\ d & e & f \\ g & h & i \end{pmatrix} \middle| a,b,c,d,e,f,g,h,i \in Z^+ \cup \{0\} \right\},$$

$V_2 = \{Z^+ \cup \{0\} \times Z^+ \cup \{0\} \times Z^+ \cup \{0\} \times Z^+ \cup \{0\}\}$,

$$V_3 = \left\{ \begin{pmatrix} a & b & c & d & e & f \\ g & h & i & j & k & l \end{pmatrix} \middle| a,b,c,d,e,f,g,h,i,j,k,l \in Z^+ \cup \{0\} \right\}$$

and

$$V_4 = \left\{ \begin{pmatrix} a & a \\ a & a \\ a & a \\ a & a \\ a & a \end{pmatrix} \middle| a \in Z^+ \cup \{0\} \right\}$$

be a semigroup. Write V as a direct sum of special semigroup set linear subalgebras. Define projections using them.

22. Obtain some interesting results about the direct sum of special semigroup set linear algebras.

23. Let $V = (V_1, V_2, V_3, V_4)$ be a special semigroup set linear algebra over the semigroup $Z^o$ where $V_1 = \{Z^o \times Z^o \times Z^o \times Z^o \times Z^o\}$ where $Z^o = Z^+ \cup \{0\}$, $V_2 = \{Z^o[x] \mid$ all polynomials of degree less than or equal to 4 with coefficients from $Z^+ \cup \{0\} = Z^o\}$.

$$V_3 = \left\{ \begin{pmatrix} a & b & c \\ d & e & f \\ g & h & i \end{pmatrix} \middle| a,b,c,d,e,f,g,h,i \in Z^o = Z^+ \cup \{0\} \right\}$$

and

$$V_4 = \left\{ \begin{pmatrix} a & a & a \\ a & a & a \\ a & a & a \\ a & a & a \\ a & a & a \\ a & a & a \end{pmatrix} \middle| a \in Z^o \right\}.$$

Find special subsemigroup linear subalgebra over some proper subsemigroup of $Z^o$. Write V also as a direct sum. Define $SHom_{Z^o}(V, V) = \{Hom_{Z^o}(V_1, V_2), Hom_{Z^o}(V_2, V_1), Hom_{Z^o}(V_3, V_3), Hom_{Z^o}(V_4, V_4)\}$. Can $SHom_{Z^o}(V, V)$ be a special semigroup linear algebra over $Z^o$? Justify your claim.

24. Find some interesting properties of special set vector spaces.

25. Find a special set vector space of dimension (1, 2, 3, 4, 5).

26. Given $V = (V_1, V_2, V_3, V_4, V_5)$ where $V_1 = \{S^o \times S^o \times S^o \mid S^o = 3Z^+ \cup \{0\}\}$,

$$V_2 = \left\{ \begin{pmatrix} a & b \\ c & d \end{pmatrix} \middle| a,b,c,d \in 5Z^+ \cup \{0\} \right\},$$

$$V_3 = \left\{ \begin{pmatrix} a_1 & a_2 & a_3 & a_4 \\ a_5 & a_6 & a_7 & a_8 \end{pmatrix} \middle| a_i \in Z^+ \cup \{0\}; 1 \le i \le 8 \right\},$$

$V_4 = \{(a\ a\ a\ a\ a)\ (a\ a\ a\ a\ a\ a\ a) \mid a_i \in Z^+ \cup \{0\}\}$

and

$$V_5 = \left\{ \begin{pmatrix} a & a & a \\ a & a & a \\ a & a & a \end{pmatrix}, \begin{pmatrix} a & b & c & d \\ 0 & e & f & g \\ 0 & 0 & h & i \\ 0 & 0 & 0 & f \end{pmatrix} \middle| a,b,e,c,d,e,f,g,h,i,j \in Z^+ \cup \{0\} \right\}$$

is special set vector space over the set $S^o = Z^+ \cup \{0\}$. Find at least two distinct special set vector subspaces of V. What is the special dimension of V? Find a special generating subset of V.

27. Obtain some interesting properties about the special set linear algebra.

28. Does there exist a special set linear algebra which has no special set linear subalgebra?

29. Let $V = (V_1, V_2, V_3, V_4, V_5, V_6)$ where $V_1 = \{(1\ 1\ 1\ 1\ 1\ 1), (0\ 0\ 0\ 0\ 0\ 0), (1\ 1\ 0\ 1\ 1\ 0), (1\ 0\ 1\ 0\ 1\ 0), (1\ 1\ 1\ 1), (0\ 0\ 0\ 0), (1\ 1\ 1\ 0), (0\ 1\ 1\ 0)\}$,

$$V_2 = \left\{ \begin{pmatrix} a_1 & a_2 & a_3 & a_4 \\ a_5 & a_6 & a_7 & a_8 \end{pmatrix} \middle| a_i \in Z_2 = \{0,1\}\ 1 \le i \le 8 \right\},$$

$$V_3 = \left\{ \begin{pmatrix} a & b \\ c & d \end{pmatrix} \middle| a,b,c,d \in Z_2\{0,1\} \right\},$$

$$V_4 = \left\{ \begin{pmatrix} a_1 & a_2 \\ a_3 & a_4 \\ a_5 & a_6 \\ a_7 & a_8 \end{pmatrix} \middle| a_i \in Z_2 = \{0,1\} \quad 1 \le i \le 8 \right\},$$

$V_5$ = {all polynomials in the variable x of degree less than or equal to 8 with coefficients from $Z_2$} and $V_6 = \{Z_2 \times Z_2 \times Z_2 \times Z_2 \times Z_2 \times Z_2\}$ be a special set vector space over the set S = {0, 1}.

a. Find a special generating subset of V.
b. What is the special dimension of V?
c. Find atleast 5 distinct special set vector subspaces of V.
d. Define a special set linear operator on V.

30. Give some interesting properties about special set n-vector spaces defined over a set S.

31. Give an example of a special set n-linear algebra which is not a special set n-vector space defined over a same set S.

32. Let $V = V_1 \cup V_2 \cup V_3 \cup V_4$ where $V_1 = (V_1^1, V_2^1)$, $V_2 = (V_1^2, V_2^2, V_3^2)$, $V_3 = (V_1^3, V_2^3, V_3^3, V_4^3, V_5^3)$ and $V_4 = (V_1^4, V_2^4)$ defined over the set $S = 3Z^+ \cup \{0\}$. Take

$$V_1^1 = \left\{ \begin{pmatrix} a & b \\ c & d \end{pmatrix} \middle| a,b,c,d \text{ in } S \right\},$$

$V_2^1 = \{S \times S \times S\}$, so that $V_1 = (V_1^1, V_2^1)$ is a special set vector space over S. $V_1^2 = \{5 \times 5$ matrices with entries from S\}. $V_2^2 = \{S \times S \times S \times S \times S\}$, $V_3^2$ = {set of all $4 \times 2$ matrices with entries from S}. $V_4^2 = \{(a\ a\ a\ a), (a\ a\ a) \mid a \in S\}$ and $V_5^2$ = {all polynomials in the variable x with coefficients from S of degree less than or equal to 4}; $V_2 = \left(V_1^2, V_2^2, V_3^2, V_4^2, V_5^2\right)$ is a special set vector space over S.

Take in $V_3$, $V_1^3 = \{S \times S \times S \times S\}$, $V_2^3 = \{3 \times 3$ matrices with entries from S}, $V_3^3 = \{4 \times 4$ upper triangular matrices with entries from 5},

$$V_4^3 = \left\{ \begin{pmatrix} a & a & a & a \\ a & a & a & a \end{pmatrix}, \begin{pmatrix} a & a \\ a & a \\ a & a \end{pmatrix} \middle| a \in S \right\}$$

and $V_5^3 = \{S\ [x]$; all polynomials in x with coefficients from S of degree less than or equal to 5}, $V_3$ is a special set vector space over S. Finally define $V_1^4 = \{S \times S \times S\}$ and $V_2^4 = \{3 \times 3$ upper triangular matrices with entries from S}; $V_4$ is a special set vector space over S. Thus $V = V_1 \cup V_2 \cup V_3 \cup V_4$ is a special set vector 4-space over S.
  a. Find 3 special set vector 4-subspaces of V.
  b. Find special set generating 4-subset of V.
  c. What is special set 4-dimension of V over S?

33. Let $V = V_1 \cup V_2 \cup V_3 \cup V_4$ be a special set 4-vector space defined over the set $\{0, 1\} = S$, where $V_1 = \left(V_1^1, V_2^1\right)$, $V_1^1 = \{2 \times 2$ matrices with entries from S}, $V_2^1 = \{S \times S \times S \times S \times S\}$, $V_2 = \left(V_1^2, V_2^2, V_3^2\right)$ with

$$V_1^2 = \left\{ \begin{pmatrix} a & b \\ c & d \\ e & f \end{pmatrix} \middle| a,b,e,c,d,e,f \in Z^+ \cup \{0\} \right\},$$

$V_2^2 = \{S[x];$ polynomials of degree less than or equal to 7 over $S = \{0, 1\}\}$, $V_3^2 = \{7 \times 3$ matrices with entries from $S = \{0, 1\}\}$, $V_3 = \left(V_1^3, V_2^3, V_3^3, V_4^3\right)$ where

$$V_1^3 = \left\{ \begin{pmatrix} a \\ b \\ c \\ d \\ e \end{pmatrix} \middle| a, b, e, c, d, e \in S = \{0,1\} \right\},$$

$V_2^3 = \{$set of all $5 \times 5$ upper triangular matrices with entries from $S = \{0, 1\}\}$,

$$V_3^3 = \left\{ \begin{pmatrix} a & b \\ a & d \\ a & f \\ a & g \end{pmatrix} \middle| a, b, d, f, g \in S = \{0,1\} \text{ and } \begin{bmatrix} a & a & a & a \\ b & d & f & g \end{bmatrix} \right\},$$

$V_4^3 = \{S \times S \times S\}$ $V_4 = \left(V_1^4, V_2^4, V_3^4\right)$ where $V_1^4 = \{$set of all $4 \times 4$ lower triangular matrices with entries from the set $\{0, 1\}\}$, $V_2^4 = \{S \times S \times S \times S\}$ and

$$V_3^4 = \left\{ \begin{bmatrix} a & a & a & a & a \\ a & a & a & a & a \end{bmatrix}, \begin{pmatrix} a & a \\ a & a \\ a & a \\ a & a \end{pmatrix} \middle| a \in Z^+ \cup \{0\} \right\}.$$

a. Give 3 distinct special set fuzzy vector 4-spaces of V.
b. Find at least 2 proper special set vector 4-subspaces of V.

c. Find at least 4 proper special set fuzzy vector 4 subspaces of V.
d. Define a special set linear 4-operator on V.
e. Define a pseudo special set linear 4-operator on V.

34. Obtain some interesting results about special set fuzzy n-vector spaces.

35. For the special set linear bialgebra $V = V_1 \cup V_2$ given by $V_1 = (V_1^1, V_2^1, V_3^1, V_4^1)$ and $V_2 = (V_1^2, V_2^2, V_3^2, V_4^2, V_5^2)$ defined over the set $S = Z^+ \cup \{0\}$ where $V_1^1 = \{5 \times 5$ matrices with entries from S$\}$,

$$V_2^1 = \left\{ \begin{bmatrix} a & a & a & a \\ a & a & a & a \end{bmatrix} \middle| a \in S \right\},$$

$$V_3^1 = \left\{ \begin{bmatrix} a & a \\ a & a \\ a & a \end{bmatrix} \middle| a \in S \right\},$$

$V_4^1 = \{S[x]$ set of all polynomials of degree less than or equal to 5$\}$, $V_1^2 = \{S \times S \times S \times S\}$, $V_2^2 = \{$set of all $4 \times 4$ lower triangular matrices with entries from S$\}$, $V_3^2 = \{$all polynomials of degree less than or equal to 3$\}$,

$$V_4^2 = \left\{ \begin{bmatrix} a \\ a \\ a \\ a \\ a \end{bmatrix} \middle| a \in Z^o \right\}$$

and $V_5^2 = \{$set of all $5 \times 7$ matrices with entries from $Z^o\}$.
Obtain 3 different special set fuzzy linear bialgebras associated with V.

Find 5 different special set fuzzy linear subbialgebra associated with a proper set linear subalgebra $W = W_1 \cup W_2$.

36. Define a special set trialgebra and illustrate it by an example.

37. Let $V = V_1 \cup V_2 \cup V_3 = \left(V_1^1, V_2^1, V_3^1\right) \cup \left(V_1^2, V_2^2, V_3^2, V_4^2\right) \cup \left(V_1^3, V_2^3, V_3^3, V_4^3, V_5^3\right)$ where $V_1^1 = \{3 \times 3$ matrices with entries from $S = Z^o = Z^+ \cup \{0\}\}$, $V_2^1 = \{2 \times 6$ matrices with entries from $S = Z^o = Z^+ \cup \{0\}\}$, $V_3^1 = \{S \times S \times S \times S\}$, $V_1^2 = \{[a\ a\ a\ a\ a] \mid a \in S\}$,

$$V_2^2 = \left\{ \begin{bmatrix} a \\ a \\ a \\ a \\ a \end{bmatrix} \middle| a \in S \right\},$$

$V_3^2 = \{4 \times 4$ lower triangular matrices with entries from $S\}$,
$V_4^2 = \{7 \times 2$ matrices with entries from $S\}$, $V_1^3 = S \times S \times S$,

$$V_2^3 = \left\{ \begin{bmatrix} a & a & a & a & a \\ b & b & b & b & b \end{bmatrix} \middle| a, b \in S \right\},$$

$$V_3^3 = \left\{ \begin{bmatrix} a & d \\ b & e \\ c & f \end{bmatrix} \middle| a, b, c, d, e, f \in S \right\},$$

$$V_4^3 = \left\{ \begin{bmatrix} a_1 & a_2 & a_3 \\ a_4 & a_5 & a_6 \end{bmatrix} \middle| a_i \in S,\ 1 \le i \le 6 \right\}$$

and

$V_5^3 = \{[a\ a\ a\ a\ a\ a\ a] \mid a \in S\}$ be a special set linear trialgebra over $S = Z^o = Z^+ \cup \{0\}\}$.

Find a special set fuzzy linear trialgebra $V\eta$.

If W is a special set linear trisubalgebra of V find the special set fuzzy linear trisubalgebra $W\eta$.

Define a special set linear trioperator on V.

Find a pseudo special set trivector subspace of P of V and find its fuzzy component, $P\eta$.

Find a pseudo special set linear trioperator on V.

38. Let $V = V_1 \cup V_2 \cup V_3 \cup V_4$ be a special set vector space over the set $S = Z^+ \cup \{0\}$ where

$$V_1 = \left\{ \begin{bmatrix} a \\ a \\ a \end{bmatrix}, [a\ a\ a\ a] \middle| a \in S \right\},$$

$$V_2 = \left\{ \begin{bmatrix} a & b \\ c & d \end{bmatrix}, \begin{pmatrix} a & a & a \\ a & a & a \end{pmatrix} \middle| a,b,c,d, \in Z^o = S = Z^+ \cup \{0\} \right\},$$

$V_3 = \{$all polynomials of degree less than or equal to 4 with coefficient from $S\}$ and $V_4 = \{5 \times 5$ upper triangular matrices with entries from $S\}$.

Find $V\eta$ the special set fuzzy vector space.

Find a special set vector subspace W of V and find $W\eta$ the special set fuzzy vector space.

Find at least 4 distinct special set fuzzy vector space.

Find 3 distinct special set fuzzy vector subspaces of W.

39. Obtain some interesting results about special set fuzzy vector spaces.

40. Given $V = V_1 \cup V_2 \cup V_3$ is a special set trilinear algebra over the set $S = \{0, 1\}$ where $V_1 = \{S \times S \times S \times S \times S\}$, $V_2 = \{$All $5 \times 5$ matrices with entries from $S\}$ and $V_3 = \{3 \times 8$ matrices with entries from $S\}$.

a. Find Vη the special set fuzzy linear trialgebra.
b. Find a special set linear trisubalgebra W of V over S = {0, 1}, Find Wη the special set fuzzy linear trisubalgebra.
c. Find a special set pseudo trivector space W of V and find Wη.
d. Define a special set linear operator on V.
e. Find a special set pseudo linear operator on V.
f. Find atleast 3 special set linear subalgebra and their fuzzy analogue.

41. Obtain some interesting properties about the special semigroup set vector n-spaces.

42. Give an example of a special semigroup set vector 6-space.

43. Let $V = V_1 \cup V_2 \cup V_3 \cup V_4$ be a special semigroup set vector 4 space over the set $S = Z^+ \cup \{0\}$ where

$$V_1^1 = \left\{ \begin{pmatrix} a & b \\ c & d \end{pmatrix} \text{ such that a, b, c, d, } \in S \right\},$$

$V_2^1 = \{2 \times 6$ matrices with entries from S$\}$, $V_3^1 = S \times S \times S \times S \times S$ and $V_4^1 = \{7 \times 2$ matrices with entries from S, i.e., $V_1 = (V_1^1, V_2^1, V_3^1, V_4^1)$ and $V_2 = (V_1^2, V_2^2, V_3^2)$ where $V_1^2 = \{4 \times 4$ upper triangular matrices with entries from S$\}$ $V_2^2 = S \times S \times S \times S \times S$ and $V_3^2 = \{7 \times 3$ matrices with entries from S$\}$. $V_3 = (V_1^3, V_2^3)$ where $V_1^3 = \{$all polynomials in the variable x such that the coefficients are from S and every polynomial is of degree less than 5$\}$, $V_2^3 = \{4 \times 4$ lower triangular matrices with entries from S$\}$, $V_4 = (V_1^4, V_2^4, V_3^4, V_4^4, V_5^4)$ where $V_1^4 = \{3 \times 3$ upper triangular matrices with entries from S$\}$, $V_2^4 = S \times S \times S \times S \times S$,

$$V_3^4 = \left\{ \begin{bmatrix} a & a & a & a & a \\ b & b & b & b & b \\ c & c & c & c & c \end{bmatrix} \middle| a, b, c \in S \right\},$$

$$V_4^4 = \left\{ \begin{bmatrix} a & b \\ a & b \\ a & b \\ a & b \\ a & b \\ a & b \end{bmatrix} \middle| a, b \in S \right\}$$

and

$$V_5^4 = \{a_0 + a_1 x + a_2 x^2 + a_3 x^3 \mid a_0, a_1, a_2, a_3 \in S\}.$$

a. Find three proper special semigroup set vector 4-subspaces of V.
b. Find a proper generating special semigroup set of V.
c. Write V as a direct sum.
d. Find SHom (V, V).
e. Define a pseudo special semigroup set linear 4-operator on V.
f. Can V be made into a special set linear n-algebra?
g. Define on V a special semigroup set linear 4- operator P such that $P^2 = P$. Is P a projection on V?

44. Prove SHom(V, V) is a special group set vector n-space over a set S if V is a special semigroup set vector n-space over the set S.

45. Is PSHom(V, V) a special semigroup set vector n-space where PSHom(V, V) denotes the collection of all pseudo special semigroup set linear n-operator of V, V a special semigroup set vector space over S.
    a. What is the difference between SHom(V,V) and PSHom (V, V)?
    b. Which is of a larger n-dimension SHom(V,V) and PSHom (V, V)?

46. Let $V = V_1 \cup V_2 \cup V_3$ and $W = W_1 \cup W_2 \cup W_3$ be two special semigroup set vector trispaces over the set $\{0, 1\}$. Here
$$V = V_1 \cup V_2 \cup V_3$$
$$= \left(V_1^1, V_2^1\right) \cup \left(V_1^2, V_2^2, V_3^2\right) \cup \left(V_1^3, V_2^3, V_3^3\right)$$
and
$$W = W_1 \cup W_2 \cup W_3$$
$$= \left(W_1^1, W_2^1\right) \cup \left(W_1^2, W_2^2, W_3^2\right) \cup \left(W_1^3, W_2^3, W_3^3\right)$$
with $V_1^1 = \{S \times S \times S \times S \,/\, S = Z^+ \cup \{0\}\}$, $V_2^1 = \{\text{all } 3 \times 3$ matrices with entries from $2S\}$, $V_1^2 = \{S \times S \times S\}$, $V_2^2 = \{\text{all } 4 \times 4 \text{ upper triangular matrices with entries from } S\}$.

$$V_3^2 = \left\{ \begin{bmatrix} a \\ b \\ c \\ d \end{bmatrix} \,\middle|\, a, b, c, d \in 5S \right\},$$

$V_1^3 = S \times S \times S \times S \times S$, $V_2^3 = \{\text{all } 3 \times 3 \text{ lower triangular matrices with entries from } S\}$.

$$V_3^3 = \left\{ \begin{bmatrix} a \\ b \\ c \end{bmatrix} \,\middle|\, a, b, c \in S \right\}, \quad W_1^1 = \left\{ \begin{pmatrix} a & b \\ c & d \end{pmatrix} \,\middle|\, a, b, c, d \in S \right\},$$

$W_1^2 = \{\text{all polynomials in x of degree less than or equal to 8 with coefficients from } S\}$,

$$W_1^2 = \left\{ \begin{pmatrix} a & b \\ 0 & c \end{pmatrix} \,\middle|\, a, b, c \in S \right\},$$

$W_2^2 = \{5 \times 2 \text{ matrices with entries from } S\}$, $W_3^2 = S \times S \times S \times S$,

$$W_1^3 = \left\{ \begin{pmatrix} a & b & c \\ 0 & 0 & d \\ 0 & 0 & f \end{pmatrix} \middle| a,b,c,d,f \in S \right\},$$

$$W_2^3 = \left\{ \begin{pmatrix} a_1 & a_2 & a_3 \\ a_4 & a_5 & a_6 \end{pmatrix} \middle| a_i \in S \ 1 \leq i \leq 6 \right\}$$

and

$$W_3^3 = \left\{ (a \ a \ a) \middle| a \in S \right\}.$$

a. Find SHom(V, W): Is SHom(V, W) a special semigroup set vector trispace over S?
b. Find SHom(V, V).
c. Find SHom(W, W).
d. Does there exist any relation between SHom(V, V) and SHom(W, W)?
e. What is the tridimension of SHom(V, V) and SHom (W, W)?
f. Find two special semigroup set vector trisubspaces $V^1$ and $V^2$ of V.
g. Find two special semigroup set vector trisubspaces $W^1$ and $W^2$ of W.
h. Find SHom($W^1$, $V^1$).
i. Find SHom($W^2$, $V^2$).
j. Is SHom($W^1$, $V^1$) same as SHom($V^1$, $W^1$)? Justify your claim.
k. Find a special semigroup set linear projection of V into $V_1$.

47. Let $V = V_1 \cup V_2 \cup V_3 \cup V_4$ where
$$V_1 = \left( V_1^1, V_2^1, V_3^1 \right), V_2 = \left( V_1^2, V_2^2, V_3^2, V_4^2 \right),$$
$$V_3 = \left( V_1^3, V_2^3 \right) \text{ and } V_4 = \left( V_1^4, V_2^4, V_3^4, V_4^4 \right)$$
be a special semigroup set linear 4-algebra over the semigroup $Z_{12}$. Here

$$V_1^1 = \left\{ \begin{pmatrix} a & b \\ c & d \end{pmatrix} \middle| a,b,c,d \in Z_{12} \right\},$$

$V_2^1 = \{Z_{12} \times Z_{12} \times Z_{12}\}$, $V_3^1 = \{S \times S \times S \times S \times S \mid S = \{0, 2, 4, 6, 8, 10\} \subseteq Z_{15}\}$,

$$V_1^2 = \left\{ \begin{pmatrix} a & 0 \\ b & c \end{pmatrix} \middle| a,b,c \in Z_{12} \right\},$$

$V_2^2 = $ {all polynomials in the variable x with coefficients from $Z_{12}$ of degree less than or equal to 4}, $V_3^2 = \{Z_{12} \times Z_{12} \times Z_{12} \times Z_{12}\}$,

$$V_4^2 = \left\{ \begin{pmatrix} a \\ b \\ c \\ d \end{pmatrix} \text{ such that } a,b,c,d \in Z_{12} \right\},$$

$V_1^3 = \{Z_{12} \times S \times Z_{12} \times S \times Z_{12} \times S\}$ where $S = \{0, 2, 4, 6, 8, 10\}\}$,

$$V_2^3 = \left\{ \begin{pmatrix} a \\ b \\ c \\ d \\ e \\ f \end{pmatrix} \text{ such that } a,b,c,d,e,f \in S \right\},$$

$$V_1^4 = \{Z_{12} \times Z_{12}\}, \quad V_2^4 = \left\{ \begin{pmatrix} a \\ a \\ a \\ a \\ a \\ a \end{pmatrix} \middle| a \in Z_{12} \right\},$$

$V_3^4 = \{Z_{12}^3[x]\}$, i.e., only polynomials of degree less than or equal to 3 is taken) and $V_4^4 = \{2 \times 7$ matrices with entries from $Z_{12}\}$.

a. Find atleast 4 distinct special semigroup set linear 4 subalgebras of V.
b. Find SHom (V, V).
c. Find PSHom (V, V).
d. Find a special semigroup set 4-generator of V over $Z_{12}$.
e. Write V as a direct sum.
f. Based on (e) find a special semigroup set linear 4 projection of V.

48. Let $V = V_1 \cup V_2 \cup V_3 \cup V_4 \cup V_5$ be a special group set vector 5 space over the group $Z_{15}$ where
$V_1 = \left(V_1^1, V_2^1, V_3^1\right)$, $V_2 = \left(V_1^2, V_2^2, V_3^2, V_4^2\right)$, $V_3 = \left(V_1^3, V_2^3\right)$,
$V_4 = \left(V_1^4, V_2^4, V_3^4, V_4^4\right)$ and $V_5 = \left(V_1^5, V_2^5, V_3^5\right)$

with

$$V_1^1 = \left\{ \begin{pmatrix} a & b \\ c & d \end{pmatrix} \middle| a,b,c,d \in Z_{15} \right\},$$

$$V_1^2 = \left\{ \begin{pmatrix} a \\ b \\ c \\ d \end{pmatrix} \middle| a,b,c,d \in Z_{15} \right\},$$

$$V_3^1 = \{Z_{15} \times Z_{15} \times Z_{15}\},$$

$$V_1^2 = \left\{ \begin{pmatrix} a_1 & a_2 & a_3 \\ a_4 & a_5 & a_6 \end{pmatrix} \middle| a_i \in Z_{15}, 1 \le i \le 6 \right\},$$

$$V_2^2 = \left\{ \begin{pmatrix} a_1 & a_2 \\ a_3 & a_4 \\ a_5 & a_6 \\ a_7 & a_8 \end{pmatrix} \middle| a_i \in Z_{15}, 1 \leq i \leq 8 \right\},$$

$V_3^2 = \{$all $3 \times 3$ lower triangular matrices with entries from $Z_{15}\}$, $V_4^2 = \{(a\ a\ a\ a\ a\ a) \mid a \in Z_{15}\}$, $V_1^3 = \{Z_{15} \times Z_{15} \times Z_{15}\}$, $V_2^3 = \{$all $4 \times 4$ upper triangular matrices with entries from $Z_{15}\}$, $V_1^4 = \{Z_{15} \times Z_{15}\}$

$$V_2^4 = \left\{ \begin{pmatrix} a_1 & a_2 & a_3 & a_4 \\ a_5 & a_6 & a_7 & a_8 \end{pmatrix} \middle| a_i \in Z_{15}, 1 \leq i \leq 15 \right\},$$

$$V_3^4 = \left\{ \begin{pmatrix} a_1 & a_2 \\ a_3 & a_4 \\ a_5 & a_6 \end{pmatrix} \middle| a_i \in Z_{15}, 1 \leq i \leq 6 \right\},$$

$V_4^4 = \{$all $8 \times 2$ matrices with entries from $Z_{15}\}$; $V_1^5 = \{Z_{15} \times Z_{15} \times Z_{15} \times Z_{15} \times Z_{15} \times Z_{15}\}$, $V_2^5 = \{$set of all $2 \times 7$ matrices with entries from $Z_{15}\}$, and $V_3^5 = \{(a\ a\ a\ a) \mid a \in Z_{15}\}$. V is also a linear algebra.

a. Find for the special group set vector 5-space. The dimension of the special generating set. Find the special group set linear 5-algebra as a generating set and find its dimension.
b. Find a special group set vector 5 subspace of V.
c. Find a special group set linear 5-algebra of V over $Z_{15}$. Find SHom(V, V). Is it a special group set linear 5 algebra over $Z_{15}$? Find PSHom(V, V). Is PSHom(V, V) also a special group set linear 5-algebra over $Z_{15}$?
d. Define for V as a special group set linear 5-algebra a pseudo special subgroup vector 5-space of V.

49. Obtain any interesting property about special group set linear n-algebra defined over a group G.

50. Find for the special group set linear n-algebra $V = V_1 \cup V_2 \cup V_3$ where
    $$V_1 = \left(V_1^1, V_2^1, V_3^1\right), V_2 = \left(V_1^2, V_2^2\right) \text{ and } V_3 = \left(V_1^3, V_2^3, V_3^3\right)$$
    with $V_1^1 = \{Z_2 \times Z_2\}$; $V_2^1 = \{Z_2[x]$; all polynomials of degree less than or equal to 4$\}$, $V_3^1 = \{$all $5 \times 5$ matrices with entries from $Z_2\}$, $V_1^2 = \{Z_2 \times Z_2 \times Z_2\}$, $V_2^2 = \{$all $4 \times 4$ upper triangular matrices with entries from $Z_2\}$, $V_1^3 = \{Z_2 \times Z_2 \times Z_2 \times Z_2\}$, $V_2^3 = \{$all $3 \times 3$ matrices with entries from $Z_2\}$ and
    $$V_3^3 = \left\{ \begin{pmatrix} a & b & c & d \\ e & f & g & h \end{pmatrix} \middle| a,b,c,d,e,f,g,h \in Z_2 \right\}$$
    over the group $Z_2 = \{0, 1\}$ under addition modulo 2.
    a. Find 2 special group set linear n-subalgebra.
    b. Is V simple?
    c. Is V pseudo simple?
    d. Find SHom (V, V).
    e. Find PSHom (V, V).

51. Give some interesting properties about special semigroup set fuzzy linear algebra (vector space).

52. Let $V = V_1 \cup V_2 \cup V_3 \cup V_4$ where
    $$V_1 = \left(V_1^1, V_2^1, V_3^1\right), V_2 = \left(V_1^2, V_2^2\right)$$
    $$V_3 = \left(V_1^3, V_2^3, V_3^3, V_4^3\right) \text{ and } V_4 = \left(V_1^4, V_2^4\right)$$
    with $V_1^1 = \{Z \times Z\}$ $V_2^1 = \{$all $3 \times 3$ matrices with entries from $Z\}$,
    $$V_3^1 = \left\{ \begin{pmatrix} a \\ b \\ c \end{pmatrix} \middle| a,b,c \in Z \right\},$$

$V_1^2 = \{Z \times Z \times Z \times Z\}$, $V_2^2 = \{$all $5 \times 5$ upper triangular matrices with entries from Z$\}$, $V_1^3 = \{Z \times Z \times Z \times Z \times Z \times Z\}$

$$V_2^3 = \left\{ \begin{pmatrix} a & b \\ c & d \\ e & f \\ g & h \\ i & j \\ k & l \end{pmatrix} \middle| a,b,c....k,l \in Z \right\},$$

$V_3^3 = \{Z[x]$ all polynomials in the variable x with coefficients from Z of degree less than or equal to 4$\}$, $V_4^3 = \{$All lower triangular $4 \times 4$ matrices with entries from Z$\}$. $V_1^4 = Z \times Z \times Z$ and

$$V_2^4 = \left\{ \begin{pmatrix} a & b \\ a & b \\ a & b \\ a & b \end{pmatrix} \middle| a,b \in Z \right\}$$

be a special group set linear 4-algebra over the group Z.
a. Define two distinct special group set fuzzy linear 4-algebras.
b. Find three proper distinct special group set fuzzy linear 4-subalgebras for three distinct special group set linear 4-subalgebras W, U and S.
c. Determine a special group set fuzzy linear 4-subalgebras η which can be extended to V?
d. Can all special group set fuzzy linear 4-subalgebras Wη be extended to Vη? Justify your claim.

53. Obtain some interesting properties about special semigroup set n-vector spaces (linear algebras).

54. Prove every special group set vector n space in general is not a special group set linear n-algebra and every special group set n-linear algebra is a special group set n-vector space.

55. Show every special group set vector n-space over a set is always a special semigroup set vector n-space and vice versa!

56. Give examples of special semigroup set vector n-spaces which are not special group set vector n-spaces.

57. Is every special group set vector n-space special semigroup set vector n space? Justify your claim.

58. Let V and W be two special group set vector 4-space defined over the same set S. Find SHom (V, W). Find a generating 4-set of SHom (V, W).

59. Let
$$V = V_1 \cup V_2 \cup V_3 \cup V_4 \cup V_5$$
$$= \left(V_1^1, V_2^1, V_3^1\right) \cup \left(V_1^2, V_2^2\right) \cup \left(V_1^3, V_2^3, V_3^3\right)$$
$$\cup \left(V_1^4, V_2^4\right) \cup \left(V_1^5, V_2^5, V_3^5, V_4^5\right)$$

be a special group set vector 5-space over the set $S = Z$, where $V_1^1 = Z \times Z$, $V_2^1 = $ {all $3 \times 3$ matrices with entries from $Z$},

$$V_3^1 = \left\{ \begin{pmatrix} a \\ b \\ c \end{pmatrix} \middle| a, b, c \in Z \right\},$$

$V_1^2 = Z \times Z \times Z \times Z$, $V_2^2 = $ {all $2 \times 2$ upper triangular matrices with entries from $Z$}, $V_1^3 = $ {all $4 \times 4$ lower triangular matrices with entries from $Z$},

$$V_2^3 = \left\{ \begin{pmatrix} a & a & a & a & a & a \\ b & b & b & b & b & b \end{pmatrix} \middle| a, b \in Z \right\},$$

$$V_3^3 = \left\{ \begin{pmatrix} a & a \\ a & a \\ a & a \\ a & a \\ b & b \\ b & b \\ b & b \end{pmatrix} \middle| a, b \in Z \right\},$$

$V_1^4 = Z \times Z \times Z \times Z \times Z$, $V_2^4 = \{$all $3 \times 3$ lower triangular matrices with entries from $Z\}$, $V_1^5 = \{(a\ a\ a\ a\ a) \mid a \in Z\}$,

$$V_2^5 = \left\{ \begin{pmatrix} a \\ b \\ c \\ d \\ e \\ f \\ g \end{pmatrix} \middle| a, b, \ldots g \in Z \right\}.$$

$V_3^5 = \{$all $2 \times 2$ matrices with entries from $Z\}$ and $V_4^5 = \{Z[x]$ all polynomials of degree less than or equal to 10$\}$.

a. Find the special group set 5-dimension of V.
b. Find a special group set generating set of V.
c. Find SHom (V, V). Is SHom (V, V) finite dimensional.
d. Find at least 3 proper distinct special group set vector 5 subspaces of V. W, U and T.
e. Define $\eta: V \to [0, 1]$ so that $V\eta$ is a special group set fuzzy vector 5-space.

f.  Will Wη, Uη, Tη be special group set fuzzy vector 5 subspaces where η on V is restricted to Wη, Uη and Tη, be special group set fuzzy vector 5-subspaces where η on V is restricted to W, U and T. Similarly if Wη$_\omega$ Uη$_u$ and Tη$_T$ be three distinct special group set fuzzy vector 5-subspaces of V; can their extensions to V be such that Wη'$_\omega$, Uη'$_u$ and Tη'$_T$ where η'$_\omega$, η'$_u$ and η'$_T$ extended to V be a special group set fuzzy vector 5 spaces.

g.  Find T ∈ SHom(V, V) so that T is a special group set 5-projecton on V.

h.  Write V = V$_1$ ∪ V$_2$ ∪ ... ∪ V$_5$ as a special group set direct 5-sum of V.

i.  Suppose the set Z is replaced by the set S = {0, 1, 2, ..., 10} ⊆ Z. Will V = V$_1$ ∪ V$_2$ ∪ ... ∪ V$_5$ be a special group set vector 5 space over S?

j.  What is (V's) special dimension when V is defined over S = {0, 1, 2, ..., 10}?

k.  Will SHom (V, V) as a special group set 5-vector space over S be finite or infinite? Justify your claim.

l.  Prove by illustrative example V has more than one special group set direct sum as special group set vector 5-subspaces?

m.  Show if T: V → V, U: V → V both T o U and U o T are defined.

n.  Find T and U from V into V i.e., in SHom (V, V) such that
    i.   T o U = U o T
    ii.  T o U ≠ U o T
    iii. T o T = T and U o U = U
    iv.  U o T = T and T o U = U
    v.   T o T = (0).

o.  Find η: V → [0, 1] so that Vη is a special group set fuzzy vector 5-space.

p.  Find for W ⊆ V two maps η$_1$ and η$_2$ (η$_1$≠η$_2$) such that Wη$_1$ and Wη$_2$ are two distinct special group set fuzzy vector 5-subspaces such that η$_1$ ∩ η$_2$: W → [0, 1] is an empty intersection on the subset of [0, 1].

60. Find some applications of these new structures to industrial problems.

61. Prove if $V = V_1 \cup V_2 \cup ... \cup V_n$ is a special group set vector space over a set S by changing the set S to another set T the space can be changed from infinite (finite) special group set n-dimension to finite (infinite) special group set n-dimension.

62. Give examples of simple special group set vector n-spaces.

63. Give examples of pseudo simple special group linear n-algebras. Hence or otherwise show the pseudo simple nature is dependent on the group over which the special group linear n-algebras are defined.

64. Is it true all simple special group set linear n-algebras are pseudo simple special group set linear n-algebras? Are these two notions entirely different / distinct?

65. Let $V = V_1 \cup V_2 = \left(V_1^1, V_2^1, V_3^1\right) \cup \left(V_1^2, V_2^2, V_3^2, V_4^2\right)$ where

$$V_1^1 = \{(a\ a\ a\ a) \mid a \in Z\},$$

$$V_2^1 = \left\{ \begin{pmatrix} a \\ a \\ a \\ a \\ a \end{pmatrix} \middle| a \in Z \right\},$$

$$V_3^1 = \left\{ \begin{pmatrix} a & a & a \\ a & a & a \\ a & a & a \end{pmatrix} \middle| a \in Z \right\},$$

$$V_1^2 = \left\{ \begin{pmatrix} a & a \\ a & a \end{pmatrix} \middle| a \in Z \right\},$$

$$V_2^2 = \left\{ \begin{pmatrix} a & a \\ a & a \\ a & a \\ a & a \end{pmatrix} \middle| a \in Z \right\},$$

$$V_3^2 = \{(a\ a\ a\ a\ a\ a) \mid a \in Z\}$$

and

$$V_4^2 = \left\{ \begin{pmatrix} a & 0 & 0 & 0 \\ 0 & a & 0 & 0 \\ 0 & 0 & a & 0 \\ 0 & 0 & 0 & a \end{pmatrix} \middle| a \in Z \right\}$$

be a special group set linear the group bialgebra over Z.
a. Is V simple special group set linear bialgebra?
b. Is V a pseudo simple special group set linear bialgebra.
c. Find SHom (V, V).
d. Is SHom (V, V) a simple special group set linear bialgebra?

66. Let $V = V_1 \cup V_2 \cup V_3 = \left( V_1^1, V_2^1, V_3^1, V_4^1 \right) \cup \left( V_1^2, V_2^2, V_3^2 \right)$ $\cup \left( V_1^3, V_2^3 \right)$ be a special group set linear trialgebra over the group $Z_{24}$, where $V_1^1 = \{Z_{24} \times Z_{24} \times Z_{24}\}$, $V_2^1 = \{$(set of all 2 × 2 matrices with entries from $Z_{24}\}$,

$$V_3^1 = \left\{ \begin{pmatrix} a \\ b \\ c \\ d \end{pmatrix} \middle| a, b, c, d \in Z_{24} \right\},$$

$V_4^1 = \{Z_{24}[x] /$ polynomials of degree less than or equal to 6 with entries from $Z_{24}\}$, $V_1^2 = \{Z_{24} \times Z_{24}\}$, $V_2^2 = \{$all $3 \times 3$ lower triangular matrices with entries from $Z_{24}\}$ and

$$V_3^2 = \left\{ \left. \begin{pmatrix} a & b \\ c & d \\ e & f \\ g & h \end{pmatrix} \right| a, b....h \in Z_{24} \right\},$$

$V_1^3 = \{Z_{24} \times Z_{24} \times Z_{24} \times Z_{24}\}$ and $V_2^3 = \{$all $5 \times 5$ upper triangular matrices with entries from $Z_{24}\}$.

a. Find two distinct special group set linear trisubalgebras W and U of V.
b. Find $V_\eta$, $W_\eta$ and $U_\eta$ the special group set fuzzy linear trialgebra and special group set fuzzy linear trisubalgebra.
c. Show SHom(V, V) is also a special group set linear trialgebra over $Z_{24}$.
d. Is V a simple special group set linear trialgebra?
e. Is V a pseudo simple special group set linear trialgebra.
f. If $Z_{24}$ is replaced by $S = \{0, 1\} = Z_2$ will V be a pseudo simple special group set linear trialgebra?
g. Will V as a special group set vector trispace over $Z_2$ be a simple special group set vector trispace?

# FURTHER READING

# INDEX

# ABOUT THE AUTHORS

**Dr.W.B.Vasantha Kandasamy** is an Associate Professor in the Department of Mathematics, Indian Institute of Technology Madras, Chennai. In the past decade she has guided 12 Ph.D. scholars in the different fields of non-associative algebras, algebraic coding theory, transportation theory, fuzzy groups, and applications of fuzzy theory of the problems faced in chemical industries and cement industries.

She has to her credit 646 research papers. She has guided over 68 M.Sc. and M.Tech. projects. She has worked in collaboration projects with the Indian Space Research Organization and with the Tamil Nadu State AIDS Control Society. This is her $44^{th}$ book.

On India's 60th Independence Day, Dr.Vasantha was conferred the Kalpana Chawla Award for Courage and Daring Enterprise by the State Government of Tamil Nadu in recognition of her sustained fight for social justice in the Indian Institute of Technology (IIT) Madras and for her contribution to mathematics. (The award, instituted in the memory of Indian-American astronaut Kalpana Chawla who died aboard Space Shuttle Columbia). The award carried a cash prize of five lakh rupees (the highest prize-money for any Indian award) and a gold medal.
She can be contacted at vasanthakandasamy@gmail.com
Web Site: http://mat.iitm.ac.in/home/wbv/public_html/

**Dr. Florentin Smarandache** is a Professor of Mathematics at the University of New Mexico in USA. He published over 75 books and 150 articles and notes in mathematics, physics, philosophy, psychology, rebus, literature.

In mathematics his research is in number theory, non-Euclidean geometry, synthetic geometry, algebraic structures, statistics, neutrosophic logic and set (generalizations of fuzzy logic and set respectively), neutrosophic probability (generalization of classical and imprecise probability). Also, small contributions to nuclear and particle physics, information fusion, neutrosophy (a generalization of dialectics), law of sensations and stimuli, etc. He can be contacted at smarand@unm.edu

**K. Ilanthenral** is the editor of The Maths Tiger, Quarterly Journal of Maths. She can be contacted at ilanthenral@gmail.com